\documentclass[reqno]{amsart}
\input epsf
\address{Simons Center for Geometry and Physics,
State University of New York, Stony Brook, NY 11794-3636 U.S.A. \&
Center for Geometry and Physics, Institute for Basic Sciences (IBS), Pohang, Korea} \email{kfukaya@scgp.stonybrook.edu}
\address{Center for Geometry and Physics, Institute for Basic Sciences (IBS), Pohang, Korea \& Department of Mathematics,
POSTECH, Pohang, Korea} \email{yongoh1@postech.ac.kr}
\address{Graduate School of Mathematics,
Nagoya University, Nagoya, Japan} \email{ohta@math.nagoya-u.ac.jp}
\address{Research Institute for Mathematical Sciences, Kyoto University, Kyoto, Japan}
\email{ono@kurims.kyoto-u.ac.jp}

\usepackage{graphicx}
\usepackage{amsmath}
\usepackage{amscd}
\usepackage{amssymb}
\usepackage{amstext}
\usepackage{amsmath}
\usepackage[all]{xy}
\usepackage{yhmath}
\usepackage{mathrsfs}
\usepackage{bbding}
\usepackage{color}
\usepackage{hyperref}
\hypersetup{
 setpagesize=false,
 bookmarksnumbered=true,%
 bookmarksopen=true,%
 colorlinks=true,%
 linkcolor=blue,
 citecolor=green,
}
\makeatletter
\def\l@section{\@tocline{1}{0pt}{3mm}{8mm}{}}
\def\l@subsection{\@tocline{2}{0pt}{6mm}{10mm}{}}
\def\l@subsubsection{\@tocline{3}{0pt}{9mm}{11mm}{}}
\makeatother

\def\E{\ifmmode{\mathbb E}\else{$\mathbb E$}\fi} 
\def\N{\ifmmode{\mathbb N}\else{$\mathbb N$}\fi} 
\def\R{\ifmmode{\mathbb R}\else{$\mathbb R$}\fi} 
\def\Q{\ifmmode{\mathbb Q}\else{$\mathbb Q$}\fi} 
\def\C{\ifmmode{\mathbb C}\else{$\mathbb C$}\fi} 
\def\H{\ifmmode{\mathbb H}\else{$\mathbb H$}\fi} 
\def\Z{\ifmmode{\mathbb Z}\else{$\mathbb Z$}\fi} 
\def\P{\ifmmode{\mathbb P}\else{$\mathbb P$}\fi} 
\def\T{\ifmmode{\mathbb T}\else{$\mathbb T$}\fi} 
\def\SS{\ifmmode{\mathbb S}\else{$\mathbb S$}\fi} 
\def\DD{\ifmmode{\mathbb D}\else{$\mathbb D$}\fi} 
\def\K{\ifmmode{\mathbb K}\else{$\mathbb K$}\fi}

\newcommand{\del}{\partial}

\theoremstyle{theorem}
\newtheorem{thm}{Theorem}[section]
\newtheorem{cor}[thm]{Corollary}
\newtheorem{lem}[thm]{Lemma}
\newtheorem{sublem}[thm]{Sublemma}

\newtheorem{prop}[thm]{Proposition}

\newtheorem{lemdef}[thm]{Lemma-Definition}

\theoremstyle{definition}
\newtheorem{defn}[thm]{Definition}
\newtheorem{rem}[thm]{Remark}

\newtheorem{exm}[thm]{Example}
\newtheorem{conds}[thm]{Condition}
\newtheorem{conven}[thm]{Convention}
\newtheorem{proper}[thm]{Property}
\newtheorem{defnlem}[thm]{Definition-Lemma}
\newtheorem{assump}[thm]{Assumption}

\newtheorem{shitu}[thm]{Situation}

\newtheorem{notation}[thm]{Notation}

\numberwithin{equation}{section}
\makeindex
\begin{document}

\title[Kuranishi structure
and virtual fundamental chain]
{Kuranishi structure,
Pseudo-holomorphic curve,
and
virtual fundamental chain: Part 2}
\author{Kenji Fukaya, Yong-Geun Oh, Hiroshi Ohta, Kaoru Ono}

\thanks{Kenji Fukaya is supported partially by JSPS Grant-in-Aid for Scientific Research
No. 23224002, NSF Grant No. 1406423 and
Simons Collaboration on Homological Mirror Symmetry,
Yong-Geun Oh by the IBS project IBS-R003-D1, Hiroshi Ohta by JSPS Grant-in-Aid
for Scientific Research Nos. 23340015, 15H02054 and Kaoru Ono by JSPS Grant-in-Aid for
Scientific Research, Nos. 26247006, 23224001.}

\begin{abstract}
This article consists the second parts of the article we promised at the
end of \cite[Section 1]{foootech}.
We discuss the foundation of the virtual fundamental
chain and cycle technique, especially
its version appeared in \cite{FO} and also in
\cite[Section A1, Section 7.5]{fooobook2}, \cite[Section 12]{fooo09},
\cite{fooo091}.
\par
This article is independent of our
earlier writing \cite{foootech}.
We also do not assume that the readers have any knowledge on
the pseudo-holomorphic curve.
\par
In this second part, we consider a system of spaces with Kuranishi structures
(abbreviated as K-system)
and its simultaneous perturbations.
\end{abstract}
\maketitle
\par\newpage
\date{March 31th, 2017}


\maketitle

\tableofcontents


\part{System of K-spaces and smooth correspondences}
\label{part2}


\section{Introduction to Part 2}
\label{subsec:intro2}
In Part 1, we have described the foundation of the theory of Kuranishi structure, good coordinate system,
CF-perturbation (also multivalued perturbation), and defined
the integration along the fiber (push out) of a strongly submersive map
with respect to a CF-perturbation and also prove Stokes' formula.
Using these ingredients,
we have established the notion of smooth correspondence and
proved the composition formula.
Especially, this provides us a way to obtain a virtual fundamental chain
for each K-space (space with Kuranishi structure).
Thus Part 1 is the story for each {\it single K-space}.
On the other hand, in Part 2 we are going to study a {\it system of K-spaces}.
In actual geometric applications,
there are the cases for which
it is not enough to study each single K-space individually, but is necessary to study
a system of K-spaces satisfying certain {\it compatibility
conditions}, especially, {\it compatibility conditions at boundary and corner}.
The compatibility conditions we describe depend on the situation we consider in their detail.
Here we have two geometric examples in mind. One is the Floer cohomology for
periodic Hamiltonian system which is established by \cite{FO}, \cite{LiuTi98}
for general closed symplectic manifold,
and the other is the $A_{\infty}$ algebra
associated to a Lagrangian submanifold and the Floer cohomology for the
Lagrangian intersection established by \cite{fooobook} \cite{fooobook2}.
Since this article intends to provide a `package' of the statements
appearing in the actual argument above,
we begin with axiomatizing the properties and conditions in a purely abstract setting,
motivated by these two geometric examples.
In this way we discuss two kinds of systems of K-spaces in Part 2:
One is a {\it linear K-system} containing the Floer theory for
periodic Hamiltonian system as a typical example,
and the other is a {\it tree-like K-system} containing the theory of the $A_{\infty}$ algebra
associated to a Lagrangian submanifold.
\par
We emphasize that we discuss those two cases as a prototype of the
applications of the results of this article.
In fact if we are interested only in defining
Floer cohomology of periodic Hamiltonian system and proving its
basic properties,
certainly there is a shorter proof than those given in
Sections \ref{sec:systemline1}-\ref{sec:systemline4}.
In this article we give a general proof which can be used in similar situations with
minimal change.
For this reason we tried to avoid using special feature of the
particular situation we work with but use only the arguments which are
general enough.
A similar argument works in most of the other cases
where the method of pseudo-holomorphic curves is applied.
(We do not know the case this method does not work.)
We confirm that it also works at least in the following situations.
\begin{enumerate}
\item
Constructing and proving basic properties of the Gromov-Witten invariant.
\item
Studying several Lagrangian submanifolds and constructing
an $A_{\infty}$ category (Fukaya catgory).
\item
The family version of (2) and equivariant versions of (1)(2).
\item
Including immersed Lagrangian submanifolds.
\item
Using the Lagrangian correspondence to construct
an $A_{\infty}$ functor.
\item
Including bulk deformations into the Lagrangian Floer theory to
define open-closed, closed-open maps, and  proving their
basic properties.
\item
Studying the moduli space of psuedo-holomorphic maps
from a bordered Riemann surface of arbitrary genus with Lagrangian
boundary condition to construct an IBL-infinity structure.
\item
Including the non-compact case in defining and studying
symplectic homology and wrapped Floer homology.
\item
In the case of symplectic manifolds with contact type boundary,
using closed Reeb orbits or Reeb chords to establish the
foundation of symplectic field theory and
its version with Legendrian submanifolds.
\end{enumerate}
\par
The purpose of Part 2 is summarized as follows:
If we are given a system of K-spaces satisfying the axiom we describe in this article
as an input, then
we prove that we can derive certain algebraic structures
from the system of K-spaces as an output.
This is a `package' producing an algebraic structure from a geometric input.
The problem is that the resulting algebraic structure itself
depends on the various choices made in the course of the construction, in general.
We will specify in which sense the algebraic structure is
invariant and prove the invariance in this article.
Although these typical examples of K-systems arise from moduli spaces of
pseudo-holomorphic curves, the authors expect that this kind of axiomatization
and framework will
be available for other problems arising from other situations in future.
\par
\subsection{Outline of the story of linear K-system}
\label{subsecintolinear}
In Sections \ref{sec:systemline1}-\ref{sec:systemline4}
we study systems of K-spaces, which axiomatize
the situation appearing during the construction of
Floer cohomology of periodic Hamiltonian systems.

\subsubsection{Floer cohomology of periodic Hamiltonian systems}
We first review the outline of the construction
of Floer cohomology of periodic Hamiltonian systems in  \cite{Flo89I}.
\par
Let
$H : S^1 \times M \to \R$ be a real valued smooth
function on the product of $S^1$ and a symplectic manifold $M$.
For each $t \in S^1$ we obtain a function
$H_t : M \to \R$ by
$$
H_t(x) = H(t,x).
$$
We denote its Hamiltonian vector field by $X_{H_t}$
defined by $dH_t=\omega(X_{H_t}, \cdot)$.
A periodic solution of the periodic Hamiltonian system
generated by $H$ is by definition a smooth map
$\ell : S^1 \to M$ that satisfies the equation:
\begin{equation}\label{Hameq}
\frac{d\ell}{dt} = X_{H_t} \circ \ell.
\end{equation}
Let ${\rm Per}(H)$ be the set of all solutions of (\ref{Hameq}).
We will work in the Bott-Morse situation so
we assume that ${\rm Per}(H)$ is a smooth manifold
satisfying certain nondegeneracy condition. (See \cite[Subsection 1.2 (2)]{fooospectr}
for example.)

\begin{rem}
To define Floer cohomology,
it is enough to discuss the Morse case, that is, the case when
$H$ satisfies certain nondegeneracy condition so that
${\rm Per}(H)$ is discrete.
However to calculate
Floer cohomology, it is useful to include
the Bott-Morse case especially the case with $H\equiv 0$.
\end{rem}

We decompose ${\rm Per}(H)$ into connected components and denote
$$
{\rm Per}(H) = \bigcup_{\overline\alpha \in \overline{\frak A}} \overline R_{\overline\alpha}.
$$
For $\overline\alpha_-,\overline\alpha_+ \in \overline{\frak A}$ we consider the set of solutions of the
following equation (Floer's equation) \index{Floer's equation}
for a map $u : \R \times S^1 \to M$
\begin{equation}\label{Floerseq}
\frac{\partial u}{\partial \tau}
+
J\left(
\frac{\partial u}{\partial t} - X_{H_t}(u)
\right)
=0
\end{equation}
together with the following asymptotic boundary condition.

\begin{conds}\label{asboundary}
There exist $\gamma_{-\infty} \in \overline R_{\overline\alpha_-}$ and
$\gamma_{+\infty} \in \overline R_{\overline\alpha_+}$
such that
\begin{equation}\label{evalatinfty}
\aligned
\lim_{\tau \to -\infty}   u(\tau,t) &= \gamma_{-\infty} (t),\\
\lim_{\tau \to +\infty}  u(\tau,t) &= \gamma_{+\infty} (t).
\endaligned
\end{equation}
\end{conds}
We denote by
$\overset{\circ}{\widetilde{\mathcal M}}(\overline R_{\overline\alpha_-},\overline R_{\overline\alpha_+})$
the set of solutions of (\ref{Floerseq}) satisfying
Condition \ref{asboundary}.
We can define an $\R$ action on $\overset{\circ}{\widetilde{\mathcal M}}(\overline R_{\overline\alpha_-},\overline R_{\overline\alpha_+})$
by
$(\tau_0 \cdot u)(\tau,t) = u(\tau+\tau_0,t)$,
and denote the associated quotient space by
$\overset{\circ}{{\mathcal M}}(\overline R_{\overline\alpha_-},\overline R_{\overline\alpha_+})$.
We decompose it into
$$
\overset{\circ}{\mathcal M}(\overline R_{\overline\alpha_-},\overline R_{\overline\alpha_+})
=
\bigcup_{\beta}\overset{\circ}{\mathcal M}(\overline R_{\overline\alpha_-},\overline R_{\overline\alpha_+};\beta),
$$
according to the homology class $\beta$ of $u$.
In place of this decomposition we proceed as follows.
We denote by $\frak A$ the set of pairs $(\overline\alpha,[w])$
where $\overline\alpha \in \overline{\frak A}$ and
$[w]$ is the homology class of a disc map $w$ bounding $\ell
\in \overline R_{\overline\alpha}$.\footnote{The equivalence class is defined by using
the symplectic area and the Maslov index.}
For each $\alpha \in \frak A$ we define $R_{\alpha}$ as the
set of pairs consisting of an element $\ell$ of
$\overline R_{\overline\alpha}$ and an equivalence class $[w]$ of
the homology class of disc $w$ bounding $\ell$.

Then let
$
\overset{\circ}{\mathcal M}(R_{\alpha_-},R_{\alpha_+})$
be the union of
$\overset{\circ}{\mathcal M}(\overline R_{\overline\alpha_-},\overline R_{\overline\alpha_+};\beta)$
over $\beta$
with $[w^-] \# [\beta] = [w^+]$
(where $\alpha_{\pm} = (\overline\alpha_{\pm},[w^{\pm}]$)
and denote the union by ${\mathcal M}^{\text{\rm reg}}(H;\alpha_-,\alpha_+)$.
Using the notion of stable map, we can compactify each of
${\mathcal M}^{\text{\rm reg}}(H;\alpha_-,\alpha_+)$ and
denote the compactification by
${\mathcal M}(H;\alpha_-,\alpha_+)$.
(See \cite[Definition 19.9]{FO}.)
\par
We define asymptotic evaluation maps
$
{\mathcal M}^{\text{\rm reg}}(H;\alpha_-,\alpha_+)
\to \overline R_{\alpha_{\pm}}
$ by
$$
\aligned
{\rm ev}_{-}([u]) = \gamma_{-\infty},
\qquad
{\rm ev}_{+}([u]) = \gamma_{+\infty},
\endaligned
$$
where $\gamma_{-\infty}$, $\gamma_{+\infty}$ are as in
(\ref{evalatinfty}).
They induce maps
$$
{\rm ev}_{\pm} :
{\mathcal M}(H;\alpha_-,\alpha_+)
\to R_{\alpha_{\pm}},
$$
which we call the {\it evaluation maps at infinity}.
Using the non degeneracy condition of
$R_{\alpha}$, we can show that
${\mathcal M}(H;\alpha_-,\alpha_+)$ carries a Kuranishi structure with corners
and the evaluation map
\begin{equation}\label{form117}
({\rm ev}_{-},{\rm ev}_{+})  :
{\mathcal M}(H;\alpha_-,\alpha_+)
\to
R_{\alpha_-} \times R_{\alpha_+}
\end{equation}
is a  strongly smooth and weakly submersive map.
Moreover we can find the following
isomorphism of K-spaces:
\begin{equation}\label{form1555}
\partial{\mathcal M}(H;\alpha_-,\alpha_+)
=
\bigcup_{\alpha}
{\mathcal M}(H;\alpha_-,\alpha)
\,\,{}_{{\rm ev}_{+}}\times_{{\rm ev}_{-}}
{\mathcal M}(H;\alpha,\alpha_+).
\end{equation}
When we regard the left hand side as the
{\it normalized} boundary, then
the right hand side becomes the {\it disjoint} union.
\par
See \cite[Part 5]{foootech} for the construction of such a K-system.

\subsubsection{Periodic Hamiltonian system and axiom of linear K-system}
As we explained in the previous subsubsection,
when we are given a time dependent Hamiltonian $H$,
we obtain a system
consisting of a set of smooth manifolds $R_{\alpha}$
and a set of K-spaces
${\mathcal M}(H;\alpha_-,\alpha_+)$,
together with evaluation maps
(\ref{form117}).
The axiom of linear K-systems, which we present in
Section \ref{sec:systemline1} Condition \ref{linsysmainconds},
spells out properties of such a system
which we need to define Floer cohomology.
\par
Condition \ref{linsysmainconds} (III)(IV) require the existence
of a set of manifolds $\{R_{\alpha}
\mid \alpha \in \frak A\}$,
a set of K-spaces
$\{{\mathcal M}(H;\alpha_-,\alpha_+)
\mid \alpha_-,\alpha_+ \in \frak A\}$
and evaluation maps (\ref{form117}) indexed by a countable set $\frak A$.
In this abstract situation
we call ${\mathcal M}(H;\alpha_-,\alpha_+)$
the {\it space of connecting orbits}.
\par
Condition \ref{linsysmainconds} (VI) (and (I))
requires that we can associate the Maslov index
$\mu(\alpha)$ to each $R_{\alpha}$
which determines the dimension of
${\mathcal M}(H;\alpha_-,\alpha_+)$.
\par
It is well-known that the energy
\begin{equation}\label{formulaenergy}
\int_{\R \times S^1}
\left\Vert\frac{\partial u}{\partial \tau}\right\Vert^2
+
\left\Vert
\frac{\partial u}{\partial t} - X_{H_t}(u)
\right\Vert^2 d\tau dt
\end{equation}
of the solution $u$ of (\ref{Floerseq})
is a difference of the value of certain action functional
at the asymptotic boundary values
$\alpha_-$, $\alpha_+$.
Moreover the energy is nonnegative and is zero
only when $\partial u/\partial \tau = 0$.
Condition \ref{linsysmainconds} (V) (and (I))
is an axiomatization of this property.
\par
We note that, in our Bott-Morse situation,
we need to introduce an appropriate $O(1)$-principal
bundle $o_{R_{\alpha}}$ to our critical submanifold to define Floer cohomology.
Namely the contribution of $R_{\alpha}$
to the Floer cohomology is the cohomology group of $R_{\alpha}$
with the coefficients
twisted by this local system.
(See \cite[Subsection 8.8]{fooobook2}.)  Then the orientation local system
of the moduli space  ${\mathcal M}(H;\alpha_-,\alpha_+)$
is related to the orientations of $R_{\alpha_{\pm}}$ and
to  $o_{R_{\alpha_{\pm}}}$ in an appropriate way.
Condition \ref{linsysmainconds} (VII)
is an axiomatization of this property.
\par
Note that an element of $R_{\alpha}$ is a pair $(\ell,[w])$ where
$\ell$ is a periodic orbit of our Hamiltonian system
and $[w]$ is a homology class of disks bounding $\ell$.
For an element $\beta \in H_2(M;\Z)$
represented by a sphere, we can glue $w$ with a
representative of $\beta$ to obtain another disk.
By this operation we obtain another
$R_{\beta \#\alpha}$. It is easy to see
${\mathcal M}(H;\alpha_-,\alpha_+)
\cong {\mathcal M}(H;\beta \#\alpha_-,\beta \#\alpha_+)$.
Condition \ref{linsysmainconds} (VIII)
is an axiomatization of this property.
\par
One important property of the moduli space
of pseudo-holomorphic curve or
the solution of Floer's equation,
is Gromov compactness. It claims the compactness of the
union of the moduli spaces whose elements have energy smaller
than a fixed number.
Condition \ref{linsysmainconds} (IX)
is an axiomatization of this property.
\par
The boundary of our moudli space
${\mathcal M}(H;\alpha_-,\alpha_+)$ is described as (\ref{form1555}).
Moreover this isomorphism is not only one as topological spaces
but also one as spaces with oriented Kuranishi structure.
Condition \ref{linsysmainconds} (X)
is an axiomatization of this property.
\par
Our moduli space ${\mathcal M}(H;\alpha_-,\alpha_+)$
has not only boundary but also corners in general.
Its codimension $k$ (normalized) corner
$\widehat S_k{\mathcal M}(H;\alpha_-,\alpha_+)$
is described as the disjoint union of the fiber products
$$
\aligned
{\mathcal M}(H;\alpha_-,\alpha_1)
&\times_{R_{\alpha_1}} {\mathcal M}(H;\alpha_1,\alpha_2)
\times_{R_{\alpha_2}} \dots \\
&\times_{R_{\alpha_{k-1}}} {\mathcal M}(H;\alpha_{k-1},\alpha_k)
\times_{R_{\alpha_k}} {\mathcal M}(H;\alpha_{k},\alpha_+),
\endaligned
$$
where $\alpha_1,\dots,\alpha_k \in \frak A$.
Condition \ref{linsysmainconds} (XI)
is an axiomatization of this property.
\par
We explain Condition \ref{linsysmainconds} (XII) in Subsubsection \ref{subsub:ccconds}.
\par
A system satisfying Condition \ref{linsysmainconds}
is called a {\it linear K-system}
(Definition \ref{linearsystemdefn} (2)).
Now the main result of Sections \ref{sec:systemline1}-\ref{sec:systemline4}
is as follows.
Suppose we are given a linear K-system.
We consider a direct sum of $\R$ vector spaces
\begin{equation}\label{new156}
\bigoplus_{\alpha \in \frak A} \Omega(R_{\alpha};o_{R_{\alpha}}).
\end{equation}
Using (\ref{new156}) and the energy filtration
we define (in Definitions \ref{Fvectspace} and \ref{defn141212}) a module
$$
CF(\mathcal C;\Lambda_{0,{\rm nov}})
$$
over the universal Novikov ring $\Lambda_{0,{\rm nov}}$.
(See Definition \ref{defn:Nov} for the definition of $\Lambda_{0,{\rm nov}}$.)
Here $\mathcal C$ denotes the totality of the part of data of our linear K-system
which is related to $\{R_{\alpha}\}$. We call it critical
submanifold data. (See Definition \ref{linearsystemdefn} (1).)

\begin{thm}\label{linesysmainth1Intro}
To each linear K-system,
we can associate a cochain complex, Floer cochain complex,
which we denote by $(CF(\mathcal C;\Lambda_{0,{\rm nov}}),d)$.
This complex is independent of the choices up to cochain homotopy equivalence.
\end{thm}
This is a slightly simplified version of Theorem \ref{linesysmainth1}.

\subsubsection{Construction of Floer cochain complex}
\label{subsub15-1-3}

We will prove Theorem \ref{linesysmainth1Intro} (or
Theorem \ref{linesysmainth1})
in detail in Section \ref{sec:systemline3}.
The proof we present there is written
in a way so that it is a prototype of
the proof of various similar results
and can be adapted easily to the proof
of similar results.
\par
The coboundary operator $d$ in Theorem \ref{linesysmainth1Intro}
is a sum of the exterior differential $d_0 : \Omega(R_{\alpha};o_{R_{\alpha}})
\to \Omega(R_{\alpha};o_{R_{\alpha}})$
and the operator $d_{\alpha_-,\alpha_+} :
\Omega(R_{\alpha_-};o_{R_{\alpha_-}}) \to \Omega(R_{\alpha_+};o_{R_{\alpha_+}})$ obtained
by the smooth correspondence
\begin{equation}\label{form15777}
\begin{CD}
R_{\alpha_-} @<{{\rm ev}_{-}}<< {\mathcal M}(H;\alpha_-,\alpha_+)
@>{{\rm ev}_{+}}>> R_{\alpha_+}.
\end{CD}
\end{equation}
Then (\ref{form1555}) together with
Stokes' formula (\cite[Theorem 9.26]{part11}) and Composition formula
(\cite[Theorem 10.20]{part11}) `imply'
$$
d_0 \circ d_{\alpha_-,\alpha_+} + d_{\alpha_-,\alpha_+} \circ d_0
+ \sum_{\alpha} d_{\alpha_-,\alpha} \circ d_{\alpha,\alpha_+}
=0.
$$
This will imply that $d = d_0 + \sum d_{\alpha_-,\alpha_+}$
satisfies $d\circ d = 0$.
Thus we obtain Floer cohomology.
\par
More precisely speaking, to define the operator
$d_{\alpha_-,\alpha_+}$ from (\ref{form15777})
as smooth correspondence we need to take and fix a CF-perturbation on
${\mathcal M}(H;\alpha_-,\alpha_+)$.
\par
To apply Stokes' formula, our CF-perturbations
of various moduli spaces must be compatible with the
isomorphism (\ref{form1555}).
Namely we need to show the next statement.
\begin{prop}{\rm (slightly imprecise statement)}\label{impreciseexistence}
For each given linear K-system and sufficiently small $\epsilon >0$, there exists a
system of CF-perturbations $\widehat{\frak S^{\epsilon}_{\alpha_-,\alpha_+}}$ on ${\mathcal M}(H;\alpha_-,\alpha_+)$
such that:
\begin{enumerate}
\item
$\widehat{\frak S^{\epsilon}_{\alpha_-,\alpha_+}}$ is transversal
to $0$ and ${\rm ev}_{+}$ is strongly submersive with respect
to this CF-perturbation.
\item
The restriction of $\widehat{\frak S^{\epsilon}_{\alpha_-,\alpha_+}}$ to
the boundary is equivalent to the fiber product
of $\widehat{\frak S^{\epsilon}_{\alpha_-,\alpha}}$
and $\widehat{\frak S^{\epsilon}_{\alpha,\alpha_+}}$
via the isomorphism (\ref{form1555}).
\end{enumerate}
\end{prop}
This is slightly imprecise statement and is {\it not}
the statement we will prove.
The precise statement we will prove is
Proposition \ref{prop161}.
The main difference between Proposition \ref{impreciseexistence}
and Proposition \ref{prop161} is the following.
\begin{enumerate}
\item[(a)]
The CF-perturbation
$\widehat{\frak S^{\epsilon}_{\alpha_-,\alpha_+}}$ is not defined on the
Kuranishi structure of  ${\mathcal M}(H;\alpha_-,\alpha_+)$ itself,
which is given by the axiom of linear K-system, but
defined on its thickening.
\item[(b)]
We replace ${\mathcal M}(H;\alpha_-,\alpha_+)$
by ${\mathcal M}(H;\alpha_-,\alpha_+)^{\boxplus\tau_0}$,
which
is a K-space obtained by putting the collar
to the space ${\mathcal M}(H;\alpha_-,\alpha_+)$ {\it outside}.
\item[(c)]
We fix $E_0$ and construct CF-perturbations for
only finitely many moduli spaces, that is, the moduli spaces consisting of
the elements of energy (\ref{formulaenergy}) $\le E_0$.
\end{enumerate}
The reason of Item (a) is that, to construct a CF-perturbation
we need first to construct a good coordinate system and
then go back to the Kuranishi structure.
We explained this point already in \cite[Subsection 1.2]{part11}.
\par
The reason of Item (b) is more technical. We perform various
operations in a neighborhood of the boundary  and corner
of our K-space.
Those constructions are easier to carry out if
the charts of our K-space have collars.
We will use it to extend the Kuranishi structure
on the boundary which is a thickening of the given one,
to the interior.
(See however Remark \ref{rem1555} (1).)
Existence of collar on a given cornered manifold or orbifold
is fairly standard fact in differential topology.
In case of Kuranishi structure or good coordinate system,
to put the collar to the all charts so that coordinate
changes preserve it is rather cumbersome.
(This is because the way to put collar to a given
cornered orbifold is not unique.)
We take a short cut and put the collar `outside'
rather than `inside'.
The process to put the collar outside is
describe in detail in Section \ref{sec:triboundary}. See
Subsection \ref{subsection:introextensionlem}
for more detailed explanation on this point.
\par
The reason of Item (c) is that it is difficult to perturb
infinitely moduli spaces simultaneourly.
We go back to this point in Subsubsection \ref{subsubsection:inductivelimit}.

\begin{rem}\label{rem1555} (1)
Recall that in this article
we start from a purely abstract setting of a K-space with boundary and corner,
which is not necessarily arising from
a particular geometric situation like the moduli space of pseudo-holomorphic curves.
In the geometric setting
studied in \cite{FO},\cite{fooobook},\cite{fooobook2},
an extension
of the Kuranishi structure
to a small neighborhood of $\partial X$ in $X$ is given from its construction.
In fact,
we start from a Kuranishi structure
$\partial\widehat{\mathcal U}$ on the boundary
(which we obtain from geometry or analysis) and construct
a good coordinate system  ${\widetriangle{\mathcal U_{\partial}}}$
and use it to find $\widehat{\mathcal U^+_{\partial}}$ and its perturbation.
We need to extend $\widehat{\mathcal U^+_{\partial}}$ and
its perturbation to a neighborhood of $\partial X$.
In this situation, the Kuranishi charts of
$\widehat{\mathcal U^+_{\partial}}$ are obtained as open subcharts
of certain Kuranishi charts of ${\partial}{\widetriangle{\mathcal U}}$.
(See the proof of \cite[Theorem 3.30, Proposition 6.44]{part11}.)
Therefore it can be indeed extended using the extension of
${\partial}\widehat{\mathcal U}$ directly.
\par
(2) In the proof of well-definedness of the virtual fundamental
chain, that corresponds to the well defined-ness of the Gromov-Witten invariant,
which is given in Part 1 of this article,
`trivialization of corner' is not necessary.
This is because
we only need to apply Stokes' formula and do not
need the chain level argument.
See \cite[Propositions 8.15,8.16]{part11} and their proofs.
\end{rem}

\subsubsection{Corner compatibility conditions}
\label{subsub:ccconds}

The proof of Proposition \ref{impreciseexistence}
(or Proposition \ref{prop161}) is
by induction on energy.
We consider the isomorphism (\ref{form1555}):
$$
\partial{\mathcal M}(H;\alpha_-,\alpha_+)
=
\bigcup_{\alpha}
{\mathcal M}(H;\alpha_-,\alpha)
\,\,{}_{{\rm ev}_{+}}\times_{{\rm ev}_{-}}
{\mathcal M}(H;\alpha,\alpha_+).
$$
We observe the energy of the moduli space
appearing in the right hand side is strictly smaller than
one appearing in the left hand side.
So by induction hypothesis the CF-perturbation
of the right hand side is already given.
Therefore the statement we need to work out
this induction is something like the following (*).
\begin{enumerate}
\item[(*)]
Let $(X,\widehat{\mathcal U})$ be a K-space with
corner. Suppose a CF-perturbation
$\widehat{\frak S^{\epsilon}_{\partial}}$ is given
on the normalized boundary $\partial(X,\widehat{\mathcal U})$,
satisfying certain transversality properties.
Then we can find a CF-perturbation
$\widehat{\frak S^{\epsilon}}$ on $(X,\widehat{\mathcal U})$
which has the same transversality property and
whose restriction to the boundary coincides with
$\widehat{\frak S^{\epsilon}_{\partial}}$.
\end{enumerate}
However, we note that the statement (*), as it is, does {\it not} hold.
In fact, since $(X,\widehat{\mathcal U})$ has not only boundary
but also corners, we need to assume certain
compatibility conditions for $\widehat{\frak S^{\epsilon}_{\partial}}$
at the corner.
Let us elaborate this point below.
\par
We remark that we use the {\it normalized} corner of an
orbifold (or Kuranishi structure) with  corners.
Typically
a point in the corner $\widehat S_2U$ of an orbifold $U$ corresponds
to two points in the {\it normalized} boundary.
In other words
we have a double cover
\begin{equation}\label{doublecover1519}
\pi : \partial\partial U \to \widehat S_2U.
\end{equation}
Suppose we are given a CF-perturbation $\widehat{\frak S_{\partial}^{\epsilon}}$
on the normalized
corner $\partial U$. The compatibility condition
we need to assume for $\widehat{\frak S_{\partial}^{\epsilon}}$
is that if $\pi(x) = \pi(y)$ then
the perturbation $\widehat{\frak S_{\partial}^{\epsilon}}$ at $x$
coincides with $\widehat{\frak S_{\partial}^{\epsilon}}$ at $y$.
Namely we need to require the next condition:
\begin{enumerate}
\item[($\star$)]
There exists a CF-perturbation $\widehat{\frak S_{S_2U}^{\epsilon}}$
on $\widehat S_2U$, whose pull-back to $\partial\partial U$
is equivalent to the restriction of
$\widehat{\frak S_{\partial}^{\epsilon}}$ to
$\partial\partial U$
\end{enumerate}
We can state a similar condition for Kuranishi structure on $X$.
We note that to state the condition ($\star$) precisely
we first need to clarify the relationship between the Kuranishi structure
on $\partial\partial X$ and one on $\widehat S_2X$.
\par
In the situation of our application
where $X ={\mathcal M}(H;\alpha_-,\alpha_+)$,
we have isomorphisms
\begin{equation}
\aligned
&\partial\partial{\mathcal M}(H;\alpha_-,\alpha_+) \\
&\cong
\partial\left(\bigcup_{\alpha}
{\mathcal M}(H;\alpha_-,\alpha)
\,\,{}_{{\rm ev}_{+}}\times_{{\rm ev}_{-}}
{\mathcal M}(H;\alpha,\alpha_+)\right) \\
&\cong
\bigcup_{\alpha}
\partial({\mathcal M}(H;\alpha_-,\alpha))
\,\,{}_{{\rm ev}_{+}}\times_{{\rm ev}_{-}}
{\mathcal M}(H;\alpha,\alpha_+) \\
&\quad\cup
\bigcup_{\alpha}
{\mathcal M}(H;\alpha_-,\alpha)
\,\,{}_{{\rm ev}_{+}}\times_{{\rm ev}_{-}}
\partial({\mathcal M}(H;\alpha,\alpha_+))
\\
&\cong
\bigcup_{\alpha_1,\alpha_2}
({\mathcal M}(H;\alpha_-,\alpha_1)
\,\,{}_{{\rm ev}_{+}}\times_{{\rm ev}_{-}}
({\mathcal M}(H;\alpha_1,\alpha_2))
\,\,{}_{{\rm ev}_{+}}\times_{{\rm ev}_{-}}
{\mathcal M}(H;\alpha_2,\alpha_+)\\
&\quad\cup
\bigcup_{\alpha_1,\alpha_2}
{\mathcal M}(H;\alpha_-,\alpha_1)
\,\,{}_{{\rm ev}_{+}}\times_{{\rm ev}_{-}}
(({\mathcal M}(H;\alpha_1,\alpha_2)
\,\,{}_{{\rm ev}_{+}}\times_{{\rm ev}_{-}}
{\mathcal M}(H;\alpha_2,\alpha_+)).
\endaligned
\nonumber
\end{equation}
On the other hand, by Condition \ref{linsysmainconds} (XI) we assumed:
$$
\aligned
&\widehat S_2{\mathcal M}(H;\alpha_-,\alpha_+) \\
&\cong
\bigcup_{\alpha_1,\alpha_2}
{\mathcal M}(H;\alpha_-,\alpha_1)
\,\,{}_{{\rm ev}_{+}}\times_{{\rm ev}_{-}}
{\mathcal M}(H;\alpha_1,\alpha_2)
\,\,{}_{{\rm ev}_{+}}\times_{{\rm ev}_{-}}
{\mathcal M}(H;\alpha_2,\alpha_+)
\endaligned
$$
By these isomorphisms we obtain a double cover
$$
\pi' : \partial\partial{\mathcal M}(H;\alpha_-,\alpha_+)
\to \widehat S_2{\mathcal M}(H;\alpha_-,\alpha_+).
$$
The Condition \ref{linsysmainconds} (XII)
(the second of corner compatibility condition)
requires that this double cover $\pi'$
coincides with the double cover $\pi$
in (\ref{doublecover1519}).
\begin{rem}
\begin{enumerate}
\item
The condition $\pi = \pi'$
(Condition \ref{linsysmainconds} (XII)) is {\it not} automatic
and we need to {\it assume} it as a part of the axiom of linear K-system.
In fact, we can define the covering map $\pi$ in a canonical way for an
arbitrary
K-space $X$. On the other hand, the covering map $\pi'$
depends on the choice of the isomorphism (\ref{form1555})
and similar isomorphisms for the corner.
(Condition \ref{linsysmainconds} (XI)).
Note that in our axiomatization
only the {\it existence} of the isomorphism (\ref{form1555})
is required. The isomorphism such as (\ref{form1555}) is not unique.
In fact, we can change it by composing any automorphism
of $\partial{\mathcal M}(H;\alpha_-,\alpha_+)$.
If we change the isomorphism (\ref{form1555}) then
the identity $\pi = \pi'$ will no longer hold.
\par
In other words, Condition \ref{linsysmainconds} (XII)
is one on the consistency between various choices of the isomorphisms
(\ref{form1555}) and similar isomorphisms for the corner.
\item
In our geometric situation, we define
the isomorphism (\ref{form1555}) using geometric
description of the boundary of our moduli space
${\mathcal M}(H;\alpha_-,\alpha_+)$.
Then the condition $\pi = \pi'$ is fairly obvious.
In this article, we need to state this condition explicitly
because our purpose here is to formulate the
precise conditions for our system of K-spaces under which
we can define the Floer cohomology
in the way independent of the geometric origin of such
a system of K-spaces.
\end{enumerate}
\end{rem}
Actually we need to require consistency at the corner of
arbitrary codimension using the covering space
\begin{equation}\label{cornercoverinintro}
\pi_{m,\ell} : \widehat S_m(\widehat S_{\ell} X)  \to \widehat S_{m+\ell} X,
\end{equation}
which exists for any K-space $X$ with corner.
(See Proposition \ref{prop2813}.)
If we assume the corner compatibility condition,
the property ($\star$) and its
analogue for higher codimensional
corner can be shown inductively.
The inductive step of this induction
can be stated as follows.
\begin{prop}{\rm (slightly imprecise statement)}
\label{propp1577}
Let $(X,\widehat{\mathcal U})$ be a
K-space with corners.
Suppose for each $k$ we have a
CF-pertubation $\widehat{\frak S_k}$ on $\widehat S_k(X,\widehat{\mathcal U})$
with the following properties.
\par
For each $m$ and $\ell$, the following two CF-perturbations
on $\widehat S_m(\widehat S_{\ell} X)$ are equivalent each other.
\begin{enumerate}
\item
The restriction of $\widehat{\frak S_{\ell}}$ to $\widehat S_m(\widehat S_{\ell} X)$.
\item
The pull-back of $\widehat{\frak S_{m+\ell}}$ by the
covering map (\ref{cornercoverinintro}).
\end{enumerate}
\par
Then there exists a CF-perturbation $\widehat{\frak S}$ on
$(X,\widehat{\mathcal U})$ such that its restriction
to $\widehat S_k(X,\widehat{\mathcal U})$ coincides with $\widehat{\frak S_k}$
for each $k$.
\end{prop}
This is a simplified statement and is {\it not} the
statement we will prove in Section \ref{sec:triboundary}.
The statement we will prove is Proposition \ref{prop529rev}.
The difference between Proposition \ref{propp1577} and
Proposition \ref{prop529rev} is the following.
\begin{enumerate}
\item[(a)]
The CF-perturbation we start with is not
given on $\widehat S_k(X,\widehat{\mathcal U})$
itself but is given on a thickening of $\widehat S_k(X,\widehat{\mathcal U})$.
The CF-perturbation we obtain is defined
on a thickening of $(X,\widehat{\mathcal U})$.
\item[(b)]
We replace $X$ by $X^{\boxplus\tau_0}$,
which is a K-space obtained from $X$ by putting the collar
outside.
\item[(c)]
We assume that $\widehat{\frak S_k}$ satisfies an appropriate transversality
property and will find a CF-perturbation $\widehat{\frak S}$
satisfying the same transversality
property.
\end{enumerate}
The reason for (a) is that we need to go once to a good coordinate system
and come back to construct a CF-perturbation.
The reason for (b) is explained in detail in Subsection \ref{subsection:introextensionlem}.
We actually need to construct a system of CF-perturbations satisfying
certain transversality properties. This is the reason for (c).

\subsubsection{Well-defined-ness of Floer cohomology
and morphism of linear K-system}
\label{subsub15-1-5}
The most important property of Floer cohomology of periodic Hamiltonian
system is its invariance under the choice of Hamiltonians.
Our story contains axiomatization of this equivalence.
For this purpose we introduce the notion of morphisms between two
linear K-systems.
To explain the relevant axiom
we consider the case of
linear K-system arising from the periodic Hamiltonian system and
the associated Floer equation.
Let $H^i : S^1 \times M \to \R$ be a periodic Hamiltonian
function for $i=1,2$.
Using the set of critical points we obtain a set of
manifolds $\{ R^i_{\alpha} \mid \alpha \in \mathfrak A_i\}$,
$i=1,2$ and we obtain a set, the compactified moduli space,
${\mathcal M}(H^i;\alpha_-,\alpha)$ of solutions of
Floer's equation (\ref{Floerseq}) for $H = H^{i}$,
$\alpha_{\pm} \in \mathfrak A_i$.
They define cochain complexes
$(CF(\mathcal C^i;\Lambda_{0,{\rm nov}}),d^i)$
and its Floer cohomologies by Theorem \ref{linesysmainth1Intro}
for $i=1,2$.
\par
The well established method to prove the independence of
the Floer cohomology of periodic Hamiltonian system under the choice of
Hamiltonian is to use the moduli space of
the next equation (\ref{Floerseqpara}).
(This method was invented by
Floer \cite{Flo89I}.)
We take a function $\mathcal H : \R \times S^1 \times M \to \R$
such that
\begin{equation}
\mathcal H(\tau,t,x) =
\begin{cases}
H^1(t,x)     &\text{if $\tau < -C$} \\
H^2(t,x)    &\text{if $\tau > C$}
\end{cases}
\end{equation}
where $C$ is a sufficiently large fixed number.
We put $H_{\tau,t}(x) = \mathcal H(\tau,t,x)$ and consider the equation
\begin{equation}\label{Floerseqpara}
\frac{\partial u}{\partial \tau}
+
J\left(
\frac{\partial u}{\partial t} - X_{H_{\tau,t}}(u)
\right)
=0
\end{equation}
with the asymptotic boundary condition for $\tau \to \pm\infty$
given by $R^1_{\alpha_-}$, $R^2_{\alpha'_+}$, respectively.
We denote the compactified moduli space of the solution of
(\ref{Floerseqpara}) with this boundary condition by
$\mathcal M(H_{\tau,t};\alpha_-,\alpha'_+)$.
We will define the notion of morphism
of linear K-systems
in Definition \ref{linearsystemmorphdefn} and
Condition \ref{morphilinsys}.
The notion of {\it interpolation space}
$\mathcal N(\alpha_-,\alpha'_+)$
appearing in the definition of morphisms is the axiomatization
of the properties of this moduli space
$\mathcal M(H_{\tau,t};\alpha_-,\alpha'_+)$.
For example, (\ref{formula1211morph}) corresponds to the
property of the boundary of $\mathcal M(H_{\tau,t};\alpha_-,\alpha'_+)$,
that is,
\begin{equation}\label{form151313}
\aligned
\partial \mathcal M(H_{\tau,t};\alpha_-,\alpha'_+)
\cong
&\bigcup_{\alpha \in \frak A_1}
\mathcal M(H^1;\alpha_-,\alpha)
\times_{R^1_{\alpha}} \mathcal M(H_{\tau,t};\alpha,\alpha'_+) \\
&\cup\bigcup_{\alpha' \in \frak A_2}
\mathcal M(H_{\tau,t};\alpha_-,\alpha')
\times_{R^2_{\alpha'}} \mathcal M(H^2;\alpha',\alpha'_+).
\endaligned
\end{equation}
Thus the set of K-spaces $\{\mathcal M(H_{\tau,t};\alpha_-,\alpha'_+)
\mid \alpha_- \in \frak A_1, \alpha'_+ \in \frak A_2\}$
together with various other data defines a
morphism from the linear K-system associated to $H^1$
to the linear K-system associated to $H^2$.
\par
In the study of Floer cohomology of periodic Hamiltonian system
the moduli space
$\mathcal M(H_{\tau,t};\alpha_-,\alpha'_+)$ is used to define
a cochain map from Floer's cochain complex associated to $H^{1}$
to Floer's cochain complex associated to $H^{2}$.
We can carry out this construction by using the properties
spelled out in Definition \ref{linearsystemmorphdefn} and
Condition \ref{morphilinsys} only and prove the next result.
\begin{thm}\label{thm158585}
If $\frak N$ is a morphism from one linear K-system $\mathcal F_1$
to another linear K-system $\mathcal F_2$, then
$\frak N$ induces a cochain map
$$
\frak N_* : (CF(\mathcal C^1;\Lambda_{{\rm nov}}),d^1)
\to (CF(\mathcal C^2;\Lambda_{{\rm nov}}),d^2).
$$
Here $(CF(\mathcal C^i;\Lambda_{\rm nov}),d^i)$
is the cochain complex associated to $\mathcal F_i$
by Theorem \ref{linesysmainth1Intro}.
\par
The cochain map $\frak N_*$ depends on various choices.
However it is independent of the choices up to
cochain homotopy.
\end{thm}
Theorem \ref{thm158585} is Theorem \ref{linesysmainth2} (1).
The proof is similar to the proof of Theorem \ref{linesysmainth1Intro}
and is given in Subsection \ref{subsec:proofsec14main2}.
\par
We define the notion of composition
of morphisms in Section \ref{section:compomorphis} and show that
$\frak N \mapsto \frak N_*$ is functorial
with respect to the composition of the morphisms in
Subsection \ref{subsec:compochain}.
When the interpolation spaces of the morphism
$\frak N_{i+1i} : \mathcal F_i \to \mathcal F_{i+1}$ is given by
$\mathcal N_{ii+1}(\alpha^i,\alpha^{i+1})$
for $i=1,2$, the interpolation space
of the composition $\frak N_{32} \circ \frak N_{21} : \mathcal F_1 \to \mathcal F_{3}$
is the K-space
$\mathcal N_{13}(\alpha^1,\alpha^{3})$
obtained, roughly speaking, by gluing the K-spaces
\begin{equation}\label{compositioninterpolation}
\mathcal N_{12}(\alpha^1,\alpha^{2})
\times_{R^2_{\alpha_2}} \mathcal N_{23}(\alpha^2,\alpha^{3})
\end{equation}
for various $\alpha^2$ along the boundaries and corners.
We need to smooth a part of corners of this fiber product
to glue. For this purpose, we need the definition of smoothing corners
of K-spaces with corners. We will discuss it in
Section \ref{section:compomorphis}.
See Subsection \ref{subsec:introcomp}
for an issue of smoothing corners
of K-spaces.
To define  smoothing corners in a canonical way
we use the collar (which was put outside).
So more precisely we use
\begin{equation}\label{form16200Intro}
\bigcup_{\alpha_2 \in \frak A_2}
\mathcal N_{12}(\alpha_1,\alpha_2) \times^{\boxplus\tau}_{R^2_{\alpha_2}}
\mathcal N_{23}(\alpha_2,\alpha_3)
\end{equation}
in place of (\ref{compositioninterpolation}).
See Definition \ref{defn1635} for the definition of (\ref{form16200Intro}).
\par
We also define the notion of homotopy and homotopy of homotopies
etc. of morphisms and show that homotopy between morphisms
$\frak N$ and $\frak N'$
induces a cochain homotopy between $\frak N_*$ and $\frak N'_*$.

\subsubsection{Identity morphism}
\label{subsub15-1-6}
To make the assignment $\frak N \mapsto \frak N_*$
functorial, we need the notion of the identity morphisms.
In the second half of Section \ref{section:compomorphis}
we define and prove a basic property of the identity morphism
of linear K-system.
\par
In the geometric situation of the linear K-system arising
from periodic Hamiltonian system,
the morphism among such linear K-systems are defined by
using the moduli space of the solutions of
equation (\ref{Floerseqpara}),
where $\mathcal H : \R \times S^1 \times M \to \R$.
To obtain the identity morphism of linear K-system
associated to $H : S^1 \times M \to \R$,
we consider the case of $\mathcal H$ such that
$\mathcal H(\tau,t,x) = H(t,x)$.
In other words, we use $\tau$ independent $\mathcal H$.
However, note that the moduli space of solutions of (\ref{Floerseqpara})
for this $\tau$ independent $\mathcal H$
is {\it different} from the moduli space of
solutions of Floer's equation (\ref{Hameq}).
Namely
$$
\mathcal M(H_{\tau,t};\alpha_-,\alpha_+)
\ne
\mathcal M(H;\alpha_-,\alpha_+)
$$
in case $H_{\tau,t} = H_t$ for all $\tau$.
Indeed, the dimensions are different.
To define $\mathcal M(H;\alpha_-,\alpha_+)$
we divide our space by the $\R$ action given by
translation on $\tau \in \R$ direction.
Since $\mathcal H$ is happen to be $\tau$ independent,
our equation (\ref{Floerseqpara}) is invariant under this
$\R$ action, too.
However, by definition,
$\mathcal M(H_{\tau,t};\alpha_-,\alpha_+)$
is a special case of general $\mathcal H$.
For general $\mathcal H$, (\ref{Floerseqpara}) is {\it not}
invariant under $\R$
action.
{\it Before compactifiation} we can identity
$$
\overset{\circ}{\mathcal M}(H;\alpha_-,\alpha_+) \times \R
=
\overset{\circ}{\mathcal M}(H_{\tau,t};\alpha_-,\alpha_+),
$$
when $H_{\tau,t} = H_t$.
However the relationship between compactified moduli spaces
${\mathcal M}(H;\alpha_-,\alpha_+)$ and
${\mathcal M}(H_{\tau,t};\alpha_-,\alpha_+)$ is not so simple.
\par
We describe in Subsection \ref{subsec:identitylinsys}
 a way to obtain ${\mathcal M}(H_{\tau,t};\alpha_-,\alpha_+)$
(the case $H_{\tau,t}$ is $\tau$ independent)
from $\overset{\circ}{\mathcal M}(H;\alpha_-,\alpha_+)$,
in an abstract setting. In other words, we start with
the spaces of connecting orbits $
\mathcal M(\alpha_-,\alpha_+) = \mathcal M(H;\alpha_-,\alpha_+)$
of a liner K-system $\mathcal F$ and define
the interpolation spaces of the identity morphism
$\mathcal{ID} : \mathcal F \to \mathcal F$.
We also show that identity morphism is a `homotopy unit'.
Namely we show in Subsection \ref{subsec:identitylinsys} that
the composition of the identity morphism $\mathcal{ID}$
with other morphism $\frak N$ is homotopic to $\frak N$.
(Proposition \ref{properidentity}.)
\par
To construct the identity morphism and prove its
homotopy-unitality, we imitate the proof of the
corresponding results in the case of periodic Hamiltonian
system, and rewrite it so that it works in the
purely abstract setting of linear K-system
without any specific geometric origin.
Although we explain the geometric origin of the
construction of Subsection \ref{subsec:identitylinsys}
in Subsection \ref{subsection:identitygeo},
the discussion of Subsection \ref{subsection:identitygeo}
is {\it not} used in Subsection \ref{subsec:identitylinsys}
or any other part to prove main results of this article.
We expect that those explanation is useful
for readers who know Floer cohomology of periodic Hamiltonian
system to understand the contents of Subsection \ref{subsec:identitylinsys}.
\par
In Section \ref{sec:systemline3}, we use the identity morphism to prove
the second half of Theorem \ref{linesysmainth1Intro},
that is, independence of the Floer cochain complex
$(CF(\mathcal C;\Lambda_{0,{\rm nov}}),d)$
of the choices up to cochain homotopy equivalence, as follows.
We first consider the case when our linear K-system is obtained from
a Hamiltonian $H : S^1 \times M \to \R$ by using Floer's
equation (\ref{Floerseq}).
We fix a choice of an almost complex structure and
the Kuranishi structure on $\mathcal M(H;\alpha_-,\alpha_+)$.
Namely we fix choices which determine a linear K-system.
We then take two different systems of CF-perturbations on it,
which we denote by
$\widehat{\frak S^{i,\epsilon}_{\alpha_-,\alpha_+}}$,
$i=1,2$.
We then obtain two different cochain complexes
which we denote by $(CF(\mathcal C;\Lambda_{0,{\rm nov}}),d^i)$,
$i=1,2$.
We want to prove that they are cochain homotopy equivalent.
\par
In this case the interpolation spaces of the
identity morphism are
$\mathcal M(H_{\tau,t};\alpha_-,\alpha_+)$
(for various $\alpha_{\pm}$)
with $H_{\tau,t} = H_t$ for all $\tau$.
Its boundary is described by (\ref{form151313}).
In our situation it becomes:
\begin{equation}\label{form1513132}
\aligned
\partial \mathcal M(H_{\tau,t};\alpha_-,\alpha_+)
\cong
&\bigcup_{\alpha \in \frak A}
\mathcal M(H;\alpha_-,\alpha)
\times_{R_{\alpha}} \mathcal M(H_{\tau,t};\alpha,\alpha_+) \\
&\bigcup_{\alpha' \in \frak A}
\mathcal M(H_{\tau,t};\alpha_-,\alpha')
\times_{R_{\alpha'}} \mathcal M(H;\alpha',\alpha_+).
\endaligned
\end{equation}
Now we consider a CF-perturbation
$\widehat{\frak S^{1,\epsilon}_{\alpha_-,\alpha}}$
on $\mathcal M(H;\alpha_-,\alpha)$, which
is the first factor of the first term of the
right hand side,
and another CF-perturbation
$\widehat{\frak S^{2,\epsilon}_{\alpha',\alpha_+}}$
on $\mathcal M(H;\alpha',\alpha_+)$, which is the second factor of
second term of the
right hand side.
We then take a system of CF-perturbations
on various
$\mathcal M(H_{\tau,t};\alpha_-,\alpha_+)$
so that these CF-perturbations
together with
$\widehat{\frak S^{1,\epsilon}_{\alpha_-,\alpha}}$,
$\widehat{\frak S^{2,\epsilon}_{\alpha',\alpha_+}}$
are compatible with the isomorphism
(\ref{form1513132}).
(To show the existence of such a system of
CF-perturbations, we need to examine all the
corners of arbitrary codimension and check
the compatibility at the corners.
We can do so by induction using a similar
argument as explained in Subsubsection \ref{subsub:ccconds}.)
Then the correspondence by
$\mathcal M(H_{\tau,t};\alpha_-,\alpha_+)$
together with this system of CF-perturbations
defines a cochain map from
$(CF(\mathcal C;\Lambda_{0,{\rm nov}}),d^1)$
to
$(CF(\mathcal C;\Lambda_{0,{\rm nov}}),d^2)$.
This is a consequence of
Stokes' formula (\cite[Proposition 9.16]{part11}) and Composition formula
(\cite[Theorem 10.20]{part11}).
\par
This cochain map is actually an isomorphism since
it is the identity map modulo $T^{\epsilon}$ for some
$\epsilon >0$.
\par
In the case of linear K-system, which may not
come from a particular geometric construction,
we can proceed in the same way using the identity
morphism, to prove that
$(CF(\mathcal C;\Lambda_{0,{\rm nov}}),d^1)$
is cochain homotopy equivalent to
$(CF(\mathcal C;\Lambda_{0,{\rm nov}}),d^2)$.

\subsubsection{Homotopy limit}
\label{subsubsection:inductivelimit}

We note that to construct a system of
CF-perturbations for all the spaces of
connecting orbits appearing in a linear K-system,
we need to find infinitely many CF-perturbations
simultaneously.
There is an issue to do so. We explained this issue
in detail in \cite[Subsection 7.2.3]{fooobook2}.
The method to resolve it is the same as \cite[Section 7.2]{fooobook2}.
The algebraic part of this method is summarized as follows.
For $E>0$ a pair $(C,d)$ of a free $\Lambda_0$ module $C$ and
$d : C \to C$ is said to be a {\it partial cochain complex of energy cut level $E$} if
$d \circ d \equiv 0 \mod T^E$.
Let $(C_1,d), (C_2,d)$ be partial cochain complexes of energy cut level  $E$.
A $\Lambda_0$ module homomorphism $\varphi : C_1 \to C_2$ is
said to be a {\it partial cochain map of energy cut level $E$} if
$\varphi \circ d \equiv d \circ \varphi \mod T^E$.
We also note that if $(C,d)$ is a partial cochain complex of energy cut level  $E'$ and
if $E < E'$ then  $(C,d)$ is a partial cochain complex of energy cut level  $E$.
\begin{lem}\label{prop159999}
Let $(C_i,d)$ be gapped partial cochain complex of energy cut level $E_i$
for $i=1,2$ with $E_1 < E_2$. Let $\varphi : C_1 \to C_2$ be
a gapped partial cochain map of energy cut level $T^{E_1}$.
We assume $\overline\varphi : C_1/\Lambda_{+,\text{\rm nov}}C_1 \to C_2/\Lambda_{+, \text{\rm nov}}C_2$
is an isomorphism.
\par
Then there exist $d^+ : C_1 \to C_1$ and $\varphi^+ : C_1 \to C_2$ such that:
\begin{enumerate}
\item
$(C_1,d^+)$ is a partial cochain complex of energy cut level ${E_2}$.
\item
$\varphi^+ : (C_1,d^+) \to (C_2,d)$ is a partial cochain map of energy cut level $E_2$.
\item
$d^+ \equiv d \mod T^{E_1}$ and $\varphi^+ \equiv \varphi \mod T^{E_1}$.
\end{enumerate}
\end{lem}
See Definition \ref{defn141212} for the definition of gapped-ness.
Lemma \ref{prop159999} is Lemma \ref{lem1623}.
We use Lemma \ref{prop159999} to construct the cochain complex
$(CF(\mathcal C;\Lambda_{0,{\rm nov}}),d)$ appearing in Theorem \ref{linesysmainth1Intro}
as follows.
We take $0 < E_1 < E_2 < \dots$ with $ E_i \to \infty$.
We use the argument outlined in Subsubsections \ref{subsub15-1-3}
-\ref{subsub:ccconds}
using the finitely many moduli spaces (consisting of elements
of energy $< E_i$) to construct
cochain complex $(CF(\mathcal C;\Lambda_{0,{\rm nov}}),d^i)$
modulo $T^{E_i}$ for each $i$.
We next use the argument  outlined in Subsubsections \ref{subsub15-1-5}
-\ref{subsub15-1-6}
to find a cochain map
$\varphi_i : (CF(\mathcal C;\Lambda_{0,{\rm nov}}),d^i)
\to (CF(\mathcal C;\Lambda_{0,{\rm nov}}),d^{i+1})$
modulo $T^{E_i}$ for each $i$.
Now we use Lemma \ref{prop159999} inductively
and to obtain $d^{i}_k : CF(\mathcal C) \to CF(\mathcal C)$
for $k > i$
and
$\varphi_{i,k} : (CF(\mathcal C),d^{i}_k) \to (CF(\mathcal C),d^{i+1}_k)$
such that:
\begin{enumerate}
\item
$(CF(\mathcal C),d^{i}_k)$
is a cochain complex modulo $T^{E_k}$.
\item
$\varphi_{i,k}$ is a cochain map module $T^{E_k}$.
\item
$d^{i}_k \equiv d^{i}_{k+1} \mod T^{E_k}$,
$\varphi_{i,k} \equiv \varphi_{i,k+1} \mod T^{E_k}$.
\end{enumerate}
Then $\lim_{k\to \infty} d^{1}_k :
CF(\mathcal C)\to CF(\mathcal C)$
becomes the required coboundary oprator.
\par
To construct a cochain map we use
a similar argument using homotopy
modulo $T^{E}$ instead of cochain map
modulo $T^{E}$.
To construct a cochain homotopy
between cochain maps,
we also use a similar argument
using homotopy of homotopies modulo $T^E$.
Algebraic lemmas we use in place of
Lemma \ref{prop159999} or Lemma \ref{lem1623}
are Propositions \ref{prop1631} and \ref{prop1639}.

\subsubsection{Story over rational coefficient}
\label{subsec:storyoverQ}

In Section \ref{sec:systemline4} we consider the case when
all the spaces $R_{\alpha}$ are 0-dimensional
and prove that we can use Novikov ring whose ground ring is $\Q$ in that case.
The proof is based on the results of
\cite[Sections
13 and 14]{part11}.

\subsection{Outline of the story of tree-like K-system}
\label{treelikeout}

In Sections \ref{sec:systemtree1}-\ref{sec:systemtree2}
we study systems of K-spaces, which axiomatize
the situation appearing during the construction of
the filtered $A_{\infty}$ algebra associated to a Lagrangian
submanifold. (\cite{fooobook,fooobook2}.)

\subsubsection{Moduli space of pseudo-holomorphic disks: review}
\label{subsec:19-1background}

In this subsubsection, we review basic properties of the
moduli space of pseudo-holomorphic disks to motivate the
definitions in later subsubsections.
\par
Let $M$ be a symplectic manifold and $L$ its Lagrangian submanifold.
We assume that $M$ is compact or
tame (i.e., carrying a tame almost complex structure)
and $L$ is
compact, oriented and relatively spin.
We have the Maslov index group homomorphism
$
\mu : H_2(M,L;\Z) \to 2\Z
$
and the energy group homomorphism
$
E : H_2(M,L;\Z) \to \R
$
defined by
$$
E(\beta) = \int_{D^2} u^*\omega
$$
with $[u] = \beta$.
For $\beta \in H_2(M,L;\Z)$ we consider the moduli space
$\overset{\circ}{\mathcal M}_{k+1}(\beta)$\footnote{Though it is better to write it
as $\overset{\circ}{\mathcal M}_{k+1}(L;\beta)$, we omit $L$ for the simplicity of notation.}
consisting of
$((D^2,\vec z),u)$ such that:
\begin{enumerate}
\item
$u : (D^2,\partial D^2) \to (M,L)$ is pseudo-holomorphic.
\item
$\vec z = (z_0,\dots,z_k)$ are $k+1$ marked points of the boundary
$\partial D^2$.
\item $z_i \ne z_j$ if $i\ne j$.
\item
$(z_0,\dots,z_k)$ respects the counter clockwise cyclic
order of $\partial D^2$.
\end{enumerate}
We define evaluation maps
$$
{\rm ev} = ({\rm ev}_0,\dots,{\rm ev}_k)
~:~ \overset{\circ}{\mathcal M}_{k+1}(\beta) \longrightarrow L^{k+1}
$$
by ${\rm ev}_i((D^2,\vec z),u) = u(z_i)$.
Then
$\overset{\circ}{\mathcal M}_{k+1}(\beta)$
has a compactification ${\mathcal M}_{k+1}(\beta)$
to which ${\rm ev}_i$ is extended.
Moreover
${\mathcal M}_{k+1}(\beta)$ has an oriented Kuranish structure with corners
of dimension
$$
\dim {\mathcal M}_{k+1}(\beta)
= \mu(\beta) + k-2.
$$
The normalized boundary of ${\mathcal M}_{k+1}(\beta)$ is a
disjoint union of the fiber products:
$$
{\mathcal M}_{k_1+1}(\beta_1) \,\, {}_{{\rm ev}_0}\times_{{\rm ev}_i}
{\mathcal M}_{k_2+1}(\beta_2)
$$
where $\beta_1 + \beta_2 = \beta$,
$k_1 + k_2 = k+1$, $i=1,\dots,k_2$.
\par
These facts are proved in \cite[Subsection 7-1]{fooo09}.

\subsubsection{Axiom of tree-like K-system and
main theorem constructing the filtered $A_{\infty}$ algebra}
\label{subsec:19-2treelike}

Axioms of the {\it tree-like K-system over $L$}
or the {\it $A_{\infty}$ correspondence over $L$} are
given as Conditions \ref{linAinfmainconds} and
Definition \ref{defn198} and
is obtained by axiomatizing the properties of the system of
the moduli spaces ${\mathcal M}_{k+1}(\beta)$ and the
evaluation maps ${\rm ev}_i$, which are described in the
previous subsubsection.
The way to axiomatize various structures are parallel
to the case of linear K-system and we do not repeat it.
The main result we obtain is the next theorem:
For a closed oriented manifold $L$ we denote by
$
\Omega(L)
$
the de Rham complex of $L$.
We put
$$
\Omega(L;\Lambda_0)
=
\Omega(L) \widehat{\otimes} \Lambda_0.
$$
See Definition \ref{defn:Nov} for the coefficient ring $\Lambda_0$.
Here $\widehat{\otimes}$ denotes the completion of the algebraic
tensor product.
Namely, an element of $\Omega(L;\Lambda_0)$ is a formal sum
$$
\sum_{i=0}^{\infty} T^{\lambda_i} h_i
$$
where $\lambda_i \in \R_{\ge 0}$ with $\lambda_1 < \lambda_2 < \dots$,
$\lim_{i\to \infty} \lambda_i = + \infty$ and $h_i \in \Omega(L)$.
\begin{thm}\label{thm11510}
Suppose $({\mathcal M}_{k+1}(\beta),{\rm ev},\mu,E)$ is a tree-like K-system
over $L$.
Then we can associate a filtered $A_{\infty}$ structure
$\{\frak m_{k} \mid k=0,1,2\dots\}$ on $\Omega(L;\Lambda_0)$.
\par
The filtered $A_{\infty}$ algebra $(\Omega(L;\Lambda_0),\{\frak m_{k} \mid k=0,1,2\dots\})$
is independent of the various choices up to homotopy equivalence.
\end{thm}
This is the de Rham version of the half of \cite[Theorem A]{fooobook}.
(\cite[Theorem A]{fooobook} also contains the part of constructing
tree-like K-system arising from a geometric situation described in
the previous subsubsection.)
It is a consequence of Theorem \ref{theorem1934} (1)
by putting
\begin{equation}\label{defnmkmk0}
\frak m_k = \sum_{\beta} T^{E(\beta)} \frak m_{k,\beta}.
\end{equation}
We recall that the filtered $A_{\infty}$ structure assigns maps
\begin{equation}\label{defnmkmk}
\frak m_k :
 \underbrace{\Omega(L;\Lambda_0)[1] \widehat\otimes \dots \widehat\otimes \Omega(L;\Lambda_0)[1]}
_{\text{$k$ times}}
\to \Omega(L;\Lambda_0)[1]
\end{equation}
$k=0,1,2,\dots$,
that satisfy the $A_{\infty}$ relations
\begin{equation}\label{Ainfinityre}
\aligned
\sum_{k_1+k_2=k+1}\sum_{i=1}^{k-k_2+1}
(-1)^*{\frak m}_{k_1}(x_1,\ldots,{\frak m}_{k_2}(x_i,\ldots,x_{i+k_2-1}),\ldots,x_{k}) = 0.
\endaligned\end{equation}
Here $[1]$ is the degree $+1$ shift functor.
When we define $\frak m_k$ by (\ref{defnmkmk0}),
the formula (\ref{Ainfinityre}) follows from (\ref{Ainfinityrelbeta}).
See \cite[Definition 3.2.20]{fooobook} for the definition
of filtered $A_{\infty}$ algebra and \cite[Definition 4.2.42]{fooobook}
for the definition of homotopy equivalence of filtered $A_{\infty}$ algebras.
Roughly speaking, the $A_{\infty}$ operation
$$
\frak m_{k,\beta} ~:~ \underbrace{\Omega(L)[1] \otimes \dots \otimes \Omega(L)[1]}
_{\text{$k$ times}}
\longrightarrow \Omega(L)[1]
$$
is defined by
\begin{equation}\label{roughAope}
\frak m_{k,\beta}(h_1,\dots,h_k) =
{\rm ev}_0!({\rm ev}_1^*h_1 \wedge \dots \wedge {\rm ev}_k^*h_k)
\end{equation}
using the correspondence
\begin{equation}\label{smoothcorrespondence}
	\xymatrix{
                              &
{\mathcal{M}_{k+1} (\beta)} \ar[dl]_{(\text{\rm ev}_1,\dots,\text{\rm ev}_k)} \ar[dr]^{\text{\rm ev}_0} &
      \\
 L^k & & L}
\end{equation}
See (\ref{form2012}) for the precise
definition.\footnote{Strictly speaking,
in \eqref{form2012} we define a {\it partial} $A_{\infty}$ algebra structure
(see Definition \ref{defn:partial})
which depends on a parameter $\epsilon >0$.
We then use a `homotopy limit' in the way similar to that of the construction of
$A_\infty$ algebra explained
in Subsubsection \ref{subsubsection:inductivelimit}.}
The integration along the fiber ${\rm ev}_0!$
in the formula (\ref{roughAope}) is defined
by using an appropriate system of CF-perturbations $\widehat{\frak S}_{k+1}(\beta)$
on the K-space ${\mathcal M}_{k+1}(\beta)$,
which is the main part of the data defining
tree-like K-system.
\par
The formula (\ref{Ainfinityre})
(or (\ref{Ainfinityrelbeta}))
is obtained from Stokes' formula (\cite[Theorem 9.26]{part11}) and Composition formula
(\cite[Theorem 10.20]{part11})
via the isomorphism
\begin{equation}\label{boudaryisointroain}
\partial {\mathcal M}_{k+1}(\beta)
=
\bigcup_{k_1+k_2=k}\bigcup_{i=1,\dots,k_2} \bigcup_{\beta_1+\beta_2=\beta}
{\mathcal M}_{k_1+1}(\beta_1) \,\, {}_{{\rm ev}_0}\times_{{\rm ev}_i}
{\mathcal M}_{k_2+1}(\beta_2),
\end{equation}
which is a part of the axiom of a tree-like K-system,
Condition \ref{linAinfmainconds} (IX).
\par
For this argument to work, we need to
choose a system of CF-perturbations $\widehat{\frak S}_{k+1}(\beta)$
so that it is compatible with the isomorphism (\ref{boudaryisointroain}).
Proposition \ref{prop203333} is the precise
statement which claims the existence of such a system.
\par
Construction of such a system of CF-perturbations
is parallel to that of a linear K-system.
We uses an induction over $k$ and $E(\beta)$ to construct
$\widehat{\frak S}_{k+1}(\beta)$.
The inductive step of this
construction uses Proposition \ref{propp1577}
(or its precise version Proposition \ref{prop529rev}.)
To inductively verify the assumptions of  Proposition \ref{propp1577}
we need to construct our system
$\widehat{\frak S}_{k+1}(\beta)$ to be compatible
not only along the boundary but also at the corners.
Therefore we need to assume the compatibility of
Kuranishi structures on ${\mathcal M}_{k+1}(\beta)$
at the corners.
This is the corner compatible conditions
Condition \ref{linAinfmainconds} (X) and (XI).

\subsubsection{Bifurcation method and pseudo-isotopy}
\label{subsec:19-3bifurcation method}

As explained in the previous subsubsection,
the construction of a filtered $A_{\infty}$
structure from the tree-like K-system
given in this article
is mostly similar to the construction
of Floer's cochain complex from a
linear K-system.
\par
The difference between two constructions
lies in the
morphism part of the construction.
In the case of linear K-system we defined a
morphism between two such K-systems
and associated to the morphism a cochain map between
their Floer's cochain complexes.
(In particular, using the identity
morphism we proved independence of
the resulting cochain complex under the various choices we make, modulo cochain homotopy equivalence.)
In the situation of tree-like K-system
we define the notion of pseudo-isotopy
between two tree-like K-systems and
of filtered $A_{\infty}$ algebras.
Then we show that the resulting
pseudo-isotopy between two tree-like K-systems
induces a pseudo-isotopy between the $A_{\infty}$ algebras.
It is easy to show that
two filtered $A_{\infty}$ algebras
are homotopy equivalent if they are
pseudo-isotopic.
(See  \cite[Theorem 8.2]{fooo091}.)
\par
In our geometric situation
of Lagrangian Floer theory, a
pseudo-isotopy of the tree-like K-system
is obtained as follows.
In the situation of Subsubsection \ref{subsec:19-1background},
we consider two
compatible almost complex structures $J_1,J_2$ on $M$.
Then we obtain the moduli spaces of pseudo-holomorhic
disks $\mathcal M_{k+1}(\beta;J_i)$ for $i=1,2$.
For each $i=1,2$ we fix some choices to define
a system of Kuranishi structures on $\mathcal M_{k+1}(\beta;J_i)$
so that it defines a tree-like K-system.
We denote these choices by $\Xi_i$
and the K-space obtained via these choices  $\mathcal M_{k+1}(\beta;J_i;\Xi_i)$.
\par
Now we consider a one parameter family of compatible almost
complex structures $\{J_t \mid t \in [1,2]\}$
which joins $J_1$ to $J_2$.
We consider the moduli space
\begin{equation}\label{eq1521}
\mathcal M_{k+1}(\beta;[1,2])
=
\bigcup_{t \in [1,2]} \mathcal M_{k+1}(\beta;J_t) \times \{t\}.
\end{equation}
Here $\mathcal M_{k+1}(\beta;J_t)$ is the moduli space
of $J_t$ holomorphic discs with boundary condition $L$,
homology class $\beta$, and $k+1$ marked points.
We can find a system of Kuranishi structures on it
such that its restriction to the part
$t =1$ (resp. $t=2$) coincides with
$\Xi_1$ (resp. $\Xi_2$.)
We have the evaluation maps
$$
{\rm ev} = ({\rm ev}_0,\dots,{\rm ev}_k) :
\mathcal M_{k+1}(\beta;[1,2])
\to
L^{k+1}
$$
and
$$
{\rm ev}_{[1,2]} : \mathcal M_{k+1}(\beta;[1,2])
\to
[1,2].
$$
We axiomatize the properties of the system consisting of
$\mathcal M_{k+1}(\beta;[1,2])$,
the evaluation maps, etc.. and define the notion of $[1,2]$-parametrized
family of $A_{\infty}$ correspondences.
(See  Condition \ref{linAinfmaincondspara} and
Definition \ref{defn1913piso}. We define more general
notion of $P$-parametrized $A_{\infty}$ correspondence
in Definition \ref{defn2013}.)
\par
We now make choices of CF-perturbations etc.
on $\mathcal M_{k+1}(\beta;J_i;\Xi_i)$ $i=1,2$ and we use
them to construct the filtered $A_{\infty}$ structures.
\par
In the same way as in
(\ref{roughAope}) we use the evaluation maps ${\rm ev}$ and
${\rm ev}_{[1,2]}$ together with our CF-perturbations to
define operators
$$
\frak m_{k,\beta} : \Omega(L \times [1,2])^{\otimes k}
\to \Omega(L \times [1,2])
$$
Then
$
\frak m_k = \sum_{k,\beta} T^{E(\beta)} \frak m_{k,\beta}
$
satisfies the $A_{\infty}$ relation
(\ref{Ainfinityre}).
The system of operators $\frak m_{k,\beta}$ on
$\Omega(L \times [1,2])$
which satisfies the $A_{\infty}$ relation and some additional
properties is called a {\it pseudo-isotopy of
filtered $A_{\infty}$ algebras}.
See Definition \ref{pisotopydef2}.
\par
Thus we will prove the following:
\begin{thm}\label{thm1511}
A pseudo-isotopy of $A_{\infty}$ correspondences
induce a pseudo-isotopy of filtered $A_{\infty}$
algebras.
\end{thm}
Theorem \ref{thm1511} is Theorem \ref{theorem1934} (3).
\par
We note that Theorem \ref{thm1511} implies
the second half of
Theorem \ref{thm11510} as follows.
Let $({\mathcal M}_{k+1}(\beta),{\rm ev},\mu,E)$ be an $A_{\infty}$
correspondence.
We can define a pseudo-isotopy
of this $A_{\infty}$
correspondence with itself
by taking
$$
{\mathcal M}_{k+1}(\beta;[1,2])
= {\mathcal M}_{k+1}(\beta) \times [1,2]
$$
etc..
Then we apply Theorem \ref{thm1511}
to show that the filtered $A_{\infty}$
algebras obtained by two different
CF-perturbations
from $({\mathcal M}_{k+1}(\beta),{\rm ev},\mu,E)$ are pseudo-isotopic to
each other.
\begin{rem}\label{rem1512}
Here we use
the construction of a pseudo-isotopy
from an $A_{\infty}$
correspondence to itself for the construction of a pseudo-isotopy of $A_\infty$ algebras
in the way similar as we use the identity morphism
in Subsubsection \ref{subsub15-1-6} for the construction of a cochain map between
Floer cochain complexes.
We like to mention that the construction of a pseudo-isotopy
from an $A_{\infty}$
correspondence to itself
is much easier than the construction of the
identity morphism.
\end{rem}

For the actual proof of Theorems \ref{thm11510}
and \ref{thm1511} we need to use `homotopy limit'
argument similar to those in Subsubsection \ref{subsubsection:inductivelimit}.
We need the notion of pseudo-isotopy of pseudo-isotopies etc.
for this purpose.
The algebraic lemma corresponding to Lemma \ref{prop159999}
is Propositions \ref{prop208} and \ref{prop2013}.

\subsubsection{Bifurcation method and self-gluing}

In this article we use morphism of K-systems
to prove independence of the
Floer's cochain complex associated to
a given linear K-system of the choices.
On the other hand, we use pseudo-isotopy to
prove independence of the
filtered $A_{\infty}$ structure
associated to
a tree-like K-system of the choices.
Actually we can also use morphism for the
tree-like K-system
and pseudo-isotopy
for the linear K-system.
We use two different methods in order
to demonstrate both of these two methods.
We may call `cobordism method'
\index{cobordism method}
instead of `the method using morphism',
and `bifurcation method'
\index{bifurcation method}
instead of
`the method using pseudo-isotopy'.
The difference of those two methods
is explained also in \cite[Subsection 7.2.14]{fooobook2}.
\par
The cobordism method is used in the Lagrangian Floer theory
in \cite{fooobook}. An axiomatization of
morphism of tree-like K-system
is given in \cite{fooo010}.
The bifurcation method is used in
Lagrangian Floer theory in, for example,
\cite{AFOOO,akahojoyce,fooo091}.
Each of these two methods has certain advantage and disadvantage.
\par
One advantage of the bifurcation method is that
usually it is shorter and simpler
to use the bifurcation method than
the cobordism method.
See for example, Remark \ref{rem1512}.
\par
On the other hand, we cannot
prove independence of Floer cohomology
of periodic Hamiltonian system
under the change of Hamiltonian
function, by bifurcation method.
This is because we need to study the
situation where the sets of critical
points are different.
\par
Micheal Hutchings \cite{H2}\footnote{Actually
Hutchings' concern is not so much
on the proof of independence of Morse-Novikov cohomology
of the choices but rather the explicit form
of the cochain homotopy equivalence,
between two Morse-Novikov complexes before and after
wall crossing.
The discussion below clarifies the way to prove the
independence of Morse-Novikov cohomology,
but to find the explicit form
of of the cochain homotopy equivalence we need to study
more. See \cite {H1}.} and Paul Seidel \cite[Remark 10.14]{Se}
mentioned some issue to prove
invariance of Floer cohomology by using
the bifurcation method.
We explain below how those issue
was resolved in our previous writings.
\par
We discuss the case of Morse-Novikov cohomology.
Let $M$ be a compact Riemannian manifold and
$h$ a closed 1-form given by the exterior derivative
of a Morse function locally.
Let $R(h)$ be the set of the critical points of $h$.
For $p,q \in R(h)$ we consider the compactified moduli space
of gradient lines of $h$ joining $p$ to $q$
and denote it by $\mathcal M(h;p,q)$.
(We identify two gradient lines $\ell, \ell'$ as an
element of $\mathcal M(h;p,q)$ if $\ell(\tau) = \ell'(\tau
+\tau_0)$ for some $\tau_0 \in \R$.)
We take its subset $\mathcal M(h;p,q;E)$ such that
an integration of $h$
along the gradient line is $E$.
The matrix element of the coboundary oprator of Morse-Novikov complex
is the sum of signed counts of the order of $\mathcal M(h;p,q;E)$
together with the weight $T^E$.
(We assume that the gradient vector field of $h$
is Morse-Smale.)
The proof that it defines a cochain complex
is the same as the case of
Morse complex, which is similar to one we
explained in
Subsubsection \ref{subsub15-1-5},
as was observed by Novikov.
\par
The issue is the way how to prove the
independence of the cohomology of the cochain complex associated to $h$
when we move $h$.
Suppose we have two closed 1-forms $h$ and $h'$ whose de Rham cohomology
classes in $H^1(M)$ coincide. For simplicity we assume $R(h) = R(h')$.
We take a one parameter family $h_t$ such that
$R(h) = R(h_t)$ and $h_0 = h$, $h_1 = h'$.
We consider
\begin{equation}\label{152221}
\mathcal M(h_*;p,q;E)
=
\bigcup_{t \in [0,1]} \mathcal M(h_t;p,q;E) \times \{t\}
\end{equation}
and try to use it to show this independence.
((\ref{152221}) is similar to (\ref{eq1521}). So
the method we explain below is a bifurcation method.)
\par
Note that the virtual dimension of $\mathcal M(h;p,p;E)$
is $-1$. Therefore the virtual dimension of $\mathcal M(h_*;p,p;E)$
is $0$. So there may be a discrete
subset $\{t_i\} \subset [0,1]$ where $\mathcal M(h_t;p,p;E)$
is nonempty.
Now since $\mathcal M(h_{t_1};p,p;E)$ is nonempty,
$\mathcal M(h_{t_1};p,p;2E)$ contains an object
obtained by concatenating an element of $\mathcal M(h_{t_1};p,p;E)$
with itself.
(In other words if $[\ell] \in \mathcal M(h_{t_1};p,p;E)$
then
$$
([\ell],[\ell]) \in
\mathcal M(h_{t_1};p,p;E) \times \mathcal M(h_{t_1};p,p;E)
\subset \mathcal M(h_{t_1};p,p;2E).)
$$
By continuing this process we obtain an element of
$\mathcal M(h_{t_1};p,p;kE)$ for any $k$.
So we have an element of the strata of arbitrary negative
dimension.
We explain three different ways to resolve this issue.
\par\medskip
\noindent
{\bf 1.}
Instead of the moduli spaces (\ref{152221}) we can
use a different moduli space below.
We take a non-decreasing function $\chi : \R \to [0,1]$ such that
$\chi(\tau) = 0$ for small $\tau$ and $\chi(\tau) = 1$ for
large $\tau$.
We consider the `nonautonomous' equation
\begin{equation}\label{cobordismwq}
\frac{d\ell}{d\tau}(\tau) = {\rm grad}\, h_{\chi(\tau)}
\end{equation}
such that $\ell(+\infty) = q$ and $\ell(-\infty) = p$.
Let $\mathcal M(h_{\chi(\tau)};p,q;E)$ be the set of
solutions of this equation with energy $E$.
Contrary to the definition of coboundary oprator
there is no translational symmetry.
By counting the order of $\mathcal M(h_{\chi(\tau)};p,q;E)$
we obtain a cochain map from the Morse-Novikov complex of
$h$ to one of $h'$. The standard argument shows that
it becomes a cochain homotopy equivalence.
\cite{Flo89I}.
\par
This is the standard approach to show the well-definedness of
Morse-Novikov cohomology. We note that
there is no self-gluing issue in this approach since
the equation (\ref{cobordismwq}) is not invariant under this
self gluing construction.
\par
In other words, when using the cobordism method
the issue of self gluing never occurs.
\par\bigskip
\noindent
{\bf 2.}
We next explain how the usage of the bifurcation method together
with the de Rham model resolves the issue of `self-gluing'.
\par
We consider the moduli space
$\mathcal M(h_*;p,q;E)$ in (\ref{152221}).
We define `evaluation maps' to
the interval $[0,1]$. Namely
we send
$\mathcal M(h_t;p,q;E) \times \{t\}$
to $t \in [0,1]$.
We usually consider the situation
where both projections ${\rm pr}_s$ and ${\rm pr}_t$ exist,
where the former is the source projection and the latter is the
target projection.
In our setting they are both the same map
defined above.
We have a diagram:
$$
\begin{CD}
[0,1]  @<{\rm pr}_s<< \mathcal M(h_*;p,q;E) @>{\rm pr}_t>> [0,1]
\end{CD}
$$
\par
The `pseudo-isotopy' of Morse-Novikov complex is
a cochain complex defined on
$$
\left(\bigoplus_{p\in R(h)} \Omega([0,1]) \otimes [p]\right) \widehat \otimes \Lambda_0
$$
(where $\Omega([0,1])$ is the de Rham complex
of the interval)
or its completion by using the Novikov ring.
We take the interval for each $p \in R(h)$ and
denote it by $[0,1]_p$.
So `pseudo-isotopy' of Morse-Novikov complex is
defined on
\begin{equation}\label{13form}
CF(h_*) = \Omega \left(
\coprod_{p \in R(h)} [0,1]_p
\right) \widehat \otimes \Lambda_0.
\end{equation}
The above diagram is regarded as
\begin{equation}
\begin{CD}
[0,1]_p  @<{\rm pr}_s<< \mathcal M(h_*;p,q;E) @>{\rm pr}_t>> [0,1]_q.
\end{CD}
\end{equation}
It `defines' a map
$
d_{p,q;E} : \Omega([0,1]_p) \to \Omega([0,1]_q)
$
by
\begin{equation}\label{dbysmocorr}
d_{p,q;E}(u) = ({\rm pr}_t)_!({\rm pr}_s^*(u)).
\end{equation}
Note that the pull back of differential form is defined under
rather mild assumption.
However the push out or integration along the fiber
$({\rm pr}_t)_!$ is harder to define.
This point is indeed related to the self-gluing issue
as follows.
\par
Suppose $\mathcal M(h_*;p,p;E)$ is transversal.
It consists of finitely many points.
Let us assume that it consists of a single point $\frak p$
and $t_0 = {\rm pr}_s(\frak p) = {\rm pr}_t(\frak p)$.
We take $1\in \Omega^0([0,1])$. Then
$$
d_{p,p;E}(1) = ({\rm pr}_t)_!({\rm pr}_s^*(1))
$$
is the delta form $\delta_{t_0} dt$ supported at $t_0$.
Now we remark that the pull back
$$
{\rm pr}_s^*(\delta_{t_0} dt)
$$
is not defined. The standard condition for distribution
to be pulled back is {\it not} satisfied in this case.
This point is related to the
fact the fiber product
$$
\mathcal M(h_*;p,p;E)\,
{}_{{\rm pr}_t}\times_{{\rm pr}_s}
\mathcal M(h_*;p,p;E)
$$
is not transversal.
\par
This discussion clarifies that the reason why the
problem of self-gluing occurs lies in the fact that ${\rm pr}_t$ is not a submersion.
We can not expect submersivity because of the dimensional reason.
We note that we are somehow in the Bott-Morse situation
here even in the case when the set of critical points of $h_t$ is a discrete set
for each $h_t$, in case we study a one parameter family of
Morse forms.
Therefore the way to resolve this issue is similar
to the way to study the Bott-Morse situation.
This issue can be taken care of both in de Rham and singular homology
models.
Let us first explain the case of de Rham model.
\par
The problem here is that ${\rm pr}_t$ is not a submersion.
The solution to this problem is to use a CF-perturbation.
As a simplified version of the CF perturbation, we take a
family of perturbations (globally) parameterized by a finite
dimensional space, say $W$.
For $w \in W$ we have perturbed moduli space
$\mathcal M(h_*;p,p;E;w)$. We put
$$
\mathcal M(h_*;p,p;E;W) = \bigcup_{w \in W} \mathcal M(h_*;p,p;E;w) \times \{w\}.
$$
By taking the dimension of $W$ sufficiently large we may assume that the map
\begin{equation}\label{15}
{\rm pr}_t : \mathcal M(h_*;p,p;E;W) \to [0,1]
\end{equation}
is a submersion. (The space $W$ depends on $p,E$.)
Let ${\rm pr}_W : \mathcal M(h_*;p,p;E;W) \to W$ be the
projection.
We take a differential form $\chi_W$ of degree $\dim W$ and
with compact support such that
$\int_W \chi_W = 1$. Now we define
$$
d_{p,p;E}(u) = ({\rm pr}_t)_!({\rm pr}_s^*(u) \wedge {\rm pr}_W^*\chi_W).
$$
Since (\ref{15}) is a submersion this is always well defined.
In this way we can define a cochain complex on (\ref{13form})
by
$$
\delta = d + \sum T^Ed_{p,p;E}
$$
where $T$ is a formal parameter. (Novikov parameter.)
We can also show $\delta \circ \delta = 0$.
We use $(CF(h_*),\delta)$ to prove
that $CF(h)$ is cochain homotopy equivalent to $CF(h')$ as follows.
\par
By considering embeddings $\{0\} \to [0,1]$ and $\{1\} \to [0,1]$
we have a
map
$$
CF(h_*) \to CF(h), \qquad CF(h_*) \to CF(h').
$$
We can show that they are cochain maps.
Moreover using the fact that de Rham cohomology of $[0,1]$ is $\R$
we can show that they are cochain homotopy equivalence.
Thus we find that
$CF(h)$ is cochain homotopic to $CF(h')$.
This is a baby version of the proof of independence of
filtered $A_{\infty}$ structure of the almost complex structure
etc. using the  pseudo-isotopy,
which we present in Sections \ref{sec:systemtree1}-\ref{sec:systemtree2}.
\par
In this formulation, we have an equality
$$
d_{p,p;E_1}\circ d_{p,p;E_2} = 0.
$$
This is because $d_{p,p;E}$ increase degree of differential
form by $1$ and de Rham complex of $[0,1]$ has elements in only
0-th and 1-st degree.
So self-gluing problem does not occur.
\par\bigskip
\noindent
{\bf 3.}
We finally prove the independence of Morse-Novikov
homology using the singular homology model by the
bifurcation method.
\par
Let $P$ be a smooth singular chain of $[0,1]_p$.
THat is, it is a pair of a simplex
and a smooth map from it to $[0,1]_p$.
\par
The analogue of (\ref{dbysmocorr}) in singular homology is
as follows.
We take fiber product
\begin{equation}\label{155517}
P \times_{{\rm pr}_s} \mathcal M(h_*;p,q;E)
\end{equation}
and take its triangulation.
Using the map ${\rm pr}_t$
we regard it as an element of singular chain complex
of $[0,1]_q$.
So we consider
$$
CF(h_*)^{s} = \bigoplus_{p \in R(h)} S([0,1]_p) \widehat \otimes \Lambda_0.
$$
Here $S([0,1]_p)$ is the smooth singular chain
complex of the interval $[0,1]_p$.
By (\ref{155517}), we `obtain'
$$
d_{p,q;E} : S([0,1]_p) \to S([0,1]_q).
$$
and boundary operator on $CF(h_*)^{s}$.
\par
The issue is the fiber product (\ref{155517}) may not be transversal.
Also there is no way to perturb
$\mathcal M(h_*;p,q;E)$
so that (\ref{155517}) is transversal for {\it all} $P$.
\par
The idea to resolve this issue, which appeared in \cite[Proposition 7.2.35 and etc.]{fooobook2},
is the following:
We first take fiber product (\ref{155517})
and then perturb it.
In other words, our perturbation
{\it depends not only on $\mathcal M(h_*;p,q;E)$ but also
on the singular chain $P$}.
\footnote{The way we explain below is a slightly
improved version of the one appeared in \cite{fooooverZ}.
In \cite{fooobook2} we took a countably generated subcomplex
of smooth singular chain complex.
(We use Baire's category theorem uncountably many
times in \cite{fooooverZ} then
we do not need to take countably generated subcomplex
as we did in \cite{fooobook2}.) Here we take singular chain complex
itself, that is the way of \cite{fooooverZ}.}
\par
Since we can make ${\rm pr}_s$ a submersion on each Kuranishi
chart, (this is the definition of weak submersivity!),
we can consider the fiber product (\ref{155517}) which carries a Kuranishi structure.
Then we define
$$
\mathcal M(h_*;p,q;P) := P \times_{{\rm pr}_s} \mathcal M(h_*;p,q;E),
$$
which is a space with Kuranishi structure.
We take a system of perturbations of all of them and triangulations of
their zero sets such that
the following holds.
\begin{enumerate}
\item
(See \cite[Compatibility Condition 7.2.38]{fooobook2}.)
On $\mathcal M(h_*;p,q;\partial P)
\subset \partial \mathcal M(h_*;p,q;P)$
the perturbations and triangulations are compatible.
\item
(See \cite[Compatibility Condition 7.2.44]{fooobook2}.)
On
\begin{equation}\label{form18}
\aligned
&P \times_{{\rm pr}_s} \mathcal M(h_*;p,r;E_1) \,\,
{}_{{\rm pr}_t}\times_{{\rm pr}_s} \mathcal M(h_*;r,q;E_2)\\
&\subset
\partial \mathcal M(h_*;p,q;E_1+E_2;P)
\endaligned
\end{equation}
the perturbations and triangulations are compatible.
\end{enumerate}
The meaning of (1) is clear. Let us explain the
meaning of (2).
We consider the fiber product
$\mathcal M(h_*;p,r;P) = P \times_{{\rm pr}_s} \mathcal M(h_*;p,r;E)$.
The perturbation (multisection) and a triangulation of the perturbed space
(the zero set of the multisection)
are given for this space. We regard the triangulated space of
the perturbed moduli space as a singular chain $\sum Q_i$
of $[0,1]_r$.
Then
the left hand side is the union of
$$
Q_i \times_{{\rm pr}_s} \mathcal M(h_*;r,q;E)
= \mathcal M(h_*;r,q;E;Q_i).
$$
The perturbation and a triangulation of its zero set
are also given.
We require that the restriction of the perturbation of
$\mathcal M(h_*;p,q;E;P)$ and the triangulation of its zero set
coincide with the ones which are combination of $\mathcal M(h_*;r,q;E;Q_i)$
and of $\mathcal M(h_*;p,r;P)$.
\par
Let us elaborate the last point more.
At each point in (\ref{form18}) the obstruction bundle
is a direct sum of ones of $\mathcal M(h_*;p,r;E_1)$
and of $\mathcal M(h_*;r,q;E_2)$.
We require that the first component of the
perturbation is one for $\mathcal M(h_*;p,r;P)$
and the second component of the perturbation is
one for $\mathcal M(h_*;r,q;E;Q_i)$.
\par
This is the meaning of the compatibility (2).
Construction of the perturbation and a triangulation satisfying
(1)(2) are give by an induction over $E$ and $\dim P$.
This is the way of obtaining
$
d_{p,q;E} : S([0,1]_p) \to S([0,1]_q)
$ and the way taken in \cite{fooobook2}.
\par
We elaborate this construction a bit more explicitly and show how it resolves
the issue of self-gluing.
We consider the case of $\mathcal M(h_*;p,p;E)$ that is zero dimensional.
Suppose for simplicity that it consists of one point and
its $t$ coordinate is $t_0$.
We consider a $0$-chain $P(t_1) = \{t_1\} \in [0,1]_p$.
If $t_1 \ne t_0$ then
$P(t_1) \times_{{\rm pr}_s} \mathcal M(h_*;p,p;E)$
is transversal and is the empty set.
If $t_1 = t_0$ then
$P(t_0) \times_{{\rm pr}_s} \mathcal M(h_*;p,p;E)$
is {\it not} transversal and we need to perturb it.
After perturbation it becomes empty again.
\par
Next we consider a 1-chain $P(a,b) = [a,b] \subset [0,1]_p$.
If $a,b \ne t_0$, then
$P(a,b) \times_{{\rm pr}_s} \mathcal M(h_*;p,p;E)$
is transversal. It is an empty set if $t_0 \notin [a,b]$
and is one point if $t_0 \in [a,b]$.
If $a = t_0$, then
$P(t_0,b) \times_{{\rm pr}_s} \mathcal M(h_*;p,p;E)$
is {\it not} transversal.
We already fixed perturbation of
$P(t_0) \times_{{\rm pr}_s} \mathcal M(h_*;p,p;E)
\subset \partial P(t_0,b) \times_{{\rm pr}_s} \mathcal M(h_*;p,p;E)$.
We extend it to obtain a perturbation of
$P(t_0,b) \times_{{\rm pr}_s} \mathcal M(h_*;p,p;E)$.
Whether it becomes an empty set or a one-point set depends on the
choice of the perturbation of
$P(t_0) \times_{{\rm pr}_s} \mathcal M(h_*;p,p;E)$.
\par
Now we consider the self-gluing.
We take the fiber product
\begin{equation}\label{1919}
[0,1]_p \times_{{\rm pr}_s}
\mathcal M(h_*;p,p;E)_{{\rm pr}_t}\times_{{\rm pr}_s} \mathcal M(h_*;p,p;E).
\end{equation}
We do not perturb
$[0,1]_p \times_{{\rm pr}_s} \mathcal M(h_*;p,p;E)$
and this consists of a single  point which is mapped to $t_0$ by ${\rm pr}_t$.
So the second fiber product is not transversal.
However we have already fixed the perturbation of
$P(t_0) \times_{{\rm pr}_s} \mathcal M(h_*;p,p;E)$
and by this perturbation (\ref{1919}) becomes
the empty set.
In other words, by this perturbation
the perturbed zero set does {\it not} hit the corner.
This is the way how the self-gluing issue is resolved.
This argument is a version of the way we handled
the Bott-Morse situation in \cite{fooobook2}.
\par
Note that Akaho-Joyce \cite{akahojoyce} used the bifurcation method
to show the well-defined-ness of the $A_{\infty}$ structure
using the singular homology.
We think the way they adopt is basically the same as we described
here.

\subsection{Outline of the appendices}
\label{introappend}

\subsubsection{Orbifolds and covering space of orbifolds/K-spaces}
\label{subsubofd}

Section \ref{sec:ofd} is a review of the notion of
orbifolds and vector bundles on them.
We consider effective orbifolds only and
use embeddings only as morphisms.
In this way we can avoid several delicate
issues arising in the discussion of orbifolds.
If we go beyond those cases, we need
to work with the framework of the 2-category to have a proper
notion of morphisms.
We also use the language of
chart and  coordinate transformation,
which is closer to the standard definition
of manifold.
It is well-know that there is an alternative
way using the language of
groupoid. (See for example \cite{ofdruan}.)
Using the groupoid language is somewhat similar
to the way taken in algebraic
geometry to define the notion of
stacks.
One advantage of using the groupoid language
is that the discussion then becomes
closer to the `coordinate free' exposition.
We remark that in our definition of
Kuranishi structure using the coordinates
and the coordinate transformations
is inevitable. No `coordinate
free' definition of Kuranishi structure is
known.
So we think that the coordinate description of orbifolds is more natural
for the study of Kuranishi structure.
In other approach to Kuranishi-like structure such as Joyce's, which is closer
to that of algebraic geometry, the groupoid description of orbifolds
seems to be more natural.
\par
To define the bundle extension data (See \cite[Definition 12.24]{part11})  which we used in
\cite[Sections 12 and 13]{part11}
we use some basic results of vector bundle in its orbifold version.
They are well known and have been established long time ago.
Since its proof is rather a straight forward modification
of the case of manifolds,
it seems that it is hard to find a reference which proves
them in the literature. We provide the proof
of those facts (Lemma \ref{lem2328}, Corollary \ref{cor2939},
Propositions \ref{prop2942} and \ref{prop2949}, etc.) by this reason.
\par
In Section \ref{sec:cover} we discuss the covering space of
an orbifold and a K-space, and define the covering space
\begin{equation}\label{form1531}
\widehat S_m(\widehat S_{\ell} X) \to \widehat S_{m+\ell}X
\end{equation}
for the formulation of the corner compatibilty condition.
See Subsubsection \ref{subsub:ccconds}.
We define the notion of covering space of orbifolds
in Subsections \ref{subsec:cover}
and generalize it to the case of K-spaces in Subsection \ref{subsec:cover}.
Then the covering space (\ref{form1531}) is defined in
Subsection \ref{subsec:coverconer}.

\subsubsection{Admissibility of orbifolds and of Kuranishi structures}

In Section \ref{sec:admKura}, we discuss the notion of
admissible orbifolds and
admissible Kuranishi structures.
Admissibility we study here is the property
of the coordinate change etc. with respect to the
coordinate normal to the boundary or corner.
Admissibility is used in the discussion of
Section \ref{sec:triboundary} to put the collar `outside'.
We explain how it is used there briefly below.
\par
We consider the case of an $n$ dimensional manifold $X$ with
boundary $\partial X$.
(For the simplicity of exposition we assume $X$ has
a boundary but no corner.)
Let $p \in \partial X$.
We take its coordinate chart and so
we have a diffeomorphism $\psi_p$ from
$\overline V_p \times [0,1)$ to
a neighborhood $U_p$ of $p$ in $X$.
Here $\overline V_p$ is an open subset of $\R^{n-1}$.
The space $X^{\boxplus 1}$
is obtained by taking
$\overline V_p \times [-1,1)$ for each $p$.
We glue them as follows.
Let $q \in \partial X$
and we take $\psi_q$,
$\overline V_q \times [0,1)$,
$U_q$ as above.
Let
$V_{pq} = \psi_q^{-1}(U_p)$
which is an open subset of
$\overline V_q \times [0,1)$.
The coordinate change is
$$
\varphi_{pq}
=  \psi_p^{-1} \circ \psi_q: V_{pq} \to \overline V_p \times [0,1).
$$
We extend it to
$$
\varphi_{pq}^{\boxplus 1} :
V^{\boxplus 1}_{pq} \to \overline V_p \times [0,1)
$$
as follows.
We put $\overline V_{pq} = V_{pq} \cap
(\overline V_q \times \{0\})$.
The restriction of $\varphi_{pq}$ to $\overline V_{pq}$
defines a map $\overline\varphi_{pq} : \overline V_{pq}
\to \overline V_p$. (Here we identify
$\overline V_p = \overline V_p \times \{0\}$.)
We put
$$
V^{\boxplus 1}_{pq} = V_{pq} \cup (\overline V_{pq} \times [-1,0])
$$
where we glue two spaces in the right hand side
at $\overline V_{pq} \times \{0\}$.
The map $\varphi_{pq}^{\boxplus 1}$ is defined by
$$
\varphi_{pq}^{\boxplus 1}(x,t)
=
\begin{cases}
\varphi_{pq}(x,t) &\text{if $(x,t) \in V_{pq}$}, \\
(\overline\varphi_{pq}(x),t) &\text{if $(x,t) \in \overline V_{pq} \times [-1,0]$}.
\end{cases}
$$
It is easy to see that these maps $\varphi_{pq}^{\boxplus 1}$
satisfy an appropriate cocycle condition.
Then we can glue various charts by these maps $\varphi_{pq}^{\boxplus 1}$
to obtain a space $X^{\boxplus 1}$.
It is also clear that $X^{\boxplus 1}$
is a topological manifold since
$\varphi_{pq}^{\boxplus 1}$ is a homeomorphism
to its image.
\par
However in general, $\varphi_{pq}^{\boxplus 1}$ is {\it not}
differentiable at $\overline V_{pq} \times \{0\}$.
So in general there is no obvious smooth structure
on $X^{\boxplus 1}$.
\par
We introduce the notion of the admissibility of manifolds (orbifolds, K-spaces)
so that if we start from an adimissible manifold then
the coordinate change $\varphi_{pq}^{\boxplus 1}$ becomes smooth.
Roughly speaking, the admissibility means that we are given
distinguished system of coordinates
so that the coordinate change $\varphi_{pq}$ among those coordinates
has the following additional properties.
\begin{enumerate}
\item[(*)]
We put
$
\varphi_{pq}(x,t) = (y(x,t),s(x,t))
$
then
$$
s(x,t) - t,
\qquad
\frac{\partial}{\partial t}y(x,t)
$$
together with all of their derivatives go to
zero as $t \to 0$.
\end{enumerate}
(More precisely, we assume certain exponential decay.)
It is easy to see that (*) implies that
$\varphi_{pq}^{\boxplus 1}$ is smooth.
So $X^{\boxplus 1}$ becomes a smooth manifold
in case $X$ is an admissible manifold and
we use admissible coordinate to define $X^{\boxplus 1}$.
\par
We can generalize the admissibility to the case of
orbifold with corners.
See Definitions \ref{def:3112} and \ref{def:3113}.
Then we can define admissibility
of various notions on an admissbile orbifold.
For example, admissibility of
a vector bundle, a section of it, and
a smooth map to another manifold (without boundary).
We can also define admissibility of an embedding of
an admissible orbifold
to another admissible orbifold.
We can use them to define admissibility of
Kuranishi charts and coordinate changes.
Then we can define the notion of {\it admissible
Kuranishi structure}. (Definition \ref{def:adKura}.)
On an admissible Kuranishi structure
we can define the notion of admissible
CF-perturbation.
It is mostly obvious that
the story of Part 1 can be worked out
in the admissible category.
One slightly nontrivial point to check
is the existence of
bundle extension data (see \cite[Definition 12.24]{part11})
in the admissible category. We used this notion in
\cite[Sections 12 and 13]{part11}.
So we need to establish the existence of
admissible bundle extension data
to prove the existence of CF-perturbation
etc.
in the admissible category.
As we mentioned in Subsubsection \ref{subsubofd},
the existence of bundle extension data
is proved by using certain standard construction
of vector bundle etc. (eg. existence of tubular
neighborhood). The proof of admissible version
of those standard results are nothing more than
obvious adaptations of the standard proofs.
Nevertheless for completeness' sake we provide
those proofs in Section \ref{sec:admKura}.
\par
Using admissibility we can extend
vector bundle $E$ on $X$ to a vector bundle
$E^{\boxplus 1}$ on $X^{\boxplus 1}$.
Also an admissible section of $E$
is canonically extended to a smooth section of
$E^{\boxplus 1}$.
\par
This is the way how we put a collar to the outside of a K-space $X$ in
Section \ref{sec:triboundary} and obtain
$X^{\boxplus 1}$.
We can also extend various admissible object of
a K-space $X$ to a collared object in $X^{\boxplus 1}$.
Thus the operation $X \mapsto X^{\boxplus 1}$
from admissible objects to collared objects is
completely canonical and functorial.
We write those constructions in Section \ref{sec:triboundary}
in detail for completeness. However we emphasize that
this construction is indeed straightfoward.
\par
We note that for {\it any} orbifold $X$ with corner there exists
a system of charts by which $X$ becomes an admissible orbifold.
In the case of manifold with boundary,
there exists a coordinate system
such that the coordinate change $
\varphi_{pq}
=  \psi_p^{-1} \circ \psi_q: V_{pq} \to \overline V_p \times [0,1)
$
preserves the second factor $[0,1]$ and the first factor
($\overline V_p$ factor) of $\varphi_{pq}(x,t)$ depends only
on $x$.
This statement is nothing but the existence of the collar
`inside' of a manifold $X$ with boundary.
The existence of the collar of any manifold or orbifold with
corner is a classical fact which is easy to prove.
So there is not much reason to put a collar
`outside' in the case of an orbifold. However, in the case of
Kuranishi structure, there is some cumbersome issue
to give a detailed proof of an existence of Kuranishi structure
so that all the coordinate changes preserve the collar.
(As we mentioned before, this is because the way to put
collar to a manifold is not canonical.)
The short cut we take is to use the admissible structure
to put the collar outside which is canonical.
\par
To apply this story to our geometric situation
such as the case of the moduli space of pseudo-holomorphic
curves, we need to establish existence of the
admissible structure for such moduli spaces.
This point is related to the exponential decay
estimate of the gluing analysis in the following
way.
\par
The boundary or corner of the moduli space of
pseudo-holomorphic curves appears typically
at the infinity of the moduli space and
the coordinate normal to the boundary or the
corner is the gluing parameter.
Let us consider the case of moduli space
of pseudo-holomorphic disks and consider the configuration of the
three disks as in Figure \ref{Figurep71} below.
\begin{figure}[h]
\centering
\includegraphics[scale=0.4]{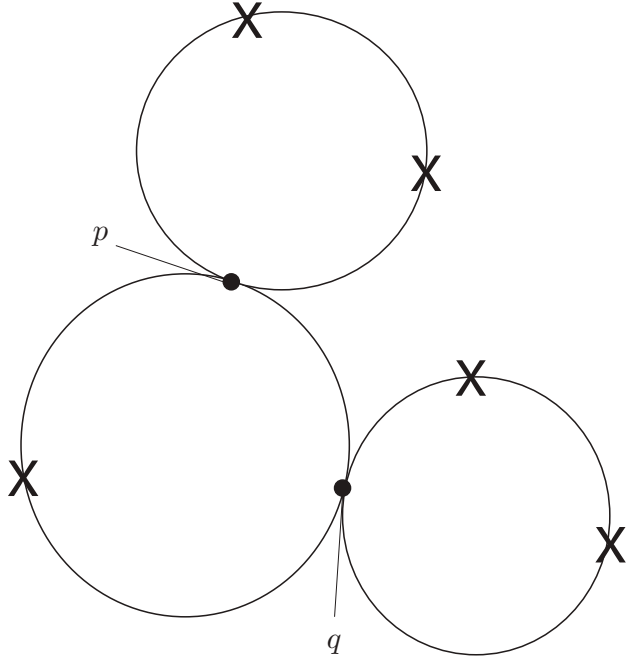}
\caption{Bordered curve consisting of three disks}
\label{Figurep71}
\end{figure}
This curve has two boundary nodes written  $p$ and $q$
in the figure.
(We add five boundary marked points so that this configuration is
stable.)
The parameter space to resolve these singularities involve two real numbers.  We write them as $T_{p}$ and $T_{q}$.
They are the length of the neck $[0,T_p]\times [0,1]$
(resp. $[0,T_q]\times [0,1]$) so
$T_p,T_q \in (C,\infty]$ for some large positive number $C$.
\par
Let us resolve these two singularities.
We can do so in one of the following three different ways.
\begin{enumerate}
\item
We first resolve the singularity at $p$ and
then at $q$.
\item
We first resolve the singularity at $q$ and
then at $p$.
\item
We resolve the two singularities at the same time.
\end{enumerate}
Let $\mathcal M_{5}(\beta)$ be the compactified moduli
space of pseudo-holomorphic discs with 5 marked points,
of homology class $\beta \in H_2(M,L)$ and boundary condition
given by a certain Lagrangian submanifold $L \subset M$.
The configuration in Figure \ref{Figurep71}
together with appropriate pseudo-holomorphic maps, gives
an element of this compactified moduli space.
Suppose this element is Fredholm regular for simplicity.
Let $\overline V$ be the intersection of its neighborhood in  $\mathcal M_{5}(\beta)$
and the stratum consisting of elements whose source curve is still singular
with 2 boundary node.
(In other words $\overline V$ is a neighborhood of this element
in the codimension two stratum
of $\mathcal M_{5}(\beta)$.)
\par
Then any one of the above three gluing constructions
gives a map from an open subset of
$\overline V \times (C,\infty]^2$ onto an open subset of $\mathcal M_{5}(\beta)$.
Let us write them as $\psi_1,\psi_2,\psi_3$, respectively.
The issue is whether $T_p$ and $T_q$ coordinates are preserved
by the coordinate change $\psi_3^{-1}\circ \psi_1$ etc.
In fact, certainly it is {\it not} preserved by this
coordinate change.
\par
This is related to the construction of the collar
of the resulting Kuranishi structure.
Namely if $T_p$, $T_q$ coordinates happen to be preserved
by the coordinate change then we can use this geometric coordinate itself
as a collaring of the corner.
In other words, the neighborhood of the corner (which we denote by $S_2\mathcal M_{5}(\alpha)$)
in $\mathcal M_{5}(\beta)$ is diffeomorphic to
$S_2\mathcal M_{5}(\beta) \times [0,\epsilon)^2$
and the coordinates of the factor $[0,\epsilon)^2$ can be taken, for example,
as $(1/T_p,1/T_q)$.
\par
However it seems not so easy to find a gluing construction
such that $\psi_3^{-1}\circ \psi_1$ etc. preserves $T_{p}$ and $T_q$
coordinates.
\par
On the other hand, there is no need at all to obtain the
collar of the corner {\it directly} by the analytic
construction of the chart.
In the case of single orbifold
existence of the collar can be proved by an easy
standard argument.
In the situation of Kuranishi structure
things are a bit more complicated,
since we need to find collars for various Kuranishi
charts which are preserved by the coordinate change.
Though we can find such a system of collars for
good coordinate system after appropriate shrinking,
its proof is a bit cumbersome to write down
in detail.
\par
Our short cut is to put collar outside,
and for this purpose we need to find an admissible
coordinate system.
\par
For this purpose it suffices to show that $\psi_3^{-1}\circ \psi_1$
preserves gluing parameter
($[0,\epsilon)^2$) modulo an error term
which is exponentially small in $T_p, T_q$.
This is easier than proving that $\psi_3^{-1}\circ \psi_1$
 exactly preserves the gluing parameter.
We can prove this property as follows.
We first observe that though the compatibility of the three different
ways of gluing (1)(2)(3) above does not hold
it is easy to construct pre-gluing for which (1)(2)(3) above are
compatible.
This is because pre-gluing is a simple process
by using partition of unity and pre-gluing on one
neck region does not affect the other neck region.
The exponential decay estimate of the
gluing construction (see \cite{foooexp} for the
detail of the proof) then implies that
actual gluing map is close to
pre-gluing modulo an error
which is exponentially small in $T_p, T_q$.
Therefore the coordinate change has the required properties.

\begin{rem}
To find a collar  related to the gluing
parameter is very much different
type of matter from a similar  issue
to put a collar
appearing in the study of
homotopy of almost complex structures etc..
The later appears, for example, when one proves
independence of the Gromov-Witten invariant of the
choice of almost complex structure.
In the later problem we can take a
homotopy $J_t$ between two almost complex structures
$J$ and $J'$ so that $J_t = J$ for $t \in [0,\epsilon]$
and $J_t = J'$ for $t \in [1-\epsilon,1]$.
So existence of collar is {\it trivial} to prove in that
situation.
The situation is different for the gluing parameter.
(We also note that to prove
independence of the Gromov-Witten invariant of the
choice of almost complex structure
we do {\it not} need to use the
collar. See the proof of
\cite[Proposition 8.16]{part11}.)
\end{rem}

\subsubsection{Stratified submersion}

When we consider a map $f$ from a manifold (an orbifold, a K-space)
with corner $X$ to another manifold $M$ without boundary
or corner,
we say that $f$ is a submersion if
its restriction to all the corners $S_kX$ are submersions.
It implies that the push out $f_{!}h$ of all the smooth forms $h$ on $X$
by $f$ is a smooth form on $M$.
\par
When we study a family of K-systems parametrized by a manifold with corner $P$,
we need to discuss the submersivity of a map
from a manifold (an orbifold, a K-space) with corner $X$
to a manifold with corner $P$.
\par
We use the notion of stratified submersion for such a purpose.
In Section \ref{subsec:ssmaptocorners} we define such a notion
and discuss push out of a differential form to a manifold with corners.

\subsubsection{Integration along the fiber and local system}

As we mentioned in certain situation
(for example when we study the Floer cohomology
of periodic Hamiltonian system in the Bott-Morse situation)
we need to introduce certain $O(1)$ principal
bundle on the space $R_{\alpha}$.
In Part 1 the integration along the fiber
is defined in the situation when our K-space is
oriented.
We need to extend it slightly to include the case
when the target and source spaces come with
$O(1)$ principal
bundle, the K-space (which
gives smooth correspondence) may not be oriened, and
the target and source spaces $R_{\alpha}$ may not be oriented,
but their orientation local systems and the
$O(1)$ principal bundles we put on $R_{\alpha}$ are related by some particular way.
In Section  \ref{sec:lcalsystemtwist} we discuss
such generalization.

\par\bigskip
\noindent
{\bf Convention on the way to use several
notations.}
\begin{enumerate}
\item
[$\widehat{}$ \, and \,$\widetriangle{}$]
\enskip We use `hat' such as $\widehat{\mathcal U}$,
$\widehat{f}$, $\widehat{\frak S}$, $\widehat{h}$
of an object defined on a Kuranishi structure $\widehat{\mathcal U}$.
We use `triangle' such as $\widetriangle{\mathcal U}$,
$\widetriangle{f}$, $\widetriangle{\frak S}$, $\widetriangle{h}$
of an object defined on a good coordinate system $\widetriangle{\mathcal U}$.
\item[$p$ and $\frak p$]
For a Kuranishi structure $\widehat{\mathcal U}$ on $Z\subseteq X$
we write $\mathcal U_p$ for its Kuranishi chart, where
$p \in Z$. (We use an italic letter $p$.)
For a good coordinate system $\widetriangle{\mathcal U}$ on $Z\subseteq X$
we write $\mathcal U_{\frak p}$ for its Kuranishi chart, where
$\frak p \in \frak P$. (We use a German character $\frak p$.) Here $\frak P$ is
a partial ordered set.
\item
[$\blacksquare$]
\enskip
The mark $\blacksquare$ indicates the end of Situation.
See \cite[Situation 6.3]{part11}, for example.
\item[$M$ and $X$] Usually we denote by $M$ a smooth manifold
and by $X$ a K-space, or an orbifold unless otherwise mentioned.
\end{enumerate}

\par
\bigskip
\noindent
{\bf List of Notations in Part 1:}

\begin{enumerate}
\item[$\circ$]
${\rm Int}\, A$, $\ring A$: Interior of a subset $A$ of a topological space.
\item[$\circ$]
$\overline A$: Closure of a subset $A$ of a topological space.
\item[$\circ$]
${\rm Perm}(k)$: The permutation group of order $k!$.
\item[$\circ$]
${\rm Supp}(h)$, ${\rm Supp}(f)$: The support of a differential form $h$,
a function $f$, etc..
\item[$\circ$]
$\varphi^{\star}\mathscr F$: Pull-back of a sheaf $\mathscr F$ by
a map $\varphi$.
\item[$\circ$]
$X$: A paracompact metrizable space. (Part I).
\item[$\circ$]
$Z$: A compact subspace of $X$. (Part I).
\item[$\circ$]
$\mathcal U =(U,\mathcal E,\psi,s)$:
A Kuranishi chart, \cite[Definition 3.1]{part11}.
\item[$\circ$]
$\mathcal U\vert_{U_0} =(U_0,\mathcal E\vert_{U_0},\psi\vert_{U_0\cap s^{-1}(0)},s\vert_{U_0})$:
open subchart of $\mathcal U =(U,\mathcal E,\psi ,s)$, \cite[Definition 3.1]{part11}.
\item[$\circ$]
$\Phi = (\varphi,\widehat\varphi)$:
Embedding of Kuranishi charts, \cite[Definition 3.2]{part11}.
\item[$\circ$]
$o_p, o_p(q)$: Points in a Kuranishi neighborhood $U_p$ of $p$.
\cite[Definition 3.4]{part11}.
\item[$\circ$]
$\Phi_{21} = (U_{21},\varphi_{21},\widehat\varphi_{21})$:
Coordinate change of Kuranishi charts from $\mathcal U_1$ to $\mathcal U_2$,
\cite[Definition 3.5]{part11}.
\item[$\circ$]
$\widehat{\mathcal U} = (\{\mathcal U_p\},\{\Phi_{pq}\})$:
Kuranishi structure, \cite[Definition 3.8]{part11}.
\item[$\circ$]
$(X,\widehat{\mathcal U})$, $(X,Z;\widehat{\mathcal U})$:
K-space, relative K-space, \cite[Definition 3.11]{part11}.
\item[$\circ$]
${\widetriangle{\mathcal U}}
= (({\frak P},\le), \{\mathcal U_{\frak p}\},
\{\Phi_{\frak p\frak q}\})$:
Good coordinate system, \cite[Definition 3.14]{part11}.
\item[$\circ$]
$\vert{\widetriangle{\mathcal U}}\vert$:
\cite[Definition 3.15]{part11}.
\item[$\circ$]
$\widehat\Phi : \widehat{\mathcal U} \to \widehat{\mathcal U'}$:
KK-embedding. An embedding of Kuranishi structures, \cite[Definition 3.20]{part11}.
\item[$\circ$]
$\widetriangle\Phi : \widetriangle{\mathcal U} \to \widetriangle{\mathcal U'}$
: GG-embedding.
An embedding of good coordinate systems, \cite[Definition 3.24]{part11}.
\item[$\circ$]
$\widehat\Phi : \widehat{\mathcal U} \to \widetriangle{\mathcal U}$:
KG-embedding,
An embedding of a Kuranishi structure to a good coordinate system,
\cite[Definition 3.29]{part11}.
\item[$\circ$]
$\widehat\Phi : \widetriangle{\mathcal U} \to\widehat{\mathcal U}$:
GK-embedding.
An embedding of good coordinate system to a Kuranishi structure,
\cite[Definition 5.9]{part11}.
\item[$\circ$]
$\widehat f : (X,Z;\widehat{\mathcal U}) \to Y$ and
$\widetriangle f : (X,Z;\widetriangle{\mathcal U}) \to Y$
:
Strongly continuous map, \cite[Definitions 3.35 and 3.38]{part11}.
\item[$\circ$]
$(X,Z;\widehat{\mathcal U})  \times_{N} M$,
$(X_1,Z_1;\widehat{\mathcal U}_1)  \times_{M}
(X_2,Z_2;\widehat{\mathcal U}_2)$:
Fiber product of Kuranishi structures,
\cite[Definition 4.9]{part11}.
\item[$\circ$]
$S_k(X,Z;\widehat{\mathcal U})$, $S_k(X,Z;\widehat{\mathcal U})$:
Corner structure stratification, \cite[Definition 4.15]{part11}.
\item[$\circ$]
$\mathcal S_{\frak d}(X,Z;\widehat{\mathcal U})$,
$\mathcal S_{\frak d}(X,Z;\widetriangle{\mathcal U})$:
Dimension stratification,
\cite[Definition 5.1]{part11}.
\item[$\circ$]
$\widehat{\mathcal U} < \widehat{\mathcal U^+}$:
$\widehat{\mathcal U^+}$ is a thickening of $\widehat{\mathcal U}$.
\cite[Definition 5.3]{part11}.
\item[$\circ$]
$\mathcal S_{\frak p}(X,Z;{\widetriangle {\mathcal U}};\mathcal K)$:
\cite[Definition 5.6 (4)]{part11}.
\item[$\circ$]
$
\mathcal K = \{\mathcal K_{\frak p}\mid {\frak p
\in \frak P}\}$:
A support system. \cite[Definition 5.6 (1)]{part11}.
\item[$\circ$]
$(\mathcal K^1,\mathcal K^2)$ or
$(\mathcal K^-,\mathcal K^+)$:
A support pair, \cite[Definition 5.6 (2)]{part11}.
\item[$\circ$]
$\mathcal K^1 < \mathcal K^2$: \cite[Definition 5.6]{part11}.
\item[$\circ$]
$\vert\mathcal K\vert$:  \cite[Definition 5.6]{part11}.
\item[$\circ$]
$B_{\delta}(A)$: Metric open ball,
\cite[(6.20)]{part11}.
\item[$\circ$]
$\mathcal S_x = (W_x,\omega_x,\{{\frak s}_x^{\epsilon}\})$:
CF-perturbation (=continuous family perturbation) on one orbifold chart.
\cite[Definition 7.3]{part11}.
\item[$\circ$]
$\mathcal S_x^{\epsilon} = (W_x,\omega_x,{\frak s}_x^{\epsilon})$ for
each $\epsilon >0$:
\cite[Definition 7.3]{part11}.
\item[$\circ$]
$\frak S = \{(\frak V_{\frak r},\mathcal S_{\frak r})\mid{\frak r\in
\frak R}\}$:
Representative of a CF-perturbation on Kuranishi chart
$\mathcal U$.
\cite[Definition 7.15]{part11}.
\par\noindent
Here $\frak V_{\frak r}=(V_{\frak r},E_{\frak r},\Gamma_{\frak r},\phi_{\frak r},\widehat{\phi}_{\frak r})$ is an orbifold chart of $(U,\mathcal E)$ and
$\mathcal S_{\frak r}  = (W_{\frak r} ,\omega_{\frak r}, \{{\frak s}_{\frak r} ^{\epsilon}\})$ is a
CF-perturbation of $\mathcal U$
on $\frak V_{\frak r}$.
\item[$\circ$]
$\frak S^{\epsilon} = \{(\frak V_{\frak r},\mathcal S_{\frak r}^{\epsilon})\mid{\frak r\in
\frak R}\}$
for each $\epsilon >0$. \cite[Definition 7.15]{part11}.
\item[$\circ$]
$\widetriangle{\frak S} = \{\frak S_{\frak p} \mid \frak p \in \frak P\}$:
CF-perturbation of good coordinate system.
\cite[Definition 7.47]{part11}.
\item[$\circ$]
$\widehat{\frak S}$:
CF-perturbation of Kuranishi structure.
\cite[Definition 9.1]{part11}.
\item[$\circ$]
$\mathscr{S}$: Sheaf of CF-perturbations.
\cite[Proposition 7.21]{part11}.
\item[$\circ$]
$
\mathscr S_{\pitchfork 0}
$,
$
\mathscr S_{f \pitchfork}
$, $
\mathscr S_{f \pitchfork g}
$:
Subsheaves of $\mathscr S$. \cite[Definition 7.25]{part11}.
\item[$\circ$]
${\widetriangle f}!(\widetriangle h;\widetriangle{{\frak S}^{\epsilon}})$:
push out or integration along the fiber of $\widetriangle{h}$ by
$(\widetriangle{f},\widetriangle{{\frak S}^{\epsilon}})$
on good coordinate system.
\cite[Definition 7.78]{part11}.

\item[$\circ$]
${\widehat f}!(\widehat h;\widehat{{\frak S}^{\epsilon}})$:
push out or integration along the fiber of $\widehat{h}$ by
$(\widehat{f},\widehat{{\frak S}^{\epsilon}})$
on Kuranishi structure.
\cite[Definition 9.13]{part11}.
\item[$\circ$]
${\rm Corr}_{(\frak X,\widetriangle{{\frak S}^{\epsilon}})}$:
Smooth correspondence associated to good coordinate system.
\cite[Definition 7.85]{part11}.
\item[$\circ$]
${\rm Corr}_{(\frak X,\widehat{\frak S^{\epsilon}})}$:
Smooth correspondence of Kuranishi structure
\cite[Definition 9.23]{part11}.
\item[$\circ$]
$(\frak s_{\frak p}^{n})^{-1}(0)$:
The zero set of multisection.
\item[$\circ$]
$\Pi((\mathfrak S^{\epsilon})^{-1}(0))$:
Support set of a CF-perturbation $\mathfrak S^{\epsilon}$.
\cite[Definition 7.72]{part11}.
\item[$\circ$]
$(V,\Gamma,\phi)$: Orbifold chart, \cite[Definitions 23.1 and 23.6]{part11}.
\item[$\circ$]
$(V,E,\Gamma,\phi,\widehat\phi)$: Orbifold chart of
a vector bundle, \cite[Definitions 23.17 and 23.22]{part11}.
\item[$\circ$]
$(X,\mathcal E)$: Orbibundle, \cite[Definition 23.20]{part11}.
\end{enumerate}

\par
\bigskip
\noindent
\par
\medskip
{\bf List of Notations in Part 2:}

\begin{enumerate}
\item[$\circ$]
${\mathcal M}(\alpha_{-}, \alpha_{+})$:
Space of connecting orbits. Condition \ref{linsysmainconds}.
\item[$\circ$]
$
\mathcal C =\Big(
\frak A, \frak G, \{ R_{\alpha}\}_{\alpha \in \frak A},
\{ o_{R_{\alpha}} \}_{\alpha \in \frak A}, E, \mu,
\{ {\rm PI}_{\beta,\alpha} \}_{\beta \in \frak G, \alpha \in \frak A}
\Big)
$:
A critical submanifold data. Definition \ref{linearsystemdefn}.
\item[$\circ$]
$
\mathcal F = \Big( {\mathcal C},
\{\mathcal M(\alpha_{-},\alpha_{+}) \}_{\alpha_{\pm} \in \frak A},
({\rm ev}_{-}, {\rm ev}_{+}),
\{{\rm OI}_{\alpha_{-}, \alpha_{+}}\}_{\alpha_{\pm} \in \frak A},
\{ {\rm PI}_{\beta;\alpha_{-},\alpha_{+}} \}_{\beta \in \frak G, \alpha_{\pm} \in \frak A}
\Big)
$:
A linear K-system. Definition \ref{linearsystemdefn}.
\item[$\circ$]
$\Lambda^R_{\rm nov}$, $\Lambda^R_{0,\rm nov}$, $\Lambda^R_{+,\rm nov}$:
Universal Novikov ring, and its ideal. Definition \ref{defn:Nov}.
When $R=\R$, we drop $R$ from these notations.
\item[$\circ$]
$\Lambda^R$, $\Lambda^R_{0}$, $\Lambda^R_{+}$:
Universal Novikov ring, and its ideal. (The version without $e$.) Definition \ref{defn:Nov}.
When $R=\R$, we drop $R$ from these notations.
\item[$\circ$]
${\mathcal N}(\alpha_1, \alpha_2)$: Interpolation space. Condition \ref{morphilinsys}.
\item[$\circ$]
$\frak N_{ii+1} ~:~ \mathcal F_i \to \mathcal F_{i+1}$:
Morphism of linear K-systems.
\item[$\circ$]
${\mathcal N}_{i+1i}$:
Interpolation space of the morphism $\frak N_{ii+1}$.
Lemma-Definition \ref{lemdef1434}.
See also Remark \ref{rem:orderindex}.
\item[$\circ$]
$\mathcal {FF}=(\{E^i \}, \{ \mathcal{F}^i \}, \{ \frak N^i \})$:
An inductive system of partial linear K-systems. Definition \ref{defn1528}.
\item[$\circ$]
$V_x^{\boxplus\tau}$: Definition \ref{defn153}.
\item[$\circ$]
$\mathcal U^{\boxplus\tau}_x=(U^{\boxplus\tau}_x,\mathcal E^{\boxplus\tau}_x,\psi^{\boxplus\tau}_x, s^{\boxplus\tau}_x)$:
Trivialization of one Kuranishi chart
$\mathcal U_x$ at $x \in \overset{\circ}S_k(U)$. Lemma-Definition \ref{lem15444}.
See also Lemma \ref{Lema1515}.
\item[$\circ$]
$\mathcal S^{\boxplus\tau}_x$:
Lemma-Definition \ref{lemdef157}.
\item[$\circ$]
$\mathcal U^{\boxplus\tau}$: $\tau$-collaring, or
$\tau$-corner trivialization of
$\mathcal U$.
Lemma-Definition \ref{lemdef157tutuki}.
\item[$\circ$]
$X^{\boxplus\tau}$: $\tau$-collaring, or
$\tau$-corner trivialization of $X$.
Definition \ref{defXenhance}.
\item[$\circ$]
$(X, \widehat{\mathcal U})^{\boxplus\tau} = (X^{\boxplus\tau},\widehat{\mathcal U^{\boxplus\tau}})$:
$\tau$-collaring, or
$\tau$-corner trivialization of K-space $(X, \widehat{\mathcal U})$.
Lemma-Definition \ref{lemdef1522}.
\item[$\circ$]
$(X, \widehat{\mathcal U})^{\boxminus\tau} =
(X^{\boxminus \tau},\widehat{\mathcal U^{\boxminus \tau}})$:
Inward $\tau$-collaring of $(X, \mathcal U)$.
Definition \ref{def:inwardcollar}.
\item[$\circ$]
$
(X,\widehat{\mathcal U})^{\frak C\boxplus\tau}
= (X^{\frak C\boxplus\tau},\widehat{{\mathcal U}^{\frak C\boxplus\tau}})
$:
$\tau$-$\frak C$-corner trivialization, or
partial trivialization of corners, of $(X,\widehat{\mathcal U})$.
Definition \ref{defn1531revrev}.
\item[$\circ$]
$\mathcal N_{12}(\alpha_1,\alpha_2) \times^{\boxplus\tau}_{R^2_{\alpha_2}}
\mathcal N_{23}(\alpha_2,\alpha_3)$:
Partially trivialized fiber product.
Definition \ref{defn1635}.

\item[$\circ$]
$\partial_{\frak C}U$: Normalized $\frak C$-partial boundary of $U$.
When we denote by $\frak C$ a decomposition of
the normalized boundary $\partial U = \partial^0U  \cup \partial^1U$
into two disjoint unions,
we write $\partial_{\frak C}U= \partial^0U$.
Situation \ref{decomporbbdr}.
\item[$\circ$]
$S^{\frak C}_k(U)$:
Definition \ref{defnSCk}.

\item[$\circ$]
$\widehat{\mathcal U_1}^{\boxplus\tau_1} < \widehat{\mathcal U_2}^{\boxplus\tau_2}$ as collared Kuranishi structures:
Proposition \ref{prop161}.

\item[$\circ$]
$\mathcal G(k+1,\beta)$:
The set of all
decorated ribbon trees $(\mathcal T,\beta(\cdot))$
with $(k+1)$ exterior vertices and
$\sum_{{\rm v} \in C_{0,{\rm int}}(\mathcal T)}(\beta({\rm v})) = \beta$.
Definition \ref{defn192}.
\item[$\circ$]
$G(\mathcal{AC})$ (resp. $G(\mathcal{AC}_P)$):
The discrete submonoid associated to an (resp. a $P$-parametrized)
$A_{\infty}$ correspondence $\mathcal{AC}$
(resp. $\mathcal{AC}_P$).
Definition \ref{defassoSubmono}.

\end{enumerate}

Throughout Part 2, an orbifold with corner means an admissible orbifold with
corner in the sense of Subsection \ref{subsec:admord}.
So all notions related to the orbifold with corner are ones in the admissible category.

\section{Linear K-system: Floer cohomology I:
statement}
\label{sec:systemline1}

\subsection{Axiom of linear K-system}
\label{subsec:defnlinesys}

We  axiomatize the properties which are satisfied by the system of moduli spaces of solutions of Floer's equation.

\begin{conds}\label{linsysmainconds}
We consider the following objects.
\par\smallskip
\noindent {\bf (I)}
$\frak G$ is an additive group, and
$E : \frak G \to \R$ and $\mu : \frak G \to \Z$ are
group homomorphisms.
We call $E(\beta)$ the {\it energy}
\index{energy!on $\frak G$} of $\beta$ and
$\mu(\beta)$ the {\it Maslov index}
\index{Maslov index ! on $\frak G$} of $\beta$.
\par\smallskip
\noindent {\bf (II)}
$\frak A$ is a set on which $\frak G$ acts freely.
We assume that the quotient set $\frak A/\frak G$ is a finite set.
$E : \frak A \to \R$ and $\mu : \frak A \to \Z$ are
maps such that
$$
E(\beta \cdot\alpha) = E(\alpha) + E(\beta),
\quad
\mu(\beta \cdot\alpha) = \mu(\alpha) + \mu(\beta)
$$
for any $\alpha \in \frak A$, $\beta \in \frak G$.
We also call $E$ the {\it energy} \index{energy!on $\frak A$} and $\mu$ the {\it Maslov index}. \index{Maslov index ! on $\frak A$}
\par\smallskip
\noindent {\bf (III) (Critical submanifold)}
For any $\alpha \in \frak A$ we have a finite dimensional compact manifold
$R_{\alpha}$ (without boundary), which we call a {\it critical submanifold}.
\par\smallskip
\noindent {\bf (IV) (Connecting orbit)}
For any $\alpha_-, \alpha_+ \in \frak A$, we have a
K-space with corners $\mathcal M(\alpha_-,\alpha_+)$
and strongly smooth maps
$$
({\rm ev}_{-},{\rm ev}_{+}) :
{\mathcal M}({\alpha_-},{\alpha_+})
\to
R_{\alpha_-} \times R_{\alpha_+}.
$$
We assume that
${\rm ev}_+$ is weakly submersive.
We call $\mathcal M(\alpha_-,\alpha_+)$
the {\it  space of connecting orbits}
\index{space of connecting orbits}
and ${\rm ev}_{\pm}$ the {\it evaluation maps at infinity}.
\index{evaluation maps at infinity ! evaluation maps at infinity}
\par\smallskip
\noindent {\bf (V) (Positivity of energy)}
We assume  $
{\mathcal M}({\alpha_-},{\alpha_+}) = \emptyset$
if $E(\alpha_-) \ge E(\alpha_+)$.
\footnote{We note that ${\mathcal M}({\alpha},{\alpha}) = \emptyset$
in particular.}
\par\smallskip
\noindent {\bf (VI) (Dimension)}
The dimension of the  space of connecting orbits is given by
\begin{equation}\label{eq:dim1}
\dim
{\mathcal M}({\alpha_-},{\alpha_+})
=
\mu(\alpha_+) - \mu(\alpha_-) - 1 + \dim R_{\alpha_+}.
\end{equation}
\par\smallskip
\noindent{\bf (VII) (Orientation)}
\begin{enumerate}
\item[{\bf (i)}]
For any $\alpha\in \frak A$, a principal $O(1)$ bundle $o_{R_{\alpha}}$
on $R_{\alpha}$
is given.
We call it an {\it orientation system of the critical submanifold}.
\index{orientation system of critical submanifold}
\item[{\bf (ii)}]
For any $\alpha_1, \alpha_2 \in \frak A$, an isomorphism
\begin{equation}\label{orientiso}
{\rm OI}_{\alpha_{-},\alpha_{+}} :
{\rm ev}_{-}^*(o_{R_{\alpha_{-}}})
\cong
{\rm ev}_{+}^*(o_{R_{\alpha_{+}}}) \otimes
{\rm ev}_{+}^*(\det TR_{\alpha_{+}})
\otimes o_{{\mathcal M}(\alpha_-,\alpha_+)}
\end{equation}
of principal $O(1)$ bundles is given\footnote{
In the case of the linear K-system obtained from periodic Hamiltonian system,
we take $o_R=\Theta_R^{-}$ for each critical submanifold $R$,
where $\Theta_R^{-}$ is defined as the determinant of the index bundle
of certain family of elliptic operators.
See \cite[Definition 8.8.2]{fooobook2} for the precise definition.
Then
\cite[Proposition 8.8.6]{fooobook2}
yields the isomorphism
\eqref{orientiso}.
On the other hand, in \cite[Proposition 8.8.7]{fooobook2} we
take $o_{R}=\det TR \otimes \Theta_R^{-}$.
Note that we use singular chains in \cite{fooobook2},
while
we use differential forms in the current manuscript.
So the choices of $o_R$ are slightly different.}.
Here $o_{{\mathcal M}(R_{\alpha_-},R_{\alpha_+})}$
is an orientation bundle
of the K-space ${\mathcal M}(\alpha_-,\alpha_+)$
in the sense of \cite[Definition 3.10]{part11}.
We call the isomorphism (\ref{orientiso}) the
{\it orientation isomorphism}\footnote{See
\cite[Section 8.8]{fooobook2}.}.
(See Section \ref{sec:lcalsystemtwist}.
Otherwise the reader may consider only the case when all the spaces
$R_{\alpha}$ and ${\mathcal M}({\alpha_-},{\alpha_+})$ are
oriented and $o_{R_{\alpha_{\pm}}}$ are trivial.)
\footnote{Note that $o_{R_{\alpha}}$ may not
coincide with the principal $O(1)$ bundle
giving an orientation of $R_{\alpha}$.
For example,
$o_{R_{\alpha}}$ may be nontrivial even in the case when $R_{\alpha}$ is orientable.
On the other hand, if  ${\mathcal M}(\alpha_-,\alpha_+)$
is orientable then $o_{{\mathcal M}(\alpha_-,\alpha_+)}$
is trivial.}
More precisely, we fix a choice of homotopy class
of the isomorphism (\ref{orientiso}).
\end{enumerate}
\par\smallskip
\noindent{\bf (VIII) (Periodicity)}
\begin{enumerate}
\item[{\bf (i)}]
For any $\beta \in \frak G$  a diffeomorphism
\begin{equation}\label{formula1477}
{\rm PI}_{\beta;\alpha} : R_{\alpha} \to R_{\beta\alpha}
\end{equation}
is given such that the equality
\begin{equation}\nonumber
{\rm PI}_{\beta_2;\beta_1\alpha}\circ {\rm PI}_{\beta_1;\alpha}
=
{\rm PI}_{\beta_2\beta_1;\alpha}
\end{equation}
holds.
\item[{\bf (ii)}]
Moreover, an isomorphism
\begin{equation}\label{form128}
{\rm PI}_{\beta;\alpha_-,\alpha_+}:{\mathcal M}({\alpha_-},{\alpha_+})
\to
{\mathcal M}({\beta\alpha_-},{\beta\alpha_+})
\end{equation}
of K-spaces in the sense of \cite[Definition 4.22]{part11} is given such that
the equality
$$
{\rm PI}_{\beta_2;\beta_1\alpha_-,\beta_1\alpha_+}\circ {\rm PI}_{\beta_1;\alpha_-,\alpha_+}
=
{\rm PI}_{\beta_2\beta_1;\alpha_-,\alpha_+}
$$
holds.
The diagram below commutes.
\end{enumerate}
\begin{equation}
\begin{CD}
{\mathcal M}(\alpha_-, \alpha_+) @ > {{\rm PI}_{\beta;\alpha_-,\alpha_+}} >>
{\mathcal M}({\beta\alpha_-},{\beta\alpha_+})  \\
@ V{({\rm ev}_{-},{\rm ev}_{+})}VV @ VV{({\rm ev}_{-},{\rm ev}_{+})}V\\
R_{\alpha_-}\times R_{\alpha_+} @ >  {({\rm PI}_{\beta;\alpha_-},{\rm PI}_{\beta;\alpha_+})} >>
R_{\beta\alpha_-}\times R_{\beta\alpha_+}
\end{CD}
\end{equation}
We call ${\rm PI}_{\beta;\alpha}$, ${\rm PI}_{\beta;\alpha_-,\alpha_+}$
the {\it periodicity isomorphisms}.
\index{periodicity isomorphism ! periodicity isomorphism}
The periodicity isomorphism preserves $o_{R_{\alpha}}$ and
commutes with the orientation isomorphism.
\par\smallskip
\noindent {\bf (IX) (Gromov compactness)}
\index{Gromov compactness ! Gromov compactness}
For any $E_0 \ge 0$ and $\alpha_- \in \frak A$ the set
\begin{equation}
\{\alpha_+ \in \frak A \mid
{\mathcal M}({\alpha_-},{\alpha_+})
\ne \emptyset,
\,\, E(\alpha_+) \le E_0 + E(\alpha_-)\}
\end{equation}
is a finite set.
\par\smallskip
\noindent{\bf (X) (Compatibility at the boundary)}
\index{compatibility at the boundary !compatibility at the boundary}
The normalized boundary of the space of connecting orbits
is decomposed into a disjoint union of fiber products:\footnote{
The right hand side is a finite union by Condition (IX).}
\begin{equation}\label{formula1211_0}
\partial{\mathcal M}({\alpha_-},{\alpha_+})
\cong \coprod_{\alpha}
\left(
{\mathcal M}({\alpha_-},{\alpha})
\,\,{}_{{\rm ev}_{+}}\times_{{\rm ev}_{-}}\,\,
{\mathcal M}({\alpha},{\alpha_+})
\right).
\end{equation}
Here $\cong$ means an isomorphism of K-spaces.
When we compare the orientations between the two sides,
we swap the order of the factors in the fiber product in the right hand side of
\eqref{formula1211_0}
and put the sign as follows:
\begin{equation}\label{formula1211}
\partial{\mathcal M}({\alpha_-},{\alpha_+})
\cong \coprod_{\alpha}
(-1)^{\dim \mathcal M (\alpha,\alpha_+)}\left(
{\mathcal M}({\alpha},{\alpha_+})
\,\,{}_{{\rm ev}_{-}}\times_{{\rm ev}_{+}}\,\,
{\mathcal M}({\alpha_-},{\alpha})
\right).
\end{equation}
See Remark \ref{rem:FiberProdOrd} for this isomorphism.
This isomorphism commutes with the periodicity and
orientation isomorphisms and
is compatible with various evaluation maps.
Namely the restriction of
${\rm ev}_-$ (resp. ${\rm ev}_+$) of the left hand side of
(\ref{cornecom1}) coincides with
${\rm ev}_-$ (resp. ${\rm ev}_+$)
of the factor ${\mathcal M}({\alpha_-},{\alpha}_1)$
(resp. ${\mathcal M}({\alpha}_k,{\alpha_+})$) of the right hand side.

\begin{rem}\label{rem:FiberProdOrd}
(1)
The sign in \eqref{formula1211} is consistent with that of
\cite[p.728]{fooobook2} (the line 6 of the proof of Theorem 8.8.10 (3)).
Here we should take into account the following two points about the notation and convention.
The first point is about notation.
In \cite[Chapter 8]{fooobook2} we write $\mathcal M (\alpha,\beta)$ for
$\lim_{\tau \to -\infty} u(\tau ,t) \in R_{\beta},
\lim_{\tau \to +\infty} u(\tau ,t) \in R_{\alpha}$ for the discussion on orientations
of moduli spaces of connecting orbits.
Actually, this notation is different from those used in
other chapters of \cite{fooobook, fooobook2}, as we note at the beginning of Section 8.7 and p.723 in \cite{fooobook2}.
The order of the position of $\alpha ,\beta$ in the notation $\mathcal M (\alpha, \beta)$
used in this article is opposite to the one used in \cite[Chapter 8]{fooobook2}.
This is just a difference of notations.
\par
The second point is the order of factors in the fiber product.
This is not only a difference in conventions but also
really affects the orientation on the fiber product.
In \cite[Chapter 8]{fooobook2}, e.g. Section 8.3,
we use the evaluation map at the $0$-th marked point of {\it the second factor}
when taking the fiber product, while we are using the evaluation map
at the $0$-th marked point of {\it the first factor} in this book.
See \eqref{formula1211_0}, for example.
(As we note in Convention 8.3.1 and in the third line of \cite[p.699]{fooobook2},
we regard the $0$-th marked point $z_0$ as $+1$ in the unit disk in $\C$ which
corresponds to $+\infty$ in the strip $\R \times [0,1]$.)
In \eqref{formula1211}
we follow the convention used in \cite[Chapter 8]{fooobook2}.
That is why we swap the factors in the right hand side in
\eqref{formula1211_0}.
Now
by taking the difference of notation mentioned here
into account,
the formula in the line 6 of the proof of \cite[Theorem 8.8.10 (3), p.728]{fooobook2}
can be rewritten with the notation in this article as
\begin{equation}\label{book2p728}
\aligned
& \partial \mathcal M_{k+\ell +2}(R_{h_1}, R_{h_3}) \supset
(-1)^{d_0} \mathcal M_{k+2} (R_{h_2}, R_{h_3}) \times_{R_{h_2}}
\mathcal M_{\ell +2}(R_{h_1}, R_{h_2}),  \\
& d_0 = k\ell +k (\mu(h_2)-\mu(h_1)) + (\mu (h_3) -\mu (h_2)) + k -1 + \dim R_{h_3}.
\endaligned
\end{equation}
Applying this formula for the case $\alpha_{-} =h_1, \alpha = h_2, \alpha_{+} =h_3$ and $k=\ell=0$, we obtain \eqref{formula1211} by
the dimension formula \eqref{eq:dim1}.
\par
{\it Whenever we have to explore the orientation,
we swap the factors in the fiber product to follow the convention used
in \cite[Chapter 8]{fooobook2}.}
We will mention it where such places appear later.
On the other hand,
when we do not need to study orientations,
e.g. at (XI) (XII) compatibility at the corner below,
we do not swap the order of factors in fiber products.
See the footnote therein.
\par
(2)
In general, since we find
\begin{equation}\label{signcommute}
X_1\times_Y X_2 =(-1)^{(\dim X_1 -\dim Y)(\dim X_2 -\dim Y)} X_2 \times_Y X_1 ,
\end{equation}
we can rewrite \eqref{formula1211} as
\begin{equation}\label{formula1211_1}
\partial{\mathcal M}({\alpha_-},{\alpha_+})
\cong \coprod_{\alpha}
(-1)^{\epsilon}\left(
{\mathcal M}({\alpha_-},{\alpha})
\,\,{}_{{\rm ev}_{+}}\times_{{\rm ev}_{-}}\,\,
{\mathcal M}({\alpha},{\alpha_+})
\right)
\end{equation}
where
$$
\aligned
\epsilon & =
\dim \mathcal M (\alpha,\alpha_+) +
(\dim \mathcal M (\alpha_-,\alpha) -\dim R_{\alpha})
(\dim \mathcal M (\alpha,\alpha_+) -\dim R_{\alpha}) \\
& =
\left(\mu(\alpha)-\mu(\alpha_-)\right)
\left(\mu(\alpha_+)-\mu(\alpha) -1 + \dim R_{\alpha_+} -\dim R_{\alpha}\right)
+ \dim R_{\alpha} \\
& = \left( \mu(\alpha) -\mu(\alpha_-) \right) \dim \mathcal M (\alpha,\alpha_+)
- \left( \mu(\alpha) -\mu(\alpha_-) -1 \right) \dim R_{\alpha}.
\endaligned
$$
\end{rem}
\par
\noindent {\bf (XI) (Compatibility at the corner I)}
\index{compatibility at the corner ! compatibility at the corner}
Let $\widehat{S}_k(\mathcal M(\alpha_-,\alpha_+))$ be the
normalized corner of the K-space $\mathcal M(\alpha_-,\alpha_+)$
in the sense of Definition \ref{norcor}.
Then we have:
\begin{equation}\label{cornecom1}
\aligned
&\widehat S_k({\mathcal M}({\alpha_-},{\alpha_+}))\\
&\cong
\coprod_{\alpha_1,\dots,\alpha_k
\in \frak A}
\left(
{\mathcal M}({\alpha_-},{\alpha}_1)
\,\,{}_{{\rm ev}_{+}}\times_{R_{\alpha_1}}\,\cdots
\,\,{}_{R_{\alpha_k}}\times_{{\rm ev}_{-}}
{\mathcal M}({\alpha}_k,{\alpha_+})
\right).
\endaligned
\end{equation}
Here the right hand side is a disjoint union.
This is an isomorphism of K-spaces.
It preserves the periodicity isomorphisms.\footnote{We do not assume any
compatibility of
the orientation isomorphism {\it at the corner}, because what we will use is Stokes' formula
where boundary but not corner appears.}
It is compatible with various evaluation maps.
Namely the restriction of
${\rm ev}_-$ (resp. ${\rm ev}_+$) of the left hand side of
(\ref{cornecom1}) coincides with
${\rm ev}_-$ (resp. ${\rm ev}_+$)
of the factor ${\mathcal M}({\alpha_-},{\alpha}_1)$
(resp. ${\mathcal M}({\alpha}_k,{\alpha_+})$) of the right hand side.
\par\smallskip
\noindent {\bf (XII) (Compatibility at the corner II)}
The isomorphism in (XI) satisfies the compatibility condition given in Condition \ref{corneraddamoc} below.
\end{conds}
To describe Condition \ref{corneraddamoc} we need some notation.
By (\ref{form3077}) and (\ref{cornecom1}) we have
\begin{equation}
\aligned
&\widehat S_{\ell}(\widehat S_k({\mathcal M}({\alpha_-},{\alpha_+})))\\
&\cong
\coprod_{\alpha_1,\dots,\alpha_k \in \frak A}
\coprod_{\ell_0+\dots+\ell_k = \ell}
\left(
\widehat S_{\ell_0}{\mathcal M}({\alpha_-},{\alpha}_1)
\,\,{}_{{\rm ev}_{+}}\times_{R_{\alpha_1}}\,\cdots
\,\,{}_{R_{\alpha_k}}\times_{{\rm ev}_{-}}
\widehat S_{\ell_k}{\mathcal M}({\alpha}_k,{\alpha_+})
\right).
\endaligned
\nonumber\end{equation}
We apply (\ref{cornecom1}) to each of the fiber product factors of the
right hand side and obtain
\begin{equation}\label{1413form}
\aligned
&\widehat S_{\ell}(\widehat S_k({\mathcal M}({\alpha_-},{\alpha_+})))\\
&\cong
\coprod_{\alpha_1,\dots,\alpha_k \in \frak A}
\coprod_{\ell_0+\dots+\ell_k = \ell}
\coprod_{\alpha_{0,1},\dots,\alpha_{0,\ell_0}}
\dots
\coprod_{\alpha_{k,1},\dots,\alpha_{k,\ell_k}}
\\
&\left(
{\mathcal M}({\alpha_-},{\alpha}_{0,1})
\,\,{}_{{\rm ev}_{+}}\times_{R_{\alpha_{0,1}}}\,\cdots
\,\,{}_{R_{\alpha_{0,\ell_0}}}\times_{{\rm ev}_{-}}
{\mathcal M}({\alpha}_{0,\ell_0},{\alpha_1})
\right) \\
&
\,\,{}_{{\rm ev}_{+}}\times_{R_{\alpha_{1}}}\,
\left(
{\mathcal M}({\alpha_1},{\alpha}_{1,1})
\,\,{}_{{\rm ev}_{+}}\times_{R_{\alpha_{1,1}}}\,\cdots
\,\,{}_{R_{1,\alpha_{\ell_k}}}\times_{{\rm ev}_{-}}
{\mathcal M}({\alpha}_{1,\ell_k},{\alpha_2})
\right)\\
&
\qquad\qquad\qquad\qquad\qquad\cdots
\\
&
\,\,{}_{{\rm ev}_{+}}\times_{R_{\alpha_{k}}}\,
\left(
{\mathcal M}({\alpha_k},{\alpha}_{k,1})
\,\,{}_{{\rm ev}_{+}}\times_{R_{\alpha_{k,1}}}\,\cdots
\,\,{}_{R_{k,\alpha_{\ell_k}}}\times_{{\rm ev}_{-}}
{\mathcal M}({\alpha}_{k,\ell_k},{\alpha_+})
\right).
\endaligned
\end{equation}
By applying (\ref{cornecom1}) to $k+\ell$ in place of $k$ we obtain
\begin{equation}\label{cornecom1revrev}
\aligned
&\widehat S_{k+\ell}({\mathcal M}({\alpha_-},{\alpha_+}))\\
&\cong
\coprod_{\alpha_1,\dots,\alpha_{k+\ell}
\in \frak A}
\left(
{\mathcal M}({\alpha_-},{\alpha}_1)
\,\,{}_{{\rm ev}_{+}}\times_{R_{\alpha_1}}\,\cdots
\,\,{}_{R_{\alpha_{k+\ell}}}\times_{{\rm ev}_{-}}
{\mathcal M}({\alpha}_{k+\ell},{\alpha_+})
\right).
\endaligned
\end{equation}
It is easy to observe that each summand of
(\ref{1413form}) appears in (\ref{cornecom1revrev})
and vice versa.
\begin{conds}\label{corneraddamoc}
The covering map $\pi_{\ell,k} : \widehat S_{\ell}(\widehat S_k({\mathcal M}({\alpha_-},{\alpha_+}))) \to \widehat S_{k+\ell}({\mathcal M}({\alpha_-},{\alpha_+}))$
in Proposition \ref{prop2813}
restricts to the identity map on each summand of (\ref{1413form}).
\end{conds}
\begin{rem}
It is easy to see that each summand in (\ref{cornecom1revrev})
appears exactly $(k+\ell)!/k!\ell!$ times in (\ref{1413form}).
This is the covering index of the map $\pi_{\ell,k}$.
\end{rem}
\begin{rem}
In general, in the case of orbifold, the map
$\overset{\circ}S_{k-1}(\partial U )\to \overset{\circ}S_{k}(U)$ is not $k$ to one map
set-theoretically. (See \cite[Remark 8.6 (3)]{part11}.) However it is so in our case since the isotropy group acts trivially
on the part $[0,1)^k$ (which is the normal direction to the stratum) in our case.
\end{rem}
\begin{defn}\label{linearsystemdefn}
\begin{enumerate}
\item
Let $\frak A$, $\frak G$, $R_{\alpha}$, $o_{R_{\alpha}}$, $E$, $\mu$,
${\rm PI}_{\beta,\alpha}$ be the objects as in
Conditions \ref{linsysmainconds} (I), (II), (III), (VII)-(i),
(VIII)-(i).\footnote{
Namely, we collect the same data as in Conditions \ref{linsysmainconds}
as far as critical submanifolds concern.}
We denote them by
$$
\mathcal C =\Big(
\frak A, \frak G, \{ R_{\alpha}\}_{\alpha \in \frak A},
\{ o_{R_{\alpha}} \}_{\alpha \in \frak A}, E, \mu,
\{ {\rm PI}_{\beta,\alpha} \}_{\beta \in \frak G, \alpha \in \frak A}
\Big)
$$
and call $\mathcal C$ a {\it critical submanifold data}
\index{critical submanifold data}.
\item
Together with a critical submanifold data $\mathcal C$,
a {\it linear system of spaces with Kuranishi structures},
abbreviated as a {\it linear K-system},
\index{K-system ! linear K-system}
\index{linear system of Kuranishi structures ! {\it see: K-system}}
$\mathcal F$\footnote{We
use $\mathcal F$ to denote a linear K-system.
Here F stands for Floer.}
consists of objects
$$
\mathcal F = \Big( {\mathcal C},
\{\mathcal M(\alpha_{-},\alpha_{+}) \}_{\alpha_{\pm} \in \frak A},
({\rm ev}_{-}, {\rm ev}_{+}),
\{{\rm OI}_{\alpha_{-}, \alpha_{+}}\}_{\alpha_{\pm} \in \frak A},
\{ {\rm PI}_{\beta;\alpha_{-},\alpha_{+}} \}_{\beta \in \frak G, \alpha_{\pm} \in \frak A}
\Big)
$$
which satisfy Condition \ref{linsysmainconds}.
\item
Let $E_0 > 0$. A {\it partial linear K-system}
\index{K-system ! partial linear K-system (of energy cut level $E_0$)}
consists of the same objects as in Condition \ref{linsysmainconds}
except the following points.
We call $E_0$ its {\it energy cut level}.
\index{energy cut level!energy cut level}
\begin{enumerate}
\item
The  K-space (the space of connecting orbits)
${\mathcal M}({\alpha_-},{\alpha_+})$ is
defined only when $E({\alpha_+}) - E({\alpha_-}) \le E_0$.
\item
Periodicity and orientation isomorphisms of the space of connecting orbits
are defined only among those satisfying $E({\alpha_+}) - E({\alpha_-}) \le E_0$.
\item
Compatibility isomorphisms of Kuranishi structures at the boundary
given in Condition \ref{linsysmainconds} (X)
are defined only when the left hand side of (\ref{formula1211}) is defined.
\item
Compatibility isomorphism of Kuranishi structures  at the corners in Condition \ref{linsysmainconds} (XI) (XII)
are defined only for ${\mathcal M}({\alpha_-},{\alpha_+})$ with  $E({\alpha_+}) - E({\alpha_-}) \le E_0$.
\end{enumerate}
\end{enumerate}
\end{defn}
\begin{rem}
\begin{enumerate}
\item
We define the notion of partial linear K-system
to take care of the `running out problem'
\index{running out problem} which we discussed in
\cite[Section 7.2.3]{fooobook2}.
(See Section \ref{sec:systemline3} for the way how we will use it.)
We can use a similar argument also to define symplectic homology
(see \cite{FlHofer}, \cite {cieFlHofer}, \cite{bourgeooancea}) in the case $X$ is noncompact but convex at infinity, in
complete generality.
One difference between the current context
and the context of symplectic homology lies in the finiteness statement given in Condition \ref{linsysmainconds} (II).
\item
To construct Floer cohomology in the situation where we have only
partial linear K-systems, we need to
define the notion of
{\it inductive systems of partial linear K-systems}.
To define the notion of such an inductive system,
we need the notion of morphism of linear K-systems.
We will define it in Definition \ref{linearsystemmorphdefn}.
\end{enumerate}
\end{rem}

\subsection{Floer cohomology associated to a linear K-system}
\label{subsec:floerhomo}

To state our main theorem on linear K-systems,
we need to prepare some notations.
\begin{defn}\label{Fvectspace}
Let $\mathcal C$ be a critical submanifold data as in
Definition \ref{linearsystemdefn}.
\begin{enumerate}
\item
We put
\begin{equation}\label{Floermodule0}
CF(\mathcal C)^0 = \bigoplus_{\alpha \in \frak A}
\Omega(R_{\alpha};o_{R_{\alpha}}),
\end{equation}
where $\Omega(R_{\alpha};
o_{R_{\alpha}})$ is the de Rham complex of $R_{\alpha}$ twisted by the principal $O(1)$ bundle $o_{R_{\alpha}}$. (See Section \ref{sec:lcalsystemtwist}.)
\item
$CF(\mathcal C)^0$ is a graded filtered $\R$ vector space.
Its grading is given so that
the degree $d$ part $CF^d(\mathcal C)$ is to be
\index{degree ! of Floer cochain complex}
$$
CF^d(\mathcal C)^0 = \bigoplus_{\alpha \in \frak A} \Omega^{d-\mu(\alpha)}(R_{\alpha};o_{R_{\alpha}}),
$$
and its filtration $\frak F CF^d(\mathcal C)^0$
\index{filtration ! of Floer cochain complex} is given by
$$
\frak F^{\lambda}CF(\mathcal C)^0 = \bigoplus_{\alpha \in \frak A
\atop E(\alpha) \ge \lambda} \Omega(R_{\alpha};o_{R_{\alpha}}).
$$
Here $\lambda \in \R$.
\item
The completion of $CF(\mathcal C)^0$
with respect to the filtration $\frak F^{\lambda}CF(\mathcal C)^0$ is
denoted by $CF(\mathcal C)$.
It is a graded filtered $\R$ vector space which is complete.
Its element corresponds to an infinite formal sum
$
\sum_{i=1}^{\infty} h_i,
$
where
\begin{enumerate}
\item $h_i \in \Omega(R_{\alpha_i};o_{R_{\alpha_i}})$,
\item $\alpha_i \in \frak A$,
\item
If  $i \le j$ then $E(\alpha_i) \le E(\alpha_j)$,
\item
$\lim_{i\to \infty}E(\alpha_i) = \infty$ unless $\{\alpha_i\}$ is a finite set.
\end{enumerate}
\item
For each $\beta\in \frak G$, the inverse of the periodicity diffeomorphism
${\rm PI}_{\beta;\alpha}$
induces an isomorphism $\Omega(R_{\alpha};o_{R_{\alpha}}) \to \Omega(R_{\beta\alpha};o_{R_{\beta\alpha}})$.
Its direct sum extends to the above mentioned completion which we call
{\it periodicity isomorphism}
\index{periodicity isomorphism ! of Floer cochain complex} and write
$$
{\rm PI}_{\beta}^* : CF(\mathcal C) \to CF(\mathcal C).
$$
It is of degree $\mu(\beta)$ and satisfies
\begin{equation}
{\rm PI}_{\beta}^*(\frak F^{\lambda}CF(\mathcal C))
\subset \frak F^{\lambda+E(\beta)}CF(\mathcal C).
\end{equation}
\item
We define $d_0 : CF(\mathcal C) \to CF(\mathcal C)$
by
\begin{equation}\label{eq:deRham0}
d_0 \left(\sum_{i=1}^{\infty} h_i\right)
=
\sum_{i=1}^{\infty} (-1)^{\dim R_{\alpha_i} + \mu(\alpha_i) +1 + \deg h_i}d_{dR}h_i
\end{equation}
where $h_i \in \Omega^{\deg h_i}(R_{\alpha_i};o_{R_{\alpha_i}})$
and $d_{dR}$ on the right hand side is the de Rham differential.
(See \cite[p.150]{fooobook} for the exponent $\dim R_{\alpha_i} + \mu(\alpha_i)$ of the sign
and \cite[Remark 3.5.8]{fooobook} for the exponent $1 + \deg h_i$.)
The map
$d_0$ has degree $1$ and preserves the filtration.
\end{enumerate}
We call $CF(\mathcal C)$ together with differential
$d$ given in Theorem \ref{linesysmainth1}, the
{\it Floer cochain complex associated to $\mathcal C$}.
\index{Floer cochain complex ! of
linear K-system}
\end{defn}
\begin{thm}\label{linesysmainth1}
Suppose we are in the situation of Definition \ref{Fvectspace}.
\begin{enumerate}
\item
For any linear K-system
we can define a map $d : CF(\mathcal C) \to CF(\mathcal C)$
such that
\begin{enumerate}
\item $d \circ d = 0$.
\item  $d$ has degree $1$ and preserves the filtration.
\item $d$ commutes with the periodicity isomorphism
${\rm PI}^{\ast}_{\beta}$.
\item
There exists $\epsilon > 0$ such that:
$$
(d-d_0)(\frak F^{\lambda}CF(\mathcal C))
\subset \frak F^{\lambda+\epsilon}CF(\mathcal C).
$$
\end{enumerate}
\item
The definition of the map $d$ in (1) involves various choices
related to the associated Kuranishi structure and
$d$ depends on them. However it is independent of such choices
in the following sense.
\par
If $d_1$, $d_2$ are obtained from two different choices,
there exists $\psi : CF(\mathcal C) \to CF(\mathcal C)$ with the
following properties.
\begin{enumerate}
\item $d_2 \circ \psi = \psi \circ d_1$.
\item $\psi$ has degree $0$ and preserves the filtration.
\item $\psi$ commutes with the periodicity isomorphism
${\rm PI}^{\ast}_{\beta}$.
\item
There exists $\epsilon > 0$ such that:
$$
(\psi-{\rm id})(\frak F^{\lambda}CF(\mathcal C))
\subset \frak F^{\lambda+\epsilon}CF(\mathcal C)
$$
where ${\rm id}$ is the identity map.
\item
In particular, $\psi$ induces an isomorphism
in cohomologies:
$$
H(CF(\mathcal C),d_1)
\to H(CF(\mathcal C),d_2).
$$
\item
The cochain map $\psi$ itself depends on various choices. However
it is independent of the choices up to cochain homotopy.
Namely if $\psi_i$ $i=1,2$ are obtained from the different choices,
there exists a map $H$ of degree $-1$ such that
$$
\psi_1 - \psi_2 = d_2\circ H + H\circ d_1.
$$
Moreover there exists $\epsilon > 0$ such that:
$$
H(\frak F^{\lambda}CF(\mathcal C))
\subset \frak F^{\lambda+\epsilon}CF(\mathcal C).
$$
\end{enumerate}
\end{enumerate}
\end{thm}
The proof will be given in Section \ref{sec:systemline3}.
\par
In the independence statement such as Theorem \ref{linesysmainth1} (2)
the critical submanifolds $R_{\alpha}$ are fixed.
For most of the important applications of Floer cohomology
we need to prove certain independence statement
in the situation where the critical submanifolds vary.
In that case we need to take an appropriate Novikov ring as the relevant coefficient
ring.
We introduce the following universal Novikov ring for this purpose.
\begin{defn}\label{defn:Nov}
Let $R$ be a commutative ring with unity:
\begin{enumerate}
\item
The set $\Lambda^R_{\rm nov}$
consists of the formal sums
\begin{equation}\label{elementofNov}
\frak x = \sum_{i=1}^{\infty} P_i(e) T^{\lambda_i}
\end{equation}
where
\begin{enumerate}
\item
$P_i(e) \in R[e^{1/2},e^{-1/2}]$.
\footnote{In Item (5) we put $\deg e=2$. The Novikov ring $\Lambda_{\rm nov}^R$ here is the same as
one introduced in \cite{fooobook}, where the indeterminate $e$ has degree $2$.}
\item $T$ is a formal parameter.
\item
$\lambda_i \in \R$ and
$\lambda_i < \lambda_j$ for $i < j$.
\item
$\lim_{i\to \infty} \lambda_i = \infty$ unless \eqref{elementofNov} is a finite sum.
\end{enumerate}
\item
We define the norm of the element (\ref{elementofNov}) by
$$
\Vert \frak x\Vert = \exp \left( - \inf \{ \lambda_i \mid P_i(e) \ne 0\}\right).
$$
We put $\Vert 0\Vert = 0$.
We define a filtration on $\Lambda^R_{\rm nov}$
by
$$
\frak F^{\lambda}\Lambda^R_{\rm nov}
= \{
\frak x \in \Lambda^R_{\rm nov}
\mid \Vert \frak x\Vert
\le e^{-\lambda}\}.
$$
\item
Let us consider the subset of $\Lambda^R_{\rm nov}$
such that $P_i(e) = 0$ except finitely many indices $i$.
We define a ring structure on it in an obvious way.
It is also an $\R$-vector space
and $\Vert \frak x\Vert$ is a norm of it.
$\Lambda^R_{\rm nov}$  is its completion with respect to
this norm.
The ring structure extends to $\Lambda^R_{\rm nov}$  and
$\Lambda^R_{\rm nov}$  becomes a complete normed ring.
We call $\Lambda^R_{\rm nov}$ the {\it universal Novikov ring}
\index{Novikov ring ! $\Lambda^R_{\rm nov}$}
with {\it ground ring} $R$.
\item
The subset consisting of elements $\frak x\in \Lambda^R_{\rm nov}$
with $\Vert \frak x\Vert \le 1$ is a subring,
which we write $\Lambda^R_{0,\rm nov}$ and call it also the
{\it universal Novikov ring}.
\index{Novikov ring ! $\Lambda^R_{0,\rm nov}$}
Moreover, we put
$\Lambda^R_{+,\rm nov} = \{\frak x\in \Lambda^R_{\rm nov} ~\vert~
\Vert \frak x\Vert < 1 \}$.
\index{Novikov ring ! $\Lambda^R_{+,\rm nov}$}
\item
We put $\deg T = 0$, $\deg e = 2$ and the degree of elements of $R$ to be $0$.
Then $\Lambda^R_{0,\rm nov}$ and $\Lambda^R_{\rm nov}$
are graded rings.
\item
In case $R = \R$ we omit $\R$ and write $\Lambda_{\rm nov}$,
$\Lambda_{0,\rm nov}$ and $\Lambda_{+,\rm nov}$.
\item
When we do not include the indeterminate $e$, we write
$\Lambda^R$,
$\Lambda^R_{0}$ and $\Lambda^R_{+}$
in place of
$\Lambda^R_{\rm nov}$,
$\Lambda^R_{0,\rm nov}$ and $\Lambda^R_{+,\rm nov}$.
When $R=\R$, we also drop $\R$ from these notations.
\end{enumerate}
\end{defn}
We now start from the cochain complex obtained in Theorem \ref{linesysmainth1}
and obtain a cochain complex over the
universal Novikov ring $\Lambda_{0,\rm nov}$.
\begin{defn}\label{defn141212}
Let $\mathcal C$ be a critical submanifold data as in Definition
\ref{linearsystemdefn}.
\begin{enumerate}
\item
We consider $CF(\mathcal C)^0$ as in (\ref{Floermodule0})
and take an algebraic tensor product with ${\Lambda_{\rm nov}}$
over $\R$, that is,
$
CF(\mathcal C)^0 \otimes_{\R}{\Lambda_{\rm nov}}
$.
It is a ${\Lambda_{\rm nov}}$ module.
The ${\Lambda_{\rm nov}}$ module $CF(\mathcal C)^0 \otimes_{\R}{\Lambda_{\rm nov}}$ is graded by
$
\deg (x\otimes \frak x) = \deg x + \deg \frak x
$,
\index{degree ! of Floer cochain complex over universal Novikov ring}
where $x \in CF(\mathcal C)^0$ and $\frak x \in \Lambda_{\rm nov}$
are elements of pure degree.
The ${\Lambda_{\rm nov}}$ module $CF(\mathcal C)^0 \otimes_{\R}{\Lambda_{\rm nov}}$ is filtered by
\index{filtration ! of Floer cochain complex over universal Novikov ring}
$$
\frak F^{\lambda}(CF(\mathcal C)^0 \otimes_{\R}{\Lambda_{\rm nov}})
=
\bigcup_{\lambda_1+ \lambda_2 \ge \lambda}
\frak F^{\lambda_1}CF(\mathcal C)^0
\otimes_{\R} \frak F^{\lambda_2}{\Lambda_{\rm nov}}.
$$
\item
For $\beta \in \frak G$ we define a ${\Lambda_{\rm nov}}$ module
homomorphism
$\widehat{\rm PI}_{\beta}^* :
CF(\mathcal C)^0 \otimes_{\R}{\Lambda_{\rm nov}} \to
CF(\mathcal C)^0 \otimes_{\R}{\Lambda_{\rm nov}}$
by
$$
\widehat{\rm PI}_{\beta}^*(x\otimes \frak x)
=
{\rm PI}_{\beta}^*(x) \otimes T^{-E(\beta)}e^{-\mu(\beta)/2}\frak x.
$$
It preserves the degree and the filtration.
\item
We consider the closed ${\Lambda_{\rm nov}}$ submodule of
$CF(\mathcal C)^0 \otimes_{\R}{\Lambda_{\rm nov}}$ generated
by elements of the form
$\widehat{\rm PI}_{\beta}^*(x\otimes \frak x) - x\otimes \frak x$
for $x\otimes \frak x \in CF(\mathcal C)^0 \otimes_{\R}{\Lambda_{\rm nov}}$
and $\beta \in \frak G$.
We denote by
$
(CF(\mathcal C)^0 \otimes_{\R}{\Lambda_{\rm nov}})/\sim
$
the quotient module by this submodule.
It is filtered and graded.
\item
We denote the completion of
$
(CF(\mathcal C)^0 \otimes_{\R}{\Lambda_{\rm nov}})/\sim
$ with respect to the filtration by
$
CF(\mathcal C;\Lambda_{\rm nov}).
$
It is a filtered and graded ${\Lambda_{\rm nov}}$ module.
We write its filtration by $\frak F^{\lambda}CF(\mathcal C;\Lambda_{\rm nov})$.
\item
We put
$CF(\mathcal C;\Lambda_{0,\rm nov}) = \frak F^{0}CF(\mathcal C;\Lambda
_{\rm nov})$.
It is a $\Lambda_{0,\rm nov}$ module.
\item
A subset $G \subset \R_{\ge 0} \times \Z$ is called a {\it discrete submonoid} if
the following holds.
Denote by $E : G  \to \R_{\ge 0}$ and $\mu : G  \to \Z$ the natural projections.
\begin{enumerate}
\item If $\beta_1,\beta_2 \in  G$, then $\beta_1 + \beta_2 \in G$.
$(0,0) \in G$.
\item The image $E(G) \subset \R_{\ge 0}$ is discrete.
\item For each $E_0 \in \R_{\ge 0}$ the inverse image
$G \cap E^{-1}([0,E_0])$ is a finite set.
\end{enumerate}
\item
A $\Lambda_{0,\rm nov}$-module homomorphism
$\psi : CF(\mathcal C_1;\Lambda_{0,\rm nov}) \to CF(\mathcal C_2;\Lambda_{0,\rm nov})$
\footnote{Here the case $\mathcal C_1=\mathcal C_2$ is also included.}
is called {\it $G$-gapped},
if $\psi$ is decomposed into the following form:
\begin{equation}\label{eq:gapped}
\psi = \sum_{\beta \in G} \psi_{\beta} T^{E(\beta)}e^{\mu(\beta)/2},
\end{equation}
where $\psi_{\beta} : CF(\mathcal C_1) \to \mathcal CF(\mathcal C_2)$ is an $\R$-linear map.
In addition, if $\psi$ is a cochain map, it is called a
{\it $G$-gapped cochain map}.
We simply say that $\psi$ is {\it gapped}
if it is $G$-gapped for some discrete submonoid $G$.
\index{gapped ! cochain map}
\item
If an operator $d : CF(\mathcal C;\Lambda_{0,\rm nov}) \to
CF(\mathcal C;\Lambda_{0,\rm nov})$ satisfying $d\circ d =0$ is gapped,
we call $(CF(\mathcal C;\Lambda_{0,\rm nov}),d)$
a {\it gapped cochain complex}.
\index{gapped ! cochain complex}
\end{enumerate}
\end{defn}
\begin{lem}
$CF(\mathcal C;\Lambda_{0,\rm nov})$ is a free
$\Lambda_{0,\rm nov}$ module.
\end{lem}
\begin{proof}
This is a consequence of Condition \ref{linsysmainconds} (II)
and its definition.
(See also \cite[Lemma 2.3]{fooospectr}.)
\end{proof}
\begin{cor}\label{corollary1216}
Let $\mathcal C$ be the critical submanifold data given in Definition
\ref{linearsystemdefn}.
\begin{enumerate}
\item
In the situation of
Theorem \ref{linesysmainth1} (1),
the operator $d$ on
$CF(\mathcal C)$ induces an operator
$: CF(\mathcal C;{\Lambda_{0,\rm nov}}) \to CF(\mathcal C;{\Lambda_{0,\rm nov}})$, which we also denote by $d$.
It has the same properties as in Theorem \ref{linesysmainth1} (1) (a)-(d).
Moreover, it is gapped
in the sense of Definition \ref{defn141212}.
\item
In the situation of Theorem \ref{linesysmainth1}(2), the cochain map $\psi$ induces a cochain map
$\psi : (CF(\mathcal C,{\Lambda_{0,\rm nov}}),d_1) \to
(CF(\mathcal C,{\Lambda_{0,\rm nov}}),d_2)$.
It has the same properties as in Theorem \ref{linesysmainth1} (2) (a)-(d).
Moreover, it is a gapped cochain map.
\end{enumerate}
\end{cor}
\begin{proof}
(1) This is immediate from
Theorem \ref{linesysmainth1} (1)  (b)(c)(d).
The gapped-ness follows from the compactness axiom (IX) in
Condition \ref{linsysmainconds}.
\par
(2) This is immediate from
Theorem \ref{linesysmainth1} (2)
(b)(c)(d).
From the the construction of $\psi$ given in
Subsection \ref{subsec:proofsec14main2},
we can see the gapped-ness. Here the compactness axiom
Condition \ref{morphilinsys} (IX) for morphisms of linear K-systems is used.
\end{proof}
\begin{defn}
In the situation of Theorem \ref{linesysmainth1} (1),
we call the cohomology of $(CF(\mathcal C;\Lambda_{\rm nov}),d)$
the {\it Floer cohomology} of linear K-system $\mathcal F$.
\index{Floer cohomology ! of linear K-system}
\index{K-system ! Floer cohomology of linear K-system}
It is independent of the choices by Corollary \ref{corollary1216}.
We call the cochain complex $(CF(\mathcal C;\Lambda_{\rm nov}),d)$
the {\it Floer cochain complex over the universal Novikov ring}.
\index{Floer cochain complex ! over universal Novikov ring}
\end{defn}

\subsection{Morphism of linear K-systems}
\label{subsec:morphlinsys}

We next define morphisms between
linear K-systems.
\par
\begin{shitu}\label{situatinmor}
\begin{enumerate}
\item
For each $i=1,2$,
let
$$
\mathcal C_i =\Big(
\frak A_i, \frak G_i, \{ R^i_{\alpha_i}\}_{\alpha \in \frak A_i},
\{ o_{R^i_{\alpha_i}} \}_{\alpha \in \frak A_i}, E, \mu,
\{ {\rm PI}^i_{\beta_i,\alpha_i} \}_{\beta_i \in \frak G_i, \alpha_i \in \frak A_i}
\Big)
$$
be a critical submanifold data
and
$$
\mathcal F_i = \Big( {\mathcal C}_i,
\{\mathcal M^i(\alpha_{i-},\alpha_{i+}) \}_{\alpha_{i\pm} \in \frak A_i},
({\rm ev}_{-}, {\rm ev}_{+}),
\{{\rm OI}^i_{\alpha_{i-}, \alpha_{i+}}\}_{\alpha_{i\pm} \in \frak A_i},
\{ {\rm PI}^i_{\beta_i;\alpha_{i-},\alpha_{i+}} \}_{\beta_i \in \frak G_i, \alpha_{i\pm} \in \frak A_i}
\Big)
$$
a linear K-system.
We assume $\frak G_1 = \frak G_2$
(together with energy $E$ and the Maslov index $\mu$ on it) and denote it by $\frak G$.
\item
We also consider the same objects as in (1) except we only suppose that they consist of
partial linear K-systems of energy cut
level $E_0$.$\blacksquare$
\end{enumerate}
\end{shitu}
\begin{conds}\label{morphilinsys}
In Situation \ref{situatinmor} (1)
we consider the following objects.
\par\smallskip
\noindent {\bf (I)(II)(III)}  Nothing to add to those from Condition \ref{linsysmainconds}.
\par\smallskip
\noindent {\bf (IV) (Interpolation space)}
For any $\alpha_1\in \frak A_1$ and $\alpha_2\in \frak A_2$,
we have a
K-space with corners $\mathcal N(\alpha_1,\alpha_2)$
and strongly smooth maps
$$
({\rm ev}_{-},{\rm ev}_{+})  :
{\mathcal N}({\alpha_1},{\alpha_2})
\to
R^1_{\alpha_1} \times R^2_{\alpha_2}.
$$
We assume that ${\rm ev}_{+}$ is weakly submersive.
We call $\mathcal N(\alpha_1,\alpha_2)$
the {\it  interpolation space}
\index{interpolation space}
and ${\rm ev}_{\pm}$ the {\it evaluation maps at infinity}.
\index{evaluation maps at infinity ! of interpolation space}
\par\smallskip
\noindent {\bf (V) (Energy loss)}
We assume  $
\mathcal N(\alpha_1,\alpha_2) = \emptyset$
if $E(\alpha_1) \ge E(\alpha_2) +c$ for some $c\ge 0$.
We call $c$ the {\it energy loss}.
\index{energy loss ! of morphism of linear K-systems}
\par
There is an exception in case $c=0$. See Definition \ref{defn1521}.
\par\smallskip
\noindent {\bf (VI) (Dimension)}
The dimension of the interpolation space is given by
\begin{equation}
\dim
{\mathcal N}({\alpha_1},{\alpha_2})
=
\mu(\alpha_2) - \mu(\alpha_1)  + \dim R_{\alpha_2}.
\end{equation}
\par\smallskip
    \noindent{\bf (VII) (Orientation)}
For any $\alpha_1\in \frak A_1$ and $\alpha_2\in \frak A_2$, an isomorphism
\begin{equation}
{\rm OI}_{\alpha_1,\alpha_2} :
{\rm ev}_{-}^*(o_{R^1_{\alpha_{1}}})
\cong
{\rm ev}_{+}^*(o_{R^2_{\alpha_{2}}})
\otimes
{\rm ev}_{+}^*(\det TR^2_{\alpha_{2}})
\otimes
o_{{\mathcal N}({\alpha_1},{\alpha_2})}
\end{equation}
of principal $O(1)$ bundles is given. Here $o_{{\mathcal N}({\alpha_1},{\alpha_2})}$ is the
orientation bundle which gives an orientation of K-space
${\mathcal N}({\alpha_1},{\alpha_2})$.
We call ${\rm OI}_{\alpha_1,\alpha_2}$ an {\it orientation isomorphism}.
\index{orientation isomorphism ! of interpolation space}
\par\smallskip
    \noindent{\bf (VIII) (Periodicity)}
For any $\beta \in \frak G$ an isomorphism
\begin{equation}\label{form128morph}
{\rm PI}_{\beta;\alpha_1,\alpha_2}
: {\mathcal N}({\alpha_1},{\alpha_2})
\to
{\mathcal N}({\beta\alpha_1},{\beta\alpha_2})
\end{equation}
of K-spaces is given such that the equality
$$
{\rm PI}_{\beta_2;\beta_1\alpha_1,\beta_1\alpha_2}\circ {\rm PI}_{\beta_1;\alpha_1,\alpha_2}
=
{\rm PI}_{\beta_2\beta_1;\alpha_1,\alpha_2}
$$
holds and
the following diagram commutes.
\begin{equation}
\begin{CD}
{\mathcal N}({\alpha_1},{\alpha_2}) @ > {{\rm PI}_{\beta;\alpha_1,\alpha_2}} >>
{\mathcal N}({\beta\alpha_1},{\beta\alpha_2})  \\
@ V{({\rm ev}_{1},{\rm ev}_{2})}VV @ VV{({\rm ev}_{1},{\rm ev}_{2})}V\\
R^1_{\alpha_1}\times R^2_{\alpha_2} @ >  {({\rm PI}_{\beta;\alpha_1},{\rm PI}_{\beta;\alpha_2})} >>
R^1_{\beta\alpha_1}\times R^2_{\beta\alpha_2}
\end{CD}
\end{equation}
We call ${\rm PI}_{\beta;\alpha_1,\alpha_2}$
the {\it periodicity isomorphism}.
\index{periodicity isomorphism ! of interpolation space}
The periodicity isomorphism commutes with the orientation
isomorphism.
\par\smallskip
\noindent {\bf (IX) (Gromov compactness)}
\index{Gromov compactness ! of interpolation space}
For any $E_0 \ge 0$ and $\alpha_1 \in \frak A_1$ the set
\begin{equation}
\{\alpha_2 \in \frak A_2 \mid
{\mathcal N}({\alpha_1},{\alpha_2})
\ne \emptyset,
\,\, E(\alpha_2) \le E_0 + E(\alpha_1)\}
\end{equation}
is a finite set.
\par\smallskip
\noindent {\bf (X) (Compatibility at the boundary)}
\index{compatibility at the boundary ! of interpolation space}
The  normalized boundary of the interpolation space
is decomposed into the fiber products as follows:\footnote{See Remark \ref{rem:FiberProdOrd} for the sign and the order of the factors in the fiber products.}
\begin{equation}\label{formula1211morph}
\aligned
&\partial{\mathcal N}({\alpha_1},{\alpha_2})
\\
\cong &\coprod_{\alpha'_1 \in \frak A_1}
(-1)^{\dim {\mathcal N}({\alpha}'_1,{\alpha_2})}\left(
{\mathcal N}({\alpha}'_1,{\alpha_2})
\,\,{}_{{\rm ev}_{-}}\times_{{\rm ev}_{+}}\,\,
{\mathcal M}^1({\alpha_1},{\alpha'_1})
\right) \\
&\sqcup \coprod_{\alpha'_2 \in \frak A_2}
(-1)^{\dim {\mathcal M}^2({\alpha}'_2,{\alpha_2})}\left(
{\mathcal M}^2({\alpha}'_2,{\alpha_2})
\,\,{}_{{\rm ev}_{-}}\times_{{\rm ev}_{+}}\,\,
{\mathcal N}({\alpha_1},{\alpha'_2})
\right),
\endaligned
\end{equation}
where the right hand side is the disjoint union.\footnote{
The right hand side is a finite sum by Condition (IX).}
Here $\cong$ means the isomorphism of K-spaces.
This isomorphism
commutes with the periodicity and
orientation isomorphisms.
It is compatible with various evaluation maps.
\par\smallskip
\noindent {\bf (XI) (Compatibility at the corner I)}
\index{compatibility at the corner ! of interpolation space}
The normalized corner
$\widehat S_{k}({\mathcal N}({\alpha_1},{\alpha_2}))$
is decomposed into a disjoint union of
\begin{equation}\label{eachtermofSkN}
\aligned
&{\mathcal M^1}({\alpha_1},\alpha_{1,1})
\times_{R^1_{\alpha_{1,1}}}  \dots
\times_{R^1_{\alpha_{1,k_1-1}}}
{\mathcal M^1}({\alpha_{1,k_1-1}},\alpha_{1,k_1}) \\
&\times_{R^1_{\alpha_{1,k_1}}}
{\mathcal N}(\alpha_{1,k_1},\alpha_{2,1}) \\
&\times_{R^1_{\alpha_{2,1}}}{\mathcal M}^2({\alpha_{2,1}},\alpha_{2,2})
\times_{R^2_{\alpha_{2,2}}}
 \dots
\times_{R^2_{\alpha_{2,k_2}}}
{\mathcal M^2}({\alpha_{2,k_2}},\alpha_{2})
\endaligned
\end{equation}
where $k_1+k_2 = k$, $\alpha_{1,i} \in \frak A_1$, $\alpha_{2,i} \in \frak A_2$.
This is an isomorphism of K-spaces.
The isomorphism preserves the
periodicity isomorphisms
and is compatible with various evaluation maps.
\par\smallskip
\noindent {\bf (XII) (Compatibility at the corner II)}
The isomorphism (XI) satisfies the compatibility condition given by Condition \ref{condsmorcompcomp} below.
\end{conds}
\begin{conds}\label{condsmorcompcomp}
Condition \ref{morphilinsys} (XI) and
Condition \ref{linsysmainconds} (XI) (together with (\ref{form3077}))
define an isomorphism from
$\widehat S_{\ell}(\widehat S_{k}({\mathcal N}({\alpha_1},{\alpha_2})))$
to the disjoint union of the summand (\ref{eachtermofSkN})
with $k$ replaced by $k+\ell$.
Then the map $\pi_{\ell,k} : \widehat S_{\ell}(\widehat S_{k}({\mathcal N}({\alpha_1},{\alpha_2})))
\to \widehat S_{\ell+k}({\mathcal N}({\alpha_1},{\alpha_2}))$
in Proposition \ref{prop2813}
becomes the identity map under those isomorphisms.
\end{conds}
\begin{defn}\label{linearsystemmorphdefn}
Let us consider Situation \ref{situatinmor} .
\begin{enumerate}
\item
A {\it morphism of linear K-systems}
\index{morphism ! of linear K-systems}
\index{K-system ! morphism of linear K-systems}
from $\mathcal F_1$ to $\mathcal F_2$
consists of the objects as in Condition \ref{morphilinsys}.
We say $c$ the {\it energy loss} of our morphism.
\index{energy loss ! energy loss}
\item
Let $E_i > 0$ and let
$\mathcal F_i$ be a partial linear K-system
with energy cut levels $E_i$ for each $i=1,2$.
Suppose $E_1 \ge E_2 + c$ with $c\ge 0$.
A {\it morphism of partial linear K-systems}
\index{morphism ! of partial linear K-systems}
\index{K-system ! morphism of partial linear K-systems}
from $\mathcal F_1$ to $\mathcal F_2$
consists of the same objects as those in Condition \ref{morphilinsys}
except the following:
\begin{enumerate}
\item
The K-space, the interpolation space,
${\mathcal N}({\alpha_1},{\alpha_2})$ is
defined only when $-c \le E({\alpha_2}) - E({\alpha_1}) \le E_2$.
\item
The periodicity isomorphisms of the interpolation spaces
are defined only among those satisfying $-c \le E({\alpha_2}) - E({\alpha_1}) \le E_2$.
\item
The compatibility of the isomorphism of K-spaces at the boundary in Condition \ref{morphilinsys} (X)
is defined only when the left hand side of (\ref{formula1211morph})
is defined.
\item
The compatibility of the isomorphism of K-spaces at the corners in Condition \ref{morphilinsys} (XI)
is required only for ${\mathcal N}({\alpha_1},{\alpha_2})$
with $-c \le E({\alpha_2}) - E({\alpha_1}) \le E_2$.
\end{enumerate}
\end{enumerate}
\end{defn}
\begin{defn}\label{defn1521}
Suppose that the index set $\frak A_1$ of the
critical submanifold of $\mathcal F_1$
is identified with the index set $\frak A_2$ of the
critical submanifold of $\mathcal F_2$.
We define the notion of a {\it morphism with energy loss $0$
and congruent to an isomorphism} by requiring the following:
\index{morphism with energy loss $0$
and congruent to the isomorphism}
\par
For $\alpha = \alpha_1 = \alpha_2 \in \frak A_1
= \frak A_2$, we require
$
\mathcal N(\alpha,\alpha) = R_{\alpha}$
instead of requiring it to be an empty set.
We also require the maps  ${\rm ev}_{\pm} : R_{\alpha} \to R_{\alpha}$ to be identity maps.
Moreover for $\alpha_1 \ne \alpha_2$ we require
$\mathcal N(\alpha_1,\alpha_2) = \emptyset$ if $E(\alpha_1)=E(\alpha_2)$.
\par
If $\mathcal F_1 = \mathcal F_2$ in addition, we call the morphism a
{\it morphism with energy loss $0$
and congruent to the identity morphism},
instead of a
morphism with energy loss $0$
and congruent to an isomorphism.
\end{defn}

\subsection{Homotopy and higher homotopy
of morphisms of linear K-systems}
\label{subsec:higherhomotopyksys}
To establish basic properties of Floer cohomology associated to a
linear K-system or an inductive system
of partial linear K-systems,
we will use the notion of homotopy between the morphisms
and also higher homotopy such as homotopy of homotopies etc..
To include the most general case, we define the notion of morphisms
parametrized by a manifold with corners $P$.
\begin{shitu}\label{situparaPmorph}
\begin{enumerate}
\item
Let $\mathcal C_i, \mathcal F_i$ ($i=1,2$)
be critical submanifold data and linear K-systems
as in Situation \ref{situatinmor}, respectively.
We assume
$\frak G_1=\frak G_2$ and denote the common group by $\frak G$.
\item
Let $P$ be a smooth manifold with boundary or corner.
$\blacksquare$
\end{enumerate}
\end{shitu}
\begin{conds}\label{Pparamorphi}
In Situation \ref{situparaPmorph} we consider the following
objects.
\par\smallskip
\noindent
{\bf (I), (II), (III)}  Nothing to add.
\par\smallskip
\noindent
{\bf (IV)
(Interpolation space)}
For any $\alpha_1 \in \frak A_1$ and $\alpha_2 \in \frak A_2$, we have a
K-space with corners  $\mathcal N(\alpha_1,\alpha_2;P)$
and a strongly smooth map
$$
({\rm ev}_{P},{\rm ev}_{-},{\rm ev}_{+})  :
\mathcal N(\alpha_1,\alpha_2;P)
\to
P\times R^1_{\alpha_1} \times R^2_{\alpha_2}.
$$
We assume that
$$
({\rm ev}_{P},{\rm ev}_{+}) :
\mathcal N(\alpha_1,\alpha_2;P)
\to
P \times R^2_{\alpha_2}
$$
is a corner stratified weak submersion.
(See Definition \ref{defn3055} for the notion of corner stratified
submersivity.)\index{submersion ! corner stratified weak submersion}
\index{corner ! corner stratified weak submersion}
We call $\mathcal N(\alpha_1,\alpha_2;P)$  a {\it $P$-parametrized interpolation space}.
\index{$P$-parametrized interpolation space}
\par\smallskip
\noindent {\bf (V) (Energy loss)}
We assume  $
\mathcal N(\alpha_1,\alpha_2;P) = \emptyset$
if $E(\alpha_1) \ge E(\alpha_2) +c$.
We call $c\ge 0$ the {\it energy loss}.
\par
Exception: In case $c=0$ and $\alpha = \alpha_1 = \alpha_2$,
we require
 $
\mathcal N(\alpha,\alpha;P) = P \times R_{\alpha}$
instead of requiring it to be an empty set.
We also require that ${\rm ev}_{\pm} : P \times R_{\alpha} \to R_{\alpha}$
and ${\rm ev}_{P} : P\times R_{\alpha} \to P$  are projections.
\par\smallskip
\noindent {\bf (VI) (Dimension)}
The dimension is given by
\begin{equation}
{\mathcal N}({\alpha_1},{\alpha_2};P)
=
\mu(\alpha_2) - \mu(\alpha_1)  + \dim R^2_{\alpha_2} + \dim P.
\end{equation}
\par\smallskip
\noindent {\bf (VII) (Orientation)}
$P$ is oriented.
For any $\alpha_1 \in \frak A_1$ and $\alpha_2 \in \frak A_2$,
we are given an {\it orientation isomorphism}
\begin{equation}
{\rm OI}_{\alpha_1,\alpha_2;P} :
{\rm ev}_{-}^*(o_{R^1_{\alpha_{1}}})
\cong
{\rm ev}_{+}^*(o_{R^2_{\alpha_{2}}})
\otimes
{\rm ev}_{+}^*(\det TR^2_{\alpha_2})
\otimes
o_{{\mathcal N}({\alpha_1},{\alpha_2};P)}.
\end{equation}
\par\smallskip
\noindent {\bf (VIII) (Periodicity)}
For any $\beta \in \frak G$  an isomorphism
\begin{equation}\label{form1282}
{\rm PI}_{\beta;\alpha_1,\alpha_2;P}
:
{\mathcal N}({\alpha_1},{\alpha_2};P)
\to
{\mathcal N}({\beta\alpha_1},{\beta\alpha_2};P)
\end{equation}
of K-spaces is given.
The equality
$$
{\rm PI}_{\beta_2;\beta_1\alpha_1,\beta_1\alpha_2;P}\circ {\rm PI}_{\beta_1;\alpha_1,\alpha_2;P}
=
{\rm PI}_{\beta_2\beta_1;\alpha_1,\alpha_2;P}
$$
holds.
The isomorphism ${\rm PI}_{\beta;\alpha_1,\alpha_2;P}$
is compatible with
$({\rm ev}_{-},{\rm ev}_{+},{\rm ev}_{P})$
and preserves the orientation isomorphism.
\par\smallskip
\noindent{\bf (IX) (Gromov compactness)}
For any $E_0\ge 0$ and $\alpha_1 \in \frak A_1$ the set
\begin{equation}
\{\alpha_2 \in \frak A_2 \mid
{\mathcal N}({\alpha_1},{\alpha_2};P)
\ne \emptyset,
\,\, E(\alpha_2) \le E_0 + E(\alpha_1)\}
\end{equation}
is a finite set.
\end{conds}
To state the next condition we note that Lemma \ref{lem304}
implies the following.
\begin{lem}
We assume ${\mathcal N}({\alpha_1},{\alpha_2};P)$
satisfies Condition \ref{Pparamorphi} (IV)-(IX).
\begin{enumerate}
\item
We can define a K-space ${\mathcal N}({\alpha_1},{\alpha_2};\partial P)$
by the fiber product:
$$
{\mathcal N}({\alpha_1},{\alpha_2};\partial P) :=
\partial P {}_P\times_{{\rm ev}_P}{\mathcal N}({\alpha_1},{\alpha_2};P).
$$
\item
The periodicity and orientation isomorphisms of ${\mathcal N}({\alpha_1},{\alpha_2};P)$
induce those of ${\mathcal N}({\alpha_1},{\alpha_2};\partial P)$.
\item
$({\rm ev}_-,{\rm ev}_+,{\rm ev}_P)$ of
${\mathcal N}({\alpha_1},{\alpha_2};P)$ induces
$({\rm ev}_-,{\rm ev}_+,{\rm ev}_{\partial P})$ of
${\mathcal N}({\alpha_1},{\alpha_2};\partial P)$.
\item
Objects defined in (1)-(3) above satisfy
Condition \ref{Pparamorphi} (I)-(IX).
\end{enumerate}
\end{lem}

\begin{conds}\label{famiboudaru22}
\noindent{\bf (X) (Compatibility at the boundary)}
The normalized boundary of the $P$-parametrized  interpolation space
is decomposed into the fiber products as follows.
\footnote{See Remark \ref{rem:FiberProdOrd}.}
\begin{equation}\label{formula1211iso}
\aligned
&\partial
{\mathcal N}({\alpha_1},{\alpha_2};P)
\\
\cong
&\coprod_{\alpha \in \frak A_2}
(-1)^{\dim {\mathcal M}^2({\alpha},{\alpha_2})}
{\mathcal M}^2({\alpha},{\alpha_2})
\,\,{}_{{\rm ev}_{-}}\times_{{\rm ev}_{+}}
{\mathcal N}({\alpha_1},{\alpha};P)
\\
&\sqcup \coprod_{\alpha \in \frak A_1}
(-1)^{\dim {\mathcal N}({\alpha},{\alpha_2};P)}
{\mathcal N}({\alpha},{\alpha_2};P)
\,\,{}_{{\rm ev}_{-}}\times_{{\rm ev}_{+}}
{\mathcal M}^1({\alpha_1},{\alpha})
\\
&\sqcup {\mathcal N}({\alpha_1},{\alpha_2};\partial P),
\endaligned
\end{equation}
where the right hand side is the disjoint union.
Here $\cong$ means the isomorphism of K-spaces.
This isomorphism preserves the
orientation isomorphism
and the periodicity isomorphism,
which in the right hand side are obtained by taking fiber products thereof respectively.
It is also compatible with various evaluation maps.
\end{conds}
To state the next condition we note that  Lemma \ref{lem304}  also implies
the following:
\begin{lem}\label{Skparafamilylem}
We assume ${\mathcal N}({\alpha_1},{\alpha_2};P)$
satisfies Condition \ref{Pparamorphi} (IV)-(IX).
\begin{enumerate}
\item
We can define a K-space by the fiber product
$$
\widehat S_k(P)
{}_P\times_{{\rm ev}_P} {\mathcal N}({\alpha_1},{\alpha_2};P).
$$
We denote it by ${\mathcal N}({\alpha_1},{\alpha_2};\widehat S_k(P))$.
\item
The periodicity
isomorphism of ${\mathcal N}({\alpha_1},{\alpha_2};P)$
induces one of ${\mathcal N}({\alpha_1},{\alpha_2};\widehat S_k(P))$.
\item
$({\rm ev}_-,{\rm ev}_+,{\rm ev}_P)$ of
${\mathcal N}({\alpha_1},{\alpha_2};P)$ induces
$({\rm ev}_-,{\rm ev}_+,{\rm ev}_{\widehat S_k(P)})$ of
${\mathcal N}({\alpha_1},{\alpha_2};\widehat S_k(P))$.
\item
Objects defined in (1)-(3) above satisfy
Condition \ref{Pparamorphi} (I)-(IX).
\end{enumerate}
\end{lem}

The next lemma is also a conclusion of Lemma \ref{lem304}.
\begin{lem}\label{deflem1425}
Suppose we are in the situation of Lemma \ref{Skparafamilylem}.
\begin{enumerate}
\item
The $(k+\ell)/k!\ell!$ fold covering map
$\widehat S_k(\widehat S_{\ell}(P))
\to \widehat S_{k+\ell}(P)$ induces a $(k+\ell)!/k!\ell!$ fold covering
map
$: {\mathcal N}({\alpha_1},{\alpha_2};\widehat S_k(\widehat S_{\ell}(P))) \to
{\mathcal N}({\alpha_1},{\alpha_2};\widehat S_{k+\ell}(P))$.
\item
The covering map (1) commutes with the
periodicity
isomorphisms. It is also compatible with various evaluation maps.
Especially the following diagram commutes.
\begin{equation}
\begin{CD}
{\mathcal N}({\alpha_1},{\alpha_2};\widehat S_k(\widehat S_{\ell}(P))) @ > {} >>
{\mathcal N}({\alpha_1},{\alpha_2};\widehat S_{k+\ell}(P))  \\
@ V{{\rm ev}_{\widehat S_k(\widehat S_{\ell}(P))}}VV
@ VV{{\rm ev}_{\widehat S_{k+\ell}(P)}}V\\
\widehat S_k(\widehat S_{\ell}(P)) @ >>>
\widehat S_{k+\ell}(P)
\end{CD}
\end{equation}
\end{enumerate}
\end{lem}

\begin{conds}\label{boundarycompPban1}
\noindent{\bf (XI) (Compatibility at the corner)}
The normalized corner $\widehat S_k \mathcal N((\alpha_-,\alpha_+);P)$
is decomposed into a disjoint union of fiber products:
\begin{equation}\label{formula143333}
\aligned
&{\mathcal M^{1}}(\alpha_-,\alpha_1)
\times_{R_{\alpha_1}}
\dots
\times_{R_{\alpha_{{k_1}-1}}}
{\mathcal M^{1}}(\alpha_{{k_1}-1},\alpha_{{k_1}}) \\
&
\times_{R_{\alpha_{k_1}}}
{\mathcal N}({\alpha_{k_1}},{\alpha_{k_1+1}};\widehat S_{k_3}(P)) \\
&
\times_{R_{\alpha_{k_1+1}}}
{\mathcal M^{2}}(\alpha_{k_1+1},\alpha_{k_1+2})
\times_{R_{\alpha_{k_1+2}}}
\dots
\times_{R_{\alpha_{k_1+k_2}}}
{\mathcal M^{2}}(\alpha_{{k_1}+k_2},\alpha_+).
\endaligned
\end{equation}
Here $k_1+k_2+k_3 = k$.
This is an isomorphism of K-spaces.
This isomorphism respects the
periodicity isomorphisms.
It is also compatible with various evaluation maps.
Moreover Condition \ref{furthercompatifiber} below holds.
\end{conds}
To state Condition \ref{furthercompatifiber}  we make a digression.
We consider $\widehat S_{\ell}((\ref{formula143333}))$,
where (\ref{formula143333}) stands for the K-space defined by
(\ref{formula143333}).
Applying (\ref{form3077}) to ({\ref{formula143333}),
we get an isomorphism from the normalized corner
$\widehat S_{\ell}(\widehat S_k( \mathcal N(\alpha_-,\alpha_+;P)))$
to the disjoint union of fiber products of the normalized corners
of the factors of ({\ref{formula143333}).
We can identify the normalized corners of various factors of
({\ref{formula143333}) by using
(\ref{cornecom1}) and Condition \ref{boundarycompPban1}.
Thus we obtain the next lemma.
\begin{lem}
The isomorphisms in Condition \ref{boundarycompPban1} and
in (\ref{cornecom1})  canonically induce an isomorphism from
$\widehat S_{\ell}(\widehat S_k(\mathcal N(\alpha_-,\alpha_+;P)))$
to a disjoint union of the fiber products of the following form:
\begin{equation}\label{formula143333rev}
\aligned
&{\mathcal M^{1}}(\alpha_-,\alpha'_1)
\times_{R_{\alpha'_1}}
\dots
\times_{R_{\alpha'_{k'_1 -1}}}
{\mathcal M^{1}}(\alpha'_{{k'_1}-1},\alpha'_{{k'_1}}) \\
&
\times_{R_{\alpha'_{k'_1}}}
{\mathcal N}({\alpha'_{k'_1}},{\alpha'_{k'_1+1}};\widehat S_{\ell'}(\widehat S_{k'_3}(P))) \\
&
\times_{R_{\alpha'_{k'_1+1}}}
{\mathcal M^{2}}(\alpha'_{k'_1+1},\alpha'_{k'_1+2})
\times_{R_{\alpha'_{k'_1+2}}}
\dots
\times_{R_{\alpha'_{k'_1+k'_2}}}
{\mathcal M^{2}}(\alpha'_{{k'_1}+k'_2},\alpha_+).
\endaligned
\end{equation}
Here $k'_1+k'_2+k'_3 + \ell' = \ell +k$.
(Note that the same copies of the form (\ref{formula143333rev})
appear several times in  $\widehat S_{\ell}(\widehat S_k(\mathcal N(\alpha_-,\alpha_+;P)))$.)
\end{lem}
Now we state the last condition.
\begin{conds}\label{furthercompatifiber}
{\rm\bf (Compatibility at the corner continued)}
Under the isomorphism between
(\ref{formula143333rev}) and
(\ref{formula143333})
(with $k$ replaced by $k+\ell$),
the $(k+\ell)!/k!\ell!$ fold covering map
$\widehat S_{\ell}(\widehat S_k (\mathcal N(\alpha_-,\alpha_+;P)))
\to \widehat S_{\ell+k}(\mathcal N(\alpha_-,\alpha_+;P))$
obtained in Proposition \ref{prop2813}
is identified with the fiber product of the identity
maps and of
$$
{\mathcal N}({\alpha'_{k'_1}},{\alpha'_{k'_1+1}};\widehat S_{\ell'}(\widehat S_{k'_3}(P)))
\to
{\mathcal N}({\alpha'_{k'_1}},{\alpha'_{k'_1+1}};\widehat S_{\ell'+k'_3}(P)),
$$
which is the map given in Lemma \ref{deflem1425}.
\end{conds}
\begin{defn}\label{defn1429}
\begin{enumerate}
\item
A {\it $P$-parametrized family of morphisms}
\index{parametrized family of morphisms ! parametrized family of morphisms} from $\mathcal F_1$
to $\mathcal F_2$ consists of objects satisfying Conditions \ref{Pparamorphi},
\ref{famiboudaru22}, \ref{boundarycompPban1}, \ref{furthercompatifiber}.
\item
If a $P$-parametrized family of morphisms
$\frak N_P$ from $\mathcal F_1$
to $\mathcal F_2$ and a $\partial P$-parametrized family of morphisms
$\frak N_{\partial P}$ from $\mathcal F_1$
to $\mathcal F_2$ are related as in Condition \ref{famiboudaru22},
we call $\frak N_{\partial P}$ the {\it boundary of $\frak N_P$}
\index{boundary ! of parametrized family of morphisms}
and write $\partial \frak N_P$.
\item
Let $E_i > 0$ and
$\mathcal F_i$ be partial linear K-systems
with energy cut levels $E_i$ ($i=1,2$).
Suppose $E_1 \ge E_2 + c$ for some $c\ge 0$.
A {\it $P$-parametrized family of morphisms of partial linear K-systems}
\index{parametrized family of morphisms
! of partial linear K-systems}
\index{K-system ! parametrized family of morphisms
of partial linear K-systems}
from $\mathcal F_1$ to $\mathcal F_2$
consists of the same objects as in Condition \ref{Pparamorphi}
except the following:
\begin{enumerate}
\item
The K-space, the $P$ parametrized interpolation space,
${\mathcal N}({\alpha_1},{\alpha_2};P)$ is
defined only when $-c \le E({\alpha_2}) - E({\alpha_1}) \le E_2$.
\item
The periodicity and orientation isomorphisms of the $P$-parametrized interpolation spaces
are defined only among those which satisfy $-c \le E({\alpha_2}) - E({\alpha_1}) \le E_2$.
\item
The compatibility of the isomorphism of K-spaces at the boundary in Condition \ref{famiboudaru22} (X)
is defined only when the left hand side of (\ref{formula1211iso}) is defined.
\item
We require Condition \ref{boundarycompPban1}, \ref{furthercompatifiber} (XI)
only when $\mathcal N(\alpha_-,\alpha_+;P)$ is
defined.
\end{enumerate}
\end{enumerate}
\end{defn}
\begin{defn}\label{def:homotopymorph}
Let $\mathcal F_1$, $\mathcal F_2$ be linear K-systems
and $\frak N_1$, $\frak N_2$ two morphisms between them.
\begin{enumerate}
\item
 {\it A homotopy from $\frak N_1$ to $\frak N_2$}
\index{homotopy ! morphisms of linear K-systems}
\index{K-system ! homotopy between morphisms of linear K-systems}
is a $P = [1,2]$ parametrized family of morphisms
from $\frak N_1$ to $\frak N_2$
such that its boundary is $\frak N_1 \sqcup -\frak N_2$.
Here $-\frak N_2$ is $\frak N_2$ with the orientation systems
of the critical submanifolds and
orientation isomorphisms inverted
and $\sqcup$ denotes the disjoint union.
\item
We say that $\frak N_1$ is {\it homotopic}
\index{homotopic ! morphisms of linear K-systems} to $\frak N_2$ if there
exists a homotopy between them.
\item
The case of morphisms between partial linear K-systems
is defined in the same way.
\end{enumerate}
\end{defn}

\begin{thm}\label{linesysmainth2}
Let $\mathcal F_1$, $\mathcal F_2$, $\mathcal F_3$ be linear K-systems.
We make choices mentioned in Theorem \ref{linesysmainth1} (2) and
obtain cochain complexes $CF(\mathcal F_i;\Lambda_{\rm nov})$, $i=1,2,3$.
\begin{enumerate}
\item
A morphism $\frak N : \mathcal F_1 \to \mathcal F_2$ induces
a $\Lambda_{\rm nov}$ module homomorphism
\begin{equation}
\psi_{\frak N} : CF(\mathcal F_1;\Lambda_{\rm nov})
\to CF(\mathcal F_2;\Lambda_{\rm nov})
\end{equation}
with the following properties.
\begin{enumerate}
\item
$\psi_{\frak N}  \circ d = d \circ \psi_{\frak N}$,
where $d$ in the left hand side (resp. right hand side)
is the coboundary oprator of
$CF(\mathcal F_1;\Lambda_{\rm nov})$ (resp.
$CF(\mathcal F_2;\Lambda_{\rm nov})$.)
\item
The homomorphism $\psi_{\frak N}$ preserves degree and satisfies
$$
\psi_{\frak N}\left(\frak F^{\lambda}CF(\mathcal F_1;\Lambda_{\rm nov})
\right) \subset
\frak F^{\lambda-c}CF(\mathcal F_1;\Lambda_{\rm nov}),
$$
where $c$ is the energy loss of $\frak N$.
In particular, if $c=0$ then $\psi_{\frak N}$ induces a
$\Lambda_{0,\rm nov}$ module homomorphism :
$$
\psi_{\frak N} : CF(\mathcal F_1;\Lambda_{0,\rm nov})
\to CF(\mathcal F_2;\Lambda_{0,\rm nov}).
$$
\end{enumerate}
\item
Let $\frak N_1 , \frak N_2 : \mathcal F_1 \to \mathcal F_2$ be morphisms.
If $\frak H$ is a homotopy from $\frak N_1$ to $\frak N_2$,
then it induces a $\Lambda_{\rm nov}$ module homomorphism:
\begin{equation}
\psi_{\frak H} : CF(\mathcal F_1;\Lambda_{\rm nov})
\to CF(\mathcal F_2;\Lambda_{\rm nov})
\end{equation}
with the following properties.
\begin{enumerate}
\item
$\psi_{\frak H}  \circ d + d \circ \psi_{\frak H} =\psi_{\frak N_1}
- \psi_{\frak N_2}$.
\item
The homomorphism $\psi_{\frak H}$ decreases degree by $1$ and
satisfies
$$
\psi_{\frak H}\left(\frak F^{\lambda}CF(\mathcal F_1;\Lambda_{\rm nov})
\right) \subset
\frak F^{\lambda-c}CF(\mathcal F_1;\Lambda_{\rm nov})
$$
where $c$ is the energy loss of $\frak H$.
In particular, if $c=0$ then $\psi_{\frak H}$ induces a
$\Lambda_{0,\rm nov}$ module homomorphism :
$$
\psi_{\frak H} : CF(\mathcal F_1;\Lambda_{0,\rm nov})
\to CF(\mathcal F_2;\Lambda_{0,\rm nov}).
$$
\end{enumerate}
\item Let $\frak N_{i+1 i} : \mathcal F_i \to \mathcal F_{i+1}$ be morphisms of linear K-systems of for $i=1,2$.
Let  $\frak N_{31} : \mathcal F_1 \to \mathcal F_{3}$ be the composition
$\frak N_{32}\circ \frak N_{21}$.
Then the composition $\psi_{\frak N_{23}} \circ \psi_{\frak N_{12}}$
is cochain homotopic to $\psi_{\frak N_{13}}$.
\item
If ${\mathcal ID}_{\mathcal F_i} : \mathcal F_i \to  \mathcal F_i$ is the identity morphism, the cochain
map $\psi_{{\mathcal ID}_{\mathcal F_i}}$ induced by it is cochain homotopic to the identity.
\end{enumerate}
\end{thm}
\begin{rem}
The morphism $\psi_{\frak N}$ in Item (1) depends on the choices made for its construction.
More precisely, we first make choices mentioned in Theorem \ref{linesysmainth1} Item (1)
to define coboundary oprators of
$CF(\mathcal{F}_i;\Lambda_{\rm nov})$, $i=1,2$.
Then we make a choice (compatible with the first choices)
to define the map $\psi_{\frak N}$ which is
a cochain map with respect to the coboundary operators
obtained from the choices we made.
Note that we can use Item (2) to show that up to cochain homotopy the cochain map
$\psi_{\frak H}$ is independent of the choices
we made to define it.
\end{rem}

We will define the notion of the identity morphism later in
Lemma-Definition \ref{lemdefidentity} of Subsection \ref{subsec:identitylinsys}.
The definition of composition of morphisms is given in Subsection \ref{subsec:compmorline} and the precise definition of the notation used there is postponed until
Section \ref{section:compomorphis} where we will discuss `smoothing corners' systematically.
The proof of Theorem \ref{linesysmainth2} is given in Section \ref{sec:systemline3}.

\begin{defn}\label{def:congmodT}
We say two morphisms $\frak N_1$ and $\frak N_2$ are
{\it congruent modulo $T^E$} if the following holds.
\begin{enumerate}
\item
$\frak N_1(\alpha_1,\alpha_2) \cong \frak N_2(\alpha_1,\alpha_2)$
for any $\alpha_1,\alpha_2$ with $E(\alpha_2) - E(\alpha_1) \le E$.
\item The above isomorphism is one between K-spaces.
It preserves the orientation, the periodicity isomorphism,
and the isomorphisms describing the compatibilities at their boundary and corners.
\end{enumerate}
\end{defn}

\subsection{Composition of morphisms
of linear K-systems}
\label{subsec:compmorline}

In this subsection we discuss composition of morphisms.
\begin{shitu}\label{composisitu}
Suppose we are in Situation \ref{situparaPmorph} for $i=1,2,3$
and let $\mathcal F_i$ be linear K-systems
for $i=1,2,3$.$\blacksquare$
\end{shitu}
\begin{lemdef}\label{lemdef1434}
Suppose we are in Situation \ref{composisitu}.
\begin{enumerate}
\item
Let $\frak N_{i+1 i} : \mathcal F_i \to \mathcal F_{i+1}$ be morphisms
and
${\mathcal N}_{i i+1}({\alpha_i},{\alpha_{i+1}})$
their interpolation spaces.
Then we can define the {\rm composition}
\index{composition ! of morphisms of
linear K-systems}
$\frak N_{32} \circ \frak N_{21} : \mathcal F_1 \to \mathcal F_3$.
The interpolation space\footnote{We denote by ${\mathcal N}_{12}$ the interpolation space of the morphism $\frak N_{21}$. See Remark \ref{rem:orderindex}.} of
$\frak N_{32} \circ \frak N_{21}$ is the union of the fiber products
\begin{equation}\label{compintspacess}
\bigcup_{\alpha_2 \in \frak A_2}
{\mathcal N}_{12}({\alpha_1},{\alpha_{2}})
\,{}_{{\rm ev}_+} \times^{\boxplus 1}_{{\rm ev}_-}
{\mathcal N}_{23}({\alpha_2},{\alpha_{3}}).
\end{equation}
The precise meaning of the union $\cup$ in (\ref{compintspacess}) is defined
during the proof in this subsection and in Section \ref{section:compomorphis}.
The symbol $\times^{\boxplus 1}$ will be defined in Definition \ref{defn1635}.
\item
Let $\frak N_{i+1 i}^{P_i} : \mathcal F_i \to \mathcal F_{i+1}$
be $P_i$-parametrized morphisms
and ${\mathcal N}_{i i+1}({\alpha_i},{\alpha_{i+1}};P_i)$
their interpolation spaces.
We can define the {\rm composition}
$\frak N_{32}^{P_2} \circ \frak N_{21}^{P_1} : \mathcal F_1 \to \mathcal F_3$
that is a $P_1\times P_2$-parametrized morphism.
\index{composition ! of parametrized morphisms of
linear K-systems}
\index{K-system ! composition of parametrized morphisms of
linear K-systems}
Its interpolation space is
\begin{equation}\label{eq143636}
\bigcup_{\alpha_2 \in \frak A_2}
{\mathcal N}_{12}({\alpha_1},{\alpha_{2}};P_1)
\,{}_{{\rm ev}_+} \times^{\boxplus 1}_{{\rm ev}_-}
{\mathcal N}_{23}({\alpha_2},{\alpha_{3}};P_2).
\end{equation}
The precise meaning of the union $\cup$ in (\ref{eq143636}) is defined
during the proof in this subsection and in Section \ref{section:compomorphis}.
The symbol $\times^{\boxplus 1}$ will be also defined in Definition \ref{defn1635}.
\item
The energy loss of the morphism $\frak N_{32} \circ \frak N_{21}$
is the sum of the energy loss of $\frak N_{32}$ and of $\frak N_{21}$.
The same holds for the $P_i$ parametrized version.
\item
We have $\frak N_{43} \circ (\frak N_{32} \circ \frak N_{21})
= (\frak N_{43} \circ \frak N_{32}) \circ \frak N_{21}$.
The same holds for the $P_i$ parametrized version.
\item We have
\begin{equation}\label{form14378}
\partial(\frak N_{32}^{P_2} \circ \frak N_{21}^{P_1})
=
(\frak N_{32}^{\partial P_2} \circ \frak N_{21}^{P_1})
\cup
(\frak N_{32}^{P_2} \circ \frak N_{21}^{\partial P_1})
\end{equation}
The boundary in the left hand side is taken in the sense of
Definition \ref{defn1429} (2).
The precise meaning of the union $\cup$ in the right hand side is defined
during the proof  in this subsection and by Definition-Lemma \ref{glueparamorphi}
\item
We can generalize (1)-(5) to the case of
partial linear K-systems.
\end{enumerate}
\end{lemdef}
\begin{proof}[Idea of the proof]
Here we explain only the basic geometric idea of
the proof.
To give a rigorous proof, there is an issue
about providing a precise definition of `smoothing corner' in the
definition of the union in (\ref{compintspacess}) and (\ref{eq143636}).
We postpone this point until Section \ref{section:compomorphis}.
\par
(1)
We first explain the meaning of the union in (\ref{compintspacess}).
Note that\footnote{See Remark \ref{rem:FiberProdOrd} for the order of factors and signs.}
\begin{equation}\label{partcompintspacess}
\aligned
&\bigcup_{\alpha_2 \in \frak A_2}
\partial\left({\mathcal N}_{23}({\alpha_2},{\alpha_{3}})
\,{}_{{\rm ev}_-} \times_{{\rm ev}_+}
{\mathcal N}_{12}({\alpha_1},{\alpha_{2}})
\right) \\
 \supseteq
& \bigcup_{\alpha_2 \in \frak A_2}
\partial {\mathcal N}_{23}({\alpha_2},{\alpha_{3}})
\,{}_{{\rm ev}_-} \times_{{\rm ev}_+}
{\mathcal N}_{12}({\alpha_1},{\alpha_{2}}) \\
&\quad
\cup
\bigcup_{\alpha_2 \in \frak A_2}
(-1)^{\dim {\mathcal N}_{23}({\alpha_2},{\alpha_{3}}) + \dim R_{\alpha_2}}
{\mathcal N}_{23}({\alpha_2},{\alpha_{3}})
\,{}_{{\rm ev}_-} \times_{{\rm ev}_+}
\partial{\mathcal N}_{12}({\alpha_1},{\alpha_{2}}).
\endaligned
\end{equation}
See \cite[Lemma 8.2.3 (1)]{fooobook2} for the sign.
We also note that the fiber product
\begin{equation}\label{3fiberproductrem}
{\mathcal N}_{12}({\alpha_1},{\alpha_{2}})
\,{}_{{\rm ev}_+} \times_{{\rm ev}_-}
{\mathcal M}^2({\alpha_2},{\alpha'_{2}})
\,{}_{{\rm ev}_+} \times_{{\rm ev}_-}
{\mathcal N}_{23}({\alpha'_2},{\alpha_{3}})
\end{equation}
appears in both of the second and third lines of
(\ref{partcompintspacess}).
To construct the union in (\ref{compintspacess}),
we glue several components along the codimension one boundaries
such as (\ref{3fiberproductrem}).
We note that
those parts we glue contain
certain (higher codimensional) corners.\footnote{So the order of factors in \eqref{3fiberproductrem} does not matter
as we note in Remark \ref{rem:FiberProdOrd}.}
See Figure \ref{Figure15-1}.
\begin{figure}[h]
\centering
\includegraphics[scale=0.5]{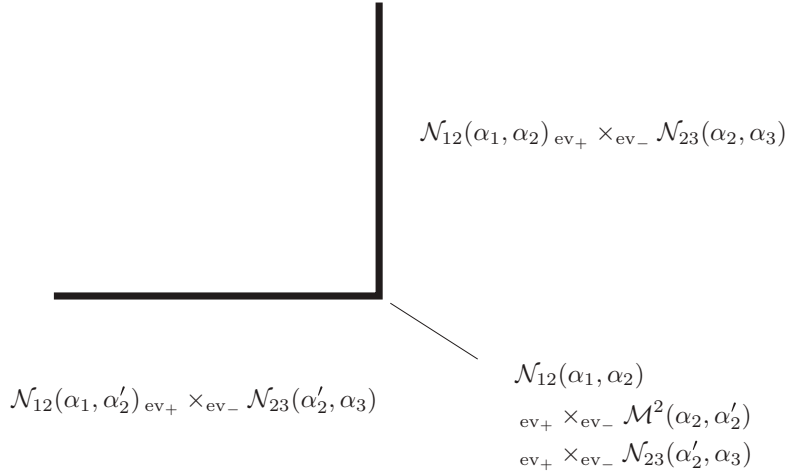}
\caption{(\ref{3fiberproductrem}) looks like a corner}
\label{Figure15-1}
\end{figure}
For this purpose we use the notion of
smoothing corners of Kuranishi structure.
We will discuss the smoothing in detail in Subsections
\ref{subsec:subsec32-1}-\ref{subsec:smoothcornerkstr}.
We also observe that there exists a certain K-space
such that the disjoint union of the summands of
(\ref{compintspacess}) appears in its boundary.
(See Proposition \ref{prop1636}.)
Combining them, we can put the
structure of K-space on the union
(\ref{compintspacess}).
See Lemma-Definition \ref{1638defken}.
Thus we obtain the interpolation space of the
composition
$$
\frak N_{32} \circ \frak N_{21}.
$$
Then it is straightforward to check that this
interpolation space satisfies the defining conditions of morphism.
\par
The proof of (2) is entirely smilar to the proof of (1).
(See Subsection \ref{subsec:paracompositions}.)
(3) is immediate from the definition.
\par
We next discuss the proof of (4). Note that both sides of the equality are a union
$$
\bigcup_{\alpha_2,\alpha_3}
{\mathcal N}_{12}({\alpha_1},{\alpha_{2}})
\,{}_{{\rm ev}_+} \times_{{\rm ev}_-}
{\mathcal N}_{23}({\alpha_2},{\alpha_{3}})
\,{}_{{\rm ev}_+} \times_{{\rm ev}_-}
{\mathcal N}_{34}({\alpha_3},{\alpha_{4}}).
$$
While we define this union from the disjoint union,
we perform the process of smoothing corner and gluing K-spaces
along the boundary twice which we described during the proof of (1).
Once at the boundary components of the form
$$
{\mathcal N}_{12}({\alpha_1},{\alpha'_{2}})
\,{}_{{\rm ev}_+} \times_{{\rm ev}_-}
{\mathcal M}^{2}({\alpha'_2},{\alpha_{2}})
\,{}_{{\rm ev}_+} \times_{{\rm ev}_-}
{\mathcal N}_{23}({\alpha_2},{\alpha_{3}})
\,{}_{{\rm ev}_+} \times_{{\rm ev}_-}
{\mathcal N}_{34}({\alpha_3},{\alpha_{4}})
$$
and once at the boundary components of the form
$$
{\mathcal N}_{12}({\alpha_1},{\alpha_{2}})
\,{}_{{\rm ev}_+} \times_{{\rm ev}_-}
{\mathcal N}_{23}({\alpha_2},{\alpha'_{3}})
\,{}_{{\rm ev}_+} \times_{{\rm ev}_-}
{\mathcal M}^{3}({\alpha'_3},{\alpha_{3}})
\,{}_{{\rm ev}_+} \times_{{\rm ev}_-}
{\mathcal N}_{34}({\alpha_3},{\alpha_{4}}).
$$
Since these two components do not intersect at the interior,
we can perform these two process independently
and can exchange the order of them.
(4) follows.
(There is an issue of showing that this isomorphism
is compatible with the smooth structure we gave
during the proof of (1). See Subsection \ref{subsec:compassoci}.)
\par
We next prove (5).
By (\ref{formula1211iso})
applied to $\mathcal N_{12}(\alpha_1,\alpha_2;{P_1})$
and to $\mathcal N_{23}(\alpha_2,\alpha_3;{P_2})$,
the normalized boundary of the
interpolation space of
$\frak N_{32}^{P_2} \circ \frak N_{21}^{P_1}$
is a disjoint union of the
following four kinds of components:
\begin{enumerate}
\item[(A)]
$\mathcal N_{12}(\alpha_1,\alpha_2;{\partial P_1})\times_{R_{\alpha_2}}
\mathcal N_{23}(\alpha_2,\alpha_3;{P_2})
$.
\item[(B)]
$\mathcal N_{12}(\alpha_1,\alpha_2;{P_1})
\times_{R_{\alpha_2}}
\mathcal N_{23}(\alpha_2,\alpha_3;{\partial P_2})
$.
\item[(C)]
$
\mathcal M^1(\alpha_1,\alpha'_1)
\times_{R_{\alpha'_1}}
\mathcal N_{12}(\alpha'_1,\alpha_2;{P_1})
\times_{R_{\alpha_2}}
\mathcal N_{23}(\alpha_2,\alpha_3;{P_2})$.
\item[(D)]
$
\mathcal N_{12}(\alpha_1,\alpha_2;{P_1})
\times_{R_{\alpha_2}}
\mathcal N_{23}(\alpha_2,\alpha'_3;{P_2})
\times_{R_{\alpha'_3}}
\mathcal M^3(\alpha'_3,\alpha'_3).
$
\end{enumerate}
Note that potentially there are components of the
form
\begin{equation}\label{formula1545}
\mathcal N_{12}(\alpha_1,\alpha'_2;{P_1})
\times_{R_{\alpha'_2}}
\mathcal M^{2}(\alpha'_2,\alpha_2)
\times_{R_{\alpha_2}}
\mathcal N_{23}(\alpha_2,\alpha_3;{P_2})
\end{equation}
in the boundary of  the
interpolation space of
$\frak N_{32}^{P_2} \circ \frak N_{21}^{P_1}$.
However (\ref{formula1545}) appears twice
in the boundary of the interpolation space of
$\frak N_{32}^{P_2} \circ \frak N_{21}^{P_1}$:
Once in
$$
\partial\mathcal N_{12}(\alpha_1,\alpha_2;{P_1})
\times_{R_{\alpha_2}}
\mathcal N_{23}(\alpha_2,\alpha_3;{P_2})
$$
and once in
$$
\mathcal N_{12}(\alpha_1,\alpha_2;{P_1})
\times_{R_{\alpha_2}}
\partial\mathcal N_{23}(\alpha_2,\alpha_3;{P_2}).
$$
When we define the interpolation space of
$\frak N_{32}^{P_2} \circ \frak N_{21}^{P_1}$,
we glue
$$
\mathcal N_{12}(\alpha_1,\alpha_2;{P_1})
\times_{R_{\alpha_2}}
\mathcal N_{23}(\alpha_2,\alpha_3;{P_2})
$$
for various $\alpha_2$ along (\ref{formula1545}).
(See Definition-Lemma \ref{glueparamorphi} for detail.)
Therefore those components of the form (\ref{formula1545})
do not appear in the boundary of
the interpolation space of
$\frak N_{32}^{P_2} \circ \frak N_{21}^{P_1}$.
Hence the disjoint union of (A)-(D)
is the normalized boundary of the
interpolation space of
$\frak N_{32}^{P_2} \circ \frak N_{21}^{P_1}$.
\par
On the other hand, the right hand side of (\ref{form14378}) is the
union of the components of types (A),(B).
Note that by Definition \ref{defn1429} (2)
and  (\ref{formula1211iso}),
we have
$$
\aligned
&(\text{The interpolation space of}\,\, \partial(\frak N_{32}^{P_2} \circ \frak N_{21}^{P_1}))
\cup {\rm (C)} \cup {\rm (D)}  \\
&= \partial\left(
\text{The interpolation space of}\,\, (\frak N_{32}^{P_2} \circ \frak N_{21}^{P_1})
\right).
\endaligned
$$
Then (5) follows from these facts.
\par
The proof of (6) is the same as the proofs of (1)-(5).
\end{proof}

\subsection{Inductive system of linear K-systems}
\begin{defn}\label{defn1528}
Let $\mathcal C$ be a critical submanifold data as in
Definition \ref{linearsystemdefn}.
\begin{enumerate}
\item
Let $E_0 < E_1$.
A partial linear K-system with
energy cut level $E_1$ induces
a partial linear K-system with
energy cut level $E_0$
by forgetting all the structures of the former
which do not appear in the definition of the latter.
The same applies to the morphism, homotopy etc.
We call these processes the {\it energy cut at $E_0$}.
\index{energy cut level!energy cut}
\item
An {\it inductive system of partial linear K-systems}
\index{K-system ! inductive system of partial linear K-systems}
$$
\mathcal F \mathcal F = (\{E^i \}, \{ \mathcal F^i \}, \{ \frak N^i \} )
$$
consists of the following objects.\footnote{Hereafter we use the superscript $i$ as the suffix of the inductive system.}
\begin{enumerate}
\item
We are given an increasing sequence $E^i$ of positive numbers
such that $\lim_{i\to \infty}E^i = \infty$.
\item
For each $i$ we are given
a partial linear K-system with
energy cut level $E^i$, which we denote by $\mathcal F^i$.
\item
The critical submanifold data
$$
\mathcal C = \Big(\frak A, \frak G, \{R_{\alpha}\}_{\alpha \in \frak A},
\{o_{R_{\alpha}} \}_{\alpha \in \frak A}, E,
\mu,
\{ {\rm PI}_{\beta,\alpha}\}_{\beta \in \frak G, \alpha \in \frak A}
\Big)
$$
that is a part of data of $\mathcal F^i$, is independent of $i$ and
is given at the beginning.
\item
For each $i$ we are given a
morphism $\frak N^i : \mathcal F^i \to \mathcal F^{i+1}\vert_{E^i}$,
where $\mathcal F^{i+1}\vert_{E^i}$
is the partial linear K-system
of energy cut level $E^i$ that is induced from $\mathcal F^{i+1}$
by the energy cut.
\item
The energy loss of $\frak N^i$ is $0$. \item
$\frak N^i$ is congruent to the identity morphism modulo $\epsilon_i$
for some $\epsilon_i > 0$. See Definition \ref{def:congmodT}.
\item
We assume the following {\it uniform Gromov compactness}.
\index{uniform Gromov compactness}
For any $E_0 \ge 0$ and $\alpha_1\in \frak A_1$ the set
\begin{equation}
\{\alpha_2 \in \frak A_2 \mid \exists i\,~
{\mathcal N}^i({\alpha_1},{\alpha_2})
\ne \emptyset,
\,\, E(\alpha_2) \le E_0 + E(\alpha_1)\}
\end{equation}
is a finite set. Moreover
\begin{equation}
\{\alpha_2 \in \frak A_2 \mid \exists i\,~
{\mathcal M}^i({\alpha_1},{\alpha_2})
\ne \emptyset,
\,\, E(\alpha_2) \le E_0 + E(\alpha_1)\}
\end{equation}
is a finite set. Here
${\mathcal N}^i({\alpha_1},{\alpha_2})$
is the interpolation space of the morphism $\frak N^i$
and ${\mathcal M}^i({\alpha_1},{\alpha_2})$ is the space
of connecting orbits of $\mathcal F^{i}$.
\end{enumerate}
\item
Let $\mathcal {FF}_j = (\{E^i_j\},\{\mathcal F_j^i\},\{\frak N^i_j\})$ $(j = 1,2)$ be
two inductive systems of partial linear K-systems.
We assume $E^i_1\ge E^i_2 - c$ for some $c\ge 0$.
A {\it morphism}
\index{morphism ! of inductive systems of partial linear K-systems}
\index{K-system ! morphism of inductive systems of partial linear K-systems}
$$
(\{\frak N^{i}_{21}\}, \{\frak H^i\})
: \mathcal {FF}_1 \to \mathcal {FF}_2
$$
with energy loss $c$
is a pair of
$\{\frak N^{i}_{21}\}$ and $\{\frak H^i\}$ with the following properties:
\begin{enumerate}
\item
$\frak N^{i}_{21}
: \mathcal F_1^i \to \mathcal F_2^i$ is a morphism of energy loss $c$.
\item
$\frak H^i$ is a homotopy between
$\frak N^{i}_{2} \circ \frak N^{i}_{21}$  and
$\frak N^{i+1}_{21} \circ \frak N^{i}_{1}$.
Here we regard them as morphisms with energy cut level
$E_2^i$ and of energy loss $c$.
\end{enumerate}
\item
In the situation of (3), let
$(\{\frak N^{i}_{(k)21}\},\{\frak H_{(k)}^i\})
: \mathcal {FF}_1 \to \mathcal {FF}_2$ be
morphisms for $k=a,b$.
A {\it homotopy}
\index{homotopy ! morphisms of inductive systems of partial linear K-systems}
\index{K-system ! homotopy of morphisms of
inductive systems of partial linear K-systems}
from $(\{\frak N^{i}_{(a)21}\},\{\frak H_{(a)}^i\})$
to $(\{\frak N^{i}_{(b)21}\},\{\frak H_{(b)}^i\})$
is the pair $(\{\frak H_{(ab)}^i\},\{\mathcal H^i\})$
with the following properties:
\begin{enumerate}
\item
$\frak H_{(ab)}^i$ is a homotopy from $\frak N^{i}_{(a)21}$ to
$\frak N^{i}_{(b)21}$.
\item
$\mathcal H^i$ is a $[0,1]^2$ parametrized morphism from
$\mathcal F_1^i$  to $\mathcal F_2^{i+1}$.
\footnote{In other words it is a homotopy of homotopies.}
\item The normalized boundary
$\partial\mathcal H^i$ is a disjoint union of the following 4 homotopies, see Figure \ref{Figure15-2}:
\begin{enumerate}
\item $\frak H_{(a)}^i$
\item $\frak N^{i}_{2} \circ \frak H_{(ab)}^i$
\item $\frak H_{(b)}^i$
\item  $\frak H_{(ab)}^{i+1} \circ \frak N^{i}_{1}$
\end{enumerate}
\begin{figure}[h]
\centering
\includegraphics[scale=0.25]{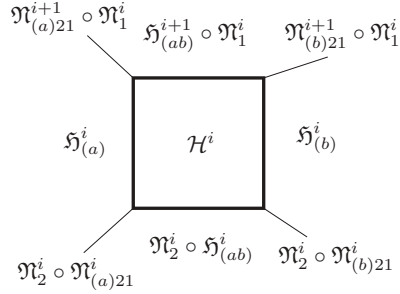}
\caption{homotopy of homotopies $\mathcal H^i$}
\label{Figure15-2}
\end{figure}
\item
For any $E_0\ge 0$ and $\alpha_1\in \frak A_1$ the set
\begin{equation}
\{\alpha_2 \in \frak A_2 \mid \exists i\,~
{\mathcal N}({\alpha_1},{\alpha_2};{\frak H_{(ab)}^i})
\ne \emptyset,
\,\, E(\alpha_2) \le E_0 + E(\alpha_1)\}
\end{equation}
is a finite set. Moreover
\begin{equation}
\{\alpha_2 \in \frak A_2 \mid \exists i\,~
{\mathcal N}({\alpha_1},{\alpha_2};{\mathcal H^i})
\ne \emptyset,
\,\, E(\alpha_2) \le E_0 + E(\alpha_1)\}
\end{equation}
is a finite set.
Here ${\mathcal N}({\alpha_1},{\alpha_2};{\frak H_{(ab)}^i})$
(resp. ${\mathcal N}({\alpha_1},{\alpha_2};{\mathcal H^i})$)
is the interpolation space which we used to define
${\frak H_{(ab)}^i}$
(resp. ${\mathcal H^i}$.).
\end{enumerate}
Here all the (parametrized) morphisms have
energy cut level $E_2^i$.
\item
Two morphisms of inductive systems of partial linear K-systems
are said to be {\it homotopic}
\index{homotopic ! morphisms of inductive systems of partial linear K-systems}
if there exists a homotopy between them.
\end{enumerate}
\end{defn}
\begin{rem}
In certain cases especially in the study of symplectic homology (see \cite{FlHofer}, \cite{cieFlHofer},  \cite{bourgeooancea})
and wrapped Floer homology (see \cite{abouzaiseidel}), we need to study the case
when the critical submanifold data $\mathcal C$
varies. We can actually study such a situation in a similar way.
However we do not try to work it out in this article because this
is an expository article
providing the technical detail of the results which had been established in the previously published
literatures.
\end{rem}
\begin{lemdef}\label{lemdef1437}
\begin{enumerate}
\item
We can compose morphisms of inductive system of partial linear K-systems.
\item
Composition of homotopic morphisms are homotopic.
\item
Homotopy between morphisms is an equivalence relation.
\end{enumerate}
\end{lemdef}
\begin{proof}
(1)
Let
$(\{\frak N^{i}_{j+1 j}\},\{\frak H_{(j+1 j)}^i\})
: \mathcal {FF}_j \to \mathcal {FF}_{j+1}$
be morphisms for $j=1,2$.
We put
$$
\frak N^{i}_{31} = \frak N^{i}_{32} \circ \frak N^{i}_{21}.
$$
Then
$\frak H_{(32)}^i \circ \frak N^{i}_{21}$
is a homotopy from
$\frak N^i_3\circ \frak N^{i}_{31}
= \frak N^i_3\circ \frak N^i_{32} \circ \frak N^i_{21}$ to
$\frak N^{i+1}_{32} \circ \frak N^{i}_{2} \circ \frak N^{i}_{21}$
and
$\frak N_{32}^{i+1} \circ \frak H^{i}_{(21)}$
is a homotopy from
$\frak N^{i+1}_{32} \circ \frak N^{i}_{2} \circ \frak N^{i}_{21}$
to $\frak N^{i+1}_{31} \circ \frak N^{i}_{1}
= \frak N^{i+1}_{32} \circ \frak N^{i+1}_{21} \circ \frak N^{i}_{1}$.
Therefore we can construct $\frak H_{(31)}^i$
by gluing
$\frak H_{(32)}^i \circ \frak N^{i}_{21}$ and
$\frak N_{32}^{i+1} \circ \frak H^{i}_{(21)}$.
(See Subsubsection
\ref{glueinterpolation} for gluing process.
.)
See Figure \ref{Figure15-3}.
\begin{figure}[h]
\centering
\includegraphics[scale=0.25]{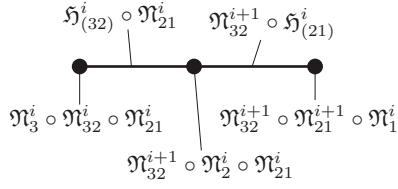}
\caption{composition of homotopies}
\label{Figure15-3}
\end{figure}
\par
(2) Let $(\{\frak N^{i}_{k,21}\},\{\frak H_{(k,21)}^i\})
: \mathcal {FF}_1 \to \mathcal {FF}_{2}$
be morphisms for $k=a,b$ and
let  $(\{\frak N^{i}_{32}\},\{\frak H_{(32)}^i\})
: \mathcal {FF}_2 \to \mathcal {FF}_{3}$
be a morphism.
Let $(\{\frak H_{(ab)}^i\},\{\mathcal H^i\})$
be a homotopy from $(\{\frak N^{i}_{a,21}\},\{\frak H_{(a,21)}^i\})$
to
$(\{\frak N^{i}_{b,21}\},\{\frak H_{(b,21)}^i\})$.
Then we have
\begin{equation}
\aligned
\partial (\frak N^{i+1}_{32} \circ \mathcal H^i)
=
&\frak N^{i+1}_{32} \circ \frak H_{(ab)}^{i+1} \circ \frak N^{i}_{1}
\cup
\frak N^{i+1}_{32} \circ \frak H_{(b,21)}^i \\
&\cup
\frak N^{i+1}_{32} \circ  \frak N^{i}_{2} \circ \frak H_{(ab)}^{i}
\cup
\frak N^{i+1}_{32} \circ \frak H_{(a,21)}^i
\endaligned
\end{equation}
and
\begin{equation}
\aligned
\partial (
\frak H_{(32)}^i \circ \frak H_{(ab)}^{i})
=
&\frak N^{i+1}_{32} \circ  \frak N^{i}_{2} \circ \frak H_{(ab)}^{i}
\cup \frak H^{i}_{32} \circ  \frak N^{i}_{b,21}
\\
&\cup
\frak N^{i}_{3} \circ \frak N^{i}_{32}  \circ \frak H_{(ab)}^{i}
\cup \frak H^{i}_{32} \circ  \frak N^{i}_{a,21}.
\endaligned
\end{equation}
We note that both are $[0,1]^2$-parametrized families of morphisms.
\smallskip
By the definition of composition we have
$$
\aligned
\frak N^{i+1}_{32} \circ \frak H_{(a,21)}^i
&\cup
\frak H^{i}_{32} \circ  \frak N^{i}_{a,21}
=
\frak H^i_{(a,31)}, \\
\frak N^{i+1}_{32} \circ \frak H_{(b,21)}^i
&\cup
\frak H^{i}_{32} \circ  \frak N^{i}_{b,21}
=
\frak H^i_{(b,31)}.
\endaligned
$$
We put
$$
\aligned
\frak H_{(ab)}^{i\prime}
&= \frak N^{i}_{32} \circ \frak H_{(ab)}^{i} \\
\mathcal H^{i\prime} &=
(\frak N^{i+1}_{32} \circ \mathcal H^i)
\cup_{\frak N^{i+1}_{32} \circ  \frak N^{i}_{2} \circ \frak H_{(ab)}^{i}}
(
\frak H_{(32)}^i \circ \frak H_{(ab)}^{i}).
\endaligned$$
See Figure \ref{Figure15-4}.
Here in the right hand side of the second formula,
we glue two $[0,1]^2$-parametrized morphisms
along one of the components of their boundaries.
We can do it in the same way as the
proof of Lemma-Definition \ref{lemdef1434}. See Subsubsection
\ref{glueinterpolation} for detail.
We can then easily see that
$$
\partial(\mathcal H^{i\prime})
=
\frak H_{(a,31)}^i \cup
\frak H_{(b,31)}^i \cup
(\frak N^{i}_{3} \circ \mathfrak H_{(ab)}^{i}) \cup
(\mathfrak H_{(ab)}^{i+1} \circ \frak N^{i}_{1}).
$$
Namely
$
(\{\frak N^{i}_{32}\},\{\frak H_{(32)}^i\})
\circ (\{\frak N^{i}_{a,21}\},\{\frak H_{(a,21)}^i\})$
is homotopic to
$
(\{\frak N^{i}_{32}\},\{\frak H_{(32)}^i\})
\circ (\{\frak N^{i}_{b,21}\},\{\frak H_{(b,21)}^i\})$.
\smallskip
\begin{figure}[h]
\centering
\includegraphics[scale=0.4,angle=-90]{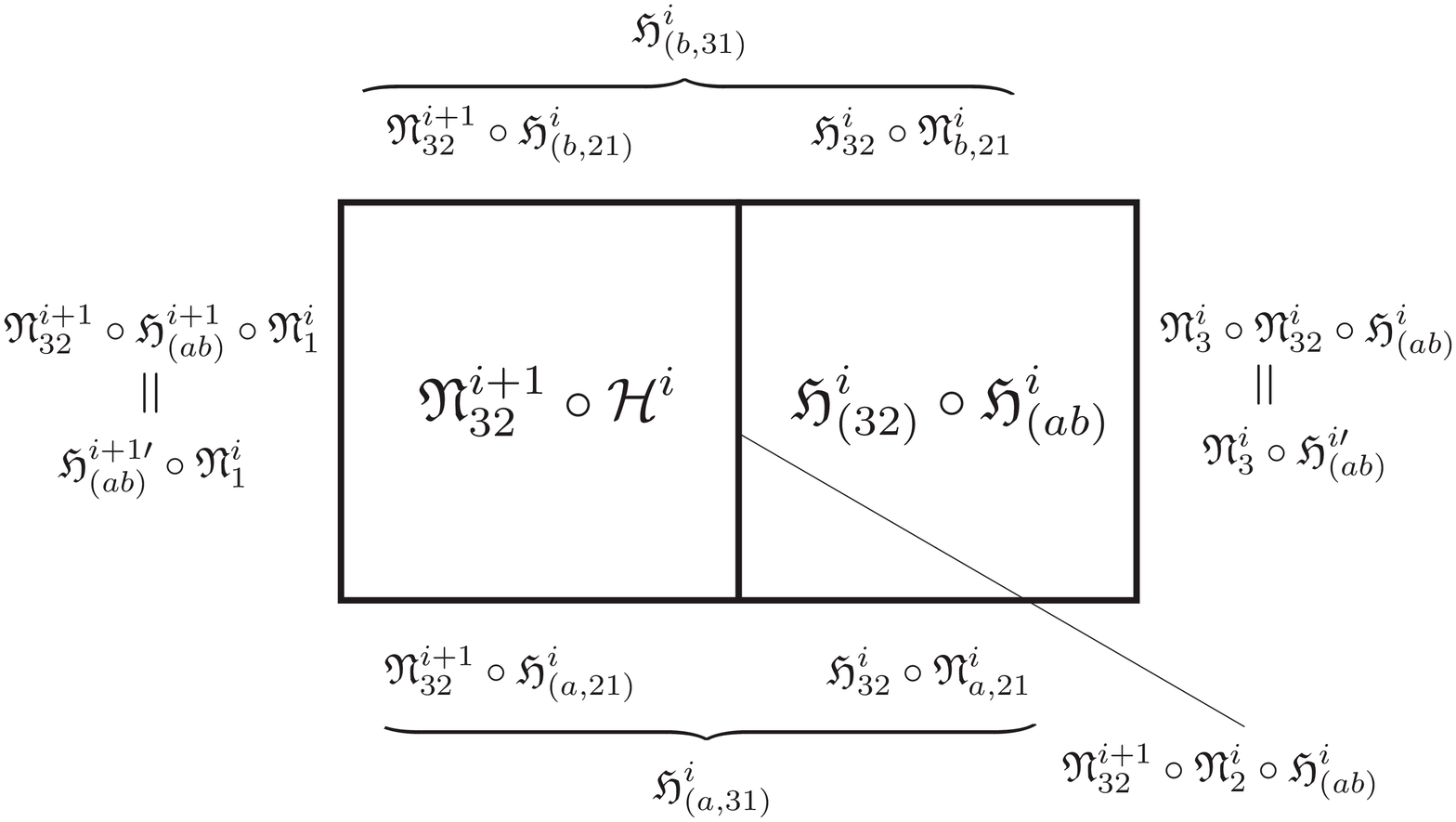}
\caption{$\mathcal H^{i\prime}$}
\label{Figure15-4}
\end{figure}
The case of the composition of homotopies and morphisms
in the opposite direction is similar.
\par
(3) will be proved in Subsubsection \ref{glueinterpolation}.
\end{proof}

\begin{thm}\label{indlinesysmainth1}
Let $\mathcal {FF} = (\{E^i\},\{\mathcal F^i\},\{\frak N^i\})$
be an inductive system of partial linear K-systems.
Note that the $\Lambda_{0,{\rm nov}}$ module
$CF(\mathcal{F}^i;\Lambda_{0,\rm nov})$ is independent of $i$,
which we denote by
$CF(\mathcal{FF};\Lambda_{0,\rm nov})$.
\begin{enumerate}
\item
To our inductive system of partial linear K-systems
$\mathcal{FF}$,
we can associate a map $d : CF(\mathcal{FF};\Lambda_{0,\rm nov}) \to CF(\mathcal{FF};\Lambda_{0,\rm nov})$
such that:
\begin{enumerate}
\item $d \circ d = 0$.
\item
There exists $\epsilon > 0$ such that:
$$
(d-d_0)(\frak F^{\lambda}CF(\mathcal{FF};\Lambda_{0,\rm nov}))
\subset \frak F^{\lambda+\epsilon}CF(\mathcal{FF};\Lambda_{0,\rm nov}).
$$
See \eqref{eq:deRham0} for the definition of $d_0$.
\end{enumerate}
\item
The definition of the map $d$ in (1) involves various choices and
$d$ depends on them. However it is independent of such choices
in the following sense:
Suppose $d_1$, $d_2$ are obtained by two different choices.
Then there exists a map
$\psi : CF(\mathcal{FF};\Lambda_{0,\rm nov}) \to CF(\mathcal{FF};\Lambda_{0,\rm nov})$
with the following properties.
\begin{enumerate}
\item $d_2 \circ \psi = \psi \circ d_1$.
\item $\psi$ is degree $0$ and preserves the filtration.
\item
There exists $\epsilon > 0$ such that:
$$
(\psi-{\rm id})(\frak F^{\lambda}CF(\mathcal{FF};\Lambda_{0,\rm nov}))
\subset \frak F^{\lambda+\epsilon}CF(\mathcal{FF};\Lambda_{0,\rm nov})
$$
where ${\rm id}$ is the identity map.
\item
In particular, $\psi$ induces an isomorphism
on cohomologies:
$$
H(CF(\mathcal{FF};\Lambda_{0,\rm nov}),d_1) \to
H(CF(\mathcal{FF};\Lambda_{0,\rm nov}),d_2).
$$
\item
$\psi$ depends on various choices but it is independent of the
choices up to cochain homotopy.
\end{enumerate}
\item
A morphism $\frak{NN} : \mathcal{FF}_1 \to \mathcal{FF}_2$ induces
a $\Lambda_{\rm nov}$ module homomorphism
\begin{equation}
\psi_{\frak{NN}} : CF(\mathcal{FF}_1;\Lambda_{0,\rm nov})
\to CF(\mathcal{FF}_2;\Lambda_{0,\rm nov})
\end{equation}
with the following properties.
\begin{enumerate}
\item
$\psi_{\frak{NN}}  \circ d = d \circ \psi_{\frak{NN}}$,
where $d$ in the left hand side (resp. right hand side)
is the coboundary oprator of
$CF(\mathcal{FF}_1;\Lambda_{0,\rm nov})$ (resp.
$CF(\mathcal{FF}_2;\Lambda_{0,\rm nov})$.)
\item
$\psi_{\frak{NN}}$ preserves degree and
$$
\psi_{\frak{NN}}\left(\frak F^{\lambda}CF(\mathcal{FF}_1;\Lambda_{\rm nov})
\right) \subset
\frak F^{\lambda-c}CF(\mathcal{FF}_1;\Lambda_{\rm nov})
$$
where $c$ is the energy loss of $\frak{NN}$.
In particular, if $c=0$ then $\psi_{\frak{NN}}$ induces a
$\Lambda_{0,\rm nov}$ module homomorphism :
$$
\psi_{\frak{NN}} : CF(\mathcal{FF}_1;\Lambda_{0,\rm nov})
\to CF(\mathcal{FF}_2;\Lambda_{0,\rm nov}).
$$
\end{enumerate}
\item
Let $\frak{NN}_a$,
$\frak{NN}_b : \mathcal{FF}_1 \to \mathcal{FF}_2$ be morphisms.
If $\frak{HH}$ is a homotopy from $\frak{NN}_a$ to $\frak{NN}_b$,
it induces a  $\Lambda_{\rm nov}$ module homomorphism:
\begin{equation}
\psi_{\frak{HH}} : CF(\mathcal{FF}_1;\Lambda_{\rm nov})
\to CF(\mathcal{FF}_2;\Lambda_{\rm nov})
\end{equation}
with the following properties.
\begin{enumerate}
\item
$\psi_{\frak{HH}}  \circ d + d \circ \psi_{\frak H} =\psi_{\frak{NN}_1}
- \psi_{\frak{NN}_2}$.
\item
$\psi_{\frak{HH}}$ decreases degree by $1$ and
$$
\psi_{\frak{HH}}\left(\frak F^{\lambda}CF(\mathcal{FF}_1;\Lambda_{\rm nov})
\right) \subset
\frak F^{\lambda-c}CF(\mathcal{FF}_1;\Lambda_{\rm nov})
$$
where $c$ is the energy loss of $\frak{HH}$.
In particular, if $c=0$ then $\psi_{\frak{HH}}$ induces a
$\Lambda_{0,\rm nov}$ module homomorphism :
$$
\psi_{\frak{HH}} : CF(\mathcal{FF}_1;\Lambda_{0,\rm nov})
\to CF(\mathcal{FF}_2;\Lambda_{0,\rm nov}).
$$
\end{enumerate}
\item
Let $\mathcal{FF}_a$, $\mathcal{FF}_b$, $\mathcal{FF}_c$
be inductive systems of partial linear K-systems
and $\frak{NN}_{ba} : \mathcal{FF}_a \to \mathcal{FF}_{b}$,
$\frak{NN}_{cb} : \mathcal{FF}_b \to \mathcal{FF}_{c}$
be morphisms. Then we have
$$
\psi_{\frak{NN}^{cb}\circ \frak{NN}^{ba}}
\sim
\psi_{\frak{NN}^{cb}} \circ \psi_{\frak{NN}^{ba}}
$$
where $\sim$ means cochain homotopic.
\item
The identity morphism induces a cochain map which is cochain homotopic to the identity map.
\end{enumerate}
\end{thm}
The proof will be given in Section \ref{sec:systemline3}.
The notion of identity morphism which appears in Item (6)
will be defined in Section \ref{subsec:identitylinsys}.
\begin{rem}
The morphism $\psi_{\frak{NM}}$ in Item (3) depends on the choices made for its construction.
More precisely, we first make choices mentioned in Item (2)
to define coboundary oprators of
$CF(\mathcal{FF}_i;\Lambda_{\rm nov})$, $i=1,2$.
Then we make a choice (compatible with the first choices)
to define the map $\psi_{\frak{NM}}$ which is
a cochain map with respect to the coboundary operators
obtained from the choices we made.
Note that we can use Item (4) to show that up to cochain homotopy the cochain map
$\psi_{\frak{NM}}$ is independent of the choices
we need to make to define it.
\end{rem}
\begin{defn}
We call the cohomology group of
$(CF(\mathcal{FF};\Lambda_{0,\rm nov}),d)$
the {\it Floer cohomology} of an inductive system
of partial linear K-systems $\mathcal{FF}$
\index{Floer cohomology ! of inductive system
of the partial linear K-systems, $HF(\mathcal{FF};\Lambda_{0,\rm nov})$}
\index{K-system ! Floer cohomology of inductive system
of the partial linear K-systems, $HF(\mathcal{FF};\Lambda_{0,\rm nov})$}
and denote the cohomology by $HF(\mathcal{FF};\Lambda_{0,\rm nov})$.
\end{defn}

\section{Extension of Kuranishi structure and its perturbation from
boundary to its neighborhood}
\label{sec:triboundary}

\subsection{Introduction to Section \ref{sec:triboundary}}
\label{subsection:introextensionlem}

In Section \ref{sec:systemline1} we formulated various versions of
corner compatibility conditions.
To prove the results stated in Section \ref{sec:systemline1}
(and which also appear in other places of this article and will appear in future)
we need to extend the Kuranishi structure given on the boundary
$\partial X$ satisfying corner compatibility conditions to one on
$X$. More precisely, we will start with the situation where we have
a Kuranishi structure $\widehat{\mathcal U}$ on $X$ and a Kuranishi structure
$\widehat {\mathcal U^+_\partial}$ on $\partial X$
such that
$$
\partial\widehat{\mathcal U}
< \widehat{\mathcal U^+_{\partial}}
$$
and want to find a Kuranishi structure $\widehat{\mathcal U^+}$ on $X$
such that
$$
\partial\widehat{\mathcal U^+}
=\widehat{\mathcal U^+_{\partial}},
\quad
\widehat{\mathcal U} <
\widehat{\mathcal U^+}.
$$
We also need a similar statement for CF-perturbations.
\par
It is rather easy to show that
if a Kuranishi structure is defined on a neighborhood $\Omega$
of a compact set $K$ then we can extend the Kuranishi structure without
changing it in a neighborhood of $K$.
(See Lemma \ref{lem1544}.)
This statement however is not enough to prove the
existence of  extension in the
above situation, since there we are given
a Kuranishi structure on $\partial X$ only and not on its
neighborhood.
In this section, we discuss the problem of extending a Kuranishi structure
and  a CF-perturbations on $\partial X$ to its neighborhood.

\begin{rem}
If we carefully examine the whole proofs of the geometric
applications appeared in previous literatures such as
\cite{FO},\cite{fooobook},\cite{fooobook2}, we will find out that, for the Kuranishi structure $\widehat{\mathcal U^+_\partial}$ we actually use,
an extension to its small neighborhood of $\partial X$ in $X$
 is given from its construction.
(This is the reason why in our earlier papers we did not elaborate on the statement
and its proof given in this section.)
In fact,
in the actual situation of applications, we start from a Kuranishi structure
$\partial\widehat{\mathcal U}$ on the boundary
(which we obtain from geometry or analysis) and construct
a good coordinate system  ${\widetriangle{\mathcal U_{\partial}}}$
and use it to find $\widehat{\mathcal U^+_{\partial}}$ and its perturbation.
We need to extend $\widehat{\mathcal U^+_{\partial}}$ and
its perturbation to a neighborhood of $\partial X$.
In this situation, the Kuranishi charts of
$\widehat{\mathcal U^+_{\partial}}$ are obtained as open subcharts
of certain Kuranishi charts of ${\partial}{\widetriangle{\mathcal U}}$.
(See the proof of \cite[Theorem 3.30, Proposition 6.44]{part11}.)
Therefore it can indeed be  extended using the extension of ${\partial}\widehat{\mathcal U}$, directly.
\par
Nevertheless, the reason why we prove these extension results in this article
is as follows.
In this article we want to present various parts of the proofs of the
whole story into a `package' as much as possible.
In other words, we want to decompose whole proofs of the geometric results
into pieces, so that each piece can be stated and proved rigorously
and independently from other pieces.
Namely, for each divided `package',
we want to state the precise assumption and conclusion
together with the proof.
In this way, one can use and quote each `package' without
referring their proofs.
We want to do so because the whole story
has now grown huge which becomes
harder to follow all at once.
\par
For this purpose, we want to specify and restrict
the information which we `remember' at each
step of the inductive construction of the
K-system and its perturbations.
When we construct a system of virtual fundamental
chains of, for example, $\mathcal M(\alpha_-,\alpha_+)$,
we use, as an induction hypothesis, a certain structure
on $\mathcal M(\alpha,\alpha')$ for
$E(\alpha_-) \le E(\alpha) < E(\alpha') \le E(\alpha_+)$
(and $E(\alpha') - E(\alpha) < E(\alpha_+) - E(\alpha_-)$)
and use only those structures during our construction.
We need to make a careful choice of the structures
we `remember' during the inductive steps,
when we proceed from one step of the induction to the
next,
so that the induction works.
Our choice in this article is that we `remember'
Kuranishi structure and its perturbation
but forget the good coordinate system (which we use to
construct the perturbation) and objects on it.
We also forget the way how those Kuranishi structures
and various objects on them
are constructed, but explicitly list up all the properties
we use in the next step of the induction.
Therefore the relation between
$\widehat{{\mathcal U}_{\partial}^+}$
and $\partial\widehat{\mathcal U}$
(except $\partial\widehat{\mathcal U}
< \widehat{{\mathcal U}_{\partial}^+}$)
is among the
data which we forget when we proceed to the next step of the induction.
\end{rem}
The simplest version of the extension result we mentioned above is the following.
\par\medskip
{\bf ($\ast$)} Suppose we are given a continuous function $f$ on $\partial([0,1)^n)$.
We assume that the restriction of $f$ to
$[0,1)^k \times \{0\} \times [0,1)^{n-k-1}$
is smooth for any $k$.
Then $f$ is extended to a smooth function on $[0,1)^n$.
\par\medskip
The statement {\bf ($\ast$)} is, of course, classical.
(See \cite[Lemma 7.2.121]{fooo09} for its proof, for example.)
We can use this statement together with various
technique of manifold theory
(such as partition of unity and induction on the number of coordinate charts) to prove the
existence of an extension of CF-perturbation
if we include enough structure into the assumptions.
(However the proof is rather cumbersome to work out in detail.
In fact, it uses triple induction.
Two of the triple are indexed by the partially ordered set appearing in the
definition of good coordinate system. Another is the codimension of the starta
of corner structure stratifiation.)
\par
In this section we take a short cut in the following way.
Suppose $M$ is a manifold with boundary (but has no corner).
Then it is well-known that a neighborhood of $\partial M$
is identified with $\partial M \times [0,\epsilon)$.
If $M$ has corners,
we can identify a neighborhood of $\overset{\circ}S_k(M)$
in $M$ with a twisted product
$$
\overset{\circ}S_k(M) \tilde\times [0,1)^{n-k}
$$
that is a fiber bundle over
$\overset{\circ}S_k(M) = S_k(M) \setminus S_{k+1}(M)$ (see \cite[Definition 4.13]{part11}),
whose structure group is a finite group of permutations of
$(n-k)$ factors.
We also require a compatibility of these structures for various $k$.
Such a structure may be regarded as a special case of the
`{\it system of tubular neighborhoods}'\index{system of tubular neighborhoods} of a stratified space introduced by Mather \cite{Math73}.
Its existence is claimed and can be proved in the same way as in \cite{Math73}.
(See also \cite[Proposition 8.1]{fooooverZ}.)
\par
We can include certain `topological' objects such as vector bundle
(especially obstruction bundle) and generalize the notion
of trivialization of the structure in a neighborhood of
boundary and corner, and prove its existence
without assuming extra conditions.
However, we can {\it not} expect the Kuranishi map
respects this trivialization. Namely the Kuranishi map
in general may not be constant in the $[0,1)^k$ factor.
So we need some discussion to use it on the way how we extend
perturbation to a neighborhood of $\partial X$.
\par
The main idea we use in this section for the short cut is summarized as follows.
In case when $M$ is a manifold with boundary,
we can attach $\partial M \times [-1,0]$ to $M$
by identifying $\partial M \times \{0\}$ with $\partial M \subset M$,
and enhance $M$ to a manifold $M^{\boxplus 1}$ with boundary
so that a neighborhood of its boundary is {\it canonically} identified with $\partial M \times [-1,0)$.
Then we can extend a vector bundle $E$ on $M$ to $M^{\boxplus 1}$
so that its restriction to this neighborhood of the boundary is
{\it canonically} identified with $E\vert_{\partial M} \times [-1,0)$.
Then its section such as Kuranishi map can be extended
to $M^{\boxplus 1}$ so that on $\partial M \times [-1,0)$ it is
constant in the $[-1,0)$ direction.
In other words, {\it we attach the collar `{\bf outside}' of $M$ in place
of constructing it `{\bf inside}'.}
\par
We can perform a similar construction in the case when $M$
has a corner, the case of orbifold, and the case of K-space.
Then we replace the moduli spaces such as $\mathcal M(\alpha_-,\alpha_+)$
by $\mathcal M(\alpha_-,\alpha_+)^{\boxplus 1}$
so that the fiber product description of their boundaries remains the same and
their Kuranishi structures have {\it canonical trivializations} near the boundary.
We can use them to extend the Kuranishi structures and
their perturbations given on the boundary to a
neighborhood of the boundary.
\par
There is a slight issue about the smoothness of the structure
at the point $\partial M \times \{0\}$.
The shortest and simplest way to resolve the issue is to use {\it exponential
decay estimate} of various objects appearing there
in geometric situations, which we proved in \cite[Lemma A1.58]{fooobook2}, \cite[Theorems 13.2 and 19.5]{foootech}, \cite{foooexp}.
In the abstract setting,
we impose certain {\it admissibility} on
orbifold, vector bundle, etc.,
to incorporate such exponential decay property.
See Section \ref{sec:admKura}.
\begin{rem}
In fact, we can work in the piecewise smooth category
for the purpose of
putting the collar outside. However,
working in the smooth category rather than the piecewise smooth category
reduces the amount of checks we need
during various constructions.
\end{rem}
We will carry out this idea in detail in this section.
It is rather lengthy since we need to write down and check various compatibilities of many
objects with our trivialization.
However, the compatibilities in almost all cases are
fairly obvious.
In other words, the proofs are nontrivial
only because they are lengthy.

\subsection{Trivialization of corners on one chart}
\label{subsec:trivonechart}
In this section, for a Kuranishi chart $\mathcal U =(U,\mathcal E, \psi ,s)$,
we will define in Lemma-Definition \ref{lem15444} a certain Kuranishi chart
$\mathcal U^{\boxplus \tau}$ at a point $x$ of corners of $U$, which is called
{\it trivialization of corners}. Geometrically speaking,
`trivialization' means `trivialization of normal bundle of corners',
which gives us certain coordinates on the corners.
We consider the following situation in this subsection.
\begin{shitu}\label{situ151}
Let $\mathcal U = (U,\mathcal E,\psi, s)$ be a Kuranishi neighborhood of $X$
and
$$
x \in \overset{\circ}S_k(U).
$$
Let $\frak V_x = (V_x,E_x,\Gamma_x,\phi_x,\widehat\phi_x)$
be an orbifold chart of $(U,\mathcal E)$ at $x$.
Recall from
Definition \ref{defn2655} (also Definition \ref{defn2613}) that
$V_x$ is a manifold and
$$
\phi_x : V_x \to X, \quad \widehat\phi_x : V_x \times E \to \mathcal E.
$$
Let $s_x$ be a representative of the Kuranishi map $s$ on
$\frak V_x$.
We assume that $V_x$ is an open subset of the direct product $\overline{V_x} \times
[0,1)^k$ where $\overline{V_x}$ is a manifold without boundary
and that $o_x \in V_x$ corresponds to $(\overline o_x,(0,\dots,0))
\in \overline{V_x} \times
[0,1)^k$.
For $y \in V_x \subset \overline{V_x} \times [0,1)^k$ we denote by
$\overline{y} \in \overline{V_x}$ the $\overline{V_x}$-component of $y$.
See Lemma \ref{lem26999} for this notation.
$\blacksquare$
\end{shitu}
\begin{defn}\label{defn153}
Under Situation \ref{situ151},
we define the {\it retraction map} \index{retraction map}
$$
\mathcal R_x : \overline{V_x} \times (-\infty,1)^k \to  \overline{V_x} \times [0,1)^k
$$
by
$$
\mathcal R_x(\overline y,(t_1,\dots,t_k))
= (\overline y,(t'_1,\dots,t'_k)),
$$
where
$$
t'_i =
\begin{cases}
t_i  &\text{if $t_i \ge 0$,}\\
0    &\text{if $t_i \le 0$.}
\end{cases}
$$
For $\tau > 0$,
we define an open subset $V_x^{\boxplus\tau}$ of
$\overline{V_x} \times [-\tau,1)^k$
by
$\mathcal R_x^{-1}(V_x) \cap (\overline{V_x} \times [-\tau,1)^k$).
Here the retraction map
$\mathcal R_x$ naturally induces a map
$$\mathcal R_x : V_x^{\boxplus \tau} \to V_x
\subseteq \overline{V_x} \times
[0,1)^k.$$
\end{defn}
Next, we define a $\Gamma_x$ action on  $V_x^{\boxplus \tau}$.
\begin{figure}[h]
\centering
\includegraphics[scale=0.4]{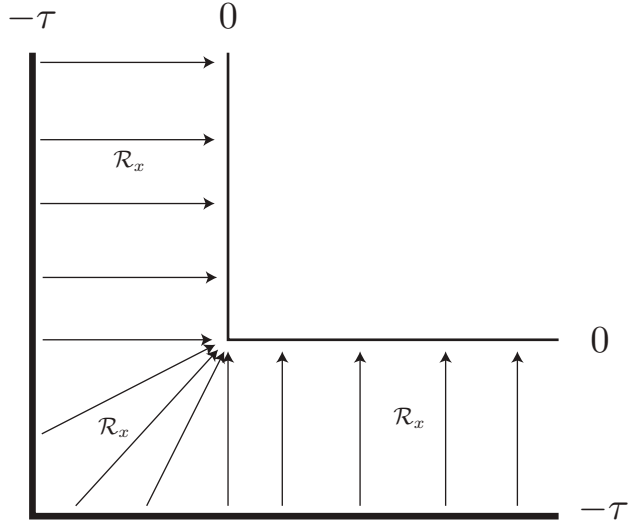}
\caption{the retraction $\mathcal R^x$}
\label{Figure16-1}
\end{figure}
\begin{defn}\label{defn154444}
Let $\gamma \in \Gamma_x$.
By the definition of admissible orbifold
(Definition \ref{defn297} (1)(b)), there exists $\sigma_{\gamma}
\in {\rm Perm}(k)$ with the following properties.
We define maps
$\varphi_{0}^{\gamma} : \overline V_x \times [0,1)^k
\to \overline V_x$ and
$\varphi_i^{\gamma} : \overline V_x \times [0,1)^k
\to [0,1)$ requiring
the components of $\gamma\cdot (\overline y,(t_1,\dots,t_k))$ to be
$$
\gamma\cdot (\overline y,(t_1,\dots,t_k))
=
\left(\varphi_{0}^{\gamma} (\overline y,(t_1,\dots,t_k)),
\left(\varphi_{\sigma_{\gamma}^{-1}(i)}^{\gamma}(\overline y,(t_1,\dots,t_k))\right)_{i=1}^k\right).
$$
Here the left hand side is the given $\Gamma_x$ action
on $V_x$. Then we find the following:
\begin{enumerate}
\item
The map $\varphi_{0}^{\gamma}$ is admissible.
\item
The function $(\overline y,(t_1,\dots,t_k)) \mapsto \varphi_i^{\gamma}
(\overline y,(t_1,\dots,t_k)) - t_i$
is exponentially small at the boundary.
\end{enumerate}
The existence of such $\varphi_{0}^{\gamma}$,
$\varphi_i^{\gamma}$ is nothing but the definition of an admissible action.
Then we define
$\widehat\varphi_{0}^{\gamma} : \overline V_x \times [-\tau,1)^k
\to \overline V_x$ and
$\widehat\varphi_i^{\gamma} : \overline V_x \times [-\tau,1)^k
\to [-\tau,1)$
as follows.
\begin{enumerate}
\item[(A)]
$\widehat\varphi_{0}^{\gamma}
= \varphi_{0}^{\gamma}\circ \mathcal R_x$.
\item[(B)]
$$
\widehat\varphi_{i}^{\gamma}(\overline y,(t_1,\dots,t_k))
=
\begin{cases}
t_i &\text{if $t_i \le 0$,}\\
(\varphi_{i}^{\gamma}\circ \mathcal R_x)(\overline y,(t_1,\dots,t_k))
&\text{if $t_i \ge 0$.}
\end{cases}
$$
\end{enumerate}
We now put
$$
\gamma\cdot (\overline y,(t_1,\dots,t_k))
=
\left(\widehat\varphi_{0}^{\gamma} (\overline y,(t_1,\dots,t_k)),
\left(\widehat\varphi_{\sigma_{\gamma}^{-1}(i)}^{\gamma}(\overline y,(t_1,\dots,t_k))\right)_{i=1}^k\right).
$$
We note that $y \mapsto \gamma \cdot y$ coincides with the given $\gamma$ action
if $y \in V_x$.
\end{defn}
We now prove that this defines a $\Gamma_x$-action on $V_x^{\boxplus \tau}$.
We first show:
\begin{lem}\label{lem155222}
We have $\widehat\varphi_{i}^{\gamma}(\overline y,(t_1,\dots,t_k)) > 0$
if and only if $t_i >0$.
Moreover $\widehat\varphi_{i}^{\gamma}(\overline y,(t_1,\dots,t_k)) = 0$
if and only if $t_i =0$.
\end{lem}
\begin{proof}
If $t_i \le 0$ then
$\widehat\varphi_{i}^{\gamma}(\overline y,(t_1,\dots,t_k)) = t_i \le 0$ by (B).
On the other hand, if $t_i > 0$ then
$\widehat\varphi_{i}^{\gamma}(\overline y,(t_1,\dots,t_k))
= (\varphi_{i}^{\gamma}\circ \mathcal R_x)(\overline y,(t_1,\dots,t_k))
> 0$ by (B). The lemma follows.
\end{proof}
\begin{lem}\label{lemma1556}
We have $\gamma \cdot (\mathcal R_x(y)) = \mathcal R_x(\gamma \cdot y)$
for any $y \in V_x^{\boxplus \tau}$.
\end{lem}
\begin{proof}
The coincidence of $\overline V_x$ coordinates is an immediate consequence
of the definition.
We will check the coincidence of the ${\sigma_{\gamma}(i)}$-th
coordinates of $[-\tau,1)^k$ factors.
\par
If $t_i \le 0$, then this coordinate is $0$ for both
$\gamma \cdot (\mathcal R_x(y))$ and $\mathcal R_x(\gamma \cdot y)$.
\par
If $t_i > 0$, they both are the ${\sigma_{\gamma}(i)}$-th coordinates
of the $[-\tau,1)^k$ factor of $\gamma \cdot (\mathcal R_x(y))$.
\par
In fact, this is obvious for  $\gamma \cdot (\mathcal R_x(y))$.
To show this claim for $\mathcal R_x(\gamma \cdot y)$,
we first observe that the ${\sigma_{\gamma}(i)}$-th coordinate
of the $[-\tau,1)^k$ factor of $\gamma \cdot y$
is the ${\sigma_{\gamma}(i)}$-th coordinate
of the $[-\tau,1)^k$ factor of $\gamma \cdot (\mathcal R_x(y))$
by (B). Then we use Lemma \ref{lem155222} to show
that $\mathcal R_x$ does not change this coordinate.
\end{proof}
\begin{lem}\label{lem15777}
We have
$
\gamma \cdot (\mu \cdot y)
= (\gamma \mu) \cdot y
$
for any $y \in V_x^{\boxplus \tau}$ and $\gamma,\mu \in \Gamma_x$.
\end{lem}
\begin{proof}
We calculate
$$
\widehat\varphi^{\gamma}_0(\mu \cdot y)
=
\varphi^{\gamma}_0(\mathcal R_x(\mu \cdot y))
=
\varphi^{\gamma}_0(\mu \cdot \mathcal R_x(y))
=
\widehat\varphi^{\gamma\mu}_0(y).
$$
Therefore the $\overline V_x$ factor coincides.
\par
Let $y = (\overline y,(t_1,\dots,t_k))$.
Suppose $t_i > 0$. Then
$\widehat\varphi_{i}^{\mu}(\overline y,(t_1,\dots,t_k))) > 0$
by Lemma \ref{lem155222}.
Therefore we obtain
$$
\widehat\varphi^{\gamma}_{\sigma_{\mu}(i)}(\mu \cdot y)
=
\varphi^{\gamma}_{\sigma_{\mu}(i)}(\mathcal R_x(\mu \cdot y))
=
\varphi^{\gamma}_{\sigma_{\mu}(i)}(\mu \cdot \mathcal R_x(y))
=
\widehat\varphi^{\gamma\mu}_i(y).
$$
Suppose  $t_i \le 0$. Then
$\widehat\varphi_{i}^{\mu}(\overline y,(t_1,\dots,t_k))) = t_i \le 0$
by Lemma \ref{lem155222}.
Thus we get
$$
\widehat\varphi^{\gamma}_{\sigma_{\mu}(i)}(\mu \cdot y)
=
t_i
=
\widehat\varphi^{\gamma\mu}_i(y).
$$
The proof of Lemma \ref{lem15777} is complete.
\end{proof}
We have thus defined a $\Gamma_x$ action on
$V_x^{\boxplus\tau}$.
\begin{lem}\label{Vplusissmooth}
\begin{enumerate}
\item
For any admissible map $f : V_x \to M$ the composition
$f \circ \mathcal R_x : V^{\boxplus\tau}_x \to M$ is smooth.
\item
If $f$ is a submersion, so is $f \circ \mathcal R_x$.
\item
The $\Gamma_x$ action on $V^{\boxplus\tau}_x$ is smooth and
$\mathcal R_x$ is $\Gamma_x$ equivariant.
\end{enumerate}
\end{lem}
\begin{proof}
(1).
Since the problem is local, it suffices to prove the case when $M = \R$.
We use Lemma \ref{canonicalformforadmi} to obtain
$f_I$ such that $f = \sum_{I}f_I$.
Here $I\subset \{1,\dots,k\}$ and $f_I : \overline V_x \times [0,1)^I \to \R$
is exponentially small near the boundary.
We define a function $\widehat f_I : V_x^{\boxplus\tau} \to \R$ by
$$
\widehat f_I(\overline y,(t_1,\dots,t_k))
=
\begin{cases}
\widehat f_I(\overline y,t_I)
&\text{if $t_i \ge 0$ for all $i \in I$,}
\\ 0
&\text{if $t_i \le 0$ for some $i \in I$.}
\end{cases}
$$
Since $f_I : \overline V_x \times [0,1)^I \to \R$
is exponentially small near the boundary,
$\widehat f_I$ is smooth.
It is easy to see that
$f\circ \mathcal R_x = \sum_I \widehat f_I$.
Therefore $f\circ \mathcal R_x$ is smooth.
\par
We note that the submersivity of $f$ by definition
implies the submersivity of the restriction of $f$ to
$\overset{\circ}S_k V_x$. (2) is immediate from this fact and
the construction.
\par
The first half of (3) follows from definition. The second half of (3) is
Lemma \ref{lemma1556}.
\end{proof}
We define
$U^{\boxplus\tau}_x$ to be the quotient $V^{\boxplus\tau}_x/\Gamma_x$ under the
$\Gamma_x$ action on $V^{\boxplus\tau}_x$ described above.
Then Lemma \ref{lem15777} yields the retraction map
$$
\mathcal R_x : U^{\boxplus\tau}_x \to U_x
$$
induced by the one on $V^{\boxplus\tau}_x$.
We denote
$$
\mathcal E_x = (E_x \times V_x)/\Gamma_x,
$$
which is a vector bundle on an orbifold $U_x =V_x/\Gamma_x$.
Then we define $\mathcal E^{\boxplus\tau}_x$ to be the pull-back
\begin{equation}\label{form151111}
\mathcal E^{\boxplus\tau}_x =
\mathcal R_x^*(\mathcal E_x)
= (E_x \times V^{\boxplus\tau}_x)/\Gamma_x,
\end{equation}
which is a smooth vector bundle on
an orbifold $U^{\boxplus\tau}_x$.
The Kuranishi map $s_x$ which is a section of $\mathcal E_x$
induces a section $s_x^{\boxplus\tau}$ of $\mathcal E^{\boxplus\tau}_x$.
In the same way as in the proof of Lemma \ref{Vplusissmooth} (1)
we can show that $s_x^{\boxplus\tau}$ defines a smooth section.
\par
We define
\begin{equation}\label{form151112}
(X \cap U_x)^{\boxplus\tau} := (s_x^{\boxplus\tau})^{-1}(0)/\Gamma_x,
\end{equation}
which is a paracompact Hausdorff space.
(We note that $X\cap U_x$ in the notation
$(X \cap U_x)^{\boxplus\tau}$ does NOT
stand for the set-theoretical intersection, but is just a notation.)
We have a map
$\psi_x^{\boxplus\tau} : (s_x^{\boxplus\tau})^{-1}(0)/\Gamma_x \to (X \cap U_x)^{\boxplus\tau}$
that is the identity map.
Then from the definition we find
\begin{lemdef}\label{lem15444}
Let $\mathcal U = (U, \mathcal E,\psi,s)$ be a Kuranishi chart of $X$  and
$x \in \overset{\circ}S_k(U)$ as in Situation \ref{situ151}. Then
$(U^{\boxplus\tau}_x=V^{\boxplus\tau}_x/\Gamma_x,\mathcal E^{\boxplus\tau}_x,\psi_x^{\boxplus\tau}, s_x^{\boxplus\tau})$
is a Kuranishi chart of $(X \cap U_x)^{\boxplus\tau}$.
We denote
$$
\mathcal U^{\boxplus\tau}_x=(U^{\boxplus\tau}_x,\mathcal E^{\boxplus\tau}_x,\psi^{\boxplus\tau}_x, s^{\boxplus\tau}_x)
$$
and call it {\rm trivialization of corners} of $\mathcal U_x$.
\index{corner ! trivialization of corners for one chart}
\index{trivialization of corners ! for one chart}
\end{lemdef}
Let $\mathcal S_x = (W_x,\omega_x,\frak s_x)$
be a CF-perturbation of $\mathcal U$ on $\frak V_x$.
(\cite[Definition 7.5]{part11}.)
We define $\frak s_x^{\boxplus\tau} : V_x^{\boxplus\tau} \times W_x \to E_x$ by
\begin{equation}
\frak s_x^{\boxplus\tau}(y,\xi) = \frak s_x(\mathcal R_x(y),\xi).
\end{equation}
\begin{lemdef}\label{lemdef157}
\begin{enumerate}
\item
The boundary $\partial(U_x^{\boxplus\tau})$ (resp. $\partial(\mathcal S_x^{\boxplus\tau})$)  is canonically diffeomorphic
to $(\partial U_x)^{\boxplus\tau}$ (resp. $(\partial \mathcal S_x)^{\boxplus\tau}$).
\item
The triple
$(W_x,\omega_x,\frak s_x^{\boxplus\tau})$ is a CF-perturbation
of $(V^{\boxplus\tau}_x/\Gamma_x,\mathcal E^{\boxplus\tau}_x,s_x^{\boxplus\tau},\psi_x^{\boxplus\tau})$.
We denote it by $\mathcal S^{\boxplus\tau}_x$.
\item
If $\mathcal S_x$ is equivalent to $\mathcal S'_x$,
then
$\mathcal S^{\boxplus\tau}_x$ is equivalent to $\mathcal S^{\prime \boxplus\tau}_x$.
\item
If $f : U \to M$ is a smooth map and strongly
submersive on $\frak V_x$ with respect to $\mathcal S_x$,
then $f \circ \mathcal R_x$ is strongly
submersive with respect to $\mathcal S^{\boxplus\tau}_x$.
We denote it by $f^{\boxplus\tau}_x : U^{\boxplus\tau}_x \to M$.
\item
The strong transversality to $N \to M$ is also preserved.
\item
The versions of (2)-(5) where
`CF-perturbation' is replaced by
`multivalued perturbation' also hold.
\item
If $h_x$ is a differential form on $U_x$, that is,
$\widetilde h_x$ is a $\Gamma_x$-invariant differential
form on $V_x$, then $\widetilde {h}_x^{\boxplus\tau}
:=\mathcal R_x^*\widetilde h_x$ defines a
differential form on $U^{\boxplus\tau}_x$.
We denote it by $h^{\boxplus\tau}_x$.
\item
In the situation of (2)(4)(7),
we assume that $h_x$ is compactly supported. Then we have:
\begin{equation}\label{form1551}
(f_{x}^{\boxplus\tau})!(h_x^{\boxplus\tau};\mathcal S_x^{\boxplus\tau})
=
f_{x}!(h_x;\mathcal S_x).
\end{equation}
\end{enumerate}
\end{lemdef}
\begin{proof}
The proofs are mostly immediate from the definition. We only prove (8) for
completeness' sake.
To prove (8), we recall the definition of the push out from
\cite[Definition 7.10]{part11},
which
results in
\begin{equation}\label{intboxsquare}
\int_M f_x^{\boxplus\tau}!(h_x^{\boxplus\tau};\mathcal S_x^{\boxplus\tau}) \wedge \rho
 =  (-1)^{|\rho||\omega|}\int_{(s_x^{\boxplus\tau})^{-1}(0)}
\pi_1^* h^{\boxplus\tau}_x \wedge \pi_1^*(f^{\boxplus\tau}_x)^* \rho \wedge \pi_2^*\omega_x
\end{equation}
\begin{equation}\label{intnobox}
\int_M f_x!(h_x;\mathcal S_x) \wedge \rho
=  (-1)^{|\rho||\omega|}\int_{(s_x)^{-1}(0)}
\pi_1^*\widetilde h_x \wedge \pi_1^*f_x^* \rho \wedge \pi_2^*\omega_x
\end{equation}
for any differential form $\rho$ on $M$. Here $\pi_1: V_x^{\boxplus\tau} \times W_x \to V_x^{\boxplus\tau}$,
$\pi_2: V_x^{\boxplus\tau} \times W_x \to W_x$ are projections,
$\widetilde h^{\boxplus\tau}$ is a $\Gamma_x$-invariant differential form
on $V_x^{\boxplus\tau}$ defined by ${\mathcal R}_x^*\widetilde h_x$ and
$f^{\boxplus\tau} = f_x \circ {\mathcal R}_x : V_x^{\boxplus\tau} \to M$ is a submersion defined by
Lemma \ref{Vplusissmooth}.
\par
We compare the right hand sides of the above two integrals.
First we note $s_x^{\boxplus\tau} = s_x \circ (\mathcal R_x \times \text{\rm id}_{W_x})$ and
$\mathcal R_x|_{V_x} = \text{\rm id}_{V_x}$ on $\del V_x \subset V_x \subset V_x^{\boxplus\tau}$.
We decompose
$$
(s_x^{\boxplus\tau})^{-1}(0) = \left((s_x^{\boxplus\tau})^{-1}(0) \cap
(V_x \times W_x)\right)
\cup \left((s_x^{\boxplus\tau})^{-1}(0) \setminus (V_x \times W_x)\right),
$$
and its associated integral and note $s_x^{\boxplus\tau} \equiv s_x$ on $V_x \subset V_x^{\boxplus\tau}$.
Then contribution of the integral \eqref{intboxsquare}
over the first part of the domain becomes \eqref{intnobox}.
\par
It remains to check the contribution of \eqref{intboxsquare} on the region
$(s_x^{\boxplus\tau})^{-1}(0) \setminus (V_x \times W_x)$.
In the rest of the proof,
we may assume that $\rho$ is chosen
so that the degree of $\pi_1^*\widetilde h_x \wedge \pi_1^*f_x^* \rho \wedge \pi_2^*\omega_x$ matches the dimension of $(s_x)^{-1}(0)$.
\par
We note that the retraction $\mathcal R_x \times \text{\rm id}_{W_x}:V_x^{\boxplus\tau} \setminus V_x \times W_x\to \del V_x \times W_x$
also induces a retraction of $(s_x^{\boxplus\tau})^{-1}(0) \setminus (V_x^{\boxplus\tau} \times W_x)$ to
$s_x^{-1}(0) \cap (\del V_x \times W_x)$ with one-dimensional fiber by definition of $s_x^{\boxplus\tau}$.
We also note
$$
\widetilde h^{\boxplus\tau} \wedge f^{\boxplus\tau} \rho = (\mathcal R_x \times \text{\rm id}_{W_x})^*\eta
$$
for some form $\eta$ defined on $\del V_x$ by the definitions of $\widetilde h^{\boxplus\tau}$ and $\wedge f^{\boxplus\tau}$
given above. We derive
$$
\aligned
& \int_{(s_x^{\boxplus\tau})^{-1}(0) \setminus (V_x \times W_x)}
\pi_1^* h^{\boxplus\tau} \wedge \pi_1^*f^{\boxplus\tau} \rho \wedge \pi_2^*\omega_x \\
=  & \int_{(s_x^{\boxplus\tau})^{-1}(0) \setminus (V_x \times W_x)}
\pi_1^* (\mathcal R_x \times \text{\rm id}_{W_x})^* \eta \wedge \pi_2^*\omega_x\\
=& \int_{(s_x^{\boxplus\tau})^{-1}(0) \setminus (V_x \times W_x)}(\mathcal R_x \times \text{\rm id}_{W_x})^*
(\pi_1^*\eta \wedge \pi_2^*\omega_x)\\
=& \int_{s_x^{-1}(0) \cap (\del V_x \times W_x)} (\pi_1^*\eta \wedge \pi_2^*\omega_x).
\endaligned
$$
For the last equality we use
$$
(s_x^{\boxplus\tau})^{-1}(0) \setminus (V_x \times W_x) = (\mathcal R_x \times \text{\rm id}_{W_x})^{-1}
(s_x^{-1}(0) \cap (\del V_x \times W_x))
$$
by definition of $s_x^{\boxplus\tau}$.
By submersion property of $s_x$, we have
$$
\dim s_x^{-1}(0) \cap (\del V_x \times W_x) = \dim s_x^{-1}(0) - 1.
$$
Therefore by the degree assumption made on $\rho$ above, the last integral vanishes.
Now the proof of (\ref{form1551}) is complete.
\end{proof}
\begin{lem}\label{lem1599}
We put
$
\overset{\circ\circ}S_k(U_x^{\boxplus\tau})
=
S_k(U_x^{\boxplus\tau})
\cap \mathcal R_x^{-1}(\overset{\circ}S_k(U_x))
$.
\begin{enumerate}
\item
The closure of $\overset{\circ\circ}S_k(U_x^{\boxplus\tau})$
in $\overset{\circ}S_k(U_x^{\boxplus\tau})$
is an orbifold with corners.
\item
The map $\mathcal R_x$ induces an orbifold diffeomorphism
from
${\rm Clos}(\overset{\circ\circ}S_k(U_x^{\boxplus\tau}))$ to
$\widehat S_k(U_x)$.
\end{enumerate}
\end{lem}
\begin{proof}
It suffices to prove the lemma for the case when $V = [0,1)^k$, which is
obvious. (See Figure \ref{Figure16-2}.)
\begin{figure}[h]
\centering
\includegraphics[scale=0.3]{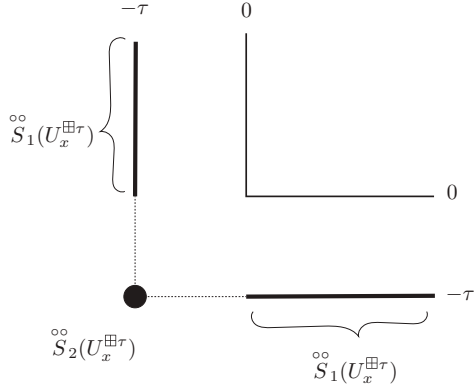}
\caption{$\overset{\circ\circ}S_k(U_x^{\boxplus\tau})$}
\label{Figure16-2}
\end{figure}
\end{proof}

\subsection{Trivialization of corners and embedding}
To shorten the discussion we study trivialization of the corner for
the case of coordinate change and of a local representative
of an embedding simultaneously.
In this subsection we consider the following situation.

\begin{shitu}\label{situ153}
Suppose we are in Situation \ref{situ151}.
Let
$x' \in \overset{\circ}S_k(U)$ and
$\frak V_{x'} = (V_{x'},E_{x'},\Gamma_{x'},\phi_{x'},\widehat\phi_{x'})$
be an orbifold chart of $(U',\mathcal E')$ at $x$.
(See Definition \ref{defn2613}.)
Let $s_{x'}$ be a representative of the Kuranishi map $s$ on
$\frak V_{x'}$.
Let
$(U',\mathcal E') \to (U,\mathcal E)$ be an embedding.
We take its local representative $(h_{xx'},\varphi_{xx'},\widehat{\varphi}_{xx'})$.
$\blacksquare$
\end{shitu}
This includes the case when $(h_{xx'},\varphi_{xx'},\widehat{\varphi}_{xx'})$
represents an isomorphism, which is nothing but
the case of coordinate change of orbifold.
\begin{lem}\label{lema1510}
In Situation \ref{situ153}, there exists a unique injective map
$\frak j : \{1,\dots,k'\} \to \{1,\dots, k\}$
with the following properties.
\begin{enumerate}
\item
For $\gamma' \in \Gamma_{x'}$,
 we have:
$
\sigma(\gamma')(\frak j(i))
=
\frak j(\sigma'(\gamma')(i))
$.
\item
If
$
\varphi_{xx'}(\overline y',(t'_1,\dots,t'_{k'}))
= (\overline y,(t_1,\dots,t_k))
$
then
$
t_{\frak j(i)} = 0
$
if and only if $t'_i = 0$.
\end{enumerate}
\end{lem}
\begin{proof}
The existence of $\frak j$ satisfying (2) is immediate from the fact
that $\varphi_{xx'}$ preserves stratification $S_n(V_{x'})$, $S_n(V_{x})$.
Such $\frak j$ is necessarily unique.
Then (1) follows from this uniqueness.
\end{proof}
For $A \subset \{1,\dots,k\}$ we put
\begin{equation}
V(A) =
\{
(\overline y,(t_1,\dots,t_k)) \in V_x
\mid \text{If $i \in A$ then
$t_i = 0$}
\}.
\end{equation}
\begin{defn}
We define
$
\varphi_{xx'}^{\boxplus\tau} : V^{\boxplus\tau}_{x'} \to V^{\boxplus\tau}_x
$
as follows.
\begin{enumerate}
\item
If $y' \in V_{x'}$ then
$
\varphi_{xx'}^{\boxplus\tau}(y') = \varphi_{xx'}(y')
$.
\item
Let $y' =  (\overline y',(t'_1,\dots,t'_{k'}))$ and
$\mathcal{R}_{x'} (y') \in V_{x'}(A')$, where $A' \subset \{1,\dots,k'\}$.
We define
$
y_0 = \varphi_{xx'}(\mathcal R_{x'}(y'))
$
and write
$
y_0 = (\overline y_0,(t_{0,1},\dots,t_{0,k}))
$.
Then we define
$$
\varphi_{xx'}^{\boxplus \tau}(y') =(\overline y_0,(t_{1},\dots,t_{k}))
$$
where
$$
t_i
=
\begin{cases}
t'_{i'}   &\text{if  $i = \frak i(i')$, $i' \in A'$,}\\
t_{0,i}  &\text{if  $i \notin \frak i(A')$.}
\end{cases}
$$
\end{enumerate}
\end{defn}
\begin{lem}\label{lem15333}
\begin{enumerate}
\item The map
$
\varphi_{xx'}^{\boxplus\tau} : V^{\boxplus\tau}_{x'} \to V^{\boxplus\tau}_x
$
is a smooth embedding of manifolds.
\item The map
$
\varphi_{xx'}^{\boxplus\tau} : V^{\boxplus\tau}_{x'} \to V^{\boxplus\tau}_x
$ is $h_{xx'}$ equivariant.
\item We have
$
\varphi_{xx'} \circ \mathcal R_{x'}
= \mathcal R_{x} \circ \varphi^{\boxplus\tau}_{xx'}.
$
\end{enumerate}
\end{lem}
\begin{proof}
Statements (2) and (3) are obvious from the definition.
We prove (1).
For simplicity of notation we consider the case $\frak i(i) = i$.
We study the smoothness at the point
$y' = (\overline y',(t'_1,\dots,t'_{k'}))$.
We may assume without loss of generality that
$$
t'_1 = \dots = t'_{\ell} = 0, \quad
t'_{\ell+1},\dots,t'_{m} < 0, \quad
t'_{m+1},\dots,t'_{k'} > 0.
$$
\par
Let
$z' = (\overline z',(s'_1,\dots,s'_{k'})) \in V_{x'}^{\boxplus \tau}$ be a
point in a neighborhood of $y'$.
By taking the neighborhood sufficiently small, we may assume that
$s'_{\ell+1},\dots,s'_{m} < 0$ and
$s'_{m+1},\dots,s'_{k'} > 0$.
We put
$$
\mathcal R_{x'}(z')
= (\overline z',(s''_1,\dots,s''_{k'})).
$$
Then we have
$s''_{\ell+1},\dots,s''_{m} = 0$,
$s''_{m+1} = s'_{m+1},\dots,s''_{k'} = s'_{k'}$.
Moreover for $i\le \ell$,
$$
s''_i =
\begin{cases}
0  &\text{if $s'_i \le 0$} \\
s'_i  &\text{if $s'_i \ge 0$.}
\end{cases}
$$
We denote $\varphi_{xx'}^{\boxplus \tau}(z') = z$
and
$z = (\overline z,(s_1,\dots,s_{k}))$.
By definition we have
$$
\overline z = \pi_0(\varphi_{xx'}(\overline z',(s''_1,\dots,s''_{k'})))
$$
where $\pi_0 : V_x \times [0,1)^k \to V_x$ is the projection.
Therefore, the smooth dependence of $\overline z$ on $z'$ can be proved
in the same way as the proof of Lemma \ref{Vplusissmooth} (1).
(Namely we use Lemma \ref{newlem319} (2).)
\par
We also have
$$
(s_{m+1},\dots,s_{k})
=
\pi_{m+1,\dots,k}(\varphi_{xx'}(\overline z',(s''_1,\dots,s''_{k'})))
$$
where  $\pi_{m+1,\dots,k} : V_x \times [0,1)^{k} \to [0,1)^{k-m}$
is the projection to the last $k-m$ factors.
Therefore, the smooth dependence of $(s_{k'+1},\dots,s_{k})$ on $z'$ can be proved
in the same way as in the proof of Lemma \ref{Vplusissmooth}.
Moreover
$
s_i = s'_i
$
for $i= \ell+1,\dots,m$. So the smooth dependence of $s_i$ on
$z'$ is obvious.
Finally
for $i=1,\dots,\ell$ we have
$$
s_i =
\begin{cases}
\pi_{i}(\varphi_{xx'}(\overline z',(s''_1,\dots,s''_{k'})))  &\text{if $s'_i \ge 0$} \\
s'_i &\text{if $s'_i \le 0$}.
\end{cases}
$$
Here $\pi_{i} : V_x \times [0,1)^{k} \to [0,1)$
is the projection to the $i$-th factor of $[0,1)^{k}$.
Then the smooth dependence of $(s_{1},\dots,s_{m})$ on $z'$ follows
from Lemma \ref{canonicalformforadmi}.
\par
The injectivity of $\varphi^{\boxplus\tau}_{xx'}$ can be proved easily.
The smoothness of local inverse of $\varphi^{\boxplus\tau}_{xx'}$
can be proved in the same way as in the proof of
smoothness of $\varphi^{\boxplus\tau}_{xx'}$.
\end{proof}
By Lemma \ref{lem15333} (3) the embedding of bundles
$\widehat\varphi_{xx'} : \mathcal E_{x'} \to \mathcal E_{x}$
over $\varphi_{xx'}$ induces a map
$$
\mathcal R^*_{x'}\mathcal E_{x'} \to \mathcal R^*_{x}\mathcal E_{x}
$$
over $\varphi_{xx'}^{\boxplus\tau}$.
Since $\mathcal E_{x'}^{\boxplus\tau} =\mathcal R^*_{x'}\mathcal E_{x'}$ and
$\mathcal E_{x}^{\boxplus\tau} =\mathcal R^*_{x}\mathcal E_{x}$
by definition,
we obtain
$$
\widehat\varphi_{xx'}^{\boxplus\tau}
: \mathcal E_{x'}^{\boxplus\tau} \to \mathcal E^{\boxplus\tau}_{x}.
$$
In the same way as in the proof of Lemma \ref{lem15333} (1),
we can show that $\widehat\varphi_{xx'}^{\boxplus\tau}$ is a smooth embedding
of vector bundles.
We have thus proved the next lemma.
\begin{lem}\label{coordinatechage+}
Under the situation above, $(h_{xx'},\varphi^{\boxplus\tau}_{xx'},\widehat\varphi_{xx'}^{\boxplus\tau})$
is an embedding of orbifold charts.
If $(h_{xx'},\varphi_{xx'},\widehat{\varphi}_{xx'})$ is a
coordinate change,
$(h_{xx'},\varphi^{\boxplus\tau}_{xx'},\widehat\varphi_{xx'}^{\boxplus\tau})$
is also a coordinate change.
\end{lem}
We also have the following:
\begin{lem}\label{1513lem}
Let $(h_{xx''},\varphi_{xx''},\widehat\varphi_{xx''})$
(resp. $(h_{x'x''},\varphi_{x'x''},\widehat\varphi_{x'x''})$)
be as in Situation \ref{situ153},
where $x,x'$ in Situation \ref{situ153}
is replaced by $x,x''$ (resp. $x',x''$).
We define $(h_{xx''},\varphi^{\boxplus\tau}_{xx''},\widehat\varphi_{xx''}^{\boxplus\tau})$
(resp. $(h_{x'x''},\varphi^{\boxplus\tau}_{x'x''},\widehat\varphi_{x'x''}^{\boxplus\tau})$)
from $(h_{xx''},\varphi_{xx''},\widehat\varphi_{xx''})$
(resp. $(h_{x'x''},\varphi_{x'x''},\widehat\varphi_{x'x''})$)
in the same way as in the proof of Lemma \ref{coordinatechage+}.
Then we have
$$
\varphi_{xx''}^{\boxplus\tau}
= \varphi_{xx'}^{\boxplus\tau}\circ \varphi_{x'x''}^{\boxplus\tau},
\quad
\widehat\varphi_{xx''}^{\boxplus\tau}
= \widehat\varphi_{xx'}^{\boxplus\tau}\circ \widehat\varphi_{x'x''}^{\boxplus\tau}, \quad
s_x^{\boxplus\tau} \circ \varphi_{xx'}^{\boxplus\tau}
=
\widehat\varphi_{xx'}^{\boxplus\tau} \circ s_{x'}^{\boxplus\tau}.
$$
\end{lem}
\begin{defn}
In the Situation \ref{situ153}
we define:
\begin{equation}
(X \cap U)^{\boxplus\tau} = \bigcup_{x \in \psi(s^{-1}(0))}
\left((s_x^{\boxplus\tau})^{-1}(0)/\Gamma_x\right)/\sim.
\end{equation}
Here the equivalence relation $\sim$ in the right hand is defined
as follows.
Let $y_i \in V^{\boxplus\tau}_{x_i}/\Gamma_{x_i}$ with $s_{x_i}(\tilde y_i) = 0$, $[\tilde y_i] = y_i$ for $i=1,2$.
Then $y_1 \sim y_2$
if and only if there exist $\frak V_x$,
$\tilde y \in V_x$, and $(h_{x_i x},\varphi_{x_ix},\widehat\varphi_{x_ix})$
as in Situation \ref{situ153} such that $s(\tilde y) = 0$
and
$$
[\varphi_{x_ix}(\tilde y)] = y_i
$$
in $V^{\boxplus\tau}_{x_i}/\Gamma_{x_i}$ for $i=1,2$.
\footnote{Lemma \ref{1513lem} implies that this is an equivalence relation.}
\par
The maps $\mathcal R_{x} : (s_x^{\boxplus\tau})^{-1}(0) \to V_{x}$ for various $x$ induce a map
$(X \cap U)^{\boxplus\tau} \to X$, which we denote by
\begin{equation}\label{defformR}
\mathcal R : (X \cap U)^{\boxplus\tau} \to X.
\end{equation}
\end{defn}
\begin{lem}\label{Lema1515}
In Situation \ref{situ151} we have a Kuranishi chart
$$
\mathcal U^{\boxplus\tau}
= (U^{\boxplus\tau},\mathcal E^{\boxplus\tau},\psi^{\boxplus\tau},s^{\boxplus\tau})
$$
of $(X \cap U)^{\boxplus\tau}$
such that
$(V^{\boxplus\tau}_x/\Gamma_x,\mathcal E^{\boxplus\tau}_x,\psi_x^{\boxplus\tau},
s_x^{\boxplus\tau})$
becomes its orbifold chart.
\end{lem}
\begin{proof}
This is a consequence of Lemmas \ref{lem15333}, \ref{coordinatechage+}, \ref{1513lem}.
\end{proof}

\begin{lemdef}\label{lemdef157tutuki}
In the situation of Lemma \ref{Lema1515},
we call $\mathcal U^{\boxplus\tau}$ the {\rm $\tau$-collaring}, or
{\rm $\tau$-corner trivialization} of
$\mathcal U$.
\index{$\tau$-collaring ! of Kuranishi chart}
\index{corner ! $\tau$-corner trivialization}
\index{corner ! trivialization of corners, {\it see: $\tau$-corner trivialization}}
\begin{enumerate}
\item
$\partial(U^{\boxplus\tau})$ is canonically diffeomorphic to
$(\partial U)^{\boxplus\tau}$.
\item
Let
$\frak S = \{(\frak V_{\frak r},\mathcal S_{\frak r}) \mid \frak r \in \frak R\}$ be a CF-perturbation
of $\mathcal U$.
Then $\{(\frak V^{\boxplus\tau}_{\frak r},\mathcal S^{\boxplus\tau}_{\frak r}) \mid \frak r \in \frak R\}$ is
a CF-perturbation
of $\mathcal U^{\boxplus\tau}$.
We denote it by $\frak S^{\boxplus\tau}$ and call it
the  {\rm $\tau$-collaring} of $\frak S$.
\index{$\tau$-collaring ! of CF-perturbation}
\item
If $\frak S$ is equivalent to $\frak S'$, then
$\frak S^{\boxplus\tau}$ is equivalent to $\frak S^{\prime \boxplus\tau}$.
\item
If $f : U \to M$ is a smooth map and strongly
submersive on $K$ with respect to $\frak S$,
then $f \circ \mathcal R$ is strongly
submersive with respect to $\frak S^{\boxplus\tau}$.
We denote it by $f^{\boxplus\tau}$ and call it {\rm $\tau$-collaring}
\index{$\tau$-collaring ! of smooth map}
of $f$.
\item
The strong transversality to a map $N \to M$ is also preserved.
\item
The versions of (2)-(5) where
`CF-perturbation' is replaced by
`multivalued perturbation' also hold.
\index{$\tau$-collaring ! of multivalued perturbation}
\item
If $h$ is a differential form on $U$, then $\widetilde h_x$ for
various $x$ are glued to define a
differential form on $U^{\boxplus\tau}$.
We denote it by $h^{\boxplus\tau}$ and call it
the {\rm $\tau$-collaring} of $h$.
\index{$\tau$-collaring ! of differential form}
\item
In the situation of (2)(3)(4)(7),
we assume that $h$ is compactly supported. Then we have
\begin{equation}
f^{\boxplus\tau}_{!}(h^{\boxplus\tau};\frak S^{\boxplus\tau})
=
f_{!}(h;\frak S).
\end{equation}
\end{enumerate}
\end{lemdef}
This follows immediately from Lemma-Definition \ref{lemdef157}.
Also the next lemma is a straightforward generalization of
Lemma \ref{lem1599} to the case of Kuranishi chart.
\begin{lem}\label{lem1519}
We put
$$
\overset{\circ\circ}S_k(U^{\boxplus\tau})
=
S_k(U^{\boxplus\tau})
\cap \mathcal R^{-1}(\overset{\circ}S_k(U)).
$$
\begin{enumerate}
\item
The closure of $\overset{\circ\circ}S_k(U^{\boxplus\tau})$
in $\overset{\circ}S_k(U^{\boxplus\tau})$ is an orbifold with
corners.
We call $\overset{\circ\circ}S_k(U^{\boxplus\tau})$
a {\rm small corner} of codimension $k$.
\index{small corner}
\item
The retraction map
$\mathcal R$ induces an orbifold diffeomorphism
from ${\rm Close}(\overset{\circ\circ}S_k(U^{\boxplus\tau}))$  onto
$\widehat S_k(U)$.
\end{enumerate}
\end{lem}
\subsection{Trivialization of corners of Kuranishi structure}
\label{subsection:triconkura}
In this subsection and the next,
we study trivialization of corners of Kuranishi structure and good coordinate system.
For a K-space $(X,\widehat{\mathcal U})$ we firstly describe the underlying topological space $X^{\boxplus\tau}$ of the trivialization of corners
of $X$ in Definition \ref{defXenhance}.
In the next subsection, we define the trivialization of corners
$(X^{\boxplus\tau},\widehat{\mathcal U^{\boxplus\tau}})$
(or sometimes called $\tau$-collaring) of the K-space $(X,\widehat{\mathcal U})$.
We first consider the following situation.

\begin{shitu}\label{situ156}
Let $\mathcal U_i = (U_i,\mathcal E_i,s_i,\psi_i)$ be Kuranishi charts of $X$
and $\Phi_{21} = (\varphi_{21},\widehat\varphi_{21})$ an embedding
of Kuranishi charts.
We may decorate $\mathcal U_i$ by some of the following in addition:
\begin{enumerate}
\item
We are given CF-perturbations $\frak S^i$ of
$\mathcal U_i$ ($i=1,2$) such that $\frak S^1$, $\frak S^2$ are compatible
with $\Phi_{21}$.
\item
We are given differential forms $h_i$ on $U_i$ ($i=1,2$) such that $h_1 = \varphi_{21}^*h_2$.
\item
We are given smooth maps $f_i : U_i \to M$ ($i=1,2$) such that $f_1 = f_2 \circ \varphi_{21}$.
\item
We are given multivalued perturbations $\frak s^i$ of
$\mathcal U_i$ ($i=1,2$) such that $\frak s^1$, $\frak s^2$ are compatible
with $\Phi_{21}$.
\item
We have another Kuranishi chart $\mathcal U_3$ and
an embedding $\Phi_{32} : \mathcal U_2 \to \mathcal U_3$.
We put $\Phi_{31} =  \Phi_{32}\circ \Phi_{21}$.$\blacksquare$
\end{enumerate}
\end{shitu}
\begin{lem}\label{lem151888}
In Situaion \ref{situ156} we have an embedding
of Kuranishi charts $\Phi^{\boxplus\tau}_{21} : \mathcal U^{\boxplus\tau}_1
\to \mathcal U^{\boxplus\tau}_2$, whose restriction to $\mathcal U_1$
coincides with $\Phi_{21}$.
Moreover we have the following.
\begin{enumerate}
\item
In case of  Situation \ref{situ156} (1), $\frak S^{1\boxplus\tau}$, $\frak S^{2\boxplus\tau}$ are compatible
with $\Phi_{21}^{\boxplus\tau}$.
\item
In case of  Situation \ref{situ156} (2),
$h^{\boxplus\tau}_1 = (\varphi_{21}^{\boxplus\tau})^* h^{\boxplus\tau}_2$.
\item
In case of Situation \ref{situ156} (3),
$f_1^{\boxplus\tau} = f_2^{\boxplus\tau} \circ \varphi_{21}^{\boxplus\tau}$.
\item
In case of Situation \ref{situ156} (4),
$\frak s^{1\boxplus\tau}$ and $\frak s^{2\boxplus\tau}$ are compatible
with $\Phi_{21}^{\boxplus\tau}$.
\item
In case of Situation \ref{situ156} (5),
we have
$\Phi^{\boxplus\tau}_{31} =  \Phi^{\boxplus\tau}_{32}\circ \Phi^{\boxplus\tau}_{21}$.
\item
$\varphi^{\boxplus\tau}_{21} \circ \mathcal R_1 = \mathcal R_2 \circ \varphi^{\boxplus\tau}_{21}$.
\end{enumerate}
\end{lem}
\begin{rem}
Both of $\mathcal U^{\boxplus\tau}_1$  and
$\mathcal U^{\boxplus\tau}_2$ are Kuranishi charts of the
topological space $(X
\cap U)^{\boxplus\tau}_2$.
\end{rem}
\begin{proof}[Proof of Lemma \ref{lem151888}]
We use Lemma \ref{lem26999}
to obtain objects
$\{\frak V^i_{\frak r} \mid \frak r \in \frak R_i\}$ and
$(h_{\frak r,21},\varphi_{\frak r,21},\hat\varphi_{\frak r,21})$
which have the properties spelled out
there.
Let $\frak r \in \frak R_1$.
Then the map
$\varphi_{\frak r,21} : V^1_{\frak r} \to V^2_{\frak r}$
is extended to
$\varphi^{\boxplus\tau}_{\frak r,21} : V^{1 \boxplus\tau}_{\frak r} \to V^{2 \boxplus\tau}_{\frak r}$
as follows.
Let
$(\overline y',(t'_1,\dots,t'_{d({\frak r})}))
\in V^1_{\frak r}$
and
$$
(\overline y'',(t''_1,\dots,t''_{d({\frak r})}))
=
\varphi_{\frak r,21}(\mathcal R_{\frak r}(\overline y',(t'_1,\dots,t'_{d({\frak r})}))).
$$
Then we put
$$
\varphi_{\frak r,21}^{\boxplus\tau}(\overline y',(t'_1,\dots,t'_{d({\frak r})}))
=
(\overline y,(t_1,\dots,t_{d({\frak r})}))
$$
where $\overline y' = \overline y''$ and
$$
t_i
=
\begin{cases}
t'_i  &\text{if $t'_i \le 0$}, \\
t''_i  &\text{if $t'_i \ge 0$}.
\end{cases}
$$
We define $\widehat\varphi^{\boxplus\tau}_{\frak r,21}$ in a similar way.
Using Lemma \ref{lem26999} (4) and Lemma \ref{lem2622}, it is easy to see that
$(h_{\frak r,21},\varphi^{\boxplus\tau}_{\frak r,21},\widehat\varphi^{\boxplus\tau}_{\frak r,21})$
is a representative of the required embedding.
\par
It is straightforward to check (1)-(6).
\end{proof}
We will glue $\mathcal U_p$ for various $p$ using
Lemma \ref{lem151888} to obtain
a Kuranishi structure $\widehat{\mathcal U^{\boxplus\tau}}$.
(See Lemma-Definition
\ref{lemdef1522}.)
Its underlying topological space is obtained as follows.

\begin{defn}\label{defXenhance}
(1) Let $(X,\widehat{\mathcal U})$ be a K-space.
We define a topological space $X^{\boxplus\tau}$
as follows.
We take a disjoint union
$$
\coprod_{p\in X}  (s_p^{\boxplus\tau})^{-1}(0)/\Gamma_{p}
$$
and define an equivalence relation
$\sim$ as follows:
Let $x_p \in (s^{\boxplus\tau}_{p})^{-1}(0)$ and
$x_q \in (s^{\boxplus\tau}_{q})^{-1}(0)$.
We define $[x_p] \sim [x_q]$ if there exist
$r \in X$ and $x_r \in (s^{\boxplus\tau}_p)^{-1}(0)
\cap U^{\boxplus\tau}_{pr} \cap U^{\boxplus\tau}_{qr}$
such that
\begin{equation}\label{equiXXcondplus}
[x_p] = \varphi_{pr}^{\boxplus\tau}([x_r]),
\quad
[x_q] = \varphi_{qr}^{\boxplus\tau}([x_r]).
\end{equation}
As we will see in Lemma \ref{lem1521},
$\sim$ is an equivalence relation.
We define $X^{\boxplus\tau}$ as the set of the equivalence classes of this
equivalence relation $\sim$
\begin{equation}\label{def:collaredX}
X^{\boxplus\tau} :=
\left(\coprod_{p\in X}  (s_p^{\boxplus\tau})^{-1}(0)/\Gamma_{p}
\right)/ \sim.
\end{equation}
(2) For $k=1,2,\dots$ we define
\begin{equation}\label{def:kcollaredX}
S_k(X^{\boxplus\tau}) :=
\left(\coprod_{p\in X}  S_k(U_p^{\boxplus\tau}) \cap
\left((s_p^{\boxplus\tau})^{-1}(0)/\Gamma_{p}\right)
\right) /\sim.
\end{equation}
\end{defn}
The relation $\sim$ is defined on the sets
$s_p^{-1}(0)/\Gamma_p$ or the enhanced sets $(s_p^{\boxplus\tau})^{-1}(0)/\Gamma_{p}$.
So the messy process to shrink the domain to ensure the
consistency of coordinate changes are not necessary here.
In fact, we show the following.
\begin{lem}\label{lem1521}
The relation $\sim$ in Definition \ref{defXenhance} is an equivalence relation.
\end{lem}
\begin{proof}
(1) We just check transitivity.
Other properties are easier to prove.
Suppose $[x_p] \sim [x_q]$ and $[x_q] \sim [x_r]$.
By Definition \ref{defXenhance}, there exist $u, v \in X$ and $x_u \in U_u^{\boxplus \tau}, x_v \in U_v^{\boxplus \tau}$ such that
$$[x_p]=\varphi_{pu}^{\boxplus \tau}([x_u]), \ [x_q]=\varphi_{qu}^{\boxplus \tau}([x_u])$$
$$[x_q]=\varphi_{qv}^{\boxplus \tau} ([x_v]), \ [x_r]=\varphi_{rv}^{\boxplus \tau} ([x_v]).$$
By Lemma \ref{lem151888} (6),
we obtain
$\psi_p(\mathcal R_p([x_p])) = \psi_q(\mathcal R_q([x_q])) = \psi_r(\mathcal R_r(x_r))$
from \eqref{equiXXcondplus}.
Denote this common point by $t \in X$.
We may take $U_t = \overline{U_t} \times [0,1)^d$
where the $ [0,1)^d$ component of $o_t$ is zero
and $\overline U_t$ has no boundary.
We note
$U^{\boxplus\tau}_t = \overline{U_t} \times [-\tau,1)^d$.
Since $t \in s_p^{-1}(0) \cap s_q^{-1}(0) \cap s_r^{-1}(0) \cap s_u^{-1}(0) \cap s_v^{-1}(0)$, we have coordinate changes $\varphi_{pt}, \varphi_{qt}, \varphi_{rt}, \varphi_{ut}, \varphi_{vt}$
so that
$$\mathcal R_p([x_p])=\varphi_{pu} \circ \varphi_{ut}(o_t) = \varphi_{pt}(o_t), $$
$$\mathcal R_q([x_q])=\varphi_{qu} \circ \varphi_{ut}(o_t)= \varphi_{qv} \circ \varphi_{vt}(o_t) = \varphi_{qt}(o_t),$$
and
$$\mathcal R_r([x_r]) = \varphi_{rv} \circ \varphi_{vt}(o_t) = \varphi_{rt}(o_t).$$
Also note that
%
%
the restrictions of $\varphi_{pt}^{\boxplus \tau}, \varphi_{qt}^{\boxplus
\tau}, \varphi_{rt}^{\boxplus \tau}, \varphi_{ut}^{\boxplus \tau}$ and
$\varphi_{vt}^{\boxplus \tau}$ to $\mathcal R_t^{-1}(o_t)$ are bijections to
$\mathcal R_p^{-1}(\mathcal R_p([x_p]))$, $\mathcal R_q^{-1}(\mathcal
R_q([x_q]))$, $\mathcal R_r^{-1}(\mathcal R_r([x_r]))$, $\mathcal
R_u^{-1}(\mathcal R_u([x_u]))$
and $\mathcal R_v^{-1}(\mathcal R_v([x_v]))$, respectively.
Hence there exists $x_t \in U_t^{\boxplus \tau}$ such that $[x_q]=\varphi_{qt}^{\boxplus \tau} ([x_t])$,
$[x_u]=\varphi_{ut}^{\boxplus \tau}([x_t])$ and $[x_v]=\varphi_{vt}^{\boxplus \tau}([x_t])$.
Therefore we have
$$[x_p]= \varphi_{pu}^{\boxplus}([x_u]) = \varphi_{pt}^{\boxplus \tau} ([x_t]) \ \text{and} \
[x_r]=\varphi_{rv}^{\boxplus \tau}([x_v]) = \varphi_{rt}^{\boxplus \tau}([[x_t]),$$
which imply that $[x_p] \sim [x_r]$.
\end{proof}
\begin{rem}
In this section we consider the case of a Kuranishi structure on a compact
metrizable space $X$.
We may also consider the case of Kuranishi structure of a pair $Z \subseteq X$
of a metrizable space $X$ and its compact subspace $Z$.
There is actually nothing new to do so, except the following
point.
To define a topological space $X^{\boxplus\tau}$
we used a Kuranishi structure on $X$.
If we are given a Kuranishi structure of $Z \subseteq X$,
we can still define $Z^{\boxplus\tau}$.
However we can define $X^{\boxplus\tau}$ only
in a neighborhood of $Z$.
The situation here is similar to the situation
we met in defining the boundary $\partial(X,Z;\widehat{\mathcal U})$.
See \cite[Remark 8.9 (2)]{part11}.
It seems unlikely that this point becomes an important issue
in the application.
In fact, in all the cases we know so far appearing in the actual
applications, it is enough to define  $X^{\boxplus\tau}$
in a neighborhood of $Z$, or there is an obvious way
to define $X^{\boxplus\tau}$
in the particular situations.
\end{rem}
\begin{defn}
Using a good coordinate system $\widetriangle{\mathcal U}$ of $X$,
we can define $X^{\boxplus\tau}$ in a similar way as the case of Kuranishi structure.
\end{defn}

\subsection{Collared Kuranishi structure}
\label{subsec:collaredKura}

For a point $p \in X$ we can define its Kuranishi neighborhood as
${\mathcal U}^{\boxplus\tau}_p$.
There is a slight issue in defining a Kuranishi neighborhood compatible
with the collar structure.
\begin{exm}
We consider an orbifold $X=[0,\infty)^2/\Z_2$, where $\Z_2$ acts by
exchanging the factors.
We want to regard it as a `$1$-collared orbifold'.
If $p = [0,0]$ we can take an obvious choice $[0,1)^2/\Z_2$
as its `collared neighborhood'.
There is an issue in case $p = [(0,0.5)]$.
We might try to take its neighborhood such as
$[0,1) \times (0.3,0.7)$.
This however does not work.
In fact $(0.4,0.6) \sim (0.6,0.4)$ but  $\Z_2$ is not contained
in the isotropy group of $(0,0.5)$.
It seems impossible to find a good `collared neighborhood' of $[(0,0.5)]$
such that the `length' of the collar is $1$.
This is a technical problem and certainly we should regard
$[0,\infty)^2/\Z_2$
to have a collar of length $\ge 1$.
\par
It seems to the authors that
the best way to define the appropriate notion of $\tau$-collared
cornered orbifold is as follows: We do not
define an orbifold chart at the points in $[0,1)^2 \setminus \{(0,0)\}$.
The points of $[0,1)^2 \setminus \{(0,0)\}$ are contained in the chart at $(0,0)$ so we do not need
an orbifold chart at those points.
We will define the notion of {\it $\tau$-collared Kuranishi structure}
along this line below.
\end{exm}
\begin{rem}\label{rem1528}
The above mentioned trouble occurs
only when the action of the isotropy group on the
normal factor $[0,1)^k$ is nontrivial.
So it does not occur in the situation
of our applications in Sections \ref{sec:systemline1}-\ref{sec:systemtree2}.
However, we present the formulation which works in more general cases.
It actually appears when we will study
the moduli space of pseudo-holomorphic curves from
a bordered Riemann surface of arbitrary genus with
arbitrary number of boundary components.
\par
There occurs no similar issue for the definition of the $\tau$-collared good coordinate system.
\end{rem}
\begin{defn}\label{defn1529}
Given $\tau >0$
let
$$
X' = X^{\boxplus\tau}
$$ be the $\tau$-collaring of certain Kuranishi structure
$\widehat{\mathcal U}$ on $X$.
We put
\begin{equation}\label{def:marumaruSk}
\overset{\circ\circ}S_k(X',\widehat{\mathcal U})
= S_k(X')
\cap \mathcal R^{-1}
(\overset{\circ}S_k(X,\widehat{\mathcal U})),
\end{equation}
where $S_k(X')$ is defined by \eqref{def:kcollaredX},
and define a subset
$B_{\tau}(\overset{\circ\circ}S_k(X',\widehat{\mathcal U})) \subset X'$
as the union of
\begin{equation}
\psi^{\boxplus\tau}_p
\left(
(s^{\boxplus\tau}_p)^{-1}(0) \cap
\{(\overline y,(t_1,\dots,t_k))
\mid t_i \le 0, \,\, i=1,\dots,k\}
\right)
\end{equation}
for $p \in \overset{\circ}S_k(X,\widehat{\mathcal U})$.
\end{defn}
We note that
if $p' \in \overset{\circ\circ}S_k(X',\widehat{\mathcal U})$ then
$p' = \psi^{\boxplus\tau}_p(\overline y,(-\tau,\dots,-\tau))$
for $p = \mathcal R(p')$.
Therefore  $\overset{\circ\circ}S_k(X',\widehat{\mathcal U})
\subset B_{\tau}(\overset{\circ\circ}S_k(X',\widehat{\mathcal U}))$.
We also note that
\begin{equation}
\aligned
B_{\tau}(\overset{\circ\circ}S_0(X',\widehat{\mathcal U})) &=
\overset{\circ\circ}S_0(X',\widehat{\mathcal U}) = \overset{\circ}S_0(X,\widehat{\mathcal U}),\\
B_{\tau}(\overset{\circ\circ}S_k(X',\widehat{\mathcal U}))
\cap X &= \overset{\circ}S_k(X,\widehat{\mathcal U}).
\endaligned
\end{equation}
See Figure \ref{Figure16-3}.
\begin{figure}[h]
\hskip-5cm
\centering
\hskip3cm
\includegraphics[scale=0.4]{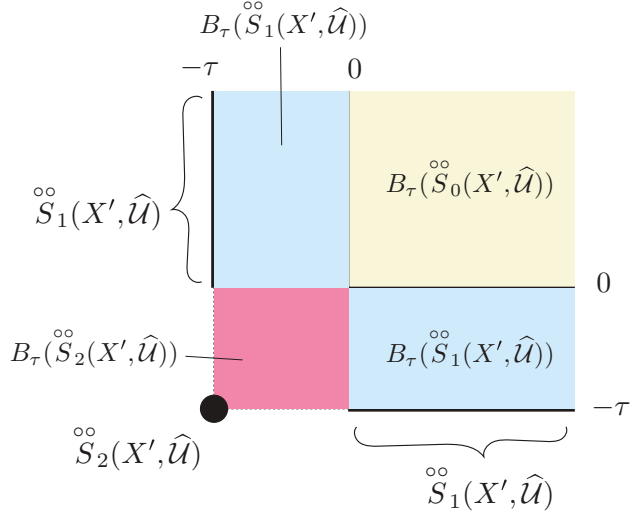}
\caption{$B_{\tau}(\overset{\circ\circ}S_k(X',\widehat{\mathcal U}))$}
\label{Figure16-3}
\end{figure}
\begin{lem}\label{lem1530}
Let $X'=X^{\boxplus\tau}$ be as above. Then it has the following decomposition:
$$
X' = \coprod_k B_{\tau}(\overset{\circ\circ}S_k(X',\widehat{\mathcal U}))
$$
where the right hand side is the disjoint union.
\end{lem}
\begin{proof}
Let $p' \in X'$ and put $p = \mathcal R(p')$.
We can write as
$p' = \psi^{\boxplus\tau}_p(\overline y,(t_1,\dots,t_k))$.
Without loss of generality, we may assume that
$t_1,\dots,t_{\ell} \le 0 < t_{\ell+1},\dots,t_k$ for some $\ell$.
We put
$q = \psi^{\boxplus\tau}_p(\overline y,(0,\dots,0,t_{\ell+1},\dots,t_k))
\in \overset{\circ}S_{\ell}(X,\widehat{\mathcal U})$.
We may choose the coordinate of $q$ so that
the map
$\frak j : \{1,\dots,\ell\} \to \{1,\dots,k\}$
appearing in the coordinate change from an orbifold chart at $q$
to an orbifold chart at $p$, which appeared in Lemma
\ref{lema1510}, is $\frak j(i) = i$.
Then we can take $\overline y'$ such that
$$
\psi^{\boxplus\tau}_p(\overline y,(0,\dots,0,t_{\ell+1},\dots,t_k))
=
\psi^{\boxplus\tau}_q(\overline y',(0,\dots,0)).
$$
Then
$p' = \psi^{\boxplus\tau}_q(\overline y',(t_{1},\dots,t_{\ell}))
\in B_{\tau}(\overset{\circ\circ}S_{\ell}(X',\widehat{\mathcal U}))$.
Moreover,
$$
B_{\tau}(\overset{\circ\circ}S_k(X',\widehat{\mathcal U}))
\cap B_{\tau}(\overset{\circ\circ}S_{\ell}(X',\widehat{\mathcal U}))
=
\emptyset
$$
for $k\ne\ell$ is obvious from definition.
\end{proof}
\begin{defn}\label{defn1531}
Suppose we are in the situation of Definition \ref{defn1529}.
In particular, $X'$ is a compact metrizable space
homeomorphic to $X^{\boxplus\tau}$ for a certain
$K$-space $(X,\widehat{\mathcal U})$.
\begin{enumerate}
\item
Let $p' \in \overset{\circ\circ}S_k(X',\widehat{\mathcal U})$.
A {\it $\tau$-collared Kuranishi neighborhood} at $p'$
\index{$\tau$-collared ! Kuranishi neighborhood}
is a Kuranishi chart $\mathcal U_{p'}$ of $X'$ such that
$\mathcal U_{p'} = (\mathcal U_{p})^{\boxplus\tau}$
for a certain Kuranishi neighborhood $\mathcal U_{p}$ of
$p =\mathcal R(p')$. (See Figure \ref{Figure16-4}.)
\item
For $p' \in \overset{\circ\circ}S_k(X',\widehat{\mathcal U})$ and
$q' \in \overset{\circ\circ}S_{\ell}(X')$,
let $\mathcal U_{p'} = \widehat{\mathcal U^{\boxplus\tau}_p}
= (\mathcal U_{p})^{\boxplus\tau}$ and $\mathcal U_{q'}
= \widehat{\mathcal U^{\boxplus\tau}_q}
$
be their $\tau$-collared Kuranishi neighborhoods, respectively.
Suppose
$q' \in \psi_{p'}(s_{p'}^{-1}(0))$.
A {\it $\tau$-collared coordinate change}
\index{$\tau$-collared ! coordinate change}
$\Phi_{p'q'}$ from  $\mathcal U_{q'}$
to $\mathcal U_{p'}$ is
$\Phi_{pq}^{\boxplus\tau}$ defined by Lemma \ref{lem151888},
where $\Phi_{pq}$ is a coordinate change
from $\mathcal U_{q}$
to $\mathcal U_{p}$.
\item
A {\it $\tau$-collared Kuranishi structure}
\index{$\tau$-collared ! Kuranishi structure}
$\widehat{\mathcal U'}$ on $X'$ consists of the
following objects:
\begin{enumerate}
\item
For each $p' \in \overset{\circ\circ}S_k(X',\widehat{\mathcal U})$,
$\widehat{\mathcal U'}$ assignes a
$\tau$-collared Kuranishi neighborhood $\mathcal U_{p'}$.
\item
For each
$p' \in \overset{\circ\circ}S_k(X',\widehat{\mathcal U})$ and
$q' \in \overset{\circ\circ}S_{\ell}(X',\widehat{\mathcal U})$
with $q' \in \psi_{p'}(s_{p'}^{-1}(0))$,
$\widehat{\mathcal U'}$ assignes a $\tau$-collared
coordinate change $\Phi_{p'q'}$.
\item
If $p' \in \overset{\circ\circ}S_k(X',\widehat{\mathcal U})$,
$q' \in \overset{\circ\circ}S_{\ell}(X',\widehat{\mathcal U})$,
$r' \in \overset{\circ\circ}S_{m}(X',\widehat{\mathcal U})$
with
$q' \in \psi_{p'}(s_{p'}^{-1}(0))$ and
$r' \in \psi_{q'}(s_{q'}^{-1}(0))$,
then we require
$$\Phi_{p'q'} \circ \Phi_{q'r'}\vert_{U_{p'q'r'}}
=
\Phi_{p'r'}\vert_{U_{p'q'r'}}
$$
where $U_{p'q'r'} = U_{p'r'} \cap \varphi_{q'r'}^{-1}(U_{p'q'})$.
\end{enumerate}
\item
A {\it $\tau$-collared K-space}
is a pair of a compact metrizable space
and its $\tau$-collared Kuranishi structure.
\index{$\tau$-collared ! K-space}
\index{K-space ! $\tau$-collared}
It is obtained from a K-space $(X,\widehat{\mathcal U})$
as in Lemma-Definition \ref{lemdef1522} below.
\item
We can define the notion
of $\tau$-collared CF-perturbation,
\index{$\tau$-collared ! CF-perturbation}
$\tau$-collared multivalued perturbation,
\index{$\tau$-collared ! multivalued perturbation}
$\tau$-collared good coordinate system,
\index{$\tau$-collared ! good coordinate system}
$\tau$-collared Kuranishi chart,
\index{$\tau$-collared ! Kuranishi chart}
$\tau$-collared vector bundle,
\index{$\tau$-collared ! vector bundle}
$\tau$-collared smooth section,
$\tau$-collared embedding of various kinds,
\index{$\tau$-collared ! embedding}
etc.
in the same way.
Actually those objects on $(X^{\boxplus\tau}
,\widehat{\mathcal U^{\boxplus\tau}})$ are obtained from the corresponding objects on  $(X,\widehat{\mathcal U})$
by applying the process of $\tau$-collaring on each chart
as in Lemma \ref{lemma1523} below.
\item
An $\mathscr A$-parametrized family of
$\tau$-collared CF-perturbations is
said to be {\it uniform}
\index{uniform family ! of $\tau$-collared CF-perturbations}
if it is of the form
$\{ \widehat{\mathfrak S^{\boxplus\tau}_{\sigma}} \mid \sigma \in \mathscr A\}$
for a certain uniform family
$\{ \widehat{\mathfrak S_{\sigma}} \mid \sigma \in \mathscr A\}$
of CF-perturbations on $(X,\widehat{\mathcal U})$.
The definition of
uniform family of $\tau$-collared multivalued perturbations
is the same.
\index{uniform family ! of $\tau$-collared multivalued perturbations}
\end{enumerate}
\end{defn}
\begin{figure}[h]
\centering
\includegraphics[scale=0.4]{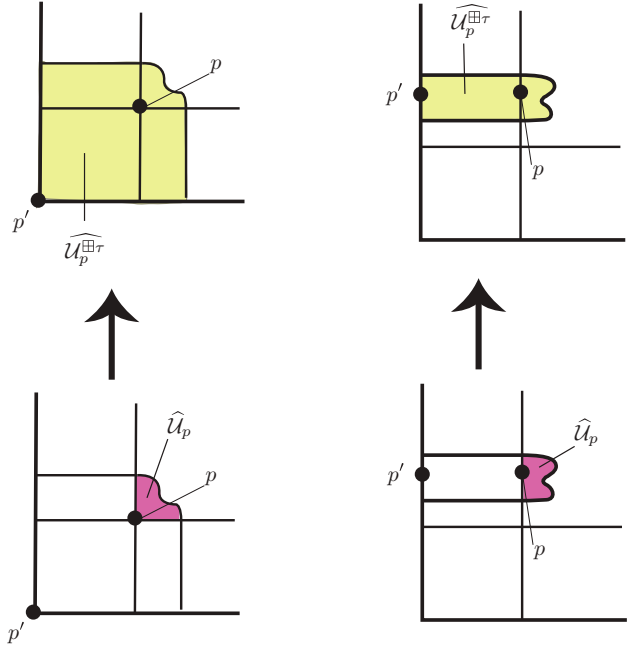}
\caption{$\widehat{\mathcal U^{\boxplus\tau}_p}$}
\label{Figure16-4}
\end{figure}

\begin{lemdef}\label{lemdef1522}
For any K-space $(X,\widehat{\mathcal U})$
we can assign a $\tau$-collared K-space $(X^{\boxplus\tau},\widehat{\mathcal U^{\boxplus\tau}})$
such that:
\begin{enumerate}
\item
Its underlying topological space $X^{\boxplus\tau}$ is as in
Definition \ref{defXenhance}.
\item
If $p \in \overset{\circ\circ}S_k(X^{\boxplus\tau})$,
its Kuranishi neighborhood is
${\mathcal U}^{\boxplus\tau}_{\mathcal R(p)}$, where
${\mathcal U}^{\boxplus\tau}_{\mathcal R(p)}$ is defined in Lemma \ref{Lema1515}.
\item
The coordinate changes are $\Phi^{\boxplus\tau}_{\mathcal R(p)\mathcal R(q)}$,
where $\Phi^{\boxplus\tau}_{\mathcal R(p)\mathcal R(q)}$ is defined in Lemma \ref{lem151888}.
\end{enumerate}
We call $(X^{\boxplus\tau},\widehat{\mathcal U^{\boxplus\tau}})$
the {\rm $\tau$-collaring}
(or {\rm trivialization of corners})
of $(X,\widehat{\mathcal U})$.
\index{$\tau$-collaring ! $\tau$-collaring}
\index{trivialization of corners ! of K-space $(X,\widehat{\mathcal U})$}
We sometimes write
$(X,\widehat{\mathcal U})^{\boxplus\tau}$
in place of
$(X^{\boxplus\tau},\widehat{\mathcal U^{\boxplus\tau}})$.
\end{lemdef}
\begin{proof}
This is immediate from the definition.
\end{proof}
For a $\tau$-collared Kuranishi structure
$\widehat{\mathcal U^{\boxplus\tau}}$ on
$X' = X^{\boxplus\tau}$
we sometimes write
$\overset{\circ\circ}S_k((X,\widehat{\mathcal U})^{\boxplus\tau})$
etc. in place of
$\overset{\circ\circ}S_k(X^{\boxplus\tau},\widehat{\mathcal U^{\boxplus\tau}})$
etc..
\begin{defn}
Let
$\widetriangle{\mathcal U}
= (({\frak P},\le), \{\mathcal U_{\frak p}\mid \frak p \in \frak P\},
\{\Phi_{\frak p\frak q}
\mid \frak q \le \frak p\})$
be a good coordinate system of $X$. We define a
good coordinate system $\widetriangle{\mathcal U^{\boxplus\tau}}$ of $X^{\boxplus\tau}$
as
$$
\widetriangle{\mathcal U^{\boxplus\tau}}
= (({\frak P},\le), \{\mathcal U^{\boxplus\tau}_{\frak p}\mid \frak p \in \frak P\},
\{\Phi^{\boxplus\tau}_{\frak p\frak q}
\mid \frak q \le \frak p\}).
$$
We call
$(X^{\boxplus\tau},\widetriangle{\mathcal U^{\boxplus\tau}})$
the {\it $\tau$-collaring}
(or {\it trivialization of corners}) of $(X,\widetriangle{\mathcal U})$.
\index{good coordinate system ! $\tau$-collaring}
\index{good coordinate system ! trivialization of corners}
\index{$\tau$-collaring ! of good coordinate system}
\index{trivialization of corners ! of good coordinate system}
\par
We call a good coordinate syetem to be {\it $\tau$-collared}
\index{$\tau$-collared ! good coordinate system}
\index{good coordinate system ! $\tau$-collared}
if it is isomorphic to $\widetriangle{\mathcal U^{\boxplus\tau}}$
for some $\widetriangle{\mathcal U}$.
\end{defn}
The next lemma says that many objects
defined on $(X,\widehat{\mathcal U})$
have corresponding collared objects
on $(X^{\boxplus\tau},\widehat{\mathcal U^{\boxplus\tau}})$.
This is actually a trivial statement to prove.

\begin{lem}\label{lemma1523}
We consider the situation of Lemma-Definition \ref{lemdef1522}.
\begin{enumerate}
\item
$(\partial(X,\widehat{\mathcal U}))^{\boxplus\tau}$
is canonically isomorphic to
$\partial(X^{\boxplus\tau},\widehat{\mathcal U^{\boxplus\tau}})$.
\item
A CF-perturbation $\widehat{\frak S}$
on $(X,\widehat{\mathcal U})$ induces a $\tau$-collared
CF-perturbation $\widehat{\frak S^{\boxplus\tau}}$
on $(X^{\boxplus\tau},\widehat{\mathcal U^{\boxplus\tau}})$.
\item
A strongly continuous
map $\widehat f$ from $(X,\widehat{\mathcal U})$
induces a  $\tau$-collared strongly continuous
map $\widehat{f^{\boxplus\tau}}$ from $(X^{\boxplus\tau},\widehat{\mathcal U^{\boxplus\tau}})$.
Strong smoothness and weak submersivity are preserved.
\item
In the situation of (2)(3),
if $\widehat f$ is strongly submersive with respect to
$\widehat{\frak S}$,
then  $\widehat{f^{\boxplus\tau}}$ is a $\tau$-collared strongly submersive
map with respect to
$\widehat{\frak S^{\boxplus\tau}}$.
\item
Transversality to a map $N\to M$ is also preserved .
\item
The versions of (2)(4) where
`CF-perturbation' is replaced by
`multivalued perturbation' also hold.
\item
A differential form $\widehat h$ on $(X,\widehat{\mathcal U})$
induces a $\tau$-collared differential form $\widehat {h^{\boxplus\tau}}$
on $(X^{\boxplus\tau},\widehat{\mathcal U^{\boxplus\tau}})$.
\item
In the situation of (2)(4)(7),
if $\widehat f : (X,\widehat{\mathcal U}) \to M$
is strongly submersive with respect to $\widehat{\frak S}$,
then we have
\begin{equation}\label{form1515}
\widehat f!(\widehat h;\widehat{\frak S})
=
\widehat{f^{\boxplus\tau}}!(\widehat{h^{\boxplus\tau}};\widehat{\frak S^{\boxplus\tau}}).
\end{equation}
\item
We put
$
\overset{\circ\circ}S_k(X^{\boxplus\tau},\widehat{\mathcal U^{\boxplus\tau}})
=
S_k(X^{\boxplus\tau},\widehat{\mathcal U^{\boxplus\tau}})
\cap \mathcal R^{-1}
(\overset{\circ}S_k(X,\widehat{\mathcal U}))
$ and call it
the {\rm small codimension $k$ corner}.
\index{small corner}
(We note that $\overset{\circ\circ}S_k(X^{\boxplus\tau},\widehat{\mathcal U^{\boxplus\tau}})
=
\overset{\circ\circ}S_k(X^{\boxplus\tau},\widehat{\mathcal U})$
in \eqref{def:marumaruSk}.)
Then the Kuranishi structure of $S_k(X^{\boxplus\tau},\widehat{\mathcal U^{\boxplus\tau}})$
induces a Kuranishi structure on
${\rm Clos}(\overset{\circ\circ}S_k(X^{\boxplus\tau},\widehat{\mathcal U^{\boxplus\tau}}))$.
\item
The restriction of the retraction map $\mathcal R$ is an
underlying homeomorphism of an isomorphism between the K-spaces
${\rm Clos}(\overset{\circ\circ}S_k(X^{\boxplus\tau},\widehat{\mathcal U^{\boxplus\tau}}))$
and
$\widehat S_k(X,\widehat{\mathcal U})$.
\item
If there exists an embedding $\widehat{\mathcal U} \to \widehat{\mathcal U^+}$
of Kuranishi structures,
then the space $X^{\boxplus\tau}$ defined by $\widehat{\mathcal U}$
is canonically homeomorphic to the one defined by
$\widehat{\mathcal U^+}$.
The same holds for various types of embeddings between Kuranishi structures
$\widehat{\mathcal U}$
and/or good coordinate systems ${\widetriangle{\mathcal U}}$.
\item
Various types of embeddings between Kuranishi structures
$\widehat{\mathcal U}$
and/or good coordinate systems ${\widetriangle{\mathcal U}}$
induce $\tau$-collared embeddings between
$\widehat{\mathcal U^{\boxplus\tau}}$
and/or
${\widetriangle{\mathcal U^{\boxplus\tau}}}$.
Compatibility among various objects on them (such as
CF-perturbation) is preserved under the operation $\boxplus$.
\item
If $\widehat{\mathcal U^+}$ is a thickening of
$\widehat{\mathcal U}$, then
$\widehat{\mathcal U^{+\boxplus\tau}}$ is a thickening of
$\widehat{\mathcal U^{\boxplus\tau}}$.
\item
If $(\widetilde{X},\widetilde{\widehat{\mathcal U}})$
is a $k$-fold covering of $(X,{\widehat{\mathcal U}})$,
then
$(\widetilde{X}^{\boxplus\tau},\widetilde{\widehat{\mathcal U^{\boxplus\tau}}})$
is a $k$-fold covering of
$(X^{\boxplus\tau},{\widehat{\mathcal U^{\boxplus\tau}}})$.
\item
The same results as (1)-(14) hold when  we replace Kuranishi structure by
good coordinate system.
\end{enumerate}
\end{lem}
\begin{proof}
(1)-(8) are consequence of Lemma-Definition \ref{lemdef157tutuki} (1)-(8),
respectively.
(9), (10) are consequences of Lemma \ref{lem1519} (1), (2),
respectively.
The proof of (11) is similar to the proof of Lemma \ref{lem1521}.
(12) follows from Lemma \ref{lem151888}.
We can prove (13) by putting
$O^{\boxplus\tau}_p = \mathcal R_{\mathcal R(p)}^{-1}(O_{\mathcal R(p)})$
and
$W_{p}(q)^{\boxplus\tau} = \mathcal R_{\mathcal R(p)}^{-1}(W_{\mathcal R(p)}(\mathcal R(q))$, where the notations
for $O_{\mathcal R(p)}, W_{\mathcal R(p)}$
are as in \cite[Definition 5.3 (2)]{part11}.
(14) is obvious from the definition.
The proof of (15) is the same as the proof of (1)-(14).
\end{proof}

\begin{lem}\label{lem153535}
If $(X',\widehat{\mathcal U'})$ is
$\tau$-collared, then for any $0 < \tau' < \tau$,
$X'$ has a
$\tau'$-collared Kuranishi structure which is determined
in a canonical way from the $\tau$-collared Kuranishi structure $(X',\widehat{\mathcal U'})$.
The same holds for CF-perturbation,
multivalued perturbation, good coordinate system and
various other objects.
\end{lem}
\begin{proof}
The lemma follows from, roughly speaking,
\begin{equation}\label{1515form}
((X,\widehat{\mathcal U})^{\boxplus\tau_1})^{\boxplus\tau_2}
\cong
(X,\widehat{\mathcal U})^{\boxplus\tau}
\end{equation}
where $\tau = \tau_1 + \tau_2$.
(Here $\tau_2$ corresponds to $\tau'$ in Lemma \ref{lem153535}.)
To be precise,
we will define a $\tau_2$-collared Kuranishi structure of
$(X,\widehat{\mathcal U})^{\boxplus\tau}$ below.
In other words, we will define a
Kuranishi structure on $X^{\boxplus\tau_1}$.
We first note that the equality
\begin{equation}
((\mathcal U_p)^{\boxplus\tau_1})^{\boxplus\tau_2}
=
\mathcal U_p^{\boxplus(\tau_1+\tau_2)}
\end{equation}
literally holds for a Kuranishi {\it chart} $\mathcal U_p$.
\par
We put $X' = X^{\boxplus\tau}$, $X'' =  X^{\boxplus\tau_1}$.
Then we have
$X' = X''^{\boxplus\tau_2}$.\footnote{Strictly speaking,
we defined the space
$X''^{\boxplus\tau_2}$ only when $X''$ has a
Kuranishi structure. We have defined a
$\tau_{1}$-collared Kuranishi structure on $X''$ but not
a Kuranishi structure on it yet.
But it is straightforward
to define a space $X''^{\boxplus\tau_2}$ when
a collared Kuranishi structure on $X''$ is given.
On the other hand,
we will also define a Kuranishi structure on $X''$.
So $X''^{\boxplus\tau_2}$ makes sense in either way.
(Indeed, the two ways to define $X''^{\boxplus\tau_2}$ coincide.)}
We note that the set $\overset{\circ\circ}S_k(X',\widehat{\mathcal U^{\boxplus\tau}})$
depends on whether we regard $X' = X^{\boxplus\tau}$
or $X' = X''^{\boxplus\tau_2}$.
So we write $\overset{\circ\circ}S_k(X';\tau)$
when we regard $X' = X^{\boxplus\tau}$,
and $\overset{\circ\circ}S_k(X';\tau_2)$
when we regard $X' = X''^{\boxplus\tau_2}$.
\par
We will define a $\tau_2$-collared Kuranishi chart
$\mathcal U'_{p}$ for each $p \in \overset{\circ\circ}S_k(X';\tau_2)$.
Let $p \in \overset{\circ\circ}S_k(X';\tau_2)$.
By Lemma \ref{lem1530} there exists a unique $n$ such that
$p \in B_{\tau}(\overset{\circ\circ}S_n(X';\tau))$.
Therefore there exists $p' \in \overset{\circ}S_n(X,\widehat{\mathcal U})$
such that
$$
p = \psi_{p'}(\overline y,(t_1,\dots,t_n))
$$
and $s_{p'}(\overline y,(t_1,\dots,t_n)) = 0$.
Such $p'$ is unique, if we require
$$
p' = \psi_{p'}(\overline y,(0,\dots,0))
$$
in addition.
In fact, $p' = \mathcal R(p)$. We will take this choice.
\par
We may take our coordinates $\overline V_{p'} \times
[0,c)^n$
of $V_{p'}$ so that $\Gamma_{p'}$ acts
by permutation of the $[0,c)^n$ factors.
\footnote{Finding such a choice so that it is compatible with various
coordinate changes is nontrivial.
However it is easy to make such a choice at
each point.
See Remark \ref{rem1551} for the ambiguity caused by the choice of such coordinates.}
By changing the enumeration of the $t_i$'s, we may assume,
without loss of generality, that
$$
t_1 = \dots = t_{k} = -\tau < -\tau_1 <
t_{k+1} \le \dots \le t_n.
$$
(It suffices to consider the case $k\ge 1$.)
In fact, we have $t_i \notin (-\tau,-\tau_1]$ by the assumption
$p \in \overset{\circ\circ}S_k(X';\tau_2)$.
\par
We take $k = a_0 < a_1 < \dots < a_m \le n$
such that
$
\{t_1,\dots,t_n\} = \{-\tau, t_{a_1},\dots,t_{a_m}\}
$
and $i\ne j \Rightarrow t_{a_i} \ne t_{a_j}$.
Then for given $t'_1,\dots,t'_n$, we define
$s_0,s'_1,s_1,s'_2,\dots,s'_m,s_m$ by
$$
s_0 = \max \{t'_1,\dots,t'_{k},-\tau_1\},
$$
and
$$
s'_j = \min \{t'_i \mid t_i = t_{a_j}\}, \qquad
s_j = \max \{t'_i \mid t_i = t_{a_j}\}
$$
for $j=1,2,\dots,m$.
Now we define
$$
V_p(\tau_2)
=
\{
(\overline y,(t'_1,\dots,t'_n)) \in V_{p'}^{\boxplus\tau}
\mid  s_0 < s'_1 \le  s_1 < s'_2 \le  \dots < s_{m-1} < s'_m\}
$$
and put
$\Gamma_p = \{\gamma \in \Gamma_{p'} \mid \gamma p = p\}$.
\begin{figure}[h]
\centering
\includegraphics[scale=0.4,angle=90]{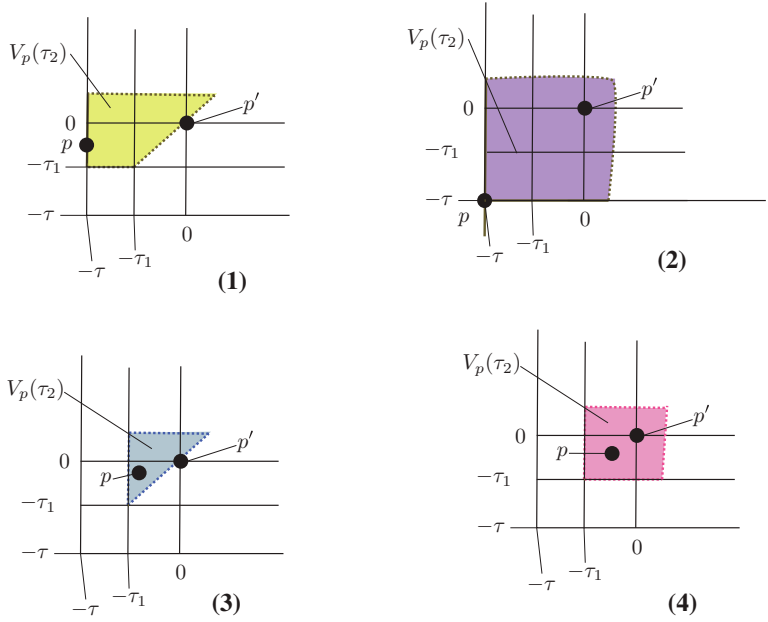}
\caption{$V_p(\tau_2)$}
\label{Figure16-5}
\end{figure}
\begin{sublem}
\begin{enumerate}
\item
If $\gamma \in \Gamma_p$, then $\gamma V_p(\tau_2) = V_p(\tau_2)$.
\item
If $\gamma \in \Gamma_{p'}$
and $\gamma V_p(\tau_2) \cap V_p(\tau_2) \ne \emptyset$,
then $\gamma \in \Gamma_p$.
\end{enumerate}
\end{sublem}
\begin{proof}
Let $A_j = \{i \in \{1,\dots,n\} \mid t_i = t_{a_j}\}$.
Then $\Gamma_{\hat p}$ induces a permutation of
$\{1,\dots,n\} $ by Definition \ref{defn297} (1)(b).
It is easy to see that
$$
\Gamma_p
=
\{\gamma \in \Gamma_{\hat p}
\mid \gamma A_j = A_j ~ \text{for all $j$} \}.
$$
The sublemma follows from this fact and the definition.
\end{proof}
We restrict $(\mathcal U_{p'})^{\boxplus\tau}$
to $V_p(\tau_2)/\Gamma_p$ to obtain
a Kuranishi chart $\mathcal  U'_p$.
We denote $U_{p}(\tau_2) = V_p(\tau_2)/\Gamma_p$
\par
We observe the following fact.
\begin{enumerate}
\item[(*)]
If $(\overline y,(t'_1,\dots,t'_n))
\in V_p(\tau_2)$ and
$t''_1,\dots,t''_{k} \in [-\tau,-\tau_1)$,
then
$$
(\overline y,(t''_1,\dots,t''_{k},t'_{\ell+1},\dots,t'_n))
\in V_p(\tau_2).
$$
\end{enumerate}
Using (*), we can prove
$$
(V_p(\tau_2) \cap V_p^{\boxplus\tau_1})^{\boxplus\tau_2}
= V_p(\tau_2).
$$
Therefore $\mathcal  U'_p$ is $\tau_2$-collared.
It is easy to construct coordinate change to obtain
a $\tau_2$-collared Kuranishi structure.
\par
The second half of the lemma follows easily from the construction of
the Kuranishi chart $\mathcal  U'_p$.
\end{proof}
\begin{rem}\label{rem1537}
In the situation of Lemma \ref{lem153535} we define a Kuranishi structure
$\widehat{\mathcal U''}$
on $X^{\boxplus\tau_1}$ such that
$(X^{\boxplus\tau_1},\widehat{\mathcal U''})^{\boxplus\tau_2}
= (X',\widehat{\mathcal U'})$ as follows.
We first consider the case when $p \in {\rm Int}\,X^{\boxplus\tau_1}$.
We put
\begin{equation}\label{formula1516}
{\mathcal U''_p}= {\mathcal U^{\boxplus\tau}_{\mathcal R(p)}}\vert_{U^{\boxplus\tau_1}_{p'} \cap U_p(\tau_2)}.
\end{equation}
Here $\mathcal R : X' \to X$ is the retraction map as in (\ref{defformR}).
\par
Let us elaborate the
right hand side of (\ref{formula1516}).
If $p \in {\rm Int} X$, then
$\mathcal U'_p = \mathcal U_p$ and
$U^{\boxplus\tau_1}_{\mathcal R(p)} = U_p$.
Therefore $\mathcal U''_p = \mathcal U'_p =  \mathcal U_p$.
\par
Suppose $p \notin  {\rm Int} X$.
(This case corresponds to (3) and (4) in Figure \ref{Figure16-5}.)
Then
$\mathcal R(p) = p' \in \overset{\circ} S_k(X,\widehat{\mathcal U})$ for $k \ge 1$.
Using the parametrization map $\psi_{p'}^{\boxplus\tau}$
of the Kuranishi chart $\mathcal U_{p'}^{\boxplus\tau}$
we can write
$$
p = \psi_{p'}^{\boxplus\tau}(\overline y,(t_1,\dots,t_k)).
$$
Since $p \in {\rm Int}\, X'$, we have $t_i > -\tau_1$. Therefore, by definition,
$$
U^{\boxplus\tau_1}_{p'} \cap U_p(\tau_2)
=
\{(\overline y',(t'_1,\dots,t'_k)) \in
U^{\boxplus\tau_1}_{p'}
\mid t'_i  > -\tau_1\}.
$$
Then
$\mathcal U''_p$ is the restriction of $\mathcal U_{p'}^{\boxplus\tau}$
to this set.
\par
We note that if $p \notin {\rm Int}\,X$ and $p \in {\rm Int}\, X^{\boxplus\tau_1}$,
then $p \notin \overset{\circ\circ} S_k(X,\tau_1)$ for any $k$.
\par
We now consider the case when
$p \in \overset{\circ}S_k(X^{\boxplus\tau_1},\widehat{\mathcal U})$, $k\ge 1$.
Consider the map $X^{\boxplus\tau} \to X^{\boxplus\tau_1}$
defined as in \eqref{defformR}.
This is the retraction map
when we regard
$X^{\boxplus\tau} = (X^{\boxplus\tau_1})^{\boxplus\tau_2}$.
We denote it by $\mathcal R'$,
which is different from the retraction map
$\mathcal R : X^{\boxplus\tau}
\to X$.
Then there exists a {\it unique} $\widehat p
\in \overset{\circ} S_k(X^{\boxplus\tau})$ such that
$\mathcal R'(\widehat p) = p$.
(See Lemma \ref{lem538} below.)
We put
$$
\mathcal U''_p = \mathcal U'_{\widehat p}\vert_{U_{p'}^{\boxplus\tau_1} \cap U_{\hat p}(\tau_2)}.
$$
See Figure \ref{Figure16-6} and compare it with Figure \ref{Figure16-5} (1)(2).
\begin{figure}[h]
\centering
\includegraphics[scale=0.4,angle=90]{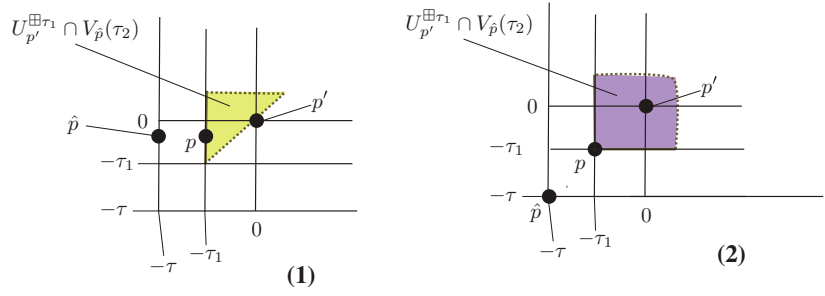}
\caption{$\mathcal U''_p$}
\label{Figure16-6}
\end{figure}
\end{rem}
We used the next lemma in the above remark.
\begin{lem}\label{lem538}
Let $p \in X^{\boxplus\tau}$ and $\mathcal R(p) \in \overset{\circ} S_k(X,\widehat{\mathcal U})$.
Then there exists a unique
$\widehat p
\in \overset{\circ} S_k(X^{\boxplus\tau},\widehat{\mathcal U^{\boxplus\tau}})$ such that
$\mathcal R'(\widehat p) = p$, $\mathcal R(\widehat p) = \mathcal R(p)$.
\end{lem}
\begin{proof}
We use the parametrization map $ \psi_{\mathcal R(p)}^{\boxplus\tau}$ to write
$
p = \psi_{\mathcal R(p)}^{\boxplus\tau}(\overline y,(t_1,\dots,t_k))
$, $t_i \le 0$.
Then
$
\widehat p = \psi_{\mathcal R(p)}^{\boxplus\tau}(\overline y,(-\tau,\dots,-\tau)).
$
\end{proof}

In the construction of this or the next section,
we need to replace a $\tau$-collared structure
by a $\tau'$-collared structure
with $\tau' < \tau$ several times.

\subsection{Extension of collared Kuranishi structure}
\label{subsec:extenonechart}
The tiresome routine works to repeat the definitions in earlier
sections are mostly over.
What we gain from these routine works is Propositions \ref{prop528}
and \ref{prop529} below which are extension theorems of a $\tau$-collared
Kuranishi structure and a $\tau$-collared CF-perturbation, respectively.
In this subsection we prove Proposition \ref{prop528} and
in the next subsection we prove Proposition \ref{prop529}.
\begin{rem}
Let $
S_k(X;\widehat{\mathcal U})
$
be a codimension $k$ stratum of a K-space $(X, \widehat{\mathcal U})$
and
$
\widehat S_k(X;\widehat{\mathcal U})
$
a normalized codimension $k$ corner of $(X, \widehat{\mathcal U})$.
(See \cite[Definition 4.15]{part11} for
$
S_k(X;\widehat{\mathcal U})
$ and
Definition \ref{norcor}
for
$
\widehat S_k(X;\widehat{\mathcal U})
$, respectively.)
If $\widehat{\mathcal U} \to \widehat{\mathcal U^+}$
is a KK-embedding (embedding of Kuranishi structures) of $X$,
the underlying topological space of
$
S_k(X;\widehat{\mathcal U})
$
 (resp. $
\widehat S_k(X;\widehat{\mathcal U})
$)
is canonically homeomorphic to the
underlying topological space of
$
S_k(X;{\widehat{\mathcal U^+}})
$
(resp.
$
\widehat S_k(X;{\widehat{\mathcal U^+}})
$.)
\par
Hereafter we write $S_k(X)$ or $\widehat S_k(X)$ in place of
$
S_k(X;\widehat{\mathcal U})
$, $
\widehat S_k(X;\widehat{\mathcal U})
$.
They stand for the underlying topological spaces.

\subsubsection{Statement}\label{subsubsec:statementextKuranishi}
To state the extension theorem (Proposition \ref{prop528}) of a $\tau$-collared Kuranishi structure
we consider the following situation.
\end{rem}
\begin{shitu}\label{sit1526}
Let $X$ be a paracompact Hausdorff space
with $\tau$-collared Kuranishi structure $\widehat{\mathcal U}$.
Let $\partial X$ be the normalized boundary of $X$.
\par
For given $\tau >0$, we are given a $\tau$-collared Kuranishi structure $\widehat {\mathcal U^+_{\partial}}$
of $\partial X$ such that
\begin{equation}\label{embeddingbetween}
\partial\widehat{\mathcal U}
<
\widehat {\mathcal U^+_{\partial}}.
\end{equation}
We assume that $\widehat {\mathcal U^+_{\partial}}$ satisfies the
following conditions:
\begin{enumerate}
\item
For each $k \ge 1$
there exists a $\tau$-collared Kuranishi structure
$\widehat{\mathcal U^+_{S_k}}$
on
$\widehat S_{k}(X)$ such that
$\widehat{\mathcal U_{S_1}^+} = \widehat {\mathcal U^+_{\partial}}$.
\item
The $\tau$-collared K-space
$\widehat S_{k}(\widehat S_{\ell}(X),\widehat{\mathcal U_{S_\ell}^+})$
is isomorphic to the $(k+\ell)!/k!\ell!$ fold covering space of
$(\widehat S_{k+\ell}(X),\widehat{\mathcal U_{S_{k+\ell}}^+})$.
\item
The following diagram of K-spaces commutes.
\begin{equation}\label{diagin26277XX}
\begin{CD}
\widehat S_{k_1}(\widehat S_{k_2}(\widehat S_{k_3}(X),\widehat{\mathcal U_{S_{k_3}}^+}))) @ >{\pi_{k_1,k_2}}>>
\widehat S_{k_1+k_2}(\widehat S_{k_3}(X),\widehat{\mathcal U^+_{S_{k_3}}})) \\
@ V{}VV @ VV{}V\\
\widehat S_{k_1}(\widehat S_{k_2+k_3}(X),\widehat{\mathcal U_{S_{k_2+k_3}}^+})) @ > {} >> (\widehat S_{k_1+k_2+k_3}(X),\widehat{\mathcal U_{S_{k_1+k_2+k_3}}^+})
\end{CD}
\end{equation}
Here $\pi_{k_1,k_2}$ is the covering map in Proposition \ref{prop2813}.
The right vertical and lower horizontal arrows are covering maps in (2).
The left vertical arrow is induced by the covering map
$\widehat S_{k_2}(\widehat S_{k_3}(X),\widehat{\mathcal U_{S_{k_3}}^+}))
\to (\widehat S_{k_2+k_3}(X),\widehat{\mathcal U_{S_{k_2+k_3}}^+})$ in (2).
\item There exists a $\tau$-collared embedding
${\widehat S_{k}(X,\widehat{\mathcal U})} \to \widehat{\mathcal U_{S_k}^+}$.
\item
The following diagram of K-spaces commutes.
\begin{equation}\label{diag15main}
\begin{CD}
\widehat S_{k}(\widehat S_{\ell}(X,\widehat{\mathcal U}))
@ > {} >>\widehat S_{k}(\widehat S_{\ell}(X),\widehat{\mathcal U^+_{S_\ell}})
\\
@ V{}VV @ VV{}V \\
\widehat S_{k+\ell}(X,\widehat{\mathcal U}) @ > {} >>(\widehat S_{k+\ell}(X),\widehat{\mathcal U^+_{S_{k+\ell}}})
\end{CD}
\end{equation}
Here  the map in the first horizontal line is induced by the embedding
$\widehat S_{\ell}(X,\widehat{\mathcal U})
\to
\widehat{\mathcal U_{S_{\ell}}^+}$.
The map in the second horizontal line is give by (4).
The map in the first vertical column is given by Proposition \ref{prop2813}.
The map in the second vertical column
is given by (2). $\blacksquare$
\end{enumerate}
\end{shitu}
\begin{rem}
\begin{enumerate}
\item
Here we used the notion of covering space of K-space we discuss
in Section \ref{sec:cover} to formulate the compatibility
condition in Situation \ref{sit1526} at the
corner of general codimension.
In our application in Sections \ref{sec:systemline1}-\ref{sec:systemtree2}, the
stratum $\widehat S_{k}(\widehat S_{\ell}(X),\widehat {\mathcal U^+_{S_k}})$
is a disjoint union of $(k+\ell)!/k!\ell!$ copies of $\widehat S_{k+\ell}(X,\widehat{\mathcal U^+_{S_{k+\ell}}})$.
So the notion of covering space of K-space
is not necessary, there.
The result in the generality stated here will become necessary to
study the case of higher genus Lagrangian Floer theory and
or symplectic field theory.
\item
If  $\widehat{\mathcal U^+_{S_{k+\ell}}}$ is a restriction
of Kuranishi structure $\widehat {\mathcal U^+}$ such
that
$
\widehat{\mathcal U}
<
\widehat{\mathcal U^{+}}
$,
then Condition (1)-(5) above follows from
Proposition \ref{prop2813}.
Proposition \ref{prop528} below may be regarded as a converse of this statement.
\end{enumerate}
\end{rem}
\begin{defn}\label{def:inwardcollar}
In Situation \ref{sit1526}, we define
\begin{equation}
X_0 := X \setminus
X^{\boxminus\tau}.
\end{equation}
Here $X^{\boxminus \tau}$ is a subset of $X$ such that $(X^{\boxminus \tau},\widehat{\mathcal U^{\boxminus \tau}})
^{\boxplus\tau} = (X,\widehat{\mathcal U})$.
We note that $X_0$ is an open neighborhood of $S_1(X)$ in $X$.
We call $(X, \widehat{\mathcal U})^{\boxminus\tau} :=
(X^{\boxminus \tau},\widehat{\mathcal U^{\boxminus \tau}})$
{\it inward $\tau$-collaring} of $(X,\widehat{\mathcal U})$.
\index{$\tau$-collaring ! inward $\tau$-collaring}
\end{defn}
Now the next proposition is our main result of this subsection.
We complete the proof at the end of this subsection.
\begin{prop}\label{prop528}
Under Situation \ref{sit1526},
for any $0<\tau'<\tau$, there
exists a $\tau'$-collared Kuranishi structure
$\widehat{\mathcal U^+}$
on $X_0$ with the following properties.
\begin{enumerate}
\item
The restriction
of $\widehat{\mathcal U^+}$ to $\widehat S_k(X)$
is isomorphic to $\widehat {\mathcal U_{S_k}^+}$
as $\tau'$-collared Kuranishi structures.
\item
There exists an embedding of $\tau'$-collared Kuranishi
structres
$
\widehat{\mathcal U}\vert_{X_0}
\to
\widehat {\mathcal U^+}$.
\item
The restriction of (2) to $\widehat S_kX$
coincides with the one induced from
(\ref{embeddingbetween}) under the identification (1).
\item
There exists an isomorphism between the
K-spaces $\widehat S_{k}(X,\widehat {\mathcal U^+})$
and $(\widehat S_{k}(X),
\widehat {\mathcal U_{S_k}^+})$
such that the following diagram of K-spaces commutes.
\begin{equation}\label{prop1519}
\begin{CD}
\widehat S_{k}(\widehat S_{\ell}(X,\widehat {\mathcal U^+}))
@>{\cong}>>\widehat S_{k}(\widehat S_{\ell}(X),
\widehat {\mathcal U^+_{S_\ell}})
\\
@ VV{}V @ VV{}V\\
\widehat S_{k+\ell}(X,\widehat {\mathcal U^+}) @>\cong>>
(\widehat S_{k+\ell}(X),
\widehat {\mathcal U^+_{S_{k+\ell}}})
\end{CD}
\end{equation}
Here the first horizontal arrow is (1), the
second horizontal arrow is one claimed in (4),
the left vertical arrow is given by
Proposition \ref{prop2813} and the right vertical arrow is
given in Situation \ref{sit1526} (2).
\item
The two embeddings
$\widehat{\mathcal U}\vert_{\widehat S_{k}(X,\widehat{\mathcal U^+})}
\to
\widehat{\mathcal U^+}\vert_{\widehat S_{k}(X,\widehat{\mathcal U^+})}$
and
$\widehat{\mathcal U}\vert_{\widehat S_{k}(X,\widehat{\mathcal U^+})}
\to
\widehat{\mathcal U_{S_k}^+}$
coincide via the isomorphism in (4).
Here the first embedding is induced by the embeddings
$\widehat{\mathcal U}
\to \widehat{\mathcal U^{+}}$ and
the second embedding is as in Situation \ref{sit1526} (4).
\end{enumerate}
\end{prop}

\subsubsection{Extension theorem
for a single collared Kuranishi chart}\label{subsubsec:extKurachart}
The main part of the proof of Proposition \ref{prop528}
is to prove the corresponding result for one Kuranishi chart.
For this purpose we consider the following situation.
\begin{shitu}\label{sit1526new}
Let ${\mathcal U}$
be a $\tau$-collared Kuranishi chart of $X$
and $ {\mathcal U^+_{\partial}}$
a $\tau$-collared Kuranishi chart of $\partial X$.
We assume that there exists an embedding
\begin{equation}\label{embeddingbetween2}
\partial{\mathcal U}
 \to
{\mathcal U^+_{\partial}}
\end{equation}
of $\tau$-collared Kuranishi charts and the
following conditions are satisfied:
\begin{enumerate}
\item
For each $k\ge 1$
there exists a $\tau$-collared Kuranishi chart
${\mathcal U}^+_{S_k}$
on
$
\widehat S_{k}(X)
$ such that
${\mathcal U}^+_{S_1}
={\mathcal U}^+_{\partial}$.
\item
The orbifold
$\widehat S_{k}(U^+_{S_{\ell}})$
is isomorphic to the $(k+\ell)!/k!\ell!$ fold covering space of $U^+_{S_{k+\ell}}$.
The restriction of the obstruction bundle and Kuranishi map of ${\mathcal U}^+_{S_{\ell}}$ to  $\widehat S_{k}(U^+_{S_{\ell}})$
are pull-backs of ones of ${\mathcal U}^+_{S_{k+\ell}}$.
\item There exists an embedding
${\mathcal U}\vert_{\widehat S_{k}(X,{\mathcal U})} \to {\mathcal U}^+_{S_k}$
of $\tau$-collared Kuranishi charts.
\item
The following diagram commutes.
\begin{equation}\label{diagin26277XXrev}
\begin{CD}
\widehat S_{k_1}(\widehat S_{k_2}({\mathcal U_{S_{k_3}}^+})) @ >{\pi_{k_1,k_2}}>>
\widehat S_{k_1+k_2}({\mathcal U^+_{S_{k_3}}}) \\
@ V{}VV @ VV{}V\\
\widehat S_{k_1}({\mathcal U_{S_{k_2+k_3}}^+}) @ > {} >>
{\mathcal U^+_{S_{k_1+k_2+k_3}}}
\end{CD}
\end{equation}
Here $\pi_{k_1,k_2}$ is the covering map in Proposition \ref{prop2813}.
The right vertical and lower horizontal arrows are the covering maps given in (2).
The left vertical arrow is induced by the covering map
$\widehat S_{k_2}({\mathcal U_{S_{k_3}}^+})
\to {\mathcal U_{S_{k_2+k_3}}^+}$ of (2).
\item There exists an embedding
${\widehat S_{k}({\mathcal U})} \to {\mathcal U_{S_k}^+}$.
\item
The following diagram commutes.
\begin{equation}\label{diag15main2}
\begin{CD}
\widehat S_{k}(\widehat S_{\ell}({\mathcal U}))
@ > {} >>\widehat S_{k}({\mathcal U^+_{S_{\ell}}})
\\
@ V{}VV @ VV{}V \\
\widehat S_{k+\ell}({\mathcal U}) @ > {} >>{\mathcal U}^+_{S_{k+\ell}}
\end{CD}
\end{equation}
The maps are as in the case of Diagram (\ref{diag15main}).
$\blacksquare$
\end{enumerate}
\end{shitu}
The following is the extension theorem
for a single collared Kuranishi chart.
\begin{lem}\label{prop528new}
Under
Situation \ref{sit1526new},
we put
$$
U_0 := U \setminus U^{\boxminus\tau}.
$$
Then for any $0< \tau'<\tau$ there
exists a $\tau'$-collared Kuranishi chart
${\mathcal U^+}$
of $X_0 = X \setminus 
X^{\boxminus\tau}$ with the following properties.
\begin{enumerate}
\item
The restriction
of ${\mathcal U^+}$ to $\widehat S_k(U)$
is isomorphic to $\mathcal U^+_{S_k}$
as $\tau'$-collared Kuranishi charts.
\item
There exists an embedding of $\tau'$-collared Kuranishi charts
$
\mathcal U\vert_{U_0}
\to
{\mathcal U}^+$.
\item
The restriction of (2) to $\widehat S_k(U)$
coincides with one induced from
Situation \ref{sit1526} (4) under the identification (1).
\item
The following diagram commutes.
\begin{equation}\label{prop1519}
\begin{CD}
\widehat S_{k}({\widehat S_{\ell}}({\mathcal U^+}))
@>{\cong}>>\widehat S_{k}({\mathcal U_{S_{\ell}}^+})
\\
@ VV{}V @ VV{}V\\
\widehat S_{k+\ell}({\mathcal U}^+)  @>\cong>>
{\mathcal U_{S_{k+\ell}}^+}
\end{CD}
\end{equation}
Here the first horizontal arrow is induced from (1). The
second horizontal arrow is (1). The left vertical arrow is given by
Proposition \ref{prop2813}. The right vertical arrow is induced by a map
given in Situation \ref{sit1526new} (2).
\item
The embeddings
${\widehat S_{k}}({\mathcal U})\vert_{U_0}
\to {\widehat S_{k}}(\mathcal U^+)$
and
${\widehat S_{k}}({\mathcal U})\vert_{U_0}
\to {\mathcal U}^+_{S_{k}}$
coincide via the
isomorphism in Situation  \ref{sit1526new}  (3).
Here
the first embedding is induced by the embedding
${\mathcal U} \vert_{U_0}
\to  {\mathcal U^{+}}$
and the second embedding is as in Situation  \ref{sit1526new}  (3).
\end{enumerate}
\end{lem}
The proof of Lemma \ref{prop528new} occupies
Subsubsections \ref{subsubsec:constructionU}--\ref{subsubsec:pfprop528new}.
\subsubsection{Construction of Kuranishi chart $\mathcal U^+$}
\label{subsubsec:constructionU}
Let $p' \in S_1(U)$.
Firstly we will construct a Kuranishi chart $\mathcal U^+_{p'}$
mentioned in Lemma \ref{prop528new}.
We use Kuranishi neighborhoods of various $\tilde p' \in \widehat S_k(U)$
with $\pi(\tilde p') = p'$.
The Kuranishi neighborhoods we use are those of the Kuranishi chart
$\mathcal U^+_{S_k}$ given in
Situation \ref{sit1526new} (1).
We will modify and glue them in a canonical way to obtain $\mathcal U^+_{p'}$.
The detail is in order.
\par
We first begin with describing the situation of
the Kuranishi chart $\mathcal U$ in Situation \ref{sit1526new}
in more detail.
Suppose $p \in \overset{\circ}S_k(U)$ for some $k$.
Let $\frak V_{\frak r}
= (V_{\frak r},\Gamma_{\frak r},\phi_{\frak r})$
be an orbifold chart of $U$, which is the
underlying orbifold of our Kuranishi chart $\mathcal U$.
(Here $\frak r$ stands for the index of orbifold charts.)
Let $p \in U_{\frak r} = V_{\frak r}/\Gamma_{\frak r}$
such that the base point $o_{\frak r}$ of the chart goes to $p$,
(Definition \ref{2661}) that is,
$p = \phi_{\frak r}({o_{\frak r}}) \in \overset{\circ}S_k(U)$.
Then we may assume that $V_{\frak r} \subset \overline V_{\frak r}
\times [0,1)^{k}$ and $o_{\frak r}
= (\overline o_{\frak r},(0,\dots,0))$.
Here $\overline V_{\frak r}$ is a manifold without boundary.
We have a representation $\sigma : \Gamma_{\frak r}
\to {\rm Perm}(k)$ by the definition of admissible orbifold (See
Definition \ref{defn297} (1)(b)),
where $ {\rm Perm}(k)$ is the group of permutations of
$\{1,\dots,k\}$.
\par
Let $A \subset \{1,\dots,k\}$.
We put
\begin{equation}\label{1514formula}
\aligned
\overset{\circ}S_A([0,1)^k)
&=
\{(t_1,\dots,t_k) \in  [0,1)^{k}
\mid i\in A \Rightarrow t_i = 0,
\,\, i \notin A \Rightarrow t_i > 0\}, \\
\overset{\circ}S_A(V_{\frak r})
&=
V_{\frak r} \cap (\overline V_{\frak r} \times \overset{\circ}S_A([0,1)^k)),
\\
\Gamma_{\frak r}^A
&=
\{ \gamma \in \Gamma_{\frak r}
\mid \gamma A = A\}.
\endaligned
\end{equation}
The subset $A$ determines a point $p(A)$ of $\widehat S_{\ell}(U)$
which goes to $p$ by the map
$\widehat{S}_{\ell}(U) \to U$.
Here $\ell =\# A$.
An orbifold chart of $\widehat S_{\ell}(U)$ at
$p(A)$ is given by
$S_A(V_{\frak r})$ which is the closure of $\overset{\circ}S_A(V_{\frak r})$,
the isotropy group $\Gamma_{\frak r}^A$ and $\psi_A$
which is a lift of the restriction of  $\psi$ to $S_A(V_{\frak r})$.
We put
\begin{equation}\label{formulaU(A)}
V_{\frak r}(p;A)
=
(S_A(V_{\frak r}))^{\boxplus\tau} \times [-\tau,0)^{A}.
\end{equation}
Let us elaborate (\ref{formulaU(A)}).
We first note
$$
S_A(V_{\frak r}) \subset \overline V_{\frak r} \times [0,1)^{A^c}
$$
where $A^c = \{1,\dots,k\} \setminus A$.
(In fact, if $i \in A$ then $t_i = 0$ for an element of
$S_A(V_{\frak r})$.)
We have a retraction map
$$
\mathcal R_{A^c} : \overline V_{\frak r} \times [-\tau,1)^{A^c}
\to \overline V_{\frak r} \times [0,1)^{A^c},
$$
which changes $t_i < 0$ to $t_i = 0$. Then we find
$$
(S_A(V_{\frak r}))^{\boxplus\tau}
= (\mathcal R_{A^c})^{-1}(S_A(V_{\frak r}))
\subset \overline V_{\frak r} \times [-\tau,1)^{A^c}.
$$
Let $\Pi_{A^c} : \overline V_{\frak r} \times [-\tau,1)^{k}
\to \overline V_{\frak r} \times [-\tau,1)^{A^c}$
be the projection.
Then we can write
\begin{equation}
\aligned
V_{\frak r}(p;A)
=
& \{
y = (\overline y,(t_1,\dots,t_k)) \in
\overline V_{\frak r} \times [-\tau,1)^{k}
~ \vert ~ \\
& \Pi_{A^c}(y) \in (S_A(V_{\frak r}))^{\boxplus\tau},
i \in A \Rightarrow t_i \in [-\tau,0) \}.
\endaligned
\end{equation}
\begin{rem}
We observe that the space $V_{\frak r}(p,A)$, when it is written as (\ref{formulaU(A)}), is defined by the data of $S_A(V_{\frak r})$ only.
In particular, it is independent of $\overset{\circ}S_0(V_{\frak r})$.
Note, for our extension $V_{\frak r}^+$, we are given only
data on $S_1(V_{\frak r}^+)$.
We will use this remark to construct $V_{\frak r}^+(p,A)$ in this situation.
\end{rem}
We next suppose that $B \supset A$ with $\# B = \ell + m$.
Then the triple $(p,A,B)$ determines a point
$p(A,B) \in \widehat S_{m}(\widehat S_{\ell}(U))$.
We consider the maps $\pi_m : \widehat S_{m}(\widehat S_{\ell}(U))
\to S_{m}(\widehat S_{\ell}(U))$
and $\pi_{m,\ell} : \widehat S_{m}(\widehat S_{\ell}(U))
\to \widehat S_{m+\ell}(U)$
defined in Proposition \ref{prop2813}.
We have $\pi_m(p(A,B)) = p(A)$ and
$\pi_{m,\ell}(p(A,B)) = p(B)$.
\par
We put $\Gamma^{A,B}_{\frak r} = \Gamma^{A}_{\frak r} \cap \Gamma^{B}_{\frak r}$.
The map
$S_A(V_{\frak r}) /\Gamma^{A,B}_{\frak r}  \to S_A(V_{\frak r})/\Gamma^{A}_{\frak r} $ is a restriction of $\pi_m$
and the map $S_A(V_{\frak r})/\Gamma^{A,B}_{\frak r}  \to S_A(V_{\frak r})/\Gamma^{B}_{\frak r} $
is a restriction of $\pi_{m,\ell}$.
We note
$$
V_{\frak r}(p;B) = \{(\overline y,(t_1,\dots,t_k)) \in V_{\frak r}(p;A)
\mid i \in B \Rightarrow t_i \in [-\tau,0)\}.
$$
Therefore we get $V_{\frak r}(p;B) \subset V_{\frak r}(p;A)$.
\par
Using the expression (\ref{formulaU(A)}),
we can rewrite the embedding
$V_{\frak r}(p;B) \subset V_{\frak r}(p;A)$ as follows.
We have
\begin{equation}\label{form1526new}
S_{B\setminus A}(S_A(V_{\frak r})) \times [0,\epsilon)^{B\setminus A} \subset   S_A(V_{\frak r}).
\end{equation}
Here $S_{B\setminus A}(S_A(V_{\frak r}))$ is a subset of
$\widehat S_m(S_A(V_{\frak r}))$. Note
$\pi_0(\widehat S_m(S_A(V_{\frak r})))$ corresponds one to one
to the set of $B$ satisfying
$\{1,\dots,k\} \supset B \supset A$ with $\#B = \ell +m$.
Then $S_{B\setminus A}(S_A(V_{\frak r}))$
is the connected component corresponding to $B$ under this
one to one correspondence.
\par
(\ref{form1526new}) implies
\begin{equation}\label{form1526}
(S_{B\setminus A}(S_A(V_{\frak r})))^{\boxplus\tau} \times [-\tau,0)^{B-A}
\subset  (S_A(V_{\frak r}))^{\boxplus\tau}.
\end{equation}
The map $\pi_{m,\ell}$ induces a map
\begin{equation}\label{form1527}
\pi_{B \setminus A,A}:  S_{B \setminus A}(S_A(V_{\frak r})) \to   S_{B}(V_{\frak r}),
\end{equation}
which is an isomorphism.
(The $(m+\ell)!/m!\ell!$ different choices of the points in $\pi^{-1}_{m,\ell}
(\text{one point})$ corresponds
to $(m+\ell)!/m!\ell!$ different choices of $A$ in the given $B$.)
Therefore composing the inverse of (\ref{form1527})
and the inclusion (\ref{form1526}), we obtain
\begin{equation}\label{form1528}
\aligned
(S_{B}(V_{\frak r}))^{\boxplus\tau} \times &[-\tau,0)^{B} \\
\overset{(\pi_{B \setminus A,A}^{\boxplus\tau})^{-1}\times {\rm id}}\longrightarrow
&(S_{B\setminus A}(S_AV_{\frak r}))^{\boxplus\tau} \times [-\tau,0)^{B\setminus A}
\times [-\tau,0)^{A}\\
\longrightarrow
&(S_A(V_{\frak r}))^{\boxplus\tau} \times [-\tau,0)^{A}.
\endaligned
\end{equation}
It is easy to see that (\ref{form1528}) coincides with the inclusion
$V_{\frak r}(p;B) \subset V_{\frak r}(p;A)$.
\par
For $A \subset B\subset C$ we have
$V_{\frak r}(p;C) \subset V_{\frak r}(p;B) \subset V_{\frak r}(p;A)$.
The composition of the two embeddings $V_{\frak r}(p;C) \subset V_{\frak r}(p;B)$ and
$V_{\frak r}(p;B) \subset V_{\frak r}(p;A)$ coincides with
$V_{\frak r}(p;C) \subset V_{\frak r}(p;A)$.
This is equivalent to the commutativity of Diagram \eqref{diagin26277XX}.
(See Sublemma \ref{lem1540}.)
We put
$$
V_{\frak r}(p) = \bigcup_{A\subseteq \{1,\dots,k\}}V_{\frak r}(p;A).
$$
This is a $\Gamma_{\frak r}$ equivariant open subset of
$\mathcal R^{-1}(\{p\}) \setminus U$.
We may take $V_{\frak r}(p)/\Gamma_{\frak r}$ (together with other data) as
a Kuranishi neighborhood of the
point in $\mathcal R^{-1}(\{p\}) \setminus U$.
\par\smallskip
To construct $V_{\frak r}^+(p)$ we imitate the above description
using only the data given on the boundary as follows.
\par
The Kuranishi chart $\mathcal U^+_{S_\ell}$ given in Situation \ref{sit1526new} induces
$V^+_{\frak r,S_A}/\Gamma_{\frak r}^A$. (It is an open subset of $U^+_{S_\ell}$. Also recall $\#A =\ell$.)
We define
\begin{equation}\label{formulaU(A)++}
V_{\frak r}^+(p;A)
=
(V^+_{\frak r,S_A})^{\boxplus\tau} \times [-\tau,0)^{A}.
\end{equation}
For $B \supset A$ we have an embedding denoted by $h_{A,B}$
\begin{equation}\label{cornercanemb1515}
h_{A,B} : (S_{B\setminus A}(V^+_{\frak r,S_A}))^{\boxplus\tau} \times [-\tau,0)^{B\setminus A}
\hookrightarrow (V^+_{\frak r,S_A})^{\boxplus\tau}.
\end{equation}
The covering map $\widehat S_{m}(\mathcal U^+_{S_\ell}) \to \mathcal U^+_{S_{\ell+m}}$
given in Situation \ref{sit1526new} (2) induces
the map
\begin{equation}\label{form1527ref}
\pi'_{B \setminus A,A}:  S_{B \setminus A}(V^+_{\frak r,S_A}) \to   V^+_{\frak r,S_B}
\end{equation}
which is an isomorphism.
We define
$$
\phi_{AB} : V_{\frak r}^+(p;B) \to V_{\frak r}^+(p;A),
$$
by
\begin{equation}\label{form1528rev}
\aligned
(V^+_{\frak r,S_{B}})^{\boxplus\tau} \times &[-\tau,0)^{B} \\
\overset{(\pi_{B \setminus A,A}^{\prime\boxplus\tau})^{-1}\times {\rm id}}\longrightarrow
&(S_{B\setminus A}(V^+_{\frak r,S_{A}}))^{\boxplus\tau} \times [-\tau,0)^{B\setminus A}
\times [-\tau,0)^{A}\\
\overset{h_{A,B}\times {\rm id}}
\longrightarrow
&(V^+_{\frak r,S_{A}})^{\boxplus\tau} \times [-\tau,0)^{A}.
\endaligned
\end{equation}
\begin{sublem}\label{lem1540}
If $C \supset B \supset A$, then $\phi_{AB}\circ \phi_{BC} = \phi_{AC}$.
\end{sublem}
\begin{proof}
The sublemma follows from the commutativity of Diagram (\ref{diagin26277XXrev})
as follows.
Recall that $h_{A,B}$ is
the inclusion map in (\ref{cornercanemb1515}).
Then the following diagram
commutes.
\begin{equation}\label{diagram1537}
\begin{CD}
(V^+_{\frak r,S_{B}})^{\boxplus\tau}
@<{\pi_{B\setminus A,A}^{\prime \boxplus\tau}}<<
(S_{B\setminus A}(V^+_{\frak r,S_{A}}))^{\boxplus\tau}
\\
@ AA{h_{B,C}}A @ AA{h_{B,C}}A\\
(S_{C\setminus B}(V^+_{\frak r,S_{B}}))^{\boxplus\tau} \times [-\tau,0)^{C\setminus B}
@<({S_{C\setminus B}(\pi'_{B\setminus A,A}))^{\boxplus\tau}\times {\rm id}}<<
(S_{C\setminus B}(S_{B\setminus A}(V^+_{\frak r,S_{A}})))^{\boxplus\tau}
\times [-\tau,0)^{C\setminus B}
\end{CD}
\end{equation}
\par\smallskip
\noindent
The commutativity of Diagram \eqref{diagram1537}
is a consequence of the commutativity of Diagram \eqref{diagin26277XXrev} and
the fact that $\pi'_{B\setminus A,A}$
is a diffeomorphism of cornered manifolds and the definition of
$h_{*,*}$.
We also have the following commutative diagram.
\begin{equation}\label{diagram1538}
\begin{CD}
(V^+_{\frak r,S_{A}})^{\boxplus\tau} \times [-\tau,0)^{A}
@<{h_{A,B}\times {\rm id}_A}<<
(S_{B\setminus A}(V^+_{\frak r,S_{A}}))^{\boxplus\tau}\times [-\tau,0)^{B}\\
@ AA{h_{A,C}\times {\rm id}_A}A @ AA{h_{B,C}\times {\rm id}_B}A\\
(S_{C\setminus A}(V^+_{\frak r,S_{A}}))^{\boxplus\tau}\times [-\tau,0)^{C}
@<(\pi_{C\setminus B,B\setminus A})^{\boxplus\tau}\times {\rm id}_C<<
(S_{C\setminus B}(S_{B\setminus A}(V^+_{\frak r,S_{A}})))^{\boxplus\tau}
\times [-\tau,0)^{C}
\end{CD}
\end{equation}
\par\smallskip
\noindent
Note that $\pi_{C\setminus B,B\setminus A}$ in Diagram \eqref{diagram1538} is the map in
Proposition \ref{prop2813}.
(The map $\pi'_{B\setminus A,A}$ in Diagram \ref{diagram1537} is one in
Situation \ref{sit1526new} (2).)
The commutativity of  Diagram \ref{diagram1538}  is an immediate
consequence of the definition of $\pi_{C\setminus B,B\setminus A}$
and $h_{*,*}$.
Therefore we have
$$
\aligned
&\phi_{AB}\circ \phi_{BC} \\
&=(h_{A,B}\times {\rm id}) \circ
((\pi_{B\setminus A,A}^{\prime\boxplus\tau})^{-1}\times {\rm id})
\circ
(h_{B,C}\times {\rm id}) \circ
((\pi_{C\setminus B,B}^{\prime\boxplus\tau})^{-1}\times {\rm id})
\\
&=
(h_{A,B}\times {\rm id}) \circ
(h_{B,C}\times {\rm id}) \circ
((S_{C\setminus B}(\pi'_{B\setminus A,A}))^{\boxplus\tau})^{-1} \times {\rm id})
\circ
((\pi_{C\setminus B,B}^{\prime\boxplus\tau})^{-1}\times {\rm id})
\\
&=(h_{A,C}\times {\rm id}) \circ
((\pi_{C\setminus B,B\setminus A})^{\boxplus\tau} \times  {\rm id})
\\
&\qquad\qquad\qquad\,\circ
((S_{C\setminus B}(\pi'_{B\setminus A,A}))^{\boxplus\tau})^{-1} \times {\rm id})
\circ
((\pi_{C\setminus B,B}^{\prime\boxplus\tau})^{-1}\times {\rm id})
\\
&= (h_{A,C}\times {\rm id}) \circ
((\pi_{C\setminus A,A}^{\prime\boxplus\tau})^{-1}\times {\rm id})
\\
&= \phi_{AC}.
\endaligned
$$
Here the first equality is the definition.
The second equality is the commutativity of  Diagram (\ref{diagram1537}).
The third equality is the commutativity of
Diagram (\ref{diagram1538}).
The fourth equality is the commutativity of Diagram (\ref{diagin26277XXrev}).
\end{proof}
We consider the disjoint union
$$
\coprod_A V_{\frak r}^+(p;A)
$$
and define $\sim$ on it as follows.
For $x \in V_{\frak r}^+(p;A)$ and $y \in V_{\frak r}^+(p;B)$ we say $x \sim y$
if and only if there exist $C$ and $z \in V_{\frak r}^+(p;C)$
such that $x = \phi_{AC}(z)$ and $y = \phi_{BC}(z)$.
\begin{sublem}
$\sim$ is an equivalence relation.
\end{sublem}
\begin{proof}
It suffices to prove the transitivity.
Let $x = (x',(t_i)_{i\in A})$ where $x' \in (V^+_{\frak r,S_A})^{\boxplus\tau}$
and $t_i \in [-\tau,0)$ for $i \in A$.
We furthermore write
$x' = (\overline x,(t_i)_{i\in A^c})$.
We observe that
$x'$ is in the image of
$(S_{C\setminus A}(V^+_{\frak r,S_{A}}))^{\boxplus\tau} \times [-\tau,0)^{C\setminus A}$
if and only if $t_i < 0$ for all $i\in C$.
Therefore for each $x$ there exists {\it unique} $C$ such that
\begin{enumerate}
\item
$x \in {\rm Im}(\phi_{AC})$.
\item
If $x \in {\rm Im}(\phi_{AD})$ then $D \subseteq C$.
\end{enumerate}
Transitivity follows from this fact, Sublemma \ref{lem1540} and the fact that $\phi_{AB}$
is injective.
\end{proof}
We define
$$
V_{\frak r}^+(p) = \coprod_A V_{\frak r}^+(p;A)/\sim.
$$
\begin{sublem}\label{lemlemHausdorff}
The quotient space $V_{\frak r}^+(p)$ is Hausdorff with respect to the quotient topology.
\end{sublem}
\begin{proof}
Let $x,y \in V_{\frak r}^+(p)$ such that $x \ne y$.
We take $C_x$ and $C_y$ as in (1)(2) above
and take the representatives $\tilde x \in V_{\frak r}^+(p;C_x)$ and
$\tilde y \in V_{\frak r}^+(p;C_y)$, respectively.
If $C_x = C_y = C$, we can find an open set $U_x, U_y$ in
$V_{\frak r}^+(p;C)$ such that $\tilde x \in U_x$,  $\tilde y \in U_y$
and $U_x\cap U_y = \emptyset$.
The images of $U_x$ and $U_y$ in $V_{\frak r}^+(p)$ separate $x,y$.
\par
Suppose $C_x \ne C_y$. We may assume that there exists $j
\in C_x \setminus C_y$.
We write
$\tilde x = (x',(t_i^0)_{i\in C_x})$, $x' = (\overline x,(t_i^0)_{i\in C_x^c})$.
We also write
$\tilde y = (y',(s_i^0)_{i\in C_y})$, $y' = (\overline y,(s_i^0)_{i\in C_y^c})$.
Then $t_j^0 <0$ and $s_j^0 \ge0$.
Let $U_x$ be the set of all points in
$V_{\frak r}^+(p;C_x)$ such that $t_j < t_j^0/2$
and $U_y$ be the set of all points in
$V_{\frak r}^+(p;C_y)$ such that $s_j > t_j^0/2$.
They induce disjoint open sets in $V_{\frak r}^+(p)$ containing $x$ and $y$ respectively.
\end{proof}
Since $\phi_{AB}$'s are open embeddings of manifolds, Sublemma \ref{lemlemHausdorff}
implies that $V_{\frak r}^+(p)$ is a smooth manifold.
\par
We next define an obstruction bundle and a Kuranishi map on it.
(\ref{formulaU(A)++}) shows that each $V_{\frak r}^+(p;A)$
comes with an obstruction bundle and a Kuranishi map on it.
We denote them by $\mathcal V^+_{{\frak r},p;A}$.
Moreover (\ref{form1528rev}) implies that each $\phi_{AB}$ is covered
by the bundle isomorphism and Kuranishi map is compatible with it.
Moreover the identity $\phi_{AB}\circ \phi_{BC} = \phi_{AC}$
is promoted to the identity among bundle maps.
Therefore we obtain an obstruction bundle and a Kuranishi map on
$V_{\frak r}^+(p)$.
(This is nothing but \cite[Lemma 3.17]{part11}.)
We denote them by $\mathcal E_{{\frak r},p}$ and $s^+_{{\frak r},p}$.
We can also define $\psi_{{\frak r}, p}^+ : (s^+_{{\frak r},p})^{-1}(0) \to X_0$
in an obvious way.
The following is immediate from the construction.
\begin{sublem}\label{gammaequivequiv}
For each $\gamma \in \Gamma_p$ and $A$ there exists
$\varphi_{\gamma,A} : \mathcal V^+_{p;A} \to \mathcal V^+_{p;\gamma A}$.
Moreover
$\varphi_{\gamma,A} \circ \phi_{AB} = \phi_{(\gamma A) (\gamma B)} \circ \varphi_{\gamma,B}$.
\end{sublem}
Then Sublemma \ref{gammaequivequiv} implies that
$$
(V_{\frak r}^+(p),
\Gamma_{{\frak r},p},\mathcal E_{{\frak r},p},s^+_{{\frak r},p},
\psi_{{\frak r},p}^+)
$$
is a Kuranishi chart at each point of $\mathcal R^{-1}(p)$.
\par\smallskip
Next we define coordinate change.
Let $p' \in \mathcal R^{-1}(p) \cap U_0$,
$q' \in \mathcal R^{-1}(q) \cap U_0$
and $q' \in \psi_p^+((s^+_p)^{-1}(0))$.
(We assume $q \in S_1(U)$.)
Then we have
$q \in \psi_{A_0}(s_{A_0}^{-1}(0))$
for some $A_0$.
Here $s_{A_0}$ is a Kuranishi map
of a Kuranishi chart  $\mathcal U_{S_{k'}}$ at $p$.
We put $k = \#A_0$.
\par
Let $p \in \overset{\circ}S_k(X)$.
We will use the same notations
as those used in the construction of the
Kuranishi chart
$(V_{\frak r}^+(p),
\Gamma_{{\frak r},p},\mathcal E_{{\frak r},p},s^+_{{\frak r},p},\psi_{{\frak r},p}^+)$.
\par
We take $k' \le k$ such that $q \in \overset{\circ}S_{k'}(X)$\
and we may assume $q \in \psi_{A_q}(s_{A_q}^{-1}(0))$
with $\#A_q = k'$.
We have $\tilde q \in S_{A_q}(V_{\frak r})$ which parametrizes $q$.
Moreover for each $A \subseteq A_{q}$ there exists $\tilde q_A \in S_{A}(V_{\frak r})$ which parametrizes $q$.
\par
For $A \subset A_q$ with
$\#A = \ell$, the Kuranishi chart $\mathcal U^+_{S_{\ell}}$
gives an orbifold chart $\frak V_{\frak o,A}^+
= (V_{\frak o,A}^+,\Gamma^A_{\frak o},\psi_{\frak o,A}^+)$
at $\tilde q_A$.
The coordinate change of the underlying orbifold $U^+_{S_{\ell}}$ of
$\mathcal U^+_{S_{\ell}}$ induces a
group homomorphism
$h_{\frak r \frak o}^A : \Gamma^A_{\frak o} \to \Gamma^A_{\frak r}$
and an $h_{\frak r \frak o}^A$ equivariant
embedding
\begin{equation}\label{defnvarphicha15}
\varphi^{A+}_{\frak r\frak o}
: V_{\frak o,A}^+ \to V_{\frak r,A}^+.
\end{equation}
Moreover the admissibility of our orbifolds implies that
we have
\begin{equation}
\aligned
V_{\frak o,A}^+
&\subset \overline V_{\frak o,A}^+ \times [0,1)^{A_q\setminus A},
\\
V_{\frak r,A}^+
&\subset \overline V_{\frak r,A}^+ \times [0,1)^{A_p\setminus A}
\endaligned
\end{equation}
and
$$
\varphi^{A+}_{\frak r\frak o}(\overline y,(t_i)_{i\in A_q\setminus A})
=
\left(\varphi^{A+}_{\frak r\frak o,0}(\overline y,(t_i)), \enskip
(\varphi^{A+}_{\frak r\frak o,j}(\overline y,(t_i))
_{j\in A_p\setminus A}\right)
$$
such that :
\begin{enumerate}
\item
$\varphi^{A+}_{\frak r\frak o,0}$ is admissible.
\item
For $j \in A_q\setminus A$,
$\varphi^{A+}_{\frak r\frak o,j} -t_j$
is exponentially small near the boundary.
\item
For $j \in A_p\setminus A_q$,
$\varphi^{A+}_{\frak r\frak o,j}$ is admissible.
\end{enumerate}
Below we will extend $\varphi^{A+}_{\frak r\frak o}$ to
$$
\varphi^{A+\boxplus\tau}_{\frak r\frak o} :
V^+(q,A) \to V^+(p,A).
$$
Note
$$
\aligned
V_{\frak o}^+(q,A) &= (V_{\frak o,A}^+)^{\boxplus\tau} \times [-\tau,0)^A, \\
V_{\frak r}^+(p,A) &= (V_{\frak r,A}^+)^{\boxplus\tau} \times [-\tau,0)^A.
\endaligned
$$
Let $y = (y',(t_i)_{i\in A}) \in V_{\frak o}^+(q,A)$.
We define
\begin{equation}
\varphi^{A+\boxplus\tau}_{\frak r\frak o}(y)
=
((\varphi^{A+}_{\frak r\frak o})^{\boxplus\tau}(y'),(t_i)_{i\in A})
\in V_{\frak r}^+(p,A).
\end{equation}
\begin{sublem}\label{coordinatechangeABcomp}
If $A \subseteq B \subseteq A_q$ then
$$
\varphi^{A+\boxplus\tau}_{\frak r\frak o}\circ \phi_{AB}
=
\phi_{AB} \circ \varphi^{B+\boxplus\tau}_{\frak r\frak o}.
$$
\end{sublem}
\begin{proof}
This is a consequence of the fact that two
maps appearing in (\ref{form1528rev}) is
compatible with the coordinate change.
\end{proof}
By Sublemma \ref{coordinatechangeABcomp} and
\cite[Lemma 3.18]{part11}
we can glue $\varphi^{A+\boxplus\tau}_{\frak r\frak o}$ for
various $A$ to obtain a map
$$
\varphi^{+\boxplus\tau}_{\frak r\frak o} :
V_{\frak o}^+(q) \to V_{\frak r}^+(p)
$$
and a bundle map
$$
\widehat\varphi^{+\boxplus\tau}_{\frak r\frak o} :
\mathcal E_{{\frak o},q}^+\to \mathcal E_{{\frak r},p}^+.
$$
They are $h_{pq} : \Gamma_q \to \Gamma_p$ equivariant
by construction. Moreover the Kuranishi maps
and parametrizations $\psi_p^+$, $\psi_q^+$ are compatible
to it.
We have thus constructed the coordinate change.
\par
The cocycle condition among the coordinate changes follows
from the cocycle condition among the coordinate
changes of various $\mathcal U^+_{S_{\ell}}$.
\footnote{Since we are constructing the space $U$ together
with its orbifold structure, we need to check the cocycle condition.
It is easy to check however.}
\subsubsection{Completion of the proof of Lemma \ref{prop528new}}
\label{subsubsec:pfprop528new}
To complete the proof of Lemma  \ref{prop528new}
it suffices to check that the Kuranishi chart $\mathcal U^+$
we obtained has the required properties.
\begin{proof}[Proof of the $\tau'$-collaredness of $\widehat{\mathcal U^+}$]
Each $V^+(p,A) = (V^+_{\frak r,A})^{\boxplus\tau} \times [-\tau,0)^A$ is
$\tau'$-collared. The open embedding $\phi_{AB}$ which we used to glue them
are $\tau'$-collared.
The coordinate change is obtained by gluing $\varphi^{A+\boxplus\tau}_{\frak r\frak o}$,
which is $\tau'$-collared.
Therefore, we can construct a $\tau'$-collared Kuranishi structure
in the same way as in the proof of Lemma \ref{lem153535}.
\end{proof}
\begin{rem}\label{rem1551}
We note that the definition of $\tau'$-collared Kuranishi structure
produced in Lemma \ref{lem153535} has slight ambiguity.
Namely it depends on the choice of the coordinate
$\overline V_{\widehat p} \times [-\tau,c)^k$ on which
isotropy group acts by permutation on the $[-\tau,c)^k$ factor.
We mentioned this point in the footnote in the proof of
Lemma \ref{lem153535}.
Note that on the region $[-\tau,0]$ there is no ambiguity like that at all.
Actually we observe that the proof of the
$\tau'$-collared-ness of $\widehat{\mathcal U^+}$
is the {\it only} place where Lemma \ref{lem153535}
is used in the application.
In the situation of  Lemma \ref{lem153535}
the coordinate $\overline V_{\widehat p} \times [-\tau,c)^k$ for which
the action of isotropy group is given by exchanging the factor
on $[-\tau,c)^k$ is given.
Namely in our situation $[-\tau,c)$ corresponds to
$[-\tau,0]$ where we shift the parameter so that
$0 \in [-\tau,c)$  corresponds $\tau' \in [-\tau,0]$.
Thus the ambiguity mentioned in the footnote in the proof of
Lemma \ref{lem153535} is not at all
an issue here.
\par
Note that Lemma \ref{lem153535} is literally correct with the proof given.
The concern of this remark is the precise meaning of the
word `canonical' in Lemma \ref{lem153535}.
\end{rem}
\begin{proof}[Proof of Lemma  \ref{prop528new} (1)]
Let $p' \in S_{\ell}(U^{\boxplus\tau} \setminus U)$ and $\mathcal R(p') = p \in S_{\ell}(U)$.
We take $k \ge \ell$ such that $p \in \overset{\circ}S_{k}(U)$.
We use the same notations used in the construction of $V^+(p)$.
Take $p'_A$  a point in the underlying topological space of $\widehat S_{\ell}(U^{\boxplus\tau} \setminus U)$
which goes to $p'$. We have a corresponding point $p_A \in \widehat S_{\ell}(U)$ which goes to $p$.
The point $p_A$ corresponds to a certain subset $A \subset \{1,\dots,k\}$ with $\# A = \ell$.
An orbifold neighborhood of $p_A$ in $U^+_{S_{\ell}}$ is $V_{\frak r,S_A}^+/\Gamma_{\frak r}^A$ by definition.
\par
Note that
$
V_{\frak r}^+(p;A)
=
V_{\frak r,S_A}^+ \times [-\tau,0)^{A}
$ and
a neighborhood of $p'_A$ in $S_A(V_{\frak r}^+(p;A))$ is
$V_{\frak r,S_A}^+ \times \{(0,\dots,0)\}$ in $V_{\frak r}^+(p;A)$.
\par
Thus we have shown that $\widehat S_{\ell}(\mathcal U^+)$ and
$\mathcal U^+_{S_{\ell}}$ are locally diffeomorphic each other as orbifolds.
This diffeomorphism is compatible with the gluing by $\phi_{AB}$
and by coordinate changes.
So the underlying orbifolds of $\widehat S_{\ell}(\mathcal U^+)$ and
$\mathcal U^+_{S_{\ell}}$ are diffeomorphic.
Moreover it is covered by the bundle isomorphism of obstruction
bundles which is compatible with coordinate change and Kuranishi map.
\end{proof}
\begin{proof}[Proof of Lemma  \ref{prop528new} (2)]
By assumption there exists an embedding
$\mathcal U\vert_{S_k(U)} \to \mathcal U^+_{S_k}$. (Situation \ref{sit1526new}
(3).)
By comparing (\ref{formulaU(A)}) and (\ref{formulaU(A)++})
it induces an embedding $V_{\frak r}(p;A)\to V_{\frak r}^+(p;A)$.
\par
By the commutativity of (\ref{diag15main2})
we have the following commutative diagram
$$
\begin{CD}
V_{\frak r}(p;B)
@>{}>> V_{\frak r}^+(p;B)
\\
@ VV\phi_{AB}{}V @ VV{\phi_{AB}}V\\
V_{\frak r}(p;A)
@>{}>> V_{\frak r}^+(p;A)
\end{CD}
$$
for $B \supset A$. Therefore we have an
embedding $V_{\frak r}(p) \to V_{\frak r}^+(p)$.
It is covered by a bundle map and is $\Gamma_p$ equivariant.
Moreover it is compatible with Kuranishi map.
Thus this embedding $V_{\frak r}(p) \to V_{\frak r}^+(p)$ is promoted to an
embedding of Kuranishi charts.
\par
On the other hand, the embeddings  $V_{\ast}(p;A)\to V_{\ast}^+(p;A)$ commute with the embeddings $\varphi^{A+}_{\frak r\frak o}$ and
$\varphi^{A}_{\frak r\frak o}$.
Therefore we can glue the embeddings $V_{\ast}(p) \to V_{\ast}^+(p)$ to obtain the required
embedding.
\end{proof}
\begin{proof}[Proof of Lemma  \ref{prop528new} (3)]
This follows from the proof of  Lemma  \ref{prop528new} (1),(2).
\end{proof}
\begin{proof}[Proof of Lemma  \ref{prop528new} (4)]
This follows from the proof of  Lemma  \ref{prop528new} (1).
\end{proof}
\begin{proof}[Proof of Lemma  \ref{prop528new} (5)]
This follows from the proof of  Lemma  \ref{prop528new} (2).
\end{proof}
Therefore the proof of Lemma  \ref{prop528new} is now complete.
\subsubsection{Proof of Proposition \ref{prop528}}
\begin{proof}
It suffices to construct the coordinate change between
Kuranishi charts produced
in Lemma \ref{prop528new} and show that the coordinate changes
are compatible with various embeddings and isomorphisms
appearing in the statement of Proposition \ref{prop528}
and of Lemma \ref{prop528new}.
This is indeed straightforward.
In fact, the Kuranishi structure in Lemma \ref{prop528new}
is constructed from $\mathcal U^+_{S_k}$,
which are Kuranishi charts of $\widehat{\mathcal U^+_{S_k}}$.
They are compatible with the coordinate change by definition.
The process to construct our Kuranishi chart from $\mathcal U^+_{S_k}$
is by trivialization of the corner, $* \mapsto *^{\boxplus\tau}$, and gluing by the map
in Situation \ref{sit1526new} (2).
The former is compatible with coordinate change
as we proved in the first half of this section.
The latter is compatible since it is induced by the corresponding
map (Situation \ref{sit1526} (2)) of Kuranishi structures.
\end{proof}
\begin{rem}
What is written as $\mathcal U$ in the notation of Proposition \ref{prop528}
corresponds to  $\mathcal U^{\boxplus\tau}$  in the notation of Lemma \ref{prop528new}.
\end{rem}

\subsection{Extension of collared CF-perturbation}
\label{subsection:extcfp}
In this section we prove Proposition \ref{prop529}.
\begin{shitu}\label{situ529}
In Situation \ref{sit1526},
let $\widehat{\frak S_{\partial}^{+}}$ be a
$\tau$-collared CF-perturbation of
$\widehat{\mathcal U^+_{\partial}}$.
We assume the following:
\begin{enumerate}
\item
For each $k \ge 1$
there exists a $\tau$-collared CF-perturbation
$\widehat{\frak S_{S_k}^+}$ of
$\widehat{\mathcal U_{S_k}^+}$ such that
$\widehat{\frak S_{S_1}^+} = \widehat{\frak S_{\partial}^{+}}$.
\item
The pull-back
of $\widehat{\frak S_{S_{k+\ell}}^+}$
by
$
\pi_{k,\ell} : \widehat S_{k}(\widehat S_{\ell}(X),\widehat{\mathcal U^+_{ S_{\ell}}})
\to (\widehat S_{k+\ell}(X),\widehat{\mathcal U^+_{ S_{k+\ell}}})$
is equivalent to the restriction of $\widehat{\frak S^+_{\ell}}$. $\blacksquare$
\end{enumerate}
\end{shitu}

\begin{prop}\label{prop529}
Suppose we are in Situation \ref{situ529}. Then for any $0 < \tau' < \tau$ there exists a
$\tau'$-collared CF-perturbation
$\widehat{\frak S^+}$ of the Kuranishi structure
$\widehat{\mathcal U^+}$ obtained in Proposition \ref{prop528}
such that the restriction of $(\widehat S_k(X),\widehat{\mathcal U_{S_k}^+})$
is equivalent to $\widehat{\frak S^+_{S_k}}$.
When $\widehat{\frak S^+_{S_k}}$ varies in a uniform family,
we may take $\widehat{\frak S^+}$ to be uniform.
\end{prop}
\begin{proof}
We first consider the situation of one chart.
We use the same notations used in the construction of
the chart $V^+(p;A)$ during the proof of Lemma  \ref{prop528new}.
By assumption (Situation \ref{situ529} (1)) we are
given a CF-perturbation $\mathcal S^+_{\frak r,A}$ on $V^+_{\frak r,S_A}/\Gamma^A_{\frak r}$.
It induces $\mathcal S^{+\boxplus\tau}_{\frak r,A}$
on $(V^+_{\frak r,S_A})^{\boxplus\tau}/\Gamma^A_{\frak r}$.
We extend it by constant on the $[-\tau,0)^{A}$ factor to obtain $\mathcal S_{p,A}$ on
$V^+(p;A)/\Gamma^A_{\frak r}$.
\par
If $A' = \gamma A$ for $\gamma \in \Gamma_{\frak r}$,
then $\mathcal S_{p,A}$ is isomorphic to $\mathcal S_{p,\gamma A}$. Therefore we obtain
$\mathcal S^+_{p,\ell}$
on
\begin{equation}\label{Vraafffff}
\left(\bigcup_{A;\#A = \ell} V_{\frak r}^+(p;A)\right) / \Gamma_{\frak r}.
\end{equation}
Note that (\ref{Vraafffff}) is diffeomorphic to
$V_{\frak r}^+(p,A)/\Gamma^A_{\frak r}$.
The open sets (\ref{Vraafffff}) for various $\ell$ cover
$V_{\frak r}^+(p)$.
\par
Let $B \supset A$. Situation \ref{situ529} (2) implies that
the restriction of  $\mathcal S_{p,A}$ by the map $\phi_{AB}$
is equivalent to $\mathcal S_{p,B}$.
Therefore $\mathcal S^+_{p,\ell}$ is equivalent to $\mathcal S^+_{p,m}$
on the intersection of the domains (\ref{Vraafffff}).
Thus we get a CF-perturbation
$\frak S^+_p$
on $V_{\frak r}^+(p)/\Gamma_{\frak r}$.
\par
If $q \in \psi^+_{p}((s_p^+)^{-1}(0))$, then we have
$\varphi^{A+}_{\frak r\frak o}
: V_{\frak o,A}^+ \to V_{\frak r,A}^+$, that is
(\ref{defnvarphicha15}).
Since $\widehat{\frak S^+_{S_{\ell}}}$ is a CF-perturbation,
the pull-back of $\mathcal S^+_{\frak r,A}$ by $\varphi^{A+}_{\frak r\frak o}$  is
equivalent to $\mathcal S^+_{\frak o,A}$.
Therefore
$\frak S^+_p$ and $\frak S^+_q$ are glued
to define a CF-perturbation on the union of domains.
We have thus constructed a CF-perturbation
on each of the Kuranishi charts.
The compatibility with the coordinate change
follows from the corresponding compatibility
of $\widehat{\mathcal U^+_{ S_{k}}}$'s.
Thus we have obtained a CF-perturbation $\widehat{\frak S^+}$.
\par
The equivalence of the restriction of $\widehat{\frak S^+}$ to
 $(\widehat S_k(X),\widehat{\mathcal U_{S_k}^+})$
and $\widehat{\frak S^+_{k}}$ is obvious from the construction.
The uniformity also follows from the construction.
\end{proof}

\subsection{Extension of Kuranishi structure and CF-perturbation
from a neighborhood of a compact set}
\label{subsec:extensionfromopen}
In this subsection
we prove extension lemmas of Kuranishi structure and of CF-perturbation
defined on a {\it neighborhood} of a compact set.
These lemmas are used in the next subsection.
Note that they are results about Kuranishi structure and CF-perturbation,
not about $\tau$-collared ones.

\begin{lem}\label{lem1544}
Let $K$ be a compact set of $X$ and $Z \subseteq X$ a compact neighborhood
of $K$
such that $K \subset {\rm Int}\, Z$.
Suppose we are given a Kuranishi structure
$\widehat{\mathcal U}$ on $X$ and
$\widehat{\mathcal U^{+}_{Z}}$ on $Z$.
We assume
$$
\widehat{\mathcal U}\vert_{Z} < \widehat{\mathcal U_{Z}^+}.
$$
Let $\Omega$ be a relatively compact neighborhood of $K$
in $Z$ such that
$$
K \subset \Omega \subset \overline{\Omega} \subset
{\rm Int}\, Z.
$$
We also assume the following:
\begin{enumerate}
\item[(i)]
Let $p \in K$ and $\mathcal U^{+}_{Z,p}
= (U^{+}_{Z,p},\mathcal E^{+}_{Z,p},s^{+}_{Z,p},\psi^{+}_{Z,p})$
be the
Kuranishi neighborhood of $\widehat{\mathcal U^{+}_{Z}}$ at $p$.
We assume
$$
\psi^{+}_{Z,p}((s^{+}_{Z,p})^{-1}(0)) \subset \Omega.
$$
\item[(ii)]
Let $p \in K$ and $\mathcal U_{p}
= (U_{p},\mathcal E_{p},
s_{p},\psi_{p})$ be the
Kuranishi neighborhood of $\widehat{\mathcal U}$ at $p$.
We assume
$$
\psi_{p}((s_{p})^{-1}(0)) \subset \Omega.
$$
\end{enumerate}
\par
Then there exist a Kuranishi structure $\widehat{\mathcal U^+}$ on $X$
and an embedding $\widehat{\mathcal U} \to \widehat{\mathcal U^+}$
with the following properties:
\begin{enumerate}
\item
For any $p\in K$
the Kuranishi neighborhood $\mathcal U^+_p$ of $\widehat{\mathcal U^+}$
at $p$ is isomorphic to the Kuranishi neighborhood  $\mathcal U^{+}_{Z,p}$ of
$\widehat{\mathcal U^{+}_{Z}}$ at $p$.
For any $p,q \in K$
the coordinate change between $\mathcal U^+_p$ and $\mathcal U^+_q$
coincides with the coordinate change between $\mathcal U^{+}_{Z,p}$ and
$\mathcal U^{+}_{Z,q}$.
\item
$\widehat{\mathcal U^+}\vert_{\Omega}$
is an open substructure of $\widehat{\mathcal U^{+}_{Z}}\vert_{\Omega}$.
\item
$\widehat{\mathcal U} < \widehat{\mathcal U^+}$.
\item
The next diagram commutes.
\begin{equation}\label{diag1544strict}
\xymatrix{
\widehat{\mathcal U}\vert_{\Omega} \ar[rr]^{\rm embedding} \ar[dr]_{\rm embedding} &  &
\widehat{\mathcal U^+}\vert_{\Omega} \ar[dl]^{\rm embedding} \\
& \widehat{\mathcal U_{Z}^+}\vert_{\Omega}  &
}
\end{equation}
Here the right down arrow is the open embedding given by (2).
\item
The embedding
$\widehat{\mathcal U}\vert_{K}
\to \widehat{\mathcal U_{Z}^+}\vert_{K}$
coincides
with the embedding
$\widehat{\mathcal U}\vert_{K}
\to \widehat{\mathcal U^+}\vert_{K}$
via the isomorphism (1).
\end{enumerate}
\end{lem}
\begin{proof}
We take an open set $\Omega_1 \subset X$ such that
$$
\overline \Omega \subset \Omega_1 \subset
\overline \Omega_1 \subset {\rm Int}\, Z.
$$
We next take an open set $W_1 \subset X$ such that
$$\overline W_1 \cap \overline\Omega = \emptyset, \quad
\Omega_1 \cup W_1 = X.
$$
We replace $\widehat{\mathcal U}$ and
$\widehat{\mathcal U_{Z}^+}$ by their open substructures
(but without changing the Kuranishi neighborhoods of the point $p\in K$)
if necessary, and may assume that the following holds.
\begin{enumerate}
\item
If
$\psi^+_{Z,p}(U^+_{Z,p} \cap (s^+_{Z,p})^{-1}(0)) \cap
\overline W_1 \ne \emptyset$,
then
$\psi^+_{Z,p}(U^+_{Z,p} \cap (s^+_{Z,p})^{-1}(0)) \cap \overline\Omega = \emptyset$.
\item
If
$\psi_p(U_p \cap s_p^{-1}(0)) \cap \overline W_1 \ne \emptyset$,
then
$\psi_p(U_p \cap s_p^{-1}(0)) \cap \overline\Omega = \emptyset$.
\end{enumerate}
We define a Kuranishi structure $\widehat{\mathcal U'}$ on $X$
as follows.
\begin{enumerate}
\item[(a)]
If $p \in \overline{\Omega}_1$, we put $\mathcal U'_p
= \mathcal U^+_{Z,p}$.
\item[(b)]
If $p \notin \overline{\Omega}_1$, we put $\mathcal U'_p
=
\mathcal U_p\vert_{U_p \setminus \psi^{-1}_p(\overline\Omega_1)}$.
\end{enumerate}
The coordinate change is defined as follows.
Let $q \in \psi'_p((s'_p)^{-1}(0))$.
If  $p,q \in \overline{\Omega}_1$, then
we define $\Phi'_{pq} = \Phi^+_{Z,pq}$.
If $p,q \notin \overline{\Omega}_1$, then we define
$\Phi'_{pq} = \Phi_{pq}\vert_{U_{pq} \setminus \psi^{-1}_q(\overline{\Omega}_1)}$.
Among the other two cases
$q \in \overline{\Omega}_1$,
$p \notin \overline{\Omega}_1$ cannot occur (by (b)).
We consider the remaining case
$q \notin \overline{\Omega}_1$,
$p \in \overline{\Omega}_1$.
We have an embedding $\Phi_q : \mathcal U_q
\to \mathcal U^+_{Z,q}$. We compose it with the embedding
of Kuranishi structure $\widehat{\mathcal U^+_{Z}}$ to
obtain
$$
\Phi^+_{Z,pq} \circ \Phi_q
: \mathcal U_q\vert_{(\varphi_{q})^{-1}(U^+_{Z,pq})}
\to \mathcal U^+_{Z,q}\vert_{U^+_{Z,pq}}
\to \mathcal U^+_{Z,p}.
$$
The composition gives the coordinate change $\Phi'_{pq}$ in this case.
\par
Note $\Phi^+_{Z,pq} \circ \Phi_q
=\Phi_p \circ  \Phi_{pq}$ on  $(\varphi_{q})^{-1}(U^+_{Z,pq})$,
by the definition of embedding of Kuranishi structures.
Using this fact, it is easy to see that they define a Kuranishi structure
on $X$.
\par
The Kuranishi structure $\widehat{\mathcal U'}$
has all the properties we need, except the property (3).
We will further modify $\widehat{\mathcal U'}$
to $\widehat{\mathcal U^+}$ by this reason.
Firstly, we use \cite[Propositions 6.44 and 6.49]{part11} to find
a Kuranishi structure $\widehat{\mathcal U''}$ such that
$$
\widehat{\mathcal U'} < \widehat{\mathcal U''}.
$$
Though there are various choices of such $\widehat{\mathcal U''}$,
we choose one of them in the proof of Lemma \ref{lem1544}.
For the later purpose, we will take more specific
$\widehat{\mathcal U''}$ in the proof of Lemma
\ref{lem1557}.
\par
Next, we modify $\widehat{\mathcal U''}$ to
$\widehat{\mathcal U^+}$ which has the required properties as follows.
We take an open set $W_2 \subset X$ such that
$$\overline W_2 \cap \overline\Omega_1 = \emptyset,
\quad
Z \cup W_2 = X.
$$
We replace various Kuranishi structures involved by their
open substructures
(but without changing the Kuranishi neighborhoods of the point $p\in K$)
and may assume the following.
\begin{enumerate}
\item[(I)]
If $\psi_p(U^+_{Z,p} \cap (s^+_{Z,p})^{-1}(0)) \cap \overline W_2 \ne \emptyset$,
then
$\psi^+_{Z,p}(U^+_{Z,p}  \cap (s^+_{Z,p})^{-1}(0)) \cap
\overline\Omega_1 = \emptyset$.
\item[(II)]
If
$\psi_p(U_p \cap s_p^{-1}(0)) \cap \overline W_2 \ne \emptyset$,
then
$\psi_p(U_p \cap s_p^{-1}(0)) \cap \overline\Omega_1 = \emptyset$.
\end{enumerate}
Now we define a Kuranishi structure $\widehat{\mathcal U^+}$ on $X$
as follows.
\begin{enumerate}
\item[(A)] If $p \in \overline W_2$, we put $\mathcal U^{+}_p
= \mathcal U''_p$.
\item[(B)]
If $p \notin \overline W_2$, we put $\mathcal U^+_p
= \mathcal U'_p\vert_{U'_p\setminus \psi^{-1}_p(\overline W_2)}$.
\end{enumerate}
The coordinate change is defined as follows.
Let $q \in \psi^+_p((s_p^+)^{-1}(0))$.
If $p,q \in \overline {W}_2$, then we define
$\Phi^+_{pq} = \Phi''_{pq}$.
If $p,q \notin \overline {W}_2$, then we define
$\Phi^+_{pq} = \Phi'_{pq}\vert_{U'_{pq} \setminus \psi^{-1}_{q}(\overline W_2)}$.
Among the other two cases
$q \in \overline W_2$,
$p \notin \overline W_2$ cannot occur.
Suppose
$p \in \overline W_2$,
$q \notin \overline W_2$.
Then there is an embedding $\Phi_q : \mathcal U'_q
\to \mathcal U''_q$.
The coordinate change of $\widehat{\mathcal U^+}$ is given by the
composition
$$
\Phi''_{pq} \circ \Phi_q
: \mathcal U'_q\vert_{(\varphi_{q})^{-1}(U''_{pq})}
\to \mathcal U''_q\vert_{U''_{pq}}
\to \mathcal U''_p.
$$
It is easy to see that this Kuranishi structure $\widehat{\mathcal U^+}$
has the required properties.
\end{proof}
We next discuss extension of CF-perturbations.

\begin{shitu}\label{Situ1556}
Suppose we are in the situation of Lemma \ref{lem1544}.
We assume the following in addition.
\begin{enumerate}
\item
We have a strongly smooth and weakly submersive map
$f : (X,\widehat{\mathcal U}) \to M$ to a manifold $M$.
(See \cite[Definition 3.38]{part11}.)
\item
We have a CF-perturbation $\widehat{\frak S^+_{Z}}$ of $\widehat{\mathcal U^+_{Z}}$.
\item
We have a strongly smooth map $f_Z : (Z,\widehat{\mathcal U_Z^+}) \to M$ which is strongly submersive with respect to
$\widehat{\frak S^+_{Z}}$.
(See \cite[Definition 9.2]{part11}.)
\item
The following diagram commutes.
\begin{equation}\nonumber
\xymatrix{
(X,\widehat{\mathcal U})\vert_{Z} \ar[rr]^{\rm embedding} \ar[dr]_{f\vert_Z} &  &
(Z,\widehat{\mathcal U_Z^+}) \ar[dl]^{f_Z} \\
& M  &
}
\end{equation}
$\blacksquare$
\end{enumerate}
\end{shitu}
\begin{lem}\label{lem1557}
In Situation \ref{Situ1556}, we may choose the Kuranishi structure $\widehat{\mathcal U^+}$
in Lemma \ref{lem1544} so that the following holds in addition.
\begin{enumerate}
\item
There exists a CF-perturbation
$\widehat{\frak S^+}$ of $\widehat{\mathcal U^+}$.
\item
The map $f$ extends to a strongly smooth map  $f^+ : (X,\widehat{\mathcal U^+}) \to M$.
\item
The extended map $f^+ : (X,\widehat{\mathcal U^+}) \to M$ is strongly submersive with respect to $\widehat{\frak S^+}$.
\item
For any $p \in K$, two CF-perturbations $\widehat{\frak S^+_{Z}}$ and
$\widehat{\frak S^+}$
assign the same CF-perturbation on the Kuranishi chart
$\mathcal U^+_{Z,p}
= \mathcal U^+_{p}$.
\item
$\widehat{\frak S^+_{Z}}\vert_{\Omega}$ is a restriction of  $\widehat{\frak S^+}\vert_{\Omega}$ to the
open substructure.
\item
When $\widehat{\frak S^+_{Z}}\vert_{\Omega}$
varies in a uniform family, we may take $\widehat{\frak S^+}$
to be uniform.
\end{enumerate}
\end{lem}
\begin{proof}
The lemma is a consequence of combination of results in
\cite[Sections 3,6,7 and 9]{part11}.
Before we start the proof,
we recall from \cite[Section 7]{part11}
that we used a good coordinate system
to define a CF-perturbation.
Thus we also need to involve an extension of
good coordinate system in the course of the proof of
Lemma \ref{lem1557}.
Indeed,
we use a good coordinate system from the given Kuranishi structure
to find an extension of
the given CF-perturbation, and come back from the good coordinate system to
a Kuranishi structure together with the CF-perturbation.
This process is described in \cite[Section 9]{part11}.
This is a rough description of the structure of the proof
of Lemma \ref{lem1557} given below.
\par
Now we start the proof.
In the proof of Lemma \ref{lem1544},
we took a relatively compact open subset $\Omega_1 \subset {\rm Int}\, Z$ such that
$$
\overline\Omega
\subset \Omega_1 \subset \overline\Omega_1 \subset {\rm Int}\, Z.
$$
As in the proof, we replace $\widehat{\mathcal U}$ and
$\widehat{\mathcal U_{Z}^+}$ by their open substructures
without changing the Kuranishi neighborhoods of the point $p\in K$
if necessary, and may assume that
the map $f_Z : (Z,\widehat{\mathcal U_Z^+}) \to M$ is strictly strongly submersive
(\cite[Definition 9.2]{part11})
with respect to $\widehat{\frak S^+_{Z}}$.
We put
$$
Z_1 = \overline{\Omega}_1, \quad
\widehat{\mathcal U^+_{Z_1}}
=
\widehat{\mathcal U^+_{Z}}\vert_{Z_1}.
$$
Then by the definition of the Kuranishi structure $\widehat{\mathcal U'}$ in
the proof of Lemma \ref{lem1544},
we note
\begin{equation}\label{eq1648}
\widehat{\mathcal U'}\vert_{Z_1}
=
\widehat{\mathcal U^+_{Z_1}}.
\end{equation}
\par
We apply \cite[Theorem 3.30]{part11} to the Kuranishi structure
$\widehat{\mathcal U^+_{Z_1}}$ to find
a good coordinate system $\widetriangle{\mathcal U^+_{Z_1}}$ on $Z_1$ and
a KG-embedding
\begin{equation}\label{KGembedd1649}
\Phi ~:~ \widehat{\mathcal U^+_{Z_1}} \longrightarrow
\widetriangle{\mathcal U^+_{Z_1}}.
\end{equation}
Since we are given a CF-perturbation $\widehat{\mathfrak S^+_{Z}}\vert_{Z_1}$
of $\widehat{\mathcal U^+_{Z_1}}$,
\cite[Lemma 9.10]{part11} shows that there exists a CF-perturbation
$\widetriangle{\mathfrak S^+_{Z_1}}$
of $\widetriangle{\mathcal U^+_{Z_1}}$ such that
$\widetriangle{\mathfrak S^+_{Z_1}}$ and
$\widehat{\mathfrak S^+_{Z}}\vert_{Z_1}$ are compatible with
the KG-embedding $\Phi$ in \eqref{KGembedd1649}.
In addition, \cite[Lemma 7.53 and Proposition 7.57]{part11}
yield that
there exists a strongly smooth map
$
\widetriangle{f_{Z_1}} : (Z_1, \widetriangle{\mathcal U^+_{Z_1}}) \to
M
$
such that it is a strongly submersive with respect to
$\widetriangle{\mathfrak S^+_{Z_1}}$ and
the following diagram commutes.
\begin{equation}\nonumber
\xymatrix{
(Z_1,\widehat{\mathcal U^+_{Z_1}}) \ar[rr]^{\Phi} \ar[dr]_{f_Z\vert_{Z_1}} &  &
(Z_1,\widetriangle{\mathcal U_{Z_1}^+}) \ar[dl]^{\widetriangle{f_{Z_1}}} \\
& M  &
}
\end{equation}
Since
$
\widehat{\mathcal U^+_{Z_1}}=\widehat{\mathcal U'}\vert_{Z_1}
$
is the restriction of the Kuranishi structure
$\widehat{\mathcal U'}$ on $X$ to $Z_1$, we can apply
\cite[Proposition 7.52]{part11} for the case $Z_1 =Z_1$ and $Z_2 =X$ to obtain
a good coordinate system $\widetriangle{\mathcal U'}$ on $X$ such that
it is an extension of $\widetriangle{\mathcal U^+_{Z_1}}$,
and
$\widetriangle{\mathcal U'}$ and $\widehat{\mathcal U'}$ are compatible
in the sense of \cite[Definition 3.32]{part11}, i.e.,
there exists a KG-embedding
$$
\widehat{\mathcal U'} \longrightarrow \widetriangle{\mathcal U'}.
$$
Moreover, the CF-perturbation
$\widetriangle{\mathfrak S^+_{Z_1}}$ is also extended to
a CF-perturbation $\widetriangle{\mathfrak S'}$ of
$\widetriangle{\mathcal U'}$
by \cite[Proposition 7.57]{part11}.
In addition, by \cite[Lemma 7.53]{part11} there exists a strongly smooth and
strongly submersive map
$$
\widetriangle{f'} ~:~
(X, \widetriangle{\mathcal U'}) \longrightarrow M
$$
with respect to $\widetriangle{\mathfrak S'}$
such that the following diagram commutes.
\begin{equation}\nonumber
\xymatrix{
(Z_1,\widehat{\mathcal U'}\vert_{Z_1}) \ar[rr]^{} \ar[dr]_{f_Z\vert_{Z_1}} &  &
(Z_1,\widetriangle{\mathcal U'}\vert_{Z_1}) \ar[dl]^{\widetriangle{f'}\vert_{Z_1}} \\
& M  &
}
\end{equation}
\par
Now we go back to Kuranishi structure from good coordinate system.
We apply \cite[Proposition 6.44]{part11} for the case $\widetriangle{\mathcal U}_0=\widetriangle{\mathcal U}=\widetriangle{\mathcal U'}$ to find a Kuranishi structure
$\widehat{\mathcal U''}$ together with a GK-embedding
$$
\widetriangle{\mathcal U'} \longrightarrow \widehat{\mathcal U''},
$$
and a strongly smooth map
$$
\widehat{f'} : (X,\widehat{\mathcal U''}) \longrightarrow
M
$$
such that $\widetriangle{f'}$ is a pull-back of
$\widehat{f'}$.
By \cite[Lemma 5.14]{part11}, $\widehat{\mathcal U''}$ is a thickening of
$\widehat{\mathcal U'}$:
$$
\widehat{\mathcal U'} < \widehat{\mathcal U''}.
$$
Moreover, \cite[Lemma 9.9]{part11} shows that there exists a CF-perturbation
$\widehat{\mathfrak S''}$ of $\widehat{\mathcal U''}$ such that
$\widehat{f'}$ is strongly submersive with respect to
$\widehat{\mathfrak S''}$.
\par
Finally in the exactly same way as in the proof of
Lemma \ref{lem1544}, we modify the Kuranishi structure $\widehat{\mathcal U''}$ obtained above to $\widehat{\mathcal U^+}$.
Accordingly, we also have the corresponding CF-perturbation
$\widehat{\mathfrak S^+}$ of $\widehat{\mathcal U^+}$ and
the corresponding map $f^+ : (X,\widehat{\mathcal U^+}) \to M$.
Then all the assertions of Lemma \ref{lem1557} follow from the construction.
\end{proof}

\subsection{Main results of Section \ref{sec:triboundary}}
\label{subsection:concltrisection}

We now combine the results of Subsections \ref{subsec:extenonechart}-\ref{subsec:extensionfromopen}
to prove results which we will use later in our applications.

\begin{prop}\label{prop1562}
In Situation \ref{sit1526}, there exists a $\tau'$-collared
Kuranishi structure $\widehat{\mathcal U^{++}}$ on $X$
for any $0< \tau'<\tau$ such that
Proposition \ref{prop528} holds,
by replacing $\widehat{\mathcal U^{+}}$ by
$\widehat{\mathcal U^{++}}$.
\end{prop}
\begin{rem}
The difference between Proposition \ref{prop528} and Proposition \ref{prop1562}
is that the  $\tau'$-collared Kuranishi structure $\widehat{\mathcal U^{++}}$
in Proposition \ref{prop1562} is defined on whole $X$,
while the  $\tau'$-collared Kuranishi structure $\widehat{\mathcal U^{+}}$
in Proposition \ref{prop528}
is defined only on a neighborhood of the boundary.
\end{rem}
\begin{proof}
We will use Lemma \ref{lem1544} to prove Proposition \ref{prop1562}.
For this purpose, we will slightly modify the $\tau'$-collared Kuranishi structure
$\widehat{\mathcal U^{+}}$
produced in Proposition \ref{prop528}
to get $\widehat{\mathcal U^+_{\Omega}}$
satisfying the assumption of Lemma \ref{lem1544}.
The detail is in order.
\par
We take $0< \tau' < \tau'' < \tau''' < \tau$.
Let $\widehat{\mathcal U^{+}}$ be the $\tau'$-collared
Kuranishi structure produced in Proposition \ref{prop528}.
We note that
$X' = X^{\boxplus\tau} = (X^{\boxplus(\tau-\tau')})^{\boxplus\tau'}$.
We take the Kuranishi structure $\widehat{\mathcal U''}$
on $X^{\boxplus(\tau-\tau')}$ such that
$(X^{\boxplus(\tau-\tau')},\widehat{\mathcal U''})^{\boxplus\tau'}
= (X^{\boxplus\tau},\widehat{\mathcal U^{+}})$.
(See Remark \ref{rem1537} for the description of the Kuranishi
structure $\widehat{\mathcal U''}$.)
We shrink the Kuranishi neighborhood $\mathcal U''_p$
of the Kuranishi structure $\widehat{\mathcal U''}$
to obtain $\mathcal U'''_p$ and
a Kuranishi structure $\widehat{\mathcal U'''}$
so that the following is satisfied.
\begin{enumerate}
\item
If $p \notin X^{\boxplus(\tau-\tau'')}$,
then
$
\psi''_p((s''_p)^{-1}(0)) \cap X^{\boxplus(\tau-\tau''')}
= \emptyset.
$
\item
If $p \in S_1(X^{\boxplus(\tau-\tau')})$,
then $\partial{\mathcal U'''_p} = \partial{\mathcal U''_p}$.
\end{enumerate}
Note that (2) above implies
\begin{equation}\label{formula1541}
(X,\widehat{\mathcal U'''})^{\boxplus(\tau-\tau')}\vert_
{\overline{X \setminus X^{\boxplus(\tau-\tau')}}}
= (X',\widehat{\mathcal U^+})\vert_
{\overline{X \setminus X^{\boxplus(\tau-\tau')}}}.
\end{equation}
We put $X'' = {\rm Int}\, (X \setminus X^{\boxplus(\tau-\tau')})$,
$K = \overline{X'' \setminus X^{\boxplus(\tau-\tau'')}}$
and
$Z = \overline{X'' \setminus X^{\boxplus(\tau-\tau''')}}$.
We now apply Lemma  \ref{lem1544}.
Here the role of $X$, $K$, $Z$, $\widehat{\mathcal U^+_{Z}}$
in Lemma \ref{lem1544} is played by
$X''$, $K$, $Z$ and $\widehat{\mathcal U'''}$, respectively.
\begin{rem}
As we mentioned at the top of Subsection \ref{subsec:extensionfromopen},
Lemma \ref{lem1544} is about genuine Kuranishi structure
and not $\tau$-collared Kuranishi structure.
So here we apply Lemma  \ref{lem1544} literally, not its
$\tau$-collared version.
\end{rem}
We thus obtain a Kuranishi structure
which we wrote $\widehat{\mathcal U^+}$
in Lemma  \ref{lem1544}.
We denote it here by $\widehat{\mathcal U^{\prime +}}$.
\par
We put $\widehat{\mathcal U^{++}}
= \widehat{\mathcal U^{\prime + \boxplus\tau'}}$.
(\ref{formula1541}), Lemma  \ref{lem1544}
(1) and the fact
that  $\widehat{\mathcal U^+}$ satisfies Proposition \ref{prop528}
imply that $\widehat{\mathcal U^{++}}$
satisfies Proposition \ref{prop528}.
The proof of Proposition \ref{prop1562} is now complete.
\end{proof}
We next include CF-perturbations.
\begin{prop}\label{prop529rev}
In Situation \ref{situ529} there exists a
$\tau'$-collared CF-perturbation $\widehat{\frak S^{++}}$ on the Kuranishi structure
$\widehat{\mathcal U^{++}}$ obtained in Proposition \ref{prop1562}
such that
\begin{enumerate}
\item
Its restriction to $(\widehat S_k(X),\widehat{\mathcal U_{S_k}^+})$
is equivalent to $\widehat{\frak S^+_{k}}$.
\item
If $f : (X,\widehat{\mathcal U}) \to M$ is weakly submersive
and $f\vert_{\widehat{\mathcal U^+_{k}}}$
is strongly submersive with respect to $\widehat{\frak S^+_{k}}$, then
we may take $\widehat{\frak S^{++}}$ such that
$f : (X,\widehat{\mathcal U^{++}}) \to M$ is strongly submersive
with respect to  $\widehat{\frak S^{++}}$.
\item
Uniformity of CF-perturbations is
preserved in this construction.
\end{enumerate}
\end{prop}
\begin{proof}
We use the notation in the proof of Proposition \ref{prop1562}.
We apply Proposition \ref{prop529} to obtain a CF-perturbation
$\widehat{\frak S^{+}}$ on $\widehat{\mathcal U^{+}}$.
Since $\widehat{\frak S^{+}}$ is $\tau'$-collared,
it is induced by a CF-perturbation $\widehat{\frak S''}$ of
$\widehat{\mathcal U''}$.
Therefore by restriction we obtain a CF-perturbation
$\widehat{\frak S^+_{Z}}$ on
$\widehat{\mathcal U^+_{Z}}$.
Here $Z$ is one in the proof of
Proposition \ref{prop1562}.
Thus we can apply Lemma \ref{lem1557}
to obtain required $\widehat{\frak S^{++}}$ and $\widehat{\mathcal U^{++}}$.
\end{proof}
\par\medskip
One minor point remains to be explained to apply the results of this subsection.
Note that a $\tau$-collared Kuranishi structure is not a Kuranishi structure,
since some points are not assigned its Kuranishi neighborhood to.
Since in Part 1 the `push out' is defined for the case of Kuranishi structure and
good coordinate system but not for the case of $\tau$-collared Kuranishi
structure,
we need some explanation to define the notion of the `push out' for
$\tau$-collared Kuranishi structure.
\par
One can define the notion of $\tau$-collared good coordinate system,
and prove the existence of  $\tau$-collared good coordinate system
compatible with each $\tau$-collared Kuranishi structure, and use it
to define the `push out'.
(We note that $\tau$-collared good coordinate system
is a special case of good coordinate system.
See the end of Remark \ref{rem1528}.)
It is certainly possible to proceed in that way.
\par
Here we take a slightly different way which seems shorter as follows.
\begin{lem}\label{lem1569}
Let $\widehat{\mathcal U}$ be a $\tau$-collared Kuranishi structure
on $X' = X^{\boxplus\tau}$.
\begin{enumerate}
\item We can associate a Kuranishi structure $\widehat{\mathcal U'}$
on $X'$ such that:
\begin{enumerate}
\item If $p \in \overset{\circ\circ}S_k(X')$ then
$\mathcal U'_p = \mathcal U_p$.
\item
For each $p \in X'$  we take the (unique) point $\widehat p \in \overset{\circ\circ}S_k(X')$
such that $\mathcal R(p) = \mathcal R(\widehat p)$ and
$\mathcal R(p) \in \overset{\circ\circ}S_k(X)$.
(See Lemma \ref{lem538}.)
Then ${\mathcal U'}_p$ is an open subchart of ${\mathcal U}_{\widehat p}$
\end{enumerate}
\item
If $\widehat{\frak S}$ is a $\tau$-collared CF-perturbation on
$\widehat{\mathcal U}$,
it induces a CF-perturbation $\widehat{\frak S'}$
on $\widehat{\mathcal U'}$.
\item
If $\widehat{\mathcal U'}$ and $\widehat{\mathcal U''}$ are
two Kuranishi structures satisfying (1) (a)(b) above,
then there exists a Kuranishi structure
$\widehat{\mathcal U'''}$ such that it satisfies (1) (a)(b) above
and it is an open substructure of both $\widehat{\mathcal U'}$ and $\widehat{\mathcal U''}$.
Similar uniqueness statement holds for the CF-perturbation in (2).
\item
When a $\tau$-collared CF-perturbation $\widehat{\frak S}$ on
$\widehat{\mathcal U}$ varies in a uniform family,
we may take the induced CF-perturbation $\widehat{\frak S'}$ in (2) to be uniform.
\end{enumerate}
\end{lem}
The proof is obvious so omitted.

\begin{defn}\label{pushoutdeftau}
Let $f : (X',\widehat{\mathcal U}) \to M$ be a strongly smooth map
which is strongly submersive with respect to $\widehat{\frak S}$.
Let $h$ be a differential form on $(X',\widehat{\mathcal U})$.
We take a Kuranishi structure $\widehat{\mathcal U}'$ as in Lemma \ref{lem1569} and define its
{\it push out}
\index{integration along the fiber (push out) ! integration along the fiber}
\index{push out ! {\it see: integration along the fiber}}
\begin{equation}\label{pushtau}
f_!(h;\widehat{\frak S}^{\epsilon}) = f_!(h;\widehat{\frak S^{\prime\epsilon}}).
\end{equation}
Here $h$ in the right hand side is the differential form on
$(X',\widehat{\mathcal U'})$ which is
induced from $h$ by Lemma \ref{lem1569} (1).
\end{defn}
\begin{lem}\label{lem17535352}
The right hand side of (\ref{pushtau}) is independent of the choices of $\widehat{\mathcal U}'$, $\widehat{\frak S'}$.
Moreover Stokes' formula (\cite[Theorem 9.26]{part11}) and
the composition formula (\cite[Theorem 10.20]{part11})
hold for the push out in Definition \ref{pushoutdeftau}.
\end{lem}
\begin{proof}
Independence follows from \cite[Theorem 9.14]{part11}.
Stokes' formula and the composition formula
are direct consequences of the corresponding results
(\cite[Theorems 9.26 and 10.20]{part11}).
\end{proof}

\section{Smoothing corners and composition of morphisms}
\label{section:compomorphis}
The goal of this section is to define composition of morphisms
of linear K-systems (Lemma-Definition \ref{1638defken}) and to show that it is associative (Proposition \ref{prop1640}).
There are two key ingredients for the construction of composition of morphisms.
One is `partial trivialization of corners' and the other is
`smoothing corners'.
The precise definition of `partial trivialization of corners' will be given
in Definition \ref{defn1531revrev}
and `smoothing corners' will be described in detail
in Subsections \ref{subsec:smoothcornermodel} and \ref{subsec:smoothcornerkstr}.
After that, we will define composition of morphisms in Subsection \ref{subsec:complinkurasmcorner}
and prove its associativity in Subsection \ref{subsec:compassoci}.
We also discuss identity morphism in Subsection \ref{subsec:identitylinsys}.

\subsection{Introduction to Section \ref{section:compomorphis}}
\label{subsec:introcomp}
In this subsection, we explain the idea of the construction of composition of morphisms and its geometric background.
\par
Firstly, we explain the reason why we need
the notion of `partial trivialization of corners',
or more generally, `partially trivialized fiber products', instead of
trivialization of corners or usual fiber products.
Let $\mathcal N_{12}(\alpha_1,\alpha_2)$ and $\mathcal N_{23}(\alpha_2,\alpha_3)$
be interpolation spaces of morphisms $\frak N_{12}$ and $\frak N_{23}$
of linear K-systems, respectively.
As we mentioned in Lemma-Definition \ref{lemdef1434},
the interpolation space of the morphism $\frak N_{32}\circ\frak N_{21}$  is a union of
fiber products
\begin{equation}\label{unionoffiber}
\bigcup_{\alpha_2} \mathcal N_{12}(\alpha_1,\alpha_2)
\times_{R_{\alpha_2}} \mathcal N_{23}(\alpha_2,\alpha_3).
\end{equation}
The union in (\ref{unionoffiber}) may not be a disjoint union, in general.
In fact, the summands corresponding to $\alpha_2$ and to $\alpha'_2$
may intersect at
\begin{equation}\label{unionoffiber2}
\mathcal N_{12}(\alpha_1,\alpha_2)
\times_{R_{\alpha_2}}
\mathcal M^{2}(\alpha_2,\alpha'_2)
\times_{R_{\alpha'_2}}
\mathcal N_{23}(\alpha'_2,\alpha_3).
\end{equation}
Moreover three of such summands may intersect at
\begin{equation}\label{unionoffiber3}
\mathcal N_{12}(\alpha_1,\alpha_2)
\times_{R_{\alpha_2}}
\mathcal M^{2}(\alpha_2,\alpha'_2)
\times_{R_{\alpha'_2}}
\mathcal M^{2}(\alpha'_2,\alpha''_2)
\times_{R_{\alpha''_2}}
\mathcal N_{23}(\alpha''_2,\alpha_3).
\end{equation}
The pattern how the summands of (\ref{unionoffiber})
intersect each other is similar to the way how the
components of the boundary of certain K-space (or of
orbifold) intersect each other. (Namely, each of the summands
of (\ref{unionoffiber})  corresponds to a codimension one boundary
and (\ref{unionoffiber2}), (\ref{unionoffiber3}) correspond to codimension 2 and 3 corners, respectively.)
\par
We can use this observation to apply a version of Proposition \ref{prop528},
that is, we
 `put the collar'  {\it outside} the union  (\ref{unionoffiber})
to obtain a collared K-space so that its boundary is
$$
\bigcup_{\alpha_2} \mathcal N_{12}(\alpha_1,\alpha_2)
\times_{R_{\alpha_2}}^{\boxplus\tau} \mathcal N_{23}(\alpha_2,\alpha_3).
$$
Here the notation $\times_{R_{\alpha_2}}^{\boxplus\tau}$
is defined in Definition \ref{defn1635}, which
is called
{\it partially trivialized fiber product}.
\index{partially trivialized fiber product}
\index{fiber product ! partially trivialized fiber product}
\index{corner ! partially trivialized fiber product}
The notion of {\it partial trivialization of corners}
\index{trivialization of corners ! partial}
\index{corner ! partial trivialization of corners}
is
introduced in Subsection \ref{subsec:prtialtrivi}.
See Definition \ref{defn1531revrev}.
The reason why we introduce the notion of partial trivialization of corner
is as follows. Note that
the boundary or corner of K-space of the summand of
(\ref{unionoffiber}) has different components from those appearing in
(\ref{unionoffiber2}), (\ref{unionoffiber3}).
In fact, a boundary of the form
$$
\mathcal M^{1}(\alpha_1,\alpha'_1)
\times_{R_{\alpha'_1}}
\mathcal N_{12}(\alpha'_1,\alpha_2)
\times_{R_{\alpha_2}} \mathcal N_{23}(\alpha_2,\alpha_3)
$$
also appears.
To study this kind of boundary components together,
we introduce the notion of partial
trivialization of corner and use it to modify
Proposition \ref{prop528} so that it
can be directly applicable to our situation.
\par
After we have done partial trivialization of corners, we will next discuss
{\it smoothing corners}.
Here we use the fact that a K-space $X$ has a collar where all the objects are
`constant' in the direction of the collar.
See Subsection \ref{subsec:subsec32-1} for the reason why
this fact is useful for our construction.
Then we can finally define composition of morphisms.
\par\medskip
In the rest of this subsection, we explain a geometric origin
of the idea that the union (\ref{unionoffiber}) looks like a
boundary of certain K-space, by considering the situation of
the Floer cohomology of periodic Hamiltonian system.
(The story of linear K-system will be applied to
define and study the Floer cohomology of periodic Hamiltonian system.
Namely, we associate
such a system to each periodic Hamiltonian
function $H$.
See Subsection
\ref{subsecintolinear}.)
Let $H^1$, $H^2$, $H^3$ be periodic Hamiltonian functions.
To define a cochain map between Floer's cochain complexes
associated to them, we use $\tau \in \R$ dependent Hamiltonian
functions interpolating them.
Namely, we take $H^{ij} : \R \times S^1 \times M \to \R$
such that
\begin{equation}\label{def:Hij}
H^{ij}(\tau,t,x)
=
\begin{cases}
H^j(t,x)
&\text{if $\tau \le -T_0$}, \\
H^i(t,x)
&\text{if $\tau \ge T_0$}.
\end{cases}
\end{equation}
We then use the moduli space of the solutions of the equation
\begin{equation}\label{equa165}
\frac{\partial u}{\partial \tau}
+
J
\left(
\frac{\partial u}{\partial t} - X_{H^{ij}_{\tau,t}}(u)
\right) = 0
\end{equation}
where $H^{ij}_{\tau,t}(x) = H^{ij}(\tau,t,x)$ and $X_{H^{ij}_{\tau,t}}$ is its
Hamiltonian vector field.
The solution space of (\ref{equa165}) with an
appropriate boundary condition becomes an interpolation
space $\mathcal N^{ji}(\alpha_i,\alpha_j)$ of the morphism
$\frak N^{ji}$.
\par
To study the relation between $\frak N^{31}$ and the composition
$\frak N^{32} \circ \frak N^{21}$ we use one parameter family of
$\tau$-dependent Hamiltonian
functions $H^{31,T}$ where
\begin{equation}\label{eq1666}
H^{31,T}(\tau,t,x)
=
\begin{cases}
H^1(t,x)
&\text{if $\tau \le -T_0-T$}\\
H^{21}(\tau+T,t,x)
&\text{if $-T_0-T\le \tau \le T_0 - T$}\\
H^2(t,x)
&\text{if $T_0 - T \le \tau \le T-T_0$}\\
H^{32}(\tau-T,t,x)
&\text{if $T-T_0\le \tau \le T+T_0 $}\\
H^3(t,x)
&\text{if $T+T_0 \le \tau$}.
\end{cases}
\end{equation}
\begin{figure}[h]
\centering
\includegraphics[scale=0.5]{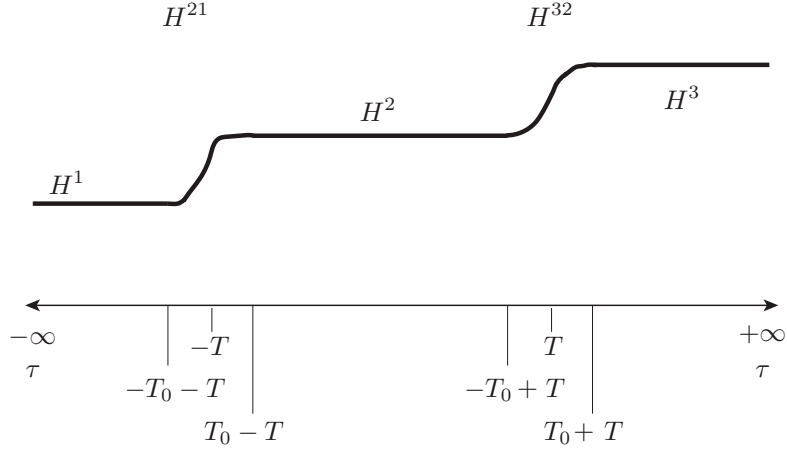}
\caption{the concatenation $H^{31}$}
\label{Figure17-1}
\end{figure}
See Figure \ref{Figure17-1}.
\par
We may choose $H^{31} = H^{31,2T_0}$, for example.
We consider the limit as $T\to +\infty$
and the set of solutions of (\ref{equa165})
with $X_{H^{ij}_{\tau,t}}$ replaced by $X_{H^{31,T}_{\tau,t}}$.
The moduli space of its solutions becomes the union of fiber products
(\ref{unionoffiber}).
Thus the union of solution spaces for $T \in [2T_0,\infty)$ and the
fiber product (\ref{unionoffiber}) gives a homotopy between
the morphism defined by $H^{31} = H^{31,2T_0}$ and
the composition whose interpolation space is (\ref{unionoffiber}).
The space (\ref{unionoffiber}) itself is a part of the boundary of this
cobordism.
\subsection{Partial trivialization of cornered K-space}
\label{subsec:prtialtrivi}

In this subsection we explain that the story of trivialization of
corner in Section \ref{sec:triboundary} can be generalized
to the case of partial trivialization in a quite straightforward way.
Because of the nature of this book,
we repeat the statements.
We believe that
the readers can go through those parts very
quickly since there is nothing new to do.

\begin{shitu}\label{decomporbbdr}
Let $U$ be an admissible orbifold with corners.
We decompose its normalized boundary $\partial U$
into a {\it disjoint} union
$$
\partial U = \partial^0U  \cup \partial^1U
$$
where both $\partial^0U$  and $\partial^1U$ are open subsets in $\partial U$.
We denote this decomposition by $\frak C$.
We also denote $\partial^0U$  by
$\partial_{\frak C}U$. $\blacksquare$
\end{shitu}

\begin{defn}\label{defnSCk}
In Situation \ref {decomporbbdr} we define a closed subset
$S^{\frak C}_k(U)$ of $U$ as follows.
Let $p \in U$. We take its orbifold chart $(V_p,\Gamma_p,\phi_p)$
where $V_p \subset \overline V_p \times [0,1)^{k'}$ and
$p = \phi(y_0,(0,\dots,0))$.
For $i=1,2,\dots,k'$ we put
$$
\partial_i V_p = V_p \cap (\overline V_p \times [0,1)^{i-1} \times \{0\}
\times [0,1)^{k'-i}).
$$
We require
\begin{equation}\label{eq1677}
\#
\{
i \in \{1,\dots,k'\} \mid \phi(\partial_i V_p) \subset \partial^0 U
\}
= k.
\end{equation}
Then
$\overset{\circ}{S^{\frak C}_{k}}(U)$ is the set of all $p \in U$ such that
(\ref{eq1677}) is satisfied.
We put
$$
{S^{\frak C}_k}(U)
=
\bigcup_{\ell \ge k} \overset{\circ}{S^{\frak C}_{\ell}}(U).
$$
\end{defn}
\begin{conven}\label{conv1633}
In case $p \in S^{\frak C}_k(U)$ as above, we take its orbifold chart
$(V_p,\Gamma_p,\phi_p)$ as above such that
$
 \phi(\partial_i V_p) \subset \partial^0U
$
if and only if $i = k'-k+1,\dots,k'$.
\end{conven}
\par
\begin{shitu}\label{decomporbbdrkura}
Let $(X,\widehat{\mathcal U})$
be an $n$-dimensional K-space.
We assume that for each $p \in \partial X$, we have a decomposition of
its Kuranishi neighborhood into a disjoint union
\begin{equation}\label{eq16888}
\partial U_p = \partial^0U_p  \cup \partial^1U_p
\end{equation}
such that for each coordinate change $\Phi_{pq}$ with
$p,q \in \partial X$ we have
\begin{equation}\label{eq16888tugi}
\varphi_{pq}(U_{pq} \cap \partial^0U_p) \subset \partial^0U_p,
\qquad
\varphi_{pq}(U_{pq} \cap \partial^1U_p) \subset \partial^1U_p.
\end{equation}
We also denote this decomposition by $\frak C$.$\blacksquare$
\end{shitu}
\begin{defn}
In Situation \ref{decomporbbdrkura} we define
$S^{\frak C}_k(X,\widehat{\mathcal U})$
as follows. If $p \notin \partial X$ then
$p \in S^{\frak C}_0(X,\widehat{\mathcal U})$ but $p \notin S^{\frak C}_k(X,\widehat{\mathcal U})$ for $k\ge 1$.
If $p \in \partial X$ then $p \in S^{\frak C}_k(X,\widehat{\mathcal U})$
if and only if $o_p \in S^{\frak C}_k(U_p)$.
We denote $\partial_{\frak C}X = S_1^{\frak C}(X,\widehat{\mathcal U})$.
We also put
$$
\overset{\circ}{S^{\frak C}_k}(X,\widehat{\mathcal U})
=
{S^{\frak C}_k}(X,\widehat{\mathcal U})
\setminus \bigcup_{\ell > k} {S^{\frak C}_{\ell}}(X,\widehat{\mathcal U}).
$$
We can define a similar notion for good coordinate system
by modifying the above definition in an obvious way.
\footnote{We do not define it in detail since it is never used in this book.}
\end{defn}
We can generalize Proposition \ref{prop2813}
without change as follows.
\begin{prop}\label{prop2813rev}
In Situation \ref {decomporbbdrkura},
for each $k$ there exist a compact
$(n-k)$-dimensional K-space $\widehat S^{\frak C}_k(X,\widehat{\mathcal U})$ with
corners, a map
$\pi_k :  \widehat S^{\frak C}_k(X,\widehat{\mathcal U})  \to
S^{\frak C}_k(X,\widehat{\mathcal U})$, and
a decomposition  of
$\partial\widehat S^{\frak C}_k(X,\widehat{\mathcal U})$ as in Situation \ref {decomporbbdrkura},
(which we also denote by $\frak C$),
and a map $\pi_{\ell,k} : \widehat S^{\frak C}_{\ell}(\widehat S^{\frak C}_k(X,\widehat{\mathcal U}))
\to \widehat S^{\frak C}_{k+\ell}(X,\widehat{\mathcal U})$
for each $\ell, k$
such that they enjoy the following properties:
\begin{enumerate}
\item
The map $\pi_k$ is a continuous map of underlying topological space.
\item The interior of $\widehat S^{\frak C}_k(X,\widehat{\mathcal U})$
is isomorphic to $\overset{\circ}{S^{\frak C}_k}(X,\widehat{\mathcal U})$.
The underlying homeomorphism of this isomorphism is the
restriction of $\pi_k$.
\item
The map $\pi_{\ell,k}$ is an $(\ell+k)!/\ell!k!$ fold covering map of K-spaces.
\item
The following objects on $(X,\widehat{\mathcal U})$ induce those on $\widehat S^{\frak C}_k(X,\widehat{\mathcal U})$.
Moreover the induced objects are compatible with the covering maps $\pi_{\ell,k}$.
\begin{enumerate}
\item
CF-perturbation.
\item
Multivalued perturbation.
\item
Differential form.
\item Strongly continuous map. Strongly smooth map.
\item Covering map.
\end{enumerate}
\item
The following diagram commutes.
\begin{equation}\label{diagin26277XXrev2}
\begin{CD}
\widehat S^{\frak C}_{k_1}(\widehat S^{\frak C}_{k_2}(\widehat S^{\frak C}_{k_3}(X,\widehat{\mathcal U}))) @ >{\pi_{k_1,k_2}}>>
\widehat S^{\frak C}_{k_1+k_2}(\widehat S^{\frak C}_{k_3}(X,\widehat{\mathcal U})) \\
@ V{\widehat S^{\frak C}_{k_1}(\pi_{k_2,k_3})}VV @ VV{\pi_{k_1+k_2,k_3}}V\\
\widehat S^{\frak C}_{k_1}(\widehat S^{\frak C}_{k_2+k_3}(X,\widehat{\mathcal U})) @ > {\pi_{k_1,k_2+k_3}} >> \widehat S^{\frak C}_{k_1+k_2+k_3}(X,\widehat{\mathcal U})
\end{CD}
\end{equation}
Here $\widehat S^{\frak C}_{k_1}(\pi_{k_2,k_3})$ is the covering map induced from $\pi_{k_2,k_3}$.
\item
For $i=1,2$ let $f_i : (X_i,\widehat{\mathcal U_i}) \to M$ be a strongly smooth map.
If $f_1$ is transversal to $f_2$, then
$$
\widehat S^{\frak C}_{k}\left((X_1,\widehat{\mathcal U_1}) \times_M (X_2,\widehat{\mathcal U_2})\right)
\cong
\coprod_{k_1+k_2=k}
\widehat S^{\frak C}_{k_1}(X_1,\widehat{\mathcal U_1}) \times_M \widehat S^{\frak C}_{k_2}(X_2,\widehat{\mathcal U_2}).
$$
Here the right hand side is the disjoint union.
The decomposition of the fiber product
$(X_1,\widehat{\mathcal U_1}) \times_M (X_2,\widehat{\mathcal U_2})$
as in Situation \ref {decomporbbdrkura}
is induced from those of $(X_1,\widehat{\mathcal U_1})$ and
$(X_2,\widehat{\mathcal U_2})$ as follows.
$$
\aligned
& \partial^0((X_1,\widehat{\mathcal U_1}) \times_M (X_2,\widehat{\mathcal U_2}))\\
& =
\partial^0(X_1,\widehat{\mathcal U_1}) \times_M (X_2,\widehat{\mathcal U_2})
\cup
(-1)^{\dim X_1 + \dim M}
(X_1,\widehat{\mathcal U_1}) \times_M \partial^0(X_2,\widehat{\mathcal U_2}).
\endaligned
$$
\item
(1)-(6) also hold when we replace `Kuranishi structure' by `good coordinate system'.
\item
Various kinds of embeddings of Kuranishi structures and/or good coordinate systems induce ones of
$\widehat S^{\frak C}_{k}(X,\widehat{\mathcal U})$.
\end{enumerate}
\end{prop}
The proof is the same as the proof of Proposition \ref{prop2813}
and so is omitted.
\par\medskip
Next we generalize the process of trivializing corner in Section
\ref{sec:triboundary}
to partial trivialization.

\begin{shitu}\label{situ16777}
Suppose that we are in Situation \ref{situ151} and
a decomposition
$
\partial U = \partial^0U  \cup \partial^1U
$
is given as in Situation \ref{decomporbbdr}. $\blacksquare$
\end{shitu}
In Situation \ref{situ16777}
We define a map
$$
\mathcal R^{\frak C}_x : \overline{V_x} \times [0,1)^{k'-k}
\times (-\infty,1)^k \to
\overline{V_x} \times [0,1)^{k'}
$$
by
$$
\mathcal R^{\frak C}_x(\overline y,(t_1,\dots,t_k))
= (\overline y,(t'_1,\dots,t'_k))
$$
where
$$
t'_i =
\begin{cases}
t_i  &\text{if $t_i \ge 0$,}\\
0    &\text{if $t_i \le 0$.}
\end{cases}
$$
For $\tau \ge 0$, we define an open subset $V_x^{\frak C\boxplus\tau}$ of
$\overline{V_x} \times  [0,1)^{k'-k}
\times [-\tau,1)^k$
to be
$$
V_x^{\frak C\boxplus\tau} =
(\mathcal R_x^{\frak C})^{-1}(V_x) \cap (\overline{V_x}
\times [0,1)^{k'-k} \times [-\tau,1)^k).
$$
Then $\mathcal R^{\frak C}_x$ induces a map
$\mathcal R^{\frak C}_x : V_x^{\frak C\boxplus \tau} \to V_x \subset
\overline{V_x} \times [0,1)^{k'}$.
\par
We can define a $\Gamma_x$ action on  $V_x^{\frak C\boxplus\tau}$
in the same way as in Definition \ref{defn154444}
and put $U_x^{\frak C\boxplus\tau}=V_x^{\frak C\boxplus\tau}/\Gamma_x$.
We define $\mathcal E^{\frak C\boxplus\tau}_x$ in the same way as in
(\ref{form151111})  by
$$
\mathcal E^{\frak C\boxplus\tau}_x =
(\mathcal R_x^{\frak C})^*(\mathcal E_x)
= (E_x \times V^{\frak C\boxplus\tau}_x)/\Gamma_x.
$$
The section $s_x$ of $\mathcal E_x$ induces a section $s_x^{\frak C\boxplus\tau}$ of
$\mathcal E^{\frak C\boxplus\tau}_x$ in an obvious way.
We define
$$
(X \cap V_x)^{\frak C\boxplus\tau} = (s_x^{\frak C\boxplus\tau})^{-1}(0)/\Gamma_x.
$$
Let
$\psi_x^{\frak C\boxplus\tau} : (s_x^{\frak C\boxplus\tau})^{-1}(0)/\Gamma_x \to (X \cap V_x)^{\frak C\boxplus\tau}$ be the identity map.
Then similarly to Lemma-Definition \ref{lem15444}, we find that
$$
\mathcal U^{\frak C\boxplus\tau} =
(V^{\frak C\boxplus\tau}_x/\Gamma_x,\mathcal E^{\frak C\boxplus\tau}_x,
\psi_x^{\frak C\boxplus\tau}, s_x^{\frak C\boxplus\tau})
$$
is a Kuranishi chart of $(X \cap V_x)^{\frak C\boxplus\tau}$.
Moreover
the following objects on $\mathcal U = (U,\mathcal E,s,\psi)$
induce the corresponding objects on $\mathcal U^{\frak C\boxplus\tau}$.
The proof is the same as that of Lemma-Definition \ref{lemdef157} so is omitted.
\par
\begin{enumerate}
\item[$\bullet$]
CF-perturbation.
\item[$\bullet$]
Strongly smooth map.
\item[$\bullet$]
Differential form.
\item[$\bullet$]
Multivalued perturbation.
\end{enumerate}
\par
We put
$
\overset{\circ\circ}S_k(V_x^{\frak C\boxplus\tau})
=
S_k(V_x^{\frak C\boxplus\tau})
\cap (\mathcal R_x^{\frak C})^{-1}(\overset{\circ}S_k(V_x^{\frak C\boxplus\tau}))
$.
Then Lemma \ref{lem1599} is generalized in an obvious way.
Furthermore we have
\begin{lem}
Suppose we are in Situation \ref{situ156} and
$
\partial U_i = \partial^0U_i  \cup \partial^1U_i
$
for $i=1,2$. We denote by $\frak C$ the decomposition of the boundary.
Then
$\Phi_{21} = (\varphi_{21},\widehat\varphi_{21})$
induces an embedding
$\Phi^{\frak C\boxplus\tau}_{21} : \mathcal U^{\frak C\boxplus\tau}_1
\to \mathcal U^{\frak C\boxplus\tau}_2$ of Kuranishi charts
whose restriction to $\mathcal U_1$
coincides with $\Phi_{21}$.
Moreover we have
\begin{enumerate}
\item
In case of  Situation \ref{situ156} (1), $\frak S^{1\frak C\boxplus\tau}$, $\frak S^{2\frak C\boxplus \tau}$ are compatible
with $\Phi_{21}^{\frak C\boxplus\tau}$.
\item
In case of  Situation \ref{situ156} (2),
$(\varphi_{21}^{\frak C\boxplus\tau})^*(h^{\frak C\boxplus\tau}_2) = (h^{\frak C\boxplus\tau}_1)$.
\item
In case of Situation \ref{situ156} (3),
$f_2^{\frak C\boxplus\tau} \circ \varphi_{21}^{\frak C\boxplus\tau} = f_1^{\frak C\boxplus\tau}$.
\item
In case of Situation \ref{situ156} (4),
$\frak s^{1\frak C\boxplus \tau}$, $\frak s^{2\frak C\boxplus\tau}$ are compatible
with $\Phi_{21}^{\frak C\boxplus\tau}$.
\item
In case of Situation \ref{situ156} (5),
we have
$\Phi^{\boxplus\tau}_{\frak C31} =  \Phi^{\frak C\boxplus\tau}_{32}\circ \Phi^{\frak C\boxplus\tau}_{21}$.
\item
$\varphi^{\frak C\boxplus\tau}_{21} \circ \mathcal R^{\frak C}_1 = \mathcal R^{\frak C}_2 \circ \varphi^{\frak C\boxplus\tau}_{21}$.
\end{enumerate}
\end{lem}
The proof is the same as the proof of Lemma \ref{lem151888}.

\begin{defn}\label{defXenhanceref}
Suppose we are in Situation \ref{decomporbbdrkura}.
Consider a disjoint union
$$
\coprod_{p\in X}  (s_p^{\frak C\boxplus\tau})^{-1}(0)/\Gamma_{p}
$$
and define an equivalence relation
$\sim$ on it as follows:
Let $x_p \in (s^{\frak C\boxplus\tau}_{p})^{-1}(0)$ and
$x_q \in (s^{\frak C\boxplus\tau}_{q})^{-1}(0)$.
Then we define
$[x_p] \sim [x_q]$ if there exist
$r \in X$ and $x_r \in (s^{\frak C\boxplus\tau}_p)^{-1}(0)
\cap U^{\frak C\boxplus\tau}_{pr} \cap U^{\frak C\boxplus\tau}_{qr}$
such that
\begin{equation}\label{equiXXcondplusrevrev}
[x_p] = \varphi_{pr}^{\frak C\boxplus\tau}([x_r]),
\quad
[x_q] = \varphi_{qr}^{\frak C\boxplus\tau}([x_r]).
\end{equation}
The same argument as in Lemma \ref{lem1521} show that $\sim$
is an equivalence relation.
Then we define a topological space $X^{\frak C\boxplus\tau}$
by the set of the equivalence classes of this
equivalence relation $\sim$ with the quotient topology:
$$
X^{\frak C\boxplus\tau} :=
\left(\coprod_{p\in X}  (s_p^{\frak C\boxplus\tau})^{-1}(0)/\Gamma_{p}\right) /\sim .
$$
\end{defn}
In the situation of Definition \ref{defXenhanceref}
we put \begin{equation}
\overset{\circ\circ}S_k(X^{\frak C\boxplus\tau})
= S_k(X^{\frak C\boxplus\tau},\widehat{\mathcal U}^{\frak C\boxplus\tau})
\cap (\mathcal R^{\frak C})^{-1}
(\overset{\circ}S_k(X,\widehat{\mathcal U})).
\end{equation}
We define $B_{\tau}(\overset{\circ\circ}S_k(X^{\frak C\boxplus\tau})) \subset X^{\frak C\boxplus\tau}$
as the union of
\begin{equation}
\psi^{\frak C\boxplus\tau}_p
\left(
(s^{\frak C\boxplus\tau}_p)^{-1}(0) \cap
\{(\overline y,(t_1,\dots,t_k))
\mid t_i \le 0, \,\, i=k'-k+1,\dots,k'\}
\right).
\end{equation}
We can show
\begin{equation}
X^{\frak C\boxplus\tau}= \coprod_k B_{\tau}(\overset{\circ\circ}S_k(X^{\frak C\boxplus\tau}))
\end{equation}
in the same way as in Lemma \ref{lem1530}.

\begin{defn}\label{defn1531revrev}
\begin{enumerate}
\item
Let $p' \in \overset{\circ\circ}S_k(X^{\frak C\boxplus\tau})$.
A {\it $\tau$-$\frak C$-collard Kuranishi neighborhood} at $p'$
\index{Kuranishi chart ! $\tau$-$\frak C$-collard Kuranishi neighborhood}
is a Kuranishi chart $\mathcal U_{p'}$ of $X^{\frak C\boxplus\tau}$ which is
$(\mathcal U_{p})^{\frak C\boxplus\tau}$
for certain Kuranishi neighborhood $\mathcal U_{p}$ of
$p =\mathcal R^{\frak C}(p')$.
\item
Let $p' \in \overset{\circ\circ}S_k(X^{\frak C\boxplus\tau})$,
$q' \in \overset{\circ\circ}S_{\ell}(X^{\frak C\boxplus\tau})$
and let $\mathcal U_{p'} = (\mathcal U_{p})^{\frak C\boxplus\tau}$, $\mathcal U_{q'}
= (\mathcal U_{q})^{\frak C\boxplus\tau}$
be their $\tau$-$\frak C$-collared Kuranishi neighborhoods, respectively.
Suppose
$q' \in \psi_{p'}(s_{p'}^{-1}(0))$.
A {\it $\tau$-$\frak C$-collared coordinate change} $\Phi_{p'q'}$
\index{Kuranishi chart ! $\tau$-$\frak C$-collared coordinate change}
from $\mathcal U_{q'}$
to $\mathcal U_{p'}$ is by definition
$\Phi_{pq}^{\frak C\boxplus\tau}$ where $\Phi_{pq}$ is a coordinate change
from $\mathcal U_{q}$
to $\mathcal U_{p}$.
\item
A {\it $\tau$-$\frak C$-collared Kuranishi structure}
\index{Kuranishi structure ! $\tau$-$\frak C$-collared Kuranishi structure}
$\widehat{\mathcal U'}$ on $X^{\frak C\boxplus\tau}$ is the
following objects.
\begin{enumerate}
\item
For each $p' \in \overset{\circ\circ}S_k(X^{\frak C\boxplus\tau})$,
$\widehat{\mathcal U'}$ assigns a
$\tau$-$\frak C$-collared Kuranishi neighborhood $\mathcal U_{p'}$.
\item
For each
$p' \in \overset{\circ\circ}S_k(X^{\frak C\boxplus\tau})$ and
$q' \in \overset{\circ\circ}S_{\ell}(X^{\frak C\boxplus\tau})$
with $q' \in \psi_{p'}(s_{p'}^{-1}(0))$,
$\widehat{\mathcal U'}$ assigns a $\tau$-$\frak C$-collared
coordinate change $\Phi_{p'q'}$.
\item
If $p' \in \overset{\circ\circ}S_k(X^{\frak C\boxplus\tau})$,
$q' \in \overset{\circ\circ}S_{\ell}(X^{\frak C\boxplus\tau})$,
$r' \in \overset{\circ\circ}S_{m}(X^{\frak C\boxplus\tau})$
with
$q' \in \psi_{p'}(s_{p'}^{-1}(0))$ and
$r' \in \psi_{q'}(s_{q'}^{-1}(0))$,
then we require
$$\Phi_{p'q'} \circ \Phi_{q'r'}\vert_{U_{p'q'r'}}
=
\Phi_{p'r'}\vert_{U_{p'q'r'}}
$$
where $U_{p'q'r'} = U_{p'r'} \cap \varphi_{q'r'}^{-1}(U_{p'q'})$.
\end{enumerate}
\item
A {\it $\tau$-$\frak C$-collared K-space}
\index{K-space ! $\tau$-$\frak C$-collared} is a pair of paracompact Hausdorff space
$X^{\frak C\boxplus\tau}$
and its $\tau$-$\frak C$-collared Kuranishi structure
$\widehat{{\mathcal U}^{\frak C\boxplus\tau}}$.
We call
$$
(X^{\frak C\boxplus\tau},\widehat{{\mathcal U}^{\frak C\boxplus\tau}})
$$
the {\it $\tau$-$\frak C$-corner trivialization}, or
{\it partial trivialization of corners}, of $(X,\widehat{\mathcal U})$.
\index{corner ! $\tau$-$\frak C$-corner trivialization}
\index{partial ! trivialization of corners}
\index{corner ! partial trivialization of corners}
\index{trivialization of corners ! partial}
We sometimes write
$$
(X,\widehat{\mathcal U})^{\frak C\boxplus\tau}
$$
in place of
$(X^{\frak C\boxplus\tau},\widehat{{\mathcal U}^{\frak C\boxplus\tau}})$.

\item
We can define the notion
of {\it $\tau$-$\frak C$-collared CF-perturbation,
\index{Kuranishi structure ! $\tau$-$\frak C$-collared CF-perturbation}
\index{CF-perturbation ! $\tau$-$\frak C$-collared}
$\tau$-$\frak C$-collared multivalued perturbation,
\index{multivalued perturbation ! $\tau$-$\frak C$-collared}
$\tau$-$\frak C$-collared good coordinate system,
\index{good coordinate system ! $\tau$-$\frak C$-collared}
$\tau$-$\frak C$-collared Kuranishi chart,
$\tau$-$\frak C$-collared vector bundle,
$\tau$-$\frak C$-collared smooth section,
$\tau$-$\frak C$-collared embedding} of various kinds,
etc.
in the same way.
\end{enumerate}
\end{defn}
The decomposition $\frak C$ induces a decomposition of the boundary
$\partial(X^{\frak C\boxplus\tau},\widehat{\mathcal U}^{\frak C\boxplus\tau})$
in an obvious way.
We also dente it by $\frak C$.
\begin{lem}
Lemma \ref{lemma1523} can be generalized to a ${\frak C}$-version
in an obvious way.
\end{lem}
\begin{lem}\label{lem153535revref}
If $(X',\widehat{\mathcal U'})$ is
$\tau$--$\frak C$-collared, then for any $0 < \tau' < \tau$,
$X'$ has a
$\tau'$-$\frak C$-collared Kuranishi structure determined
in a canonical way from the $\tau$--$\frak C$-collared Kuranishi structure $(X',\widehat{\mathcal U'})$.
The same holds for CF-perturbation,
multivalued perturbation and good coordinate system.
\end{lem}
The proof is the same as the proof of Lemma \ref{lem153535}.
\begin{shitu}\label{situ1613}
Let $(X,\widehat{\mathcal U})$ be a K-space.
Suppose $\frak C_1$, $\frak C_2$, $\frak C$ are
decompositions of $\partial(X,\widehat{\mathcal U})$ as in
Situation \ref{decomporbbdrkura}.
We assume the following two conditions.
\begin{enumerate}
\item
$\partial_{\frak C_1}U_p \cap \partial_{\frak C_2}U_p
= \emptyset$.
\item
$\partial_{\frak C_1}U_p \cup \partial_{\frak C_2}U_p
= \partial_{\frak C}U_p$.$\blacksquare$
\end{enumerate}
\end{shitu}
\begin{lem}\label{lem1614}
In Situation \ref{situ1613} we have the following canonical isomorphisms.
$$
((X,\widehat{\mathcal U})^{\frak C_1\boxplus\tau})^{\frak C_2\boxplus\tau}
\cong
((X,\widehat{\mathcal U})^{\frak C_2\boxplus\tau})^{\frak C_1\boxplus\tau}
\cong
(X,\widehat{\mathcal U})^{\frak C\boxplus\tau}.
$$
\end{lem}
\begin{rem}
The decomposition $\frak C_2$ of
$\partial(X,\widehat{\mathcal U})$ induces one on
$\partial((X,\widehat{\mathcal U})^{\frak C_1\boxplus\tau})$,
which we denote by the same symbol.
We used this fact in the statement of Lemma \ref{lem1614}.
\end{rem}
\begin{proof}[Proof of Lemma \ref{lem1614}]
Suppose a Kuranishi chart of $(X,\widehat{\mathcal U})$
is given as a quotient open subset $V_p$ of
$\overline V_p \times [0,1)^{k_1}\times [0,1)^{k_2}\times [0,1)^{k_3}$
where
$$
\aligned
\partial_{\frak C_1}V_p
&= V_p \cap \left(\overline V_p \times [0,1)^{k_1}\times\partial [0,1)^{k_2}\times [0,1)^{k_3}\right),
\\
\partial_{\frak C_2}V_p
&= V_p \cap \left(\overline V_p \times [0,1)^{k_1}\times [0,1)^{k_2}\times\partial [0,1)^{k_3}\right).
\endaligned
$$
Then
$$
\partial_{\frak C}V_p
= V_p \cap \left(\overline V_p \times [0,1)^{k_1}\times\partial([0,1)^{k_2}\times [0,1)^{k_3})\right).
$$
Therefore we have
$$
V_p^{\frak C_1\boxplus\tau}
=
\mathcal R^{-1}(V_p) \cap \left(\overline V_p \times [0,1)^{k_1}\times [-\tau,1)^{k_2}\times [0,1)^{k_3}\right)
$$
and
$$
V_p^{\frak C\boxplus\tau}
=
\mathcal R^{-1}(V_p) \cap\left(\overline V_p \times [0,1)^{k_1}\times [-\tau,1)^{k_2}\times [-\tau,1)^{k_3})\right).
$$
Lemma \ref{lem1614} follows easily.
\end{proof}

\subsection{In which sense smoothing corner is canonical?}
\label{subsec:subsec32-1}

We next discuss smoothing corner.
Smoothing corner of manifolds is a standard process
and its generalization to orbifolds is also straightforward.
An issue to generalize the smoothing corner to Kuranishi structures
lies in the way of how we fix a smooth structure around the corners
and how much we can make the smooth structure canonical.
This technicalities can be go around fairly
nicely especially when we use {\it trivialization of
corner} introduced in Section \ref{sec:triboundary},
that is exactly the case we use in our story.
We discuss those issues in Subsections
\ref{subsec:subsec32-1}-\ref{subsec:smoothcornerkstr}.

We begin with the following remark.
Let $M$ be a manifold (or an orbifold) with corners.
We have a structure of manifold with boundary
(but without corner) on the {\it same}
underlying topological space.
We write this manifold with boundary (but without corner) as $M'$.
We denote by $i' : M \to M'$ the identity map.
Then it has the following properties.
\begin{enumerate}
\item[(*)]
$i'$ induces  a smooth embedding  $\widehat S_k(M) \to M'$.
\end{enumerate}
Moreover if $M$ is admissible, we have the following:
\begin{enumerate}
\item If $f : M\to \R$ is an admissible function, then
$f \circ (i')^{-1}$ is smooth.
\item
If $\mathcal E \to M$ is an admissible vector bundle, the underlying continuous
map $\mathcal E \to M$ has a structure of $C^{\infty}$-vector bundle on $M'$.
We write it $\mathcal E' \to M'$.
\item
If $s$ is an admissible section of $\mathcal E$,
the same (set-theoretical) map $M' \to \mathcal E'$
is a smooth section.
\end{enumerate}
The proofs of (1)-(3) above are easy. Since we do not use them, we do not prove them.
We next see how much the smooth structure of $M'$ is canonical.
The following statement is also standard.
\begin{lem}\label{canosmooth}
We can construct $M'$ from $M$ in such a way that the differential manifold
$M'$ is well-defined modulo diffeomorphism.
More precisely, we have the following:
Suppose we obtain another $M''$ from $M$.
Let $i'' : M \to M''$ be the identity map.
Then we have a diffeomorphism $f : M' \to M''$ such that
$
f(i'(S_k(M)) = i''(S_k(M)),
$ for $k=0,1,2,\dots$.
\end{lem}
In other words, $M'$ is well-defined modulo stratification preserving diffeomorphism. Here the stratification means the corner structure
stratification of $M$.
This lemma is fairly standard and its proof is omitted.
It seems that it is more nontrivial to find a `canonical' way
so that the above diffeomorphism $f$ can be taken to be the identity map.
In this article, we do not try to find such a way in general situation
of (admissible) orbifold, but will do so in the case of {\it collared orbifold}.
\par
Before doing so, we explain a reason why the uniqueness in the sense of Lemma \ref{canosmooth}
is not enough for our purpose.
When we generalize the process of smoothing corner to that for Kuranishi structures,
we need to study the situation where we have an embedding $N \to M$ of
cornered orbifolds. When we smooth the corners of $N$ and $M$,
we want the {\it same} map $N' \to M'$ to be a {\it smooth} embedding.
This is not obvious because of the non-uniqueness
of the smooth structure we put on $M'$ and $N'$.
However,
it is still true and not difficult to prove that we {\it can find} a smooth structures of $M'$ and $N'$
so that $N' \to M'$ is a smooth embedding.
\par
On the other hand, in order to smooth the corner of Kuranishi structure we need to smooth the
corner of  all the Kuranishi charts,
simultaneously. This now becomes a nontrivial problem.
If we try to use the uniqueness in Lemma \ref{canosmooth},
we should include the diffeomorphism $f$ as a part of data
in the construction.
Then the
compatibility of the
coordinate change might be broken.\footnote
{It seems that we can still prove that for a given good coordinate system
we can construct a smoothing corner compatible with the
coordinate change for any proper open substructure of it.
We need to work out rather cumbersome induction to prove it
at the level of detailed-ness we intend to achieve in this article.}
\par
We go around this issue by using the collar.
When we use the collar, the way how we smooth the corner still involves
choices.
However, we can make the choice to smooth
the model $[0,1)^k$ only once and then use that particular
choice to smooth all the collared orbifolds, simultaneously.
This way is canonical enough so that all the embeddings
of Kuranishi charts become smooth embeddings automatically after
smoothing the corners.

\subsection{Smoothing corner of $[0,\infty)^k$}
\label{subsec:smoothcornermodel}
\index{corner ! smoothing corner}
\index{smoothing corner ! smoothing corner}

In this subsection, we fix data which we need to smooth the corner of a partially collared orbifold and a partially collared Kuranishi structure.
Namely we fix a way to smooth the local model $[0,\infty)^k$ so that it is compatible
with various $k$ and also with the symmetry exchanging the factors.
The latter is important to study the case of orbifolds.
Let ${\rm Perm}({k+1})$ be a group of permutations of $\{1,\dots,k+1\}$.
\begin{defn}\label{def:permaction}
We define a ${\rm Perm}({k+1})$ action on $\R^{k}$ as follows.
We regard
$$
\R^{k} = \left\{ \left.(t_0,\dots,t_k) \in \R^{k+1}
~\right\vert~
\sum_{i=0}^k t_i   = 0\right\}.
$$
The group ${\rm Perm}({k+1})$ acts on $\R^{k+1}$ by exchanging the factors.
It restricts to an action of ${\rm Perm}({k+1})$ to $\R^{k}$.
In this section, the ${\rm Perm}({k+1})$ action on $\R^{k}$ always means
this particular action.
\end{defn}

Below we will show the existence of
a set of homeomorphisms
\begin{equation}
\Phi_{k} : [0,\infty)^{k} \to \R^{k-1} \times [0,\infty)
\end{equation}
and smooth structures $\frak{sm}_{k}$ on $ [0,\infty)^{k}$,
simultaneously, for any $k \in \Z_{+}$, which satisfy
Condition \ref{cornersmoexi}.

\begin{conds}\label{cornersmoexi}
We require $\Phi_{k}$ and $\frak{sm}_k$ to satisfy the following conditions.
\begin{enumerate}
\item
$\Phi_{k}$ is a diffeomorphism from $( [0,\infty)^{k},\frak{sm}_k)$
to $\R^{k-1} \times [0,\infty)$. Here we use the standard smooth structure on $\R^{k-1} \times [0,\infty)$.
\item
Let $\Phi_{k}(t_1,\dots,t_k) = (x,t)$. Then
$$
\Phi_{k}(ct_1,\dots,ct_k) = (x,ct)
$$
for $c \in \R_+$.
\item
For ${\bf t} \in [0,\infty)^{k} \setminus 0$
there exists an open neighborhood of ${\bf t}$ that is
isometric to $V \times [0,\epsilon)^{\ell}$ where $V \subset \R^{k-\ell}$.
(Here $V$ is an open set.
When we say `isometric', we use the
Euclidian metrics
on $[0,\infty)^{k}$, $\R^{k-\ell}$, $[0,\epsilon)^{\ell}$.)
Then the map
$$
{\rm id} \times  \Phi_{\ell} : V \times [0,\epsilon)^{\ell} \to V \times \R^{\ell-1}  \times [0,\infty)
$$
is a diffeomorphism onto its image.
Here we put the restriction of the smooth structure $\frak{sm}_k$ to
$V \times [0,\epsilon)^{\ell}$. (The space $V \times [0,\epsilon)^{\ell}$ is
identified with an open subset of $[0,\infty)^{k} \setminus 0$ by the isometry.)
\item
The map $
\Phi_{k} : [0,\infty)^{k} \to \R^{k-1} \times [0,\infty)
$
is ${\rm Perm}({k})$ equivariant.
Here ${\rm Perm}({k})$ acts on $[0,\infty)^{k}$ by permutation of factors,
and acts on $\R^{k-1}$ by
Definition \ref{def:permaction}.
On the last factor $[0,\infty)$ the action is trivial.
\end{enumerate}
\end{conds}
\begin{rem}
In (3) above we require a neighborhood of ${\bf t}$ to
be isometric to $V \times [0,\epsilon)^{\ell}$.
The reason why we require such a rather restrictive assumption that they
are isometric is that we want to specify the diffeomorphism
from a neighborhood of ${\bf t}$ to $V \times [0,\epsilon)^{\ell}$.
We use the Euclidian metric on $\R^n$ here  only to make
the choice of diffeomorphism (that is isometry) as
canonical as possible.
\end{rem}
\begin{lem}
For any $k \in \Z_+$
there exist $\Phi_{k}$ and $\frak{sm}_k$ satisfying Condition \ref{cornersmoexi}.
\end{lem}
\begin{proof}
The proof is by induction on $k$. For $k=1$,
$\Phi_1$ is the identity map and  $\frak{sm}_1$ is the standard
smooth structure on $[0,\infty)$.
\par
Suppose we have $\Phi_{i}$, $\frak{sm}_i$ for $i < k$.
We observe that Condition (3) determines a smooth structure $\frak{sm}_k$
on $[0,\infty)^{k} \setminus 0$ uniquely.
Indeed, well-defined-ness of this structure can be checked by
Condition (3) itself inductively.
Moreover, by the definition of the smooth structure
$\frak{sm}_k$,
we find that the ${\rm Perm}({k})$-action is smooth with respect
to this smooth structure.
The map $(t_1,\dots,t_k) \mapsto (ct_1,\dots,ct_k)$
is also a diffeomorphism for this smooth structure
if $c \in \R_+$.
\par
Next we will construct a homeomorphism $\Phi_k$ and
extend the smooth structure $\frak{sm}_k$ to
$[0,\infty)^{k}$.
We choose a compact subset $S \subset [0,\infty)^{k} \setminus 0$ which is a smooth $(k-1)$-dimensional
submanifold with corners (with respect to the standard
structure of manifold with corners of $ [0,\infty)^{k+1}$) such that:
\begin{enumerate}
\item[(a)]
$S$ is a slice of the multiplicative $\R_+$ action on $[0,\infty)^{k} \setminus 0$.
\item[(b)]
$S$ is perpendicular to all the strata $\overset{\circ}S_{\ell}([0,\infty)^{k})$. We can take
a tubular neighborhood of $\overset{\circ}S_{\ell}([0,\infty)^{k}) \cap S$ in $S$
such that the fiber of the projection to $\overset{\circ}S_{\ell}([0,\infty)^{k}) \cap S$
is flat with respect to the Euclidean metric of $[0,\infty)^k$.
\item[(c)]
$S$ is invariant under the ${\rm Perm}({k})$-action on $[0,\infty)^{k}$.
\end{enumerate}
\begin{figure}[h]
\centering
\includegraphics[scale=0.3]{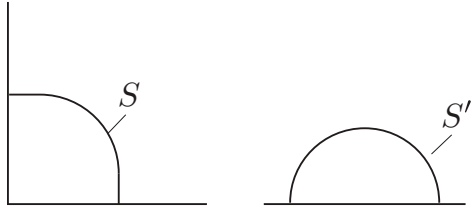}
\caption{submanifolds $S$ and $S'$}
\label{Figure17-2}
\end{figure}
We can find such $S$ by ${\rm Perm}({k})$-equivariantly
modifying the intersection of unit ball and  $[0,\infty)^{k}$
around the boundary a bit.
\par
By Condition (b)  $S$ is a smooth submanifold with boundary
of $([0,\infty)^{k} \setminus 0,\frak{sm}_{k})$.
Since we can construct $S$ by modifying the intersection of unit ball and  $[0,\infty)^{k}$
around the boundary a bit,  there exists a ${\rm Perm}({k})$ equivariant diffeomorphism from
$S$ to
$$
S' = \{ (x,t) \in \R^{k-1} \times \R_{\ge 0} \mid \Vert x\Vert^2 + t^2 = 1\}.
$$
Here we use the smooth structure of $S$ induced from
$([0,\infty)^{k} \setminus 0,\frak{sm}_{k})$ and the standard
smooth structure on $S'$.
See Figure \ref{Figure17-2}.
(Note that
$S$ becomes a manifold with boundary and without corner
with respect to this smooth structure.)
We fix this diffeomorphism. Then we can define
$\Phi_k$ by extending the diffeomorphism $S \to S'$ so that Condition (2) is satisfied.
By construction, $\Phi_k$ is a diffeomorphism outside the origin.
\par
Then we extend the smooth structure $\frak{sm}_{k}$ to the origin
so that $\Phi_k$ also becomes a diffeomorphism at the origin.
The proof is now complete by induction.
\end{proof}
\begin{rem}
During the proof we made choices of $S$ and a diffeomorphism
between $S$ and $S'$ for each $k$.
The resulting smooth structure $\frak{sm}_{k}$ {\it depends}
on these choices in the sense that the identity map is not a diffeomorphism
when we use two smooth structures obtained by different choices
for the source and the target.
However since two different choices of $S$ and the diffeomorphism $S \to S'$
are isotopic to each other,
the resulting smooth structure  $\frak{sm}_{k}$
is independent of the choices in the sense of diffeomorphism.
(This is a proof of Lemma \ref{canosmooth} in this case.)
\end{rem}

When we apply the construction of smoothing corner of
Kuranishi structure, we sometimes need to put collars
to the smoothed K-space.
We use Lemma \ref{lem328col} below for this purpose.
Let $\Phi_k$ and $\frak{sm}_k$ be as in Condition \ref{cornersmoexi}.

\begin{conds}\label{cond327}
For any $k \in \Z_+$
we consider $\frak{Trans}_{k-1}$ and $\Psi_k$ with the following
properties.
\begin{enumerate}
\item
$\frak{Trans}_{k-1}$ is a smooth $(k-1)$-dimensional submanifold of $[0,\infty)^k$
and is contained in $(0,\infty)^k \setminus (1,\infty)^k$.
\item
$\frak{Trans}_{k-1}$ is invariant under the ${\rm Perm}(k)$ action on
$[0,\infty)^k$.
\item
$
\frak{Trans}_{k-1} \cap ([0,\infty)^{k-1} \times [1,\infty))
=
\frak{Trans}_{k-2} \times [1,\infty).
$
This is an equality
as subsets of $[0,\infty)^k
= [0,\infty)^{k-1} \times [0,\infty)$.
\item
$$
\Psi_k : [0,1] \times \frak{Trans}_{k-1}
\to [0,\infty)^k
$$
is a homeomorphism onto its image.
Let $\frak U_k$ be its image.
\item
Using the smooth structure $\frak{sm}_k$ on $[0,\infty)^k$,
the subset
$\frak U_k \subset [0,\infty)^k$ is a smooth $k$-dimensional submanifold with boundary
and $\Psi_k$ is a diffeomorphism. Moreover
$$
\partial \frak U_k
= \partial([0,\infty)^k) \cup \frak{Trans}_{k-1}
$$
and the restriction of $\Psi_k$ to $\{0\}\times \frak{Trans}_{k-1}$
is a diffeomorphism onto $\partial([0,\infty)^k)$.
The restriction of $\Psi_k$ to $\{1\}\times \frak{Trans}_{k-1}$
is the identity map.
\item
$\Psi_k$ is equivariant under the ${\rm Perm}(k)$ action.
(The ${\rm Perm}(k)$
action on $\frak{Trans}_{k-1}$ is defined in Item (2) and
the action
on $[0,\infty)^k$ is by permutation of factors.
\item
If $s \ge 1$, $t \in [0,1]$ and $(x_1,\dots,x_{k-1}) \in \frak{Trans}_{k-2}$,
then
$$
\Psi_k(t,(x_1,\dots,x_{k-1},s)) = (\Psi_{k-1}(t,(x_1,\dots,x_{k-1})),s).
$$
Here we use the identification in Item (3) to define the left hand side.
\end{enumerate}
\end{conds}
\begin{figure}[h]
\centering
\includegraphics[scale=0.3]{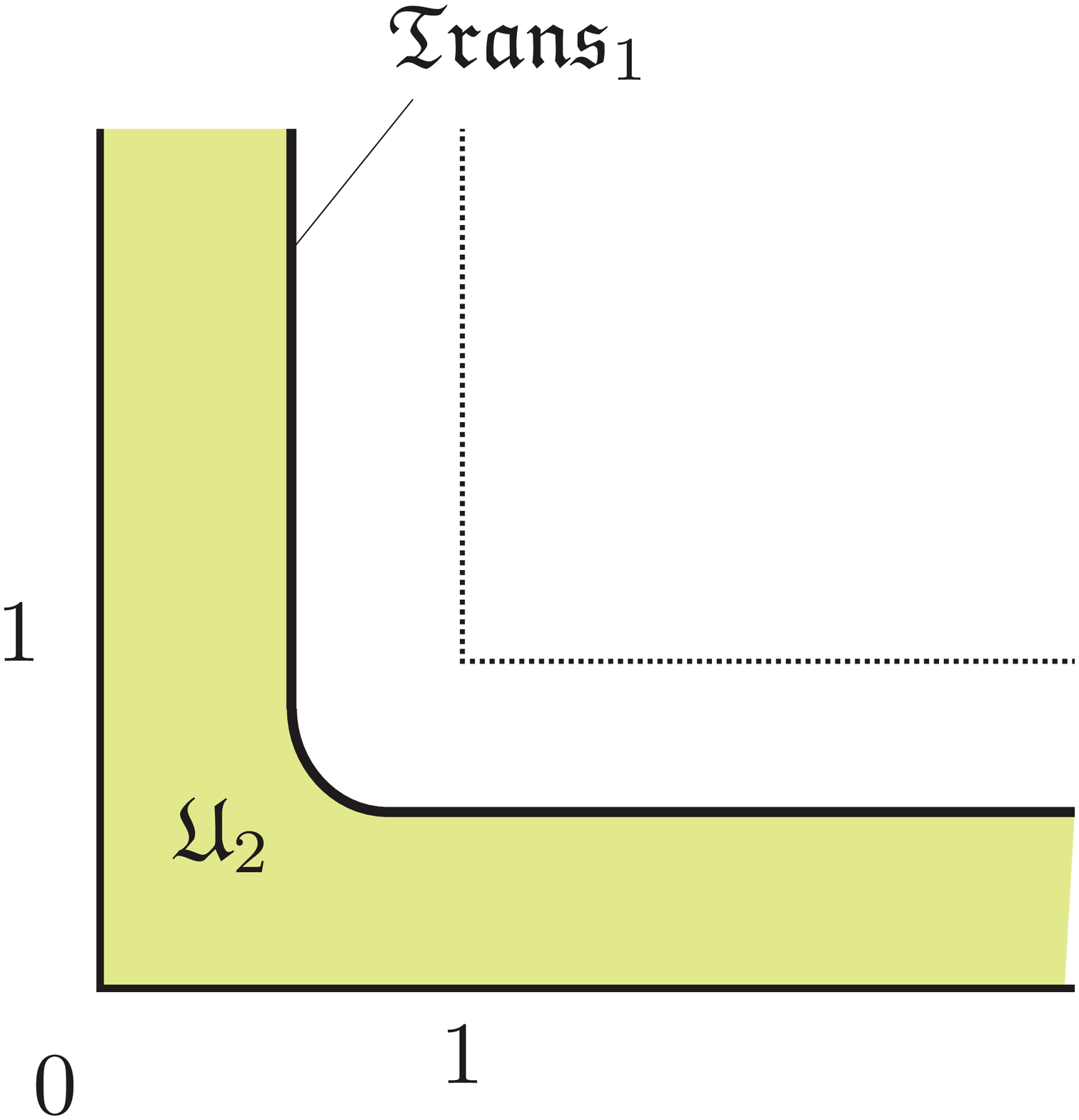}
\caption{$\frak{Trans}_{k-1}$ and $\frak U_k$}
\label{Figure17-3}
\end{figure}
\begin{lem}\label{lem328col}
For any $k \in \Z_+$ there exist
$\frak{Trans}_{k-1}$ and $\Psi_k$ satisfying Condition \ref{cond327}.
Moreover, for each given $\delta > 0$, we may take them so that
$\frak U_k$ contains $[0,\infty)^k \setminus [1-\delta,\infty)^k$.
\end{lem}
\begin{proof}
The proof is by induction.
If $k=1$, we put $\frak{Trans}_{0} = \{ 1 - \delta/2\}$ and
there is nothing to prove.
Suppose we have $\frak{Trans}_{k'-1}$, $\Psi_{k'}$ for $k' < k$.
Conditions \ref{cond327} (2) (3) determine
$\frak{Trans}_{k-1}$ outside $[0,1]^k$.
Conditions \ref{cond327} (6) (7) determine $\Psi_k$ outside
$[0,1]^k$.
It is easy to see that we can extend them to $[0,1]^k$ and obtain $\frak{Trans}_{k-1}$, $\Psi_k$.
\end{proof}
We note that $\Psi_k$ defines a collar of $([0,\infty)^k,\frak{sm}_k)$,
which is a manifold with boundary (but without corner).

\subsection{Smoothing corner of collared orbifolds and of Kuranishi structures}
\label{subsec:smoothcornerkstr}

In this subsection we combine the story of partial trivialization of
Kuranishi structure in Subsection \ref{subsec:prtialtrivi} with the story of smoothing corner
in Subsection \ref{subsec:subsec32-1}-\ref{subsec:smoothcornermodel}.

\begin{shitu}\label{sit327}
Let $U$ be an orbifold.
We consider its normalized boundary $\partial U$.
Let $\frak C$ be a decomposition of $\partial U$ as in Situation \ref{decomporbbdr}.
We denote by $U^{\frak C\boxplus\tau}$ the partial trivialization
of corners.
$\blacksquare$
\end{shitu}

We will define smoothing corners of $U^{\frak C\boxplus\tau}$ below.
\par
Let $p \in \overset{\circ\circ}S_k(U^{\frak C\boxplus\tau})$
and put $\overline p = \mathcal R^{\frak C}(p)$.
The point $p$ has an orbifold neighborhood
\begin{equation}\label{form3222}
\frak V_p = (\overline V_{\overline p} \times [-\tau,0)^k,
\Gamma_{\overline p},\phi_p).
\end{equation}
Here
$(\overline V_{\overline p} \times [0,1)^k,
\Gamma_{\overline p},\phi_{\overline p})$
is an orbifold neighborhood of $\overline p$ in $U$ with
$\overline p \in \overset{\circ}{S^{\frak C}_{k}}(U)$, and
$\overline V_{\overline p}$ may have boundary or corners
but the boundary of $\overline V_{\overline p}$ does not
correspond to a boundary component in $\frak C$.

\begin{defn}\label{defn3277}
We define a smooth structure of $\frak V_p$ in (\ref{form3222}) as follows.
We identify $[-\tau,0)^k \cong [0,\tau)^k$ by the diffeomorphism
$(t_1,\dots,t_k) \mapsto (t_1+\tau,\dots,t_k+\tau)$.
We use the smooth structure $\frak{sm}_k$ on $[0,\tau)^k$
to obtain a smooth structure on $[-\tau,0)^k$.
Then by taking direct product we obtain
a smooth structure on $\overline V_{\overline p} \times [-\tau,0)^k$.
\end{defn}

\begin{lem}\label{lem3288}
Consider the smooth structure on $\overline V_{\overline p} \times [-\tau,0)^k$ as
in Definition \ref{defn3277}. Then we have:
\begin{enumerate}
\item
The $\Gamma_{\overline p}$ action on
$\overline V_{\overline p} \times [-\tau,0)^k$ is smooth with respect to this
smooth structure.
\item
If $q \in \psi_p(\overline V_{\overline p} \times [-\tau,0)^k)$
and $q \in \overset{\circ\circ}S_{\ell}(U^{\frak C\boxplus\tau})$,
then the coordinate change from $\frak V_q$ to $\frak V_p$ is
smooth with respect to the above smooth structure.
\end{enumerate}
\end{lem}
\begin{proof}
Item (1) is a consequence of Condition \ref{cornersmoexi} (4) and
Item (2) follows from Condition \ref{cornersmoexi} (3).
\end{proof}

\begin{defn}\label{def:partialsmooth}
We obtain an atlas of an orbifold structure on $U^{\frak C\boxplus\tau}$
by Lemma \ref{lem3288} and Definition \ref{defn3277}.
We call $U^{\frak C\boxplus\tau}$ with this smooth structure
the {\it orbifold with corner obtained by $\tau$-$\frak C$-partial smoothing of
corners},
\index{corner ! $\tau$-$\frak C$-partial smoothing of
corners}
\index{smoothing corner ! $\tau$-$\frak C$-partial smoothing of
corners}
and when no confusion can occur,
we simply call the smoothing {\it $\tau$-partial smoothing of
corners}.
\index{corner ! $\tau$-partial smoothing of
corners}
\index{smoothing corner ! $\tau$-partial smoothing of
corners}
We denote it by
$$
U^{{\rm sm}\frak C\boxplus\tau}.
$$
In case $\frak C$ is the whole set of all the components of the boundary,
this smooth structure has no corner.
\end{defn}

\begin{lem}
If $p \in \overset{\circ} S_{k+\ell}(U^{\frak C\boxplus\tau})
\cap  \overset{\circ}{S^{\frak C}_{k}}(U^{\frak C\boxplus\tau})$,
then $p \in \overset{\circ}S_{\ell+1}(U^{{\rm sm}\frak C\boxplus\tau})$.
\end{lem}

The proof is obvious.
We also have the following:

\begin{lem}
If $0< \tau' < \tau$, the orbifold obtained by $\tau$-partial smoothing of
corners is $\tau'$-collared.
\end{lem}
\begin{proof}
Let $(V,\Gamma,\phi)$ be an orbifold chart of $U$.
We may chose $V$ so that it is an open subset of
$\overline V \times [0,1)^k \times [0,1)^{k'-k}$
and Convention \ref{conv1633} is satisfied.
Then the corresponding orbifold chart of $U^{\frak C\boxplus\tau}$
is
$$
V =
\mathcal R^{-1}(V) \cap
\left(
\overline V \times [-\tau,1)^k \times [0,1)^{k'-k}
\right).
$$
We use the smooth structure $\frak{sm}_k$ on the
$ [-\tau,1)^k$ factor (which we identity with $[0,1+\tau)^k$)
and obtain the orbifold chart of
$U^{{\rm sm}\frak C\boxplus\tau}$.
\par
Now we use $\frak{Trans}_{k-1}$ and $\Psi_k$ produced in Lemma \ref{lem328col}.
The map $\Psi_k$  is a diffeomorphism
$$
\Psi_k : \frak{Trans}_{k-1} \times [0,1) \to [-\tau,0)^k
$$
to the image. The image is a neighborhood of
$\partial [-\tau.0)^k$.
Then we have a smooth embedding
$$
\overline V \times \frak{Trans}_{k-1} \times [0,1) \times [0,1)^{k'-k}
\supseteq V'
\to V
$$
where $V'$ is an open subset. Its image is a neighborhood of
$\partial V$.
This embedding is a diffeomorphism to its image
if we use the differential structure after smoothing corners.
This gives a collar of $U$ on this chart.
Using Conditions \ref{cond327} (6)(7), we can show that this collar is
compatible with the coordinate change.
\par
We can take $\tau' >0$ for any $\tau >0$ by choosing $\delta$ small
in Lemma \ref{lem328col}.
\end{proof}

In the next lemma we summarize the properties of (partial) smoothing corner.

\begin{lem}\label{lem32222}
Suppose we are in Situation \ref{sit327}.
Let $0< \tau' <\tau$.
\begin{enumerate}
\item
If $E$ is an admissible vector bundle on $U$,
$E^{{\rm sm}\frak C\boxplus\tau}$ becomes an
admissible vector bundle on $U^{{\rm sm}\frak C\boxplus\tau}$.
It is $\tau' $-collared.
\item
If $f : U \to M$ is an admissible map, $f^{\frak C \boxplus\tau} : U^{{\rm sm}\frak C\boxplus\tau} \to M$
is an admissible map from $U^{{\rm sm}\frak C\boxplus\tau}$.
When $M$ has corner, we take the decomposition $\frak C$ of the normalized boundary $\partial U$
so that $\partial_{\frak C}U$ is contained in the horizontal boundary
as in Definition \ref{defn3288}. Then
the same assertion also holds.
\item
If $f : U_1 \to U_2$ is an admissible embedding and $f(\partial_{\frak C_1}U_1)\subset
\partial_{\frak C_2}U_2 \cap f(U_1)$, then
$f^{\frak C_1\boxplus\tau} : U_1^{{\rm sm}\frak C_1\boxplus\tau} \to U_2^{{\rm sm}\frak C_2\boxplus\tau}$
is an admissible embedding.
It is $\tau' $-collared.
\item
In the situation of (3) suppose $E_1 \to U_1$, $E_2 \to U_2$ are admissible vector bundles and
$f$ is covered by an admissible embedding of vector bundles $\widehat f : E_1 \to E_2$.
Then $\widehat f^{\frak C_1\boxplus\tau} : E_1^{{\rm sm}\frak C_1\boxplus\tau} \to E_2^{{\rm sm}\frak C_2\boxplus\tau}$
is an admissible embedding of vector bundles which covers $f^{\frak C_1\boxplus\tau}$.
It is $\tau' $-collared.
\item
In the situation of (1) if $s$ is an admissible section of $E$,
$s^{\frak C\boxplus\tau}$ is an admissible section of $E^{{\rm sm}\frak C\boxplus\tau}$.
It is $\tau' $-collared.
\end{enumerate}
\end{lem}
The proof is obvious.
\par
Now we consider the case of Kuranishi structure.
This generalization is quite straightforward.
\begin{shitu}\label{sit328}
Let $(X,\widehat{\mathcal U})$ be a K-space.
Suppose for each $p \in X$ we have a decomposition of
$\partial U_p$ into two unions of connected components.
Let $\frak C_p$ be the first union. We assume that
for each $q \in \psi_p(s_p^{-1}(0))$ we have
$$
\varphi_{pq}(\partial_{\frak C_q}U_q\cap U_{pq})
=
\varphi_{pq}(U_{pq}) \cap \partial_{\frak C_p}U_p.
$$
$\blacksquare$
\end{shitu}

\begin{defn}\label{defsmoothing}
In Situation \ref{sit328}, let
$
(X^{\frak C\boxplus\tau},\widehat{{\mathcal U}^{\frak C\boxplus\tau}})
$
be a $\tau$-$\frak C$-collared K-space
as in Definition \ref{defn1531revrev}.
We define its {\it partial smoothing of corners}
as follows:
\par
If $p \in \overset{\circ\circ}S_k(X^{\frak C\boxplus\tau})$
and its partially $\tau$-collared Kuranishi neighborhood is
$\mathcal U_{\overline p}^{\frak C\boxplus\tau}
= (U_{\frak C\overline p}^{\frak C\boxplus\tau},E_{\frak C\overline p}^{\boxplus\tau},s_{\frak C\overline p}^{\frak C\boxplus\tau},\psi_{\overline p}^{\frak C\boxplus\tau})$,
then we change its smooth structure by smoothing corner.
We obtain a Kuranishi chart denoted by
$\mathcal U_{\overline p}^{{\rm sm}\frak C\boxplus\tau}$.
Then the coordinate change is induced by one of $\widehat{{\mathcal U}^{\frak C\boxplus\tau}}$
because of Lemma \ref{lem32222}.
It is $\tau' $-collared for any $0< \tau' < \tau$.
We call the Kuranishi structure obtained above the
{\it Kuranishi structure obtained by partial smoothing of corners}
\index{Kuranishi structure ! partial smoothing of corners}
\index{smoothing ! Kuranishi structure obtained by partial smoothing of corners}
and denote it by
$$
\widehat{{\mathcal U}^{{\rm sm}\frak C\boxplus\tau}}.
$$
\index{$
\widehat{{\mathcal U}^{{\rm sm}\frak C\boxplus\tau}}
$}
\end{defn}
\begin{lem}\label{lem1633}
In the situation of Definition \ref{defsmoothing} the following ($\tau$-collared)
object on $(X,\widehat{\mathcal U})$ induces the corresponding ($\tau'$-collared)
object
of the Kuranishi structure obtained by partial
smoothing of corners.
\begin{enumerate}
\item
Strongly smooth map.
\item
CF-perturbation.
\item
Multivalued perturbation.
\item
Differential form.
\end{enumerate}
\end{lem}
The proof is immediate from construction.
\begin{lem}\label{lem1634lemlem}
In the situation of Lemma \ref{lem1633}, we put
$$
\partial_{\frak C}(X,\widehat{\mathcal U})
=
\coprod_{c\in \frak C} \partial_c(X,\widehat{\mathcal U}).
$$
Let $\widehat{{\mathcal U}^{\frak C}}$ be the
Kuranishi structure of
$\partial_{\frak C}(X,\widehat{\mathcal U})$ with smoothed corners.
We consider the K-space
$(\partial_{\frak C}(X),\widehat{{\mathcal U}^{\frak C}})$.
Let $h$ be a differential form on $X$ and
$f : X \to M$ a strongly smooth map which is
strongly suvbmersive with respect to a
CF-perturbation $\frak S^{\epsilon}$.
The CF-perturbation $\frak S^{\epsilon}$
induces one on $(\partial_{\frak C}(X),\widehat{{\mathcal U}^{\frak C}})$ and on
$\partial_c(X,\widehat{\mathcal U})$.
We also denote them by $\frak S^{\epsilon}$.
Then we have
$$
\sum_c f!(h\vert_{\partial_c(X)};\frak S^{\epsilon})
=
f!(h\vert_{\partial_{\frak C}(X)};\frak S^{\epsilon}).
$$
\end{lem}
The proof is obvious and so omitted.

\subsection{Composition of morphisms of linear K-systems}
\label{subsec:complinkurasmcorner}

We refer Conditions \ref{linsysmainconds}, \ref{morphilinsys},
\ref{Pparamorphi} and
Definitions \ref{linearsystemdefn}, \ref{linearsystemmorphdefn},
\ref{def:homotopymorph}
for the various definitions and notations concerning linear K-system,
morphisms and homotopy.
\begin{shitu}\label{shitu1634}
Suppose we are in Situation \ref{composisitu} and
$\frak N_{i+1 i}$ is a morphism from $\mathcal F_i$ to
$\mathcal F_{i+1}$.
We denote by
$$\mathcal M^i(\alpha_-,\alpha_+)$$
the space of connecting orbits of $\mathcal F_i$
and by
$$\mathcal N_{i i+1}(\alpha_-,\alpha_+)$$
the interpolation space of $\frak N_{i+1 i}$.
Let $R_{\alpha_i}^i$ be a critical submanifold
of $\mathcal F_i$ and $\alpha_i \in \frak A_i$.
$\blacksquare$
\end{shitu}
\begin{rem}\label{rem:orderindex}
In the above definition we denote the interpolation space of a morphism from $\mathcal F_i$ to
$\mathcal F_{i+1}$ by $\mathcal N_{i i+1}(\ast, \ast)$,
while we denote the corresponding morphism by $\frak N_{i+1 i}$
to be compatible with algebraic formulas of compositions of morphisms.
\end{rem}

\begin{defn}\label{defn1635}
In Situation \ref{shitu1634} we define the
{\it partially trivialized fiber product}
\index{partially trivialized fiber product}
\index{corner ! partially trivialized fiber product}
\index{fiber product ! partially trivialized fiber product}
\begin{equation}\label{form1616}
\mathcal N_{12}(\alpha_1,\alpha_2) \times^{\boxplus\tau}_{R^2_{\alpha_2}}
\mathcal N_{23}(\alpha_2,\alpha_3)
\end{equation}
as follows.
We consider the
fiber product
$
\mathcal N_{12}(\alpha_1,\alpha_2) \times_{R^2_{\alpha_2}}
\mathcal N_{23}(\alpha_2,\alpha_3)
$
equipped with fiber product Kuranishi structure.
Then we consider the decomposition $\frak C$ of its boundary
which consists of the following two kinds of
components of its normalized boundary.
$$
\aligned
&\mathcal N_{12}(\alpha_1,\alpha'_2) \times_{R^2_{\alpha'_2}}
\mathcal M^2(\alpha_{\alpha'_2},\alpha_{\alpha_2})
\times_{R^2_{\alpha_2}}
\mathcal N_{23}(\alpha_2,\alpha_3), \\
&\mathcal N_{12}(\alpha_1,\alpha_2) \times_{R^2_{\alpha_2}}
\mathcal M^2(\alpha_{\alpha_2},\alpha_{\alpha'_2})
\times_{R^2_{\alpha'_2}}
\mathcal N_{23}(\alpha'_2,\alpha_3).
\endaligned
$$
The first line is contained in
$
\partial\mathcal N_{12}(\alpha_1,\alpha_2) \times_{R^2_{\alpha_2}}
\mathcal N_{23}(\alpha_2,\alpha_3)
$
and the second line is contained in
$
\mathcal N_{12}(\alpha_1,\alpha_2) \times_{R^2_{\alpha_2}}
\partial\mathcal N_{23}(\alpha_2,\alpha_3)
$.
Now we define (\ref{form1616}) by
$$
(\mathcal N_{12}(\alpha_1,\alpha_2) \times_{R^2_{\alpha_2}}
\mathcal N_{23}(\alpha_2,\alpha_3))^{\frak C\boxplus\tau}.
$$
\end{defn}
\begin{prop}\label{prop1636}
In Situation \ref{shitu1634} there exists a
$\tau$-$\frak C$-collared K-space
$$
\mathcal{N}_{123}(\alpha_1,\alpha_3)
$$
with the following properties:
\begin{enumerate}
\item
Its normalized boundary in $\frak C$ is isomorphic to the disjoint union of
\begin{equation}\label{prop1636form1}
\mathcal N_{12}(\alpha_1,\alpha_2) \times^{\boxplus\tau}_{R^2_{\alpha_2}}
\mathcal N_{23}(\alpha_2,\alpha_3)
\end{equation}
over various $\alpha_2$.
\item
Its normalized boundary which is not in $\frak C$  is isomorphic to  the disjoint union of
\begin{equation}\label{prop1636form2}
\mathcal M^1({\alpha_1,\alpha'_1})
\times_{R^1_{\alpha'_1}} \mathcal{N}_{123}(\alpha'_1,\alpha_3)
\end{equation}
over various $\alpha'_1$ and of
\begin{equation}\label{prop1636form3}
\mathcal{N}_{123}(\alpha_1,\alpha'_3)
\times_{R_{\alpha'_3}} \mathcal M^3({\alpha'_3,\alpha_3})
\end{equation}
over various $\alpha'_3$.
\item
The isomorphisms (1)(2) satisfy the compatibility
conditions, Condition \ref{conds1637}, at the corners.
\item
The evaluation maps, periodicity and orientation isomorphisms
are defined on $\mathcal{N}_{123}(\alpha_1,\alpha_3)$
and commute with the isomorphisms in (1)(2)(3) above.
\item
Similar statements of Conditions \ref{morphilinsys} (V)(IX) hold.
\end{enumerate}
\end{prop}
See Figure \ref{Figure17-40}.
\begin{figure}[h]
\centering
\includegraphics[scale=0.3,angle=-90]{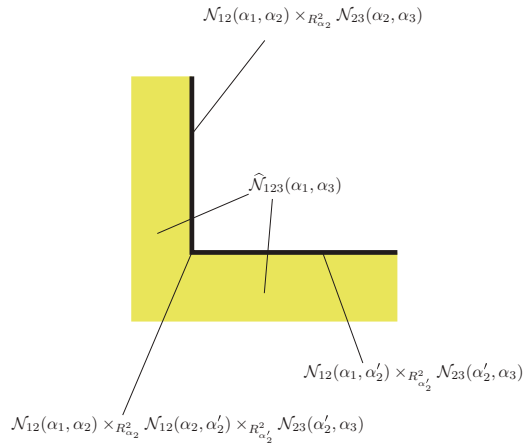}
\caption{$\widehat{\mathcal{N}}_{123}(\alpha_1,\alpha_3)$}
\label{Figure17-40}
\end{figure}
We omit describing the precise ways to modify Conditions \ref{morphilinsys} (V)(IX)
to our situation.
Since they are not hard, we leave them to the readers.
\begin{conds}\label{conds1637}
\begin{enumerate}
\item
The $k$ dimensional normalized corner
$\widehat{S}_k(\mathcal{N}_{123}(\alpha_1,\alpha_3))$ of
$\mathcal{N}_{123}(\alpha_1,\alpha_3)$ is
isomorphic to a disjoint union of the
one of the two types (\ref{eq1617}) and  (\ref{eq1618}):
\begin{equation}\label{eq1617}
\aligned
&\mathcal M^1({\alpha_1,\alpha_1^2}) \times_{R^1_{\alpha_1^2}}
\dots \times_{R^1_{\alpha_1^{k_1-1}}}
\mathcal M^1({\alpha_1^{k_1-1},\alpha_1^{k_1}})
\\
&\times_{R^1_{\alpha_1^{k_1}}}\widehat{S}_m^{\frak C}\left(\mathcal N_{12}
(\alpha_1^{k_1},\alpha_2) \times^{\boxplus\tau}_{R^2_{\alpha_2}}
\mathcal N_{23}(\alpha_2,\alpha_3^{1})\right) \\
&\times_{R^3_{\alpha_3^{1}}}
\mathcal M^3({\alpha_3^{1},\alpha_3^{2}}) \times_{R^3_{\alpha_3^{2}}}
\dots \times_{R^3_{\alpha_3^{k_3-1}}}
\mathcal M^3({\alpha_3^{k_3-1},\alpha_3}),
\endaligned
\end{equation}
and
\begin{equation}\label{eq1618}
\aligned
&\mathcal M^1({\alpha_1,\alpha_1^2}) \times_{R^1_{\alpha_1^2}}
\dots \times_{R^1_{\alpha_1^{k_1-1}}}
\mathcal M^1({\alpha_1^{k_1-1},\alpha_1^{k_1}})
\\
&\times_{R^1_{\alpha_1^{k_1}}}\mathcal{N}_{123}(\alpha_1^{k_1},\alpha_3^1) \\
&\times_{R^3_{\alpha_3^{1}}}
\mathcal M^3({\alpha_3^{1},\alpha_3^{2}}) \times_{R^3_{\alpha_3^{2}}}
\dots \times_{R^3_{\alpha_3^{k_3-1}}}
\mathcal M^3({\alpha_3^{k_3-1},\alpha_3}).
\endaligned
\end{equation}
Here in (\ref{eq1617}), $k_1 + m - 1 + k_3 = k$
and in (\ref{eq1618}), $k_1  + k_3 -2 = k$.
Note $k_1, k_3 \in \Z_{\ge 1}$.
In case $k_1 =1$ the first line of (\ref{eq1617}) or  (\ref{eq1618})
is void.
In case $k_3 =1$ the third line of (\ref{eq1617}) or  (\ref{eq1618})
is void.
We also note that (1) (2) of Proposition \ref{prop1636}
is the case $k=1$ of this condition.
\item
(1) implies that
$\widehat{S}_{\ell}(\widehat{S}_k(\mathcal{N}_{123}(\alpha_1,\alpha_3)))$
is isomorphic to the disjoint union of the spaces of type
(\ref{eq1617}),  (\ref{eq1618}) or
\begin{equation}\label{eq1619}
\aligned
&\mathcal M^1({\alpha_1,\alpha_1^2}) \times_{R^1_{\alpha_1^2}}
\dots \times_{R^1_{\alpha_1^{k_1-1}}}
\mathcal M^1({\alpha_1^{k_1-1},\alpha_1^{k_1}})
\\
&\times_{R^1_{\alpha_1^{k_1}}}\widehat{S}^{\frak C}_{n}
\left(\widehat{S}_m^{\frak C}\left(\mathcal N_{12}
(\alpha_1^{k_1},\alpha_2) \times^{\boxplus\tau}_{R^2_{\alpha_2}}
\mathcal N_{23}(\alpha_2,\alpha_3^{1})\right)\right) \\
&\times_{R^3_{\alpha_3^{1}}}
\mathcal M^3({\alpha_3^{1},\alpha_3^{2}}) \times_{R^3_{\alpha_3^{2}}}
\dots \times_{R^3_{\alpha_3^{k_3-1}}}
\mathcal M^3({\alpha_3^{k_3-1},\alpha_3}).
\endaligned
\end{equation}
The covering map
$\widehat{S}_{\ell}(\widehat{S}_k(\mathcal{N}_{123}(\alpha_1,\alpha_3)))
\to \widehat{S}_{\ell+k}(\mathcal{N}_{123}(\alpha_1,\alpha_3)))$
is the identity map on the component of type
(\ref{eq1617}),  (\ref{eq1618}) and is induced by
$$
\aligned
&\widehat{S}^{\frak C}_{n}\left(\widehat{S}_m^{\frak C}\left(\mathcal N_{12}
(\alpha_1^{k_1},\alpha_2) \times^{\boxplus\tau}_{R^2_{\alpha_2}}
\mathcal N_{23}(\alpha_2,\alpha_3^{1})\right)\right)
\\
&\longrightarrow
\widehat{S}_{n+m}^{\frak C}\left(\mathcal N_{12}
(\alpha_1^{k_1},\alpha_2) \times^{\boxplus\tau}_{R^2_{\alpha_2}}
\mathcal N_{23}(\alpha_2,\alpha_3^{1})\right)
\endaligned$$
on the component of type
(\ref{eq1619}).
\end{enumerate}
\end{conds}
\begin{lemdef}\label{1638defken}
By Proposition \ref{prop1636} we can partially smooth the corner
of $\mathcal{N}_{123}(\alpha_1,\alpha_3)$.
Then we obtain a K-space which is the union of
$
\mathcal N_{12}(\alpha_1,\alpha_2) \times^{\boxplus\tau}_{R^2_{\alpha_2}}
\mathcal N_{23}(\alpha_2,\alpha_3)
$
over various $\alpha_2$:
\begin{equation}\label{form16200}
\bigcup_{\alpha_2}
\mathcal N_{12}(\alpha_1,\alpha_2) \times^{\boxplus\tau}_{R^2_{\alpha_2}}
\mathcal N_{23}(\alpha_2,\alpha_3).
\end{equation}
There is a morphism from $\mathcal F_1$ to $\mathcal F_3$
whose interpolation space is
(\ref{form16200}).
We call this morphism the {\rm composition of }
$\frak N_{21}$ and $\frak N_{32}$ and write
$\frak N_{32} \circ \frak N_{21}$.
\end{lemdef}
This follows immediately from Proposition \ref{prop1636}.
\begin{rem}
Before proving Proposition \ref{prop1636},
we note the following.
Using the moduli space of the solutions of
the equation (\ref{equa165}), which uses one parameter family
of homotopies of Hamiltonians
(\ref{eq1666}), we can find a space
$\widehat{\mathcal{N}}_{123}(\alpha_1,\alpha_3)$
whose normalized boundary is the union of
\begin{equation}\label{prop1636form1rev}
\mathcal N_{12}(\alpha_1,\alpha_2) \times_{R^2_{\alpha_2}}
\mathcal N_{23}(\alpha_2,\alpha_3),
\end{equation}
and
\begin{equation}\label{prop1636form2rev}
\mathcal M^1({\alpha_1,\alpha'_1})
\times_{R^1_{\alpha'_1}} \widehat{\mathcal{N}}_{123}(\alpha'_1,\alpha_3)
\end{equation}
over various $\alpha'_1 \in \frak A_1$ and
\begin{equation}\label{prop1636form3rev}
\widehat{\mathcal{N}}_{123}(\alpha_1,\alpha'_3)
\times_{R_{\alpha'_3}} \mathcal M^3({\alpha'_3,\alpha_3})
\end{equation}
over various $\alpha'_3 \in \frak A_3$.
\par
We take $\frak C$ as the boundary of type (\ref{prop1636form1rev}).
Then we can take
$$
{\mathcal{N}}_{123}(\alpha_1,\alpha_3)
=
\widehat{\mathcal{N}}_{123}(\alpha_1,\alpha_3)^{\frak C\boxplus\tau}
\setminus
\widehat{\mathcal{N}}_{123}(\alpha_1,\alpha_3).
$$
\end{rem}
In the proof of Proposition \ref{prop1636} below we will
construct this space without using the whole
${\mathcal{N}}_{123}(\alpha_1,\alpha_3)$ but
using only its `boundary' which is (\ref{prop1636form1rev}),
(\ref{prop1636form2}) and (\ref{prop1636form3}).
This argument is very much similar to the proof of
Proposition \ref{prop528}.

\begin{proof}[Proof of Proposition \ref{prop1636}]
We begin with proving Lemma \ref{lem1639}.
For $A \subset \{1,\dots,k\}$
we regard
$$
[0,1)^A =
\{(t_1,\dots,t_k) \in [0,1)^k \mid t_i = 0
\,\,\, \text{for $i \notin A$}\}.
$$
\begin{lem}\label{lem1639}
For $A \subset B \subseteq \{1,\dots,k\}$
we have a smooth embedding
$$
\phi_{AB} :
([0,1)^{B{^c}})^{\boxplus\tau} \times [-\tau,0)^{\# B}
\to ([0,1)^{A^c})^{\boxplus\tau}\times [-\tau,0)^{\# A},
$$
where $A^c, B^c$ are the complements of $A,B$ in the set $\{1,\dots,k\}$ respectively.
Moreover, if $A \subset B \subset C \subseteq \{1,\dots,k\}$,
we have
\begin{equation}\label{form162727}
\phi_{AB} \circ \phi_{BC}
= \phi_{AC}.
\end{equation}
\end{lem}
\begin{proof}
Regrading $V^+_{\frak r,S_A}, V^+_{\frak r,S_B}$ in \eqref{form1528rev}
as $[0,1)^{A^c}, [0,1)^{B^c}$, respectively,
we can adopt \eqref{form1528rev} to define the map $\phi_{AB}$ above.
(Use the obvious inclusion map and the local inverse of the covering map.)
For example, if $A = \{1,\dots,a\}$,
$B = \{1,\dots,b\}$,
then we have
$$
\phi_{AB}((t_{b+1},\dots,t_k),(s_1,\dots,s_{b}))
=
(s_{a+1},\dots,s_b,t_{b+1},\dots,t_k,s_1,\dots,s_{a}).
$$
The formula (\ref{form162727}) is easy to prove
from this definition.
\end{proof}
For $\alpha_2^1,\dots,\alpha_2^k \in \frak A_2$,
we put
\begin{equation}\label{form1628}
\aligned
&\mathcal N(\alpha_1,\alpha_3;\alpha_2^1,\dots,\alpha_2^k)\\
&=
\left(
\mathcal N_{12}(\alpha_1,\alpha_2^1)
\times_{R^2_{\alpha_2^1}} \mathcal M^2(\alpha_2^1,\alpha_2^2)
\times_{R^2_{\alpha_2^{2}}} \right.\\
&\quad\left.\dots \times_{R^2_{\alpha_2^{k-1}}}
\mathcal M^2(\alpha_2^{k-1},\alpha_2^k)
\times_{R^2_{\alpha_2^{k}}} \mathcal N_{23}(\alpha_2^{k},\alpha_3)\right)^{\frak C
\boxplus\tau} \times [-\tau,0)^k.
\endaligned
\end{equation}
Here $\frak C$ is the decomposition
of the boundary such that its {\it complement} consists of
$$
\aligned
&\mathcal M^1(\alpha_1,\alpha'_1)
\times_{R^1_{\alpha'_1}} \mathcal N_{12}(\alpha'_1,\alpha_2^1)
\times_{R^2_{\alpha_2^1}} \mathcal M^2(\alpha_2^1,\alpha_2^2)
\times_{R^2_{\alpha_2^{2}}}\\
&\quad\quad\quad \dots \times_{R^2_{\alpha_2^{k-1}}}
\mathcal M^2(\alpha_2^{k-1},\alpha_2^k)
\times_{R^2_{\alpha_2^{k}}} \mathcal N_{23}(\alpha_2^{k},\alpha_3)
\endaligned
$$
and
$$
\aligned
&\mathcal N_{12}(\alpha_1,\alpha_2^1)
\times_{R^2_{\alpha_2^1}} \mathcal M^2(\alpha_2^1,\alpha_2^2)
\times_{R^2_{\alpha_2^{2}}}\\
&\quad \dots \times_{R^2_{\alpha_2^{k-1}}}
\mathcal M^2(\alpha_2^{k-1},\alpha_2^k)
\times_{R^2_{\alpha_2^{k}}} \mathcal N_{23}(\alpha_2^{k},\alpha'_3)
\times_{R^3_{\alpha'_3}} \mathcal M^3(\alpha'_3,\alpha_3).
\endaligned
$$
\par
Let $A = \{i_1,\dots,i_a\} \subset \{1,\dots,k\}$ with
$i_1 < i_2 < \dots < i_a$.
We define an embedding
$$
\hat\phi_{A\{1,\dots,k\}} :
\mathcal N(\alpha_1,\alpha_3;\alpha_2^{i_1},\dots,\alpha_2^{i_a})
\to
\mathcal N(\alpha_1,\alpha_3;\alpha_2^1,\dots,\alpha_2^k)
$$
as follows.
For any $p \in \mathcal N(\alpha_1,\alpha_3;\alpha_2^1,\dots,\alpha_2^k)$ we have
a Kuranishi neighborhood $U_p$ of
$\mathcal N(\alpha_1,\alpha_3;\alpha_2^1,\dots,\alpha_2^k)$
such that $U_p$ has an
orbifold chart at $p$ of the form
\begin{equation}\label{16282828}
V_p^{\frak C\boxplus\tau} \times ([-\tau,0)^k).
\end{equation}
On the other hand, $p$ has a Kuranishi neighborhood $U'_p$
of $\mathcal N(\alpha_1,\alpha_3;\alpha_2^{i_1},\dots,\alpha_2^{i_a})$
such that $U'_p$ has an orbifold chart at $p$ of the form
\begin{equation}\label{16292929}
V_p^{\frak C\boxplus\tau} \times ([0,1)^{k-\#A})^{\boxplus\tau}
\times  [-\tau,0)^{\#A}.
\end{equation}
Therefore the map $\hat\phi_{A\{1,\dots,k\}}$ in Lemma \ref{lem1639} defines an
embedding from (\ref{16282828}) to (\ref{16292929}).
This is compatible with the coordinate change
of Kuranishi structure and defines a required embedding.
We can glue various
$\mathcal N(\alpha_1,\alpha_3;\alpha_2^1,\dots,\alpha_2^k)$
by the embeddings  $\hat\phi_{A\{1,\dots,k\}}$.
More precisely, we will glue charts of their Kuranishi structures.
We can do so in the same way as the proof of
Proposition \ref{prop528} using the compatibility of
interpolation spaces and (\ref{form162727}).
\par
We have thus obtained a $\frak C$-collared K-space
$\mathcal{N}_{123}(\alpha_1,\alpha_3)$.
It is easy to see from the construction that
$\mathcal{N}_{123}(\alpha_1,\alpha_3)$ has the
required properties.
(See Figure \ref{Figure17-40}.)
\end{proof}
By construction we have the following:
\begin{lem}\label{lemma1642}
Suppose that there exists a K-space
$$
\mathcal{N}'_{123}(\alpha_1,\alpha_3)
$$
and its boundary components $\frak C$ such that the following
holds.
\begin{enumerate}
\item
$\partial_{\frak C}\mathcal{N}'_{123}(\alpha_1,\alpha_3)$
is a disjoint union of
$\mathcal N_{12}(\alpha_1,\alpha_2) \times_{R^2_{\alpha_2}}
\mathcal N_{23}(\alpha_2,\alpha_3)$
over $\alpha_2$.
\item
The $k$-th normalized corner
$\widehat S^{\frak C}_k(\mathcal{N}'_{123}(\alpha_1,\alpha_3))$
is identified with the disjoint union of
$$
\aligned
&\mathcal N_{12}(\alpha_1,\alpha_2^1) \times_{R^2_{\alpha_2^1}}
\mathcal M^2(\alpha_2^1,\alpha_2^2) \times_{R^2_{\alpha_2^2}}\\
&\dots
\times_{R^2_{\alpha_2^{k-1}}}
\mathcal M^2(\alpha_2^{k-1},\alpha_2^k )
\times_{R^2_{\alpha_2^k}}
\mathcal N_{23}(\alpha_2^k,\alpha_3)
\endaligned
$$
over $\alpha_2^1, \dots , \alpha_2^{k}$.
\item
The map $\hat S_{k}(\hat S_{\ell}(\mathcal{N}'_{123}(\alpha_1,\alpha_2)))
\to \hat S_{k+\ell}(\mathcal{N}'_{123}(\alpha_1,\alpha_2))$
is compatible with the isomorphisms (1)(2) in the sense similar to Condition \ref{conds1637} (2).
\end{enumerate}
Then there exists an isomorphism
\begin{equation}
\mathcal{N}'_{123}(\alpha_1,\alpha_2)^{\frak C\boxplus\tau}
\setminus \mathcal{N}'_{123}(\alpha_1,\alpha_3)
\cong \mathcal{N}_{123}(\alpha_1,\alpha_3).
\end{equation}
This isomorphism is compatible with the isomorphisms in (1)(2) above and
Condition \ref{conds1637} (2) at the corners.
\end{lem}

\subsection{Associativity of the composition}
\label{subsec:compassoci}

In this subsection we present the detail of the proof of the
associativity of the composition.

\begin{prop}\label{prop1640}
Suppose we are in Situation \ref{shitu1634} for $i=1,2,3$.
Then we have the following identity.
\begin{equation}\label{form1630form}
(\frak N_{43} \circ \frak N_{32}) \circ \frak N_{21}
=
\frak N_{43} \circ (\frak N_{32} \circ \frak N_{21}).
\end{equation}
\end{prop}
\begin{proof}
As we already discussed in Subsection \ref{subsec:compmorline},
this equality is mostly obvious.
Namely the interpolation space of the left hand
and the right hand sides
are isomorphic to each other by the associativity of the fiber product.
The only issue is about the way how we smooth the corner.
Since smoothing corners is not completely canonical,
this point is non-trivial.
To clarify this tiny technicality we will use the
following.\footnote{The fact that the left and right hand sides induce the
same cochain map in the de Rham complex by the smooth correspondence
follows without comparing the smooth structures near the corner.
Since the part we smooth the corner
lies in the collar, it does not contribute to the integration along the fiber
(see (\ref{form1515}).)
So we never need the part of the proof of this proposition given in
this subsection for applications.}
\begin{lem}\label{lem1641}
There exists a K-space
$$
\mathcal{N}_{1234}(\alpha_1,\alpha_4)
$$
with the following properties.
\begin{enumerate}
\item
Its normalized boundary is isomorphic to the union of the following four types of
fiber products.
\begin{equation}\label{form1631}
\mathcal{N}_{123}(\alpha_1,\alpha_3)
\times_{R^3_{\alpha_3}}
\mathcal{N}_{34}(\alpha_3,\alpha_4),
\end{equation}
\begin{equation}\label{form1632}
\mathcal{N}_{12}(\alpha_1,\alpha_2)
\times_{R^2_{\alpha_2}}
\mathcal{N}_{234}(\alpha_2,\alpha_4),
\end{equation}
\begin{equation}\label{form1633}
\mathcal M^1(\alpha_1,\alpha'_1)
\times_{R^1_{\alpha'_1}}
\mathcal{N}_{1234}(\alpha_1',\alpha_4),
\end{equation}
\begin{equation}
\mathcal{N}_{1234}(\alpha_1,\alpha'_4)
\times_{R^4_{\alpha'_4}}
\mathcal M^4(\alpha'_4,\alpha'_1).
\end{equation}
\item
A similar compatibility condition as Condition \ref{conds1637}
holds at the corners.
\item
The evaluation maps, periodicity and orientation isomorphisms
are defined on $\mathcal{N}_{1234}(\alpha_1,\alpha_4)$
and commute with the isomorphisms in (1)(2) above.
\item
An appropriate version of Condition \ref{morphilinsys} (V)(IX) holds.
\end{enumerate}
\end{lem}
We omit the precise definition of the compatibility condition
in Lemma \ref{lem1641} (2).
Since it is not difficult, we leave it to the readers.
We will prove Lemma \ref{lem1641} later in this subsection.
We continue the proof of Proposition \ref{prop1640}.
\par
We consider two decompositions
$\frak C_1$, $\frak C_2$ of the boundary of
$\mathcal{N}_{1234}(\alpha_1,\alpha_4)$, where
$\frak C_1$ consists of the boundary of type
(\ref{form1631}) and
$\frak C_2$ consists of the boundary of type
(\ref{form1632}).
Thus we are in the situation of Lemma \ref{lem1614}.
Here we use the next lemma.
\begin{lem}\label{lem1642}
Suppose we are in the situation of Lemma \ref{lem1614}.
Then the following constructions described in {\rm(A)} and {\rm (B)}
give the same K-space.
\begin{enumerate}
\item[(A)]
\begin{enumerate}
\item[(1)]
We first take $\tau$-$\frak C_1$-trivialization of the
corner.
\item[(2)]
We smooth the
boundaries in $\frak C_1$.
\item[(3)]
We next take $\tau$-$\frak C_2$-trivialization of the
corner.
\item[(4)]
We smooth the
boundaries in $\frak C_2$.
\end{enumerate}
\smallskip
\item[(B)]
\begin{enumerate}
\item[(1)]
We first take $\tau$-$\frak C_2$-trivialization of the
corner.
\item[(2)]
We smooth the
boundaries in $\frak C_2$.
\item[(3)]
We next take $\tau$-$\frak C_1$-trivialization of the
corner.
\item[(4)]
We smooth the
boundaries in $\frak C_1$.
\end{enumerate}
\end{enumerate}
\end{lem}
\begin{proof}
It suffices to consider the case of $[0,1)^{k_1} \times [0,1)^{k_2}$
where $\frak C_1$ corresponds to $\partial[0,1)^{k_1} \times [0,1)^{k_2}$
and $\frak C_2$ corresponds to $[0,1)^{k_1} \times \partial[0,1)^{k_2}$.
However, the proof for this case is straightforward from the definition.
\end{proof}
We note that the intersection of
(\ref{form1631}) and (\ref{form1632}) is
\begin{equation}\label{1635mod}
\mathcal{N}_{12}(\alpha_1,\alpha_2)
\times_{R^2_{\alpha_2}}
\mathcal{N}_{23}(\alpha_2,\alpha_3)
\times_{R^3_{\alpha_3}}
\mathcal{N}_{34}(\alpha_3,\alpha_4).
\end{equation}
The interpolation spaces of
both sides of (\ref{form1630form})
are obtained by appropriately
modifying the union of (\ref{1635mod})
for $\alpha_2$, $\alpha_3$.
\par
We study how the processes (A) and (B) of
Lemma \ref{lem1642} affect the subspace
(\ref{1635mod}).
By (A) (1), (\ref{1635mod}) becomes
$$
\mathcal{N}_{12}(\alpha_1,\alpha_2)
\times_{R^2_{\alpha_2}}
\mathcal{N}_{23}(\alpha_2,\alpha_3)
\times^{\boxplus\tau}_{R^3_{\alpha_3}}
\mathcal{N}_{34}(\alpha_3,\alpha_4).
$$
The step (A) (2) smoothes the corners of the
(union over $\alpha_3$ of) the second fiber product
factor. So it becomes
$$
\mathcal{N}_{12}(\alpha_1,\alpha_2)
\times_{R^2_{\alpha_2}}
\mathcal{N}_{24}(\alpha_2,\alpha_4).
$$
Here $\mathcal{N}_{24}(\alpha_2,\alpha_4)$
is the interpolation space of the composition
$\frak N_{43} \circ \frak N_{32}$.
Then (A)(3) changes it to
$$
\mathcal{N}_{12}(\alpha_1,\alpha_2)
\times_{R^2_{\alpha_2}}^{\boxplus\tau}
\mathcal{N}_{24}(\alpha_2,\alpha_4).
$$
Thus when (A) completed
we obtain an interpolation space of
$$
(\frak N_{43} \circ \frak N_{32}) \circ \frak N_{21}.
$$
In the same way, (B) gives an interpolation space of
$$
\frak N_{43} \circ (\frak N_{32} \circ \frak N_{21}).
$$
Thus Proposition \ref{prop1640} follows from
Lemma \ref{lem1642}.
\begin{proof}[Proof of Lemma \ref{lem1641}]
We fix $\alpha_1, \alpha_4$.
Let $\alpha_2^1,\dots,\alpha_2^{k_2} \in \frak A_2$ and
$\alpha_3^1,\dots,\alpha_3^{k_3} \in \frak A_3$.
We put
$$
\aligned
&\mathcal N(\alpha_1,\alpha_3;\alpha_2^1,\dots,\alpha_2^{k_2};
\alpha_3^1,\dots,\alpha_3^{k_3})\\
&=
\left(
\mathcal N_{12}(\alpha_1,\alpha_2^1)
\times_{R^2_{\alpha_2^1}} \mathcal M^2(\alpha_2^1,\alpha_2^2)
\times_{R^2_{\alpha_2^{2}}} \right.\\
&\quad\dots \times_{R^2_{\alpha_2^{k_2-1}}}
\mathcal M^2(\alpha_2^{k_2-1},\alpha_2^{k_2})
\times_{R^2_{\alpha_2^{k_2}}} \mathcal N_{23}(\alpha_2^{k_2},\alpha_3^1)
\times_{R^3_{\alpha_3^{1}}}
\mathcal M^3(\alpha_3^1,\alpha_3^2)\times_{R^3_{\alpha_3^{2}}}\\
&\quad\left.
\dots \times_{R^3_{\alpha_3^{k_3-1}}}
\mathcal M^3(\alpha_3^{k_3-1},\alpha_3^{k_3})
\times_{R^3_{\alpha_3^{k_3}}} \mathcal N_{34}(\alpha_3^{k_3},\alpha_4)
\right)^{\frak C
\boxplus\tau} \times [-\tau,0)^{k_2+k_3}.
\endaligned
$$
Here $\frak C$ is defined in a similar way as in (\ref{form1628}).
We can glue them in a similar way as in the proof of Proposition
\ref{prop1636} to obtain the required
$\mathcal{N}_{1234}(\alpha_1,\alpha_4)$.
\end{proof}
The proof of Proposition \ref{prop1640} is now complete.
\end{proof}
We note that in the geometric situation of
periodic Hamiltonian system
we can prove Lemma \ref{lem1641}
using the 2-parameter family of moduli space of solutions of the
equation
\begin{equation}\label{equa165rev}
\frac{\partial u}{\partial \tau}
+
J
\left(
\frac{\partial u}{\partial t} - X_{H^{ST}_{\tau,t}}(u)
\right)
= 0,
\end{equation}
where
\begin{equation}\label{eq1666rev}
H^{ST}(\tau,t,x)
=
\begin{cases}
H^1(t,x)
&\text{if $\tau \le -T-T_0$}\\
H^{21}(\tau+T,t,x)
&\text{if $-T-T_0\le \tau \le -T+T_0$}\\
H^2(t,x)
&\text{if $-T+T_0 \le \tau \le  -T_0$}\\
H^{32}(\tau,t,x)
&\text{if $-T_0\le \tau \le T_0$}\\
H^3(t,x)
&\text{if $T_0 \le \tau \le  S-T_0$}\\
H^{43}(\tau -S,t,x)
&\text{if $S-T_0 \le \tau \le S+T_0$} \\
H^4(t,x)
&\text{if $S+T_0 \le \tau$},
\end{cases}
\end{equation}
and
$H^{ij}(\tau,t,x) : \R \times S^1 \times M \to \R$ $(i,j =1,\dots,4)$ is defined
by \eqref{def:Hij}.

\subsection{Parametrized version of morphism: composition and gluing}
\label{subsec:paracompositions}

\subsubsection{Compositions of parametrized morphisms}
\label{compparamor}

\begin{shitu}\label{shitu1634para}
Suppose we are in Situation \ref{composisitu} and
$\frak N_{i+1 i}^{P_i}$ is a $P_i$-parametrized morphism from $\mathcal F_i$ to
$\mathcal F_{i+1}$.
We denote by
$$\mathcal M^i(\alpha_-,\alpha_+)$$
the space of connecting orbits of $\mathcal F_i$
and by
$$\mathcal N_{i i+1}(\alpha_-,\alpha_+;P_i)$$
the interpolation space of $\frak N_{i+1 i}^{P_i}$.
Let $R_{\alpha_i}^i$ be a critical submanifold
of $\mathcal F_i$ and $\alpha_i \in \frak A_i$.
$\blacksquare$
\end{shitu}
\begin{defn}\label{defn1635para}
In Situation \ref{shitu1634para} we define a K-space
\begin{equation}\label{form1616para}
\mathcal N_{12}(\alpha_1,\alpha_2;P_1) \times^{\boxplus\tau}_{R^2_{\alpha_2}}
\mathcal N_{23}(\alpha_2,\alpha_3;P_2)
\end{equation}
as follows.
We take the
fiber product
$
\mathcal N_{12}(\alpha_1,\alpha_2;P_1) \times_{R^2_{\alpha_2}}
\mathcal N_{23}(\alpha_2,\alpha_3;P_2)
$
with fiber product Kuranishi structure and
consider the decomposition $\frak C$ of its boundary
consisting of the following two kinds of
components of its normalized boundary:
$$
\aligned
&\mathcal N_{12}(\alpha_1,\alpha'_2;P_1) \times_{R^2_{\alpha'_2}}
\mathcal M^2(\alpha_{\alpha'_2},\alpha_{\alpha_2})
\times_{R^2_{\alpha_2}}
\mathcal N_{23}(\alpha_2,\alpha_3;P_2), \\
&\mathcal N_{12}(\alpha_1,\alpha_2;P_1) \times_{R^2_{\alpha_2}}
\mathcal M^2(\alpha_{\alpha_2},\alpha_{\alpha'_2})
\times_{R^2_{\alpha'_2}}
\mathcal N_{23}(\alpha'_2,\alpha_3;P_2).
\endaligned
$$
The first line is contained in
$
\partial\mathcal N_{12}(\alpha_1,\alpha_2;P_1) \times_{R^2_{\alpha_2}}
\mathcal N_{23}(\alpha_2,\alpha_3;P_2)
$
and the second line is contained in
$
\mathcal N_{12}(\alpha_1,\alpha_2;P_1) \times_{R^2_{\alpha_2}}
\partial\mathcal N_{23}(\alpha_2,\alpha_3;P_2)
$.
Now we define (\ref{form1616para}) by
$$
(\mathcal N_{12}(\alpha_1,\alpha_2;P_1) \times_{R^2_{\alpha_2}}
\mathcal N_{23}(\alpha_2,\alpha_3);P_2)^{\frak C\boxplus\tau}.
$$
Note that besides the evaluation maps ${\rm ev}_{-}$ and ${\rm ev}_{+}$,
(\ref{form1616para}) carries an evaluation map to $P_1 \times P_2$,
which is stratumwise weakly submersive.
\end{defn}
\begin{prop}\label{prop1636para}
In Situation \ref{shitu1634para} there exists a
$\tau$-$\frak C$-collared K-space
$$
\mathcal{N}_{123}(\alpha_1,\alpha_3;P_1\times P_2)
$$
with the following properties:
\begin{enumerate}
\item
Its normalized boundary in $\frak C$ is isomorphic to the disjoint union of
\begin{equation}\label{prop1636form1para}
\mathcal N_{12}(\alpha_1,\alpha_2;P_1) \times^{\boxplus\tau}_{R^2_{\alpha_2}}
\mathcal N_{23}(\alpha_2,\alpha_3;P_2)
\end{equation}
over various $\alpha_2$.
\item
Its normalized boundary not in $\frak C$ is isomorphic to the disjoint union of
\begin{equation}\label{prop1636form2para}
\mathcal M^1({\alpha_1,\alpha'_1})
\times_{R^1_{\alpha'_1}} \mathcal{N}_{123}(\alpha'_1,\alpha_3;P_1\times P_2)
\end{equation}
over various $\alpha'_1$ and
\begin{equation}\label{prop1636form3para}
\mathcal{N}_{123}(\alpha_1,\alpha'_3;P_1\times P_2)
\times_{R_{\alpha'_3}} \mathcal M^3({\alpha'_3,\alpha_3})
\end{equation}
over various $\alpha'_3$.
\item
Conditions similar to
(3)(4)(5) of Proposition \ref{prop1636} hold.
\item
We have an evaluation map
$
\mathcal{N}_{123}(\alpha_1,\alpha_3;P_1\times P_2) \to P_1 \times P_2
$
which is a stratumwise weakly submersive map.
It is compatible with the boundary description given in (1)(2) above.
\end{enumerate}
\end{prop}
The proof is the same as the proof of Proposition \ref{prop1636} so omitted.
\begin{defnlem}\label{1638defkenpara}
By Proposition \ref{prop1636para} we can partially smooth the corner
of $\mathcal{N}_{123}(\alpha_1,\alpha_3;P_1\times P_2)$
to obtain a K-space given by
\begin{equation}\label{form16200para}
\bigcup_{\alpha_2}
\mathcal N_{12}(\alpha_1,\alpha_2;P_1) \times^{\boxplus\tau}_{R^2_{\alpha_2}}
\mathcal N_{23}(\alpha_2,\alpha_3;P_2).
\end{equation}
There is a $(P_1\times P_2)$-parametrized morphism from $\mathcal F_1$ to $\mathcal F_3$
whose interpolation space is given by
(\ref{form16200para}).
We call this morphism the {\it parametrized composition of }
$\frak N_{21}^{P_1}$ and $\frak N_{32}^{P_2}$ and write
$\frak N_{32}^{P_1} \circ \frak N_{21}^{P_2}$.
\end{defnlem}
\begin{lem}
The parametrized composition is associative in the same sense as in
Proposition \ref{prop1640}.
Moreover the boundary of the parametrized composition
$\frak N_{32}^{P_1} \circ \frak N_{21}^{P_2}$ is
the disjoint union of
$$
\frak N_{32}^{\partial P_1} \circ \frak N_{21}^{P_2}
\qquad \text{and} \qquad
\frak N_{32}^{P_1} \circ \frak N_{21}^{\partial P_2}.
$$
\end{lem}
\begin{proof}
The proof of the first half is similar to the proof of Proposition
\ref{prop1640}. The second half is obvious from the construction.
\end{proof}

\subsubsection{Gluing parametrized morphisms}
\label{glueinterpolation}

We first review well known obvious facts of gluing cornered manifolds
along their boundaries.
\begin{defnlem}\label{gluemfdaroungvoundary}
Let $P_1$, $P_2$ be two admissible manifolds with corners,
and for each $i=1,2$
let $\partial _{\frak C_i}P_i \subset \partial P_i$ be an open and closed subset of the
normalized boundary $\partial P_i$ and
assume $S_2^{\frak C_i}(P_i) =\emptyset$.
(See Definition \ref{defnSCk} for the notation.)
Let $I : \partial _{\frak C_1}P_1 \to \partial _{\frak C_2}P_2$
be an admissible diffeomorphism of manifolds with corners.
\begin{enumerate}
\item
We can define a structure of admissible manifold with corners
on
\begin{equation}\label{glueP1andP2}
P_1  \,_{\frak C_1} \cup_{\frak C_2} P_2
=
(P_1 \cup P_2)/\sim .
\end{equation}
Here $\sim$ is defined by $x \sim I(x)$ for $x \in \partial _{\frak C_1}P_1$.
\item
The boundary
$
\partial(P_1  \,_{\frak C_1} \cup_{\frak C_2} P_2)
$
is described as follows.
Let $\partial _{\frak C^c_i}P_i$ be the complement of $\partial _{\frak C_i}P_i$
in $\partial P_i$. Then
$\frak C_i$ induces
a decomposition
$$
\partial(\partial_{\frak C^c_i}P_i)
=
\partial_{\frak C_i}(\partial_{\frak C^c_i}P_i)
\cup
\partial_{\frak C^c_i}(\partial_{\frak C^c_i}P_i).
$$
Moreover $I$ induces a diffeomorphism
$$
\partial_{\frak C_1}(\partial_{\frak C^c_1}P_1)
\cong
\partial_{\frak C_2}(\partial_{\frak C^c_2}P_2).
$$
Now we have, see Figure \ref{Figure17-35},
\begin{equation}
\partial(P_1  \,_{\frak C_1} \cup_{\frak C_2} P_2)
=
\partial_{\frak C^c_1}P_1
\,\,{}_{\frak C_1} \cup_{\frak C_2} \,\, \partial _{\frak C^c_2}P_2.
\end{equation}
\item
Suppose that we have one of the following objects on $P_1$ and on $P_2$
such that their restrictions to $\partial_{\frak C_1}P_1$ and
to $\partial_{\frak C_2}P_2$ are isomorphic each other
(or coincide with each other)
when we identify  $\partial_{\frak C_1}P_1$ with
$\partial_{\frak C_2}P_2$ under the diffeomorphism $I$.
Then we obtain a glued object on
$P_1
\,\,{}_{\frak C_1} \cup_{\frak C_2} P_2$.
\begin{enumerate}
\item
Vector bundle.
\item
Section of vector bundle.
\item
Differential form.
\item
Smooth map to a manifold.
\end{enumerate}
\item
In addition, suppose that we have admissible manifolds with corners $P'_1$, $P'_2$ and
for each $i=1,2$ let $\partial_{\frak C'_i}P'_i$ be an open and closed subset of
the normalized boundary $\partial P'_i$. Let
$I' : \partial_{\frak C'_1}P'_1 \cong \partial_{\frak C'_2}P'_2$ be
an admissible diffemorphism as above.
We also assume that for each $i=1,2$ we have
an admissible smooth embedding
$f_i : P'_i \to P_i$ of cornered orbifolds such that
\begin{enumerate}
\item
$f_i^{-1}(\partial_{\frak C_i}P_i) = \partial_{\frak C'_i}P'_i$.
\item
$$
I \circ f_1 = f_2 \circ I'
$$
on $\partial_{\frak C'_1}P'_1$.
\end{enumerate}
Then $f_1,f_2$ induce a smooth admissible embedding of cornered orbifolds
$$
f_1 \,\,{}_{\frak C_1} \cup_{\frak C_2} \,\, f_2
:
P'_1\,\,{}_{\frak C'_1} \cup_{\frak C'_2} \,\,P'_2
\to
P_1\,\,{}_{\frak C_1} \cup_{\frak C_2} \,\,P_2.
$$
\item
In the situation of (4), suppose that we are given vector bundles $E_i$ on $P_i$ and
$E'_i$ on $P'_i$. We assume $E_1\vert_{\partial_{\frak C_1}P_1}
\cong E_2\vert_{\partial_{\frak C_2}P_2}$ and the isomorphism covers $I$.
We assume the same condition for $E'_1$ and $E'_2$.
Moreover we assume that there exist embeddings of vector bundles $\hat f_i : E'_i \to E_i$
which cover $f_i$ for $i=1,2$, and
assume that $\hat f_1$ and $\hat f_2$ are compatible with
the isomorphisms $E_1\vert_{\partial_{\frak C_1}P_1}
\cong E_2\vert_{\partial_{\frak C_2}P_2}$
and $E'_1\vert_{\partial_{\frak C'_1}P'_1}
\cong E'_2\vert_{\partial_{\frak C'_2}P'_2}$.
\par
Then we obtain an embedding of vector bundles:
$$
\hat f_1 \,\,{}_{\frak C_1} \cup_{\frak C_2} \,\, \hat f_2
:
E'_1\,\,{}_{\frak C'_1} \cup_{\frak C'_2} \,\,E'_2
\to
E_1\,\,{}_{\frak C_1} \cup_{\frak C_2} \,\,E_2.
$$
\item
In addition, suppose that we have $P'_1$, $P'_2$ and
$\partial_{\frak C'_i}P'_i$ (for $i=1,2$),
$I' : \partial_{\frak C'_1}P'_1 \cong \partial_{\frak C'_2}P'_2$ as above.
We also assume that for each $i=1,2$ we have an
admissible smooth stratumwise submersion
$\pi_i : P'_i \to P_i$ of cornered orbifolds such that
\begin{enumerate}
\item
$\pi_i^{-1}(\partial_{\frak C_i}P_i) = \partial_{\frak C'_i}P'_i$.
\item
$$
I \circ \pi_1 = \pi_2 \circ I'
$$
on $\partial_{\frak C'_1}P'_1$.
\end{enumerate}
Then $\pi_1,\pi_2$ induce an admissible smooth stratumwise submersion of cornered orbifolds
$$
\pi_1 \,\,{}_{\frak C_1} \cup_{\frak C_2} \,\, \pi_2
:
P'_1\,\,{}_{\frak C'_1} \cup_{\frak C'_2} \,\,P'_2
\to
P_1\,\,{}_{\frak C_1} \cup_{\frak C_2} \,\,P_2.
$$
\end{enumerate}
\end{defnlem}
\begin{figure}[h]
\centering
\includegraphics[scale=0.4]{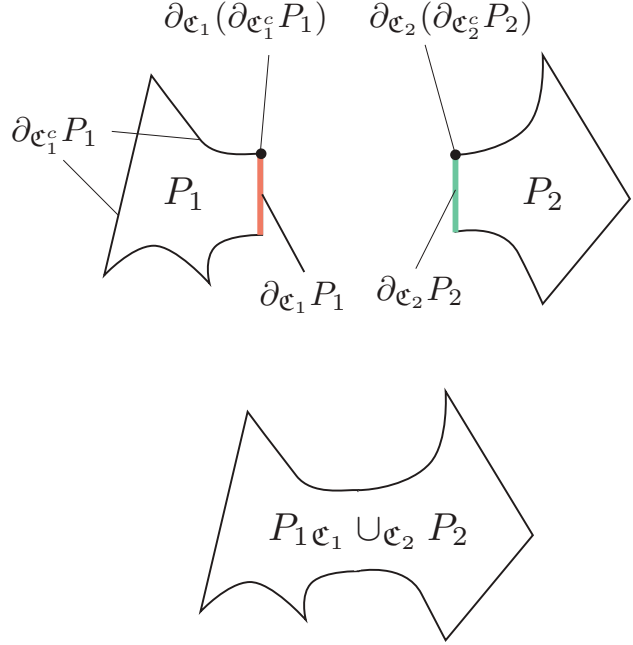}
\caption{gluing cornered manifolds}
\label{Figure17-35}
\end{figure}
\begin{proof}
The proof is mostly obvious except the smoothness of
various objects (such as coordinate change) at
the boundary where we glue $P_1$ and $P_2$.
The smoothness can be proved by using admissibility and
Lemma \ref{newlem319} (2).
In fact, the differential in the direction normal to the boundary of all the
objects vanishes in infinite order.
So gluing along the boundary gives smooth objects.
\end{proof}
\begin{shitu}\label{situ17533}
Let
$P_1$, $P_2$, $\partial _{\frak C_i}P_i \subset \partial P_i$
and $I : \partial _{\frak C_1}P_1 \to \partial _{\frak C_2}P_2$
be given as in Definition-Lemma \ref{gluemfdaroungvoundary}.
For $i=1,2$ let $(X_i,\widehat{\mathcal U_i})$ be a K-space and
${\rm ev}_{P_i} : (X_i,\widehat{\mathcal U_i}) \to P_i$
be a strongly smooth and stratumwise weakly submersive map.
Via the map ${\rm ev}_{P_i} : (X_i,\widehat{\mathcal U_i}) \to P_i$,
the decomposition
$\partial P_i = \partial _{\frak C_i}P_i \cup \partial_{\frak C_i^c}P_i$
induces a decomposition
$\partial(X_i,\widehat{\mathcal U_i})
= \partial_{\frak C_i}(X_i,\widehat{\mathcal U_i})
\cup \partial_{\frak C_i^c}(X_i,\widehat{\mathcal U_i})$
as in Situation \ref{decomporbbdrkura}.
Namely,
$\partial_{\frak C_i}(X_i,\widehat{\mathcal U_i})$
is defined by taking Kuranishi charts so that they are mapped
to $\partial _{\frak C_i}P_i$ via ${\rm ev}_{P_i}$ and
the coordinate change satisfies the condition as in Situation \ref{decomporbbdrkura}.
We assume that there exists an isomorphism
\begin{equation}
\hat I :
\partial_{\frak C_1}(X_1,\widehat{\mathcal U_1})
\cong
\partial_{\frak C_2}(X_2,\widehat{\mathcal U_2})
\end{equation}
which is compatible with the admissible diffeomorphism
$I :
\partial _{\frak C_1}P_1 \to \partial _{\frak C_2}P_2$.
\par
Then we glue $X_1$ and $X_2$ by the underlying homeomorphism
of $\hat I$ to obtain
$$
X = X_1 \,\,{}_{\frak C_1} \cup {}_{\frak C_2}\,\, X_2
$$
in the same way as (\ref{glueP1andP2}).
$\blacksquare$
\end{shitu}
\begin{defnlem}\label{glueKuraaroungvoundary}
In Situation \ref{situ17533} we can glue $\widehat{\mathcal U_1}$ and
$\widehat{\mathcal U_2}$ to obtain a
Kuranishi structure $\widehat{\mathcal U}$ of $X$.
We write
$$
(X,\widehat{\mathcal  U}) =(X_1,\widehat{\mathcal U_1}) \,\,{}_{\frak C_1} \cup {}_{\frak C_2}\,\,
(X_2,\widehat{\mathcal U_2}), \qquad
\widehat{\mathcal U} = \widehat{\mathcal U_1} \,\,{}_{\frak C_1} \cup {}_{\frak C_2}\,\,
\widehat{\mathcal U_2}.
$$
It has the following properties.
\begin{enumerate}
\item
We have
$$
\partial(X,\widehat{\mathcal U})
=\partial _{\frak C^c_1}(X_1,\widehat{\mathcal U_1})
\,\,{}_{\frak C_1} \cup_{\frak C_2} \,\, \partial_{\frak C^c_2}(X_2,\widehat{\mathcal U_2}).
$$
\item
Suppose that we have one of the following objects on $(X_1,\widehat{\mathcal U_1})$ and on $(X_2,\widehat{\mathcal U_2})$
such that their restrictions to $\partial_{\frak C_1}(X_1,\widehat{\mathcal U_1})$ and
to $\partial_{\frak C_2}(X_2,\widehat{\mathcal U_2})$
are isomorphic each other (or coincide each other)
when we identify  $\partial_{\frak C_1}(X_1,\widehat{\mathcal U_1})$ with
$\partial_{\frak C_2}(X_2,\widehat{\mathcal U_2})$ under the diffeomorphism $\hat I$.
Then we obtain the corresponding object on
$(X,\widehat{\mathcal U})$.
\begin{enumerate}
\item
Differential form.
\item
Strongly smooth map to a manifold.
\item
CF-perturbation.
\item
Multi-valued perturbation.
\end{enumerate}
\item
Let $h_i$ ($i=1,2$) be differential forms on $(X_i,\widehat{\mathcal U_i})$,
$f_i : (X_i,\widehat{\mathcal U_i}) \to M$
strongly smooth maps, and $\widehat{\frak S_i}$ CF-perturbations of
$(X_i,\widehat{\mathcal U_i})$ such that $f_i$ is strongly submersive with respect to
$\widehat{\frak S_i}$ respectively. Suppose we can glue them in the sense of Item (2)
and obtain $h$, $f$ and $\widehat{\frak S}$. Then we have
$$
f!(h;\widehat{\frak S^{\epsilon}})
=
f_1!(h_1;\widehat{\frak S_1^{\epsilon}})
+
f_2!(h_2;\widehat{\frak S_2^{\epsilon}}).
$$
\item
In addition, suppose that we have $(X_i,\widehat{\mathcal U'_i})$ ($i=1,2$) and
$\hat I' : \partial_{\frak C_1}(X_1,\widehat{\mathcal U'_1}) \cong \partial_{\frak C_2}(X_2,\widehat{\mathcal U'_2})$ as above.
(Note that both of $\widehat{\mathcal U'_i}$ and $\widehat{\mathcal U_i}$ are
Kuranishi structures of $X_i$ for each $i=1,2$.)
\par
We also assume that we have
embeddings
$(X_i,\widehat{\mathcal U'_i}) \to (X_i,\widehat{\mathcal U_i})$ of Kuranishi
structures for $i=1,2$
that are compatible with $\widehat I$ and $\widehat I'$
in an obvious sense.
Then the embeddings induce an embedding of Kuranishi structures
$
:
\widehat{\mathcal U'_1}\,{}_{\frak C_1} \cup_{\frak C_2} \,\widehat{\mathcal U'_2}
\to
\widehat{\mathcal U_1}\,{}_{\frak C_1} \cup_{\frak C_2} \,\widehat{\mathcal U_2}.
$
\par
The compatibility of objects (a)(b)(c)(d) in Item (2) is preserved by this process.
\item
We have a stratumwise submersive map:
$$
(X_1,\widehat{\mathcal U_1}) \,{}_{\frak C_1} \cup {}_{\frak C_2}\,
(X_2,\widehat{\mathcal U_2})
\to P_1\,{}_{\frak C_1} \cup_{\frak C_2} \,P_2.
$$
\end{enumerate}
\end{defnlem}
The proof is immediate from Definition-Lemma \ref{gluemfdaroungvoundary}.
\begin{shitu}\label{situ1753322}
Suppose that we are in Situation \ref{situparaPmorph}.
Let
$P_1$, $P_2$, $\partial _{\frak C_i}P_i \subset \partial P_i$
and $I : \partial _{\frak C_1}P_1 \to \partial _{\frak C_2}P_2$
be as in Definition-Lemma \ref{gluemfdaroungvoundary}.
Let $\frak N^{P_i}$ is a $P_i$-parametrized morphism from $\mathcal F_1$ to
$\mathcal F_{2}$.
We denote by
$$
\mathcal N_{i}(\alpha_-,\alpha_+;P_i)
$$
the interpolation space of $\frak N_{i}^{P_i}$.
We obtain
a $\partial_{\frak C_i}P_i$-morphism
$\frak N^{\partial_{\frak C_i}P_i}$ from $\mathcal F_1$ to
$\mathcal F_{2}$,
whose interpolation space is
$$
\mathcal N_{i}(\alpha_-,\alpha_+;\partial_{\frak C_i}P_i)
=
\partial_{\frak C_i}P_i \times_{P_i}
\mathcal N_{i}(\alpha_-,\alpha_+;P_i).
$$
We assume that
$\frak N^{\partial_{\frak C_1}P_1}$ is isomorphic to
$\frak N^{\partial_{\frak C_2}P_2}$ under an isomorphism
which covers $I : \partial _{\frak C_1}P_1 \to \partial _{\frak C_2}P_2$
in an obvious sense.
$\blacksquare$
\end{shitu}
\begin{defnlem}\label{glueparamorphi}
In Situation \ref{situ1753322},
we can glue $\frak N^{P_1}$ and $\frak N^{P_2}$
to define a
$P$-parametrized morphism $\frak N^{P}$
for $P=P_1\,{}_{\frak C_1} \cup_{\frak C_2} P_2$
from
$\mathcal F_1$ to
$\mathcal F_{2}$
with the following properties.
\begin{enumerate}
\item
The interpolation space of $\frak N^{P}$ is
$$
\mathcal N_{1}(\alpha_-,\alpha_+;P_1)
\,{}_{\frak C_1} \cup {}_{\frak C_2}\,
\mathcal N_{2}(\alpha_-,\alpha_+;P_2).
$$
\item
The boundary of $\frak N^{P}$ is
$$
\frak N^{P_1}\vert_{\partial_{\frak C^c_1} P_1}
\,{}_{\frak C_1} \cup {}_{\frak C_2}\,
\frak N^{P_1}\vert_{\partial_{\frak C^c_2} P_2}.
$$
\end{enumerate}
\end{defnlem}
The proof is immediate from Definition-Lemma \ref{glueKuraaroungvoundary}.
\par\smallskip
Now we return to the situation of Lemma-Definition \ref{lemdef1437}.
We recall that to prove Item (1) we glued two $[0,1]$-parametrized
morphisms $\frak H_{(32)}^i \circ \frak N^{i}_{21}$ and
$\frak N_{32}^i \circ \frak H^{i}_{(21)}$
along a part of their boundaries, and
to prove Item (2) we glued two
$[0,1]^2$-parametrized morphisms
$\frak N^{i+1}_{32} \circ \mathcal H^i$
and $\frak H_{(32)}^i \circ \frak H_{(ab)}^{i}$
along a part of their boundaries.
Using Definition-Lemma \ref{glueparamorphi}, we can perform such gluing.
\begin{proof}[Proof of Lemma-Definition \ref{lemdef1437} (3)]
We can prove transitivity by gluing two homotopies
by using  Definition-Lemma \ref{glueparamorphi}.
Other properties are easier to prove so omitted.
\end{proof}

\subsection{Identity morphism}
\label{subsec:identitylinsys}
In this subsection, we define the identity morphism.
\begin{rem}
Note that we can {\it not} define the identity morphism by
defining its interpolation space such as
\begin{equation}\label{formrem1747}
\mathcal N(\alpha_-,\alpha_+)
=
\begin{cases}
R_{\alpha_-}   &\alpha_- = \alpha_+, \\
\emptyset & \alpha_- \ne \alpha_+.
\end{cases}
\end{equation}
In fact, for $\alpha_- \ne \alpha_+$, (\ref{formula1211morph}) and
(\ref{formrem1747}) yield
\begin{equation}\label{form1548}
\aligned
&\partial \mathcal N(\alpha_-,\alpha_+)  \\
&\supseteq  (\mathcal M(\alpha_-,\alpha_+) \times_{R_{\alpha_+}} \mathcal N(\alpha_+,\alpha_+))
 \cup (\mathcal N(\alpha_-,\alpha_-) \times_{R_{\alpha_-}} \mathcal M(\alpha_-,\alpha_+)) \\
& = \mathcal M(\alpha_-,\alpha_+) \sqcup \mathcal M(\alpha_-,\alpha_+)
\endaligned
\end{equation}
which is not an empty set in general.
\par
Suppose that we put two different systems of perturbations
on $\mathcal M(\alpha,\beta)$.
They give two different coboundary oprations of the Floer
cochain complex on $CF(\mathcal C)$.
We can use the identity morphism and a system of perturbations on
its interpolation spaces
to define a cochain map between them.
This is the way how we proceed in Section \ref{sec:systemline3}.
While working out this proof, we put two different perturbations on the
two copies on $\mathcal M(\alpha_-,\alpha_+)$
appearing in the right hand side of (\ref{form1548})
and extend them to $\mathcal N(\alpha_-,\alpha_+)$.
Thus the fact that $\mathcal N(\alpha_-,\alpha_+)$ is nonempty is important
for proving that Floer cohomology is independent
of perturbations by using the identity morphism.
\end{rem}
Let
$\mathcal F$ be a linear K-system whose space of connecting orbits is given by
$\mathcal M(\alpha_-,\alpha_+)$.
We will define
the identity morphism from $\mathcal F$ to itself.
The interpolation spaces are defined as follows.
\begin{defn}\label{identitystrati}
\begin{enumerate}
\item
We define $\overset{\circ}{\mathcal N}(\alpha_-,\alpha_+)$
as follows.
\begin{enumerate}
\item
If $\alpha_- \ne \alpha_+$, then
$$
\overset{\circ}{\mathcal N}(\alpha_-,\alpha_+)
=
\overset{\circ}{\mathcal M}(\alpha_-,\alpha_+) \times \R.
$$
\item
If $\alpha = \alpha_- = \alpha_+$, then
$$
\overset{\circ}{\mathcal N}(\alpha,\alpha)
=
R_{\alpha}.
$$
\end{enumerate}
\item
We compactify $\overset{\circ}{\mathcal N}(\alpha_-,\alpha_+)$ as follows.
In case $\alpha = \alpha_- = \alpha_+$, we put
${\mathcal N}(\alpha,\alpha) =
\overset{\circ}{\mathcal N}(\alpha,\alpha)$.
\par
In case $\alpha_- \ne \alpha_+$, a stratum
$\overset{\circ}S_k({\mathcal N}(\alpha_-,\alpha_+))$
of the compactification
${\mathcal N}(\alpha_-,\alpha_+)$
is the union of the following two types of fiber products:
\begin{enumerate}
\item
\begin{equation}\label{1442form}
\aligned
&\overset{\circ}{\mathcal M}(\alpha_-,\alpha_1)
\times_{R_{\alpha_1}}
\overset{\circ}{\mathcal M}(\alpha_1,\alpha_2)
\times_{R_{\alpha_{2}}}
\dots
\times_{R_{\alpha_{k'-1}}}
\overset{\circ}{\mathcal M}(\alpha_{k'-1},\alpha_{k'})  \\
&\times_{R_{\alpha_{k'}}}
(\overset{\circ}{\mathcal M}(\alpha_{k'},\alpha_{k'+1}) \times \R)
\\
&\times_{R_{\alpha_{k'+1}}}
\overset{\circ}{\mathcal M}(\alpha_{k'+1},\alpha_{k'+2})
\times_{R_{\alpha_{k'+2}}}
\dots
\times_{R_{\alpha_{k}}}
\overset{\circ}{\mathcal M}(\alpha_{k},\alpha_{+}).
\endaligned
\end{equation}
Here we take all possible choices of $k'$,
$\alpha_1,\dots,\alpha_k$ with $0 \le k' \le k+1$ and
$$
E(\alpha_-) < E(\alpha_1) < \dots < E(\alpha_{k'}) <
\dots < E(\alpha_k) < E(\alpha_+).
$$
\item
\begin{equation}\label{1443form}
\aligned
&\overset{\circ}{\mathcal M}(\alpha_-,\alpha_1)
\times_{R_{\alpha_1}}
\overset{\circ}{\mathcal M}(\alpha_1,\alpha_2)
\times_{R_{\alpha_{2}}}
\dots
\times_{R_{\alpha_{k'-1}}}
\overset{\circ}{\mathcal M}(\alpha_{k'-1},\alpha_{k'})  \\
&\times_{R_{\alpha_{k'}}}
R_{\alpha_{k'}}
\\
&\times_{R_{\alpha_{k'}}}
\overset{\circ}{\mathcal M}(\alpha_{k'},\alpha_{k'+1})
\times_{R_{\alpha_{k'+1}}}
\dots
\times_{R_{\alpha_{k-1}}}
\overset{\circ}{\mathcal M}(\alpha_{k-1},\alpha_{+}).
\endaligned
\end{equation}
Here we take all possible choices of $k'$,
$\alpha_1,\dots,\alpha_{k-1}$ with $0 \le k' \le k$ and
$$
E(\alpha_-) < E(\alpha_1) < \dots < E(\alpha_{k'}) <
\dots < E(\alpha_{k-1}) < E(\alpha_+).
$$
\end{enumerate}
\end{enumerate}
\end{defn}
\begin{rem}
\begin{enumerate}
\item
In (2)(a) above, we regard $\alpha_{0} = \alpha_-$ and
$\alpha_{k+1} = \alpha_+$.
For example, in the case $k' = 0$, the stratum (\ref{1442form})
becomes
$$
(\overset{\circ}{\mathcal M}(\alpha_{-},\alpha_{1}) \times \R)
\times_{R_{\alpha_{1}}}
\overset{\circ}{\mathcal M}(\alpha_{1},\alpha_{2})
\times_{R_{\alpha_{2}}}
\dots
\times_{R_{\alpha_{k}}}
\overset{\circ}{\mathcal M}(\alpha_{k},\alpha_{+}).
$$
\item
Similarly in (2)(b) above, we regard $\alpha_{0} = \alpha_-$ and
$\alpha_{k} = \alpha_+$.
For example, in the case $k' = k$, the stratum (\ref{1443form})
becomes
\begin{equation}\label{1449form}
\overset{\circ}{\mathcal M}(\alpha_{-},\alpha_{1})
\times_{R_{\alpha_{1}}}
\overset{\circ}{\mathcal M}(\alpha_{1},\alpha_{2})
\times_{R_{\alpha_{2}}}
\dots
\times_{R_{\alpha_{k-1}}}
\overset{\circ}{\mathcal M}(\alpha_{k-1},\alpha_{+})
\times_{R_{\alpha_+}} R_{\alpha_+}.
\end{equation}
\item
Indeed, (\ref{1443form}) is isomorphic to
$$
\aligned
&\overset{\circ}{\mathcal M}(\alpha_-,\alpha_1)
\times_{R_{\alpha_1}}
\overset{\circ}{\mathcal M}(\alpha_1,\alpha_2)
\times_{R_{\alpha_{2}}}
\dots
\times_{R_{\alpha_{k'-1}}}
\overset{\circ}{\mathcal M}(\alpha_{k'-1},\alpha_{k'})  \\
&\times_{R_{\alpha_{k'}}}
\overset{\circ}{\mathcal M}(\alpha_{k'},\alpha_{k'+1})
\times_{R_{\alpha_{k'+1}}}
\dots
\times_{R_{\alpha_{k-1}}}
\overset{\circ}{\mathcal M}(\alpha_{k-1},\alpha_{+}).
\endaligned
$$
We write it as in (\ref{1443form}) to distinguish
the various components, that is actually the same space.
However,
in the definition of $\overset{\circ}S_k({\mathcal N}(\alpha_-,\alpha_+))$ above,
we regard those various components written as (\ref{1443form}) as {\it different spaces}.
\item
Similarly, the space (\ref{1442form}) is independent of $k' = 0,\dots,k+1$
up to isomorphism.
However, we regard them as {\it different spaces} in the definition of
$\overset{\circ}S_k({\mathcal N}(\alpha_-,\alpha_+))$.
\end{enumerate}
\end{rem}
\begin{lemdef}\label{lemdefidentity}
There exists a morphism from $\mathcal F$ to itself
whose interpolation space has a stratification
described in Definition \ref{identitystrati}.
We call the morphism the {\rm identity morphism}\index{identity morphism ! of $\mathcal F$}.
We can also define the identity morphism in the case of partial
linear K-systems in the same way.
\end{lemdef}
Before proving the lemma we give an example.
\begin{exm}\label{exam1443}
Let us consider the case when there is exactly
one $\alpha$ such that
$E(\alpha_-) < E(\alpha) < E(\alpha_+)$.
Then $\mathcal N(\alpha_-,\alpha_+)$
is stratified as follows:
$$
\overset{\circ}S_0(\mathcal N(\alpha_-,\alpha_+))
= \overset{\circ}{\mathcal M}(\alpha_-,\alpha_+) \times \R.
$$
\begin{equation}\label{144444}
\aligned
\overset{\circ}S_1(\mathcal N(\alpha_-,\alpha_+))
=&
(\overset{\circ}{\mathcal M}(\alpha_-,\alpha) \times \R)
\times_{R_{\alpha}}
\overset{\circ}{\mathcal M}(\alpha,\alpha_+) \\
&\cup
\overset{\circ}{\mathcal M}(\alpha_-,\alpha)
\times_{R_{\alpha}}
(\overset{\circ}{\mathcal M}(\alpha,\alpha_+) \times \R)\\
&\cup R_{\alpha_-} \times_{R_{\alpha_-}}\overset{\circ}{\mathcal M}(\alpha_-,\alpha_+)
\\
&\cup \overset{\circ}{\mathcal M}(\alpha_-,\alpha_+)
\times_{R_{\alpha_+}}  R_{\alpha_+}.
\endaligned
\end{equation}
\begin{equation}\label{1445}
\aligned
\overset{\circ}S_2(\mathcal N(\alpha_-,\alpha_+))
=&
R_{\alpha_-} \times_{R_{\alpha_-}}\overset{\circ}{\mathcal M}(\alpha_-,\alpha)
\times_{R_{\alpha}}
\overset{\circ}{\mathcal M}(\alpha,\alpha_+) \\
&\cup
\overset{\circ}{\mathcal M}(\alpha_-,\alpha)
\times_{R_{\alpha}}  R_{\alpha} \times_{R_{\alpha}}
\overset{\circ}{\mathcal M}(\alpha,\alpha_+)
\\
&
\cup
\overset{\circ}{\mathcal M}(\alpha_-,\alpha)
\times_{R_{\alpha}}
\overset{\circ}{\mathcal M}(\alpha,\alpha_+)
\times_{R_{\alpha_+}}  R_{\alpha_+}.
\endaligned
\end{equation}
See Figure \ref{Figure17-4}.
Note that in this case
$\mathcal M(\alpha_-,\alpha_+)$
has a Kuranishi structure with boundary
$\overset{\circ}{\mathcal M}(\alpha_-,\alpha)
\times_{R_{\alpha}}
\overset{\circ}{\mathcal M}(\alpha,\alpha_+)$
and without corner.
The K-space $\mathcal N(\alpha_-,\alpha_+)$
may be regarded as
$\mathcal M(\alpha_-,\alpha_+) \times [0,1]$.
However the stratification of $\mathcal N(\alpha_-,\alpha_+)$
is different from that of $\mathcal M(\alpha_-,\alpha_+) \times [0,1]$.
Namely we stratify $\mathcal M(\alpha_-,\alpha_+) \times [0,1]$
as follows:
\par\noindent
Its codimension $0$ stratum is
\begin{enumerate}
\item[(0)] $\mathcal M(\alpha_-,\alpha_+) \times [0,1]$.
\end{enumerate}
Its codimension $1$ strata are
\begin{enumerate}
\item[(1-1)] $\mathcal M(\alpha_-,\alpha_+) \times \{0\}$,
\item[
(1-2)] $\mathcal M(\alpha_-,\alpha_+) \times \{1\}$,
\item[
(1-3)] $\partial\mathcal M(\alpha_-,\alpha_+) \times [0,1/2]$,
\item[
(1-4)] $\partial\mathcal M(\alpha_-,\alpha_+) \times [1/2,1]$.
\end{enumerate}
Its codimension $2$ strata are
\begin{enumerate}
\item[(2-1)] $\partial\mathcal M(\alpha_-,\alpha_+) \times \{0\}$,
\item[
(2-2)] $\partial\mathcal M(\alpha_-,\alpha_+) \times \{1/2\}$,
\item[
(2-3)] $\partial\mathcal M(\alpha_-,\alpha_+) \times \{1\}$.
\end{enumerate}
The strata (1-1), (1-2), (1-3), (1-4)
correspond to the 1st, 2nd, 3rd, 4th terms of (\ref{144444}),
respectively.
The strata (2-1), (2-2), (2-3) correspond to the 1st, 2nd, 3rd
terms of
(\ref{1445}), respectively.
\begin{figure}[h]
\centering
\includegraphics[scale=0.5]{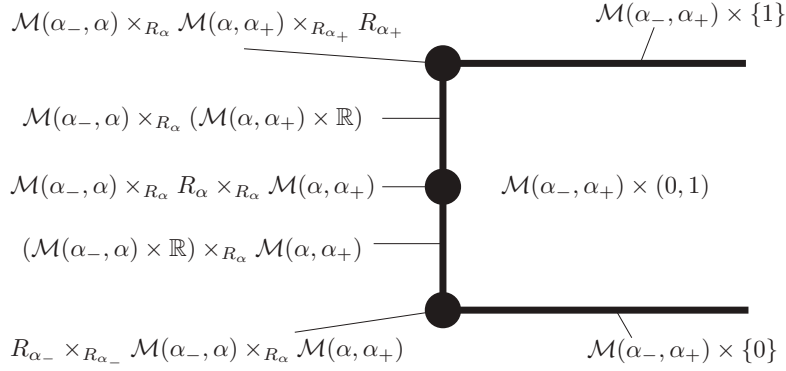}
\caption{$\mathcal N(\alpha_-,\alpha_+)$: Case 1}
\label{Figure17-4}
\end{figure}
Note that $\mathcal N(\alpha_-,\alpha_+)$ has an $\R$-action
preserving the stratification.
Namely, on the stratum where there is a $(0,1)$ factor,
we identify it with $\R$ and the $\R$ action is the
translation on it.
The other strata are in the fixed point set of this action.
\par
The quotient space
$\mathcal N(\alpha_-,\alpha_+)/\R$
is similar to $\mathcal M(\alpha_-,\alpha_+)$
but is different therefrom.
In the case of our example,
$\overset{\circ}S_0(\mathcal N(\alpha_-,\alpha_+))/\R
= \overset{\circ}S_0(\mathcal M(\alpha_-,\alpha_+))
= \overset{\circ}{\mathcal M}(\alpha_-,\alpha_+)$.
However $\overset{\circ}S_1(\mathcal N(\alpha_-,\alpha_+))/\R$
is the union of the disjoint union of
$2$ copies of $\overset{\circ}S_0(\mathcal M(\alpha_-,\alpha_+))$
and the disjoint union of two copies of
$\overset{\circ}S_1(\mathcal M(\alpha_-,\alpha_+))$.
Moreover $\overset{\circ}S_2(\mathcal N(\alpha_-,\alpha_+))/\R$
is the disjoint union of three copies of
$\overset{\circ}S_1(\mathcal M(\alpha_-,\alpha_+))$.
Then the quotient space $\mathcal N(\alpha_-,\alpha_+)/\R$
with quotient topology is non-Hausdorff.
At any rate, we do not use the quotient space in this article at all.
\par
When there exist exactly two $\alpha$ and $\alpha'$
with
$
E(\alpha_-) < E(\alpha) < E(\alpha') < E(\alpha_+),
$
the stratification of $\mathcal N(\alpha_-,\alpha_+)$
is drawn in Figure \ref{Figure17-5}.
We note that in this figure all the codimension two strata
(the vertices in Figure \ref{Figure17-5}) are contained in exactly three edges.
So its neighborhood can be identified with a corner point.
\begin{figure}[h]
\centering
\includegraphics[scale=0.55,angle=-90]{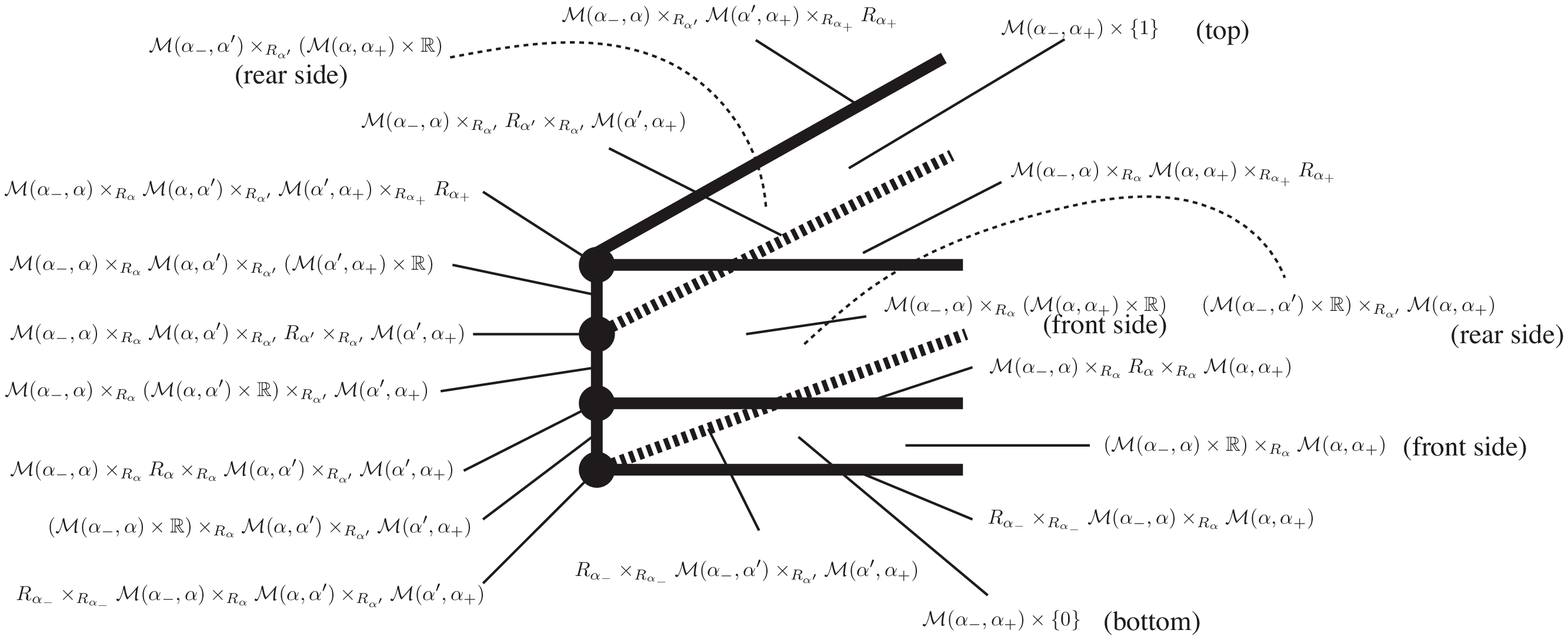}
\caption{$\mathcal N(\alpha_-,\alpha_+)$: Case 2}
\label{Figure17-5}
\end{figure}
\end{exm}
\begin{proof}[Proof of Lemma-Definition \ref{lemdefidentity}]
We first explain the way how we glue the strata and define
the topology on $\mathcal N(\alpha_-,\alpha_+)$.
(See Subsubsection \ref{productidentity} for a geometric origin of the identification
(\ref{form1452}).)
\par
Let $E$ be the energy homomorphism.
We will identity
\begin{equation}\label{form1452}
\mathcal N(\alpha_-,\alpha_+)
\cong
\mathcal M(\alpha_-,\alpha_+) \times [E(\alpha_-),E(\alpha_+)]
\end{equation}
as follows.
We identify the factor $\R$ in (\ref{1442form})
with the open interval $(E(\alpha_{k'}),E(\alpha_{k'+1}))$.
The fiber product of other factors in (\ref{1442form})
is
\begin{equation}\label{form14533}
\overset{\circ}{\mathcal M}(\alpha_-,\alpha_1)
\times_{R_{\alpha_1}}
\overset{\circ}{\mathcal M}(\alpha_1,\alpha_2)
\times_{R_{\alpha_{2}}}
\dots
\times_{R_{\alpha_{k}}}
\overset{\circ}{\mathcal M}(\alpha_{k},\alpha_{+})
\end{equation}
and is a component of $\overset{\circ}S_k({\mathcal M}(\alpha_{-},\alpha_{+}))$.
Thus (\ref{1442form}) is identified  with a subset of
$$
\overset{\circ}S_k({\mathcal M}(\alpha_{-},\alpha_{+}))
\times (E(\alpha_{k'}),E(\alpha_{k'+1}))
$$
that is a subset of (\ref{form1452}).
\par
We next consider the fiber product (\ref{1443form}).
This space is isomorphic to
(\ref{form14533}) and so is a component of
$\overset{\circ}S_k({\mathcal M}(\alpha_{-},\alpha_{+}))$.
We set its  $(E(\alpha_{-}),E(\alpha_{+}))$ factor as
$E(\alpha_{k'})$.
Thus (\ref{1443form}) is a subset of
$$
\overset{\circ}S_k({\mathcal M}(\alpha_{-},\alpha_{+}))
\times \{E(\alpha_{k'})\}
$$
that is a subset of (\ref{form1452}).
\par
We have thus specified the way how we embed all the strata of
$\mathcal N(\alpha_-,\alpha_+)$ into
(\ref{form1452}).
It is easy to see that they are disjoint from one another and their
union gives (\ref{form1452}).
Thus we identify
$\mathcal N(\alpha_-,\alpha_+)$ with
(\ref{form1452}) and
define the topology on $\mathcal N(\alpha_-,\alpha_+)$
using this identification.
\par
As we observed in Example \ref{exam1443},
the space  $\mathcal N(\alpha_-,\alpha_+)$  is
homeomorphic to
${\mathcal M}(\alpha_{-},\alpha_{+}) \times [0,1]$
but its corner structure stratification is
different from the direct product.
We next define a Kuranishi structure with corners
on $\mathcal N(\alpha_-,\alpha_+)$
compatible with the stratification
by (\ref{1442form}) and (\ref{1443form}).
We will define a Kuranishi neighborhood of each point of
$\mathcal N(\alpha_-,\alpha_+)$ below.
\par
Firstly let $\hat p$ be a point of (\ref{1442form}).
We write $\hat p = (p,s_0)$ where
$p$ is an element of (\ref{form14533}) and
$s_0 \in (E(\alpha_{k'}),E(\alpha_{k'+1}))$.
Let ${\mathcal U}_p$ be a Kuranishi neighborhood of $p$.
For simplicity of the discussion, we assume that $U_p$ consists
of a single orbifold chart $V_p/\Gamma_p$.
(In the general case we can perform the construction
below for each orbifold chart.)
Then $V_p$ is an open set of $\overline V_p \times [0,1)^k$.
We put ${\mathcal U}_p \times (s_0-\epsilon,s_0+\epsilon)$
as our Kuranishi neighborhood of $\hat p$.
In fact, the parametrization map
$s_{\hat p}^{-1}(0) \to \mathcal N(\alpha_-,\alpha_+)$
is defined as follows.
Let $(q,s) \in s_{\hat p}^{-1}(0)$.
Then $q$ parametrizes certain point of
${\mathcal M}(\alpha_{-},\alpha_{+})$.
Hence $(q,s)$ defines a point in
(\ref{form1452}).
Therefore we obtain $\psi_{\hat p}(q,s)
\in \mathcal N(\alpha_-,\alpha_+)$
by our identification.
\par
Next we consider the case $\hat p$ is in (\ref{1443form}).
First we suppose $k' = k$.
Then (\ref{1443form}) coincides with (\ref{1449form}) that is
identified with a subset of
$$
\overset{\circ}S_k(\alpha_-,\alpha_+) \times \{E(\alpha_+)\}
\subset
\overset{\circ}S_k(\alpha_-,\alpha_+) \times
\partial[E(\alpha_-),E(\alpha_+)].
$$
We take $\mathcal U_p \times (E(\alpha_+)-\epsilon,E(\alpha_+)]$
as our Kuranishi neighborhood $\mathcal U_{\hat p}$.
Here $\mathcal U_p$ is a Kuranishi neighborhood of $p$ in
${\mathcal M}(\alpha_{-},\alpha_{+})$.
The definition of the parametrization map is similar to the first case.
The case $k' = 0$ is similar to the above case.
\par
The case $k' \ne k$ and $k'\ne 0$ is the most involved, which we discuss now.
We will construct a Kuranishi neighborhood of $\hat p = (p,E(\alpha_{k'}))$.
Let $\mathcal U_p$ be the Kuranishi neighborhood of $p$ in
(\ref{form14533}).
Put $p = (p_0,p_1,\dots,p_k)$ according to this fiber product and
take a Kuranishi neighborhood $\mathcal U_{p_i}$
of $p_i$ in $\overset{\circ}{\mathcal M}(\alpha_{i},\alpha_{i+1})$.
Then the orbifold chart of $\mathcal U_p$
is given by
$$
U_{p_0} \times_{R_{\alpha_1}} \times \dots \times_{R_{\alpha_k}} U_{p_k}
\times [0,\epsilon)^k.
$$
We write its element as $(x_0,\dots,x_k;(t_1,\dots,t_k))$.
We note that the triple $(x_{k'-1},x_{k'},t_{k'})$
parametrizes a Kuranishi neighborhood of $(p_{k'-1},p_{k'})$
in $\mathcal M(\alpha_{k'-1},\alpha_{k'+1})$.
\par
For $s' \in (E(\alpha_{k'})-\epsilon,E(\alpha_{k'})+\epsilon)$
we write $(t,s) = (t_{k'},s'-E(\alpha_{k'}))$ and
change variables from $(t,s)$ to $(\rho,\sigma)  \in [0,1)^2$
by
$$
t = \rho\sigma, \qquad s = \rho - \sigma.
$$
In other words,
\begin{equation}\label{form1454}
\rho = \frac{s + \sqrt{s^2+4t}}{2},
\qquad
\sigma = \frac{-s + \sqrt{s^2+4t}}{2}.
\end{equation}
We use $x_0,\dots,x_k,t_1,\dots,t_{k'-1},t_{k'+1},\dots,t_k$ and $\rho,\sigma$
as the coordinates of
\begin{equation}\label{form1454plus}
U_{p_0} \times_{R_{\alpha_1}} \times \dots \times_{R_{\alpha_k}} U_{p_k}
\times [0,\epsilon)^k \times  (E(\alpha_{k'})-\epsilon,E(\alpha_{k'})+\epsilon).
\end{equation}
This gives a smooth structure on (\ref{form1454plus}).
In fact, the point $\hat p$ lies in the
codimension $k+1$ corner in {\it this smooth structure}
which is different from the {\it direct product smooth structure} on
(\ref{form1454plus}).
(We note that $\hat p$ lies in the codimension $k$ corner
with respect to the direct product smooth structure.)
\par
To complete the definition of the Kuranishi structure
of $\overset{\circ}{\mathcal N}(\alpha_-,\alpha_+)
$
it suffices to show that the coordinate change is admissible
in this last case.
(Admissibility of the coordinate change is trivial in other cases.)
We will check it below.
\par
We denote by $y$ the totality of the coordinates $x_0,\dots,x_k,t_1,\dots,t_{k'-1},t_{k'+1},\dots,t_k$.
Then $y,t,s$ or $y,\rho,\sigma$ are the coordinates of our
Kuranishi neighborhood.
Let $y',t',s'$ or $y',\rho',\sigma'$ be the other coordinates.
Then the coordinate change among them is given by
$$
y' = \varphi(y,t), \quad t' = \psi(y,t), \quad s' = s,
$$
where $\varphi$ and $\psi$ are smooth and satisfy
$$
\left\Vert \frac{\partial y'}{\partial t}\right\Vert_{C^k}
\le C_k e^{-c_k/t},
\quad
\Vert t' - t \Vert_{C^k}
\le C_k e^{-c_k/t}
$$
for some $c_k >0, C_k >0$.
(See Definition \ref{defn297}, Lemma \ref{newlem319} (2)
and (\ref{eq:CCexpS}).)
We note that $y'$ and $t'$ are independent of $s$.
Then using (\ref{form1454}) it is easy to check
$$
\left\Vert \frac{\partial y'}{\partial \rho}\right\Vert_{C^k}
\le C_k e^{-c_k/\rho},
\quad
\left\Vert \frac{\partial y'}{\partial \sigma}\right\Vert_{C^k}
\le C_k e^{-c_k/\sigma}
$$
and
$$
\Vert \rho' - \rho \Vert_{C^k}
\le C_k e^{-c_k/\rho},
\quad
\Vert \sigma' - \sigma \Vert_{C^k}
\le C_k e^{-c_k/\sigma}.
$$
This implies admissibility of the coordinate change
with respect to our new coordinates.
Thus we have constructed a Kuranishi structure on
${\mathcal N}(\alpha_-,\alpha_+)$.
\par
We define evaluation maps ${\rm ev}_{\pm}$ on
${\mathcal N}(\alpha_-,\alpha_+)$ by taking
one on the factor where $\alpha_-$ or $\alpha_+$ appears.
The periodicity isomorphism and orientation isomorphism
are induced by ones on ${\mathcal M}(\alpha_-,\alpha_+)$ in an obvious way.
\par
We note that in (\ref{1442form}) the factor
$\overset{\circ}{\mathcal M}(\alpha_{k'},\alpha_{k'+1}) \times \R$
can be identified with
$\overset{\circ}{\mathcal N}(\alpha_{k'},\alpha_{k'+1})
$
and in (\ref{1443form}) the
factor $R_{\alpha_{k'}}$ can be identified with
$\overset{\circ}{\mathcal N}(\alpha_-,\alpha_+)
$.
The fact that
${\mathcal N}(\alpha_-,\alpha_+)$
satisfies Conditions \ref {famiboudaru22} and \ref{furthercompatifiber}
is a consequence of the compatibility condition at the boundary and the corner
of ${\mathcal M}(\alpha_-,\alpha_+)$.
The proof of Lemma-Definition \ref{lemdefidentity} is now complete.
\end{proof}
\begin{rem}\label{geooriginid}
In the situation where we obtain a linear K-system
from the Floer equation of periodic Hamiltonian system,
morphisms among them are obtained by studying the two parameter
family of Hamiltonians $H : \R \times S^1 \times M \to \R$,
$H_{\tau,t}(x) = H(\tau,t,x)$.
In that case the interpolation space is the compactified moduli space
of solutions of the equation
\begin{equation}\label{1751form}
\frac{\partial u}{\partial \tau}
+
J\left(
\frac{\partial u}{\partial t} -  X_{H_{\tau,t}} (u)
\right) = 0.
\end{equation}
See \cite[Section 9]{fooospectr}.
In the case when $H_{\tau,t} = H_t$ is $\tau$ independent, the morphism we obtain
becomes the identity morphism defined above.
(See Subsubsection \ref{Intidentity}.)
\end{rem}
The next proposition shows
the identity morphism $\frak{ID}$ is a homotopy unit with respect to the
composition of morphisms.

\begin{prop}\label{properidentity}
For any morphism $\frak N$, the compositions
$\frak N \circ \frak{ID}$ and $\frak{ID} \circ \frak N$
are both homotopic to $\frak N$.
\end{prop}

\begin{proof}
We will prove $\frak N \circ \frak{ID} \sim \frak N$.
The proof of $\frak{ID} \circ \frak N \sim \frak N$ is similar.
Let $\mathcal N(\alpha_-,\alpha_+)$ be an interpolation space
of $\frak N$.
By the definition of the composition $\frak N \circ \frak{ID}$,
its interpolation space is decomposed into the following
two types of fiber products:
\begin{enumerate}
\item
\begin{equation}\label{1442formrev}
\aligned
&\overset{\circ}{\mathcal M^1}(\alpha_-,\alpha_1)
\times_{R_{\alpha_1}}
\overset{\circ}{\mathcal M^1}(\alpha_1,\alpha_2)
\times_{R_{\alpha_{2}}}
\dots
\times_{R_{\alpha_{k'_1-1}}}
\overset{\circ}{\mathcal M^1}(\alpha_{k'_1-1},\alpha_{k'_1})  \\
&\times_{R_{\alpha_{k'_1}}}
(\overset{\circ}{\mathcal M^1}(\alpha_{k'_1},\alpha_{k'_1+1}) \times \R)
\\
&\times_{R_{\alpha_{k'_1+1}}}
\overset{\circ}{\mathcal M^1}(\alpha_{k'_1+1},\alpha_{k'_1+2})
\times_{R_{\alpha_{k'_1+2}}}
\dots
\times_{R_{\alpha_{k_1-1}}}
\overset{\circ}{\mathcal M^1}(\alpha_{k_1-1},\alpha_{k_1})\\
&\times_{R_{\alpha_{k_1}}}
\overset{\circ}{\mathcal N}(\alpha_{k_1},\alpha'_{1})\\
&
\times_{R_{\alpha'_1}}
\overset{\circ}{\mathcal M^2}(\alpha'_1,\alpha'_2)
\times_{R_{\alpha'_{2}}}
\dots
\times_{R_{\alpha'_{k_2}}}
\overset{\circ}{\mathcal M^2}(\alpha'_{k_2},\alpha'_{+}).
\endaligned
\end{equation}
\item
\begin{equation}\label{1443formrev00}
\aligned
&\overset{\circ}{\mathcal M^1}(\alpha_-,\alpha_1)
\times_{R_{\alpha_1}}
\overset{\circ}{\mathcal M^1}(\alpha_1,\alpha_2)
\times_{R_{\alpha_{2}}}
\dots
\times_{R_{\alpha_{k'_1-1}}}
\overset{\circ}{\mathcal M^1}(\alpha_{k'_1-1},\alpha_{k'_1})  \\
&\times_{R_{\alpha_{k'_1}}}
R_{\alpha_{k'_1}}
\\
&\times_{R_{\alpha_{k'_1}}}
\overset{\circ}{\mathcal M^1}(\alpha_{k'_1},\alpha_{k'_1+1})
\times_{R_{\alpha_{k'_1+1}}}
\dots
\times_{R_{\alpha_{k_1-1}}}
\overset{\circ}{\mathcal M^1}(\alpha_{k_1-1},\alpha_{k_1})\\
&\times_{R_{\alpha_{k_1}}}
\overset{\circ}{\mathcal N}(\alpha_{k_1},\alpha'_{1})\\
&
\times_{R_{\alpha'_1}}
\overset{\circ}{\mathcal M^2}(\alpha'_1,\alpha'_2)
\times_{R_{\alpha'_{2}}}
\dots
\times_{R_{\alpha'_{k_2}}}
\overset{\circ}{\mathcal M^2}(\alpha'_{k_2},\alpha'_{+}).
\endaligned
\end{equation}
\end{enumerate}
Note that the process of gluing, described in the proof of
Lemma-Definition \ref{lemdef1434}, is included in this description.
Namely, the fiber products (\ref{1442formrev}), (\ref{1443formrev00})
appear only once in this decomposition.
We also note that ${\mathcal N}(\alpha_{-},\alpha'_{+})$,
that is an interpolation space of  $\frak N$,
is decomposed into
\begin{equation}\label{1442formrevrev}
\aligned
&\overset{\circ}{\mathcal M^1}(\alpha_-,\alpha_1)
\times_{R_{\alpha_1}}
\dots
\times_{R_{\alpha_{k_1-1}}}
\overset{\circ}{\mathcal M^1}(\alpha_{k_1-1},\alpha_{k_1})\\
&\times_{R_{\alpha_{k_1}}}
\overset{\circ}{\mathcal N}(\alpha_{k_1},\alpha'_{1})\\
&
\times_{R_{\alpha'_1}}
\overset{\circ}{\mathcal M^2}(\alpha'_1,\alpha'_2)
\times_{R_{\alpha'_{2}}}
\dots
\times_{R_{\alpha'_{k_2}}}
\overset{\circ}{\mathcal M^2}(\alpha'_{k_2},\alpha'_{+}).
\endaligned
\end{equation}
By Definition \ref{def:homotopymorph} of homotopy,
it suffices to find a K-space whose boundary is a union of
(\ref{1442formrev}), (\ref{1443formrev00}), (\ref{1442formrevrev})
and the obvious component $\mathcal N (\alpha_-, \alpha_+;\partial [0,1])$
corresponding to
the last line in \eqref{formula1211iso}.
(We omit the last obvious component here.)
\par
The interpolation space of our homotopy between
$\frak N \circ \frak{ID}$ and $\frak N$
is stratified so that its strata are (\ref{1442formrev}), (\ref{1443formrev00}), (\ref{1442formrevrev})
and the following spaces (\ref{nkake1tasu}), (\ref{nkake1tasu2}):

\begin{equation}\label{nkake1tasu}
\aligned
&\overset{\circ}{\mathcal M^1}(\alpha_-,\alpha_1)
\times_{R_{\alpha_1}}
\dots
\times_{R_{\alpha_{k_1-1}}}
\overset{\circ}{\mathcal M^1}(\alpha_{k_1-1},\alpha_{k_1})\\
&\times_{R_{\alpha_{k_1}}}
\overset{\circ}{\mathcal N}(\alpha_{k_1},\alpha'_{1}) \times (0,1)\\
&
\times_{R_{\alpha'_1}}
\overset{\circ}{\mathcal M^2}(\alpha'_1,\alpha'_2)
\times_{R_{\alpha'_{2}}}
\dots
\times_{R_{\alpha'_{k_2}}}
\overset{\circ}{\mathcal M^2}(\alpha'_{k_2},\alpha'_{+}),
\endaligned
\end{equation}
\begin{equation}\label{nkake1tasu2}
\aligned
&\overset{\circ}{\mathcal M^1}(\alpha_-,\alpha_1)
\times_{R_{\alpha_1}}
\dots
\times_{R_{\alpha_{k_1-1}}}
\overset{\circ}{\mathcal M^1}(\alpha_{k_1-1},\alpha_{k_1})\\
&\times_{R_{\alpha_{k_1}}}
\overset{\circ}{\mathcal N}(\alpha_{k_1},\alpha'_{1}) \times \{0,1\}\\
&
\times_{R_{\alpha'_1}}
\overset{\circ}{\mathcal M^2}(\alpha'_1,\alpha'_2)
\times_{R_{\alpha'_{2}}}
\dots
\times_{R_{\alpha'_{k_2}}}
\overset{\circ}{\mathcal M^2}(\alpha'_{k_2},\alpha'_{+}).
\endaligned
\end{equation}
We claim that the union of
(\ref{1442formrev}), (\ref{1443formrev00}), (\ref{1442formrevrev}) and
(\ref{nkake1tasu}), (\ref{nkake1tasu2}) above has a Kuranishi structure with corner.
The proof follows.
\par
We first define a topology on the disjoint union
$\mathcal N(\alpha_-,\alpha'_+;[0,1])$
of the spaces (\ref{1442formrev}) - (\ref{nkake1tasu2}).
Let $c$ be the energy loss of $\frak N$.
We identify the underlying topological space of
$\mathcal N(\alpha_-,\alpha'_+;[0,1])$
with
\begin{equation}\label{formm1460}
\mathcal N(\alpha_-,\alpha'_+) \times [E(\alpha_-),E(\alpha'_+)+ c].
\end{equation}
See Subsubsection \ref{productInthomot} for the geometric origin
of this identification.
We identify
(\ref{1442formrev}) - (\ref{nkake1tasu2})
to a subset of (\ref{formm1460}) as follows.
\par\smallskip
\noindent{\bf (The case of (\ref{1442formrev}))}:
We identify the $\R$ factor with $(E(\alpha_{k'_1}),E(\alpha_{k'_1+1}))$.
The fiber product of the other factors is
\begin{equation}\label{formm146022}
\aligned
&\overset{\circ}{\mathcal M^1}(\alpha_-,\alpha_1)
\times_{R_{\alpha_1}}
\overset{\circ}{\mathcal M^1}(\alpha_1,\alpha_2)
\times_{R_{\alpha_{2}}}
\dots
\times_{R_{\alpha_{k_1-1}}}
\overset{\circ}{\mathcal M^1}(\alpha_{k'_1-1},\alpha_{k'_1})  \\
&\times_{R_{\alpha_{k_1}}}
\overset{\circ}{\mathcal N}(\alpha_{k_1},\alpha'_{1})\\
&
\times_{R_{\alpha'_1}}
\overset{\circ}{\mathcal M^2}(\alpha'_1,\alpha'_2)
\times_{R_{\alpha'_{2}}}
\dots
\times_{R_{\alpha'_{k_2}}}
\overset{\circ}{\mathcal M^2}(\alpha'_{k_2},\alpha'_{+}),
\endaligned
\end{equation}
which is a subset of $\mathcal N(\alpha_-,\alpha'_+)$.
Thus (\ref{1442formrev}) is identified with the
subset of (\ref{formm1460}).
\par
\noindent{\bf (The case of (\ref{1443formrev00}))}:
(\ref{1443formrev00}) is the same as
(\ref{formm146022}) and so is a subset of $\mathcal N(\alpha_-,\alpha'_+)$.
We take $E(\alpha'_{k_1})$ as the $[E(\alpha),E(\alpha')+ c]$ factor
and regard  (\ref{1443formrev00})
as a subset of (\ref{formm1460}).
\par
\noindent{\bf (The case of (\ref{1442formrevrev}))}:
(\ref{1442formrevrev}) is again
the same as (\ref{formm146022}) and so is a subset of $\mathcal N(\alpha_-,\alpha'_+)$. We take $E(\alpha')+ c$ as the $[E(\alpha),E(\alpha')+ c]$ factor
and regard  (\ref{1442formrevrev})
as a subset of (\ref{formm1460}).
\par
\noindent{\bf (The case of (\ref{nkake1tasu}) $\cup$ (\ref{nkake1tasu2}))}:
We identify $[0,1]$ with $[E(\alpha_{k_1}),E(\alpha')+ c]$.
The fiber product of the other factors is the same as
(\ref{formm146022}).
Thus (\ref{nkake1tasu}) is identified with a
subset of (\ref{formm1460}).
\par\smallskip
Thus we have identified the union
$\mathcal N(\alpha_-,\alpha'_+;[0,1])$ of (\ref{1442formrev}) - (\ref{nkake1tasu2})
with (\ref{formm1460}).
Using this identification we define the topology of
$\mathcal N(\alpha_-,\alpha'_+;[0,1])$.
\par
In a way similar to the proof of Lemma-Definition \ref{lemdefidentity}
we define a Kuranishi structure on $\mathcal N(\alpha_-,\alpha'_+;[0,1])$
as follows.
We first observe that codimension 1 strata of
$\mathcal N(\alpha_-,\alpha'_+;[0,1])$ are one of
the following four types of fiber products:
\begin{equation}\label{3tunouchino1}
(\overset{\circ}{\mathcal M^1}(\alpha_-,\alpha_{1}) \times \R)
\\
\times_{R_{\alpha_{1}}}
\overset{\circ}{\mathcal N}(\alpha_{1},\alpha'_{+}),
\end{equation}
\begin{equation}\label{3tunouchino15}
\overset{\circ}{\mathcal M^1}(\alpha_-,\alpha_{1})
\times_{R_{\alpha_{1}}}
(\overset{\circ}{\mathcal N}(\alpha_{1},\alpha'_{+}) \times (0,1)),
\end{equation}
\begin{equation}\label{3tunouchino2}
\overset{\circ}{\mathcal N}(\alpha_{-},\alpha'_{+}) \times \{0\},
\end{equation}
\begin{equation}\label{3tunouchino3}
\overset{\circ}{\mathcal N}(\alpha_{-},\alpha'_{+}) \times \{1\}.
\end{equation}
The $\R$ factor in (\ref{3tunouchino1}) is identified with
$(E(\alpha_-),E(\alpha_1))$.
The $(0,1)$ factor in (\ref{3tunouchino15})
is identified with $(E(\alpha_1),E(\alpha'_+)+c)$.
They are glued where these factors are identified with
$E(\alpha_1)$.
The union is identified with
\begin{equation}\label{3tunouchino4}
(\overset{\circ}{\mathcal M^1}(\alpha_-,\alpha_{1})
\times_{R_{\alpha_{1}}}
\overset{\circ}{\mathcal N}(\alpha_{1},\alpha'_{+})) \times (E(\alpha_-),E(\alpha'_+)+c)).
\end{equation}
However we regard
the subset
\begin{equation}\label{3tunouchino42}
(\overset{\circ}{\mathcal M^1}(\alpha_-,\alpha_{1})
\times_{R_{\alpha_{1}}}
\overset{\circ}{\mathcal N}(\alpha_{1},\alpha'_{+})) \times \{E(\alpha_1)\}
\end{equation}
of (\ref{3tunouchino4}) not as a codimension 2 stratum, that is, a corner.
In other words, we bend the space (\ref{3tunouchino4})
at (\ref{3tunouchino42}).
This bending is performed in the same way as the argument
of case $k' \ne 0,k$ in the proof of
Lemma-Definition \ref{lemdefidentity}.
\par
The point $\{0\}$ in (\ref{3tunouchino2}) is identified with
$E(\alpha_-)$.
Therefore the closures of (\ref{3tunouchino2}) and of
(\ref{3tunouchino4})  intersect at the
codimension 2 stratum
$$
(\overset{\circ}{\mathcal M^1}(\alpha_-,\alpha_{1})
\times_{R_{\alpha_{1}}}
\overset{\circ}{\mathcal N}(\alpha_{1},\alpha'_{+})) \times \{E(\alpha_-)\}.
$$
We smooth this corner in the way we explain Subsection \ref{subsec:complinkurasmcorner}.
In fact, this smoothing occurs during the definition of the
composition $\frak N \circ \frak{ID}$.
(See Definition \ref{defn1635}.)
\par
The point $\{1\}$ in (\ref{3tunouchino3}) is identified with $E(\alpha'_+)$.
Therefore, the closures of (\ref{3tunouchino3}) and of
(\ref{3tunouchino4})  intersect at the
codimension 2 stratum
$$
(\overset{\circ}{\mathcal M^1}(\alpha_-,\alpha_{1})
\times_{R_{\alpha_{1}}}
\overset{\circ}{\mathcal N}(\alpha_{1},\alpha'_{+})) \times \{E(\alpha'_+)\}.
$$
We do not smooth this corner.
\par
We can perform an appropriate bending or smoothing and
check the consistency at the corner of higher codimension
by using the compatibility conditions of
$\mathcal M^i(\alpha_1,\alpha_2)$, $(i=1,2)$ and of
${\mathcal N}(\alpha,\alpha')$.
\par
We define a map $\mathcal N(\alpha_-,\alpha'_+;[0,1])
\to [0,1]$ appearing in Condition \ref{situparaPmorph} (IV),
so that the inverse image of $\{0\}$ is the closure of the union of
the strata (\ref{3tunouchino1}) and (\ref{3tunouchino2}).
The inverse image of $\{1\}$ is the closure of the union of
the strata (\ref{3tunouchino3}).
We note that this map is {\it not} diffeomorphic to the
projection
$$
\mathcal N(\alpha_-,\alpha'_+;[0,1])
\cong \mathcal N(\alpha_-,\alpha'_+) \times [E(\alpha_-),E(\alpha'_+)+ c]
\to [E(\alpha_-),E(\alpha'_+)+ c]
\to [0,1],
$$
where the first identification is one we used to define the
topology of $\mathcal N(\alpha_-,\alpha'_+;[0,1])$ by identifying it with
(\ref{formm1460}).
\par
We are now ready to wrap up the proof of
Proposition \ref{properidentity}.
So far we have defined a $K$-space
$\mathcal N(\alpha_-,\alpha'_+;[0,1])$.
It has most of the properties of the
interpolation space of the homotopy.
In other words, its boundary has mostly the required properties,
except the boundary is related to the composition
$\frak N \circ \frak{ID}$ {\it before} trivialization of the
boundary and smoothing the corner took place.
Therefore to obtain the required interpolation space of
homotopy we perform partial trivialization of the corner
and smoothing corner at those trivialized corners.
(See Subsection \ref{subsec:complinkurasmcorner}.)
Namely the required interpolation space of the homotopy
is
$$
\mathcal N(\alpha_-,\alpha'_+;[0,1])^{\frak C\boxplus\tau}
$$
where $\frak C$ is the part of the boundary corresponding to
$\{1\} \in [0,1]$ in the factor $[0,1]$.
Thus we obtain an interpolation space of the required $[0,1]$ parametrized morphism
and complete the proof of Proposition \ref{properidentity}.
\end{proof}
\begin{exm}
We consider the case when
$\partial \mathcal N(\alpha_-,\alpha_+)$
has only one component $\mathcal M(\alpha_-,\alpha) \times_{R_{\alpha}}
\mathcal N(\alpha,\alpha_+)$.
Then the top stratum of our homotopy is
$\mathcal N(\alpha_-,\alpha_+) \times [0,1]$.
There are 4 codimension 1 strata that are
$\mathcal N(\alpha_-,\alpha_+)$,
$R_{\alpha_-} \times_{R_{\alpha_-}} \mathcal N(\alpha_-,\alpha_+) $,
$(\mathcal M(\alpha_-,\alpha) \times\R) \times_{R_{\alpha}}\mathcal N(\alpha,\alpha_+)$,
and
$\mathcal M(\alpha_-,\alpha) \times_{R_{\alpha}} (\mathcal N(\alpha,\alpha_+)
\times [0,1])$.
The first one corresponds to the morphism $\frak N$.
The union of the second and the third corresponds to the morphism $\frak N \circ \frak{ID}$.
(See Figure \ref{Figure17-6}.)
\begin{figure}[h]
\centering
\includegraphics[scale=0.55]{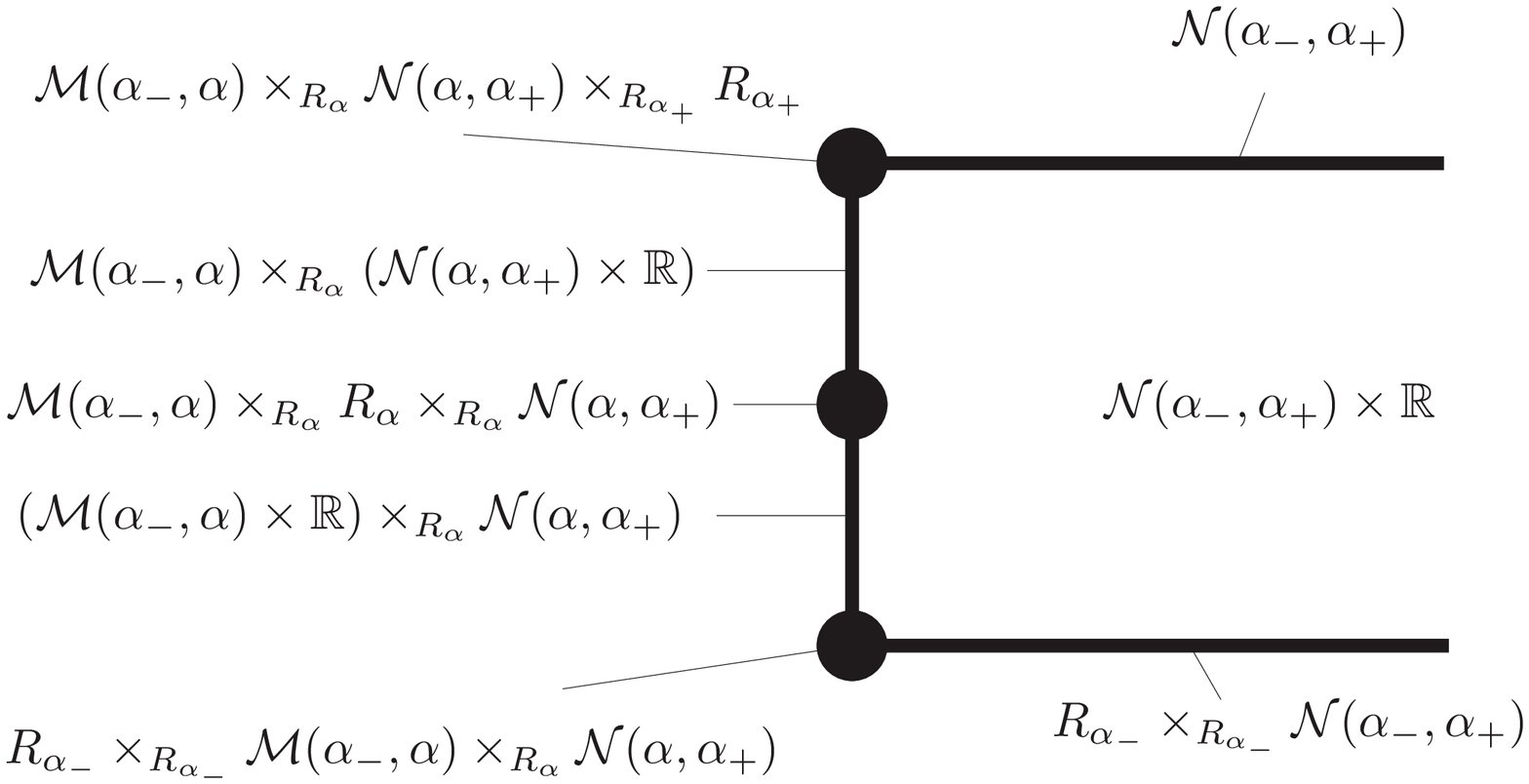}
\caption{Homotopy between $\frak N \circ \frak{ID}$ and $\frak N$: Case 1}
\label{Figure17-6}
\end{figure}
\end{exm}
\begin{exm}
Let us consider the case when there are $\alpha_1$ and $\alpha_2$ such that
$$
\partial \mathcal N(\alpha_-,\alpha_+)
=
\mathcal M(\alpha_-,\alpha_1) \times_{R_{\alpha_1}}
\mathcal N(\alpha_1,\alpha_+)
\cup
\mathcal M(\alpha_-,\alpha_2) \times_{R_{\alpha_2}}
\mathcal N(\alpha_2,\alpha_+)
$$
and the two components in the right hand side
are glued at the corner
$
\mathcal M(\alpha_-,\alpha_1) \times_{R_{\alpha_1}}
\times \mathcal N(\alpha_1,\alpha_2) \times_{R_{\alpha_2}}
\mathcal N(\alpha_2,\alpha_+)
$.
Then our homotopy has one top stratum $\mathcal N(\alpha,\alpha_+) \times (0,1)$
and 6 codimension 1 strata. Those 6 strata are
$$
\aligned
&R_{\alpha_-} \times_{R_{\alpha_-}} \mathcal N(\alpha_-,\alpha_+), \\
&(\mathcal M(\alpha_-,\alpha_1)  \times \R) \times_{R_{\alpha_1}}
\mathcal N(\alpha_1,\alpha_+) ,\\
&(\mathcal M(\alpha_-,\alpha_2)  \times \R) \times_{R_{\alpha_2}}
\mathcal N(\alpha_2,\alpha_+), \\
&\mathcal M(\alpha_-,\alpha_1) \times_{R_{\alpha_1}}
(\mathcal N(\alpha_1,\alpha_+)   \times (0,1)), \\
&\mathcal M(\alpha_-,\alpha_2) \times_{R_{\alpha_1}}
(\mathcal N(\alpha_2,\alpha_+)   \times (0,1)), \\
&\mathcal N(\alpha_-,\alpha_+).
\endaligned
$$
The first 3 strata consist $\frak N \circ \frak{ID}$
and the last stratum is $\frak N$.
See Figure \ref{Figure17-62}.
We note that each vertex in the figure
is contained in
exactly three edges.
So this configuration is one of manifold with boundary.
(However, note that the boundaries in the
interior of 3 strata consisting $\frak N \circ \frak{ID}$
are smooth during the definition of the
composition.)
\begin{figure}[h]
\centering
\includegraphics[scale=0.55,angle=-90]{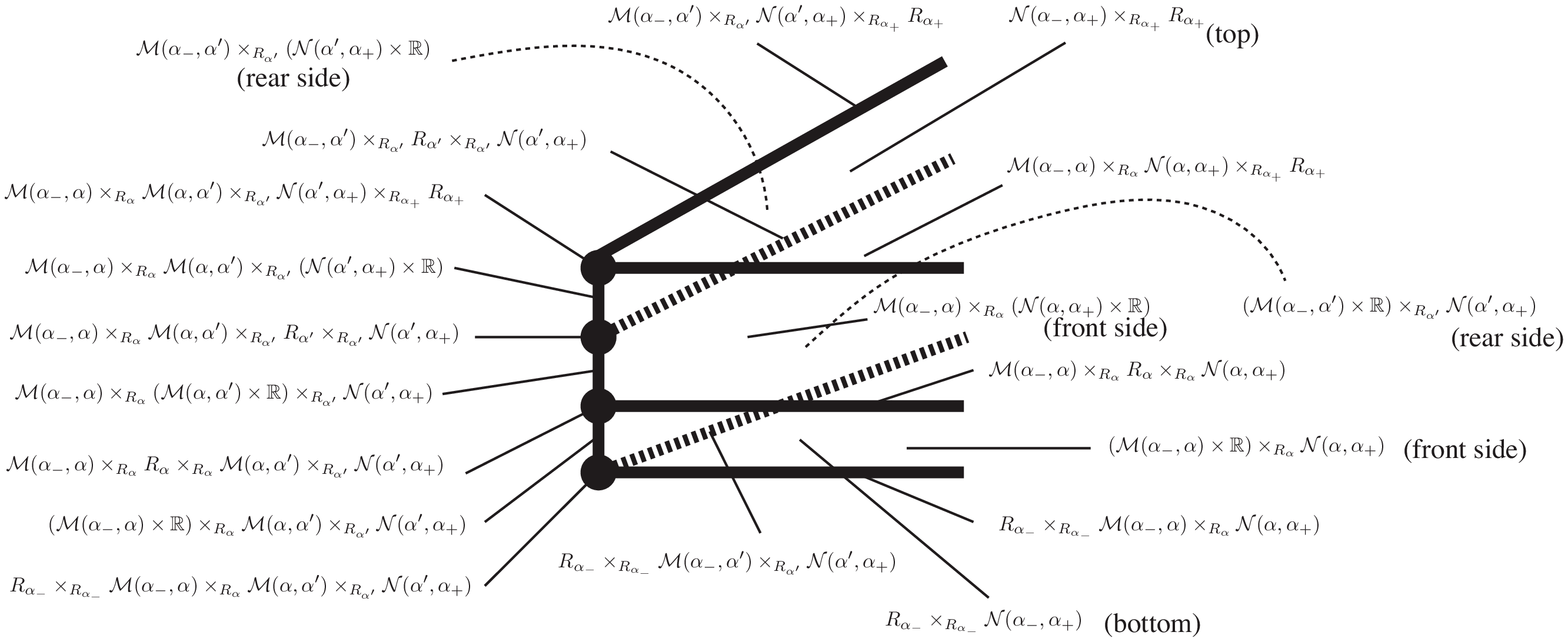}
\caption{Homotopy between $\frak N \circ \frak{ID}$ and $\frak N$: Case 2}
\label{Figure17-62}
\end{figure}
\end{exm}

\begin{rem}
The proof of Proposition \ref{properidentity} seems rather complicated.
We can prove Proposition \ref{properidentity} for the case
appearing in our geometric application, in a more intuitive way.
See Subsection \ref{subsection:identitygeo}.
\end{rem}
In the rest of this subsection, we describe an alternative
way to define the identity morphism $\frak{ID}$ and a homotopy
between $\frak N \circ \frak{ID}$ and $\frak N$.
The method we explain below is similar
to that of \cite[Subsection 4.6.1]{fooobook} in defining an
$A_{\infty}$ homomorphism,
where we used {\it time ordered fiber products}.
We put
$$
C_k(\R)
=
\{(\tau_1,\dots,\tau_k) \in (\R \cup\{\pm \infty\})^k
\mid \tau_1 \le \dots \le \tau_k\}.
$$
When $\alpha_{-} \ne \alpha_{+}$,
we consider the union of
\begin{equation}\label{topmoduli}
{\mathcal M}(\alpha_-,\alpha_1)
\times_{R_{\alpha_1}}
{\mathcal M}(\alpha_1,\alpha_2)
\times_{R_{\alpha_{2}}}
\dots
\times_{R_{\alpha_{k-1}}}
{\mathcal M}(\alpha_{k-1},\alpha_{+})
\times C_k(\R)
\end{equation}
over various $k$.
When $\alpha_- = \alpha_+ = \alpha$, we consider $R_{\alpha}$
instead of \eqref{topmoduli}.
We glue them as follows.
For the point $(\tau_1,\dots,\tau_k) \in C_k(\R)$
with $\tau_i = \tau_{i+1}$,
we put
$$
(\tau_1,\dots,\tau_i,\tau_{i+2},\dots,\tau_k)
= (\tau'_1,\dots,\tau'_{k-1}).
$$
We also consider the embedding
$$
\aligned
I_i &
:
{\mathcal M}(\alpha_-,\alpha_1)
\times_{R_{\alpha_1}}
\dots
\times_{R_{\alpha_{k-1}}}
{\mathcal M}(\alpha_{k-1},\alpha_{+})
\\
&\subset
{\mathcal M}(\alpha_-,\alpha_1)
\times_{R_{\alpha_1}}
\cdots
{\mathcal M}(\alpha_{i},\alpha_{i+2})
\dots
\times_{R_{\alpha_{k-1}}}
{\mathcal M}(\alpha_{k-1},\alpha_{+})
\endaligned
$$
which is induced by the inclusion map
$$
{\mathcal M}(\alpha_i,\alpha_{i+1})
\times_{R_{\alpha_{i+1}}}
{\mathcal M}(\alpha_{i+1},\alpha_{i+2})
\subset
\partial {\mathcal M}(\alpha_{i},\alpha_{i+2})
\subset {\mathcal M}(\alpha_{i},\alpha_{i+2}).
$$
We now identify
$$
(\frak x,(t_1,\dots,t_k))
\sim
(I_{i}(\frak x),(t'_1,\dots,t'_{k-1})).
$$
Under this identification we obtain
a K-space with corners,
which we take as an interpolation space
$\mathcal{ID}(\alpha_-,\alpha_+)$ of our
morphism $\frak{ID}$.
(We actually need to smooth the corner for this purpose.)
\par
Note that we consider the part where $\tau_i = -\infty$, $i=1,\dots,k$ to find
$$
\aligned
&\widehat S_k(\mathcal{ID}(\alpha_-,\alpha_+)) \\
&\supset
{\mathcal M}(\alpha_-,\alpha_1)
\times_{R_{\alpha_1}}
\dots
\times_{R_{\alpha_{k-1}}}
{\mathcal M}(\alpha_{k-1},\alpha_k)
\times_{R_{\alpha_{k}}}
\mathcal{ID}(\alpha_k,\alpha_+).
\endaligned
$$
We also have a similar embedding in the case when $\tau_i = +\infty$ and
consider both of them.
Then we find that $\frak{ID}$ has the required structure at the corners.
\par
To construct a homotopy from $\frak N \circ \frak{ID}$ to $\frak N$ we consider
the union of
\begin{equation}\label{topmoduli2}
{\mathcal M}(\alpha_-,\alpha_1)
\times_{R_{\alpha_1}}
{\mathcal M}(\alpha_1,\alpha_2)
\times_{R_{\alpha_{2}}}
\dots
\times_{R_{\alpha_{k-1}}}
{\mathcal N}(\alpha_{k-1},\alpha_{+})
\times C_k(\R)
\end{equation}
over various $k$ and identify them in a similar way.
Here $\tau_k \in \R \cup\{\pm \infty\} \cong [1,2]$ defines a map
to $[1,2]$.
We can prove that this space is an interpolation space of a
$[1,2]$-parametrized family of morphisms
and gives a homotopy from  $\frak N$ to $\frak N \circ \frak{ID}$.
\begin{rem}
The second method explained above seems to be shorter.
(To work out the detail of the second map, we need to glue several
K-spaces to obtain the required K-space.
We omit this detail.)
Here however we present the first one in detail since it is directly tied to the geometric situation
appearing in the Floer theory of periodic Hamiltonian system
as we mentioned in Remark \ref{geooriginid} and explained in Subsection \ref{subsection:identitygeo}.
\par
The identity morphism obtained by the second method also appears in geometric
situation in a different case. For example, when we consider two Lagrangian
submanifolds and Floer cohomology of their intersection.
Let us assume there is no pseudoholomorphic disk bounding one of
those two Lagrangian submanifolds.
Then we can define Lagrangian intersection Floer theory as
Floer did. (\cite{Flo88}. See also \cite[Chapter 2]{fooobook}.)
When we prove the independence of the choice of compatible almost complex structures,
we consider a one-parameter family of moduli spaces
${\mathcal M}(\alpha_1,\alpha_2;J_{\tau})$ where $\alpha_1, \alpha_2$ are
connected components of the intersections.
From this we can construct a similar moduli space as (\ref{topmoduli})
which provides a required cochain homotopy.
In case if the family $J_{\tau}$ is the trivial family, it boils down to the
identity morphism constructed by the second method.
\end{rem}
\begin{rem}
There is an alternative way to prove well-definedness of
Floer cohomology in the situation when $R_{\alpha}$ does {\it not}
vary. We call it the {\it bifurcation method}
in \cite[Subsection 7.2.14]{fooobook2}.
The method to use morphisms which is taken in this section
corresponds to the one which we called the {\it cobordism method} there.
There is a way to translate the bifurcation method to the cobordism
method, which is explained in \cite[Proposition 9.1, Remark 12.3]{fooo091}.
When we start from the bifurcation method and
translate it to the cobordism method, the time ordered fiber product
appears.
Then we end up with the same isomorphism as one we obtain
by the alternative proof.
\end{rem}

\subsection{Geometric origin of the definition of
the identity morphism}
\label{subsection:identitygeo}

In this subsection we explain the geometric background
of our definition of the identity morphism given in Subsection \ref{subsec:identitylinsys}.
We assume the reader is familiar with the construction of
Floer cohomology of periodic Hamiltonian system
such as the one given in \cite{Flo89I}.
Since the content of this subsection is never used in the proof of
the results of this article,
the reader can skip this subsection
if he/she prefers.
\subsubsection{Interpolation space of the identity morphism}
\label{Intidentity}
As we mentioned in Remark \ref{geooriginid},
when we study Floer cohomology of periodic Hamiltonian
system on a symplectic manifold $M$, we define a morphism between two linear K-systems
associated to $H : S^1 \times M \to \R$ and
to $H' : S^1 \times M \to \R$ by considering the
homotopy
$$
H_{\tau,t}(x) = H(\tau,t,x) : \R \times S^1 \times M \to \R,
$$
where $H(\tau,t,x) = H(t,x)$ for $\tau$ sufficiently small
and  $H(\tau,t,x) = H'(t,x)$ for $\tau$ sufficiently large,
and the interpolation space of the morphism
is the compactified moduli space
of the solutions of the equation
(\ref{1751form}).
Note that (\ref{1751form}) is {\it not} invariant under translation
of the $\tau \in \R$ direction.
We denote by
$\mathcal M(H_{\tau,t};[\gamma_-,w_-],[\gamma_-,w_+])$
the moduli space of solutions of (\ref{1751form})
with asymptotic boundary condition
$$
\lim_{\tau\to -\infty} u(\tau,t) = \gamma_-(t),
\quad
\lim_{\tau\to +\infty} u(\tau,t) = \gamma_+(t).
$$
(See Condition \ref{asboundary}.)
Moreover we assume $[w_- \# u] = [w_+]$.
Here $\gamma_{-}, \gamma_{+}$ are the periodic orbits of
the periodic Hamiltonian system associated to
$H$ and $H'$ respectively, and
$w_- : (D^2,\partial D^2) \to (M,\gamma_-)$,
$w_+ : (D^2,\partial D^2) \to (M,\gamma_+)$
are disks bounding them. We denote by
$[\gamma_{\pm},w_{\pm}]$ its homology class.
\par
We consider the Bott-Morse case. Namely the set of
$[\gamma_-,w_-]$ etc. consists of (possibly positive dimensional) smooth manifolds
$\{ R_{\alpha}\}_{\alpha}$.
We put
$$
\mathcal M(H_{\tau,t};\alpha_-,\alpha'_+)
:=
\bigcup_{[\gamma_-,w_-] \in R_{\alpha_-}, [\gamma_+,w_+] \in R_{\alpha'_+}}
\mathcal M(H_{\tau,t};[\gamma_-,w_-],[\gamma_-,w_+]),
$$
that is the interpolation space of the morphisms.
Then
the codimension $k$ corner of
$\mathcal M(H_{\tau,t};\alpha_-,\alpha'_+)$
is described by the union of
\begin{equation}\label{nkake1tasu222}
\aligned
&\mathcal M(H;\alpha_{-},\alpha_{1})
\times_{R_{\alpha_{1}}}
\dots
\times_{R_{\alpha_{k_1 -1}}}
\mathcal M(H;\alpha_{k_1-1},\alpha_{k_1})\\
&\times_{R_{\alpha_{k_1}}}
\mathcal M(H_{\tau,t};\alpha_{k_1},\alpha'_{1})\\
&
\times_{R_{\alpha'_{1}}}
\mathcal M(H';\alpha'_{1},\alpha'_{2})
\times_{R_{\alpha'_{2}}}
\dots
\times_{R_{\alpha'_{k_2}}}
\mathcal M(H';\alpha'_{k_2},\alpha'_{+})
\endaligned
\end{equation}
with $k_1 + k_2 = k$. Here
the moduli space
$\mathcal M(H;\alpha_{i},\alpha_{i+1})$ is the
set of solutions of the Floer equation
\begin{equation}\label{form1768}
\frac{\partial u}{\partial \tau}
+
J
\left(
\frac{\partial u}{\partial t} - X_{H_{t}}(u)
\right)
= 0
\end{equation}
with asymptotic boundary condition given by $R_{\alpha_{i}},R_{\alpha_{i+1}}$.
We note that we divide the set of solutions by the $\R$ action
defined by translation in the $\tau$ direction to obtain $\mathcal M(H;\alpha_i,\alpha_{i+1})$.
The space $\mathcal M(H';\alpha'_{i},\alpha'_{i+1})$ is defined in a similar
way replacing $H$ by $H'$.
\par
To define the identity morphism, we consider the case when
$H = H'$ and take the trivial homotopy.
Namely,
$$
H(\tau,t,x) \equiv H(t,x).
$$
Then the equation (\ref{1751form}) is exactly the same as (\ref{form1768}).
However, there is one important difference between them:
When we consider the interpolation space
of the morphism
$\mathcal M(H_{\tau,t};\alpha_{-},\alpha_{+})$
with $H_{\tau,t} = H_t$,
we do {\it not} divide the moduli space by $\R$ action,
even though in the case there is an $\R$ action.
In other words, we have an isomorphism:
$$
\overset{\circ}{\mathcal M}(H_{\tau,t};\alpha_{-},\alpha_{+})
= \overset{\circ}{\mathcal M}(H;\alpha_{-},\alpha_{+}) \times \R.
$$
Thus (\ref{nkake1tasu222}) becomes the closure of
(\ref{1442form}).
\par
Another point we need to
consider on $\mathcal M(H_{\tau,t};R_{\alpha_-},R_{\alpha_+})$
with $H_{\tau,t} = H_t$ is that
the case when $R = R_{\alpha_-} = R_{\alpha_+}$ can occur.
In this case, the set of solutions of
(\ref{1751form}) = (\ref{form1768})
consists of maps
each of which is constant in the $\tau \in \R$ direction and
is an element of $R$ in the $t \in S^1$ direction.
When we consider $\mathcal M(H_{\tau,t};R,R)$,
we did not regard them as the elements.
Hence this space is empty. This is because
these elements are unstable because of the $\R$ invariance.
\par
However, when we consider
$\mathcal M(H_{\tau,t};R_{\alpha_-},R_{\alpha_+})$ with
$H_{\tau,t} = H_t$ and $R = R_{\alpha_-} = R_{\alpha_+}$,
this element $u$ is included.
This is because we do not regard the $\R$ translation as
a symmetry.
In other words, we have an isomorphism
$$
\mathcal M(H_{\tau,t};R,R) = R
$$
in this case. Thus (\ref{nkake1tasu222}) becomes the closure of (\ref{1443form}).
\subsubsection{Identification of interpolation space of the identity morphism
with direct product}
\label{productidentity}
Next we explain a geometric origin of the identification
\begin{equation}\label{iso1770}
\mathcal N(\alpha_-,\alpha_+)
\cong
\mathcal M(\alpha_-,\alpha_+) \times [E(\alpha_-),E(\alpha_+)],
\end{equation}
that is (\ref{form1452}).
In our geometric situation
$$
\mathcal M(\alpha_-,\alpha_+) = \mathcal M(R_{\alpha_-},R_{\alpha_+})
$$
and $E(\alpha)$ is the value of the action functional
$$
\mathcal A_{H}([\gamma,w])
=
\int w^*\omega + \int_{S^1} H(t,\gamma(t))
$$
at $[\gamma,w] \in R_{\alpha}$.
(Note that the sign here is opposite to one \cite[(1.2)]{fooospectr}.)
Then our interpolation space $\mathcal N(\alpha_-,\alpha_+)$
is a compactification of
$\overset{\circ}{\mathcal N}(\alpha_-,\alpha_+)
=
\overset{\circ}{\mathcal M}(R_{\alpha_-},R_{\alpha_+}) \times \R$.
An element of
$\overset{\circ}{\mathcal M}(R_{\alpha_-},R_{\alpha_+})$ that is the interior
of ${\mathcal M}(R_{\alpha_-},R_{\alpha_+})$ may be regarded as a
map $u : \R\times S^1 \to M$ satisfying
(\ref{form1768}) with the asymptotic boundary condition specified
by $R_{\alpha_-}$ and $R_{\alpha_+}$.
Since we put an $\R$ factor in the definition of
$\mathcal N(\alpha_-,\alpha_+)$, its element is a map $u$ itself
but not an equivalence class by the $\R$ action
induced by the translation on the $\R$ direction.
Therefor the loop $\gamma_{\tau} : S^1 \to M$,
$\gamma_{\tau}(t) = u(\tau,t)$ is well-defined for each $\tau$.
Let $[\gamma_-,w_-]$ be an element of $R_{\alpha_-}$
such that $\gamma_-(t) = \lim_{\tau\to -\infty}u(\tau,t)$.
We consider the concatenation of
$w_-$ and the restriction of $u$ to $(-\infty,\tau] \times S^1$
and denote it by $w_{\tau}$. We put
\begin{equation}
E(u) = \mathcal A_{H}([\gamma_0,w_0]).
\end{equation}
Using the $\R$ action $(\tau_0\cdot u)(\tau,t) = u(\tau+\tau_0,t)$,
we have
$$
E(\tau_0\cdot u) = \mathcal A_{H}([\gamma_{\tau_0},w_{\tau_0}]).
$$
Therefore $E(\tau_0\cdot u)$ is an increasing function of $\tau_0 \in \R$
and
$$
\lim_{\tau_0 \to -\infty} E(\tau_0\cdot u) = E(\alpha_-),
\qquad
\lim_{\tau_0 \to +\infty} E(\tau_0\cdot u) = E(\alpha_+).
$$
Therefore using this function $E$ we obtain a homeomorphism
$$
\overset{\circ}{\mathcal N}(\alpha_-,\alpha_+)
\cong
\overset{\circ}{\mathcal M}(R_{\alpha_-},R_{\alpha_+}) \times (E(\alpha_-),
E(\alpha_+)).
$$
It is easy to see that this homeomorphism extends
to (\ref{iso1770}).
\subsubsection{Interpolation space of the homotopy}
\label{Inthomot}
We next describe a geometric origin of the
homotopy $\frak N \circ \frak{ID} \sim \frak N$
given in the proof of Proposition \ref{properidentity}.
\par
We consider the case when the interpolation space
of
$\frak N$ is  $\mathcal M(H_{\tau,t};R_{-},R_{+})$
where $H_{\tau,t}$ is a homotopy from $H$ to $H'$.
This moduli space is the set of solutions of
the equation (\ref{1751form})
with asymptotic boundary condition given by $R_-$ and $R_+$.
Note that $H_{\tau,t}$ is $\tau$ dependent.
We write this $H_{\tau,t}$ as $H^{32}_{\tau,t}$.
\par
On the other hand, the interpolation space of the morphism
$\frak{ID}$ is by definition the
set of solutions of (\ref{1751form})
with $H_{\tau,t} = H_t$, which is $\tau$ independent.
We write this $H_{\tau,t} = H_t$ as $H^{21}_{\tau,t}$.
\par
Then the composition $\frak N \circ \frak{ID}$
is obtained by using the two parameter family of
Hamiltonians concatenating $H_{\tau,t}$ and $H_t$
as in (\ref{eq1666}).
Let $H^{31,T} = H^{31,T}_{\tau,t}$ be obtained by
this concatenation.
More precisely, the interpolation space of
$\frak N \circ \frak{ID}$
appears when we take
the limit of $H^{31,T}$
as $T \to \infty$.
Note that in our case $H^{31,T}$ is
\begin{equation}\label{eq1666rev}
H^{31,T}(\tau,t,x)
=
\begin{cases}
H(t,x)
&\text{if $\tau \le -T_0-T$}\\
H^{21}(\tau+T,t,x) =H(t,x)
&\text{if $-T_0-T\le \tau \le T_0 - T$}\\
H^2(t,x) = H(t,x)
&\text{if $T_0 - T \le \tau \le T-T_0$}\\
H^{32}(\tau-T,t,x) = H_{\tau-T,t}(x)
&\text{if $T-T_0\le \tau \le T+T_0 $}\\
H^3(t,x) = H'(t,x)
&\text{if $T+T_0 \le \tau$}.
\end{cases}
\end{equation}
\begin{figure}[h]
\centering
\includegraphics[scale=0.4]{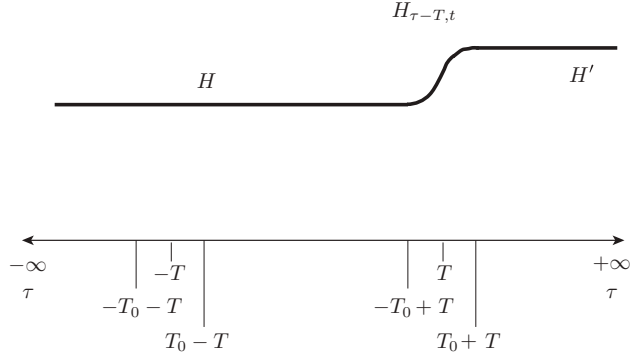}
\caption{$H^{31,T}(\tau,t,x)$}
\label{Figure17-7}
\end{figure}
See Figure \ref{Figure17-7}.
Thus actually we have
$$
H^{31,T}(\tau,t,x) = H_{\tau-T,t}(x).
$$
Therefore the set of solutions of (\ref{equa165}),
with respect to $\tau$, $t$ dependent Hamiltonian
$H^{31,T}$,
that is
\begin{equation}\label{eq1573}
\frac{\partial u}{\partial \tau}
+
J
\left(
\frac{\partial u}{\partial t} - X_{H_{\tau-T,t}}(u)
\right) = 0
\end{equation}
is indeed independent of $T$
up to the canonical isomorphism.
Moreover if $T=0$, this
equation is exactly the same as one we used
to define $\frak N$.
\par
We consider the union of the set of solutions of (\ref{eq1573})
for $T \in [0,\infty)$.
Namely
\begin{equation}\label{eq15731}
\bigcup_{T \in [0,\infty)}\overset{\circ}{\mathcal M}(H^{31,T};\alpha_-,\alpha'_+)
\times \{T\}.
\end{equation}
The interpolation space
of the homotopy is a compactification of (\ref{eq15731}).
Except the part $T = \infty$, the compactification is a
product of $[0,\infty)$ and the compactification
of the moduli space of solutions of (\ref{eq1573})
for fixed $T=0$, which is
$\mathcal N(\alpha_-,\alpha'_+) \times [0,\infty)$.
\footnote{Note $\mathcal N(\alpha_-,\alpha'_+) = {\mathcal M}(H^{31,T};\alpha_-,\alpha'_+)$.}
\par
However there is rather a delicate issue related to
the (source) $\R$ action at the limit as $T \to \infty$.
Recall that $H^{21}$ is $\tau$ independent.
As we already mentioned several times, we do not use this symmetry
to divide our moduli space.
In order to clarify this point, we take and fix a marked point
and regard our moduli space as a map from a marked cylinder.
(We take some marked point $(-T,1/2) \in (-\infty,0] \times S^1$.)
The coordinate $T \in [0,\infty)$ in turn becomes
the $[0,\infty)$ factor of $\mathcal N(\alpha_-,\alpha'_+) \times [0,\infty)$.
We note that $\mathcal N(\alpha_-,\alpha'_+)$ is the set of
solutions of
\begin{equation}\label{eq15732}
\frac{\partial u}{\partial \tau}
+
J
\left(
\frac{\partial u}{\partial t} - X_{H_{\tau,t}}(u)
\right) = 0.
\end{equation}
Then we can identify $\mathcal N(\alpha_-,\alpha'_+) \times [0,\infty)$
with the set of $(u,(-T,1/2))$ where $u$ solves the equation (\ref{eq15732}) and
$T \in [0,\infty)$.
In other words,
putting the marked point $(-T,1/2)$ corresponds to shifting
$H^{31,0}$ to $H^{31,T}$.
Then by a standard gluing analysis we can describe its limit as $T \to \infty$
by the union of
\begin{equation}\label{nkake1tasu222333}
\aligned
&\mathcal M(H;\alpha_{-},\alpha_{1})
\times_{R_{\alpha_1}}
\dots
\times_{R_{\alpha_{k'_1-2}}}
\mathcal M(H;\alpha_{k'_1-2},\alpha_{k'_1-1})\\
&\times_{R_{\alpha_{k'_1-1}}}
\mathcal M(H;\alpha_{k'_1-1},\alpha_{k'_1}) \times \R\\
&\times_{R_{\alpha_{k'_1}}}\mathcal M(H;\alpha_{k'_1},\alpha_{k'_1+1})
\times_{R_{\alpha_{k'_1+1}}}
\dots
\times_{R_{k_1-1}}
\mathcal M(H;\alpha_{k_1-1},\alpha_{k_1})\\
&\times_{R_{\alpha_{k_1}}}
\mathcal M(H_{\tau,t};\alpha_{k_1},\alpha'_{1})\\
&
\times_{R_{\alpha'_1}}
\mathcal M(H';\alpha'_{1},\alpha'_{2})
\times_{R_{\alpha'_2}}
\dots
\times_{\alpha'_{k_2}}
\mathcal M(H';\alpha'_{k_2},\alpha'_{+}).
\endaligned
\end{equation}
\par\smallskip
\noindent
See Figure \ref{Figure17-8}, that is the case
$k_1 = 1$, $k'_2 = 0$.
An object in the compactification of (\ref{eq15731}) corresponding
to the limit
as $T\to \infty$ is obtained by gluing maps
appearing in (\ref{nkake1tasu222333}).
In fact,
since an element of $\mathcal M(H^{31,T};\alpha_{-},\alpha'_{+})$
comes with a marked point, the limit
is assigned with a marked point on one of the factors.
We put the $\R$ factor to the factor on which the marked point lies.
In (\ref{nkake1tasu222333}) it lies on
the map representing an element of $\mathcal M(H;\alpha_{k'_1-1},\alpha_{k'_1})$.
\par
The space (\ref{nkake1tasu222333}) is the closure of (\ref{1442formrev}).
\begin{figure}[h]
\centering
\includegraphics[scale=0.3,angle=-90]{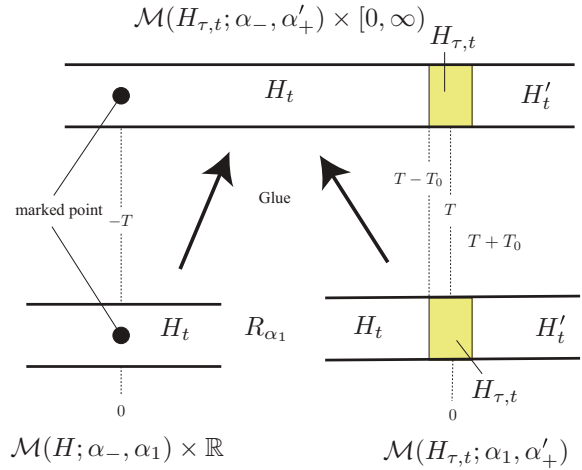}
\caption{Boundary of $\mathcal M(H_{\tau,t};\alpha_{-},\alpha'_{+}) \times [0,\infty)$}
\label{Figure17-8}
\end{figure}
By the same reason as we discussed in Subsebsection \ref{productidentity},
the limit
as $T\to \infty$ also contains a component
\begin{equation}\label{nkake1tasu22233322}
\aligned
&\mathcal M(H;\alpha_{-},\alpha_{1})
\times_{R_{\alpha_1}}
\dots
\times_{\alpha_{k_1-1}}
\mathcal M(H;\alpha_{k'_1-1},\alpha_{k'_1})\\
&\times_{R_{\alpha_{k'_1}}}
R_{\alpha_{k'_1}}\\
&\times_{R_{\alpha_{k'_1}}}\mathcal M(H;\alpha_{k'_1},\alpha_{k'_1+1})
\times_{\alpha_{k_1+1}}
\dots
\times_{R_{\alpha_{k_1-1}}}
\mathcal M(H;\alpha_{k_1-1},\alpha_{k_1})\\
&\times_{R_{\alpha_{k_1}}}
\mathcal M(H_{\tau,t};\alpha_{k_1},\alpha'_{1})\\
&
\times_{R_{\alpha'_1}}
\mathcal M(H';\alpha'_{1},\alpha'_{2})
\times_{R_{\alpha'_2}}
\dots
\times_{R_{\alpha'_{k_2}}}
\mathcal M(H';\alpha'_{k_2},\alpha'_{+}).
\endaligned
\end{equation}
Namely, this space (\ref{nkake1tasu22233322}) corresponds to the case when the marked point
lies in the component corresponding to a map $u : \R \times S^1 \to M$
which is constant in the $\R$ direction and represents an element of
$R_{\alpha_{k'_1}}$ on each $\{\tau\} \times S^1$.
The space (\ref{nkake1tasu22233322}) coincides with the closure of (\ref{1443formrev00}).
\par
Therefore the compactification of (\ref{eq15731})
coincides with the interpolation space of the homotopy constructed
in the proof of Proposition \ref{properidentity}.
\subsubsection{Identification of the interpolation space of the homotopy
with direct product}
\label{productInthomot}
Finally, we explain a geometric origin of the identification
\begin{equation}\label{form17781}
\mathcal N(\alpha_-,\alpha'_+;[0,1])
\cong \mathcal N(\alpha_-,\alpha'_+) \times [E(\alpha_-),E(\alpha'_+)+ c].
\end{equation}
The discussion is similar to that in Subsection \ref{productidentity} but is slightly
more involved.
We denote by
$\overset{\circ}{\mathcal N}(\alpha_-,\alpha'_+;[0,1])$
the set of pairs $(u,-T)$ where
$u : \R \times S^1  \to M$ is a map solving the equation (\ref{eq15732})
and satisfying the asymptotic boundary conditions as $\tau
\to \pm\infty$ given by $R_{\alpha_-}$ and $R_{\alpha'_+}$,
and $(-T,1/2)$ is the marked point in $(-\infty,0] \times S^1$.
As we discussed in Subsubsection \ref{Inthomot},
$\mathcal N(\alpha_-,\alpha'_+;[0,1])$ is a compactification
of this space.
\par
We first define a map
$$
E : \overset{\circ}{\mathcal N}(\alpha_-,\alpha'_+;[0,1]) \to \R
$$
and modify it to $E'$ later.
For $(u,T) \in \overset{\circ}{\mathcal N}(\alpha_-,\alpha'_+;[0,1])$
we put $\gamma_{\tau}(t) = u(\tau,t)$ and
$$
\gamma_- = \lim_{\tau \to -\infty} \gamma_{\tau}.
$$
The asymptotic boundary condition we assumed for
$\overset{\circ}{\mathcal N}(\alpha_-,\alpha'_+;[0,1])$
implies that there exists $w_-$ such that $[\gamma_-,w_-] \in R_{\alpha_-}$.
We denote by $w_{T}$ the concatenation of $w_-$ and the restriction of $u$
to $(-\infty,-T) \times S^1$ and define
\begin{equation}
E(u,-T)
= -\int_{D^2}  w_{T}^* \omega - \int_{t \in S^1} H_{\tau,t}(\gamma_T(t))
\in \R.
\end{equation}
We note that
\begin{equation}
\lim_{T\to \infty} E(u,-T) = \mathcal A_H([\gamma_-,w_-]) = E(\alpha_-).
\end{equation}
\par
On the other hand, contrary to the situation of
Subsection \ref{productidentity}, the map
$$
(-\infty,0] \ni -T \mapsto E(u,-T) \in \R
$$
may not be an increasing function in general, although
it is an increasing function for sufficiently large $T$.
Moreover there is no obvious relation between $E(u,0)$ and
$E(\alpha'_+)$.
Since the energy loss is $c$, the inequality
$E(\alpha_-) < E(\alpha'_+)+c$
does hold, but $E(\alpha_-) < E(\alpha'_+)$
may not hold in general.
Thus we modify $E(u,T)$ to $E'(u,T)$ with the following properties.
\begin{enumerate}
\item
$
E'(u,-T) = E(u,-T)
$ if $T$ is sufficiently large.
\item
$-T \mapsto E'(u,-T)$ is strictly increasing.
\item
$E'(u,0) = E(\alpha'_+)+c$.
\end{enumerate}
More explicitly, we can define $E'$ as follows.
We first take a sufficiently large number $T_1 > 0$ with the following properties:
\begin{enumerate}
\item[(a)]
$
E(u,-T) \le E(\alpha_+) + c
$
holds for each $u \in {\mathcal N}(\alpha_-, \alpha'_+)$ and $T\ge T_1$.
\item[(b)]
$H_{\tau,t} = H_t$ if $\tau < -T_1$.
\end{enumerate}
Then we define
\begin{equation}
E'(u,-T) =
\begin{cases}
E(u,-T)  &\text{if $T \ge T_1$} \\
\frac{T_1-T}{T_1} (E(\alpha_+) + c - E(u,-T_1)) + E(u,-T_1)
&\text{if $T \in [0,T_1].$}
\end{cases}
\end{equation}
Item (b) implies that the function $-T \mapsto E(u,-T)$ is increasing
for $T>T_1$. Item (a) implies that $-T \mapsto E'(u,-T)$ is
an increasing function for $T \in [0,T_1]$.
We have thus verified (2). (1) and (3) are obvious from definition.
We can easily extend it to ${\mathcal N}(\alpha_-,\alpha'_+;[0,1])$.
In fact, the value $E'$ is in the interval $(E(\alpha_{k'_1-1}),E(\alpha_{k'_1}))$
on (\ref{nkake1tasu222333}) and is $E(\alpha_{k'_1})$ on
(\ref{nkake1tasu22233322}).
This is a consequence of Item (1).
Then using $E'$ as the second factor, we obtain the identification (\ref{form17781}).

\section{Linear K-system: Floer cohomology II:
proof}
\label{sec:systemline3}

The purpose of this section is to prove the theorems
we claimed in Section \ref{sec:systemline1}.

\subsection{Construction of cochain complexes}
\label{subsec:constchaincpx}

Let us start with a linear K-system
or a partial linear K-system of energy cut level
$E_0$ as in Definition \ref{linearsystemdefn}.
We set $E_0 = + \infty$ for the case of linear K-system
and $E_0$ to be the energy cut level for the case of partial linear K-system.
We consider the set
\begin{equation}\label{16111}
\frak E_{\le E_0}
=
\{E(\alpha_+) - E(\alpha_-)  \mid
\mathcal M(\alpha_-,\alpha_+) \ne \emptyset, ~E(\alpha_+) - E(\alpha_-) \le E_0 \}.
\end{equation}
This is a discrete set by Condition \ref{linsysmainconds} (IX).
We put
$$
\frak E_{\le E_0} = \{E_{\frak E}^1,E_{\frak E}^2,\dots\}
$$ with
$E_{\frak E}^1 < E_{\frak E}^2 < \dots$.
We will use the results of Section \ref{sec:triboundary}
to prove the next proposition
by induction on $k$ of the energy cut level $E_{\frak E}^k$.
Put $\tau_0 =1$.
\begin{prop}\label{prop161}
For each $0<\tau < 1$,
there exist a $\tau$-collared Kuranishi structure
$\widehat{\mathcal U^+}(\alpha_-,\alpha_+)$ of
$\mathcal M(\alpha_-,\alpha_+)^{\boxplus\tau_0}$
and a CF-perturbation
$\widehat{\frak S^+}(\alpha_-,\alpha_+)$ of $\widehat{\mathcal U^+}(\alpha_-,\alpha_+)$
for every $\alpha_-$, $\alpha_+$ with $E(\alpha_+) - E(\alpha_-) \le E_{\frak E}^k$
and
they enjoy the following properties:
\begin{enumerate}
\item
Let $\widehat{\mathcal U}(\alpha_-,\alpha_+)$ be the Kuranishi structure
on $\mathcal M(\alpha_-,\alpha_+)$ given in Condition \ref{linsysmainconds}.
Then  $\widehat{\mathcal U}(\alpha_-,\alpha_+)^{\boxplus\tau_0}
< \widehat{\mathcal U^+}(\alpha_-,\alpha_+)$
as collared Kuranishi structures.\footnote{Let
$\widehat{\mathcal U_1}^{\boxplus\tau_1}, \widehat{\mathcal U_2}^{\boxplus\tau_2}$ be
two collared Kuranishi structures. We define
{\it $\widehat{\mathcal U}_1^{\boxplus\tau_1} < \widehat{\mathcal U}_2^{\boxplus\tau_2}$
as collared Kuranishi structures}
if there exist Kuranishi structures $\widehat{\mathcal U_1'}, \widehat{\mathcal U_2'}$ such that
$\widehat{\mathcal U_1'} < \widehat{\mathcal U_2'}$ and
$\widehat{\mathcal U_1'}^{\boxplus\tau_1'}=\widehat{\mathcal U_1}^{\boxplus\tau_1}$,
$\widehat{\mathcal U_2'}^{\boxplus\tau_2'}=\widehat{\mathcal U_2}^{\boxplus\tau_2}$
for some $\tau_1', \tau_2'$.}
\item
The CF-perturbation $\widehat{\frak S^+}(\alpha_-,\alpha_+)$ is transversal to $0$.
Moreover\footnote{According to Lemma \ref{lemma1523} (3), we should write
${\rm ev}^{\boxplus\tau_0}_{+}$ but not ${\rm ev}_{+}$. However,
to simplify the notation we drop ${}^{\boxplus\tau_0}$ from the notation of evaluation maps if no confusion can occur.}
\index{${\rm ev}^{\boxplus\tau_0}$}
${\rm ev}_+ : \mathcal M(\alpha_-,\alpha_+)^{\boxplus\tau_0}
\to M$ is strongly submersive with respect to $\widehat{\frak S^+}(\alpha_-,\alpha_+)$.
\item We have the following isomorphism, called the periodicity isomorphism
$$
\widehat{\mathcal U^+}(\alpha_-,\alpha_+)
\longrightarrow \widehat{\mathcal U^+}(\beta\alpha_-,\beta\alpha_+)
$$
for any $\beta \in \frak G$.
It is compatible with the periodicity isomorphism
in Condition \ref{linsysmainconds} (VIII) via the embedding
in (1).
The pull-back of $\widehat{\frak S^+}(\beta\alpha_-,\beta\alpha_+)$
by this isomorphism is equivalent to $\widehat{\frak S^+}(\alpha_-,\alpha_+)$.
\item
There exists an isomorphism of $\tau$-collared
K-spaces.\footnote{See Remark \ref{rem:FiberProdOrd} for the sign and the order of the fiber products.}
\begin{equation}\label{formula162}
\aligned
&\partial(\mathcal M(\alpha_-,\alpha_+)^{\boxplus\tau_0},\widehat{\mathcal U^+}(\alpha_-,\alpha_+)) \\
&=
\coprod_{\alpha}
(-1)^{\dim \mathcal M(\alpha,\alpha_+)}
(\mathcal M(\alpha,\alpha_+)^{\boxplus\tau_0},\widehat{\mathcal U^+}(\alpha,\alpha_+))
{}_{{\rm ev}_-}\times_{{\rm ev}_+}
(\mathcal M(\alpha_-,\alpha)^{\boxplus\tau_0},\widehat{\mathcal U^+}(\alpha_-,\alpha))
\endaligned
\end{equation}
The isomorphism (\ref{formula162}) is compatible with the isomorphism
$$
\partial\widehat{\mathcal U}(\alpha_-,\alpha_+)^{\boxplus\tau_0}
=
\coprod_{\alpha}
(-1)^{\dim \widehat{\mathcal U}(\alpha,\alpha_+)}
(\widehat{\mathcal U}(\alpha,\alpha_+)^{\boxplus\tau_0})
{}_{{\rm ev}_-}\times_{{\rm ev}_+}
(\widehat{\mathcal U}(\alpha_-,\alpha)^{\boxplus\tau_0})
$$
which is induced by Condition \ref{linsysmainconds} (X)
via the embedding in (1).
\item
The pull-back of
$\widehat{\frak S^+}(\alpha_-,\alpha_+)$
by the isomorphism
(\ref{formula162}) is equivalent to
the fiber product
$$
\widehat{\frak S^+}(\alpha,\alpha_+)
{}_{{\rm ev}_-}\times_{{\rm ev}_+}
\widehat{\frak S^+}(\alpha_-,\alpha).
$$
This fiber product is well-defined by (2).
\item
There exists an isomorphism
of $\tau$-collared K-spaces
\begin{equation}\label{cornecominduction}
\aligned
&\widehat S_k(\mathcal M(\alpha_-,\alpha_+)^{\boxplus\tau_0},\widehat{\mathcal U^+}(\alpha_-,\alpha_+))\\
&\cong
\coprod_{\alpha_1,\dots,\alpha_k
\in \frak A}
\left(
(\mathcal M(\alpha_-,\alpha_1)^{\boxplus\tau_0},\widehat{\mathcal U^+}(\alpha_-,\alpha_1))
\,\,{}_{{\rm ev}_{+}}\times_{R_{\alpha_1}}\,\cdots\right.\\
&\qquad \qquad\qquad \,\,\cdots \left.{}_{R_{\alpha_k}}\times_{{\rm ev}_{-}}
(\mathcal M(\alpha_k,\alpha_+)^{\boxplus\tau_0},\widehat{\mathcal U^+}(\alpha_k,\alpha_+))
\right).
\endaligned
\end{equation}
The isomorphism (\ref{cornecominduction}) is
compatible with the isomorphism
$$
\aligned
&\widehat S_k(\widehat{\mathcal U}(\alpha_-,\alpha_+)^{\boxplus\tau_0})\\
&\cong
\coprod_{\alpha_1,\dots,\alpha_k
\in \frak A}
\left(
\widehat{\mathcal U}(\alpha_-,\alpha_1)^{\boxplus\tau_0}
{}_{{\rm ev}_{+}}\times_{R_{\alpha_1}}
\dots {}_{R_{\alpha_k}}\times_{{\rm ev}_{-}}
\widehat{\mathcal U}(\alpha_k,\alpha_+)^{\boxplus\tau_0}
\right)
\endaligned
$$
which is induced by Condition \ref{linsysmainconds} (XI)
via the embedding in (1).
\item
The isomorphism (\ref{cornecominduction}) implies that
$\widehat S_{\ell}(\widehat S_k(\mathcal M(\alpha_-,\alpha_+)^{\boxplus\tau_0},\widehat{\mathcal U^+}(\alpha_-,\alpha_+)))$
is decomposed into similar fiber products as (\ref{cornecominduction})
where $k$ is replaced by $\ell+k$.
On the other hand, $\widehat S_{\ell+k}(\mathcal M(\alpha_-,\alpha_+)^{\boxplus\tau_0},\widehat{\mathcal U^+}(\alpha_-,\alpha_+))$
is decomposed into similar fiber products as (\ref{cornecominduction})
where $k$ is replaced by $\ell+k$.
The map
$$
\pi_{\ell,k} : \widehat S_{\ell}(\widehat S_k(\mathcal M(\alpha_-,\alpha_+)^{\boxplus\tau_0},\widehat{\mathcal U^+}(\alpha_-,\alpha_+)))
\to \widehat{S}_{\ell+k}(\mathcal M(\alpha_-,\alpha_+)^{\boxplus\tau_0},\widehat{\mathcal U^+}(\alpha_-,\alpha_+))
$$
in Proposition \ref{prop2813} becomes an identity map on each component under this identification.
\item
The pull-back of $\widehat{\frak S^+}(\alpha_-,\alpha_+)$
by the isomorphism (\ref{cornecominduction})
is the fiber product:
\begin{equation}\label{fiberproductSSinduction}
\widehat{\frak S^+}(\alpha_-,\alpha_1)\,\,
{}_{{\rm ev}_{+}}\times_{R_{\alpha_1}}
\dots {}_{R_{\alpha_k}}\times_{{\rm ev}_{-}}
\widehat{\frak S^+}(\alpha_k,\alpha_+).
\end{equation}
This fiber product is well-defined by (2).
This isomorphism is compatible with the
covering map $\pi_{\ell,k}$.
\end{enumerate}
\end{prop}
\begin{proof}
Using the results of Subsection \ref{subsection:concltrisection},
we can prove the proposition in a
straightforward way as follows.
We take $\tau^+$ with
$$
\tau < \tau^+ < 1=\tau_0.
$$
Suppose we constructed $\widehat{\mathcal U^+}(\alpha_-,\alpha_+)$
and $\widehat{\frak S^+}(\alpha_-,\alpha_+)$
satisfying (1)-(8) above for any $\alpha_-$, $\alpha_+$
with $E(\alpha_+) - E(\alpha_-) < E_{\frak E}^k$
and for $\tau$ replaced by $\tau^+$.
We consider $\alpha_-$, $\alpha_+$
with $E(\alpha_+) - E(\alpha_-) = E_{\frak E}^k$ and
$\mathcal M(\alpha_-,\alpha_+)^{\boxplus\tau_0}$.
\par
To apply Propositions \ref{prop528} and \ref{prop1562},
we check that we are in Situation \ref{sit1526}.
The space $X$ in Situation \ref{sit1526}
is $\mathcal M(\alpha_-,\alpha_+)^{\boxplus\tau_0}$ and
$\tau$ in Situation \ref{sit1526}  is $\tau^+$ here.
The $\tau$-collared Kuranishi structure $\widehat{\mathcal U}$
in  Situation \ref{sit1526}  is the Kuranishi structure $\widehat{\mathcal U}(\alpha_-,\alpha_+)^{\boxplus\tau_0}$ here.
The Kuranishi structure
$\widehat{\mathcal U^+_{S_k}}$
in Situation \ref{sit1526}
is the Kuranishi structure in the right hand side of
(\ref{cornecominduction}), which is given by the induction hypothesis.
We write it as $\widehat{\mathcal U^+_{S_k}}$.
Then it is easy to see from definition that
$\widehat S_{k}(\widehat{\mathcal U^+_{S_{\ell}}})$
is $(k+\ell)!/k!\ell!$ disjoint union of
$\widehat{\mathcal U^+_{S_{k+\ell}}}$.
Therefore  Situation \ref{sit1526}  (2) holds.
The commutativity of Diagram (\ref{diagin26277XX})
holds because all the maps in (\ref{diagin26277XX}) is the identity
map via the isomorphism (\ref{cornecominduction}).
We next construct the embedding
$\widehat S_k(\mathcal M(\alpha_-,\alpha_+)^{\boxplus\tau_0},
\widehat{\mathcal U}(\alpha_-,\alpha_+)^{\boxplus\tau_0})
\to \widehat S_{k}(\widehat{\mathcal U^+_{S_{\ell}}})$.
By induction hypothesis the embedding
$\widehat{\mathcal U}(\alpha,\alpha')^{\boxplus\tau_0}
\to \widehat{\mathcal U^+}(\alpha,\alpha')$
for $E(\alpha') - E(\alpha) < E_{\frak E}^k$ is given.
By using Condition \ref{linsysmainconds} (XI)
we have
$$
\aligned
&\widehat S_k(\widehat{\mathcal U}(\alpha_-,\alpha_1)^{\boxplus\tau_0})\\
&\cong
\coprod_{\alpha_1,\dots,\alpha_k
\in \frak A}
\left(
\widehat{\mathcal U}(\alpha_-,\alpha_1)^{\boxplus\tau_0}
{}_{{\rm ev}_{+}}\times_{R_{\alpha_1}}
\cdots {}_{R_{\alpha_k}}\times_{{\rm ev}_{-}}
\widehat{\mathcal U}(\alpha_k,\alpha_+)^{\boxplus\tau_0}
\right).
\endaligned
$$
Therefore, the right hand side is embedded in
$$
\aligned
\widehat{\mathcal U^+_{S_k}} &\cong
\coprod_{\alpha_1,\dots,\alpha_k
\in \frak A}
\left(
(\mathcal M(\alpha_-,\alpha_1)^{\boxplus\tau_0},\widehat{\mathcal U^+}(\alpha_-,\alpha_1))
\,\,{}_{{\rm ev}_{+}}\times_{R_{\alpha_1}}\,\cdots\right.\\
&\qquad \qquad\qquad \,\,\cdots \left.{}_{R_{\alpha_k}}\times_{{\rm ev}_{-}}
(\mathcal M(\alpha_k,\alpha_+)^{\boxplus\tau_0},\widehat{\mathcal U^+}(\alpha_k,\alpha_+))
\right).
\endaligned
$$
So Situation \ref{sit1526}  (4) is satisfied.
The commutativity of Diagram (\ref{diagin26277XX}) and Diagram (\ref{diag15main}) follows from
the fact that all the maps in of Diagrams  (\ref{diagin26277XX}), (\ref{diag15main}) become
the identity maps via the isomorphism (\ref{cornecominduction}),
which is a part of induction hypothesis
(Proposition \ref{prop161} (7)).
We have thus checked the assumption of Propositions \ref{prop528} and
\ref{prop1562}.
The Kuranishi structure we obtain in
Proposition \ref{prop1562} is our
Kuranishi structure $\widehat{\mathcal U^+}(\alpha_-,\alpha_+)$.
\par
We next consider CF-perturbations.
We check that we are in Situation \ref{situ529}.
We define $\widehat{\frak S^+_{S_k}}$ in
Situation \ref{situ529} by the right hand side of
(\ref{fiberproductSSinduction}).
Situation \ref{situ529} (2) can be checked easily in our case
by using inductive hypothesis.
We can thus apply Propositions \ref{prop529} and
\ref{prop529rev}. The CF-perturbations we obtain
by
Proposition \ref{prop529rev}
is $\widehat{\frak S^+}(\alpha_-,\alpha_+)$.
\par
Now we will check that $\widehat{\mathcal U^+}(\alpha_-,\alpha_+)$
and $\widehat{\frak S^+}(\alpha_-,\alpha_+)$
satisfy Proposition \ref{prop161} (1)-(8).
\par
(1)  Since $\tau < \tau^+$, the $\tau$-collared-ness of
$\widehat{\mathcal U^+}(\alpha_-,\alpha_+)$ is a consequence of
Proposition \ref{prop528} and
$\widehat{\mathcal U}(\alpha_-,\alpha_+)^{\boxplus\tau}<
\widehat{\mathcal U^+}(\alpha_-,\alpha_+)$
also follows from Proposition \ref{prop528}.
\par
(2) is a consequence of Proposition \ref{prop529rev}.
\par
(3) We apply Proposition \ref{prop528} to each
$\frak G$ equivalence class of the pair $(\alpha_-,\alpha_+)$.
Then for other $(\beta\alpha_-,\beta\alpha_+)$ we
define $\widehat{\mathcal U^+}(\alpha_-,\alpha_+)$ so that it is identified with $(\alpha_-,\alpha_+)$ by using
existence of the periodicity isomorphism on the boundary.
Then existence of the periodicity isomorphism
for $\alpha_-,\alpha_+$ with $E(\alpha_+) - E(\alpha_-) = E_{\frak E}^k$
is immediate from the definition.
The compatibility of the periodicity isomorphism with
CF-perturbations can be proved in the
same way.
\par
(4),(6),(7) This is a consequence of Proposition \ref{prop528}.
Namely it is a consequence of Proposition \ref{prop528} (1) and
the induction hypothesis.
\par
(5),(8) This is a consequence of Proposition \ref{prop529rev} (1) and the
induction hypothesis.
Hence the proof of Proposition \ref{prop161} is now complete.
\end{proof}
\begin{rem}
In the case of partial linear K-system,
the induction to prove Proposition \ref{prop161}
stops in finite steps.
In the case of linear K-system,
the number of inductive steps
is countably infinitely many.
\end{rem}
We next rewrite the geometric conclusion of Proposition \ref{prop161}
into algebraic structures.
\begin{defn}
In the situation of Proposition \ref{prop161}, we define
\begin{equation}
\frak m^{\epsilon}_{1;\alpha_+,\alpha_-} :
\Omega(R_{\alpha_-};o_{R_{\alpha_-}}) \longrightarrow \Omega(R_{\alpha_+};o_{R_{\alpha_+}})
\end{equation}
by
\begin{equation}\label{form16ten6m}
\frak m^{\epsilon}_{1;\alpha_+,\alpha_-}(h)
=
{\rm ev}_{+}!({\rm ev}_{-}^* h;\widehat{\frak S^{+ \epsilon}}(\alpha_-,\alpha_+)).
\end{equation}
Here the right hand side is defined by
Definition \ref{pushoutdeftau} on
the K-space
$$
(\mathcal M(\alpha_-,\alpha_+)^{\boxplus\tau_0},\widehat{\mathcal U^+}(\alpha_-,\alpha_+)).
$$
See Theorem \ref{thm:261} for the correspondence coupled with local systems.
Note that Condition \ref{linsysmainconds} (VII) is compatible with the conventions
\eqref{eq:261}, \eqref{eq:262}.
\end{defn}
The degree of the map $\frak m^{\epsilon}_{1;\alpha_+,\alpha_-}$ is
$$
\dim R_{\alpha_+} - \dim \mathcal M(\alpha_-,\alpha_+)
=
1 - \mu(\alpha_+) + \mu(\alpha_-)
$$
by Condition \ref{linsysmainconds} (VI) and  \cite[Definition 7.78]{part11}.
Therefore after the degree shift in Definition \ref{Fvectspace} (2)
its degree becomes $+1$.
\begin{rem}\label{runningoutsharp}
In the situation of linear K-system
where infinitely many different K-spaces are involved,
we need to make a careful choice of $\epsilon >0$.
Note that the well-definedness of push out (\cite[Theorem 9.14]{part11})
says that for each K-space and its
CF-perturbation the push out is well-defined
for sufficiently small $\epsilon >0$.
In the situation of linear K-system
infimum of such $\epsilon$ over all
$\alpha_{\pm}$ of
$\widehat{\frak S^{+ \epsilon}}(\alpha_-,\alpha_+)$ may
be $0$.
\footnote{This is a version of the `running out' problem
discussed in \cite[Subsection 7.2.8]{fooobook2}.
The way we resolve it in Subsection \ref{subsec:constchaincmps} of this book
is the same as one in \cite{fooobook2}.}
(It might be possible to prove that we can take the same $\epsilon$ for
all $\widehat{\frak S^{+ \epsilon}}(\alpha_-,\alpha_+)$
at the same time.
However, it is cumbersome to formulate the condition
without referring the construction itself.
By this reason, we do not try to prove or use such a uniformity in this book.)
More precisely speaking,
we have the following:
\begin{enumerate}
\item[($\flat$)]
For any energy cut level $E_0$ there exists $\epsilon_0(E_0)$ such that
the operator $\frak m^{\epsilon}_{1;\alpha_+,\alpha_-}$
is defined when $0 <E(\alpha_+) - E(\alpha_-) \le E_0$
and $\epsilon < \epsilon_0(E_0)$.
\end{enumerate}
Hereafter when the relationship between energy cut level and $\epsilon$
(which is the parameter of the approximation) appears in
this way, we write
{\it in the sense of $(\flat)$}\index{in the sense of $(\flat)$} in place of repeating
the above sentence over again.
\end{rem}
\begin{lem}\label{lem165}
The operators $\frak m^{\epsilon}_{1;\alpha_+,\alpha_-}$
satisfy the following equality in the sense of $(\flat)$.
\begin{equation}
d_0\circ \frak m^{\epsilon}_{1;\alpha_+,\alpha_-}
+ \frak m^{\epsilon}_{1;\alpha_+,\alpha_-} \circ d_0
+
\sum_{\alpha; E(\alpha_-) < E(\alpha)
< E(\alpha_+)}
\frak m^{\epsilon}_{1;\alpha_+,\alpha}
\circ
\frak m^{\epsilon}_{1;\alpha,\alpha_-} = 0.
\end{equation}
\end{lem}
Here and hereafter we denote
\begin{equation}\label{eq:deRham}
d_0 (h)=
(-1)^{\dim R_{\alpha} + \mu(\alpha) +1+\deg h} d_{dR} (h)
\end{equation}
for a differential form $h \in \Omega (R_{\alpha})$,
where $d_{dR}$ denotes the de Rham differential on $R_{\alpha}$.
This sign arises from $d_0$ being the classical part of $\frak m_1$ in a filtered $A_{\infty}$ structure. See \eqref{eq:deRham0} and Remark 3.5.8 \cite{fooobook} for the sign.
\begin{proof}
By Stokes' formula and the definition (\ref{form16ten6m}) we have
\begin{equation}\label{parmformind}
(d_0\circ \frak m^{\epsilon}_{1;\alpha_+,\alpha_-}
+ \frak m^{\epsilon}_{1;\alpha_+,\alpha_-} \circ d_0)(h)
=
{\rm ev}_{+}!({\rm ev}_{-}^* h;\partial\widehat{\frak S^{+ \epsilon}}(\alpha_-,\alpha_+)).
\end{equation}
By  Proposition \ref{prop161} (5) and the composition formula, the
right hand side of (\ref{parmformind}) is
$$
\sum_{\alpha; E(\alpha_-) < E(\alpha)
< E(\alpha_+)}
(\frak m^{\epsilon}_{1;\alpha_+,\alpha}
\circ
\frak m^{\epsilon}_{1;\alpha,\alpha_-})(h).
$$
The lemma follows.
\end{proof}
Lemma \ref{lem165} implies that
$$
d_0 + \sum \frak m^{\epsilon}_{1;\alpha_+,\alpha_-}
$$
is a `coboundary oprator modulo higher order term'.
Since we need to take care of the point mentioned in
Remark \ref{runningoutsharp},
we have to stop at some energy cut level.
So we still need to do some more work to
prove Theorem \ref{linesysmainth1} (1).

\subsection{Construction of cochain maps}
\label{subsec:constchaincmps}

We next consider morphism of
linear K-systems or of partial linear K-systems.

\begin{shitu}\label{situ16716}
Suppose we are in Situation \ref{situatinmor}.
We assume that for each $i=1,2$ we have a CF-perturbation
$\widehat{\frak S^+}(i;\alpha_{i-},\alpha_{i+})$ of a $\tau$-collared Kuranishi structure
$\widehat{\mathcal U^+}(i;\alpha_{i-},\alpha_{i+})$
on $\mathcal M^i(\alpha_{i-},\alpha_{i+})^{\boxplus\tau_0}$ such that they satisfy
the conclusion of
Proposition \ref{prop161}.
Here $\mathcal M^i(\alpha_{i-},\alpha_{i+})$ is as in Condition \ref{morphilinsys}.
(From now on, we write $\alpha_{\pm}$ in place of $\alpha_{i\pm}$ if no confusion can occur.)
\par
Let $0 < \tau' < \tau < \tau_0=1$ where $\tau_0, \tau$ are as in Proposition \ref{prop161}.
$\blacksquare$
\end{shitu}

\begin{prop}\label{prop1618}
Suppose we are in Situation \ref{situ16716} and are given a morphism
of partial linear K-systems of energy cut level $E_0$
($\in \R_{\ge 0} \cup \{\infty\}$)  and
energy loss $c$ from $\mathcal F_1$ to $\mathcal F_2$.
Let $\mathcal N(\alpha_1,\alpha_2)$ be its interpolation space
where $E(\alpha_2) - E(\alpha_1) \le E_0$.
\par
Then for each $\tau'$ with $0 < \tau' < \tau < \tau_0=1$ as above,
there exists a $\tau'$-collared Kuranishi structure
$\widehat{\mathcal U^+}({\rm mor};\alpha_1,\alpha_2)$
on $\mathcal N(\alpha_1,\alpha_2)^{\boxplus\tau_0}$
and a  $\tau'$-collared CF-perturbation
$\widehat{\frak S^+}({\rm mor};\alpha_1,\alpha_2)$
of $\widehat{\mathcal U^+}({\rm mor};\alpha_1,\alpha_2)$ such that
they have the following properties:
\begin{enumerate}
\item
Let $\widehat{\mathcal U}({\rm mor};\alpha_1,\alpha_2)$ be the
Kuranishi structure of $\mathcal N(\alpha_1,\alpha_2)$ in Situation \ref{morphilinsys}
(IV).
Then $\widehat{\mathcal U}({\rm mor};\alpha_1,\alpha_2)^{\boxplus\tau_0}
< \widehat{\mathcal U^+}({\rm mor};\alpha_1,\alpha_2)$
as collared Kuranishi structures.
\item
The CF-perturbation $\widehat{\frak S^+}({\rm mor};\alpha_1,\alpha_2)$ is transversal to $0$.
Moreover\footnote{As we note in the footnote of Proposition \ref{prop161} (2), we
simply write ${\rm ev}_+$ in place of ${\rm ev}_+^{\boxplus\tau_0}$.}
${\rm ev}_+ : \mathcal N(\alpha_1,\alpha_2)^{\boxplus\tau_0}
\to R_{\alpha_2}$ is strongly submersive with respect to $\widehat{\frak S^+}({\rm mor};\alpha_1,\alpha_2)$.
\item
We have the following periodicity isomorphism
$$
\widehat{\mathcal U^+}({\rm mor};\alpha_1,\alpha_2)
\longrightarrow \widehat{\mathcal U^+}({\rm mor};\beta\alpha_1,\beta\alpha_2),
$$
which is compatible with the isomorphism in Condition \ref{morphilinsys} (VIII)
via the embedding given in (1). The pull-back of
$\widehat{\frak S^+}({\rm mor};\beta\alpha_1,\beta\alpha_2)$ by this
isomorphism is $\widehat{\frak S^+}({\rm mor};\alpha_1,\alpha_2)$.
\item
There is an isomorphism of $\tau'$-collared K-spaces
\footnote{See Remark \ref{rem:FiberProdOrd} for the sign and the order of the fiber products.}
\begin{equation}\label{formula1211morphrev}
\aligned
&\partial({\mathcal N}({\alpha_1},{\alpha_2})^{\boxplus\tau_0},\widehat{\mathcal U^+}({\rm mor};\alpha_1,\alpha_2))
\\
\cong &\coprod_{\alpha'_1 \in \frak A_1}
(-1)^{\dim {\mathcal N}({\alpha}'_1,{\alpha_2})}
\left(
({\mathcal N}({\alpha}'_1,{\alpha_2})^{\boxplus\tau_0},\widehat{\mathcal U^+}({\rm mor};\alpha'_1,\alpha_2))\right.
\\
&\qquad\qquad\left.{}_{{\rm ev}_{-}}\times_{{\rm ev}_{+}}\,\,
({\mathcal M}^1({\alpha_1},{\alpha'_1})^{\boxplus\tau_0},\widehat{\mathcal U^+}(1;\alpha_1,\alpha'_1))
\right) \\
&\sqcup \coprod_{\alpha'_2 \in \frak A_2}
(-1)^{\dim {\mathcal M}^2({\alpha}'_2,{\alpha_2})}
\left(
({\mathcal M}^2({\alpha}'_2,{\alpha_2})^{\boxplus\tau_0},\widehat{\mathcal U^+}(2;\alpha'_2,\alpha_2))
\right.
\\
&\qquad\qquad\quad\left.
{}_{{\rm ev}_{-}}\times_{{\rm ev}_{+}}\,\,
({\mathcal N}({\alpha_1},{\alpha'_2})^{\boxplus\tau_0},
\widehat{\mathcal U^+}({\rm mor};\alpha_1,\alpha'_2))
\right).
\endaligned
\end{equation}
Here the first union is taken over
$\alpha'_1 \in \frak A_1$ with $E(\alpha_1) < E(\alpha'_1) \le E(\alpha_2) + c$
and the second union is taken over
$\alpha'_2 \in \frak A_2$ with $ E(\alpha_1)- c \le  E(\alpha'_2) < E(\alpha_2)$.
The number $c$ is the energy loss of the given morphism.
\par
The isomorphism (\ref{formula1211morphrev}) is compatible with the isomorphism
induced by (\ref{formula1211morph}) in Condition \ref{morphilinsys} (X) via the embedding  in (1).
\item
The pull-back of $\widehat{\frak S^+}({\rm mor};\alpha_1,\alpha_2)$ by the isomorphism
(\ref{formula1211morphrev})
is equivalent to the fiber product of $\widehat{\frak S^+}(1;\alpha_1,\alpha'_1)$,
$\widehat{\frak S^+}({\rm mor};\alpha'_1,\alpha_2)$
and of
$\widehat{\frak S^+}({\rm mor};\alpha_1,\alpha'_2)$,
$\widehat{\frak S^+}(1;\alpha'_2,\alpha_2)$.
\item
The normalized corner
$\widehat S_k({\mathcal N}({\alpha_1},{\alpha_2})^{\boxplus\tau_0},\widehat{\mathcal U^+}({\rm mor};\alpha_1,\alpha_2))$
is the disjoint union of
\begin{equation}\label{eachtermofSkNrev}
\aligned
&\widehat{\mathcal U^+}(1;\alpha_1,\alpha_{1,1})
\times_{R^1_{\alpha_{1,1}}}  \dots
\times_{R^1_{\alpha_{1,k_1-1}}}
\widehat{\mathcal U^+}(1;{\alpha_{1,k_1-1}},\alpha_{1,k_1}) \\
&\times_{R^1_{\alpha_{1,k_1}}}
\widehat{\mathcal U^+}({\rm mor};\alpha_{1,k_1},\alpha_{2,1}) \\
&\times_{R^1_{\alpha_{2,1}}}\widehat{\mathcal U^+}(2;{\alpha_2},\alpha_{2,1})
\times_{R^2_{\alpha_{2,1}}}
 \dots
\times_{R^2_{\alpha_{2,k_2-1}}}
\widehat{\mathcal U^+}(2;{\alpha_{2,k_2-1}},\alpha_{2,k_2})
\endaligned
\end{equation}
where $k_1+k_2 = k$, $\alpha_{1,i} \in \frak A_1$, $\alpha_{2,i} \in \frak A_2$.
\par
This isomorphism is compatible with the isomorphism in Condition \ref{morphilinsys} (XI)  via the the embedding in (1).
\item
(\ref{cornecominduction}) and
(\ref{eachtermofSkNrev}) imply that
$\widehat S_{\ell}(\widehat S_k({\mathcal N}({\alpha_1},{\alpha_2})^{\boxplus\tau_0},\widehat{\mathcal U^+}({\rm mor};\alpha_1,\alpha_2)))$
is a fiber product similar to (\ref{eachtermofSkNrev})
with $k$ replaced by $k+\ell$.
\par
Moreover
$\widehat S_{k+\ell}({\mathcal N}({\alpha_1},{\alpha_2})^{\boxplus\tau_0},\widehat{\mathcal U^+}({\rm mor};\alpha_1,\alpha_2))$
is also a fiber product similar to (\ref{eachtermofSkNrev})
with $k$ replaced by $k+\ell$.
The map
$$
\pi_{\ell,k} : \widehat S_{\ell}(\widehat S_k({\mathcal N}({\alpha_1},{\alpha_2})^{\boxplus\tau_0},\widehat{\mathcal U^+}({\rm mor};\alpha_1,\alpha_2)))
\to \widehat S_{k+\ell}({\mathcal N}({\alpha_1},{\alpha_2})^{\boxplus\tau_0},\widehat{\mathcal U^+}({\rm mor};\alpha_1,\alpha_2))
$$
in Proposition \ref{prop2813}  becomes the identity map via those identifications.
\item
The pull-back of the restriction of $\widehat{\frak S^+}({\rm mor};\alpha_1,\alpha_2)$
to $\widehat S_k(\widehat{\mathcal U^+}({\rm mor};\alpha_1,\alpha_2))$
by the isomorphism in (6)
is equivalent to the fiber product of $\widehat{\frak S^+}({\rm mor};*,*)$,  $\widehat{\frak S^+}(1;*,*)$,
 $\widehat{\frak S^+}(2;*,*)$ along (\ref{eachtermofSkNrev}).
\item
Suppose we are given a uniform family
$\widehat{\frak S^+}(i;\alpha_-,\alpha_+)^{\sigma}$
of CF-perturbations as in Condition \ref{situ16716}.
Then the $\sigma$-parametrized
family of CF-perturbations
$\widehat{\frak S^+}({\rm mor};\alpha_1,\alpha_2)^{\sigma}$
satisfying (5)(6)(7) for each $\sigma$ is also uniform.
\end{enumerate}
\end{prop}
\begin{rem}
Note in Formula (\ref{eachtermofSkNrev}) we simplify the notation and omit the
symbol of the underlying topological space.
Namely we write $\widehat{\mathcal U^+}(1;\alpha_1,\alpha_{1,1})$
in place of $({\mathcal M}^1({\alpha_1},{\alpha_{1,1}})^{\boxplus\tau_0},\widehat{\mathcal U^+}(1;\alpha_1,\alpha_{1,1}))$.
From now on, we use this kinds of simplified notation when no confusion can occur.
\end{rem}
\begin{proof}
The proof is by induction on $k$. It is entirely similar
to the proof of Proposition \ref{prop161}. So we omit it.
\end{proof}
We next rewrite the geometric conclusion of Proposition \ref{prop1618}
to algebraic language.
\begin{defn}\label{chainmapobtainedfromcor}
In the situation of Proposition \ref{prop1618}, we define
\begin{equation}
\psi^{\epsilon}_{\alpha_2,\alpha_1} :
\Omega(R_{\alpha_1};o_{R_{\alpha_1}}) \longrightarrow
\Omega(R_{\alpha_2};o_{R_{\alpha_2}})
\end{equation}
by
\begin{equation}\label{form16ten6}
\psi^{\epsilon}_{\alpha_2,\alpha_1}(h)
=
{\rm ev}_{+}!({\rm ev}_{-}^* h;\widehat{\frak S^{+ \epsilon}}({\rm mor};\alpha_1,\alpha_2)).
\end{equation}
Here the right hand side is defined by
Definition \ref{pushoutdeftau} on
the K-space
$$
({\mathcal N}({\alpha_1},{\alpha_2})^{\boxplus\tau_0},
\widehat{\mathcal U^+}({\rm mor};\alpha_1,\alpha_2)).
$$
\end{defn}
By Condition \ref{morphilinsys} (VI) and \cite[Definition 7.78]{part11}, the
degree of $\psi^{\epsilon}_{\alpha_1,\alpha_2}$ is
$\mu(\alpha_1) - \mu(\alpha_2)$.
Therefore, after the degree shift as in Definition \ref{Fvectspace} (2),
its degree becomes $0$.
If the energy loss of our morphism $\frak N$ is
$c$, the family $\{ \psi^{\epsilon}_{\alpha_1,\alpha_2}\}$ of maps induces
$$
\frak F^{\lambda}CF(\mathcal F_1) \longrightarrow
\frak F^{\lambda-c}CF(\mathcal F_2)
$$
where the filtration $\frak F^{\lambda}$ is defined in
Definition \ref{Fvectspace} (2)(3).
\begin{lem}\label{lem165rev}
The operators $\{ \psi^{\epsilon}_{\alpha_2,\alpha_1}\}$
satisfy the following equality in the sense of $(\flat)$.
\begin{equation}\label{form1613}
\aligned
&d_0\circ \psi^{\epsilon}_{\alpha_2,\alpha_1}
- \psi^{\epsilon}_{\alpha_2,\alpha_1} \circ d_0
\\
&+
\sum_{\alpha'_2}
\frak m^{2,\epsilon}_{1;\alpha_2,\alpha'_2}
\circ
\psi^{\epsilon}_{\alpha'_2,\alpha_1}
-
\sum_{\alpha'_1}
\psi^{\epsilon}_{\alpha_2,\alpha'_1}
\circ
\frak m^{1,\epsilon}_{1;\alpha'_1,\alpha_1}
= 0.
\endaligned
\end{equation}
Here the first sum in the second line is taken over
$\alpha'_2 \in \frak A_2$ with $ E(\alpha_1)- c \le  E(\alpha'_2) < E(\alpha_2)$
and the second sum in the second line is taken over
$\alpha'_1 \in \frak A_1$ with $E(\alpha_1) < E(\alpha'_1) \le E(\alpha_2) + c$.
Here $c\ge 0$ is the energy loss of our morphism.
\end{lem}
\begin{proof}
By Stokes' formula the sum
$d_0\circ \psi^{\epsilon}_{\alpha_2,\alpha_1}
- \psi^{\epsilon}_{\alpha_2,\alpha_1} \circ d_0$
is equal to the correspondence induced by
the boundary of
$\widehat{\frak S^{+}}({\rm mor};\alpha_1,\alpha_2)$.
By Proposition \ref{prop1618} (4) and
the composition formula this is equal to the
second line of (\ref{form1613}).
\end{proof}

\subsection{Proof of Theorem \ref{linesysmainth1} (1)
and Theorem \ref{indlinesysmainth1} (1)}
\label{subsec:proofsec14main1}

In this subsection we will prove Theorem \ref{linesysmainth1} (1).
We will also prove Theorem \ref{indlinesysmainth1} (1) at the same time.
\begin{shitu}\label{situ1612}
We study a partial linear K-system $\mathcal F^i$.
\par
We define $
\frak E$
as a set of $E \in \R$ such that one of the following holds.
\begin{enumerate}
\item
There exist $\alpha_-,\alpha_+$ and $i$ such that
$\mathcal N^i(\alpha_-,\alpha_+) \ne \emptyset$
and $E(\alpha_+) - E(\alpha_-) = E$.
\item
There exist $\alpha_-,\alpha_+$ and $i$ such that
$\mathcal M^i(\alpha_-,\alpha_+) \ne \emptyset$
and $E(\alpha_+) - E(\alpha_-) = E$.
\end{enumerate}
Then $\frak E$ is a discrete set by Definition \ref{defn1528} (2)(g).
We put
$$
\frak E = \{E_{\frak E}^{1},E_{\frak E}^{2},\dots\}
$$ such that $0 < E_{\frak E}^{1} < E_{\frak E}^{2} < \dots$.
Note that the graded and filtered vector space $CF(\mathcal F^i)$ as in Definition \ref{Fvectspace} is independent of $i$. We denote it by $CF(\mathcal F)$.
$\blacksquare$
\end{shitu}
We observe that for each fixed $i$ the set $\frak E \cap [0,E_0]$ can be
strictly bigger than the set $\frak E_{\le E_0}$ in (\ref{16111}).
Nevertheless we can replace $\frak E_{\le E_0}$
by $\frak E \cap [0,E_0]$ in Proposition \ref{prop161} etc.,
by putting $\frak m^{\epsilon}_{1;\alpha_+,\alpha_-} = 0$
when $\mathcal M^i(\alpha_-,\alpha_+) = \emptyset$.
\par
In Situation \ref{situ1612}, we
take $i \in \Z_+$ and study the relationship between the operators
$\frak m^{i,\epsilon}_{1;\alpha_+,\alpha_-}$ and
$\frak m^{i+1,\epsilon}_{1;\alpha_+,\alpha_-}$
defined by applying Proposition \ref{prop161} and
Lemma \ref{lem165} to partial linear K-systems
$\mathcal F^i$ and $\mathcal F^{i+1}$.
The next definition will be applied
with
$\frak m^1_{1;\alpha_+,\alpha_-}$ and $\frak m^2_{1;\alpha_+,\alpha_-}$ replaced by
$\frak m^{i,\epsilon}_{1;\alpha_+,\alpha_-}$ and
$\frak m^{i+1,\epsilon}_{1;\alpha_+,\alpha_-}$
respectively.
\begin{defn}\label{defn1613}
Suppose we are in Situation \ref{situ1612}.
\begin{enumerate}
\item
A {\it partial cochain complex structure on $CF(\mathcal F)$ of energy cut level $E_0$}\index{energy cut level!partial cochain complex structure}
\index{partial ! cochain complex} assigns
$\frak m_{1;\alpha_+,\alpha_-}$ to each $\alpha_+, \alpha_-$ with  $0 <E(\alpha_+) - E(\alpha_-) \le E_0$
such that (\ref{lem165}) is satisfied.
\item
Let
$(CF(\mathcal F),\{\frak m^j_{1;\alpha_+,\alpha_-}\})$ be a partial cochain complex structure on
$CF(\mathcal F)$ of energy cut level $E_0$
for $j=1,2$.
A {\it partial cochain map of energy cut level $E$ with energy loss $0$}\index{energy cut level!partial cochain map}
\index{partial ! cochain map}
from $(CF(\mathcal F),\{\frak m^1_{1;\alpha_+,\alpha_-}\})$ to
$(CF(\mathcal F),\{\frak m^2_{1;\alpha_+,\alpha_-}\})$ assigns
a map $\psi_{\alpha_2,\alpha_1}$ to each $\alpha_1,\alpha_2$ with
$0 \le E(\alpha_2) - E(\alpha_1) \le E_0$ such that (\ref{form1613})
is satisfied.
For $\alpha_1 \ne \alpha_2$ with $E(\alpha_1) =E(\alpha_2)$, we assign $\psi_{\alpha_2,\alpha_1}=0$.
For $\alpha_1 = \alpha_2$, we assign $\psi_{\alpha_2,\alpha_1}={\rm identity}$.
\item
If $(CF(\mathcal F),\{\frak m_{1;\alpha_+,\alpha_-}\})$ is a partial cochain complex structure on $CF(\mathcal F)$ of energy cut level $E_0$
and $E'_0 < E_0$, then by forgetting a part of the operations $\frak m_{1;\alpha_+,\alpha_-}$ we may regard
$(CF(\mathcal F),\{\frak m_{1;\alpha_+,\alpha_-}\})$ as a
partial cochain complex structure on $CF(\mathcal F)$ of energy cut level $E'_0$.
We call it the {\it reduction by energy cut at $E'_0$}.\index{energy cut level!reduction by energy cut}
\par
We can define a reduction by energy cut at $E'_0$ of a
partial cochain map  of energy cut level $E_0$ with energy loss $0$ in the same way.
\item
Let $(CF(\mathcal F),\{\frak m_{1;\alpha_+,\alpha_-}\})$ be a partial cochain complex structure of energy cut level $E'$.
A partial cochain complex structure of energy cut level $E$ is said to be its
{\it promotion}\index{energy cut level!promotion}\index{promotion}
if its reduction by energy cut at $E_0$ is $(CF(\mathcal F),\{\frak m_{1;\alpha_+,\alpha_-}\})$.
\par
A {\it promotion} of partial cochain map is defined in the same way.
\end{enumerate}
\end{defn}

The next lemma is a baby version of \cite[Lemma 7.2.72]{fooobook2}.

\begin{lem}\label{lem1623}
Let $(CF(\mathcal F),\{\frak m^j_{1;\alpha_+,\alpha_-}\})$ be partial cochain complexes of energy cut level $E_{\frak E}^{k_j}$
for $j=1,2$.
Suppose $k_1 < k_2$.
Let $\{\psi_{\alpha_2,\alpha_1}\}$ be a partial cochain map of energy cut level $E_{\frak E}^{k_1}$
from $(CF(\mathcal F),\{\frak m^1_{1;\alpha_+,\alpha_-}\})$ to a reduction by energy cut at
$E_{\frak E}^{k_1}$ of
$(CF(\mathcal F),\{\frak m^2_{1;\alpha_+,\alpha_-}\})$.
Then there exists a promotion of  $(CF(\mathcal F),\{\frak m^1_{1;\alpha_+,\alpha_-}\})$
to energy cut level $E_{\frak E}^{k_2}$
and a promotion of $\{\psi_{\alpha_2,\alpha_1}\}$
to the energy cut level  $E_{\frak E}^{k_2}$ from this promoted partial cochain complex structure.
\end{lem}
\begin{proof}
By an obvious induction argument it suffices to prove the case $k_2 = k_1+1$.
Suppose $E(\alpha_+) - E(\alpha_-) = E_{\frak E}^{k_2}$.
We define a linear map
$$
o(\alpha_+,\alpha_-) : \Omega(R_{\alpha_-};o_{R_{\alpha_-}})
\longrightarrow \Omega(R_{\alpha_+};o_{R_{\alpha_+}})
$$
by
$$
o(\alpha_+,\alpha_-)
=
\sum_{\alpha; E(\alpha_-) < E(\alpha)
< E(\alpha_+)}
(\frak m^{1}_{1;\alpha_+,\alpha}
\circ
\frak m^{1}_{1;\alpha,\alpha_-}).
$$
We will prove that $o(\alpha_+,\alpha_-)$ is a $d_0$-coboundary.
\begin{notation}
We use the following notation.
For an $\R$ linear map $F : CF(\mathcal F) \to CF(\mathcal F)$,
$F_{\alpha_+,\alpha_-}$ denotes the
${\rm Hom}_{\R}(\Omega(R_{\alpha_-},o_{R_{\alpha_-}}),\Omega(R_{\alpha_+},o_{R_{\alpha_+}})
)$ component of $F$.
\end{notation}
We put $\psi_{\alpha,\alpha} = {\rm id}$.
We define
$\hat d^j : CF(\mathcal F) \to CF(\mathcal F)$
($j=1,2$)
and  $\widehat\psi : CF(\mathcal F) \to CF(\mathcal F)$ by
\begin{equation}\label{defnhatd}
\hat d^j = d_0 \oplus \bigoplus_{\alpha_1,\alpha_2
\atop
E({\alpha_2}) - E({\alpha_1}) \le E_{k_1}} \frak m^{j}_{1;\alpha_2,\alpha_1},
\quad \widehat\psi = \bigoplus_{\alpha_1,\alpha_2
\atop
E({\alpha_2}) - E({\alpha_1}) \le E_{k_1}} \psi^{}_{\alpha_2,\alpha_1}.
\end{equation}
We have
\begin{equation}\label{redefo}
(\hat d^1\circ \hat d^1)_{\alpha_+,\alpha_-}
= o(\alpha_+,\alpha_-),
\end{equation}
if $E({\alpha_+}) - E({\alpha_-}) = E_{\frak E}^{k_2}$.
On the other hand,
we have
\begin{equation}\label{dddis0}
(\hat d^1\circ \hat d^1)_{\alpha_2,\alpha_1}
= 0
\end{equation}
if $E({\alpha_2}) - E({\alpha_1}) \le E_{\frak E}^{k_1}$.
We note
\begin{equation}\label{trivialidentityddd}
(\hat d^1\circ \hat d^1) \circ \hat d^1
- \hat d^1 \circ(\hat d^1\circ \hat d^1)
=0.
\end{equation}
Then (\ref{dddis0}) implies
$$
\aligned
((\hat d^1\circ \hat d^1) \circ \hat d^1)_{\alpha_+,\alpha_-}
&=
(\hat d^1\circ \hat d^1)_{\alpha_+,\alpha_-} \circ d_0, \\
(\hat d^1 \circ (\hat d^1\circ \hat d^1))_{\alpha_+,\alpha_-}
&=
d_0 \circ (\hat d^1\circ \hat d^1)_{\alpha_+,\alpha_-}.
\endaligned$$
Therefore (\ref{redefo}) and (\ref{trivialidentityddd}) imply
$$
d_0\circ o(\alpha_+,\alpha_-)
-o(\alpha_+,\alpha_-) \circ d_0 =0.
$$
Namely $o(\alpha_+,\alpha_-)$ is a cocycle.
\par
We continue to show that $o(\alpha_+,\alpha_-)$ is a coboundary.
When
$E({\alpha_+}) - E({\alpha_-}) = E_{\frak E}^{k_2}$,
we put
$$
b(\alpha_+,\alpha_-)
=
(\hat d^2\circ \widehat\psi - \widehat\psi\circ \hat d^1)_{\alpha_+,\alpha_-}.
$$
Since $\widehat\psi$ is assumed to be a cochain map of energy cut level
$E_{\frak E}^{k_1}$, we have
$$
(\hat d^2\circ \widehat\psi - \widehat\psi\circ \hat d^1)_{\alpha_2,\alpha_1}
= 0
$$
if $E({\alpha_2}) - E({\alpha_1}) \le E_{\frak E}^{k_1}$.
Therefore we find
\begin{equation}\label{form1618}
\aligned
&(\hat d^2 \circ (\hat d^2\circ \widehat\psi - \widehat\psi\circ \hat d^1)
+
(\hat d^2\circ \widehat\psi - \widehat\psi\circ \hat d^1) \circ \hat d^1)_{\alpha_+,\alpha_-}\\
&= d_0 \circ b(\alpha_+,\alpha_-) + b(\alpha_+,\alpha_-)\circ d_0.
\endaligned
\end{equation}
On the other hand, an obvious calculation leads to
\begin{equation}\label{form1619}
\aligned
&(\hat d^2 \circ (\hat d^2\circ \widehat\psi - \widehat\psi\circ \hat d^1)
+
(\hat d^2\circ \widehat\psi - \widehat\psi\circ \hat d^1) \circ \hat d^1)_{\alpha_+,\alpha_-}\\
&= ((\hat d^2 \circ \hat d^2) \circ \widehat\psi - \widehat\psi\circ (\hat d^1
\circ \hat d^1))_{\alpha_+,\alpha_-} \\
&=
(\hat d^2 \circ \hat d^2)_{\alpha_+,\alpha_-} \circ \widehat\psi_{\alpha_-,\alpha_-}
-
\widehat\psi_{\alpha_+,\alpha_+} \circ (\hat d^1 \circ \hat d^1)_{\alpha_+,\alpha_-}
\\
&=
(\hat d^2 \circ \hat d^2)_{\alpha_+,\alpha_-}
- o(\alpha_+,\alpha_-).
\endaligned
\end{equation}
We observe that
the equation $\hat d^2 \circ \hat d^2 = 0$ gives rise to
\begin{equation}\label{form1620}
d_0 \circ \frak m^{2}_{1;\alpha_+,\alpha_-}
+ \frak m^{2}_{1;\alpha_+,\alpha_-} \circ d_0
+ (\hat d^2 \circ \hat d^2)_{\alpha_+,\alpha_-}
= 0,
\end{equation}
since the energy cut level of
$(CF(\mathcal F),\{\frak m^2_{1;\alpha_+,\alpha_-}\})$
is $E_{\frak E}^{k_2}$.
Combination of (\ref{form1618}), (\ref{form1619}), (\ref{form1620})
implies
\begin{equation}\label{form1621}
o(\alpha_+,\alpha_-)
=
d_0 \circ (\frak m^{2}_{1;\alpha_+,\alpha_-} - b(\alpha_+,\alpha_-))
- (\frak m^{2}_{1;\alpha_+,\alpha_-} - b(\alpha_+,\alpha_-)) \circ d_0.
\end{equation}
Namely $o(\alpha_+,\alpha_-)$ is a $d_0$-coboundary.
Therefore there exists $\frak m^{1}_{1;\alpha_+,\alpha_-}$ such that
\begin{equation}\label{form1622}
d_0 \circ \frak m^{1}_{1;\alpha_+,\alpha_-}
+ \frak m^{1}_{1;\alpha_+,\alpha_-}\circ d_0
+ o(\alpha_+,\alpha_-)
= 0.
\end{equation}
Hence we have a promotion of
$(CF(\mathcal F),\{\frak m^1_{1;\alpha_+,\alpha_-}\})$
to the energy level $E_{\frak E}^{k_2}$.
\par
We note that the choice of
$\frak m^{1}_{1;\alpha_+,\alpha_-}$
satisfying (\ref{form1622})
is not unique and
can be changed by adding a $d_0$-cocycle.
We will use this freedom in the next part of the proof.
\par
We next promote $\psi$.
For this purpose it suffices to find $\psi_{\alpha_+,\alpha_-}$
for $E(\alpha_+) - E(\alpha_-) = E_{\frak E}^{k_2} - c$
such that
\begin{equation}\label{form1623}
\aligned
d_0 \circ \psi_{\alpha_+,\alpha_-}
- \psi_{\alpha_+,\alpha_-} \circ d_0
&+
\frak m^{2}_{1;\alpha_+,\alpha_-} \circ \psi_{\alpha_-,\alpha_-}
\\
&-
\psi_{\alpha_+,\alpha_+}
\circ
\frak m^{1}_{1;\alpha_+,\alpha_-}
+
b(\alpha_+,\alpha_-)
= 0.
\endaligned
\end{equation}
We will prove it below.
We put
$$
\aligned
o'(\alpha_+,\alpha_-)
&=
\frak m^{2}_{1;\alpha_+,\alpha_-} \circ \psi_{\alpha_-,\alpha_-}
-
\psi_{\alpha_+,\alpha_+}
\circ
\frak m^{1}_{1;\alpha_+,\alpha_-} +
b(\alpha_+,\alpha_-)
\\
&=\frak m^{2}_{1;\alpha_+,\alpha_-}
- \frak m^{1}_{1;\alpha_+,\alpha_-} +
b(\alpha_+,\alpha_-).
\endaligned
$$
The identities (\ref{form1621}) and (\ref{form1622}) imply
$$
d_0\circ o'(\alpha_+,\alpha_-) + o'(\alpha_+,\alpha_-) \circ d_0
=0.
$$
Thus $o'(\alpha_+,\alpha_-)$ is a $d_0$-cocycle.
Using the freedom of the choice of
$\frak m^{1}_{1;\alpha_+,\alpha_-}$ we mentioned above,
we may assume that $o'(\alpha_+,\alpha_-)$
is a $d_0$-coboundary.
Hence we can find $\psi_{\alpha_+,\alpha_-}$
satisfying (\ref{form1623}).
The proof of Lemma \ref{lem1623} is complete.
\end{proof}
\begin{proof}[Proof of Theorem \ref{linesysmainth1}  (1)
and Theorem \ref{indlinesysmainth1} (1)]
Below we will prove Theorem \ref{indlinesysmainth1} (1)
and indicate modifications needed to prove
Theorem \ref{linesysmainth1} (1).
\par
Suppose we are in Situation \ref{situ1612}.
(We note that
for the proof of
Theorem \ref{linesysmainth1} (1)
we consider the situation where $\mathcal F^i =\mathcal F$
and the morphisms $\frak N^i$ appearing in Definition \ref{defn1528}
(2)(d) are the identity morphisms for all $i$.)
\par
We assume that the energy cut level of $\mathcal F^i$  is $E_{\frak E}^{k_i}$
such that $k_i < k_{i+1}$.
(It implies $\lim_{i\to\infty}E_{\frak E}^{k_i} = \infty$ by the discreteness of $\frak E$.)
We apply Proposition \ref{prop161} and Lemma \ref{lem165}
to $\mathcal F^i$ for each $i$.
Then we obtain $\epsilon_{0,i}$
such that for each $\epsilon < \epsilon_{0,i}$ and $i$ we obtain
a partial cochain complex of energy cut level $E_{\frak E}^{k_i}$,
which we denote by
$(CF(\mathcal F),\{\frak m^{i,\epsilon}_{1;\alpha_+,\alpha_-}\})$.
Note that the operator $\frak m^{i,\epsilon}_{1;\alpha_+,\alpha_-}$
is defined as the correspondence of a Kuranishi structure
and a CF-perturbation.
We write them as
$(\widehat{\mathcal U^{+}}(i;\alpha_-,\alpha_+),\widehat{\frak S^{+}}(i;\alpha_-,\alpha_+))$ respectively.
We denote
\begin{equation}\label{form1624}
CF(\mathcal F^i;\epsilon) = (CF(\mathcal F),\{\frak m^{i,\epsilon}_{1;\alpha_+,\alpha_-}\}).
\end{equation}
This is well-defined if $\epsilon < \epsilon_{0,i}$,
which is a partial cochain complex of energy cut level $E_{\frak E}^{k_i}$
by Lemma \ref{lem165}.
\par
We next consider the morphisms
$\frak N^i : \mathcal F^i \to \mathcal F^{i+1}$.
Their interpolation spaces are denoted by $\mathcal N(i;\alpha_-,\alpha_+)$.
We apply Proposition \ref{prop1618}
to obtain a $\tau_{i+1}$-collared Kuranishi structure
on it, which we write $\widehat{\mathcal U^{+}}({\rm mor},i;\alpha_-,\alpha_+)$.
We note that here we regard
$\widehat{\mathcal U^{+}}(i;\alpha_-,\alpha_+)$
as a $\tau_{i+1}$-collared Kuranishi structure.
This is possible since $\tau_{i+1} < \tau_i$ and
$\widehat{\mathcal U^{+}}(i;\alpha_-,\alpha_+)$ is
$\tau_i$-collared.
\par
We next take $\rho_i > 0$ such that
\begin{equation}\label{choiceofrho}
\rho_i \le \min\{\epsilon_{k,i+1}/\epsilon_{k,i} \mid k=0,1,2,\dots\}
\end{equation}
where $\epsilon_{k,i+1}$, $\epsilon_{k,i}$,
$k=1,2,\dots$  will be defined later.
\begin{rem}\label{rem:1815}
There appear only a finite number of $k$'s. Note that
$\epsilon_{0,i}$ is already defined above and $\epsilon_{1,i}$, $\epsilon_{2,i}$ and
$\epsilon_{3,i}$
will be chosen after
Lemma \ref{lem1615}, Lemma \ref{17222lem} and Lemma \ref{lem165rev1} below, respectively.
Other $\epsilon_{4,i}$ etc will be taken during the
proof of Theorem \ref{indlinesysmainth1} (3) in
Subsection \ref{subsec:proofsec14main2} for studying homotopy of homotopies.
\end{rem}
We consider a CF-perturbation
$\epsilon \mapsto \widehat{\frak S^{+ \rho_i\epsilon}} (i+1;\alpha_-,\alpha_+)$
of $\widehat{\mathcal U^{+}}(i+1;\alpha_-,\alpha_+)$,
which is defined as follows.
Note $\widehat{\frak S^{+}}(i;\alpha_-,\alpha_+)$
is a CF-perturbation of
$\widehat{\mathcal U^{+}}(i+1;\alpha_-,\alpha_+)$.
Its local representative on the Kuranishi charts is
$(W_{\frak r},\omega_{\frak r},\frak s_{\frak r}^{\epsilon})$.
(See  \cite[Definition 7.3]{part11}.)
We replace $\frak s_{\frak r}^{\epsilon}$ by $\frak s_{\frak r}^{\epsilon\rho_i}$
but do not change anything else. It is easy to find that it is
compatible with the coordinate change etc. and defines a
CF-perturbation, which we denoted by
$\epsilon \mapsto \widehat{\frak S^{+ \rho_i\epsilon}} (i+1;\alpha_-,\alpha_+)$
above.
Hereafter we write it as
$\widehat{\frak S^{+ \rho_i \cdot}} (i+1;\alpha_-,\alpha_+)$.
By definition we have
\begin{equation}\label{form1626}
f!(h;(\widehat{\frak S^{+ \rho_i \cdot}} (i+1;\alpha_-,\alpha_+)^{\epsilon})
=
f!(h;(\widehat{\frak S^{+ \rho_i \epsilon}} (i+1;\alpha_-,\alpha_+))
\end{equation}
when the integrations along the fibers of the both hand sides are defined.
\par
Now we apply Proposition \ref{prop1618} to obtain
$\widehat{\frak S^{+i}}$ on
$\widehat{\mathcal U^{+}}({\rm mor};\alpha_-,\alpha_+)$.
Here we use the CF-perturbations
$\widehat{\frak S^{+ }} (i;\alpha_-,\alpha_+)$
on the space of connecting orbits of $\mathcal F^i$
and
$\widehat{\frak S^{+ \rho_i \cdot}} (i+1;\alpha_-,\alpha_+)$
on the space of connecting orbits of $\mathcal F^{i+1}$.
(We put $\rho_i\cdot$ for the second one.)
Then we obtain from Proposition \ref{prop1618} a
CF-perturbation on
$\widehat{\mathcal U^{+}}({\rm mor};\alpha_-,\alpha_+)$.
We denote it by
$\widehat{\frak S^{+i}}
({\rm mor},i,\sigma_i;\alpha_-,\alpha_+))$.
\begin{lem}\label{lem1615}
The family $\widehat{\frak S^{+i}}
({\rm mor},i,\rho_i;\alpha_-,\alpha_+))$
for $\rho_i\in (0,\epsilon_{1,i+1}/\epsilon_{1,i}]$
is a uniform family in the sense of \cite[Definition 9.28]{part11}.
\end{lem}
\begin{proof}
This is immediate from Proposition \ref{prop1618} (9).
\end{proof}
By \cite[Proposition 7.88]{part11} and Lemma \ref{lem1615} we can find $\epsilon_{1,i}$ such that
if $\epsilon < \epsilon_{1,i}$
the integration along the fiber is defined by using the
CF-perturbation
$\widehat{\frak S^{+i}}
({\rm mor},i,\rho_i;\alpha_-,\alpha_+))^{\epsilon}$
for $\epsilon < \epsilon_{1,i+1}$.
\begin{rem}\label{rem171717}
The fact that we can take $\epsilon$ in a way independent of $\rho_i$
is crucial here. Otherwise, the process of defining those numbers
would become circular.
Namely while working in this step of induction,
we do not know in advance how small $\epsilon_{0,i+1}$, $\epsilon_{1,i+1}$
must be. So we cannot estimate $\rho_i$ {\it from below} at this
stage.
\end{rem}
Thus for $\epsilon < \min(\epsilon_{0,i},\epsilon_{1,i})$
we obtain $\psi^{i}_{\alpha_+,\alpha_-}$ by (\ref{form16ten6}).
We use (\ref{form1626})
together with Lemma \ref{lem165rev}
to prove that $\{\psi^{i}_{\alpha_+,\alpha_-}\} :
CF(\mathcal F^i;\epsilon)  \to CF(\mathcal F^{i+1};\epsilon)$
is a partial cochain map of energy cut level $E_{\frak E}^{k_i}$.
(Here we reduce the energy cut level of
$CF(\mathcal F^{i};\epsilon)$ to $E_{\frak E}^{k_i}$.)
\par Now we use Lemma \ref{lem1623} inductively
to promote $CF(\mathcal F^{i+1};\epsilon)$
to a partial cochain complex of energy cut level $E_{\frak E}^{k_n}$
for each $n$ and $\{\psi^{i}_{\alpha_+,\alpha_-}\}$ to
a partial cochain map of energy cut level $E_{\frak E}^{k_n}$ for each $n$.
\par
To prove Theorem \ref{linesysmainth1}  (1),
we regard $\mathcal F = \mathcal F^i$
and use the identity morphism $\mathcal F \to \mathcal F$
as $\frak N^i : \mathcal F^i \to \mathcal F^{i+1}$
for all $i$.
Then Theorem \ref{indlinesysmainth1} (1)
implies Theorem \ref{linesysmainth1}  (1).
\par
The proofs of Theorem \ref{linesysmainth1}  (1)
and Theorem \ref{indlinesysmainth1} (1) are now complete.
\end{proof}

\subsection{Composition of morphisms and of induced cochain maps}
\label{subsec:compochain}

In this subsection we show that the composition of morphisms
(defined in Lemma-Definition \ref{lemdef1434}
and in Subsection \ref{subsec:complinkurasmcorner})
induces the composition of the partial cochain maps
given by Definition \ref{chainmapobtainedfromcor}.
Since the partial cochain map in Definition \ref{chainmapobtainedfromcor}
depends on the choice of the perturbation, we need to
state it a bit carefully.

We consider a situation similar but slightly different from
Situation \ref{situatinmor}.
Namely;

\begin{shitu}\label{situatinmorrev}
\begin{enumerate}
\item
For $j=1,2,3$,
let
$$
\mathcal C_j =\Big(
\frak A_j, \frak G_j, \{ R^j_{\alpha_j}\}_{\alpha \in \frak A_j},
\{ o_{R^j_{\alpha_j}} \}_{\alpha \in \frak A_j}, E, \mu,
\{ {\rm PI}^j_{\beta_j,\alpha_j} \}_{\beta_j \in \frak G_j, \alpha_j \in \frak A_j}
\Big)
$$
be critical submanifold data
and
$$
\mathcal F_j = \Big( {\mathcal C}_j,
\{\mathcal M^j(\alpha_{j-},\alpha_{j+}) \}_{\alpha_{j\pm} \in \frak A_j},
({\rm ev}_{-}, {\rm ev}_{+}),
\{{\rm OI}^j_{\alpha_{j-}, \alpha_{j+}}\}_{\alpha_{j\pm} \in \frak A_j},
\{ {\rm PI}^j_{\beta_j;\alpha_{j-},\alpha_{j+}} \}_{\beta_j \in \frak G_j,
\alpha_{j\pm} \in \frak A_j}
\Big)
$$
linear K-systems.
We assume $\frak G_1 = \frak G_2 = \frak G_3$
(together with energy $E$ and the Maslov index $\mu$ on it) and denote it by $\frak G$.
\item
The same as (1) except the assumption
that they consist of
partial linear K-systems of energy cut
level $E_0$. $\blacksquare$
\end{enumerate}
\end{shitu}

\begin{shitu}\label{situ1671699}
Suppose we are in Situation \ref{situatinmorrev} (2)
\par
We assume that for each $j=1,2,3$ we have
$\widehat{\frak S^+}(j;\alpha_{j-},\alpha_{j+})$ of
$\widehat{\mathcal U^+}(j;\alpha_{j-},\alpha_{j+})$
on $\mathcal M^j(\alpha_{j-},\alpha_{j+})^{\boxplus\tau_0}$ satisfying
the conclusions of
Proposition \ref{prop161}.
Here $\mathcal M^j(\alpha_{j-},\alpha_{j+})$ is as in Situation  \ref{situatinmorrev}.
(From now on, we write $\alpha_{\pm}$ in place of $\alpha_{j\pm}$ if no confusion can occur.)
We assume that we have a partial morphism
$\frak N_{j+1 j} : \mathcal F_j \to \mathcal F_{j+1}$ for $j=1,2$,
whose interpolation space is
$\mathcal N_{j j+1}(\alpha_j,\alpha_{j+1})$.
\par
Furthermore, we have $\widehat{\frak S^+}({\rm mor};j+1,j;\alpha_j,\alpha_{j+1})$
and  $\widehat{\mathcal  U^+}({\rm mor};j+1,j;\alpha_{j},\alpha_{j+1})$
on $\mathcal N_{j j+1}(\alpha_{j},\alpha_{j+1})$
satisfying the conclusions of Proposition \ref{prop1618}.
$\blacksquare$
\end{shitu}
We obtain a composition
$
\frak N_{31} = \frak N_{32}\circ \frak N_{21}
$
by Lemma-Definition \ref{lemdef1434} and
Lemma-Definition \ref{1638defken}.
By (\ref{form1620})
we have
\begin{equation}\label{1730form}
\aligned
&\widehat{\mathcal  U}({\rm mor};3,1;\alpha_1,\alpha_3) \\
&=
\bigcup_{\alpha_2 \in \frak A_2}
\widehat{\mathcal  U}({\rm mor};2,1;\alpha_1,\alpha_2)
\times_{R_{\alpha_2}}^{\boxplus\tau}
\widehat{\mathcal  U}({\rm mor};3,2;\alpha_2,\alpha_3).
\endaligned
\end{equation}
Here the summand of the right hand side is defined by
Definition \ref{defn1635}.
Therefore
\begin{equation}\label{1730form2}
\aligned
&\widehat{\mathcal  U}({\rm mor};3,1;\alpha_1,\alpha_3)^{\boxplus\tau} \\
&=
\bigcup_{\alpha_2 \in \frak A_2}
\widehat{\mathcal  U}({\rm mor};2,1;\alpha_1,\alpha_2)^{\boxplus\tau}
\times_{R_{\alpha_2}}
\widehat{\mathcal  U}({\rm mor};3,2;\alpha_2,\alpha_3)^{\boxplus\tau}.
\endaligned
\end{equation}
(Note we take an appropriate smoothing of corners in the right hand
side.)
On its underlying topological space
\begin{equation}\label{form1728new}
\bigcup_{\alpha_2 \in \frak A_2}
{\mathcal  N}_{12}(\alpha_1,\alpha_2)^{\boxplus\tau}
\times_{R_{\alpha_2}}
{\mathcal  N}_{23}(\alpha_2,\alpha_3)^{\boxplus\tau},
\end{equation}
we will define a Kuranishi structure
\begin{equation}\label{form1728}
\aligned
&\widehat{\mathcal  U^+}({\rm mor};3,1;\alpha_1,\alpha_3)\\
=
&\bigcup_{\alpha_2 \in \frak A_2}
\widehat{\mathcal  U^+}({\rm mor};2,1;\alpha_1,\alpha_2)
\times_{R_{\alpha_2}}
\widehat{\mathcal  U^+}({\rm mor};3,2;\alpha_2,\alpha_3)
\endaligned
\end{equation}
as follows.
\begin{rem}
The fiber product appearing in each summand of
(\ref{form1728}) obviously
gives a Kuranishi structure on each
summand of (\ref{form1728new}).
Below we explain how we smooth the corner to
obtain a Kuranishi structure of the union.
\end{rem}
Note that the Kuranishi structure $\widehat{\mathcal  U^+}({\rm mor};2,1;\alpha_1,\alpha_2)$
is $\tau'$-collared.
We take a Kuranishi structure
$\widehat{\mathcal  U^+}({\rm mor};2,1;\alpha_1,\alpha_2)^{\boxminus(\tau
-\tau')}$  such that
$$
\left(\widehat{\mathcal  U^+}({\rm mor};2,1;\alpha_1,\alpha_2)^{\boxminus(\tau
-\tau')}
\right)^{\boxplus\tau'}
=
\widehat{\mathcal  U^+}({\rm mor};2,1;\alpha_1,\alpha_2).
$$
We define
$\widehat{\mathcal  U^+}({\rm mor};3,2;\alpha_2,\alpha_3)^{\boxminus(\tau
-\tau')}$ in the same way.
We consider a Kuranishi structure
\begin{equation}\label{form1731}
\widehat{\mathcal  U^+}({\rm mor};2,1;\alpha_1,\alpha_2)^{\boxminus(\tau
-\tau')}
\times_{R_{\alpha_2}}
\widehat{\mathcal  U^+}({\rm mor};3,2;\alpha_2,\alpha_3)^{\boxminus(\tau
-\tau')}
\end{equation}
on
$$
{\mathcal  N}({\rm mor};2,1;\alpha_1,\alpha_2)^{\boxplus\tau'}
\times_{R_{\alpha_2}}
{\mathcal  N}({\rm mor};3,2;\alpha_2,\alpha_3)^{\boxplus\tau'}.
$$
See Figure \ref{Figure18-1}.
\begin{figure}[h]
\centering
\includegraphics[scale=0.5,angle=-90]{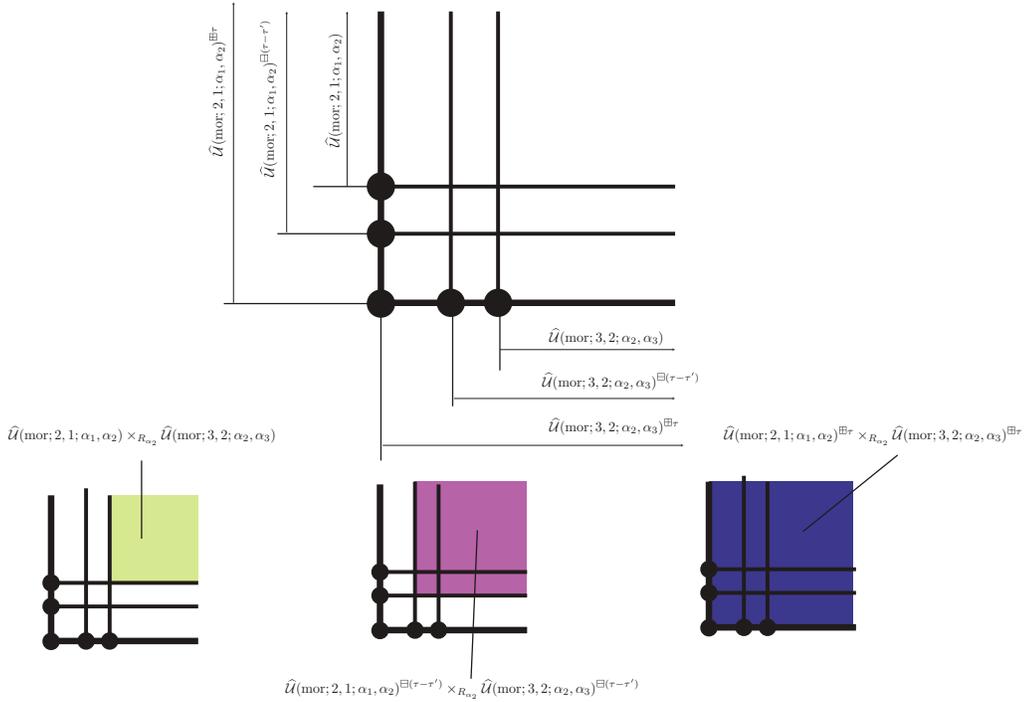}
\caption{$
{\mathcal  N}({\rm mor};2,1;\alpha_1,\alpha_2)^{\boxplus\tau}
\times_{R_{\alpha_2}}
{\mathcal  N}({\rm mor};3,2;\alpha_2,\alpha_3)^{\boxplus\tau}.
$}
\label{Figure18-1}
\end{figure}
\par
On the other hand, let
$\mathcal N_{123}(\alpha_1,\alpha_3)$ be as in
Proposition \ref{prop1636}.
That is, after smoothing corners, its boundary contains
$(
\mathcal N(\alpha_1,\alpha_3),\widehat{\mathcal  U}({\rm mor};3,1;\alpha_1,\alpha_3))$.
We denote
the Kuranishi structure of $\mathcal N_{123}(\alpha_1,\alpha_3)$ by $\widehat{\mathcal U}_{123}(\alpha_1,\alpha_3)$.
\par
We consider the complement $\frak C^c$ of the boundary components
$\frak C$ of  $\mathcal N_{123}(\alpha_1,\alpha_3)$ appearing in Proposition \ref{prop1636}.
Namely $\frak C$ consists of
$$
\bigcup
{\mathcal  N}_{12}(\alpha_1,\alpha_2)^{\boxplus\tau}
\times_{R_{\alpha_2}}
{\mathcal  N}_{23}(\alpha_2,\alpha_3)^{\boxplus\tau}.
$$
The K-space $\mathcal N_{123}(\alpha_1,\alpha_3)$
is $\tau'$-$\frak C$-collared.
We take $\mathcal N_{123}(\alpha_1,\alpha_3)^{\frak C\boxminus(\tau-\tau')}$
and consider the topological space
(see Figure \ref{Figure18-2})
\begin{equation}\label{form173333}
\left(\mathcal N_{123}(\alpha_1,\alpha_3)
\setminus \mathcal N_{123}(\alpha_1,\alpha_3)^{\frak C\boxminus(\tau-\tau')}
\right)^{\frak C^c\boxplus\tau'}.
\end{equation}
Then $\widehat{\mathcal U}_{123}(\alpha_1,\alpha_3)$
induces a Kuranishi structure on it.
We write
$\widehat{\mathcal U}_{123}(\alpha_1,\alpha_3)^{\frak C\boxminus(\tau-\tau')}$
by an abuse of notation.
\par
Starting from (\ref{form1731}) we apply Proposition \ref{prop1636}
to obtain a $\tau'$-$\frak C$-collared Kuranishi structure
on the  topological space
(\ref{form173333}),
which we denote by $\widehat{\mathcal U^+_{123}}(\alpha_1,\alpha_3)$.
The Kuranishi structure
$\widehat{\mathcal U}_{123}(\alpha_1,\alpha_3)^{\frak C\boxminus(\tau-\tau')}$ is
embedded in $\widehat{\mathcal U^+_{123}}(\alpha_1,\alpha_3)$.
\begin{figure}[h]
\centering
\includegraphics[scale=0.3,angle=-90]{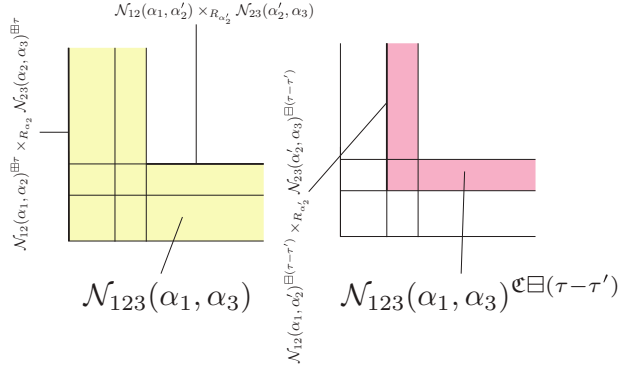}
\caption{$
\mathcal N_{123}(\alpha_1,\alpha_3)^{\frak C\boxminus(\tau-\tau')}
$}
\label{Figure18-2}
\end{figure}
\par
The $\frak C$-partial smoothing of corners
of $\widehat{\mathcal U^+_{123}}(\alpha_1,\alpha_3)$
defines our Kuranishi structure
$\widehat{\mathcal  U^+}({\rm mor};3,1;\alpha_1,\alpha_3)$
of
$\mathcal N_{13}(\alpha_1,\alpha_3)^{\boxplus\tau}$
appeared in (\ref{form1728}),
where
$\mathcal N_{13}(\alpha_1,\alpha_3)$ is the interpolation space of the morphism
$\frak N_{31}$.
We put
\begin{equation}\label{form1729}
\aligned
&\widehat{\frak S^+}({\rm mor};3,1;\alpha_1,\alpha_3)
\\
=&
\bigcup_{\alpha_2 \in \frak A_2}
\widehat{\frak S^+}({\rm mor};2,1;\alpha_1,\alpha_2)
\times_{R_{\alpha_2}}
\widehat{\frak S^+}({\rm mor};3,2;\alpha_2,\alpha_3).
\endaligned
\end{equation}
Here each summand of the right hand side
of (\ref{form1729})
gives a CF-perturbation of each
summand in (\ref{form1728}).
Using the collared-ness,
we can glue them to obtain a CF-perturbation on $\widehat{\mathcal  U^+}({\rm mor};3,1;\alpha_1,\alpha_3)$.
(In other words, they are automatically smooth on the part where they are glued.)
\begin{lem}\label{lem1721}
The Kuranishi structure
$\widehat{\mathcal  U^+}({\rm mor};3,1;\alpha_1,\alpha_3)$
and the CF-perturbation  $\widehat{\frak S^+}({\rm mor};3,1;\alpha_1,\alpha_3)$
satisfy the conclusion of  Proposition \ref{prop1618}.
\end{lem}
\begin{proof}
By construction, the corner contained in
$$
\widehat{\mathcal U^+}({\rm mor};2,1;\alpha_1,\alpha_2)
\times_{R_{\alpha_2}}
\widehat{\mathcal U^+}(2;\alpha_2,\alpha'_2)
\times_{R_{\alpha'_2}}
\widehat{\mathcal U^+}({\rm mor};3,2;\alpha'_2,\alpha_3)
$$
and
$$
\widehat{\mathcal U}({\rm mor};2,1;\alpha_1,\alpha_2)^{\boxplus\tau}
\times_{R_{\alpha_2}}
\widehat{\mathcal U}(2;\alpha_2,\alpha'_2)^{\boxplus\tau}
\times_{R_{\alpha'_2}}
\widehat{\mathcal U}({\rm mor};3,2;\alpha'_2,\alpha_3)^{\boxplus\tau}
$$
are smooth. Thereofore the boundary
of $\widehat{\mathcal  U^+}({\rm mor};3,1;\alpha_1,\alpha_3)$
consists of
$$
\widehat{\mathcal U^+}(1;\alpha_1,\alpha'_1)
\times_{R_{\alpha'_1}}
\widehat{\mathcal U^+}({\rm mor};2,1;\alpha'_1,\alpha_2)
\times_{R_{\alpha_2}}
\widehat{\mathcal U^+}({\rm mor};3,2;\alpha_2,\alpha_3)
$$
and
$$
\widehat{\mathcal U^+}({\rm mor};2,1;\alpha_1,\alpha_2)
\times_{R_{\alpha_2}}
\widehat{\mathcal U^+}({\rm mor};3,2;\alpha_2,\alpha'_3)
\times_{R_{\alpha'_3}}
\widehat{\mathcal U^+}(3;\alpha'_3,\alpha_3).
$$
This is Proposition \ref{prop1618} (4).
Since $\widehat{\mathcal U}_{123}(\alpha_1,\alpha_3)^{\frak C\boxminus(\tau-\tau')}$ is
embedded in $\widehat{\mathcal U^+_{123}}(\alpha_1,\alpha_3)$,
the K-space $\widehat{\mathcal  U}({\rm mor};3,1;\alpha_1,\alpha_3)^{\boxplus\tau}$
is embedded in
$\widehat{\mathcal  U^+}({\rm mor};3,1;\alpha_1,\alpha_3)$.
This is Proposition \ref{prop1618} (1).
The rest of the proof is obvious.
\end{proof}
\par
Let $(i,i')$ be one of $(3,2)$, $(2,1)$, $(3,1)$.
We use the pair
$\widehat{\mathcal  U^+}({\rm mor};i',i;\alpha_i,\alpha_{i'})$
and $\widehat{\frak S^+}({\rm mor};i',i;\alpha_i,\alpha_{i'})$
to apply Definition \ref{chainmapobtainedfromcor}.
We then obtain:
$$
\psi^{i'i;\epsilon}_{\alpha_{i'},\alpha_i} :
\Omega(R_{\alpha_i};o_{R_{\alpha_i}}) \longrightarrow \Omega(R_{\alpha_{i'}};o_{R_{\alpha_{i'}}}).
$$
\begin{lem}\label{17222lem}
\begin{equation}\label{lem1721form}
\psi^{31;\epsilon}_{\alpha_{3},\alpha_1}
=
\sum_{\alpha_2 \in \frak A_2}
\psi^{32;\epsilon}_{\alpha_{3},\alpha_2} \circ
\psi^{21;\epsilon}_{\alpha_{2},\alpha_1}
\end{equation}
in the sense of $(\flat)$. (See Remark \ref{runningoutsharp} for the
meaning of $(\flat)$.)
\end{lem}
\begin{proof}
This is immediate from the composition formula and
Lemma \ref{lem1634lemlem}.
\end{proof}
We will choose $\epsilon_{2,i}$ so that (\ref{lem1721form}) holds
for $\epsilon < \epsilon_{2,i}$.

\begin{cor}
Suppose we are in the Situation of Theorem \ref{linesysmainth2} (3).
Let $E$ be an arbitrary positive number.
Then we can make the choices to define $\frak N_{13}$ so that
$\frak N_{13} \equiv \frak N_{23} \circ \frak N_{12} \mod T^E$ holds.
\end{cor}

\subsection{Construction of homotopy}
\label{subsec:consthomotopy}
In this subsection we start from a homotopy of (partial) morphisms
of linear K-systems and
construct a (partial) cochain homotopy.
We can study higher homotopy in the same way.
Since the definition of parametrized morphism
is a bit heavy, we discuss the case of homotopy first
in this subsection.
The general case of higher homotopy will be
discussed in Subsection \ref{subsec:consthigerhomotopy}.

\begin{shitu}\label{situ16716rev}
Suppose we are in Situation \ref{situparaPmorph} (1).
Suppose also we are given partial linear K-systems
$\mathcal F_i$ of energy cut level
$E_0$ for $i=1,2$ and
a $\tau$-collared Kuranishi structure
$\widehat{\mathcal U^+}(i;\alpha_{i-},\alpha_{i+})$
on $\mathcal M^i(\alpha_{i-},\alpha_{i+})^{\boxplus\tau_0}$
equipped with CF-perturbations
$\widehat{\frak S^+}(i;\alpha_{i-},\alpha_{i+})$ for $i=1,2$ which satisfy
the conclusion of Proposition \ref{prop161}.
Here $0 < \tau < \tau_0=1$ as in Proposition \ref{prop161} and
$\mathcal M^i(\alpha_{i-},\alpha_{i+})$ are as in Condition \ref{morphilinsys}.
(From now on, we write $\alpha_{\pm}$ in place of $\alpha_{i\pm}$ if no confusion can occur.)
\par
Let $0 < \tau' < \tau < \tau_0=1$.
$\blacksquare$
\end{shitu}
\begin{shitu}\label{situ16716revcont}
Suppose we are in Situation \ref{situ16716rev}.
\begin{enumerate}
\item
For $j=1,2$, we are given
partial morphisms $\frak N_j : \mathcal F_1 \to \mathcal F_2$
of energy cut level $E_0$ and energy loss $c$.
We denote by $\mathcal N(j;\alpha_1,\alpha_2)$ its interpolation
space.
\item
Suppose we are given a homotopy $\frak H$ between partial morphisms
$\frak N_{1}$ and $\frak N_{2}$.
By its definition, it is a $[1,2]$-parametrized family of
partial morphisms from $\mathcal F_1$ to $\mathcal F_2$.
Suppose its energy cut level is $E_0$ and energy loss is $c$.
We denote its interpolation
space by $\mathcal N(\alpha_1,\alpha_2;[1,2])$.
$\blacksquare$
\end{enumerate}
\end{shitu}
\begin{shitu}\label{situ16716revcontcont}
Suppose we are in Situation \ref{situ16716revcont}.
Suppose also that, for $j=1,2$, we are given a $\tau'$-collared Kuranishi
structure and a $\tau'$-collared CF-perturbation on $\mathcal N(j;\alpha_1,\alpha_2)$
which satisfy the conclusions of Proposition \ref{prop1618}.
We denote them by
$\widehat{\mathcal U^+}({\rm mor},j;\alpha_1,\alpha_2)$,
$\widehat{\frak S^+}({\rm mor},j;\alpha_1,\alpha_2)$.
$\blacksquare$
\end{shitu}
Before we state the main result, we explicitly write
the boundary and corner compatibility conditions
for the case of $[1,2]$-parametrized family below.
The compatibility condition at the boundary
(Condition \ref{famiboudaru22}) is as follows.
\footnote{See Remark \ref{rem:FiberProdOrd} for the sign and the order of the fiber products.}
\begin{equation}\label{formula1627}
\aligned
&\partial \mathcal N(\alpha_1,\alpha_2;[1,2])
\\=&
\coprod_{\alpha'_1 \in \frak A_1}
(-1)^{\dim \mathcal N(\alpha'_1,\alpha_2;[1,2])}
\mathcal N(\alpha'_1,\alpha_2;[1,2])
\times_{R_{\alpha'_1}}
\mathcal M^1(\alpha_1,\alpha'_1)
\\
&\sqcup
\coprod_{\alpha'_2 \in \frak A_2}
(-1)^{\dim \mathcal M^2(\alpha'_2,\alpha_2)}
\mathcal M^2(\alpha'_2,\alpha_2)
\times_{R_{\alpha'_2}}
\mathcal N(\alpha_1,\alpha'_2;[1,2])
\\
&\sqcup \mathcal N(2;\alpha_1,\alpha_2)
\sqcup -\mathcal N(1;\alpha_1,\alpha_2).
\endaligned
\end{equation}
The first of the corner compatibility conditions
(Condition \ref{boundarycompPban1}) says that the normalized corner
$\widehat S_k(\mathcal N(\alpha_1,\alpha_2;[1,2]))$
is the disjoint union of the following two types of
fiber products:
\begin{equation}\label{formula143333revrev1}
\aligned
&{\mathcal M^{1}}(\alpha_-,\alpha_1)
\times_{R_{\alpha_1}}
\dots
\times_{R_{\alpha_{k_1}-1}}
{\mathcal M^{1}}(\alpha_{{k_1}-1},\alpha_{{k_1}}) \\
&
\times_{R_{\alpha_{k_1}}}
{\mathcal N}({\alpha_{k_1}},{\alpha_{k_1+1}};[1,2]) \\
&
\times_{R_{\alpha_{k_1+1}}}
{\mathcal M^{2}}(\alpha_{k_1+1},\alpha_{k_1+2})
\times_{R_{\alpha_{k_1+2}}}
\dots
\times_{R_{\alpha_{k_1+k_2}}}
{\mathcal M^{2}}(\alpha_{{k_1}+k_2},\alpha_+),
\endaligned
\end{equation}
with $k_1+k_2 = k$, and
\begin{equation}\label{formula143333revrev2}
\aligned
&{\mathcal M^{1}}(\alpha_-,\alpha_1)
\times_{R_{\alpha_1}}
\dots
\times_{R_{\alpha_{k_1}-1}}
{\mathcal M^{1}}(\alpha_{{k_1}-1},\alpha_{{k_1}}) \\
&
\times_{R_{\alpha_{k_1}}}
{\mathcal N}(j;{\alpha_{k_1}},{\alpha_{k_1+1}}) \\
&
\times_{R_{\alpha_{k_1+1}}}
{\mathcal M^{2}}(\alpha_{k_1+1},\alpha_{k_1+2})
\times_{R_{\alpha_{k_1+2}}}
\dots
\times_{R_{\alpha_{k_1+k_2}}}
{\mathcal M^{2}}(\alpha_{{k_1}+k_2},\alpha_+),
\endaligned
\end{equation}
with $k_1+k_2 = k-1$, $j=1,2$.
\par
The second of the corner compatibility conditions
(Condition \ref{furthercompatifiber}) says the following:
Consider the $\ell$-th normalized corner of the
K-spaces (\ref{formula143333revrev1}).
According to the descriptions of
$\widehat S_n(\mathcal N(\alpha_1,\alpha_2;[1,2]))$,
$\mathcal N(j;\alpha_1,\alpha_2)$
and $\mathcal M^i(\alpha_1,\alpha_2)$ we gave above,
we see that
the $\ell$-th normalized corner of (\ref{formula143333revrev1})
or (\ref{formula143333revrev2})
is given by the disjoint union of the same type of fiber products as
 (\ref{formula143333revrev1})
or (\ref{formula143333revrev2}).
Condition \ref{furthercompatifiber} requires that this description coincides with
the description of
$\widehat S_{k+\ell}(\mathcal N(\alpha_1,\alpha_2;[1,2]))$.
\begin{defn}\label{defn1726}
We consider the evaluation map
$$
{\rm ev}_{[1,2]} : \mathcal N(\alpha_1,\alpha_2;[1,2])
\longrightarrow [1,2].
$$
The inverse image ${\rm ev}_{[1,2]}^{-1}(\partial[1,2])$
is a part of the boundary
$\partial\mathcal N(\alpha_1,\alpha_2;[1,2])$.
We denote it by $\partial_{\frak C^v}$ and call it the {\it vertical
boundary}.\index{vertical boundary}\index{boundary ! vertical boundary}
We put
$$
\partial_{\frak C^h} \mathcal N(\alpha_1,\alpha_2;[1,2]) :=
\partial \mathcal N(\alpha_1,\alpha_2;[1,2]) \setminus
\partial_{\frak C^v}\mathcal N(\alpha_1,\alpha_2;[1,2])
$$
and call it the {\it horizontal boundary}.
\index{horizontal boundary}\index{boundary ! horizontal boundary}
\end{defn}
Note that we have the map
$$
{\rm ev}_{[1,2]} : \mathcal N(\alpha_1,\alpha_2;[1,2])^{\frak C^h\boxplus\tau}
\longrightarrow [1,2].
$$
\begin{prop}\label{homotopyexistpertstatement}
Suppose we are in Situations
\ref{situ16716rev}, \ref{situ16716revcont}, \ref{situ16716revcontcont}
and $\tau'' < \tau'$.
Then for any $\alpha_1,\alpha_2$ with $E(\alpha_2) - E(\alpha_1) \le E_0 -c$,
there exist $\widehat{\mathcal U^+}(\alpha_1,\alpha_2;[1,2])$ and
$\widehat{\frak S^+}(\alpha_1,\alpha_2;[1,2])$ such that
they enjoy the following properties.
\begin{enumerate}
\item
$\widehat{\mathcal U^+}(\alpha_1,\alpha_2;[1,2])$ is a
$\tau''$-${\frak C^h}$-collared Kuranishi structure of $\mathcal N(\alpha_1,\alpha_2;[1,2])^{\frak C^h\boxplus\tau}$
and $\widehat{\frak S^+}(\alpha_1,\alpha_2;[1,2])$ is its
$\tau''$-${\frak C^h}$-collared CF-perturbation.
\item
$\widehat{\frak S^+}(\alpha_1,\alpha_2;[1,2])$ is transversal to $0$.
Moreover the map
$$
({\rm ev}_+, {\rm ev}_{[1,2]}) : \mathcal N(\alpha_1,\alpha_2;[1,2])^{\frak C^h\boxplus\tau}
\longrightarrow R_{\alpha_2} \times [1,2]
$$
is strongly stratumwise submersive with respect to $\widehat{\frak S^+}(\alpha_1,\alpha_2;[1,2])$.
\item
We have periodicity isomorphisms among $\widehat{\mathcal U^+}(\alpha_1,\alpha_2;[1,2])$'s
that are compatible with $\widehat{\frak S^+}(\alpha_1,\alpha_2;[1,2])$.
\item
There exists an embedding of $\tau''$-collared Kuranishi structures
from
$$
\mathcal N(\alpha_1,\alpha_2;[1,2])^{{\frak C^h}\boxplus\tau}
$$
to
the  $\tau''$- ${\frak C^h}$-collared Kuranishi structure
$
\widehat{\mathcal U^+}(\alpha_1,\alpha_2;[1,2]).
$
Here for the source $\mathcal N(\alpha_1,\alpha_2;[1,2])^{{\frak C^h}\boxplus\tau}$ we use the
$\tau''$-${\frak C^h}$-collared Kuranishi structure induced by that of
$\mathcal N(\alpha_1,\alpha_2;[1,2])$
which is given by the definition of $[1,2]$-parametrized interpolation space.
\item
There is an isomorphism of $\tau''$-${\frak C^h}$-collared K-spaces
\footnote{See Remark \ref{rem:FiberProdOrd} for the sign and the order of the fiber products.}
\begin{equation}\label{formula1630}
\aligned
&\partial \widehat{\mathcal U^+}(\alpha_1,\alpha_2;[1,2])
\\=&
\coprod_{\alpha'_1\in \frak A_1}
(-1)^{\dim \widehat{\mathcal U^+}(\alpha'_1,\alpha_2;[1,2])}
\widehat{\mathcal U^+}(\alpha'_1,\alpha_2;[1,2])
\times_{R_{\alpha'_1}}
\widehat{\mathcal U^+}(1;\alpha_1,\alpha'_1)
\\
&\sqcup
\coprod_{\alpha'_2 \in \frak A_2}
(-1)^{\dim \widehat{\mathcal U^+}(2;\alpha'_2,\alpha_2)}
\widehat{\mathcal U^+}(2;\alpha'_2,\alpha_2)
\times_{R_{\alpha'_2}}
\widehat{\mathcal U^+}(\alpha_1,\alpha'_2;[1,2])
\\
&\sqcup \widehat{\mathcal U^+}({\rm mor},2;\alpha_1,\alpha_2)
\sqcup -\widehat{\mathcal U^+}({\rm mor};1;\alpha_1,\alpha_2).
\endaligned
\end{equation}
The isomorphism (\ref{formula1630}) is compatible with the isomorphism
(\ref{formula1627}) via the embedding (4).
It is also compatible with the periodicity isomorphism and the evaluation maps.
\item
The pull-back of $ \widehat{\mathcal U^+}(\alpha_1,\alpha_2;[1,2])$ by the isomorphism
(\ref{formula1630})
is equivalent to the fiber product of $\widehat{\frak S^+}(j;\alpha,\alpha')$,
$\widehat{\frak S^+}({\rm mor},j;\alpha,\alpha')$
and
$ \widehat{\frak S^+}(\alpha,\alpha';[1,2])$.
\item
The normalized corner
$\widehat S_k(\widehat{\mathcal U^+}(\alpha_1,\alpha_2;[1,2]))$
is a disjoint union of the following two types of fiber products.
\begin{equation}\label{formula143333revrev1prop}
\aligned
&\widehat{\mathcal U^+}(1;\alpha_-,\alpha_1)
\times_{R_{\alpha_1}}
\dots
\times_{R_{\alpha_{k_1}-1}}
\widehat{\mathcal U^+}(1;\alpha_{{k_1}-1},\alpha_{{k_1}}) \\
&
\times_{R_{\alpha_{k_1}}}
\widehat{\mathcal U^+}({\alpha_{k_1}},{\alpha_{k_1+1}};[1,2]) \\
&
\times_{R_{\alpha_{k_1+1}}}
\widehat{\mathcal U^+}(2;\alpha_{k_1+1},\alpha_{k_1+2})
\times_{R_{\alpha_{k_1+2}}}
\dots
\times_{R_{\alpha_{k_1+k_2}}}
\widehat{\mathcal U^+}(2;\alpha_{{k_1}+k_2},\alpha_+)
\endaligned
\end{equation}
and
\begin{equation}\label{formula143333revrev2prop}
\aligned
&\widehat{\mathcal U^+}(1;\alpha_-,\alpha_1)
\times_{R_{\alpha_1}}
\dots
\times_{R_{\alpha_{k_1}-1}}
\widehat{\mathcal U^+}(1;\alpha_{{k_1}-1},\alpha_{{k_1}}) \\
&
\times_{R_{\alpha_{k_1}}}
\widehat{\mathcal U^+}({\rm mor};j;{\alpha_{k_1}},{\alpha_{k_1+1}}) \\
&
\times_{R_{\alpha_{k_1+1}}}
\widehat{\mathcal U^+}(2;\alpha_{k_1+1},\alpha_{k_1+2})
\times_{R_{\alpha_{k_1+2}}}
\dots
\times_{R_{\alpha_{k_1+k_2}}}
\widehat{\mathcal U^+}(2;\alpha_{{k_1}+k_2},\alpha_+).
\endaligned
\end{equation}
This isomorphism is compatible with the isomorphism (\ref{formula143333revrev1}), (\ref{formula143333revrev2}) via the embedding (4).
It is also compatible with the periodicity isomorphism and the evaluation maps.
\item
The pull-back of the restriction of $\widehat{\frak S^+}(\alpha_1,\alpha_2;[1,2])$
to $\widehat S_k(\widehat{\mathcal U^+}(\alpha_1,\alpha_2;[1,2]))$
by the isomorphism in (7)
is equivalent to the fiber product of $\widehat{\frak S^+}({\rm mor};j;*,*)$,  $\widehat{\frak S^+}(1;*,*)$,
$\widehat{\frak S^+}(2;*,*)$.
\item
The isomorphism of (7) is compatible with the covering map
$$
\widehat S_{\ell}(\widehat S_k(\widehat{\mathcal U^+}(\alpha_1,\alpha_2;[1,2])))
\longrightarrow \widehat S_{k+\ell}(\widehat{\mathcal U^+}(\alpha_1,\alpha_2;[1,2])).
$$
(The precise meaning of this compatibility is the same as the case of
$$
\mathcal N(\alpha_1,\alpha_2;[1,2]),
$$ which we explained right before this proposition.)
\item
If we start from a uniform family of $\widehat{\frak S^+}(i;\alpha_-,\alpha_+)$ and
$\widehat{\frak S^+}({\rm mor},j;\alpha_1,\alpha_2)$,
then we can take the family of $\widehat{\frak S^+}(\alpha_1,\alpha_2;[1,2])$
to be uniform.
\end{enumerate}
\end{prop}
\begin{proof}
We prove the Proposition \ref{homotopyexistpertstatement}
with $E_0$ replaced by $E_{\frak E}^n$,
by induction on $n$.
The corner compatibility conditions
Conditions \ref{boundarycompPban1} and
Conditions \ref{furthercompatifiber} are written in such a way
that they immediately imply the assumptions of Proposition \ref{prop528},
which is Situation \ref{sit1526} (especially its (1)(2)).
\end{proof}
We rewrite the geometric conclusion of Proposition \ref{homotopyexistpertstatement}
into algebraic language.
\begin{defn}
In the situation of Proposition \ref{homotopyexistpertstatement}, we define
\begin{equation}
\frak h^{\epsilon}_{\alpha_2,\alpha_1} :
\Omega(R_{\alpha_1};o_{R_{\alpha_1}}) \longrightarrow \Omega(R_{\alpha_2};o_{R_{\alpha_2}})
\end{equation}
by
\begin{equation}\label{form16ten6rev}
\frak h^{\epsilon}_{\alpha_2,\alpha_1}(h)
=
{\rm ev}_{+}!({\rm ev}_{-}^* h;\widehat{\frak S^{+ \epsilon}}(\alpha_1,\alpha_2;[1,2])).
\end{equation}
Here the right hand side is defined by
Definition \ref{pushoutdeftau} on
the K-space
$$
({\mathcal N}({\alpha_1},{\alpha_2};[1,2])^{\frak C^h\boxplus\tau_0},
\widehat{\mathcal U^+}(\alpha_1,\alpha_2;[1,2])).
$$
\end{defn}
By Condition \ref{Pparamorphi} (VI) and \cite[Definition 7.78]{part11},
the degree of $\frak h^{\epsilon}_{\alpha_2,\alpha_1}$
is $\eta(\alpha_2) -\eta(\alpha_1) - 1$.
Therefore after degree shift  as in Definition \ref{Fvectspace} (2)
its degree becomes $-1$.
\par
If the energy loss of our homotopy is
$c$, the family $\{\frak h^{\epsilon}_{\alpha_2,\alpha_1}\}$ of maps induces
$$
\frak F^{\lambda}CF(\mathcal F_1) \to \frak F^{\lambda-c}CF(\mathcal F_2)
$$
where the filtration $\frak F^{\lambda}$ is defined in
Definition \ref{Fvectspace} (2)(3).
\par
Suppose we are in the situation of Proposition \ref{homotopyexistpertstatement}.
We have two morphisms of partial linear K-systems
equipped with CF-perturbations. Namely we have
$\mathcal N(j;\alpha_1,\alpha_2)$,
$\widehat{\mathcal U^+}({\rm mor},j;\alpha_1,\alpha_2)$,
$\widehat{\frak S^+}({\rm mor},j;\alpha_1,\alpha_2)$
for $j=1,2$.
We use Definition \ref{chainmapobtainedfromcor}
for $j=1,2$ to obtain maps
$$
\psi^{j,\epsilon}_{\alpha_2,\alpha_1} : \Omega(R_{\alpha_1};o_{R_{\alpha_1}}) \to \Omega(R_{\alpha_2};o_{R_{\alpha_2}}).
$$
\begin{lem}\label{lem165rev1}
The linear maps $\{\frak h^{\epsilon}_{\alpha_1,\alpha_2}\}$
satisfy the following equality in the sense of $(\flat)$ in
Remark \ref{runningoutsharp}:
\begin{equation}\label{form16353}
\aligned
&d_0\circ \frak h^{\epsilon}_{\alpha_2,\alpha_1}
+ \frak h_{\alpha_2,\alpha_1} \circ d_0
\\
=&
-\sum_{\alpha'_1}
\frak h^{\epsilon}_{\alpha_2,\alpha'_1}
\circ
\frak m^{1,\epsilon}_{1;\alpha'_1,\alpha_1} -
\sum_{\alpha'_2}
\frak m^{2,\epsilon}_{1;\alpha_2,\alpha'_2}
\circ
\frak h^{2,\epsilon}_{\alpha'_2,\alpha_1}
+\psi^{2,\epsilon}_{\alpha_2,\alpha_1} - \psi^{1,\epsilon}_{\alpha_2,\alpha_1}.
\endaligned
\end{equation}
Here the first sum in the second line is taken over
$\alpha'_1 \in \frak A_1$ with $E(\alpha_1) < E(\alpha'_1) \le E(\alpha_2) + c$
and the second sum in the second line is taken over
$\alpha'_2 \in \frak A_2$ with $ E(\alpha_1)- c \le  E(\alpha'_2) < E(\alpha_2)$.
The number $c$ is the energy loss of our morphism.
\end{lem}
We will define $\epsilon_{3,i}$ so that
(\ref{form16353}) holds for $\epsilon < \epsilon_{3,i}$.
\begin{proof}
The proof is similar to the proof of Lemmas \ref{lem165}, \ref{lem165rev}.
By Stokes' formula the left hand side is obtained from
$\partial\widehat{\frak S^{+ \epsilon}}(\alpha_1,\alpha_2;[1,2])$
in the same way as (\ref{form16ten6rev}).
We can decompose the boundary
$\partial\widehat{\frak S^{+ \epsilon}}(\alpha_1,\alpha_2;[1,2])$
into a disjoint union by Proposition \ref{homotopyexistpertstatement}.
Then we use the composition formula to obtain the right hand side of
(\ref{form16353}).
Namely the 1,2,3,4-th union of (\ref{formula1630})
correspond to the 1,2,3,4-th term of (\ref{form16353}), respectively.
\end{proof}

\subsection{Proof of Theorem \ref{linesysmainth1} (2)(except (f)),
Theorem \ref{linesysmainth2} (1) and Theorem \ref{indlinesysmainth1} (2)(except (e)), (3)}
\label{subsec:proofsec14main2}

In this subsection we use the result of Subsection
\ref{subsec:consthomotopy}
to prove  Theorem \ref{linesysmainth1} (2)
(a)-(e).
We also prove Theorem \ref{linesysmainth2} (1) and
Theorem \ref{indlinesysmainth1} (2)(3)
(except (2)(e)) at the same time.
\par
To prove Theorem \ref{linesysmainth2} (1) we need to
take the `projective limit' in a similar way as the argument of
Subsection \ref{subsec:proofsec14main1}.
Proposition \ref{prop1631} below is its algebraic part.
It is similar to Lemma \ref{lem1623} and is a baby version of
\cite[Lemma 7.2.129]{fooobook2}.
\begin{defn}\label{defn1628}
Suppose we are given two partial linear K-systems
$\mathcal F_i$ ($i=1,2$).
Also suppose we are given
$(CF(\mathcal F_i),\{\frak m^i_{1;\alpha_+,\alpha_-}\})$, a partial cochain complex structure on $CF(\mathcal F_i)$ of energy cut level $E_{(i)}$
for $i=1,2$.
We assume $0 \le c < E_0 \le E_{(1)}$ and $E_0 \le E_{(2)} - c$.
\begin{enumerate}
\item
A {\it partial cochain map $\mathcal F_1 \to \mathcal F_2$\index{partial ! cochain map} of
energy cut level $E_0$ and energy loss $c$}\index{energy cut level!partial cochain map}
\index{energy loss ! of partial cochain map} is
a family $\widehat\psi = \{\psi_{\alpha_2,\alpha_1}\}$ consisting of
the following objects $\psi_{\alpha_2,\alpha_1}$:
\par
If $\alpha_i \in \frak A_i$ and
$E(\alpha_2) \le E(\alpha_1) +E_0 -c$,
we have an $\R$ linear map
$$
\psi_{\alpha_2,\alpha_1}
~:~
\Omega(R_{\alpha_1};o_{R_{\alpha_1}}) \longrightarrow \Omega(R_{\alpha_2};o_{R_{\alpha_2}}).
$$
We require that it satisfies (\ref{form1613})
for $E(\alpha_2) \le E(\alpha_1) +E_0 -c$.
\item
For $j=1,2$, let $\{\psi_{j;\alpha_2,\alpha_1}\}$ be partial cochain maps
of energy cut level $E_0$ with energy loss $c$.
A {\it partial cochain homotopy of energy cut level $E_0$ with energy loss $c$
from $\widehat\psi_1 = \{\psi_{1;\alpha_2,\alpha_1}\}$\index{cochain homotopy ! partial cochain homotopy}\index{partial ! cochain homotopy}\index{energy cut level!partial cochain homotopy}\index{energy loss ! of partial cochain homotopy} to
$\widehat\psi_2 = \{\psi_{2;\alpha_2,\alpha_1}\}$}
is a family $\widehat{\frak h} = \{\frak h_{\alpha_2,\alpha_1}\}$ consisting of the following objects
$\frak h_{\alpha_2,\alpha_1}$:
\par
If $\alpha_i \in \frak A_i$ and
$E(\alpha_2) \le E(\alpha_1) +E_0 -c$,
we have an $\R$ linear map
$$
\frak h_{\alpha_2,\alpha_1}
~:~
\Omega(R_{\alpha_1};o_{R_{\alpha_1}}) \longrightarrow
\Omega(R_{\alpha_2};o_{R_{\alpha_2}}).
$$
We require that it satisfies (\ref{form16353})
for $E(\alpha_2) \le E(\alpha_1) +E_0 -c$.
\item
For $i=1,2$, let $\widehat\psi_{i+1 i} = \{\psi_{\alpha_{i+1},\alpha_i}\}$
be partial cochain maps of energy loss $c_i$.
Its energy cut level is $E_0$ for $i=1$ and  $E_0-c_1$ for $i=2$.
We define the {\it composition}
$
\widehat\psi_{3 1} = \widehat\psi_{3 2} \circ \widehat\psi_{2 1}
$
by
$$
(\psi_{31})_{\alpha_{3}\alpha_1}
=
\sum_{\alpha_2 \in \frak A_2}
(\psi_{32})_{\alpha_{3}\alpha_2} \circ
(\psi_{21})_{\alpha_{2}\alpha_1}.
$$
Then $\widehat\psi_{3 1}$ is a partial cochain map of energy cut level $E_0$
and energy loss $c_1 + c_2$.
\item
Consider $E'_0$ such that $c < E'_0 < E_0$.
In the situation of (1) we forget all of
$\psi_{\alpha_2,\alpha_1}$ for
$E(\alpha_2) \le E(\alpha_1) +E'_0 -c$.
We then obtain a
partial cochain map $\mathcal F_1 \to \mathcal F_2$ of
energy cut level $E'_0$.
We call it the {\it energy cut of $\widehat\psi$ at energy
cut level $E'_0$.}\index{energy cut level!energy cut}
\par
The {\it energy cut of $\widehat{\frak h}$ at energy
cut level $E'_0$} is defined in the same way.
\item
Let $\widehat\psi : \mathcal F_1 \to \mathcal F_2$ be a partial cochain map of
energy cut level $E_0$ and energy loss $c$ and let $c < E'_0 < E_0$.
If $\widehat\psi'$ is an energy cut of $\widehat\psi$ at energy
cut level $E'_0$,
we call $\widehat\psi$ a {\it promotion of $\widehat\psi'$
to the energy cut level $E_0$.}
\index{energy cut level!promotion}\index{promotion}
A promotion of a partial cochain homotopy is defined in the same way.
\end{enumerate}
\end{defn}
\begin{lemdef}
Two partial cochain maps are said to be {\it cochain homotopic}\index{cochain homotopy ! cochain homotopic}
if there exists a partial cochain homotopy between them.
This is an equivalence relation.
\end{lemdef}
\begin{proof}
If $\widehat{\frak h}^j = \{\frak h^j_{\alpha_2\alpha_1}\}$
is a partial cochain homotopy from $\widehat\psi_j$ to $\widehat\psi_{j+1}$
for $j=1,2$, then
$\{\frak h^1_{\alpha_2\alpha_1} + \frak h^2_{\alpha_2\alpha_1}\}$
is a partial cochain homotopy from $\widehat\psi_1$ to $\widehat\psi_{3}$.
The other part of the proof is obvious.
\end{proof}
\begin{prop}\label{prop1631}
Let $(CF(\mathcal F_i),\{\frak m^i_{1;\alpha_+\alpha_-}\})$
and $(CF(\mathcal F'_i),\{\frak m^{i\prime}_{1;\alpha_+\alpha_-}\})$ be
partial cochain complexes of energy cut level $E_2$
for $i=1,2$.
We take $E_1 < E_2$.
Suppose we have a diagram
\begin{equation}\label{diag1637}
\begin{CD}
(CF(\mathcal F'_1),\{\frak m^{1\prime}_{1;\alpha_+,\alpha_-}\})
@ > {\widehat\psi'_{21}} >>(CF(\mathcal F'_2),\{\frak m^{2\prime}_{1;\alpha_+,\alpha_-}\})
\\
@ A{\widehat\psi_{1}}AA @ AA{\widehat\psi_{2}}A \\
(CF(\mathcal F_1),\{\frak m^1_{1;\alpha_+,\alpha_-}\})
@ > {\widehat\psi_{21}} >>(CF(\mathcal F_2),\{\frak m^2_{1;\alpha_+,\alpha_-}\})
\end{CD}
\end{equation}
such that:
\begin{enumerate}
\item[(i)]
$\widehat\psi_{21}$, $\widehat\psi'_{21}$
are partial cochain maps of energy cut level $E_2$ and energy loss $0$.
We assume that $\widehat\psi_{21}$, $\widehat\psi'_{21}$ induces isomorphisms
modulo $T^{\epsilon}$ for a sufficiently small $\epsilon >0$.
\item[(ii)]
$\widehat\psi_{2}$ is a partial cochain map of energy cut level $E_2$ and energy loss $c$.
\item[(iii)]
$\widehat\psi_{1}$ is a partial cochain map of energy cut level $E_1$ and energy loss $c$.
\item[(iv)]
The Diagram \ref{diag1637} is homotopy commutative as partial cochain maps of energy cut level $E_1$ and energy loss $c$.
Namely there exists a partial cochain homotopy
$\widehat{\frak h} = \{\frak h_{\alpha_2,\alpha_1}\}$
where $\frak h_{\alpha_2,\alpha_1}$ is defined when
$E(\alpha_2) \le E(\alpha_1) +E_2 -c$,
and satisfies
\begin{equation}\label{form16353rev}
\aligned
&d_0\circ \frak h_{\alpha'_2,\alpha_1}
+ \frak h_{\alpha'_2,\alpha_1} \circ d_0
\\
&=
-\sum_{\hat\alpha_1 \in \frak A_1}
\frak h_{\alpha'_2,\hat\alpha_1}
\circ
\frak m^{1}_{1;\hat\alpha_1,\alpha_1}
-\sum_{\hat\alpha'_2 \in \frak A'_2}
\frak m^{2}_{1;\alpha'_2,\hat\alpha'_2}
\circ
\frak h^{2}_{\hat\alpha'_2,\alpha_1}
\\
&\quad +\sum_{\hat\alpha_2 \in \frak A_2}
(\hat\psi_{2})_{\alpha'_2\hat\alpha_2} \circ
(\hat\psi_{21})_{\hat\alpha_2\alpha_1} -
\sum_{\hat\alpha'_1 \in \frak A'_1}
(\hat\psi'_{21})_{\alpha'_2\hat\alpha'_1} \circ (\widehat\psi_{1})_{\hat\alpha'_1,\alpha_1},
\endaligned
\end{equation}
if $E(\alpha'_2) \le E(\alpha_1) +E_1 -c$.
\end{enumerate}
Then we can promote $\widehat\psi_{1}$ to a partial cochain map of energy cut level $E_2$
and energy loss $c$
and $\widehat{\frak h}$ to a partial cochain homotopy of energy cut level $E_2$
and energy loss $c$.
That is, (\ref{form16353rev})
holds if $E(\alpha'_2) \le E(\alpha_1) +E_2 -c$.
\end{prop}
\begin{proof}
By using induction and the discreteness
of the set of energies (which follows from
uniform Gromov compactness
Definition \ref{defn1528} (2) (g)), it suffices to prove
the statement for the case when
\begin{equation}\label{gapconsthomotopro}
\aligned
E(\alpha'_2) - E(\alpha_1) \notin (E_1-c,E_2-c)
\qquad
&\text{for all $(\alpha_1,\alpha'_2) \in \frak A_1 \times \frak A'_2$.}\\
E(\alpha'_i) - E(\alpha_i) \notin (E_1-c,E_2-c)
\qquad
&\text{for all $(\alpha_i,\alpha'_i) \in \frak A_i \times \frak A'_i$,
$i =1,2$.}\\
E(\alpha'_1) - E(\alpha_1) \notin (E_1,E_2)
\qquad
&\text{for all $(\alpha_1,\alpha'_1) \in \frak A_1 \times \frak A'_1$.}
\\
E(\alpha'_2) - E(\alpha_2) \notin (E_1-c,E_2-c)
\qquad
&\text{for all $(\alpha_2,\alpha'_2) \in \frak A_2 \times \frak A'_2$.}
\endaligned
\end{equation}
We will prove the proposition for this case below.
\par
We first promote $\widehat\psi_{1}$ to the energy cut level $E_2$.
We regard
{\it both} $\widehat\psi_{j}$ $(j=1,2)$ as partial cochain maps
of energy cut level $E_1$ and energy loss $c$.
We use them to define
$
\widehat\psi_{j} : CF(\mathcal F_j) \to CF(\mathcal F'_j)
$
as
$$
\widehat\psi_{j} = \bigoplus_{\alpha_j,\alpha'_j}
(\psi_{j})_{\alpha'_j\alpha_j}
$$
where we consider $\alpha_j$, $\alpha'_j$
with $E(\alpha'_{j}) - E(\alpha_j) \le E_1$.
Note that $(\psi_{2})_{\alpha'_2\alpha_2}$ is
defined for $\alpha'_2, \alpha_2$ satisfying $E(\alpha'_{2}) - E(\alpha_2) \le E_2$.
To clarify the energy cut we did for $\widehat\psi_{2}$
we write $\widehat\psi_{2}\vert_{E_1}$.
We define $\widehat{\frak h}
: CF(\mathcal F_1) \to CF(\mathcal F'_2)$
in the same way.
We also drop the energy cut level $E_2$ of $\widehat\psi'_{21}$, $\widehat\psi_{21}$
to $E_1$,
which we write $\widehat\psi'_{21}\vert_{E_1}$,
$\widehat\psi_{21}\vert_{E_1}$ respectively.
We regard them as
$$
\widehat\psi'_{21}\vert_{E_1} : CF(\mathcal F'_1) \to CF(\mathcal F'_2),
\quad
\widehat\psi_{21}\vert_{E_1} : CF(\mathcal F_1) \to CF(\mathcal F_2)
$$
in that sense.
We define
$$
\widehat d_{j}\vert_{E_1} : CF(\mathcal F_j) \to CF(\mathcal F_j),
\qquad
\widehat d'_{j}\vert_{E_1} : CF(\mathcal F'_j) \to CF(\mathcal F'_j)
$$
in the same way.
\par
Let $E(\alpha'_1) - E(\alpha_1) = E_2 - c$.
We define
$o(\alpha_1,\alpha'_1) \in {\rm Hom}_{\R}(o_{R_{\alpha_1}},o_{R_{\alpha'_1}})$ by
\begin{equation}\label{defoooaaa}
o(\alpha_1,\alpha'_1) = (\widehat d'_{1} \circ \widehat\psi_{1} - \widehat\psi_{1} \circ \widehat d_{1})_{\alpha'_1\alpha_1}.
\end{equation}
As in the proof of Lemma \ref{lem1623},
we can show
$$
d_0 \circ o(\alpha_1,\alpha'_1) + o(\alpha_1,\alpha'_1) \circ d_0 = 0.
$$
\par
We will prove that $o(\alpha_1,\alpha'_1)$ is a $d_0$-coboundary.
We consider
$$
(\widehat\psi'_{21}\vert_{E_1} \circ (\widehat d'_{1} \circ \widehat\psi_{1}
- \widehat\psi_{1} \circ \widehat d_{1}))_{\alpha'_1\alpha_1}.
$$
Using the fact that (\ref{defoooaaa}) is zero for $E(\alpha'_1) - E(\alpha_1) < E_0 - c$
and $\widehat\psi_{**}$ has energy loss 0 (Definition \ref{defn1628} (6)), we find that
$$
(\widehat\psi'_{21} \circ (\widehat d'_{1} \circ \widehat\psi_{1} - \widehat\psi_{1} \circ \widehat d_{1}) )_{\alpha_1\alpha'_1}
= o(\alpha_1,\alpha'_1).
$$
On the other hand, we have
\begin{equation}\label{form1642}
\aligned
&(\widehat\psi'_{21}\vert_{E_1} \circ (\widehat d'_{1}\vert_{E_1} \circ \widehat\psi_{1} - \widehat\psi_{1} \circ \widehat d_{1}\vert_{E_1}) )_{\alpha_1\alpha'_1}
\\
& \equiv
\left(\widehat d'_{2}\vert_{E_1} \circ
\widehat\psi'_{21}\vert_{E_1} \circ \widehat\psi_{1}
-
\widehat\psi'_{21}\vert_{E_1} \circ \widehat\psi_{1}
 \circ \widehat d_{1}\vert_{E_1} \right)_{\alpha'_1\alpha_1}.
\endaligned
\end{equation}
Here $\equiv$ means modulo $d_0$-coboundary. Such $d_0$-coboundary appears because
$$
(\widehat d'_{2}\vert_{E_1} \circ \widehat\psi'_{21}\vert_{E_1}
-
\widehat\psi'_{21}\vert_{E_1} \circ \widehat d'_{2}\vert_{E_1})_{\alpha'_1\alpha_1}
= d_0\circ (\widehat\psi'_{21})_{\alpha'_1\alpha_1} -
(\widehat\psi'_{21})_{\alpha'_1\alpha_1}
\circ d_0.
$$
Then we use
$$
(\widehat\psi'_{21}\vert_{E_1} \circ \widehat\psi_{1}
-
\widehat\psi_{2}\vert_{E_1} \circ  \widehat\psi_{21}\vert_{E_1})_{\alpha'_1\alpha_1}
=
d_0 \circ \frak h_{\alpha_1,\alpha'_1} +
\frak h_{\alpha'_1\alpha_1} \circ d_0
$$
to show that
\begin{equation}\label{form1643}
{\rm (\ref{form1642})}
\equiv
(\widehat d'_{2}\vert_{E_1} \circ
\widehat\psi_{2}\vert_{E_1} \circ \widehat\psi_{21}\vert_{E_1}
-
\widehat\psi_{2}\vert_{E_1} \circ  \widehat\psi_{21}\vert_{E_1}
 \circ \widehat d_{2}\vert_{E_1} )_{\alpha'_1\alpha_1}.
\end{equation}
The right hand side of (\ref{form1643}) is a $d_0$-coboundary since the
partial cochain maps
$\widehat\psi_{21}\vert_{E_1}$,
$\widehat \psi_{2}\vert_{E_1}$ appearing here can be promoted to the energy cut level $E_2$ by assumption.
Thus $o(\alpha_1,\alpha'_1)$ is a $d_0$-coboundary. Therefore we can
define $(\widehat\psi_{1})_{\alpha'_1\alpha_1}$ which bounds
$-o(\alpha_1,\alpha'_1)$ to promote $\widehat\psi_{1}$ to the
energy cut level $E_2$.
We note that we can still change  $(\widehat\psi_{1})_{\alpha'_1\alpha_1}$
by $d_0$-cocycle. We will use this freedom in the next step.
\par
Next we promote the homotopy $\widehat{\frak h}$.
Hereafter we denote by $\widehat\psi_{1}$ its promotion to the energy cut level $E_2$, which
we have just done above.
We now change our notation and regard $\widehat\psi_{j}$, $\widehat\psi_{21}$ etc.
as partial cochain maps of energy cut level $E_2$.
So when we regard them as $: CF(*) \to CF(*)$,
we include the ${\alpha\alpha'}$-component where
$E(\alpha') - E(\alpha) = E_2 - c$.
\par
We consider $(\alpha_1,\alpha'_2) \in \frak A_1 \times \frak A'_2$
with
$E(\alpha'_2) - E(\alpha_1) = E_2-c$.
We define
\begin{equation}\label{form1644}
o(\alpha_1,\alpha'_2)
=
\left((\widehat\psi'_{21} \circ \widehat\psi_{1}
- \widehat\psi_{2} \circ \widehat\psi'_{21})
-
(\widehat d'_{2} \circ \widehat{\frak h}
+ \widehat{\frak h}\circ \widehat d_{1})\right)_{\alpha'_2\alpha_1}.
\end{equation}
We note that the right hand side is $0$ if
$E(\alpha'_2) - E(\alpha_1) < E_2-c$ by
(\ref{gapconsthomotopro}) and the assumption.
Moreover we have
\begin{equation}\label{formula1747}
\aligned
0 = &\left(\widehat d'_2 \circ (\widehat\psi'_{21} \circ \widehat\psi_{1}
- \widehat\psi_{2} \circ \widehat\psi'_{21})
-
(\widehat\psi'_{21} \circ \widehat\psi_{1}
- \widehat\psi_{2} \circ \widehat\psi'_{21})  \circ \widehat d_1
\right.
\\
&\qquad\qquad\quad\left.-
\widehat d'_2 \circ (\widehat d'_{2} \circ \widehat{\frak h}
+ \widehat{\frak h}\circ \widehat d_{1})
+
(\widehat d'_{2} \circ \widehat{\frak h}
+ \widehat{\frak h}\circ \widehat d_{1})
\circ \widehat d_1\right)_{\alpha'_2\alpha_1}.
\endaligned
\end{equation}
In fact, we already promoted $\widehat d'_i$, $\widehat d_i$,
$\widehat{\psi}'_{21}$, $\widehat{\psi}'_{i}$ to the energy cut level
$E_2$ and in our notation (which we changed at the beginning of the
construction of $\widehat{\frak h}$) $\widehat{\psi}'_{21}$, $\widehat{\psi}'_{i}$
contain the components up to the energy cut level $E_2$.
Therefore the equalities $\widehat d'_{2} \circ \widehat{\psi}'_{21}
= \widehat{\psi}'_{21} \circ \widehat d'_{1}$ etc. hold
up to the  energy cut level $E_2$.
It implies that the first line of (\ref{formula1747}) vanishes.
The second line vanishes obviously.
We use them to obtain
\begin{equation}
d_0 \circ o(\alpha_1,\alpha'_2)
+ o(\alpha_1,\alpha'_2) \circ d_0 = 0.
\end{equation}
We now use the freedom to change $(\widehat\psi_{1})_{\alpha_1\alpha'_1}$
by $d_0$-cocycle.
Then we may assume (\ref{form1644}) is a $d_0$-coboundary.
Thus we can promote the cochain homotopy $\widehat{\frak h}$.
\end{proof}
\begin{proof}
[Proof of Theorem \ref{indlinesysmainth1} (3)]
Let $\mathcal F\mathcal F_j$ ($j=a,b$) be
inductive systems of partial linear K-systems which consist of
partial linear K-systems
$\mathcal F_{j}^{i}$, $i=1,2,\dots$ and
partial morphisms
$\frak N_j^{i+1 i} : \mathcal F_{j}^{i} \to \mathcal F_{j}^{i+1}$.
The morphism $\mathcal F\mathcal F_a \to \mathcal F\mathcal F_b$ consists of
morphisms $\frak N\frak N^{i} : \mathcal F_{a}^{i} \to \mathcal F_{b}^{i}$
and homotopy.
Namely, we have a homotopy commutative diagram of partial morphisms:
\begin{equation}\label{diagmorphiindcommute}
\begin{CD}
\mathcal F_{b}^{i}
@ > {\frak N^{i+1 i}_b} >> \mathcal F_{b}^{i+1}\\
@ A{\frak N\frak N^{i}}AA @ AA{\frak N\frak N^{i+1}}A \\
\mathcal F_{a}^{i}
@ >{\frak N^{i+1 i}_a}>> \mathcal F_{a}^{i+1}
\end{CD}
\end{equation}
\par
In the proof of Theorem \ref{indlinesysmainth1} (1) given in
Subsection \ref{subsec:proofsec14main1},
we made a choice of
$$
\widehat{\mathcal U^{+}}(j,i;\alpha_-,\alpha_+),
\quad
\widehat{\frak S^{+}}(j,i;\alpha_-,\alpha_+)
$$
that are Kuranishi structures of the spaces of connecting orbits and
their CF-perturbations to define partial
cochain complex $(CF(\mathcal F_{j}^{i}),\{\frak m^{ji;\epsilon}_{1;\alpha_+\alpha_-}\})$ and this partial
cochain complex is defined from $ \mathcal F_{j}^{i}$.
We also made a choice of
$$
\widehat{\mathcal U^{+}}({\rm mor};j,i,i+1;\alpha_-,\alpha_+),
\quad
\widehat{\frak S^{+}}({\rm mor};j,i,i+1;\alpha_-,\alpha_+;\rho_{j}^{ii+1})
$$
that are Kuranishi structures of the spaces of interpolation spaces of $\frak N^{i+1 i}_j$ and their
CF-perturbations.
Here the parameter $\rho = \rho_{j}^{ii+1} \in (0,1]$ is taken as follows.
We need to fix CF-perturbations of
$\widehat{\mathcal U^{+}}(j,i;\alpha_-,\alpha_+)$ and of
$\widehat{\mathcal U^{+}}(j,i+1;\alpha_-,\alpha_+)$ with which our
CF-perturbations on the interpolation spaces are compatible.
On one of those Kuranishi structures,
we take $\widehat{\frak S^{+}}(j,i;\alpha_-,\alpha_+)$.
On the other Kuranishi structure,  we take
$\epsilon \mapsto \widehat{\frak S^{+}}(j,i+1;\alpha_-,\alpha_+)^{\rho\epsilon}$.
The parameter $\rho$ appears here.
Using these choices, we obtain partial cochain maps at the horizontal arrows in
Diagram \eqref{diagmorphiindcommute}.
\par
Next we apply Proposition \ref{prop1618} to the vertical arrows of Diagram (\ref{diagmorphiindcommute}).
Then we can choose
$$
\widehat{\mathcal U^{+}}({\rm mor};ab,i;\alpha_-,\alpha_+),
\quad
\widehat{\frak S^{+}}({\rm mor};ab,i;\alpha_-,\alpha_+;\rho_{ab}^{i}),
$$
and
$$
\widehat{\mathcal U^{+}}({\rm mor};ab,i+1;\alpha_-,\alpha_+),
\quad
\widehat{\frak S^{+}}({\rm mor};ab,i+1;\alpha_-,\alpha_+;\rho_{ab}^{i+1}),
$$
that are Kuranishi structures and
CF-perturbations of the interpolation spaces of $\frak{NN}^{i}$,
$\frak{NN}^{i+1}$, respectively.
\par
For each $i$ we take $\epsilon_{bi} \le \epsilon_{ai}  \le \epsilon_{4,i}$
such that $\epsilon_{a i+1} \le \epsilon_{ai}$, $\epsilon_{b i+1} \le \epsilon_{bi}$.
(We specify our choice of $\epsilon_{4,i}$ later.)
We take $\rho_{j}^{ii+1} = \epsilon_{j i+1}/\epsilon_{j i}$.
We put
\begin{equation}\label{CFjinoichi}
CF_{j}^{i}:= \left(CF(\mathcal F_{j}^{i}),\{\frak m^{ji;\epsilon}_{1;\alpha_+\alpha_-}\}\right).
\end{equation}
Then we have a diagram of partial cochain complexes.
\begin{equation}\label{diagmorphiindcommutealg}
\begin{CD}
CF_{b}^{i}
@ > >> CF_{b}^{i+1}\\
@ AAA @ AAA \\
CF_{a}^{i}
@ >>> CF_{a}^{i+1}
\end{CD}
\end{equation}
By construction the vertical arrows are partial cochain maps of energy cut level
$E^i$ and $E^{i+1}$, respectively.
\par
We now use Proposition \ref{homotopyexistpertstatement} to obtain
Kuranishi structures and CF-perturbations
on the interpolation spaces of the homotopy of Diagram (\ref{diagmorphiindcommute}) which are compatible
with the choices we had made for the interpolation spaces of the arrows of Diagram (\ref{diagmorphiindcommute})
and the space of connecting orbits. (We had made the choice of them already as explained above.)
We can take $\epsilon_{4,i}$ so small that this choice of Kuranishi structures and
CF-perturbations
determines a cochain homotopy which makes Diagram (\ref{diagmorphiindcommutealg}) commutative
up to cochain homotopy.
(This is a consequence of Lemma \ref{lem165rev1}.)
Recall from Remark \ref{rem:1815} the numbers
$\epsilon_{0,i},\epsilon_{1,i},\epsilon_{2,i}, \epsilon_{3,i}$ in (\ref{choiceofrho})
are already chosen.
We will need $\epsilon_{4,i}, \epsilon_{5,i}$ etc. for homotopy of homotopies etc.
(However we need only finitely many of them.)
\par
Therefore by Proposition \ref{prop1631} we can promote Diagram (\ref{diagmorphiindcommutealg})
to homotopy commutative diagram of the partial cochain maps of energy cut level $E_{i+1}$.
\par
The rest of the proof is purely algebraic.
We now consider the following diagram:
\begin{equation}\label{diagmorphiindcommutealgseri}
\begin{CD}
CF_{b}^{1}@>>>\dots @ > >>CF_{b}^{i}
@ > >> CF_{b}^{i+1}@> >>  \dots\\
@ AAA && @ AAA @ AAA &&   \\
CF_{a}^{1} @ > >> \dots @ > >>CF_{a}^{i}
@ > >> CF_{a}^{i+1}@> >>   \dots
\end{CD}
\end{equation}
Note that our construction of cochain complex $CF(\mathcal F_a)$ is done by
promoting the horizontal lines inductively to one of energy cut level $E_k$ and
taking limit $k\to \infty$.
We do so for both of the horizontal lines.
Then
we continue to promote the vertical lines so that the whole diagram becomes
homotopy commutative.
Thus we obtain a cochain map $CF(\mathcal F_a) \to CF(\mathcal F_b)$.
This finishes the proof of Theorem \ref{indlinesysmainth1} (3).
\end{proof}
\begin{proof}
[Proof of Theorem \ref{indlinesysmainth1} (2) (a)-(d)]
When we proved Theorem \ref{indlinesysmainth1} (1)
in Subsection \ref{subsec:proofsec14main1},
we made the following choices.
\begin{enumerate}
\item
For each $i$ the pairs $\left(\widehat{\mathcal U^{+}}(i;\alpha_-,\alpha_+),
\widehat{\frak S^{+}}(i;\alpha_-,\alpha_+)\right)$ for various $\alpha_-, \alpha_+$
that are a system of Kuranishi structures of the spaces of connecting orbits and their CF-perturbations to define a partial
cochain complex
$$
(CF(\mathcal F^{i}),\{\frak m^{i;\epsilon}_{1;\alpha_+\alpha_-}\})
$$
and this partial
cochain complex is defined from the partial linear K-system $\mathcal F^{i}$ with
energy cut level $E_{\frak E}^{k_i}$.
\item
For each $i$ the pairs
$$
\left(\widehat{\mathcal U^{+}}({\rm mor};i,i+1;\alpha_-,\alpha_+),
\widehat{\frak S^{+}}({\rm mor};i,i+1;\alpha_-,\alpha_+;\rho_{j}^{ii+1})\right)
$$
for various $\alpha_-, \alpha_+$ are a system of
Kuranishi structures of the interpolation spaces of $\frak N^{i+1 i}$ and their
CF-perturbations. Here the parameter $\rho = \rho_{j}^{ii+1} \in (0,1]$ is as explained during the proof of Theorem \ref{indlinesysmainth1} (3).
\item
The small numbers $\epsilon_i$.
\end{enumerate}
We then defined $\left(CF(\mathcal F^{i}),\{\frak m^{i;\epsilon_i}_{1;\alpha_+\alpha_-}\}\right)$
by using the choices (1),(3).
Moreover using the choices (2), (3), we defined
partial cochain maps
\begin{equation}\label{formula16.49}
\widehat\psi^{i} ~:~
\left(CF(\mathcal F^{i}),\{\frak m^{i;\epsilon_i}_{1;\alpha_+\alpha_-}\}\right)
\longrightarrow
\left(CF(\mathcal F^{i+1}),\{\frak m^{i+1;\epsilon_{i+1}}_{1;\alpha_+\alpha_-}\}\right).
\end{equation}
Finally using (\ref{formula16.49}) and an algebraic result (Lemma \ref{lem1623}),
we promoted partial cochain complexes
$\left(CF(\mathcal F^{i}),\{\frak m^{i;\epsilon_i}_{1;\alpha_+\alpha_-}\}\right)$
and  partial cochain maps $\widehat\psi^{i} $ to cochain complexes and cochain maps.
This algebraic process also involves choices.
\par
We will prove that the resulting cochain complex is independent of
those choices up to cochain homotopy equivalence.
Let
$$
\left(\widehat{\mathcal U^{+}}(ji;\alpha_-,\alpha_+),
\widehat{\frak S^{+}}(ji;\alpha_-,\alpha_+)\right)
$$
and
$$
\left(\widehat{\mathcal U^{+}}({\rm mor},j;i,i+1;\alpha_-,\alpha_+),
\widehat{\frak S^{+}}({\rm mor},j;i,i+1;\alpha_-,\alpha_+;\rho_{j}^{ii+1})\right),
$$
and
$\epsilon_{ji}$
be two choices, where $j=a,b$.
\par
We consider the next diagram, which defines a morphism  $\mathcal{FF}
\to \mathcal{FF}$.

\begin{equation}\label{shiftmorphism}
\begin{CD}
\mathcal F^{1} @>{\frak N^{21}}>>\dots @ > {\frak N^{i i-1}}>>\mathcal F^{i}
@ > {\frak N^{i+1 i}}>> \mathcal F^{i+1} @> {\frak N^{i+2 i+1}}>>  \dots\\
@ AA{\frak{ID}}A && @ AA {\frak{ID}}A @ AA {\frak{ID}}A &&   \\
\mathcal F^{1} @ >{\frak N^{21}}>> \dots @ > {\frak N^{i i-1}}>>\mathcal F^{i}
@ >{\frak N^{i+1 i}} >> \mathcal F^{i+1} @>{\frak N^{i+2 i+1}} >>   \dots
\end{CD}
\end{equation}
The homotopy commutativity of Diagram (\ref{shiftmorphism}) follows
from Proposition \ref{properidentity}.
\par
Now we apply Theorem \ref{indlinesysmainth1} (3).
Namely we apply the choice we made for $j=b$ for the first line
and the choice we made for $j=a$ for the second line.
We also use the particular way to promote the inductive
system of partial cochain complexes which we used for choices $j=b$ and
$j=a$. (We obtain those cochain complexes
from the first and second lines of Diagram (\ref{shiftmorphism}).)
\par
Now we apply Proposition \ref{prop1631}
and obtain a cochain map
$
CF_{a} \to CF_{b}.
$
Here $CF_{a}$ (resp. $CF_{b}$) is the cochain complex we obtain
by this promotion of the second line (resp. first line).
The proof of Theorem \ref{indlinesysmainth1} (2) (except (e)) is complete.
\end{proof}
\begin{proof}
[Proof of Theorem \ref{linesysmainth1} (2) (a)-(e)]
This is nothing but a special case of Theorem \ref{indlinesysmainth1} (2)
where $\frak N^{i+1 i}$ is the identity morphism.
\end{proof}
\begin{proof}[Proof of Theorem \ref{linesysmainth2} (1)]
We can use Proposition \ref{properidentity}
to obtain the following homotopy commutative diagram:
\begin{equation}\label{inductivesystemwithid}
\begin{CD}
\mathcal F_2 @>{\frak{ID}}>>\dots @ > {\frak{ID}}>>\mathcal F_2
@ > {\frak{ID}}>> \mathcal F_2 @> {\frak{ID}}>>  \dots\\
@ AA{\frak N}A && @ AA {\frak N}A @ AA {\frak N}A &&   \\
\mathcal F_{1} @ >{\frak{ID}}>> \dots @ >{\frak{ID}}>>\mathcal F_{1}
@ >{\frak{ID}} >> \mathcal F_{1} @>{\frak{ID}}>>   \dots
\end{CD}
\end{equation}
In fact, we have
$
{\frak{ID}} \circ \frak N \sim \frak N \sim \frak N \circ  {\frak{ID}},
$
by Proposition \ref{properidentity}.
\par
Following the proof of Theorem \ref{linesysmainth1} (1),
we make the following choices.
\begin{enumerate}
\item
For each $j=1,2$ we regard $\mathcal F_j$ as a partial linear K-system
of energy cut level $E^i$.
We write it as $\mathcal F_{j}^{i}$.
\item
We take $\widehat{\mathcal U^+}(j,i;\alpha_-,\alpha_+)$
and $\widehat{\frak S^+}(j,i;\alpha_-,\alpha_+)$ that
are a Kuranishi structure of the space of connection orbits of $\mathcal F_{j}^{i}$
and its CF-perturbation, respectively.
They satisfy the conclusion of Proposition \ref{prop161}.
\item
We take
$\widehat{\mathcal U^+}({\rm mor},j,i i+1;\alpha_-,\alpha_+)$
and $\widehat{\frak S^+}({\rm mor},j,i i+1;\alpha_-,\alpha_+)$
that are a Kuranishi structure of the interpolation space of
$\frak{ID} : \mathcal F_{j}^{i} \to \mathcal F_{j}^{i+1}$
and its CF-perturbation, respectively.
They satisfy the conclusion of Proposition \ref{prop1618}.
\item
We also take
$\widehat{\mathcal U^+}({\rm mor},12,i;\alpha_1,\alpha_2)$
and
$\widehat{\frak S^+}({\rm mor},12,i;\alpha_1,\alpha_2)$
that are
a Kuranishi structure of the interpolation space of
$\frak{N} : \mathcal F_{1}^{i} \to \mathcal F_{2}^{i}$
and its CF-perturbation, respectively.
They satisfy the conclusion of Proposition \ref{prop1618}.
\item We take
$\widehat{\mathcal U^+}({\rm mor},i;\alpha_-,\alpha_+;[0,1])$
and $\widehat{\frak S^+}({\rm mor},i;\alpha_-,\alpha_+;[0,1])$
that are a Kuranishi structure of the interpolation space of the homotopy
$
{\frak{ID}} \circ \frak N \sim \frak N \circ  {\frak{ID}}
: \mathcal F_{1}^i \to \mathcal F_{2}^{i+1}$
and its CF-perturbation, respectively.
The Kuranishi structure is compatible with ones in
Item (1)(2)(3)(4) at the boundary.
\end{enumerate}
By these choices the geometric Diagram (\ref{inductivesystemwithid})
is converted to the algebraic diagram below:
\begin{equation}\label{diag1637rev}
\begin{CD}
\left(CF(\mathcal F_2),\{\frak m^{2 i \epsilon_{2i}}_{1;\alpha_+,\alpha_-}\}\right)
@ > {\widehat\psi_2^{i+1 i}} >>\left(CF(\mathcal F_2),\{\frak m^{2i \epsilon_{2i+1}}_{1;\alpha_+,\alpha_-}\}\right)
\\
@ A{\frak n^{i}}AA @ AA{\frak n^{i+1}}A \\
\left(CF(\mathcal F_1),\{\frak m^{1i,\epsilon_{1i}}_{1;\alpha_+,\alpha_-}\}\right)
@ > {\widehat\psi_1^{i+1 i}} >>\left(CF(\mathcal F_1),\{\frak m^{1i+1,\epsilon_{1i+1}}_{1;\alpha_+,\alpha_-}\}\right)
\end{CD}
\end{equation}
Here
\begin{enumerate}
\item[$\bullet$]
$\left(CF(\mathcal F_j),\{\frak m^{ji,\epsilon_{ji}}_{1;\alpha_+,\alpha_-}\}\right)$
is obtained by the choice (2) above.
\item[$\bullet$]
$\widehat\psi_j^{i+1 i}$, $j=1,2$ is obtained by the choice (3) above.
\item[$\bullet$]
$\frak n^{i}$ is obtained by the choice (4) above.
\item[$\bullet$]
By the choice (5) above we obtain a cochain homotopy between
$\widehat\psi_2^{i+1 i} \circ \frak n^{i}$ and $\frak n^{i+1} \circ \widehat\psi_1^{i+1 i}$.
We denote it by $\frak h_i$.
\end{enumerate}
Here we use Lemmas \ref{lem1721} and \ref{17222lem}
to prove that the composition $\widehat\psi_2^{i+1 i} \circ \frak n^{i}$
(resp. $\frak n^{i+1} \circ \widehat\psi_1^{i+1 i}$)
is the cochain map associated to the morphism
that is the composition of the identity morphism
and $\frak N^i$
(resp. $\frak N^i$ and the identity morphism).
\par
We note that at this stage the energy cut level of the objects
(partial cochain maps and partial cochain complexes) in the
right vertical line is $E^{i+1}$ and the energy cut
level of all the other objects in (\ref{diag1637rev})
(and the partial cochain homotopy $\frak h_i$) is $E^{i}$.
\par
Now we start promoting the objects in Diagram (\ref{diag1637rev}).
First according to the proof of Theorem \ref{linesysmainth1} (1),
we promote all the objects in the upper and lower
horizontal lines in Diagram (\ref{diag1637rev})
to a cochain complex and cochain map.
(Namely we promote them to the energy cut level $\infty$.)
\par
Thus all the objects other than $\frak n^i$, $\frak n^{i+1}$
and the cochain homotopy
in Diagram (\ref{diag1637rev}) are promoted to the energy cut
level $\infty$.
We now use the fact that the energy cut level of $\frak n^{i+1}$
is $E^{i+1}$ to promote
$\frak n^i$ and $\frak h^i$
to the energy cut level $E^{i+1}$.
We use Proposition \ref{prop1631} to do so.
\par
Then by induction using Diagrams (\ref{diag1637rev}) for all $i$,
we can promote $\frak n^i$, $\frak n^{i+1}$
to the energy cut level $\infty$ so that
Diagram (\ref{diag1637rev}) commutes up to cochain homotopy.
The proof of Theorem \ref{linesysmainth2} (1) is complete.
\end{proof}

\subsection{Construction of higher homotopy}
\label{subsec:consthigerhomotopy}

In this subsection, we generalize Proposition \ref{homotopyexistpertstatement}
to the case of $P$-parametrized morphisms.

\begin{shitu}\label{situation1633}
Suppose we are given a $P$-parametrized family
of morphisms of (partial) linear K-system.
Let $E_0$ be its energy cut level and $c$ be its energy loss.
Let $\mathcal N(\alpha_1,\alpha_2;P)$
be its interpolation spaces.
\par
Let
$\alpha_i \in \frak A_i$ with
$\alpha_2 - \alpha_1 \le E_0 - c$.
\begin{enumerate}
\item
We are given a $\tau$-$\frak C^h$-collared
Kuranishi structure
$\widehat{\mathcal U^{+}}(\alpha_1,\alpha_2;\widehat S_k(P))$
on the underlying topological space of
$
\mathcal N(\alpha_1,\alpha_2;S_k(P))^{\frak C^h\boxplus\tau_0}
$.
(Here $\frak C^h$ is the
{\it horizontal boundary}
\index{horizontal boundary}\index{boundary ! horizontal boundary}
component,
that is the complement of
the inverse image ${\rm ev}_{S_k (P)}^{-1}(\partial S_k (P))$ of
the evaluation map
$$
{\rm ev}_{S_k (P)} : \mathcal N(\alpha_1,\alpha_2;S_k (P))
\longrightarrow \widehat{S}_k (P)
$$
in $\partial \mathcal N(\alpha_1,\alpha_2;S_k (P))$.
See Definition \ref{defn1726} for the case $P=[1,2]$.)
\par
We assume that there exists an embedding to
$\widehat{\mathcal U^{+}}(\alpha_1,\alpha_2;\widehat S_k(P))$
from
the corner trivialization of the
Kuranishi structure of  $\mathcal N(\alpha_1,\alpha_2;S_k(P))$
in Lemma \ref{Skparafamilylem},
by which $\widehat{\mathcal U^{+}}(\alpha_1,\alpha_2;\widehat S_k(P))$
becomes a thickening of
$$
\mathcal N(\alpha_1,\alpha_2;S_k(P))^{\frak C^h\boxplus\tau_0}.
$$
\item
Evaluation maps ${\rm ev}_{\pm}$ and ${\rm ev}_{S_k(P)}$
extend to
$\widehat{\mathcal U^{+}}(\alpha_1,\alpha_2;\widehat S_k(P))$
and the above mentioned
embedding is compatible with the evaluation maps.
There exists a periodicity isomorphism which commutes with the embedding.
\item
We are given a $\tau$-collared CF-perturbation
$\widehat{\frak S^{+}}(\alpha_1,\alpha_2;\widehat S_k(P))$
of the Kuranishi structure in (1).
It is compatible with the periodicity isomorphism.
\item
The $\tau$-collared CF-perturbation $\widehat{\frak S^{+}}(\alpha_1,\alpha_2;\widehat S_k(P))$
is transversal to zero.
Moreover the evaluation map
$$
({\rm ev}_{+},{\rm ev}_{S_k(P)}) :
\widehat{\mathcal U^{+}}(\alpha_1,\alpha_2;\widehat S_k(P))
\longrightarrow R_{\alpha_2} \times \widehat S_k(P)
$$
is stratified strongly submersive with respect
to $\widehat{\frak S^{+}}(\alpha_1,\alpha_2;\widehat S_k(P))$.
\item
The Kuranishi structure
$\widehat{\mathcal U^{+}}(\alpha_1,\alpha_2;\widehat S_k(P))$
and
the CF-perturbation
$$
\widehat{\frak S^{+}}(\alpha_1,\alpha_2;\widehat S_k(P))
$$
are compatible with the isomorphisms in
Conditions \ref{famiboudaru22}, \ref{boundarycompPban1}, \ref{furthercompatifiber}.
\end{enumerate}
$\blacksquare$
\end{shitu}
\begin{prop}\label{prop1618rev}
Suppose we are in Situation \ref{situation1633}.
Then there exists a $\tau'$-$\frak C^h$-collared Kuranishi structure
$\widehat{\mathcal U^+}(\alpha_1,\alpha_2;P)$
on the underlying topological space of
$\mathcal N(\alpha_1,\alpha_2;P)^{\frak C^h\boxplus\tau_0}$
and a  $\tau'$- $\frak C^h$-collared CF-perturbation
$\widehat{\frak S^+}(\alpha_1,\alpha_2;P)$
of $\widehat{\mathcal U^+}(\alpha_1,\alpha_2;P)$ such that
they have the following properties.
\begin{enumerate}
\item
There exists an embedding from
$\mathcal N(\alpha_1,\alpha_2;P)^{\frak C^h\boxplus\tau_0}$
to $\widehat{\mathcal U^+}(\alpha_1,\alpha_2;P)$
by which the later becomes a thickening of the former.
\item
Evaluation maps ${\rm ev}_{\pm}$ and ${\rm ev}_{P}$
extend to $\widehat{\mathcal U^+}(\alpha_1,\alpha_2;P)$ and the above mentioned
embedding is compatible with them.
\item
There exists a periodicity isomorphism
$$
\widehat{\mathcal U^+}(\alpha_1,\alpha_2;P)
\longrightarrow
\widehat{\mathcal U^+}(\beta\alpha_1,\beta\alpha_2;P)
$$
which commutes with the embedding.
The pull-back of
$\widehat{\frak S^+}(\beta\alpha_1,\beta\alpha_2;P)$
by the periodicity isomorphism
is isomorphic to
$\widehat{\frak S^+}(\alpha_1,\alpha_2;P)$.
\item
The CF-perturbation
$\widehat{\frak S^+}(\alpha_1,\alpha_2;P)$  is transversal to $0$.
Moreover the evaluation map
$$
({\rm ev}_+,{\rm ev}_P) : \mathcal N(\alpha_1,\alpha_2;P)^{\frak C^h\boxplus\tau_0}
\longrightarrow R_{\alpha_2} \times P
$$
is strongly submersive with respect to
$\widehat{\frak S^+}(\alpha_1,\alpha_2;P)$.
\item
The fiber product
$$
\widehat{\mathcal U^+}(\alpha_1,\alpha_2;P)
\times_{P} \widehat S_k(P)
$$
is isomorphic to
$\widehat{\mathcal U^{+}}(\alpha_1,\alpha_2;\widehat S_k(P))$.
This isomorphism is compatible with the periodicity isomorphism
and evaluation maps.
It is also compatible with the embedding in (1) and one in
Situation \ref{situation1633} (1).
\item
The fiber product
$$
\widehat{\frak S^+}(\alpha_1,\alpha_2;P)
\times_{P} \widehat S_k(P)
$$
is equivalent to $\widehat{\frak S^{+}}(\alpha_1,\alpha_2;\widehat S_k(P))$.
\item
There exists an isomorphism:
\footnote{See Remark \ref{rem:FiberProdOrd}).}
\begin{equation}\label{boundaryiso1643}
\aligned
&\partial
\left(
\widehat{\mathcal U^+}(\alpha_1,\alpha_2;P)
\right)
\\
\cong
& \sqcup \coprod_{\alpha \in \frak A_1}
(-1)^{\dim \widehat{\mathcal U^+}(\alpha,\alpha_2;P)}
\left(
\widehat{\mathcal U^+}(\alpha,\alpha_2;P)
\,\,{}_{{\rm ev}_{-}}\times_{{\rm ev}_{+}}
\widehat{\mathcal U^+}(1;{\alpha_1},{\alpha})
\right)\\
& \coprod_{\alpha \in \frak A_2}
(-1)^{\dim \widehat{\mathcal U^+}(2;{\alpha},{\alpha_2})}\left(
\widehat{\mathcal U^+}(2;{\alpha},{\alpha_2})
\,\,{}_{{\rm ev}_{-}}\times_{{\rm ev}_{+}}
\widehat{\mathcal U^+}(\alpha_1,\alpha;P)
\right)\\
&\sqcup\widehat{\mathcal U^+}(\alpha_1,\alpha;\partial P).
\endaligned
\end{equation}
(Recall that $\widehat{\mathcal U^+}(1;{\alpha_1},{\alpha})$,
$\widehat{\mathcal U^+}(2;{\alpha},{\alpha_2})$ are
Kuranishi structures on $\mathcal M^1({\alpha_1},{\alpha})$,
$\mathcal M^2({\alpha},{\alpha_2})$, respectively,
given in Proposition \ref{prop161}.)
\par
This isomorphism is compatible with the isomorphism in Condition
\ref{famiboudaru22} via the embedding in (1).
It is also compatible with the orientation isomorphism, the periodicity isomorphism and the evaluation
maps.
\item
The normalized corner
$\widehat S_k(\widehat{\mathcal U^+}(\alpha_1,\alpha_2;P))$
is decomposed to the disjoint union
of the following fiber products.
\begin{equation}\label{cornerisostateconst}
\aligned
&\widehat{\mathcal U^+}(1;\alpha_-,\alpha_1)
\times_{R_{\alpha_1}}
\dots
\times_{R_{\alpha_{k_1}-1}}
\widehat{\mathcal U^+}(1;\alpha_{{k_1}-1},\alpha_{{k_1}}) \\
&
\times_{R_{\alpha_{k_1}}}
\widehat{\mathcal U^+}({\alpha_{k_1}},{\alpha_{k_1+1}};\widehat S_{k_3}(P)) \\
&
\times_{R_{\alpha_{k_1+1}}}
\widehat{\mathcal U^+}(2;\alpha_{k_1+1},\alpha_{k_1+2})
\times_{R_{\alpha_{k_1+2}}}
\dots
\times_{R_{\alpha_{k_1+k_2}}}
\widehat{\mathcal U^+}(2;\alpha_{{k_1}+k_2},\alpha_+).
\endaligned
\end{equation}
This isomorphism is compatible with the isomorphism in Condition
\ref{boundarycompPban1} via the embedding in (1).
It is also compatible with the periodicity isomorphism and the evaluation
maps.
\item
The isomorphism (\ref{cornerisostateconst})
satisfies the same compatibility conditions claimed for
$\mathcal N(\alpha_1,\alpha_2;P)$ in
Condition \ref{furthercompatifiber}.
\item
The CF-perturbation
$\widehat{\frak S^+}(\alpha_1,\alpha_2;P)$  is
compatible with the isomorphisms (\ref{boundaryiso1643})
and
(\ref{cornerisostateconst}).
\item
When the $\tau$-collared CF-perturbation
$\widehat{\frak S^{+}}(\alpha_1,\alpha_2;\widehat S_k(P))$ given in
Situation \ref{situation1633} varies in a uniform family,
we can take
$\widehat{\frak S^+}(\alpha_1,\alpha_2;P)$ to be uniform.
\end{enumerate}
\end{prop}
\begin{proof}
The proof is entirely the same as the proof of Proposition \ref{prop161}.
\end{proof}
The translation to algebra is fairly immediate.
\begin{defn}\label{chainmaphigherhomoto}
In the situation of Proposition \ref{prop1618rev}, we define
\begin{equation}
\psi^{P,\epsilon}_{\alpha_2,\alpha_1} ~:~
\Omega(R_{\alpha_1};o_{R_{\alpha_1}}) \longrightarrow
 \Omega(R_{\alpha_2};o_{R_{\alpha_2}})
\end{equation}
by
\begin{equation}\label{form1655}
\psi^{P,\epsilon}_{\alpha_2,\alpha_1}(h)
=
{\rm ev}_{+}!({\rm ev}_{-}^* h;\widehat{\frak S^{+ \epsilon}}(\alpha_1,\alpha_2;P)).
\end{equation}
Here the right hand side is defined by
Definition \ref{pushoutdeftau} on
the K-spaces
$$
({\mathcal N}({\alpha_1},{\alpha_2};P)^{\boxplus\tau_0},
\widehat{\mathcal U^+}({\rm mor};\alpha_1,\alpha_2;P)).
$$
\end{defn}
The degree of $\psi^{P,\epsilon}_{\alpha_2,\alpha_1}$ after
shifted is $-\dim P$.
If the energy loss of our parametrized family of morphisms is
$c$, the map $\psi^{P,\epsilon}_{\alpha_2,\alpha_1}$ induces
$$
\frak F^{\lambda}CF(\mathcal F_1) \longrightarrow
 \frak F^{\lambda-c}CF(\mathcal F_2)
$$
where the filtration $\frak F^{\lambda}$ is defined in
Definition \ref{Fvectspace} (2)(3).
\begin{lem}\label{lem165rev2}
The operators $\psi^{P,\epsilon}_{\alpha_2,\alpha_1}$
satisfy the following equality in the sense of $(\flat)$.
\begin{equation}\label{form1657}
\aligned
&d_0\circ \psi^{P,\epsilon}_{\alpha_2,\alpha_1}
- (-1)^{\deg P}\psi^{P,\epsilon}_{\alpha_2,\alpha_1} \circ d_0
+ \psi^{\partial P,\epsilon}_{\alpha_2,\alpha_1}
\\
& + (-1)^{\deg P}
\sum_{\alpha'_1}
\psi^{P,\epsilon}_{\alpha_2,\alpha'_1}
\circ
\frak m^{1,\epsilon}_{1;\alpha'_1,\alpha_1}
-
\sum_{\alpha'_2}
\frak m^{2,\epsilon}_{1;\alpha_2,\alpha'_2}
\circ
\psi^{P,\epsilon}_{\alpha'_2,\alpha_1}
= 0.
\endaligned
\end{equation}
Here the first sum in the second line is taken over
$\alpha'_1 \in \frak A_1$ with $E(\alpha_1) < E(\alpha'_1) \le E(\alpha_2) + c$
and the second sum in the second line is taken over
$\alpha'_2 \in \frak A_2$ with $E(\alpha_1)- c \le  E(\alpha'_2) < E(\alpha_2)$.
The number $c$ is the energy loss of our morphism.
\end{lem}
\begin{proof}
The proof is by Stokes' formula and the composition formula and
is entirely similar to the proof of Lemma \ref{lem165rev}.
\end{proof}

\subsection{Proof of Theorem \ref{indlinesysmainth1} (2)(e), (4)-(6) and Theorem \ref{linesysmainth1} (2)(f)}
\label{subsec:proofcomplete}

We begin with an algebraic result
that is similar to Proposition \ref{prop1631}
and is a baby version of
\cite[Theorem 7.2.212]{fooobook2}.

\begin{shitu}\label{situ1639}
For $j=1,2$, let
$(CF(\mathcal F_{j}^{i}),\hat d_{j}^{i})$,
$(CF(\mathcal F_{j}^{i+1}),\hat d_{j}^{i+1})$ be partial cochain complexes
of energy cut level  $E^{i+1}$.
\begin{enumerate}
\item
For $j=1,2$,
$\psi_j^{i+1 i} : CF(\mathcal F_{j}^{i}) \to CF(\mathcal F_{j}^{i+1})$
is a partial cochain map of energy cut level $E^{i+1}$
and energy loss $0$.
Moreover, we assume $\psi_j^{i+1 i}$ induces an isomorphism
modulo $T^{\epsilon}$ for small $\epsilon > 0$.
\item
For $k=a,b$ and $\ell =i,i+1$,
$\frak n_k^{\ell} :  CF(\mathcal F_{1}^{\ell}) \to CF(\mathcal F_{2}^{\ell})$
is a partial cochain map of energy cut level $E^{i+1}$
and energy loss $c$.
\item
$\frak h_{ab}^{i} :  CF(\mathcal F_{1}^{i}) \to CF(\mathcal F_{2}^{i})$
is a cochain homotopy between  $\frak n_a^{i}$ and $\frak n_b^{i}$
of energy cut level $E^{i}$
and energy loss $c$.
$\frak h_{ab}^{i+1} :  CF(\mathcal F_{1}^{i+1}) \to CF(\mathcal F_{2}^{i+1})$
is a cochain homotopy between  $\frak n_a^{i+1}$ and $\frak n_b^{i+1}$
of energy cut level $E^{i+1}$
and energy loss $c$. See Diagram \eqref{eq:htphtp}.
\item
For $k=a,b$,
$\frak h_{k}^{i+1 i} :  CF(\mathcal F_{1}^{i}) \to CF(\mathcal F_{2}^{i+1})$
is a partial cochain homotopy between $\frak n_k^{i+1}
\circ \psi_1^{i+1 i}$ and $\psi_2^{i+1 i} \circ \frak n_k^{i}$
of energy cut level $E^{i+1}$
and energy loss $c$.
\end{enumerate}
$\blacksquare$
\end{shitu}

\begin{equation}\label{eq:htphtp}
\vcenter{
\xymatrix{
CF(\mathcal F_{1}^{i}) \ar[rr]^-{\psi^{i+1 i}_{1}} \ar@/_10pt/[dd]_{\frak n^i_a} \ar@/^10pt/[dd]^{\frak n^i_b} & {} &
CF(\mathcal F_{1}^{i+1}) \ar@/_10pt/[dd]_{\frak n^{i+1}_a} \ar@/^10pt/[dd]^{\frak n^{i+1}_b}  \\
\overset{\frak h^{i}_{ab}}{\longrightarrow} & {} &
\overset{\frak h^{i+1}_{ab}}{\longrightarrow} \\
CF(\mathcal F_{2}^{i}) \ar[rr]^-{\psi^{i+1 i}_{2}} & {} & CF(\mathcal F_{2}^{i+1})
}
}
\end{equation}

\begin{prop}\label{prop1639}
Suppose in Situation \ref{situ1639} there exists (homotopy of homotopies)
$$
\frak H_{ab}^{i+1 i} ~:~ CF(\mathcal F_{1}^{i}) \longrightarrow
 CF(\mathcal F_{2}^{i+1})
$$
which satisfies
\begin{equation}\label{diagram1663}
\aligned
&\hat d_{2}^{i+1} \circ \frak H_{ab}^{i+1 i}
-
\frak H_{ab}^{i+1 i} \circ \hat d_{1}^{i}
\\
&= \frak h_b^{i+1 i} - \frak h_a^{i+1 i}
+ \frak h_{ab}^{i+1} \circ \psi_1^{i+1 i}
-\psi_2^{i+1 i} \circ \frak h_{ab}^{i}
\endaligned
\end{equation}
as equality of maps of energy cut level $E^i$ and energy loss $c$.
\par
Then we can promote $\frak h_{ab}^{i}$ and $\frak H_{ab}^{i+1 i}$
to the energy cut level $E^{i+1}$ so that Formula (\ref{diagram1663})
holds as an equality of maps of energy cut level $E^{i+1}$ and energy loss $c$.
\end{prop}
\begin{proof}
We first promote $\frak h_{ab}^{i}$.
Let us consider $\alpha_1,\alpha_2$ with $E(\alpha_2) - E(\alpha_1)
= E^{i+1} - c$.
We will find
$$
(\frak h_{ab}^{i})_{\alpha_2\alpha_1} ~:~
 CF(\mathcal F_{1}^{i}) \longrightarrow CF(\mathcal F_{2}^{i})
$$
such that
\begin{equation}\label{form1664}
\aligned
&d_0 \circ (\frak h_{ab}^{i})_{\alpha_2\alpha_1}
+ (\frak h_{ab}^{i})_{\alpha_2\alpha_1} \circ d_0 \\
=&
- \sum_{\alpha'_2} \frak m^2_{1,\alpha_2\alpha'_2}\circ (\frak  h_{ab}^{i})_{\alpha'_2\alpha_1}
- \sum_{\alpha'_1}
(\frak h_{ab}^{i})_{\alpha_2\alpha'_1}
\circ \frak m^1_{1,\alpha'_1\alpha_1}\\
&+ (\frak n_b^i)_{\alpha_2\alpha_1} - (\frak n_a^i)_{\alpha_2\alpha_1}.
\endaligned
\end{equation}
Let $o(\alpha_1,\alpha_2)$ be the right hand side of (\ref{form1664}).
Then we have
$$
o(\alpha_1,\alpha_2)
=
\left(-\hat d_{2}^{i} \circ {\frak h}_{ab}^{i}
-{\frak h}_{ab}^{i} \circ \hat d_{1}^{i}
+ \frak n_b^i - \frak n_a^i\right)_{\alpha_2\alpha_1}.
$$
\begin{lem}
$$
d_0 \circ o(\alpha_1,\alpha_2) - o(\alpha_1,\alpha_2)\circ d_0 = 0.
$$
\end{lem}
\begin{proof}
We observe
$$
\aligned
&\left(\hat d_{2}^{i}\circ \left(-\hat d_{2}^{i} \circ {\frak h}_{ab}^{i}
-{\frak h}_{ab}^{i} \circ \hat d_{1}^{i}
+ \frak n_b^i - \frak n_a^i\right)\right.\\
&-
\left.\left(-\hat d_{2}^{i} \circ {\frak h}_{ab}^{i}
-{\frak h}_{ab}^{i} \circ \hat d_{1}^{i}
+ \frak n_b^i - \frak n_a^i\right) \circ \hat d_{1}^{i}
\right)_{\alpha_2\alpha_1}
=0.
\endaligned
$$
This is a consequence of
$\hat d_{2}^{i} \circ \hat d_{2}^{i} = \hat d_{1}^{i} \circ \hat d_{1}^{i} = 0$
and  Situation \ref{situ1639} (2).
On the other hand we have
$$
\left(-\hat d_{2}^{i} \circ {\frak h}_{ab}^{i}
-{\frak h}_{ab}^{i} \circ \hat d_{1}^{i}
+ \frak n_b^i - \frak n_a^i
\right)_{\alpha'_2\alpha'_1} = 0
$$
if $E(\alpha'_2) - E(\alpha'_1)
\le E^{i} - c$
by Situation \ref{situ1639} (4).
The lemma follows.
\end{proof}
We will show that $o(\alpha_1,\alpha_2)$
is a coboundary.
Below $\equiv$ stands for modulo $d_0$ coboundary.
Let $o : CF(\mathcal F_{1}^{i}) \to CF(\mathcal F_{2}^{i})$
be a homomorphism such that $o_{\alpha_2\alpha_1}
= o(\alpha_1,\alpha_2)$ for $E(\alpha_2) - E(\alpha_1)
= E^{i+1} - c$ and $o_{\alpha_2\alpha_1} = 0$ otherwise.
\par
First by definition we find
\begin{equation}\label{form1771before}
\aligned
&(\psi_2^{i+1 i}\circ o)_{\alpha'_2\alpha_1} \\
& =
\left( \psi_2^{i+1 i}\circ (-\hat d_{2}^{i} \circ {\frak h}_{ab}^{i}
-{\frak h}_{ab}^{i} \circ \hat d_{1}^{i}
+ \frak n_b^i - \frak n_a^i)\right)_{\alpha'_2\alpha_1} \\
& =
\left( -\hat d_{2}^{i+1}\circ \psi_2^{i+1 i} \circ {\frak h}_{ab}^{i}
 \right.\\
&\left.\qquad -\psi_2^{i+1 i} \circ{\frak h}_{ab}^{i} \circ \hat d_{1}^{i}+ \psi_2^{i+1 i} \circ\frak n_b^i - \psi_2^{i+1 i} \circ\frak n_a^i)\right)_{\alpha'_2\alpha_1}.
\endaligned
\end{equation}
Here and hereafter $E(\alpha'_2) - E(\alpha_1)
= E^{i+1} - c$.
We used Situation \ref{situ1639} (1) in (\ref{form1771before}).
We observe
\begin{equation}\label{form1771after}
\aligned
&(-\hat d_{2}^{i+1}\circ \psi_2^{i+1 i} \circ {\frak h}_{ab}^{i}
-\psi_2^{i+1i} \circ{\frak h}_{ab}^{i} \circ \hat d_{1}^{i}) _{\alpha'_2\alpha_1}\\
& \equiv
(-\hat d_{2}^{i+1}\circ (\psi_2^{i+1 i} \circ {\frak h}_{ab}^{i})\vert_{E^i}
-(\psi_2^{i+1 i} \circ{\frak h}_{ab}^{i})\vert_{E^i} \circ \hat d_{1}^{i}) _{\alpha'_2\alpha_1}.
\endaligned
\end{equation}
Here we used the following equality
\begin{equation}\label{intform1773}
\aligned
&(\hat d_{2}^{i+1}\circ A - (-1)^{\deg A}A\circ \hat d_{1}^{i})_{\alpha'_2\alpha_1}\\
&\quad-
(\hat d_{2}^{i+1}\circ A\vert_{E^i} -
(-1)^{\deg A}A\vert_{E^i}\circ \hat d_{1}^{i})_{\alpha'_2\alpha_1} \\
&=
d_0 \circ A_{\alpha'_2\alpha_1}
- (-1)^{\deg A} A_{\alpha'_2\alpha_1} \circ d_0.
\endaligned
\end{equation}
By induction hypothesis we have
\begin{equation}\label{17715form}
\aligned
&-(\psi_2^{i+1 i} \circ \frak h_{ab}^i )\vert_{E^i}
+ (\frak h_{ab}^{i+1} \circ  \psi_1^{i+1 i})\vert_{E^i}
- (\frak h_{a}^{i+1 i})\vert_{E^i} + (\frak h_{b}^{i+1 i})\vert_{E^i}
\\
&=
\left(\hat d_{2}^{i+1} \circ \frak H_{ab}^{i+1 i}
-
\frak H_{ab}^{i+1 i} \circ \hat d_{1}^{i} \right)\vert_{E^i}.
\endaligned
\end{equation}
Therefore we obtain
\begin{equation}\label{form1772before}
\aligned
&(-\hat d_{2}^{i+1}\circ(\psi_2^{i+1 i} \circ {\frak h}_{ab}^{i})\vert_{E^i} )_{\alpha'_2\alpha_1}
\\
& =
\left(\hat d_{2}^{i+1}\circ
\left(-(\frak h_{ab}^{i+1} \circ  \psi_1^{i+1 i})\vert_{E^i}
+ (\frak h_{a}^{i+1 i})\vert_{E^i} - (\frak h_{b}^{i+1 i})\vert_{E^i} \right))\right)_{\alpha'_2\alpha_1}\\
&\quad +\left(\hat d_{2}^{i+1}\circ(\hat d_{2}^{i+1} \circ \frak H_{ab}^{i+1 i}
-
\frak H_{ab}^{i+1 i} \circ \hat d_{1}^{i})\vert_{E^i} \right)_{\alpha'_2\alpha_1}.
\endaligned
\end{equation}
By (\ref{17715form}) we have
\begin{equation}\label{form177266}
\aligned
&(-(\psi_2^{i+1 i}\circ {\frak h}_{ab}^{i})\vert_{E^i} \circ \hat d_{1}^{i})_{\alpha'_2\alpha_1}
\\
& =
\left(\left(-(\frak h_b^{i+1} \circ \psi_1^{ii+1})\vert_{E^i}
+ (\frak h_a^{i+1 i})\vert_{E^i}
- (\frak h_b^{i+1 i})\vert_{E^i}\right) \circ \hat d_{1}^{i}\right)_{\alpha'_2\alpha_1}
\\
&\quad +\left((\hat d_{2}^{i+1} \circ \frak H_{ab}^{i+1 i}
-
\frak H_{ab}^{i+1 i} \circ \hat d_{1}^{i})\vert_{E^i}\circ \hat d_{1}^{i}
\right)_{\alpha'_2\alpha_1}.
\endaligned
\end{equation}
It follows that
\begin{equation}\label{combine1775}
\aligned
&(\ref{form1772before}) + (\ref{form177266})
\\
& =
\left(\hat d_{2}^{i+1}\circ\left(
(-\frak h_{ab}^{i+1} \circ  \psi_1^{ii+1})\vert_{E^i}
+ (\frak h_{a}^{i+1 i})\vert_{E^i} - (\frak h_{b}^{i+1 i})\vert_{E^i} \right)\right)_{\alpha'_2\alpha_1}
\\
&\quad+\left(\left(-(\frak h_{ab}^{i+1} \circ \psi_1^{ii+1})\vert_{E^i}
+ (\frak h_a^{i+1 i})\vert_{E^i}
- (\frak h_b^{i+1 i})\vert_{E^i}\right) \circ \hat d_{1}^{i}\right)_{\alpha'_2\alpha_1}.
\endaligned
\end{equation}
Here we used
$\hat d_{1}^{i}\circ \hat d_{1}^{i} =
\hat d_{2}^{i}\circ \hat d_{2}^{i} = 0$
and (\ref{intform1773}) with $A = \hat d_{2}^{i+1} \circ \frak H_{ab}^{i+1 i}
-
\frak H_{ab}^{i+1 i} \circ \hat d_{1}^{i}$.
Using (\ref{intform1773}) again we have
\begin{equation}\label{form1778new}
\aligned
(\ref{combine1775})
=
&\left(\hat d_{2}^{i+1}\circ\left(
-\frak h_{ab}^{i+1} \circ  \psi_1^{i+1 i}
+ \frak h_{a}^{i+1 i} - \frak h_{b}^{i+1 i} \right)\right)_{\alpha'_2\alpha_1}
\\
&\quad+\left(\left(-\frak h_{ab}^{i+1} \circ \psi_1^{i+1 i}
+ \frak h_a^{i+1 i}
- \frak h_b^{i+1 i}\right) \circ \hat d_{1}^{i}\right)_{\alpha'_2\alpha_1}.
\endaligned
\end{equation}
By Situation \ref{situ1639} (4) we have
\begin{equation}\label{form1773rev}
(\psi_2^{i+1 i} \circ \frak n^i_k
-
\frak n^{i+1}_k\circ \psi_1^{i+1 i})_{\alpha'_2\alpha_1}
=
(\hat d_{2}^{i+1}\circ \frak h^{i+1 i}_k
+ \frak h^{i+1 i}_k \circ \hat d_{1}^{i+1})_{\alpha'_2\alpha_1}
\end{equation}
for $k=a,b$.
Thus (\ref{form1771before}),
(\ref{form1771after}), (\ref{combine1775}), (\ref{form1778new}), (\ref{form1773rev}) imply
\begin{equation}\label{form17777}
\aligned
&(\psi_2^{i+1 i}\circ o)_{\alpha'_2\alpha_1} \\
& \equiv
\left(
\hat d_{2}^{i+1} \circ (-\frak h_b^{i+1 i} + \frak h_a^{i+1 i} -
\frak h_{ab}^{i+1}\circ \psi_1^{i+1 i})\right.\\
&\qquad+
(-\frak h_b^{i+1 i} + \frak h_a^{i+1 i}
-\frak h_{ab}^{i+1}\circ \psi_1^{i+1 i})
\circ \hat d_{1}^{i}
\\
&\qquad+(\frak n_b^{i+1} \circ \psi_1^{i+1 i}
+ \hat d_{2}^{i+1} \circ \frak h_b^{i+1 i}
+ \frak h_b^{i+1 i} \circ \hat d_{1}^{i})\\
&\qquad + \left.
(-\frak n_a^{i+1} \circ \psi_1^{i+1 i}
- \hat d_{2}^{i+1} \circ \frak h_a^{i+1 i}
- \frak h_a^{i+1 i} \circ \hat d_{1}^{i})
\right)_{\alpha'_2\alpha_1}\\
& \equiv
\left(
-\hat d_{2}^{i+1} \circ \frak h_{ab}^{i+1}\circ \psi_1^{i+1 i}
- \frak h_{ab}^{i+1}\circ \psi_1^{i+1 i}
\circ \hat d_{1}^{i} \right.
\\
&\left.\qquad\qquad
+\frak n_b^{i+1} \circ \psi_1^{i+1 i}
- \frak n_a^{i+1} \circ \psi_1^{i+1 i}
\right)_{\alpha'_2\alpha_1}.
\endaligned
\end{equation}
By Situation \ref{situ1639} (3) we have
\begin{equation}\label{formla1777}
(\hat d_{2}^{i+1} \circ \frak h_{ab}^{i+1}
+
\frak h_{ab}^{i+1} \circ \hat d_{1}^{i+1})_{\alpha'_2\alpha_1}
=
(\frak n_b^{i+1} - \frak n_a^{i+1})_{\alpha'_2\alpha_1}.
\end{equation}
Moreover, by Situation \ref{situ1639} (1) we find
$$
(\hat d_{2}^{i+1} \circ \frak h_{ab}^{i+1}\circ \psi_1^{i+1 i})_{\alpha'_2\alpha_1}
=
(\hat d_{2}^{i+1} \circ \psi_1^{i+1 i} \circ \frak h_{ab}^{i+1})_{\alpha'_2\alpha_1}.
$$
Therefore the right hand side of (\ref{form17777}) vanishes.
Namely $(\psi_2^{i+1 i}\circ o)_{\alpha'_2\alpha_1}$
is a $d_0$-coboundary. Since $\psi_2^{i+1 i}$ induces an isomorphism
on $d_0$-cohomology in energy level $0$
(this follows from Definition \ref{defn1628} (6)), $o(\alpha_1,\alpha_2)$ is also a
$d_0$-coboundary.
\par
Thus we can find $(\frak h_{ab}^i)_{\alpha_2\alpha_1}$
and promote $\frak h_{ab}^i$ to the energy level $E^{i+1}$.
Note that we have a freedom to change $(\frak h_{ab}^i)_{\alpha_2\alpha_1}$
by $d_0$-cocycle.
\begin{lem}
\begin{equation}\label{tobeeliminatedybyH}
\aligned
\left(\hat d_{2}^{i+1} \circ \frak H_{ab}^{i+1 i}
-
\frak H_{ab}^{i+1 i} \circ \hat d_{1}^{i} \right.
&+ \frak h_b^{i+1 i} - \frak h_a^{i+1 i}
\\
&\left.- \frak h_{ab}^{i+1} \circ \psi_1^{i+1 i}
+\psi_2^{i+1 i} \circ \frak h_{ab}^{i}\right)_{\alpha'_2\alpha_1}
\endaligned
\end{equation}
is a $d_0$-cocycle for $E(\alpha'_2) - E(\alpha_1)
= E^{i+1} - c$.
\end{lem}
\begin{proof}
We first note that
\begin{equation}\label{form1779}
\aligned
&(\hat d_{2}^{i+1} \circ(\hat d_{2}^{i+1} \circ \frak H_{ab}^{i+1 i}
-
\frak H_{ab}^{i+1 i} \circ \hat d_{1}^{i})\\
&+
(\hat d_{2}^{i+1} \circ \frak H_{ab}^{i i+1}
-
\frak H_{ab}^{i+1 i} \circ \hat d_{1}^{i})\circ \hat d_{1}^{i})_{\alpha'_2\alpha_1}
= 0.
\endaligned
\end{equation}
Next we observe that Situation \ref{situ1639} (4) implies
\begin{equation}
\aligned
&(\hat d_{2}^{i+1} \circ( \frak h_b^{i i+1} -\frak h_a^{i+1 i})
+ (\frak h_b^{i+1 i} - \frak h_a^{i i+1})\circ \hat d_{1}^{i})_{\alpha'_2\alpha_1}\\
&= (\frak n_b^{i+1} \circ \psi_1^{i+1 i} -  \psi_2^{ii+1}\circ \frak n_b^{i}
-\frak n_a^{i+1} \circ \psi_1^{i+1 i} + \psi_2^{ii+1}\circ \frak n_a^{i}
)_{\alpha'_2\alpha_1}.
\endaligned
\end{equation}
Moreover Situation \ref{situ1639} (1)(3) imply
\begin{equation}
\aligned
&(\hat d_{2}^{i+1} \circ\frak h_{ab}^{i+1} \circ \psi_1^{i+1 i}
+
\frak h_{ab}^{i+1} \circ \psi_1^{ii+1} \circ\hat d_{1}^{i}
)_{\alpha'_2\alpha_1} \\
&=
(\hat d_{2}^{i+1} \circ\frak h_{ab}^{i+1} \circ \psi_1^{i+1 i}
+
\frak h_{ab}^{i+1} \circ\hat d_{1}^{i+1} \circ\psi_1^{i+1 i}
)_{\alpha'_2\alpha_1}
\\
&=
((\frak n_b^{i+1} - \frak n_a^{i+1})\circ\psi_1^{i+1 i}
)_{\alpha'_2\alpha_1}.
\endaligned
\end{equation}
In a similar way we have
\begin{equation}\label{form1783}
\aligned
&(\hat d_{2}^{i+1} \circ\psi_2^{i+1 i} \circ \frak h_{ab}^{i}
+
\psi_2^{i+1 i} \circ \frak h_{ab}^{i} \circ\hat d_{1}^{i}
)_{\alpha'_2\alpha_1} \\
&=
(\psi_2^{i+1 i} \circ \hat d_{2}^{i}\circ \frak h_{ab}^{i}
+
\psi_2^{i+1 i} \circ \frak h_{ab}^{i} \circ\hat d_{1}^{i}
)_{\alpha'_2\alpha_1}
\\
&=
(\psi_2^{i+1 i} \circ(\frak n_b^{i} - \frak n_b^{i})
)_{\alpha'_2\alpha_1}.
\endaligned
\end{equation}
The sum of (\ref{form1779})-(\ref{form1783}) is $0$.
\par
On the other hand,
$$
\hat d_{2}^{i+1} \circ \frak H_{ab}^{i+1 i}
-
\frak H_{ab}^{i+1 i} \circ \hat d_{1}^{i}
+ \frak h_b^{i+1 i} - \frak h_a^{i i+1}
- \frak h_{ab}^{i+1} \circ \psi_1^{i+1 i}
+\psi_2^{i+1 i} \circ \frak h_{ab}^{i}
$$
is zero up to energy level $E^i$ by assumption.
Therefore the sum of (\ref{form1779})-(\ref{form1783})
is $d_0$ differential of (\ref{tobeeliminatedybyH}).
\end{proof}
Therefore we can choose
$(\frak h_{ab}^i)_{\alpha_2\alpha_1}$ so that
(\ref{tobeeliminatedybyH}) is a $d_0$-coboundary for
$\alpha_1,\alpha_2$ with $E(\alpha_2) - E(\alpha_1) = E^{i+1} - c$.
We now choose $(\frak H_{ab}^{i+1 i})_{\alpha'_2\alpha_1}$
which bounds (\ref{tobeeliminatedybyH}).
The proof of Proposition \ref{prop1639} is now complete.
\end{proof}
\begin{rem}
In Situation \ref{situ1639} (1) we assumed that the energy zero
part of the partial cochain map of energy loss $0$ is the identity map.
Actually we only need a milder assumption that
the energy $0$ part of  the partial cochain map of energy loss $0$
induces an injection on $d_0$ cohomology.
\end{rem}
\begin{proof}
[Proof of Theorem \ref{indlinesysmainth1} (4)]
We are given $\frak N\frak N_a$, $\frak N\frak N_b$,
that are morphisms of inductive systems $\mathcal F\mathcal F_1 \to
\mathcal F\mathcal F_2$.
Each of them consists of partial morphisms
$$
\frak N^{i}_k :
\mathcal F_{1}^{i} \to \mathcal F_{2}^{i}, \quad k=a,b,
$$
respectively.
They induce partial cochain maps
$$
\frak n^{i}_k :CF(\mathcal F_{1}^{i}) \to CF(\mathcal F_{2}^{i}), \quad k=a,b
$$
for each $i$.
We are also given a homotopy $\frak H\frak H$ from
$\frak N\frak N_a$ to $\frak N\frak N_b$.
It consists of partial homotopies from $\frak N^{i}_a$ to $\frak N^{i}_b$.
By Proposition \ref{homotopyexistpertstatement}
and Lemma \ref{lem165rev2},
it induces a partial homotopy
$$
\frak h_{ab}^i : CF(\mathcal F_{1}^{i}) \to CF(\mathcal F_{2}^{i})
$$
from $\frak n^{i}_a$ to $\frak n^{i}_b$.
\par
Since $\frak N\frak N_a$, $\frak N\frak N_b$ are morphisms of
partial linear K-system,
we are given a partial homotopy from
$\frak N_2^{i+1 i} \circ \frak N^{i}_k$ to $\frak N^{i+1}_k \circ \frak N_1^{i+1 i}$.
Again by Proposition  \ref{homotopyexistpertstatement}
and Lemma \ref{lem165rev2},
it induces a partial homotopy
$$
\frak h^{i+1 i}_k : CF(\mathcal F_{1}^{i}) \to CF(\mathcal F_{2}^{i+1})
$$
from $\psi_2^{i+1 i} \circ \frak n^{i}_k$ to $\frak n^{i+1}_k \circ \psi_1^{i+1 i}$.
(Note $\psi_j^{i+1 i} : CF(\mathcal F_{j}^{i}) \to CF(\mathcal F_{j}^{i+1})$
is a cochain map induced by $\frak N_j^{i+1 i}$.)
Thus we are in Situation \ref{situ1639}.
\par
Now the existence of homotopy of homotopy $\mathcal H^i$
which is Definition \ref{defn1528} (4) (b), (c),
together with Proposition \ref{prop1618rev} and
Lemma \ref{lem165rev2}
implies that there exists
$$
\frak H_{ab}^{i+1 i} : CF(\mathcal F_{1}^{i}) \to CF(\mathcal F_{2}^{i+1})
$$
that satisfies (\ref{diagram1663}).
Thus we can apply Proposition \ref{prop1639}
to promote our partial homotopy $\frak h_{ab}^i$ to one
of energy cut level $\infty$.
Therefore $\frak n^{i}_a$ is cochain homotopic to $\frak n^{i}_b$.
\end{proof}
\begin{proof}[Proof of Theorem \ref{linesysmainth2} (2) ]
This is a special case of Theorem \ref{indlinesysmainth1} (4)
where $\frak N^{i+1 i}$ is the identity morphism.
\end{proof}
\begin{proof}[Proof of
Theorem \ref{indlinesysmainth1} (5) and Theorem \ref{linesysmainth2} (3)]
We will prove Theorem \ref{indlinesysmainth1} (5).
Theorem \ref{linesysmainth2} (3) is its special case.
\par
We recall that $\frak N_{ba} : \mathcal F\mathcal F_a \to \mathcal F\mathcal F_b$ consists of
partial morphisms
$$
\frak N_{ba}^{i} : \mathcal F_{a}^{i} \to \mathcal F_{b}^{i}
$$
and partial homotopies $\frak H_{ba}^{i}$ between
$\frak N_b^{i+1 i} \circ \frak N_{ba}^{i}$ and
$\frak N_{ba}^{i+1}  \circ \frak N_a^{i+1 i}$.
(See Diagram \ref{dd1766} below.)
\begin{equation}\label{dd1766}
\begin{CD}
\mathcal F_{b}^{i}
@ > {\frak N_b^{i+1 i}} >> \mathcal F_{b}^{i+1}
\\
@ A{\frak N_{ba}^{i}}AA @ AA{\frak N_{ba}^{i+1}}A \\
\mathcal F_{a}^{i}
@ > {\frak N_a^{i+1 i}} >> \mathcal F_{a}^{i+1}\end{CD}
\end{equation}
Also $\frak N_{cb} : \mathcal F\mathcal F_b \to \mathcal F\mathcal F_c$ consists of
partial morphisms $\frak N_{cb}^{i} : \mathcal F_{b}^{i} \to \mathcal F_{c}^{i}$
and partial homotopies $\frak H_{cb}^{i}$ in a similar way.
The definition of the composition $\frak N_{ca} = \frak N_{cb} \circ \frak N_{ba}$ is given in Lemma-Definition \ref{lemdef1437}.
\par
Let $\mathcal N_{ab}^{i}(\alpha_a,\alpha_b)$ and $\mathcal N_{bc}^{i}(\alpha_b,\alpha_c)$  be interpolation spaces of
$\frak N_{ba}^{i}$ and $\frak N_{cb}^{i}$, respectively.
For $k=a,b,c$, let $\mathcal N^{i+1 i}_k (\alpha_k,\alpha'_k)$ be an interpolation space of
$\frak N^{i+1 i}_k$.
We denote by
$\mathcal N_{ab}^{i}(\alpha_a,\alpha_b;[0,1])$ and $\mathcal N_{bc}^{i}(\alpha_b,\alpha_c;[0,1])$
interpolation spaces of
$\frak H_{ba}^{i}$ and $\frak H_{cb}^{i}$, respectively.
By definition we have
\begin{equation}
\aligned
&\partial \mathcal N_{ab}^{i}(\alpha_a,\alpha_b;[0,1]) \\
=&
\bigcup_{\alpha'_a}
\mathcal N_{a}^{i+1 i}(\alpha_a,\alpha'_a)
\times^{\boxplus\tau}_{R_{\alpha'_a}}
\mathcal N_{ab}^{i}(\alpha'_a,\alpha_b)
\\
&\cup
\bigcup_{\alpha'_b}
\mathcal N_{ab}^{i}(\alpha_a,\alpha'_b)
\times^{\boxplus\tau}_{R_{\alpha'_b}}
\mathcal N_{b}^{i+1 i}(\alpha'_b,\alpha_b).
\endaligned
\end{equation}
Here $\times^{\boxplus\tau}_{R_{\alpha'_a}}$
and
$\times^{\boxplus\tau}_{R_{\alpha'_b}}$
are as in Definition \ref{defn1635}.
Similarly we have
\begin{equation}
\aligned
&\partial \mathcal N_{bc}^{i}(\alpha_b,\alpha_c;[0,1]) \\
=&
\bigcup_{\alpha'_b}
\mathcal N_{b}^{i+1 i}(\alpha_b,\alpha'_b)
\times^{\boxplus\tau}_{R_{\alpha'_b}}
\mathcal N_{bc}^{i}(\alpha'_b,\alpha_c)
\\
&\cup
\bigcup_{\alpha'_c}
\mathcal N_{bc}^{i}(\alpha_a,\alpha'_c)
\times^{\boxplus\tau}_{R_{\alpha'_c}}
\mathcal N_{c}^{i+1 i}(\alpha'_c,\alpha_c).
\endaligned
\end{equation}
The composition $\frak N_{ca} = \frak N_{cb} \circ \frak N_{ba}$
consists of $\frak N_{ca}^{i}$ and $\frak H_{ca}^{i}$.
Here the interpolation space
$\mathcal N_{ac}^{i}(\alpha_a,\alpha_c)$ of $\frak N_{ca}$
is given by
$$
\mathcal N_{ac}^{i}(\alpha_a,\alpha_c)
=
\bigcup_{\alpha_b}
\mathcal N_{ab}^{i}(\alpha_a,\alpha_b)
\times^{\boxplus\tau}_{R_{\alpha_b}}
\mathcal N_{bc}^{i}(\alpha_b,\alpha_c).
$$
Note that the union in the above formula
is different from the disjoint union and is defined
as in Lemma-Definition \ref{1638defken}.
\par
The homotopy $\frak H_{ca}^{i}$ is obtained
by gluing $\frak H_{cb}^{i} \circ \frak N_{ba}^{i}$
and $\frak N_{cb}^{i} \circ \frak H_{ba}^{i}$
as follows.
The interpolation space of
$\frak H_{cb}^{i} \circ \frak N_{ba}^{i}$ is
\begin{equation}\label{form1769}
\bigcup_{\alpha_b}
\mathcal N_{ab}^{i}(\alpha_a,\alpha_b;[0,1])
\times^{\boxplus\tau}_{R_{\alpha_b}}
\mathcal N_{bc}^{i}(\alpha_b,\alpha_c).
\end{equation}
The interpolation space of
$\frak N_{cb}^{i} \circ \frak H_{ba}^{i}$  is
\begin{equation}\label{form1770}
\bigcup_{\alpha_b}
\mathcal N_{ab}^{i}(\alpha_a,\alpha_b)
\times^{\boxplus\tau}_{R_{\alpha_b}}
\mathcal N_{bc}^{i}(\alpha_b,\alpha_c;[0,1]).
\end{equation}
We observe that both (\ref{form1769}) and
(\ref{form1770}) contain
\begin{equation}\label{form1771}
\bigcup_{\alpha_b} \bigcup_{\alpha'_b}
\mathcal N_{ab}^{i}(\alpha_a,\alpha_b)
\times^{\boxplus\tau}_{R_{\alpha_b}}
\mathcal N_{b}^{i+1 i}(\alpha_b,\alpha'_b)
\times^{\boxplus\tau}_{R_{\alpha'_b}}
\mathcal N_{bc}^{i}(\alpha_b,\alpha_c)
\end{equation}
in its boundary. We smooth the corners
contained in (\ref{form1771}).
See Section \ref{section:compomorphis}.
\par
Next we recall that while we constructed a partial cochain map
$$
\psi_{ba}^{i} : CF_{a}^{i} \to CF_{b}^{i}
$$
(see \eqref{CFjinoichi} for the notation
$CF_{\ast}^{i}$),
we took a thickening
$\widehat{\mathcal U^{i,+}_{ab}}(\alpha_a,\alpha_b)$
of
$\widehat{\mathcal U^{i}_{ab}}(\alpha_a,\alpha_b)$
(note they are Kuranishi structures of
$\mathcal N_{ab}^{i}(\alpha_a,\alpha_b)^{\frak C^h\boxplus\tau}$)
and a CF-perturbation
$\widehat{{\frak S}_{ab}^{i}}(\alpha_a,\alpha_b)$ of $\widehat{\mathcal U^{i,+}_{ab}}(\alpha_a,\alpha_b)$.
(Here $\frak C^h$ stands for the horizontal boundary as in Definition \ref{defn1726}.)
During the construction of a partial cochain map
$$\psi_{cb}^{i} : CF_{b}^{i} \to CF_{c}^{i},$$
we took a thickening
$\widehat{\mathcal U^{i,+}_{bc}}(\alpha_b,\alpha_c)$
(of $\widehat{\mathcal U^{i}_{bc}}(\alpha_b,\alpha_c)$) and
its CF-perturbation
$\widehat{{\frak S}_{bc}^{i}}(\alpha_b,\alpha_c)$.
During the construction of partial cochain maps
$$
\psi_{k}^{i+1 i} : CF_{k}^{i} \to CF_{k}^{i+1}, \quad k=a,b,c,
$$
we took thickenings
$\widehat{\mathcal U^{i+1 i, +}_{k}}(\alpha_k,\alpha'_k)$
of $\widehat{\mathcal U^{i+1 i}_{k}}(\alpha_k,\alpha'_k)$
which are Kuranishi structures on
$\mathcal N_{k}^{i+1 i}(\alpha_k,\alpha'_k)^{\boxplus\tau}$),
and their CF-perturbations
$\widehat{{\frak S}_{k}^{i+1 i}}(\alpha_k,\alpha'_k)$.
\par
Furthermore, in the course of our construction of a partial cochain homotopy
$$
\frak h_{ba}^{i}
\quad \text {between $\frak n_{b}^{i+1 i} \circ \psi_{ba}^{i}$
and $\psi_{ba}^{i+1}  \circ \frak n_{a}^{i+1 i}$}
$$
(where $\frak n_{b}^{i+1 i} \circ \psi_{ba}^{i}$
and $\psi_{ba}^{i+1}  \circ \frak n_{a}^{i+1 i}$
are cochain maps
$: CF_a^i \to CF_b^{i+1}$),
we took a thickening
$\widehat{\mathcal U_{ab}^{i,+}}(\alpha_a,\alpha_b;[0,1])$
of  $\widehat{\mathcal U_{ab}^{i}}(\alpha_a,\alpha_b;[0,1])$
which is a Kuranishi structure on
$\mathcal N_{ab}^{i}(\alpha_a,\alpha_b;[0,1])^{\frak C^h\boxplus\tau}$, and
its CF-perturbation
$\widehat{{\frak S}_{ab}^{i}}(\alpha_a,\alpha_b;[0,1])$.
\begin{equation}\label{dd1766}
\begin{CD}
CF_{a}^{i}
@ > {\psi_{ba}^{i}} >> \mathcal CF_{b}^{i}
@ > {\psi_{cb}^{i}} >> \mathcal CF_{c}^{i}
\\
@ V{\frak n_{a}^{i+1 i}}VV  @ V{\frak n_{b}^{i+1 i}}VV
@ V{\frak n_{c}^{i+1 i}}VV\\
CF_{a}^{i+1}
@ > {\psi_{ba}^{i+1}} >> \mathcal CF_{b}^{i+1}
@ > {\psi_{cb}^{i+1}} >> \mathcal CF_{c}^{i+1}
\end{CD}
\end{equation}
During our construction of a partial cochain homotopy
$$
\frak h_{cb}^{i}
\quad \text {between $\frak n_{c}^{i+1 i} \circ \psi_{cb}^{i}$
and $\psi_{cb}^{i+1}  \circ \frak n_{b}^{i+1 i}$,}
$$
we took a thickening
$\widehat{\mathcal U^{i,+}_{bc}}(\alpha_b,\alpha_c;[0,1])$
of $\widehat{\mathcal U^{i}_{bc}}(\alpha_b,\alpha_c;[0,1])$
which is a Kuranishi structure on
$\mathcal N_{bc}^{i}(\alpha_b,\alpha_c;[0,1])^{\frak C^h\boxplus\tau}$)
and its CF-perturbation
$\widehat{{\frak S}_{bc}^{i}}(\alpha_b,\alpha_c;[0,1])$.
\par
Note that by Lemma \ref{lem1721} we may use
\begin{equation}\label{form1772}
\bigcup_{\alpha_b}
\widehat{\mathcal U^{i,+}_{ab}}(\alpha_a,\alpha_b)
\times_{R_{\alpha_b}}
\widehat{\mathcal U^{i,+}_{bc}}(\alpha_b,\alpha_c)
\end{equation}
and
\begin{equation}\label{form1773}
\bigcup_{\alpha_b}
\widehat{{\frak S}_{ab}^{i}}(\alpha_a,\alpha_b)
\times_{R_{\alpha_b}}
\widehat{{\frak S}_{bc}^{i}}(\alpha_b,\alpha_c)
\end{equation}
to define a partial cochain map
$$
\psi_{ca}^{i} : CF_{a}^{i} \to CF_{c}^{i}.
$$
Here (\ref{form1772})
is a thickening of
$
\bigcup_{\alpha_b}
\widehat{\mathcal U^{i}_{ab}}(\alpha_a,\alpha_b)
\times_{R_{\alpha_b}}
\widehat{\mathcal U^{i}_{bc}}(\alpha_b,\alpha_c)
$
which is a Kuranishi structure on
$\mathcal N_{ac}^{i}(\alpha_a,\alpha_c)^{\frak C^h\boxplus\tau}$,
and (\ref{form1773}) is its CF-perturbation.
Therefore composition formula and Lemma \ref{lem1634lemlem} yield
\begin{equation}\label{form1774}
\psi_{ca}^{i}
= \psi_{cb}^{i} \circ \psi_{ba}^{i},
\end{equation}
if we define $\psi_{ca}^{i}$ by this particular choice.
(Lemma \ref{17222lem}.)
\par
Next we consider $\frak h_{ca}^{i}$.
We take
\begin{equation}\label{form1775}
\aligned
&\left(\bigcup_{\alpha_b}
\widehat{\mathcal U^{i,+}_{ab}}(\alpha_a,\alpha_b;[0,1])
\times_{R_{\alpha_b}}
\widehat{\mathcal U^{i,+}_{bc}}(\alpha_b,\alpha_c)
\right)
\\
&\cup
\left(\bigcup_{\alpha_b}
\widehat{\mathcal U^{i,+}_{ab}}(\alpha_a,\alpha_b)
\times_{R_{\alpha_b}}
\widehat{\mathcal U^{i,+}_{bc}}(\alpha_b,\alpha_c;[0,1])
\right).
\endaligned
\end{equation}
Here we take a partial smoothing of corner of the right hand side
and glue them.
Then (\ref{form1775}) is a thickening of
$\mathcal N_{ac}^{i}(\alpha_a,\alpha_c)^{\frak C^h\boxplus\tau}$ and
\begin{equation}\label{form1776}
\aligned
&\left(\bigcup_{\alpha_b}
\widehat{{\frak S}_{ab}^{i}}(\alpha_a,\alpha_b;[0,1])
\times_{R_{\alpha_b}}
\widehat{{\frak S}_{bc}^{i}}(\alpha_b,\alpha_c)
\right)
\\
&\cup
\left(\bigcup_{\alpha_b}
\widehat{{\frak S}_{ab}^{i}}(\alpha_a,\alpha_b)
\times_{R_{\alpha_b}}
\widehat{{\frak S}_{bc}^{i}}(\alpha_b,\alpha_c;[0,1])
\right)
\endaligned
\end{equation}
is a CF-perturbation of (\ref{form1775}).
We use (\ref{form1775})  and (\ref{form1776})
to define $\frak h_{ca}^{i}$.
Then by composition formula and Lemma \ref{lem1634lemlem} again we find
\begin{equation}\label{form1777}
\frak h_{ca}^{i}
=
\psi_{cb}^{i+1} \circ \frak h_{ba}^{i}
+
\frak h_{cb}^{i} \circ \psi_{ba}^{i}.
\end{equation}
\begin{lem}
\begin{equation}\label{form1778}
\hat d \circ \frak h_{ca}^{i} + \frak h_{ca}^{i} \circ \hat d
=
\psi_{ca}^{i+1} \circ \frak n_{a}^{i+1 i}
-
\frak n_{c}^{i+1 i} \circ \psi_{ca}^{i}.
\end{equation}
\end{lem}
\begin{proof}
$$
\aligned
&\hat d \circ \frak h_{ca}^{i} + \frak h_{ca}^{i} \circ \hat d \\
=
&\psi_{cb}^{i+1}  \circ \hat d \circ \frak h_{ba}^{i} + \psi_{cb}^{i+1}  \circ
\frak h_{ba}^{i} \circ \hat d \\
&+
\frak h_{cb}^{i} \circ \hat d  \circ \psi_{ba}^{i}
+
\hat d \circ \frak h_{cb}^{i} \circ \psi_{ba}^{i}
\\
= &
\psi_{cb}^{i+1} \circ \psi_{ba}^{i+1} \circ  \frak n_{a}^{i+1 i}
-
\psi_{cb}^{i+1} \circ \frak n_{b}^{i+1 i} \circ \psi_{ba}^{i}
\\
&+
\psi_{cb}^{i+1} \circ  \frak n_{b}^{i+1 i} \circ \psi_{ba}^{i}
-
\frak n_{c}^{i+1 i} \circ \psi_{cb}^{i+1} \circ \psi_{ba}^{i}
\\
= &
\psi_{ca}^{i+1} \circ \frak n_{a}^{i+1 i}
-
\frak n_{c}^{i+1 i} \circ \psi_{ca}^{i}.
\endaligned
$$
\end{proof}
We note that the equalities
(\ref{form1774}), (\ref{form1777}), (\ref{form1778})
are ones of energy cut level $E^i$.
We recall that while we constructed
$\psi_{ba}$ and $\psi_{cb}$ we promoted
$\frak n_{k}^{i+1 i}$ ($k =a,b,c$),
$\psi_{ba}^{i}$, $\psi_{cb}^{i}$,
$\frak h_{ba}^{i}$ and $\frak h_{cb}^{i}$
to energy cut level $\infty$.
We now promote $\psi_{ca}^{i}$ and $\frak h_{ca}^{i}$
to energy cut level $\infty$
so that (\ref{form1774}), (\ref{form1777})
hold as the equalities at the energy cut level $\infty$.
Then (\ref{form1778})
holds as an equality at the energy cut level $\infty$.
\par
Thus for this particular choice of promotion, the equality
$\psi_{ca}^{i}
= \psi_{cb}^{i} \circ \psi_{ba}^{i}$ holds
not only up to cochain homotopy but also
as a strict identity.
Since $\psi_{ca}^{i} $ is independent of various choices
up to cochain homotopy (Theorem \ref{indlinesysmainth1},
which we proved in Subsection \ref{subsec:proofsec14main2}),
$
\psi_{ca}^{i}
= \psi_{cb}^{i} \circ \psi_{ba}^{i}$ holds
for any choice of $\psi_{ca}^{i}$ up to cochain homotopy.
The proof of Theorem \ref{indlinesysmainth1} (5)
is complete.
\end{proof}
We next prove a partial cochain map version of Theorem \ref{linesysmainth2} (4).
\begin{lem}
In the situation of Theorem \ref{linesysmainth2} (4) we have
$\mathcal{ID}_* \sim {\rm id} \mod T^E$ for any $E$.
Here $\sim$ means cochain homotopic.
\end{lem}
\begin{proof}
This is a consequence of Theorem \ref{linesysmainth2} (2)
and Proposition \ref{properidentity}.
In fact, Theorem \ref{linesysmainth2} (2) implies that
$\mathcal{ID}_* \circ \mathcal{ID}_* \sim \mathcal{ID}_*$.
(Here ${\frak ID}_*$ is the cochain map induced by the identity
morphism and $\sim$ is a chan homotopy).
On the other hand, Lemma \ref{lema1631} below implies
that ${\frak ID}_*$ is an isomorphism.
Therefore ${\frak ID}_* \sim {\rm id}$.
\end{proof}
\begin{lem}\label{lema1631}
Let $\widehat\psi_{21} = \{\psi_{\alpha_2\alpha_1}\} :
\mathcal F_1 \to \mathcal F_2$
be a partial cochain map of energy cut level $E_0$ and
energy loss $0$.
We assume it is congruent to the isomorphism
in the sense of Definition \ref{defn1521}. Then
there exists a partial cochain map $\widehat\psi_{12} = \{(\psi_{12})_{\alpha_1\alpha_2}\} : \mathcal F_2 \to \mathcal F_1$
of energy cut level $E_0$ and
energy loss $0$ such that $\widehat\psi_{12} \circ \widehat\psi_{21}$
and $\widehat\psi_{21} \circ \widehat\psi_{12}$ are identity maps.
\end{lem}
\begin{proof}
We construct $(\psi_{12})_{\alpha_1\alpha_2}$ by induction on
$E(\alpha_1) - E(\alpha_2)$.
This induction is possible because the set of values of $E(\alpha_1) - E(\alpha_2)$ is a discrete set
by uniform Gromov compactness Definition \ref{defn1528} (2)(g).
\par
We start with the case when $E(\alpha_1) - E(\alpha_2) = 0$.
By definition of partial cochain map
of energy loss $0$ congruent to the isomorphism
(Definition \ref{defn1521}), we have
$$
(\psi_{12})_{\alpha_1\alpha_2}
=
\begin{cases}
0  &\text{if $\alpha_1 \ne \alpha_2$} \\
{\rm id} &\text{if $\alpha_1 = \alpha_2$}.
\end{cases}
$$
In fact, the interpolation space $\mathcal N(\alpha_2,\alpha_1)$
is an empty set if $\alpha_1 \ne \alpha_2$ and $E(\alpha_1) - E(\alpha_2) = 0$
by Condition \ref{morphilinsys} (V).
This implies the first equality.
If $\alpha_1 = \alpha_2 = \alpha$,
we have $\mathcal N(\alpha,\alpha) = R_{\alpha}$
by Definition \ref{defn1521}.
Moreover, ${\rm ev}_{\pm} : \mathcal N(\alpha,\alpha) \to R_{\alpha}$
is the identity map.
This implies the second equality.
\par
Suppose
$E(\alpha_1) - E(\alpha_2) = E_0$ and
we have defined $(\psi_{12})_{\alpha'_1\alpha'_2}$
for $E(\alpha_1) - E(\alpha_2) < E_0$.
Then the condition that $\widehat\psi_{21}$ is a right
inverse of $\widehat\psi_{12}$ at the energy cut level $E_0$ is
written as
$$
(\psi_{12})_{\alpha_1\alpha_2}
+
\sum_{\alpha'_2} (\psi_{12})_{\alpha_1\alpha'_2}
\circ (\psi_{21})_{\alpha'_2\alpha_1}
=0.
$$
Since the second term is already defined
we can find $(\psi_{12})_{\alpha_1\alpha_2}$ uniquely
so that this condition holds.
Thus we have found the left inverse by induction.
We can find the right inverse in the same way.
A standard fact in group theory yields that the right and left inverse coincide.
It is also easy to see that partial inverse of a partial cochain map
is a partial cochain map.
\end{proof}

\begin{proof}[Proof of Theorem \ref{linesysmainth1} (2)(f) and
Theorem \ref{indlinesysmainth1} (2)(e)]
Theorem \ref{indlinesysmainth1} (2)(e)
follows from Theorem \ref{indlinesysmainth1} (4) applied
to the identity morphism.
Theorem \ref{linesysmainth1} (2)(f) is a special case
of Theorem \ref{indlinesysmainth1} (2)(e) when
$\frak N^{i+1 i}$ are the identity morphisms.
\end{proof}

\begin{proof}[Proof of Theorem \ref{indlinesysmainth1} (6) and Theorem \ref{linesysmainth2} (4)]
We will prove Theorem \ref{indlinesysmainth1} (6).
Theorem \ref{linesysmainth2} (4) is its special case.
\par
We first define the identity morphism of an inductive system of
linear K-systems.
Let $\mathcal{FF} = (\{\mathcal F^i\},\{\frak N^i\})$
be an inductive system of
linear K-systems.
(Definition \ref{defn1528} (2).)
We put $\mathcal{FF}_k = \mathcal{FF}$ for $k=a,b$
and $\frak N^i_{ba} = \mathcal{ID}_{\mathcal F^i}$ the
identity morphism of $\mathcal F^i$.
By Proposition \ref{properidentity} we have
$$
\frak N^i \circ \mathcal {ID}_{\mathcal F^i}
\sim \frak N^i
\sim \mathcal {ID}_{\mathcal F^i} \circ \frak N^i.
$$
Let $\frak H_{ba}^i$ be this homotopy.
(We take the particular choice of the homotopy which we gave
during the proof of Proposition \ref{properidentity}.)
\begin{defn}
We define $\mathcal{ID}_{\mathcal{FF}} = (\{\mathcal{ID}_{\mathcal F^i}\},\{
\frak H_{ba}^i\})$
the {\it identity morphism}
\index{identity morphism ! of inductive system $\mathcal{FF}$}
\index{K-system ! identity morphism of inductive system $\mathcal{FF}$}
from $\mathcal{FF}$ to itself.
\end{defn}
Theorem \ref{indlinesysmainth1} (6)
claims that the cochain map induced by $\mathcal{ID}_{\mathcal{FF}}$
is cochain homotopic to the identity.
To prove this it suffices to show the following lemma.
\begin{lem}\label{lem1745}
If $\frak N_{cb} : \mathcal{FF}_b \to \mathcal{FF}_c$
be a morphism of inductive systems of linear K-systems, then the composition
$\frak N_{cb} \circ \mathcal{ID}_{\mathcal{FF}}$ is homotopic
to $\frak N_{cb}$.
The same holds for $\mathcal{ID}_{\mathcal{FF}} \circ \frak N_{cb}$.
\end{lem}
\begin{proof}
Write $\mathcal{FF}_c = (\{\mathcal F_c^i\},\{\frak N_c^{i+1 i}\})$
and
$\frak N_{cb}=
(\{\frak N_{cb}^{i}\},\{\frak H_{cb}^i\}) : \mathcal{FF}_b \to \mathcal{FF}_c$.
Here $\frak N_{k}^{i+1 i} : \mathcal F_k^i \to \mathcal F_k^{i+1}$
($k=b,c$),
$\frak N_{cb}^{i} : \mathcal F_b^i \to \mathcal F_c^i$
are partial morphisms of partial linear K-systems, and
$\frak H_{cb}^i$ is a homotopy between
$\frak N_{cb}^{i+1}\circ \frak N_{b}^{i+1 i}$
and $\frak N_{c}^{i+1 i}\circ \frak N_{cb}^{i}$.
(They are partial morphisms $: \mathcal F_b^i \to \mathcal F_c^{i+1}$.)
\par
We denote by $\mathcal M(ki;\alpha_-,\alpha_+)$
the moduli space of connecting orbits for
$\mathcal F_k^i$. $(k=a,b,c)$.
Note that
$$
\mathcal M(ai;\alpha_-,\alpha_+) = \mathcal M(bi;\alpha_-,\alpha_+).
$$
Let $\mathcal N(k,i i+1;\alpha_-,\alpha_+)$
and $\mathcal N(bc,i;\alpha_-,\alpha'_+)$ be
interpolation spaces of $\frak N^{i+1 i}_k$ and
$\frak N^{i}_{cb}$, respectively.
Let $\mathcal N(bc,i i+1;\alpha_-,\alpha'_+;[1,2])$
be the interpolation space of $\frak H_{cb}^i$.
\begin{equation}\label{3tutakenozusiki}
\begin{CD}
\mathcal F_c^i
@ >>{\frak N^{i+1 i}_c}>\mathcal F_c^{i+1} \\
@A{\frak N_{cb}^{i}}AA   @AA{\frak N_{cb}^{i+1}}A  \\
\mathcal F_b^i
@ >>{\frak N^{i+1 i}_b}>\mathcal F_b^{i+1}
\\
@A{\mathcal{ID}}AA   @AA\mathcal{ID}A   \\
\mathcal F_a^i
@ >>{\frak N^{i+1 i}_a}>\mathcal F_c^{i+1}
\end{CD}
\end{equation}
Here we also note that $\mathcal F_a^i = \mathcal F_b^i$
and $\frak N^{i+1 i}_a = \frak N^{i+1 i}_b$.
\par
We put
$\frak N^i \circ \mathcal{ID}_{\mathcal F^i}
= (\{\frak N_{cb}^{i}\circ \mathcal{ID}\},\{\frak H_{ca}^i\})$.
By definition the homotopy $\frak H_{ca}^i$
is obtained as
\begin{equation}
\frak H_{ca}^i =
(\frak N_{cb}^{i+1}\circ\frak H_{ba}^i) \cup
(\frak H_{cb}^i\circ \mathcal{ID}).
\end{equation}
Note that
$\frak N_{cb}^{i+1}\circ\frak H_{ba}^i$
is a homotopy from $\frak N_{cb}^{i+1}\circ \mathcal{ID}
\circ \frak N^{i+1 i}_a$
to
$\frak N_{cb}^{i+1}\circ \frak N^{i+1 i}_b
\circ \mathcal{ID}$
and
$\frak H_{cb}^i\circ \mathcal{ID}$
is a homotopy from
$\frak N_{cb}^{i+1}\circ \frak N^{i+1 i}_b
\circ \mathcal{ID}$
to
$\frak N^{i+1 i}_c\circ \frak N_{cb}^{i}
\circ \mathcal{ID}$.
(See Diagram (\ref{3tutakenozusiki}).)
\par
Now we start the construction of the homotopy between
$\frak N_{cb} \circ \mathcal{ID}_{\mathcal{FF}}$ and
$\frak N_{cb}$.
By Proposition \ref{properidentity} we have
$\frak N_{cb}^{i} \circ \mathcal{ID} \sim \frak N_{cb}^{i}$.
Let $\frak H^i$ be this homotopy.
Note that
as $\frak H^i$ we take the specific homotopy we constructed during the proof of
Proposition \ref{properidentity}.
To prove Lemma \ref{lem1745} it suffices to construct a homotopy
of homotopies $\mathcal H^i$ appearing in Definition \ref{defn1528} (4).
Recall that $\mathcal H^i$ is a $[0,1]^2$-parametrized partial morphism
from $\mathcal F_a^i$ to $\mathcal F_c^{i+1}$
such that its normalized boundary $\partial\mathcal H^i$
is a disjoint union of the following 4 homotopies.
\begin{enumerate}
\item[(i)]
$\frak H_{cb}^i$.
\item[(ii)]
$\frak H^{i+1} \circ  \frak N^{i+1 i}_a$.
\item[(iii)]
$\frak H_{ca}^i =
(\frak N_{cb}^{i+1}\circ\frak H_{ba}^i) \cup
(\frak H_{cb}^i\circ \mathcal{ID})$.
\item[(iv)]
$ \frak N^{i+1 i}_c \circ \frak H^{i}$.
\end{enumerate}
We will construct an interpolation space
$$
\mathcal N(ac,i i+1,\alpha_-,\alpha_+;[1,2]^2)
$$ of the homotopy of homotopies $\mathcal H^i$
by modifying the interpolation space
$$
\mathcal N(bc,i i+1,\alpha_-,\alpha_+;[1,2])
$$
of $\frak H_{cb}^i$
in a way similar to the proof of
Proposition \ref{properidentity} as follows.
We note that
the restriction of $\mathcal N(bc,i i+1;\alpha_-,\alpha_+;[1,2])$
to $1 \in \partial [1,2]$ and to $2 \in \partial [1,2]$
is the union of the following fiber products,
respectively.
\begin{enumerate}
\item[(I)]
$\mathcal N(b,ii+1;\alpha_-,\alpha)
\times_{R_{\alpha}}^{\boxplus\tau}
\mathcal N(bc,i+1;\alpha,\alpha'_+)$.
\item[(II)]
$
\mathcal N(bc,i;\alpha_-,\alpha')
\times_{R_{\alpha'}}^{\boxplus\tau}
\mathcal N(c,ii+1;\alpha',\alpha'_{+})
$.
\end{enumerate}
There are two other kinds of boundary of
$\mathcal N(bc,i i+1;\alpha_-,\alpha_+;[1,2])$
as follows.
\begin{enumerate}
\item[(III)]
$\mathcal M(b,i;\alpha_-,\alpha)
\times_{R_{\alpha}}\mathcal N(bc,i i+1;\alpha,\alpha_+;[1,2])$.
\item[(IV)]
$
\mathcal N(bc,i i+1;\alpha,\alpha';[1,2])
\times_{R_{\alpha'}}
\mathcal M(b,i;\alpha',\alpha'_+)$.
\end{enumerate}
Let $C$ be a sufficiently large positive number.
We assume that it is large enough compared to the energy loss of
$\frak N^{cb}_{i}$.
The top dimensional stratum of our
interpolation space
$\mathcal N(ac,i i+1;\alpha_-,\alpha_+;[1,2])$
is
\begin{enumerate}
\item[(1)]
$\overset{\circ}{\mathcal N}(bc,i i+1;\alpha_-,\alpha_+;[1,2])
\times (E(\alpha_-),E(\alpha'_+) + C)$.
\end{enumerate}
This is the only stratum of top dimension.
Below we list up the codimension one strata:
\begin{enumerate}
\item[(2)]
$R_{\alpha_-} \times \{E(\alpha_-)\}$
\newline
$\times_{R_{\alpha}} \overset{\circ}{\mathcal N}(bc,i i+1;\alpha_-,\alpha_+;[1,2])
$
\vskip0.1cm
\item[(3)]
$\overset{\circ}{\mathcal M}(b,i;\alpha_-,\alpha) \times
(E(\alpha_-),E(\alpha))$
\newline
$\times_{R_{\alpha}}^{\boxplus\tau}
\overset{\circ}{\mathcal N}(bc,i i+1;\alpha,\alpha_+;[1,2])
$
\vskip0.1cm
\item[(4)]
$\overset{\circ}{\mathcal M}(b,i;\alpha_-,\alpha) $
\newline
$\times_{R_{\alpha}}
\overset{\circ}{\mathcal N}(bc,i i+1;\alpha,\alpha'_+;[1,2])
\times
(E(\alpha),E(\alpha'_+)+C)
$
\vskip0.1cm
\item[(5)]
$
\overset{\circ}{\mathcal N}(b,i i+1;\alpha_-,\alpha)
\times (E(\alpha_-),E(\alpha))
$
\newline
$
\times_{R_{\alpha}}^{\boxplus\tau}
\overset{\circ}{\mathcal N}(bc,i+1;\alpha,\alpha'_+)
$
\vskip0.1cm
\item[(6)]
$
\overset{\circ}{\mathcal N}(b,i i+1;\alpha_-,\alpha)
$
\newline
$
\times_{R_{\alpha}}^{\boxplus\tau}
\overset{\circ}{\mathcal N}(bc,i+1;\alpha,\alpha'_+)
\times (E(\alpha),E(\alpha'_+))
$
\vskip0.1cm
\item[(7)]
$R_{\alpha_-} \times \{E(\alpha'_+) + C\}$
\newline
$\times_{R_{\alpha}} \overset{\circ}{\mathcal N}(bc,i i+1;\alpha_-,\alpha_+;[1,2])
$
\vskip0.1cm
\item[(8)]
$
\overset{\circ}{\mathcal N}(bc,i;\alpha_-,\alpha')
\times  (E(\alpha_-),E(\alpha'_+)+C)
$
\newline
$
\times^{\boxplus\tau} _{R_{\alpha'}}
\overset{\circ}{\mathcal N}(c,i i+1;\alpha',\alpha'_+)
$
\vskip0.1cm
\item[(9)]
$
\overset{\circ}{\mathcal N}(bc,i i+1;\alpha_-,\alpha';[1,2])
\times  (E(\alpha_-),E(\alpha'_+)+C)
$
\newline
$
\times_{R_{\alpha'}}
\overset{\circ}{\mathcal M}(c,i;\alpha',\alpha'_+)
$
\end{enumerate}
Observe that
(3) $\cup$ (4), (5) $\cup$ (6), (8), (9) are product of (III), (I), (II), (IV)
and the interval $(E(\alpha_-),E(\alpha'_+)+C)$,
respectively.
\par
We also find that
$(2) \cup (3) = \frak H^{cb}_i\circ \mathcal{ID}$ and
$(5) = \frak N^{cb}_{i+1}\circ\frak H^{ba}_i$.
Therefore $(2) \cup (3) \cup (5)= {\rm (iii)}$.
We also have $(6) =  {\rm (ii)}$, $(7) =  {\rm (i)}$.
\par
We can regard $(8) =  {\rm (iv)}$.
Here the interval $(E(\alpha_-),E(\alpha'_+)+C)$ appearing in (8)
is not
$(E(\alpha_-),E(\alpha')+C)$, but they are diffeomorphic.
\par
Moreover, (4) is nothing but
$$
\overset{\circ}{\mathcal M}(a,i;\alpha_-,\alpha)
\times_{R_{\alpha}}
\overset{\circ}{\mathcal N}(ac,i i+1,\alpha,\alpha_+;[1,2]^2)
$$
and we can regard (9) as
$$
\overset{\circ}{\mathcal N}(ac,i i+1,\alpha_-,\alpha';[1,2]^2)
\times_{R_{\alpha'}}
\overset{\circ}{\mathcal M}(c,i+1;\alpha',\alpha).
$$
Thus the boundary of
$\mathcal N(ac,i i+1,\alpha_-,\alpha_+;[1,2]^2)$
has the required properties.
\par
We note that we need to smooth corners and bend a part of the
boundary so that corner structure stratification of
$\mathcal N(ac,i i+1,\alpha_-,\alpha_+;[1,2]^2)$
becomes a correct one.
In fact, since the union of (2)(3)(5) is (iv), we need to
smooth the corner where they intersect.
For example, the intersection of (3) and (5) is
$$
\aligned
&\overset{\circ}{\mathcal M}(b,i;\alpha_-,\alpha_1) \times^{\boxplus\tau}_{R_{\alpha_1}}
\overset{\circ}{\mathcal N}(b,i i+1;\alpha_1,\alpha_2)
\times_{R_{\alpha}}^{\boxplus\tau}
\overset{\circ}{\mathcal N}(bc,i i+1;\alpha_2,\alpha_+;[1,2])
\\
&
\times (E(\alpha_1),E(\alpha_2)).
\endaligned
$$
So we smooth the corner here.
On the other hand
$$
\aligned
&\overset{\circ}{\mathcal M}(b,i;\alpha_-,\alpha_1) \times^{\boxplus\tau}_{R_{\alpha_1}}
\overset{\circ}{\mathcal N}(b,i i+1;\alpha_1,\alpha_2)
\times_{R_{\alpha}}^{\boxplus\tau}
\overset{\circ}{\mathcal N}(bc,i i+1;\alpha_2,\alpha_+;[1,2])
\\
&
\times (E(\alpha_-),E(\alpha_1))
\endaligned
$$
is a part of the boundary of (3) which is not contained in
the boundary of (5).
We bend the boundary component
$$
\aligned
&\overset{\circ}{\mathcal M}(b,i;\alpha_-,\alpha_1) \times^{\boxplus\tau}_{R_{\alpha_1}}
\overset{\circ}{\mathcal N}(b,i i+1;\alpha_1,\alpha_2)
\times_{R_{\alpha}}^{\boxplus\tau}
\overset{\circ}{\mathcal N}(bc,i i+1;\alpha_2,\alpha_+;[1,2])
\\
&
\times (E(\alpha_-),E(\alpha_2))
\endaligned
$$
at
$$
\aligned
&\overset{\circ}{\mathcal M}(b,i;\alpha_-,\alpha_1) \times^{\boxplus\tau}_{R_{\alpha_1}}
\overset{\circ}{\mathcal N}(b,i i+1;\alpha_1,\alpha_2)
\times_{R_{\alpha}}^{\boxplus\tau}
\overset{\circ}{\mathcal N}(bc,i i+1;\alpha_2,\alpha_+;[1,2])
\\
&
\times \{E(\alpha_1)\}.
\endaligned
$$
This bending is included in the process we constructed
the homotopy $\frak H^i$ between
$\frak N_{cb}^{i}$ and $\mathcal{ID} \sim \frak N_{cb}^{i}$.
(See the proof of Proposition \ref{properidentity}.)
\par
We also note that in (4)(9) we take $\times_{R_{\alpha}}$
etc. but in other places we take
$\times_{R_{\alpha}}^{\boxplus}$.
This does {\it not} cause inconsistency at the intersection of (3) and (4)
by the following reason.
We first take the fiber product
$$
\overset{\circ}{\mathcal M}(b,i;\alpha_-,\alpha)
\times_{R_\alpha}
\overset{\circ}{\mathcal N}(bc,i i+1;\alpha,\alpha_+;[1,2])
\times
(E(\alpha_-),E(\alpha'_+))
$$
and bend it at
$$
\overset{\circ}{\mathcal M}(b,i;\alpha_-,\alpha)
\times_{R_\alpha}
\overset{\circ}{\mathcal N}(bc,i i+1;\alpha,\alpha_+;[1,2])
\times
\{ E(\alpha) \}.
$$
We then partially trivialize by choosing $\frak C$ =
union of components
$$
\overset{\circ}{\mathcal M}(b,i;\alpha_-,\alpha)
\times_{R_\alpha}
\overset{\circ}{\mathcal N}(bc,i i+1;\alpha,\alpha_+;[1,2])
\times
(E(\alpha_-),E(\alpha)).
$$
Then when we consider the corner
$$
\overset{\circ}{\mathcal M}(b,i;\alpha_-,\alpha_1)
\times_{R_{\alpha_1}}
\mathcal M(b,i;\alpha_1,\alpha_2)
\times_{R_{\alpha_2}}
\overset{\circ}{\mathcal N}(bc,i i+1;\alpha_2,\alpha_+;[1,2])
\times
(E(\alpha),E(\alpha'_+)),
$$
the part where the $\frak C$-trivialization is performed is
nothing but
$$
\overset{\circ}{\mathcal M}(b,i;\alpha_-,\alpha_1)
\times_{R_{\alpha_1}}
\mathcal M(b,i;\alpha_1,\alpha_2)
\times_{R_{\alpha_2}}
\overset{\circ}{\mathcal N}(bc,i i+1;\alpha_2,\alpha_+;[1,2])
\times
(E(\alpha_2),E(\alpha'_+)).
$$
This is because the $\frak C$-trivialization is performed
at the intersection of the two
boundary components which belongs to $\frak C$.
We smooth those corners to obtain
the union of (3) and (4).
\par
In the same way as in the proof of Proposition \ref{properidentity},
we can describe the higher codimension boundary of our
K-space $\mathcal N(ac,i i+1,\alpha_-,\alpha_+;[1,2]^2)$,
and can check the consistency in a straight forward way.
This finishes the proof of Lemma \ref{lem1745}.
\end{proof}
The proof of Theorem \ref{indlinesysmainth1} (6) is now complete.
\end{proof}
Thus we have completed the proof of all the results claimed in Section
\ref{sec:systemline1}.

\begin{rem}
In this section we use the homotopy of homotopies version of algebraic lemmas
(Proposition \ref{prop1639}) for the promotion of homotopies.
We can avoid it if we use the next lemma instead.
\begin{lem}\label{lemshortcut}
Let $\psi_1, \psi_2 : CF(\mathcal F) \to CF(\mathcal F')$ be two gapped cochain maps.
(See Definition Definition \ref{defn141212} for gapped-ness.)
Suppose $\psi_1\vert_{E}$ is partially cochain homotopic to $\psi_2\vert_{E}$ for any $E$,
up to energy cut level $E$. Then $\psi_1$ is cochain homotopic to $\psi_2$.
\end{lem}
We can prove Lemma \ref{lemshortcut} in the same way as in the proof of
\cite[Lemma 7.2.177]{fooobook2}.
Note that Lemma \ref{lemshortcut} itself is valid for any ground ring $R$,
while \cite[Lemma 7.2.177]{fooobook2} is proved only for the case when
the ground ring is $\C$ or a finite field. (In fact, \cite[Remark 7.2.181]{fooobook2} gives a counter example to \cite[Lemma 7.2.177]{fooobook2} for the case when the ground ring is $\Q$.)
The reason why Lemma \ref{lemshortcut} is valid for any ground ring $R$ is that
the equation for a map to be a cochain map or a cochain homotopy is linear.
On the other hand, the equation for a map to be an $A_{\infty}$ map is nonlinear.
\par
We choose to prove Proposition \ref{prop1639} rather than using Lemma \ref{lemshortcut},
it is more direct to apply Proposition \ref{prop1639} to other situations.
\end{rem}

\section{Linear K-system: Floer cohomology III:
Morse case by multisection}
\label{sec:systemline4}

In this section we provide the technical detail of the way how we
associate Floer cohomology with {\it ground ring $\Q$}
to a linear K-system when the critical submanifolds
are $0$ dimensional.

\begin{defn}
We say a (partial) linear K-system as in Condition \ref{linsysmainconds} and
Definition \ref{linearsystemdefn} is
{\it of Morse type}\index{K-system ! linear K-system of Morse type}
if all the critical submanifolds $R_{\alpha}$
consist of finite set.
\end{defn}

\begin{thm}\label{thm182222}
For (partial) linear K-system of Morse type,
Theorems \ref{linesysmainth1}, \ref{linesysmainth2}, \ref{indlinesysmainth1}
hold with $\R$ replaced by $\Q$.
\end{thm}

The proof is done by replacing {\it CF-perturbations} used in Section \ref{sec:systemline3}
by {\it multivalued perturbations} (or multisections) in this section.
We recall that we defined the notion of fiber product of CF-perturbations
in \cite[Section 10]{part11}.  There is certain trouble to define the notion of fiber product
of multivalued perturbations. (In fact, it is not a correct way to assume
that the evaluation map restricted to the zero set of multisection is submersive,
because this assumption is not satisfied
even in the case of generic perturbations by an obvious dimensional reason.)
On the other hand, in the case of direct product or fiber product over
$0$ dimensional spaces, we can define the notion of (fiber) product of
multivalued perturbations in an obvious way.
This is the main idea of the proof.
To work it out in detail, we need to check carefully that the whole proof
in Section \ref{sec:systemline3} (and ones in Section \ref{sec:triboundary}
which are used in Section \ref{sec:systemline3})
can be carried out using multivalued perturbation
in place of CF-perturbations.
Indeed, once we correctly state a series of lemmas, their proofs are
either automatic or straightforward analogue of the corresponding results in Sections \ref{sec:triboundary} or \ref{sec:systemline3}.
\par
The contents of this section are not used in the other part of this article.
So many of the readers may skip this section and directly go to
Section \ref{sec:systemtree1}.

\subsection{Bundle extension data revisited}
\label{subsec:18-1bundleextension}

In \cite[Subsections 13.1, 13.2, 13.4]{part11}
we proved existence of multisection or multivalued perturbation
on good coordinate system.
To resolve a technical problem mentioned
in \cite[Subsection 13.5]{part11},
we used the notion of bundle extension data (\cite[Definition 12.24]{part11}).
In our situation where we have a system of Kuranishi structures (K-system),
we first define a similar notion for a Kuranishi structure and
then discuss the compatibility of them with the
fiber product description of boundaries.
We will study this point in detail in this subsection.

\begin{defn}\label{coorchangebed}
Let $\mathcal U_i = (U_i,\mathcal E_i,\psi_i,s_i)$ $(i=1,2)$ be Kuranishi charts of $X$ and
$\Phi_{21} =
(U_{21},\varphi_{21},\widehat{\varphi}_{21}) : \mathcal U_1 \to \mathcal U_2$
a coordinate change.
A {\it bundle extension data}
\index{bundle extension data ! of Kuranishi chart}
associated with $\Phi_{21}$ is
${\frak O}_{12} = (\pi_{12},\tilde\varphi_{21},\Omega_{12})$
satisfying the following:
\begin{enumerate}
\item
$\pi_{12} : \Omega_{12} \to U_1$ is a continuous map,
where $\Omega_{12}$ is a neighborhood of $\varphi_{21}(U_{21})$
in $U_2$.
\item
$\pi_{12}$ is diffeomorphic to the projection
of the normal bundle. (See \cite[Definition 12.23]{part11}.)
\item $\tilde\varphi_{21} : \pi_{12}^*\mathcal E_1 \to \mathcal E_2$
is an embedding of vector bundle.
(See Definition \ref{def26222}.)
\item
The map
$
\varphi_{21}^* \pi^*_{12} \mathcal E_{1} \to \varphi_{21}^*\mathcal E_{2}
$
that is induced from $\tilde\varphi_{21}$
and $\varphi_{21}$ coincides with the bundle
map
$
\widehat{\varphi}_{21} : \mathcal E_{1} \to \mathcal E_{2}
$
which covers $\varphi_{21}$.
\end{enumerate}
\end{defn}
Definition \ref{coorchangebed} is mostly the same as \cite[Definition 12.24]{part11}.
In \cite[Definition 12.24]{part11} we take a compact subset $Z$ in $U_1$ and
$\Omega_{12}$ as a neighborhood of $\varphi (Z)$, while
in Definition \ref{coorchangebed} we do not take a compact subset $Z \subset U_1$
but take $\Omega_{12}$ as a neighborhood of $\varphi (U_1)$ itself.
This difference is not essential since when we use such structures we may restrict them
to compact sets.

\begin{defn}\label{compcoorchangebed}
Let $\mathcal U_i = (U_i,\mathcal E_i,\psi_i,s_i)$ $(i=1,2,3)$ be Kuranishi charts of $X$,
$\Phi_{j i} =
(U_{j i},\varphi_{j i},\widehat{\varphi}_{j i}) : \mathcal U_i
\to \mathcal U_{j}$ $(i,j \in \{1,2,3\}, i<j)$ coordinate
changes and let
$\frak O_{ji} = (\pi_{i j},\tilde\varphi_{j i},\Omega_{i j})$ be bundle extension data associated to  $\Phi_{j i}$.
\begin{enumerate}
\item
We say $\Phi_{21}$, $\Phi_{32}$ are {\it compatible with $\Phi_{31}$} if
\begin{equation}
\Phi_{31}\vert_{U_{321}} = \Phi_{32}\circ
\Phi_{21}\vert_{U_{321}}.
\end{equation}
holds on $U_{321} = \varphi_{21}^{-1}(U_{32}) \cap U_{31}$.
(Note this is the same as \cite[Definition 3.8]{part11}.)
\item
In the situation of (1),
we say $\frak O_{21}$, $\frak O_{32}$ are {\it compatible with $\frak O_{31}$} if
the following holds.
\begin{enumerate}
\item
$\pi_{23} \circ \pi_{12} = \pi_{13}$
on
$\pi^{-1}_{12}(\Omega_{23})
\cap \Omega_{12}
\cap \Omega_{13}$.
\item
$\tilde\varphi_{32} \circ \tilde\varphi_{21} = \tilde\varphi_{31}$
on
$\pi^{-1}_{12}(\Omega_{23})
\cap \Omega_{12}
\cap \Omega_{13}$.
\end{enumerate}
This is mostly the same as the compatibility in \cite[Definition 13.7]{part11}.
\end{enumerate}
\end{defn}

\begin{defn}\label{defn18555}
\begin{enumerate}
\item
Let $\widehat{\mathcal U}$ be a Kuranishi structure of $X$.
A {\it bundle extension data}
\index{bundle extension data ! of Kuranishi structure}
of $\widehat{\mathcal U}$
associates $\frak O_{pq}$ to each coordinate change
$\Phi_{pq}$ of $\widehat{\mathcal U}$
so that they are compatible in the sense of Definition \ref{compcoorchangebed} (2).
\item
Let $\widetriangle{\mathcal U}$ be a good coordinate system of $X$.
A {\it bundle extension data} of $\widetriangle{\mathcal U}$
associates $\frak O_{\frak p \frak q}$ to each coordinate change
$\Phi_{\frak p \frak q}$ of $\widetriangle{\mathcal U}$ so that they are compatible in the sense of Definition \ref{compcoorchangebed} (2).
\item
We can define a bundle extension data of various kinds of embeddings between Kuranishi structures and/or
good coordinate systems, in the same way.
\end{enumerate}
\end{defn}

\begin{lem}
Let $\widetriangle{\mathcal U}$ be a good coordinate system of $X$.
Then there exists
a bundle extension data
for any proper open substructure of $\widetriangle{\mathcal U}$.
\end{lem}
\begin{proof}
This is an immediate consequence of
\cite[Proposition 13.9]{part11}.
\end{proof}
In general, for a given Kuranishi structure,
it seems hard to prove existence of bundle extension data.
In fact, the proof of \cite[Proposition 13.9]{part11} is done by induction on the charts.
The definition of good coordinate system is designed so that such an induction
does work.
On the other hand,
since there may be infinitely many charts in Kuranishi structure,
the definition of Kuranishi structure is not suitable to work out such an induction.
This point is the same as for the construction of multisection.
The way to resolve this trouble is also the same.
Namely we first take a good coordinate system and
construct a bundle extension data associated to the good coordinate system.
Then we restrict it
to obtain one on (slightly different) Kuranishi structure.
Namely we have the next lemma.

\begin{lem}\label{lem187777}
Suppose we are in the situation of \cite[Proposition 6.44]{part11}.
We assume that there exists a bundle extension data $\widetriangle{\frak O_0}$
on good coordinate system $\widetriangle{\mathcal U_0}$.
Then we can take a Kuranishi structure $\widehat{\mathcal U}$ so that there exists a bundle extension data
$\widehat{\frak O}$ on it such that $\widetriangle{\frak O_0}$ and
$\widehat{\frak O}$ are compatible with the embedding
$\widetriangle{\mathcal U_0} \to \widehat{\mathcal U}$.
\end{lem}
\begin{proof}
By the construction of ${\widehat{\mathcal U}}$ given in the proof of
\cite[Proposition 6.44]{part11},
all the coordinate changes appearing in ${\widehat{\mathcal U}}$ are restrictions of coordinate changes of
$\widetriangle{\mathcal U_0}$.
Therefore we can obtain the bundle extension data
$\widehat{\frak O}$ by restricting $\widetriangle{\frak O_0}$.
The compatibility among members of $\widehat{\frak O}$ follows from the
compatibility among members of $\widetriangle{\frak O_0}$.
The embedding $\widetriangle{\mathcal U_0} \to \widehat{\mathcal U}$ also
consists of restriction of  coordinate changes of $\widetriangle{\frak O_0}$.
Therefore the compatibility of $\widetriangle{\frak O_0}$ and
$\widehat{\frak O}$ with the embedding also follows from the
compatibility among members of $\widetriangle{\frak O_0}$.
\end{proof}
We can prove a similar result in the situation of \cite[Theorem 3.30]{part11}.
\begin{lem}\label{lem187888}
Suppose we are in the situation of \cite[Theorem 3.30]{part11}.
We assume that we have a bundle extension data $\widehat{\frak O}$ of $\widehat{\mathcal U}$.
Then there exists a bundle extension data $\widetriangle{\frak O}$ of
the good coordinate system $\widetriangle{\mathcal U}$
such that $\widehat{\frak O}$ and $\widetriangle{\frak O}$ are compatible with
respect to the embedding $\widehat{\mathcal U} \to \widetriangle{\mathcal U}$.
\end{lem}
\begin{proof}
The proof of \cite[Theorem 3.30]{part11} is done by induction via the statement
formulated as \cite[Proposition 11.3]{part11}.
We can amplify \cite[Proposition 11.3]{part11} by adding the statement
that all the embeddings and coordinate changes are associated with
bundle extension data and they are all compatible with respect to
compositions appearing there.
\par
To prove this amplified statement,
we proceed in the same way as in the proof of \cite[Proposition 11.3]{part11}.
All we need to note are the following two points:
\begin{enumerate}
\item
We glued two coordinate changes sometimes on open sets
where they coincide each other.
\item
We invert coordinate changes which are isomorphisms.
\end{enumerate}
\par
As for point (1)
we note that whenever two coordinate changes coincide during our construction,
the bundle extension data coincide as well. Therefore we can also glue the
bundle extension data in an obvious way.
\par
As for point (2), we note that in the case when coordinate changes or
embeddings of Kuranishi charts are local isomorphisms,
bundle extension data exists and is unique up to the
process of restricting $\Omega_i$ etc. to its open subsets. So we can obviously
invert them.
Thus we obtain the required
bundle obstruction data $\widetriangle{\frak O}$ of $\widetriangle{\mathcal U}$.
\end{proof}
There are various relative versions etc. in Part I which we can prove together with bundle
extension data. The way to generalize them is straightforward so we omit them
and mention only when we use them.
\par
Next we discuss compatibility of multisection and multivalued perturbation
with bundle extension data.
\begin{defn}\label{defb18189}
For each $i=1,2$ let $\mathcal U_i$, $\Phi_{21}$, ${\frak O}_{21}$ be as in Definition \ref{coorchangebed} and
$\frak s_i$ (resp. $\widehat{\frak s_i} = \{\frak s_i^{\epsilon}\}$) a multisection (resp. multivalued perturbation)
of $\mathcal U_i$.
We say $\frak s_1$ is {\it compatible} with $\frak s_2$ with respect to
$\Phi_{21}$, ${\frak O}_{21}$
if they satisfy
\begin{equation}\label{form1226rev}
(y,\frak s_{2}(y))
=
\tilde\varphi_{21}(y,\frak s_{1}(\pi_{12}(y))),
\end{equation}
(resp.
\begin{equation}\label{form1226revrev}
(y,\frak s^{\epsilon}_{2}(y))
=
\tilde\varphi_{21}(y,\frak s^{\epsilon}_{1}(\pi_{12}(y)))),
\end{equation}
holds for all $y \in \Omega_{21}$.
\end{defn}
This definition is mostly the same as \cite[Definition 13.15]{part11}.
\begin{rem}
The equalities (\ref{form1226rev}), (\ref{form1226revrev}) are point-wise equalities
which are supposed
to hold branch-wise.
So the ambiguity of the notion of branches
explained in
\cite[Subsection 13.5]{part11} does not cause any trouble here.
\end{rem}
\begin{defn}
Let ${\widehat{\mathcal U}}$ be a Kuranishi structure and $\widehat{\frak O}$ its bundle extension data.
A multivalued perturbation $\widehat{\frak s} = \{\widehat{\frak s^{\epsilon}_p} \}$ of ${\widehat{\mathcal U}}$ is said to be {\it compatible}
\index{bundle extension data ! compatibility with multivalued perturbation}
\index{multivalued perturbation ! compatibility with bundle extension data}
\index{compatibility ! of multivalued perturbation with bundle extension data}
with $\widehat{\frak O}$ if it is compatible with the pair
$(\Phi_{pq},\frak O_{pq})$ of the coordinate change
and a member of $\widehat{\frak O}$ in the sense of Definition \ref{defb18189}.
\par
For the case of good coordinate system
with bundle extension data, the compatibility
of multivalued perturbation with extension data is defined in the same way.
\end{defn}
\begin{lem}\label{lemexistscompatiext}
Let $\widetriangle{\mathcal U}$ be a good coordinate system of $X$
and $\widetriangle{\frak O}$  its bundle extension data.
Let $\widetriangle{\mathcal U_0}$ be a proper open substructure of
$\widetriangle{\mathcal U}$.
We restrict $\widetriangle{\frak O}$ to it.
Then there exists a multivalued perturbation
$\widetriangle{\frak s} = \{\widetriangle{\frak s^{\epsilon}_{\frak p}} \}$ thereof with the following properties:
\begin{enumerate}
\item $\widetriangle{\frak s}$ is transversal to $0$.
\item It is compatible with
$\widetriangle{\frak O}\vert_{\widetriangle{\mathcal U_0}}$.
\end{enumerate}
\end{lem}
\begin{proof}
In the proof of the \cite[Theorem 6.37]{part11} given in
\cite[Seciton 13]{part11} we used induction to prove \cite[Proposition 13.23]{part11}.
The support system $\mathcal K$  which
appears in the statement of \cite[Proposition 13.23]{part11}
can be chosen so that $\mathcal U_{0,\frak p} \subset \mathcal K_{\frak p}$.
Then Lemma \ref{lemexistscompatiext} is an immediate consequence of \cite[Proposition 13.23]{part11}.
\end{proof}
We can also include multivalued perturbations in Lemmas \ref{lem187777} and \ref{lem187888}.

\begin{lem}\label{lem187777add}
Suppose we are in the situation of Lemma \ref{lem187777}.
We assume in addition that there exists a multivalued perturbation
$\widetriangle{\frak s_0}$
of $\widetriangle{\mathcal U_0}$
compatible with the bundle extension data
$\widetriangle{\frak O_0}$.
Then there exists a multivalued perturbation $\widehat{\frak s}$
of $\widehat{\mathcal U}$ compatible with the
bundle extension data
$\widehat{\frak O}$.
Moreover, $\widetriangle{\frak s_0}$ and $\widehat{\frak s}$
are compatible with the embedding
$\widetriangle{\mathcal U_0} \to \widehat{\mathcal U}$ with respect to
$\widetriangle{\frak O_0}$, $\widehat{\frak O}$.
If $\widetriangle{\frak s_0}$ is transversal to $0$,
so is $\widehat{\frak s}$.
\end{lem}
\begin{proof}
The proof is the same as the proof of Lemma \ref{lem187777}.
\end{proof}
\begin{lem}\label{lem187888add}
Suppose we are in the situation of Lemma \ref{lem187888}.
We assume in addition that we have a multivalued perturbation $\widehat{\frak s}$
of $\widehat{\mathcal U}$ compatible with
the bundle extension data $\widehat{\frak O}$.
Then there exists a multivalued perturbation $\widetriangle{\frak s_0}$
of $\widetriangle{\mathcal U}$ compatible with the
bundle extension data $\widetriangle{\frak O}$.
Moreover, $\widehat{\frak s}$ and $\widetriangle{\frak s_0}$ are
compatible with
the embedding $\widehat{\mathcal U} \to \widetriangle{\mathcal U}$
with respect to $\widehat{\frak O}$, $\widetriangle{\frak O}$.
If $\widehat{\frak s}$ is transversal to $0$,
so is $\widetriangle{\frak s_0}$.
\end{lem}
\begin{proof}
The proof is the same as the proof of Lemma \ref{lem187888}.
\end{proof}

\subsection{Virtual fundamental chain of $0$ dimensional K-space}
\label{subsec:18-2direcproduct}

\begin{defn}\label{defn1815}
Suppose
$X_1$, $X_2$ have Kuranishi structures
$\widehat{\mathcal U_1}$, $\widehat{\mathcal U_2}$ respectively,
and $\widehat f_i : (X_i,\widehat{\mathcal U_i}) \to R$
are strongly smooth.
We assume that $R$ is a $0$-dimensional compact manifold,
that is nothing but a finite set.
Then they are automatically transversal and
we have a fiber product
$
(X_1 \times_R X_2, \widehat{\mathcal U_1} \times_R \widehat{\mathcal U_2}).
$
We call this fiber product the {\it direct-like product}.
\end{defn}
The next lemma is trivial to prove.
\begin{lemdef}
Suppose we are in the situation of Definition \ref{defn1815}.
\begin{enumerate}
\item
If $\widehat{\frak O_i}$ are bundle extension data of $\widehat{\mathcal U_i}$
for $i=1,2$, then they induce bundle extension data
of their direct-like product in a canonical way.
We call it the {\rm direct-like product}\index{direct-like product}
$\widehat{\frak O_1} \times_R \widehat{\frak O_2}$.
\item
If $\widehat{\frak s_i}$ (resp. $\frak s_i$) are multivalued perturbations
(resp. multisections) of $\widehat{\mathcal U_i}$
for $i=1,2$, then they  induce a multivalued perturbation
(resp. multisection) of their fiber product in a canonical way.
\footnote{We do not need to assume that our fiber product is
direct-like in Item (2).}
We call it fiber product of the multivalued perturbation (resp.
multisection) and write $\widehat{\frak s_1} \times_R \widehat{\frak s_2}$
(resp. $\frak s_1 \times_R \frak s_2$).
\item
If $\widehat{\frak s_i}$ (resp. $\frak s_i$) are transversal to $0$,
then its direct-like product is transversal to $0$.
\item
Let $\widehat{\frak O_i}$ be bundle extension data of $\widehat{\mathcal U_i}$
and $\widehat{\frak s_i}$ (resp. $\frak s_i$) multivalued perturbations
(resp. multisections) of $\widehat{\mathcal U_i}$ for $i=1,2$.
We assume $\widehat{\frak s_i}$ (resp. $\frak s_i$) are compatible with
$\widehat{\frak O_i}$.
Then the direct-like product
$\widehat{\frak s_1} \times_R  \widehat{\frak s_2}$
(resp. $\frak s_1 \times_R \frak s_2$)
is compatible with $\widehat{\frak O_1} \times_R \widehat{\frak O_2}$.
\end{enumerate}
\end{lemdef}
Next we define
the virtual fundamental chain,
which is identified with a rational number, of a $0$ dimensional
K-space by using multivalued
perturbation.
We recall that in \cite[Definition 14.6]{part11}
we defined a virtual fundamental chain ($\in \Q$)
of $(X,\widetriangle{\mathcal U},\widetriangle{\frak s})$
where $\widetriangle{\mathcal U}$ is a
good coordinate system of $X$ and
$\widetriangle{\frak s}$ is its multivalued perturbation
transversal to $0$.
We adopt this story and proceed in the same way as in
\cite[Subsection 9.2]{part11} to define a virtual fundamental chain ($\in \Q$)
for a $0$ dimensional K-space.
\begin{shitu}\label{Shitu1817}
\begin{enumerate}
\item
We consider a quadruple
$(X,\widehat{\mathcal U},\widehat{\frak O},\widehat{\frak s})$
such that:
\begin{enumerate}
\item
$(X,\widehat{\mathcal U})$ is a K-space.
\item
$\widehat{\frak O}$ is a bundle extension data of $\widehat{\mathcal U}$.
\item
$\widehat{\frak s}$ is a multivalued perturbation of
$\widehat{\mathcal U}$ compatible with $\widehat{\frak O}$.
\end{enumerate}
\item
We consider a quadruple
$(X,\widetriangle{\mathcal U},\widetriangle{\frak O},\widetriangle{\frak s})$
such that:
\begin{enumerate}
\item
$(X,\widetriangle{\mathcal U})$ is a good coordinate system.
\item
$\widetriangle{\frak O}$ is a bundle extension data of
$\widetriangle{\mathcal U}$.
\item
$\widetriangle{\frak s}$ is a multivalued perturbation of
$\widetriangle{\mathcal U}$ compatible with $\widetriangle{\frak O}$.
\end{enumerate}
\end{enumerate}
$\blacksquare$
\end{shitu}
\begin{defn}\label{defn1818}
\begin{enumerate}
\item
Let $(X,\frak X_i) = (X,\widehat{\mathcal U_i},\widehat{\frak O_i},\widehat{\frak s_i})$ $(i=1,2)$
be as in Situation \ref{Shitu1817} (1).
We say $\frak X_2$ is a
{\it thickening}
\index{thickening ! of Kuranishi structure (with bundle extension data and multivalued perturbation)} of
$\frak X_1$
and write
$(X,\widehat{\mathcal U_1},\widehat{\frak O_1},\widehat{\frak s_1})
<
(X,\widehat{\mathcal U_2},\widehat{\frak O_2},\widehat{\frak s_2})$
if the following holds.
\begin{enumerate}
\item
$(X,\widehat{\mathcal U_1}) \to (X,\widehat{\mathcal U_2})$
by which $(X,\widehat{\mathcal U_2})$ is a thickening of
$(X,\widehat{\mathcal U_1})$.
(See \cite[Definition 5.3]{part11}.)
\item
$\widehat{\frak O_2}$, $\widehat{\frak O_1}$ are compatible with respect to this embedding.
\item
$\widehat{\frak s_2}$, $\widehat{\frak s_1}$
are compatible with this embedding
with respect to $\widehat{\frak O_2}$, $\widehat{\frak O_1}$.
\end{enumerate}
\item
In the case one or both of the objects are good coordinate systems,
we can define the notion of thickening such as
$$
\aligned
&(X,\widetriangle{\mathcal U_1},\widetriangle{\frak O_1},
\widetriangle{\frak s_1})
<
(X,\widehat{\mathcal U_2},\widehat{\frak O_2},\widehat{\frak s_2}),
\\
&(X,\widehat{\mathcal U_1},\widehat{\frak O_1},\widehat{\frak s_1})
<
(X,\widetriangle{\mathcal U_2},\widetriangle{\frak O_2},
\widetriangle{\frak s_2}),
\\
&(X,\widetriangle{\mathcal U_1},\widetriangle{\frak O_1},
\widetriangle{\frak s_1})
<
(X,\widetriangle{\mathcal U_2},\widetriangle{\frak O_2},
\widetriangle{\frak s_2}),
\endaligned
$$
in the same way.
\end{enumerate}
\end{defn}
Below we will define virtual fundamental chain of
a quadruple $(X,\widehat{\mathcal U},\widehat{\frak O},\widehat{\frak s})$
as in Situation \ref{Shitu1817} (1)
provided that $\widehat{\frak s}$ is transversal to $0$.
We begin with the following lemma.

\begin{lem}\label{1819lem}
Let
$(X,\widetriangle{\mathcal U_i},\widetriangle{\frak O_i},
\widetriangle{\frak s_i})$
$i=1,2$ be as in Situation \ref{Shitu1817} (2)
such that $\dim (X,\widetriangle{\mathcal U_i}) = 0$.
We assume
$$
(X,\widetriangle{\mathcal U_1},\widetriangle{\frak O_1},
\widetriangle{\frak s_1})
<
(X,\widetriangle{\mathcal U_2},\widetriangle{\frak O_2},
\widetriangle{\frak s_2}).
$$
Then there exists $\epsilon_0 >0$ such that the following holds for
$0 < \epsilon < \epsilon_0$ :
\par
If $\widetriangle{\frak s_1}$ and $\widetriangle{\frak s_2}$
are transversal to $0$,
then
$$
[(X,\widetriangle{\mathcal U_1},\widetriangle{{\frak s_1^{\epsilon}}})]
=
[(X,\widetriangle{\mathcal U_2},\widetriangle{{\frak s_2^{\epsilon}}})].
$$
\end{lem}
\begin{proof}
The proof is entirely the same as the proof of
\cite[Proposition 9.16]{part11}.
\end{proof}
\begin{defn}\label{defn1822}
Let $(X,\widehat{\mathcal U},\widehat{\frak O},\widehat{\frak s})$ be as
in Situation \ref{Shitu1817} (1) and
$\dim (X,\widehat{\mathcal U}) = 0$.
We assume that $\widehat{\frak s}$ is
transversal to $0$ as a family in the sense of
\cite[Remark 14.11 (2)]{part11}.
Using Lemma \ref{lem187888add}, we take
$(X,\widetriangle{\mathcal U},\widetriangle{\frak O'},\widetriangle{\frak s'})$
such that
$$
(X,\widehat{\mathcal U},\widehat{\frak O},\widehat{\frak s})
<
(X,\widetriangle{\mathcal U},\widetriangle{\frak O'},\widetriangle{\frak s'}).
$$
Then we define the {\it virtual fundamental chain of}
\index{virtual fundamental chain}
$(X,\widehat{\mathcal U},\widehat{\frak O},\widehat{\frak s})$ (at $\epsilon$) by
\begin{equation}\label{form18444}
[(X,\widehat{\mathcal U},\widehat{\frak O},\widehat{\frak s^{\epsilon}})]
=
[(X,\widetriangle{\mathcal U},\widetriangle{\frak O'},
\widetriangle{\frak s^{\prime \epsilon}})]
\in \Q,
\end{equation}
for sufficiently small $\epsilon >0$ such that $\widehat{\frak s^{\epsilon}}$ is
transversal to $0$ at $\epsilon$.
\end{defn}
\begin{lem}
The right hand side of (\ref{form18444}) is independent of
$(X,\widetriangle{\mathcal U},\widetriangle{\frak O'},
\widetriangle{\frak s^{\prime}})$
but depends only on
$(X,\widehat{\mathcal U},\widehat{\frak O},\widehat{\frak s^{\epsilon}})$
and $\epsilon$,
if $\epsilon >0$ is sufficiently small.
\end{lem}
\begin{proof}
The proof is entirely the same as the proof of
\cite[Theorem 9.14]{part11}.
\end{proof}
Now we state the multivalued perturbation versions of Stokes' formula
and composition formula.
\begin{prop}\label{Stokesmulti18}
Let $(X,\widehat{\mathcal U},\widehat{\frak O},\widehat{\frak s^{\epsilon}})$ be as
in Situation \ref{Shitu1817} (1) and
$\dim (X,\widehat{\mathcal U}) = 1$.
We assume that $\widehat{\frak s^{\epsilon}}$ is
transversal to $0$ as a family in the sense of
\cite[Remark 14.1 (2)]{part11}.
Then we have
$$
[\partial (X,\widehat{\mathcal U},\widehat{\frak O},\widehat{\frak s^{\epsilon}})]
= 0.
$$
Here
$\partial (X,\widehat{\mathcal U},\widehat{\frak O},\widehat{\frak s^{\epsilon}})$
is the normalized boundary $\partial (X,\widehat{\mathcal U})$ together with
the restrictions of $\widehat{\frak O}$ and
$\widehat{\frak s^{\epsilon}}$ to the boundary.
\end{prop}
\begin{proof}
This immediately follows from
\cite[Proposition 14.10]{part11} and the definition.
\end{proof}
To state the analogue of the composition formula we begin with
explaining the situation and introducing the notation.

\begin{shitu}\label{situ1825}
Suppose we are in Situation
\ref{defn1815}.
We assume that we are given the quadruples
$(X_i,\widehat{{\mathcal U_i}},\widehat{\frak O_i},\widehat{\frak s_i})$
for $i=1,2$ of
bundle extension data
$\widehat{\frak O_i}$ and multivalued perturbations $\widehat{\frak s_i}$.
\par
For $r \in R$ we put
$
f_i^{-1}(r) = \{x \in X_i \mid f_i(x) = r\}
$.
They carry various structures such as Kuranishi structures, which are induced by
$\widehat{{\mathcal U_i}},\widehat{\frak O_i},\widehat{\frak s_i}$.
Thus we obtain quadruples
$f_i^{-1}(r) \cap (X_i,{\widehat{\mathcal U_i}},\widehat{\frak O_i},
\widehat{\frak s_i})$ as in Situation \ref{Shitu1817} (1).
$\blacksquare$
\end{shitu}
\begin{prop}\label{prop1826}
In Situation \ref{situ1825} we assume in addition that for $i=1,2$
\begin{enumerate}
\item $\dim (X_i,\widehat{{\mathcal U_i}}) = 0$,
\item
$\widehat{\frak s_i}$ are transversal to $0$ as a family in the sense of
\cite[Remark 14.11 (2)]{part11}.
\end{enumerate}
Then the direct-like product
$$
(X_1,{\widehat{\mathcal U_1}},\widehat{\frak O_1},\widehat{\frak s_1})
\times_R
 (X_2,{\widehat{\mathcal U_2}},\widehat{\frak O_2},\widehat{\frak s_2})
$$
is $0$ dimensional and transversal to $0$ as a family.
Moreover we have
\begin{equation}
\aligned
&[(X_1,{\widehat{\mathcal U_1}},\widehat{\frak O_1},\widehat{\frak s_1^{\epsilon}})
\times_R
(X_2,{\widehat{\mathcal U_2}},\widehat{\frak O_2},\widehat{\frak s_2^{\epsilon}})] \\
&=
\sum_{r \in R}
[f_1^{-1}(r) \cap
(X_1,{\widehat{\mathcal U_1}},\widehat{\frak O_1},\widehat{\frak s_1^{\epsilon}})]
[f_2^{-1}(r) \cap
(X_2,{\widehat{\mathcal U_2}},\widehat{\frak O_2},\widehat{\frak s_2^{\epsilon}})].
\endaligned
\end{equation}
\end{prop}
\begin{proof}
The proof is the same as the proof of
\cite[Proposition 10.23]{part11}.
\footnote{It is actually easier than that.}
\end{proof}
We also note the following.
\begin{lem}
In the situation of Proposition \ref{prop1826} we replace (1) by the
following assumption.
\begin{enumerate}
\item[(1)']
~$\dim (X_1,\widehat{{\mathcal U_1}})
=
- \dim (X_2,\widehat{{\mathcal U_2}}) \ne 0.
$
\end{enumerate}
Except this point, we assume that the same condition as in
Proposition \ref{prop1826}.
Then we have
$$
[(X_1,{\widehat{\mathcal U_1}},\widehat{\frak O_1},
\widehat{\frak s_1^{\epsilon}})
\times_R
(X_2,{\widehat{\mathcal U_2}},\widehat{\frak O_2},
\widehat{\frak s_2^{\epsilon}})]
= 0.
$$
\end{lem}
\begin{proof}
We may assume $\dim (X_1,\widehat{{\mathcal U_1}}) < 0$.
Then by \cite[Lemma 14.1]{part11}
$(\frak s_1^{\epsilon})^{-1}(0)$ is an empty set.
Therefore $(\frak s_1^{\epsilon} \times_R \frak s_2^{\epsilon})^{-1}(0)$
is an empty set. The lemma follows.
\end{proof}

\subsection{Extension of multisection from boundary
to its neighborhood}
\label{subsec:18-413mane}

In this subsection we discuss an analogue of the story in Section
\ref{sec:triboundary} for multivalued perturbations.

\begin{defn}
Let $(X',\widehat{\mathcal U'})$ be a $\tau$-collared Kuranishi structure
as in Definition \ref{defn1531} (1).
A {\it $\tau$-collared bundle extension data}
\index{bundle extension data ! of $\tau$-collared Kuranishi structure}
of this $\tau$-collared Kuranishi structure
assigns a bundle extension data $\frak O_{p'q'}$ to each embedding
$\Phi_{p'q'}$ given in Definition \ref{defn1531} (1) (b)
such that in the situation of
Definition \ref{defn1531} (1) (c), the bundle extension data
$\frak O_{p'q'}$, $\frak O_{q'r'}$ are compatible with
$\frak O_{p'r'}$ in the sense of Definition \ref{compcoorchangebed}.
\par
The compatibility of $\tau$-collared multivalued perturbations with
$\tau$-collared bundle extension data
is defined in the same way as in Definition \ref{defb18189}.
\end{defn}
The next lemma is obvious from the definition.
\begin{lem}
We consider the situation of Lemma-Definition \ref{lemdef1522}.
\begin{enumerate}
\item
Let $\widehat{\frak O}$ be
a bundle extension data of $\widehat{\mathcal U}$.
Then it induces a $\tau$-collared bundle extension data of $\widehat{\mathcal U^{\boxplus\tau}}$.
We denote it by $\widehat{\frak O^{\boxplus\tau}}$.
\item
Suppose we are in the situation of (1).
Let $\widehat{\frak s} = \{\frak s^{\epsilon}\}$ be a multivalued perturbation of $(X,\widehat{\mathcal U})$
compatible with $\widehat{\frak O}$.
Then $\widehat{\frak s^{\boxplus\tau}}$ obtained in
Lemma \ref{lemma1523} (6)
is compatible with $\widehat{\frak O^{\boxplus\tau}}$.
\end{enumerate}
\end{lem}
Next we study the situation in Subsection \ref{subsec:extenonechart}.
We note that bundle extension data can be restricted to normalized
corner in an obvious way. We can also pull it back by a covering map.
Moreover, the compatibility of multivalued perturbation with bundle extension data
is preserved by restriction to normalized corner and pull-back by a covering map.

\begin{shitu}\label{shitu1830}
Suppose we are in Situation \ref{sit1526}.
We assume the following in addition.
\begin{enumerate}
\item
We are given a $\tau$-collared bundle extension data
$\widehat{\frak O^+_{S_k}}$ of a $\tau$-collared Kuranishi structure
$\widehat{\mathcal U^+_{S_k}}$.
\item
The restriction of $\widehat{\frak O^+_{S_{\ell}}}$ to $\widehat S_{k}(\widehat S_{\ell}(X),\widehat{\mathcal U_{S_\ell}^+})$
coincides with the pull-back of $\widehat{\frak O^+_{S_{\ell+k}}}$ by the isomorphism in Situation \ref{sit1526} (2). $\blacksquare$
\end{enumerate}
\end{shitu}
\begin{rem}
Here and in several similar places, we can say two bundle extension data or multisections
are `being the same' or `coincide'.
In the case of CF-perturbations, for example,
we say two such objects `being equivalent' or `isomorphic' etc., instead.
Recall that bundle extension data consists of maps and subsets.
The maps are ones between the spaces
common in the two bundle extension data in question
and the subsets are ones of the sets which are common in the two bundle extension data in question. On the other hand
in the case of CF-perturbations, for example, we also have
a parameter space $W$ of perturbations
as a part of the object. So we can say two of them are equivalent but it does not make so much sense to say that they are
the same.
The fact that we can say two bundle extension data are the same
slightly simplifies discussion here.
\end{rem}
\begin{lem}\label{lem1831}
Suppose we are in Situation \ref{shitu1830}.
Then the $\tau'$-collared Kuranishi structure
$\widehat{\mathcal U^+}$ obtained in Proposition \ref{prop528}
carries a $\tau'$-collared bundle extension data
$\widehat{\frak O^+}$ such that
the restriction
of $\widehat{\frak O^+}$ to $\widehat S_k(X)$
coincides with $\widehat{\frak O^+_{S_k}}$ under the isomorphism
of Proposition \ref{prop528} (1).
\end{lem}
\begin{proof}
The proof is immediate from the proof of Proposition \ref{prop528}.
\end{proof}
\begin{rem}
We do not assume that we have bundle extension data of
$(X,\widehat{\mathcal U})$. In our application,
the Kuranishi structure $(X,\widehat{\mathcal U})$ is one we obtain
from {\it geometry}.
As we explained in Subsection \ref{subsec:18-1bundleextension},
it seems hard to find bundle extension data, in general.
We take a good coordinate system compatible with the Kuranishi structure $\widehat{\mathcal U}$  and use it to find a multivalued
perturbation and bundle extension data.
Then $\widehat{\mathcal U^+}$ is a Kuranishi structure obtained from this good
coordinate system. We will construct bundle extension data and
multisections compatible with the direct-like product description of the corners
inductively. The results we are explaining here will be used for this purpose.
\par
Nevertheless, we can also prove the following type of statements.
If we are given bundle extension data of $\widehat{\mathcal U}$
such that this bundle extension data and
$\widehat{\frak O^+_{S_k}}$ etc.
are compatible with the embeddings, covering maps etc.
appearing in Situation \ref{sit1526} (3)(4)(5), then
the $\tau'$-collared bundle extension data $\widehat{\frak O^+}$
in Lemma \ref{lem1831}
can be taken so that it is compatible with
the embeddings and covering maps appearing in Proposition \ref{prop528}.
We do not state or prove it here since we do not use it.
\end{rem}
\begin{lem}\label{lem183333}
Suppose we are in the situation of Situation \ref{shitu1830}.
We assume in addition that we are given a $\tau$-collared multivalued perturbation
$\widehat{\frak s_{S_k}^+}$ of $\widehat{\mathcal U_{S_k}^+}$
such that
\begin{enumerate}
\item[(i)] $\widehat{\frak s_{S_k}^+}$ is compatible with
$\widehat{\frak O^+_{S_k}}$.
\item[(ii)]
The pull-back
of $\widehat{\frak s_{S_{k+\ell}}^+}$
by the map
$
\pi_{k,\ell} : \widehat S_{k}(\widehat S_{\ell}(X),\widehat{\mathcal U^+_{ S_{\ell}}})
\to (\widehat S_{k+\ell}(X),\widehat{\mathcal U^+_{ S_{k+\ell}}})$
coincides with the restriction of  $\widehat{\frak s_{S_{\ell}}^+}$.
\end{enumerate}
Then for any $0< \tau' < \tau$ there exists a $\tau'$-collared multivalued perturbation
$\widehat{\frak s^+}$ on the Kuranishi structure
$\widehat{\mathcal U^+}$ obtained in Proposition \ref{prop528}
with the following properties.
\begin{enumerate}
\item
$\widehat{\frak s^+}$ is compatible with the $\tau'$-collared bundle extension data $\widehat{\frak O^+}$
obtained in Lemma \ref{lem1831}.
\item Its restriction to $(\widehat S_k(X),\widehat{\mathcal U_{S_k}^+})$
coincides with $\widehat{\frak s_{S_k}^+}$.
\end{enumerate}
\end{lem}
\begin{proof}
The proof is the same as the proof of Proposition \ref{prop529}.
\end{proof}
\par\medskip
Now we discuss results corresponding to those in Subsection \ref{subsection:concltrisection}.
\begin{lem}\label{lem183434}
Suppose we are in the situation of Lemma \ref{lem183333}.
We assume that $\widehat{\frak s^+_{k}}$ are transversal to $0$ as a family.
Then there exists a
$\tau'$-collared multivalued perturbations $\widehat{\frak s^{++}}$ on the Kuranishi structure
$\widehat{\mathcal U^{++}}$ obtained in Proposition \ref{prop1562}
such that
\begin{enumerate}
\item
Its restriction to $(\widehat S_k(X),\widehat{\mathcal U_{S_k}^+})$
coincides with $\widehat{\frak s^+_{k}}$.
\item
$\widehat{\frak s^{++}}$ is transversal to $0$ as a family.
\item
The bundle extension data $\widehat{\frak O^+}$ obtained in Lemma \ref{lem1831}
can be extended to a collared bundle extension data
$\widehat{\frak O^{++}}$ of $\widehat{\mathcal U^{++}}$
which coincides with $\widehat{\frak O^+}$ in a neighborhood of the boundary.
\item
$\widehat{\frak s^{++}}$ is compatible with $\widehat{\frak O^{++}}$.
\end{enumerate}
\end{lem}
\begin{proof}
The proof is the same as the proof of Proposition \ref{prop529rev}.
\end{proof}

\subsection{Completion of the proof of Theorem \ref{thm182222}.}
\label{subsec:18-414comp}

\begin{proof}[Proof of Theorem \ref{thm182222}]
We first show a version of Proposition \ref{prop161}.
\begin{prop}\label{prop183939}
Suppose we are in the situation of Proposition \ref{prop161}.
Moreover we assume that the linear K-system is of Morse type.
\par
Then for any $0< \tau < 1$
there exists a $\tau$-collared Kuranishi structure
$\widehat{\mathcal U^+}(\alpha_-,\alpha_+)$ of
$\mathcal M(\alpha_-,\alpha_+)^{\boxplus\tau_0}$,
its bundle extension data  $\widehat{\frak O}(\alpha_-,\alpha_+)$,
multivalued perturbations
$\widehat{\frak s^+}(\alpha_-,\alpha_+)$ of $\widehat{\mathcal U^+}(\alpha_-,\alpha_+)$
for every $\alpha_-$, $\alpha_+$ with $E(\alpha_+) - E(\alpha_-) \le E_{\frak E}^k$,
such that $\widehat{\frak s^+}(\alpha_-,\alpha_+)$ is compatible with
$\widehat{\frak O}(\alpha_-,\alpha_+)$ and
enjoy the following properties.
\begin{enumerate}
\item
The same as Proposition \ref{prop161} (1).
\item
$\widehat{\frak s^+}(\alpha_-,\alpha_+)$ is transversal to $0$.
\item
The pull-back of $\widehat{\frak s^+}(\beta\alpha_-,\beta\alpha_+)$
by the periodicity isomorphism coincides with $\widehat{\frak s^+}(\alpha_-,\alpha_+)$.
\item Proposition \ref{prop161} (4) holds .
Moreover $\widehat{\frak O}(\alpha_-,\alpha_+)$ is preserved by the periodicity isomorphism.
\item
The pull-back of
$\widehat{\frak s^+}(\alpha_-,\alpha_+)$
by the isomorphism
(\ref{formula162}) coincides with
the fiber product
$$
\widehat{\frak s^+}(\alpha_-,\alpha)
{}_{{\rm ev}_+}\times_{{\rm ev}_-}
\widehat{\frak s^+}(\alpha,\alpha_+).
$$
The pull-back of
$\widehat{\frak O}(\alpha_-,\alpha_+)$
by the isomorphism
(\ref{formula162}) coincides with
the fiber product
$$
\widehat{\frak O}(\alpha_-,\alpha)
{}_{{\rm ev}_+}\times_{{\rm ev}_-}
\widehat{\frak O}(\alpha,\alpha_+).
$$
These fiber products are well-defined since they are direct-like products.
\item
Proposition \ref{prop161} (6) holds.
Moreover (\ref{cornecominduction}) preserves
bundle extension data.
\item
Proposition \ref{prop161} (7) holds.
Moreover $\pi_{\ell,k}$ preserves bundle extension data.
\item
The pull-back of $\widehat{\frak s^+}(\alpha_-,\alpha_+)$
by the isomorphism (\ref{cornecominduction})
is the fiber product
\begin{equation}\label{fiberproductSinduction}
\widehat{\frak s^+}(\alpha_-,\alpha_1)\,\,
{}_{{\rm ev}_{+}}\times_{R_{\alpha_1}}
\dots {}_{R_{\alpha_k}}\times_{{\rm ev}_{-}}
\widehat{\frak s^+}(\alpha_k,\alpha_+).
\end{equation}
This fiber product is well-defined because it is direct-like.
\end{enumerate}
\end{prop}
\begin{proof}
Using the results of previous subsections,
the proof is the same as the proof of Proposition \ref{prop161}.
\end{proof}
We next rewrite Proposition \ref{prop183939} in the algebraic language.
In our situation (linear K-system of Morse type)
we define
\begin{equation}
\Omega(R_{\alpha}) =  \bigoplus_{r \in \Omega(R_{\alpha})} \Q [r].
\end{equation}
We next define
\begin{equation}
\frak m^{\epsilon}_{1;\alpha_+,\alpha_-} :
\Omega(R_{\alpha_-}) \to \Omega(R_{\alpha_+})
\end{equation}
in the case when $\dim \mathcal M(\alpha_-,\alpha_+) = 0$
by the formula
\begin{equation}\label{form16ten6revrev}
\aligned
&\frak m^{\epsilon}_{1;\alpha_+,\alpha_-}([r_-]) \\
&=
\sum_{r_+ \in R_{\alpha_+}} [(ev_-,ev_+)^{-1}((r_-,r_+)) \\
&\qquad\qquad \cap
(\mathcal M(\alpha_-,\alpha_+)^{\boxplus\tau_0},
\widehat{\mathcal U^+}(\alpha_-,\alpha_+),\widehat{\frak O}(\alpha_-,\alpha_+),
\widehat{{\frak s^{+ \epsilon}}}(\alpha_-,\alpha_+))][r_+].
\endaligned
\end{equation}
Here the coefficient of $[r_+]$ in the right hand side is the virtual fundamental chain as in Definition \ref{defn1822}, which is a rational number, and $r_- \in R_{\alpha_-}$.
\par
We modify ($\flat$) appearing in Remark \ref{runningoutsharp} as follows.
\begin{enumerate}
\item[($\flat'$)]
\index{in the sense of ($\flat'$)}
For any energy cut level $E_0$ there exists $\epsilon_0(E_0)>0$ such that
the operator $\frak m^{\epsilon}_{1;\alpha_+,\alpha_-}$
is defined when $0 <E(\alpha_+) - E(\alpha_-) \le E_0$
and $\epsilon$ is in a dense open subset of $\{\epsilon \mid 0 < \epsilon < \epsilon_0(E_0)\}$.
\end{enumerate}
\begin{lem}\label{lem165rev3}
The operators $\frak m^{\epsilon}_{1;\alpha_+,\alpha_-}$ in (\ref{form16ten6revrev})
satisfy the following equality in the sense of $(\flat')$:
\begin{equation}
\sum_{\alpha; E(\alpha_-) < E(\alpha)
< E(\alpha_+)}
\frak m^{\epsilon}_{1;\alpha_+,\alpha}
\circ
\frak m^{\epsilon}_{1;\alpha,\alpha_-} = 0.
\end{equation}
\end{lem}
\begin{proof}
Using Propositions \ref{Stokesmulti18}, \ref{prop1826} and \ref{prop183939}, the proof is the same as the proof of
Lemma \ref{lem165}.
\end{proof}
We have thus rewritten Subsection \ref{subsec:constchaincpx}
by using multivalued perturbations
in the case of linear K-system of Morse type.
It is now obvious that we can rewrite Subsections \ref{subsec:constchaincmps}-\ref{subsec:proofcomplete}
in the same way and complete the proof of Theorem \ref{thm182222}.
\end{proof}

\section{Tree-like K-system\index{K-system ! tree-like K-system}: $A_{\infty}$ structure I: statement}
\label{sec:systemtree1}

In Sections \ref{sec:systemtree1}-\ref{sec:systemtree2} we discuss
construction of a filtered $A_{\infty}$ structure associated
to a relatively spin Lagrangian submanifold $L$
of a symplectic manifold.
This construction had been written in great detail in the article
\cite{fooo08},\cite{fooo09} based on singular homology.
Furthermore its de Rham version was given in
\cite[ Section 12]{fooo09}, \cite{foootoric3},
\cite{fooo091}.
\par
In this article, based on de Rham cohomology, we give its detail again.
We also provide a package so that the construction part of various
Kuranishi structures and their usage part to obtain a filtered $A_{\infty}$
structure are clearly separated from each other.
For this purpose, we take an axiomatic approach as in the case of
linear K-system developed up to the previous sections.
A similar axiomatic treatment was written in
\cite{fooo010}.
The axiom we give here is slightly different from
that of \cite{fooo010}.
In \cite{fooo010} we used a geometric operad
(the Stasheff operad) to formulate Kuranishi $A_{\infty}$ correspondence.
In this article
we use the {\it K-system
organized by a tree},
sometimes called
{\it tree-like K-system} in short, (plus certain additional data) rather than
the Stasheff operad.
In fact, in \cite{fooo010} the Stasheff operad was used only to
describe the combinatorial data on the way how various
strata (which are fiber products of moduli spaces of pseudo-holomorphic
disks) are glued. Since the members of the Stasheff operad are cells,
they do not carry a nontrivial homology class.
So in  \cite{fooo010} we actually did not use the geometric
data of the Stasheff operad but used only its combinatorial
structure, (that is, the way how various strata intersect).

\subsection{Axiom of tree-like K-system: $A_{\infty}$ correspondence}
\label{subsec:19-2}
In the rest of Part 2, we assume that
$L$ is a smooth oriented closed manifold.
\begin{shitu}\label{situ191}
Let $\frak G$ be an additive group and $\mu : \frak G \to \Z$,
$E : \frak G \to \R$ group homomorphisms.
We call $\mu(\beta)$ the
{\it Maslov index}\index{Maslov index ! on $\frak G$}
of $\beta$ and
$E(\beta)$ the {\it energy}\index{energy!on $\frak G$}
of $\beta$.
$\blacksquare$
\end{shitu}

\begin{defn}\label{defn192}
A {\it decorated rooted metric ribbon tree}
\index{decorated rooted metric ribbon tree}
\index{ribbon tree ! {\it see: decorated rooted metric ribbon tree}}
is $(\mathcal T,\beta(\cdot))$
such that:
\begin{enumerate}
\item $\mathcal T$ is a connected tree.
Let $C_0(\mathcal T)$, $C_1(\mathcal T)$ be the set of
all vertices and edges of $\mathcal T$, respectively.
\item
For each ${\rm v} \in C_0(\mathcal T)$
we fix a cyclic order
of the set of edges containing ${\rm v}$.
This is equivalent to fixing an isotopy type of an embedding
of $\mathcal T$ to the plane $\R^2$.
(Namely, the cyclic order of the edges is given by the orientation
of the plane so that the edges are enumerated according to the counter clockwise orientation.
We call it a {\it ribbon structure} at the vertex ${\rm v}$.)
\item
$C_0(\mathcal T)$ is divided into the set of {\it exterior vertices}
$C_{0,{\rm ext}}(\mathcal T)$ and
the set of {\it interior vertices} $C_{0,{\rm int}}(\mathcal T)$.
\item
We fix one element of $C_{0,{\rm ext}}(\mathcal T)$, which we call the {\it root}.
\item
The valency of all the exterior vertices are $1$.
\item
$\beta() : C_{0,{\rm int}}(\mathcal T) \to \frak G$ is a map.
We require $E(\beta({\rm v})) \ge 0$.
Moreover if $E(\beta({\rm v})) = 0$ then $\beta({\rm v})$ is
required to be the unit.
\item {\bf (Stability)}
For each ${\rm v} \in C_{0,{\rm int}}(\mathcal T)$
we assume that one of the following holds.
\begin{enumerate}
\item
$E(\beta({\rm v})) > 0$.
\item
The valency of ${\rm v}$ is not smaller than $3$.
\end{enumerate}
\end{enumerate}
We denote by $\mathcal G(k+1,\beta)$ the set of all
decorated ribbon trees $(\mathcal T,\beta(\cdot))$ such that:
\begin{enumerate}
\item[(I)]
$\# C_{0,{\rm ext}}(\mathcal T) = k+1$.
\item[(II)]
$\sum_{{\rm v} \in C_{0,{\rm int}}(\mathcal T)}(\beta({\rm v})) = \beta.$
\end{enumerate}
We decompose the set of edges $C_{1}(\mathcal T)$ as follows.
If an edge ${\rm e}$ contains an exterior vertex,
we call ${\rm e}$ an {\it exterior edge}. Otherwise we call
${\rm e}$ an {\it interior edge}.
We denote by $C_{1,{\rm int}}(\mathcal T)$,
(resp. $C_{1,{\rm ext}}(\mathcal T)$) the set of all interior (resp.
exterior) edges.
\end{defn}
Next we define fiber product of K-spaces along an element
$(\mathcal T,\beta(\cdot))$ of $\mathcal G(k+1,\beta)$.
Suppose we are given
K-spaces $\mathcal M_{k+1}(\beta)$
and maps
$$
{\rm ev} = ({\rm ev}_0,\dots,{\rm ev}_k)
: \overset{\circ}{\mathcal M}_{k+1}(\beta) \to L^{k+1}
$$
for each $\beta$ and $k \in \Z_{\ge 0}$.
Let
$(\mathcal T,\beta(\cdot)) \in\mathcal G(k+1,\beta)$.
We define the K-space
\begin{equation}\label{form191}
\prod_{(\mathcal T,\beta(\cdot))} {\mathcal M}_{k_{\rm v}+1}(\beta({\rm v}))
\end{equation}
as follows.
We consider the direct product
$
\prod_{{\rm v} \in C_{0,{\rm int}}(\mathcal T)}  {\mathcal M}_{k_{\rm v}+1}(\beta({\rm v})),
$
where $k_{\rm v}+1$ is the valency of the vertex ${\rm v}$.
We take two copies of $L$ for each interior
edge ${\rm e} \in C_{1,{\rm int}}(\mathcal T)$.
We define a map
\begin{equation}\label{edegeevaluationmap}
{\rm ev} :
\prod_{{\rm v} \in C_{0,{\rm int}}(\mathcal T)}  {\mathcal M}_{k_{\rm v}+1}(\beta({\rm v}))
\to \prod_{{\rm e} \in C_{1,{\rm int}}(\mathcal T)} L^2
\end{equation}
as follows.
For each ${\rm v} \in C_{0,{\rm int}}(\mathcal T)$ we enumerate the edges
containing
${\rm v}$ as
$$
{\rm e}_{{\rm v},0},\dots,{\rm e}_{{\rm v},k_{\rm v}}
$$
such that the following conditions are satisfied.
\begin{conds}
\begin{enumerate}
\item
The edge ${\rm e}_{{\rm v},0}$ is contained in the connected component of $\mathcal T \setminus \{{\rm v}\}$
which contains the root.
\item
$({\rm e}_{{\rm v},0},\dots,{\rm e}_{{\rm v},k_{\rm v}})$ respects the cyclic ordering
of the edges given by the ribbon structure at ${\rm v}$.
\end{enumerate}
\end{conds}
Such an enumeration is unique.
Each edge ${\rm e}$ contains two vertices.
For one of them ${\rm v}_-$ we have ${\rm e} = {\rm e}_{{\rm v}_-,0}$.
For the other vertex ${\rm v}_+$ contained in ${\rm e}$, we have ${\rm e} = {\rm e}_{{\rm v}_+,i}$ for some
$i\in \{1,\dots,k_{{\rm v}_+}\}$.
We define the ${\rm e}$ component of ${\rm ev}(({\bf x}_{\rm v})_{{\rm v} \in C_{0,{\rm int}}(\mathcal T)})$ as
$({\rm ev}_0( {\bf x}_{{\rm v}_-}),{\rm ev}_i( {\bf x}_{{\rm v}_+}))$,
where ${\rm e} = {\rm e}_{{\rm v}_+,i}$ and ${\bf x}_{\rm v}
\in \mathcal M_{k_{\rm v}+1}(\beta(\rm v))$. (Here ${\rm ev}$ is the map in (\ref{edegeevaluationmap}).)
\begin{defn}\label{defn1933}
The fiber product (\ref{form191}) is defined by
\begin{equation}\label{fiberproducttree}
\left(\prod_{{\rm v} \in C_{0,{\rm int}}(\mathcal T)}  {\mathcal M}_{k_{\rm v}+1}(\beta({\rm v}))\right)
\,\, {}_{{\rm ev}}\times_{\prod_{{\rm e} \in C_{1,{\rm int}}(\mathcal T)} L^2}
\left(\prod_{{\rm e} \in C_{1,{\rm int}}(\mathcal T)} L\right).
\end{equation}
Here
${\rm ev}$ is as in \eqref{edegeevaluationmap}
and
$\prod_{{\rm e} \in C_{1,{\rm int}}(\mathcal T)} L$ is the product of the
diagonal $L \subset L^2$ and is contained in
$\prod_{{\rm e} \in C_{1,{\rm int}}(\mathcal T)} L^2$.
We call (\ref{fiberproducttree}) {\it the fiber product of $\mathcal M_{k+1}(\beta)$
along $(\mathcal T,\beta(\cdot))$}.
\end{defn}
\begin{rem}
The fiber product (\ref{fiberproducttree})
in the sense of K-spaces may not be defined because of
the transversality problem.
It is defined if the following Condition \ref{zerossubmersive} is satisfied.
\end{rem}
\begin{conds}\label{zerossubmersive}
The map ${\rm ev}_0 : \mathcal M_{k+1}(\beta) \to L$ is weakly submersive.
\end{conds}

\begin{conds}\label{linAinfmainconds}
We consider the following objects.
\par\smallskip
\noindent {\bf (I)}
$\frak G$ is an additive group.
(We denote the unit $0 \in \frak G$ by $\beta_0$.)
$E : \frak G \to \R$ and $\mu : \frak G \to \Z$ are
group homomorphisms.
We call $E(\beta)$ the {\it energy} of $\beta$ and
$\mu(\beta)$ the {\it Maslov index} of $\beta$.
\par\smallskip
\noindent {\bf (II)}
$L$ is a smooth oriented manifold without boundary.
\par\smallskip
\noindent {\bf (III) (Moduli space)}
For each $\beta \in \frak G$ and $k\in \Z_{\ge 0}$ we have a
K-space with corners $\mathcal M_{k+1}(\beta)$
and a strongly smooth map
$$
{\rm ev} = ({\rm ev}_0,\dots,{\rm ev}_k)
: {\mathcal M}_{k+1}(\beta) \to L^{k+1}.
$$
We assume that
${\rm ev}_0$ is weakly submersive.
We call ${\mathcal M}_{k+1}(\beta)$
the {\it  moduli space of $A_{\infty}$ operations}
\index{$A_{\infty}$ operations ! moduli space of $A_{\infty}$ operations ${\mathcal M}_{k+1}(\beta)$}.
\par\smallskip
\noindent {\bf (IV) (Positivity of energy)}
We assume
${\mathcal M}_{k+1}(\beta) = \emptyset$
if $E(\beta) < 0$.
\par\smallskip
\noindent {\bf (V) (Energy zero part)}
In case $E(\beta) = 0$, we have $
{\mathcal M}_{k+1}(\beta) = \emptyset$ unless $\beta = 0$
and $k \ge 2$.
If $\beta = \beta_0 = 0$ then
${\mathcal M}_{k+1}(\beta_0) = L \times D^{k-2}$ and
${\rm ev}_i :  {\mathcal M}_{k+1}(\beta_0) \to L$ is the projection.
Here we regard $D^{k-2}$ as a Stasheff cell which is a
manifold with corners. (See \cite[Section 10]{foh}, for example.)
\par\smallskip
\noindent {\bf (VI) (Dimension)}
The dimension of the moduli space of $A_{\infty}$ operations is given by
\begin{equation}\label{eq:dimM}
\dim
{\mathcal M}_{k+1}(\beta)
=
\mu(\beta) + \dim L + k -2.
\end{equation}
\par\smallskip
\noindent {\bf (VII) (Orientation)}
${\mathcal M}_{k+1}(\beta)$ is oriented.
\par\smallskip
\noindent {\bf (VIII) (Gromov compactness)}
For any $E_0$ the set
\begin{equation}
\{\beta \in \frak G \mid \exists k\,\,
{\mathcal M}_{k+1}(\beta)
\ne \emptyset,
\,\, E(\beta) \le E_0\}
\end{equation}
is a finite set.
\par\smallskip
\noindent {\bf (IX) (Compatibility at the boundary)}
The normalized boundary of the moduli space of $A_{\infty}$ operations
is decomposed into the disjoint union of fiber products as follows.
\footnote{See Remark \ref{rem:FiberProdOrd} for the sign and the order of the fiber products.}
\begin{equation}\label{formula1660}
\partial {\mathcal M}_{k+1}(\beta)
\cong
\coprod_{\beta_1,\beta_2,k_1,k_2,i}
(-1)^{\epsilon}{\mathcal M}_{k_1+1}(\beta_1) \,\, {}_{{\rm ev}_i}\times_{{\rm ev}_0}
{\mathcal M}_{k_2+1}(\beta_2)
\end{equation}
where
\begin{equation}\label{eq:boundarysign}
\epsilon = (k_1 -1)(k_2 -1) + \dim L + k_1 +
(i-1)\Big(1+(\mu(\beta_2)+k_2) \dim L\Big)
\end{equation}
and
the union is taken over
$\beta_1,\beta_2,k_1,k_2,i$ such that $\beta_1 + \beta_2 = \beta$,
$k_1 + k_2 = k+1$, $i=1,\dots,k_2$ and $i$.
This isomorphism is compatible with
orientation and
is compatible with evaluation maps in the following sense.
Let ${\bf x}_1 \in {\mathcal M}_{k_1+1}(\beta_1)$ and ${\bf x}_2
 \in {\mathcal M}_{k_2+1}(\beta_2)$. We denote by
${\bf x}$ the element of $\partial {\mathcal M}_{k+1}(\beta)$
by the isomorphism (\ref{formula1660}).
Then
\begin{equation}\label{form1977}
{\rm ev}_j({\bf x})
=
\begin{cases}
{\rm ev}_j({\bf x}_1)   &\text{if $j =0,\dots, i-1$,} \\
{\rm ev}_{j-i+1}({\bf x}_2)   &\text{if $j = i,\dots,i+k_2-1$,} \\
{\rm ev}_{j-k_2+1}({\bf x}_1)   &\text{if $j = i+k_2,\dots,k$.}
\end{cases}
\end{equation}
See Remark \ref{rem:formula1660} below for
the sign \eqref{eq:boundarysign} and \eqref{form1977}.
\par
In case $\beta = 0$, we require that
(\ref{formula166}) coincides with the standard decomposition appearing
at the boundary of Stasheff cell. (See \cite[Section 10]{foh}.)
\par\smallskip
\noindent {\bf (X) (Compatibility at the corner I)}
Let $\widehat{S}_m( {\mathcal M}_{k+1}(\beta))$ be the
normalized corner of the K-space $ {\mathcal M}_{k+1}(\beta)$
in the sense of Definition \ref{norcor}.
Then it is isomorphic to the
disjoint union of
\begin{equation}\label{cornecomAinf1}
\prod_{(\mathcal T,\beta(\cdot))} {\mathcal M}_{k_{\rm v}+1}(\beta({\rm v})).
\end{equation}
Here the union is taken over all
$(\mathcal T,\beta(\cdot))
\in \mathcal G(k+1,\beta)$
such that $\#  C_{1,{\rm int}}(\mathcal T) = m$.
This isomorphism is compatible with the evaluation map
in the following sense.
\par
Let ${\rm v}_i$ be the $i$-th exterior vertex of $\mathcal T$.
The (unique) edge ${\rm e}$ containing ${\rm v}_i$ contains
one interior vertex denoted by ${\rm v}$.
Suppose ${\rm e}$ is the $j$-th edge of ${\rm v}$.
Let an element $({\bf x}_{\rm v})_{{\rm v}\in C_{0,{\rm int}}(\mathcal T)}$
of (\ref{cornecomAinf1}) correspond to an element ${\bf x}$
in $\widehat{S}_m( {\mathcal M}_{k+1}(\beta))$.
Then we require
\begin{equation}\label{form19999}
{\rm ev}_i({\bf x}) = {\rm ev}_j({\bf x}_{\rm v}).
\end{equation}
When $\beta = 0$, we require that
(\ref{cornecomAinf1}) coincides with the standard decomposition appearing
at the corner of the Stasheff cell.
\par\smallskip
\noindent {\bf (XI) (Compatibility at the corner II)}
Condition (X) implies that
$$
\widehat{S}_{\ell}(\widehat{S}_m( {\mathcal M}_{k+1}(\beta)))
$$
is a disjoint union of $(m+\ell)!/m!\ell!$ copies
of (\ref{cornecomAinf1}),
where the union is taken over all
$(\mathcal T,\beta(\cdot))
\in \mathcal G(k+1,\beta)$
such that $\#  C_{1,{\rm int}}(\mathcal T) = m+\ell$.
\par
The map
$\widehat{S}_{\ell}(\widehat{S}_m( {\mathcal M}_{k+1}(\beta)))
\to \widehat{S}_{m+\ell}( {\mathcal M}_{k+1}(\beta))$
is identified with the identity map
on each of the component (\ref{cornecomAinf1}).
\end{conds}
\begin{rem}\label{rem:formula1660}
The sign in \eqref{formula1660} is consistent with
our conventions adopted in \cite{fooobook2}.
\cite[Proposition 8.3.3]{fooobook2} is the same as
\eqref{formula1660} for the case $i=1$.
Also we can derive the sign in \eqref{formula1660} by using
\cite[Proposition 8.3.3]{fooobook2} and \cite[(8.4.5)]{fooobook2}.
Indeed, the formula \cite[(8.4.5)]{fooobook2} yields
\begin{equation}\label{eq:845}
{\mathcal M}_{k_1+1}(\beta_1) \,\, {}_{{\rm ev}_i}\times_{{\rm ev}_0}
{\mathcal M}_{k_2+1}(\beta_2)
=
(-1)^{\delta}{\mathcal M}_{k_1+1}(\beta_1) \,\, {}_{{\rm ev}_1}\times_{{\rm ev}_0}
{\mathcal M}_{k_2+1}(\beta_2)
\end{equation}
where
\begin{equation}\label{eq:845_1}
\aligned
\delta & = (i-1)(1+\dim L \dim {\mathcal M}_{k_2+1}(\beta_2) + \dim L) \\
& \equiv (i-1)\Big(1+(\mu (\beta_2)+k_2) \dim L\Big) \quad \mod 2,
\endaligned
\end{equation}
by taking the dimension formula \eqref{eq:dimM} into account.
Moreover,
the convention of the order of boundary marked points after gluing described in \cite[Remark 8.3.4]{fooobook2} is nothing but
the convention \eqref{form1977} for the case $i=1$.
\end{rem}
\begin{defn}\label{defn198}
A {\it tree-like K-system}\index{tree-like K-system ! {\it see: $A_{\infty}$ correspondence}}, or sometimes called
an {\it $A_{\infty}$ correspondence}
\index{$A_{\infty}$ correspondence ! $A_{\infty}$ correspondence}
over $L$, is a system of
(${\mathcal M}_{k+1}(\beta), {\rm ev}, \mu, E$)
satisfying Condition \ref{linAinfmainconds}.
\end{defn}
We next define a notion of {\it partial $A_{\infty}$ correspondence}.
The moduli space $\mathcal M_{k+1}(\beta)$ depends on $k$ and $\beta$.
In the version of partial $A_{\infty}$ correspondence which
we use in this article, we include the Kuranishi structure on $\mathcal M_{k+1}(\beta)$
for only a finite number of the pairs $(\beta,k)$.
This coincides with the way taken in \cite[Section 7]{fooobook2},
where we used the
notion of $A_{n,K}$ structure.
In \cite{fooo091} the notion of $A_{\infty}$ structure modulo $E_0$ was used.
It includes only a finite number of $\beta$'s but infinitely many $k$'s are
included. In \cite{fooo091} the Kuranishi structure of $\mathcal M_{k+1}(\beta)$
such that it is compatible with the forgetful map
$\mathcal M_{k+1}(\beta) \to \mathcal M_{1}(\beta)$ was used.
This is the reason
why infinitely many of $k$'s were included in \cite{fooo091}.
\par
Since we postpone technical detail concerning the forgetful map
to \cite{foootech3}, we use the formulation where only a finitely many
$k$'s are included in the partial structure.
Since we are working in de Rham theory, it is certainly possible
to include infinitely many $k$'s at this stage. However,
it seems easier to use only a finite number of
moduli spaces at each step of the construction.
\begin{defn}
A {\it partial $A_{\infty}$ correspondence of energy cut level $E_0$
and minimal energy $e_0$ over $L$}
\index{$A_{\infty}$ correspondence ! partial $A_{\infty}$ correspondence}
\index{partial !  $A_{\infty}$ correspondence}
is defined in the same way as $A_{\infty}$ correspondence except the following:
\begin{enumerate}
\item
The moduli space ${\mathcal M}_{k+1}(\beta)$ of $A_{\infty}$ operations
is defined only when $E(\beta) + ke_0 \le E_0$.
\item
The compatibility conditions that are Conditions \ref{linAinfmainconds}
(IX)(X)(XI) are assumed only when  $E(\beta) + ke_0\le E_0$.
\item
We assume that $\mathcal M_{k+1}(\beta) = \emptyset$ if $0 < E(\beta) < e_0$.
\end{enumerate}\end{defn}
Hereafter we say ${\mathcal M}_{k+1}(\beta)$
is a (partial) $A_{\infty}$ correspondence,
(and omit ${\rm ev}$, $\mu$, $E$) for simplicity.
\par
We next describe a parametrized version of Condition \ref{linAinfmainconds}.
\begin{conds}\label{linAinfmaincondspara}
We consider the following objects.
\par\smallskip
\noindent {\bf (I)}
$\frak G$ is an additive group.
(We denote the unit $0$ by $\beta_0$.)
$E : \frak G \to \R$ and $\mu : \frak G \to \Z$ are
group homomorphisms.
We call $E(\beta)$ the {\it energy} of $\beta$ and
$\mu(\beta)$ the {\it Maslov index} of $\beta$.
\par\smallskip
\noindent {\bf (II)}
$L$ is a smooth oriented manifold without boundary.
$P$ is a smooth oriented manifold with corners.
\par\smallskip
\noindent {\bf (III) (Moduli space)}
For each $\beta \in \frak G$ and $k\in \Z_{\ge 0}$ we have a
K-space with corners $\mathcal M_{k+1}(\beta;P)$
and a strongly smooth map
$${\rm ev} = ({\rm ev}_P, {\rm ev}_0,\dots,{\rm ev}_k)
: {\mathcal M}_{k+1}(\beta;P) \to P\times L^{k+1}.
$$
We assume that
$({\rm ev}_P,{\rm ev}_0)$ is weakly submersive stratumwisely.
We call ${\mathcal M}_{k+1}(\beta;P)$
the {\it  moduli space of $P$-parametrized $A_{\infty}$ operations}.
\index{$A_{\infty}$ operations ! moduli space of $P$-parametrized $A_{\infty}$ operations ${\mathcal M}_{k+1}(\beta;P)$}
\par\smallskip
\noindent {\bf (IV) (Positivity of energy)}
We assume  $
{\mathcal M}_{k+1}(\beta;P) = \emptyset$
if $E(\beta) < 0$.
\par\smallskip
\noindent {\bf (V) (Energy zero part)}
In case $E(\beta) = 0$, we have $
{\mathcal M}_{k+1}(\beta;P) = \emptyset$ unless $\beta = 0$
and $k \ge 2$.  If $\beta = \beta_0 = 0$, then
${\mathcal M}_{k+1}(\beta_0;P) = P\times L \times D^{k-2}$ and
${\rm ev}_i :  {\mathcal M}_{k+1}(\beta_0;P) \to L$ is the projection.
Also ${\rm ev}_P :  {\mathcal M}_{k+1}(\beta_0;P) \to P$ is the projection.
Here we again identify $D^{k-2}$ with the Stasheff cell.
\par\smallskip
\noindent {\bf (VI) (Dimension)}
The dimension of the moduli space of $P$-parametrized $A_{\infty}$ operations is given by
\begin{equation}
\dim
{\mathcal M}_{k+1}(\beta;P)
=
\mu(\beta) + \dim L + k -2+\dim P.
\end{equation}
\par\smallskip
\noindent {\bf (VII) (Orientation)}
${\mathcal M}_{k+1}(\beta;P)$ is oriented.
\par\smallskip
\noindent {\bf (VIII) (Gromov compactness)}
For any $E_0$ the set
\begin{equation}
\{\beta \in \frak G \mid \exists k\,\,
{\mathcal M}_{k+1}(\beta;P)
\ne \emptyset,
\,\, E(\beta) \le E_0\}
\end{equation}
is a finite set.
\par\smallskip
\noindent {\bf(IX) (Compatibility at the boundary)}
The  normalized boundary of the moduli space of $A_{\infty}$ operations
is decomposed into the fiber products as follows.
\footnote{Following our convention in \cite[(8.9.1)]{fooobook2},
that the parameter space $P$ is put on the first factor in the fiber product.
So there is no extra sign contribution from the parameter space in the second line on
the right hand side of \eqref{formula166}.}
\begin{equation}\label{formula166}
\aligned
&\partial {\mathcal M}_{k+1}(\beta;P)\\
\cong
&\coprod_{\beta_1,\beta_2,k_1,k_2,i}
(-1)^{\epsilon}
{\mathcal M}_{k_1+1}(\beta_1;P) \,\, {}_{({\rm ev}_P,{\rm ev}_i)}
\times_{({\rm ev}_P,{\rm ev}_0)}
{\mathcal M}_{k_2+1}(\beta_2;P) \\
& \sqcup \left(\partial P \,{}_P\times_{{\rm ev}_P} {\mathcal M}_{k+1}(\beta;P)\right).
\endaligned
\end{equation}
where
\begin{equation}\label{eq:boundarysign2}
\epsilon = (k_1 -1)(k_2 -1) + \dim L + k_1 +
(i-1)\Big(1+(\mu(\beta_2)+k_2 + \dim P) \dim L\Big)
\end{equation}
and the union in the second line is taken over
$\beta_1,\beta_2,k_1,k_2,i$ such that $\beta_1 + \beta_2 = \beta$,
$k_1 + k_2 = k+1$, $i=1,\dots,k_2$.
The fiber product in the second line is taken over $P \times L$.
This isomorphism is compatible with
orientation.
It is compatible with evaluation maps,
in the same sense as (\ref{form1977}).
In case $\beta = 0$, we require that
(\ref{formula166}) coincides with the
decomposition induced from the standard decomposition appearing
at the boundary of Stasheff cell.
\par\smallskip
\noindent {\bf (X) (Compatibility at the corner I)}
Let $\widehat{S}_m( {\mathcal M}_{k+1}(\beta;P))$ be the
normalized corner of the K-space $ {\mathcal M}_{k+1}(\beta;P)$
in the sense of Definition \ref{norcor}.
Then it is isomorphic to the
disjoint union of
\begin{equation}\label{cornecomAinf122}
\prod_{(\mathcal T,\beta(\cdot))} {\mathcal M}_{k_{\rm v}+1}(\beta({\rm v});
\widehat{S}_{m'}(P)).
\end{equation}
(We will explain the fiber product (\ref{cornecomAinf122})
right after Condition \ref{linAinfmaincondspara}.)
Here the union is taken over all
$(\mathcal T,\beta(\cdot))
\in \mathcal G(k+1,\beta)$ and $m' \in \Z_{\ge 0}$
such that $\#  C_{1,{\rm int}}(\mathcal T)
+ m'= m$ and
we put
$$
 {\mathcal M}_{k_{\rm v}+1}(\beta({\rm v});
\widehat{S}_{m'}(P))
=
 \widehat{S}_{m'}(P) \,\,_P\times_{{\rm ev}_P} {\mathcal M}_{k_{\rm v}+1}(\beta({\rm v});
P).
$$
This isomorphism is compatible with the evaluation map
in the same sense as (\ref{form19999}).
In case $\beta = 0$, we require that
(\ref{cornecomAinf122}) coincides with the decomposition induced from the standard decomposition appearing
at the corner of the Stasheff cell.
\par\smallskip
\noindent {\bf (XI) (Compatibility at the corner II)}
Condition (X) implies that
$\widehat{S}_{\ell}(\widehat{S}_m( {\mathcal M}_{k+1}(\beta;P)))$
is a disjoint union of copies
\begin{equation}\label{cornecomAinf122333}
\prod_{(\mathcal T,\beta(\cdot))} {\mathcal M}_{k_{\rm v}+1}(\beta({\rm v});
\widehat{S}_{\ell'}(\widehat{S}_{m'}(P)))
\end{equation}
where the union is taken over all
$(\mathcal T,\beta(\cdot))
\in \mathcal G(k+1,\beta)$,
and $m',\ell'$
such that $\#  C_{1,{\rm int}}(\mathcal T) +
m' + \ell'= m+\ell$
and
we put
$$
 {\mathcal M}_{k_{\rm v}+1}(\beta({\rm v});
\widehat{S}_{\ell'}(\widehat{S}_{m'}(P)))
=
 \widehat{S}_{\ell'}(\widehat{S}_{m'}(P)) \,\,_P\times_{{\rm ev}_P}
 {\mathcal M}_{k_{\rm v}+1}(\beta({\rm v});
P).
$$
\par
The covering map
$\widehat{S}_{\ell}(\widehat{S}_m( {\mathcal M}_{k+1}(\beta;P)))
\to \widehat{S}_{m+\ell}( {\mathcal M}_{k+1}(\beta;P))$
is identified with the map induced from the covering map
$\widehat{S}_{\ell'}(\widehat{S}_{m'}(P))
\to \widehat{S}_{\ell'+m'}(P)$
on each of the component (\ref{cornecomAinf122333}).
\end{conds}
Now we define the fiber product (\ref{cornecomAinf122}) as follows.
We consider the direct product
$\prod_{{\rm v} \in C_{0,{\rm int}}(\mathcal T)}  {\mathcal M}_{k_{\rm v}+1}(\beta({\rm v});P)$
as in (\ref{edegeevaluationmap}). We have
$$
{\rm ev} :
\prod_{{\rm v} \in C_{0,{\rm int}}(\mathcal T)}  {\mathcal M}_{k_{\rm v}+1}(\beta({\rm v});P)
\to \prod_{{\rm e} \in C_{1,{\rm int}}(\mathcal T)} L^2.
$$
Using ${\rm ev}_P :  {\mathcal M}_{k_{\rm v}+1}(\beta({\rm v});P) \to P$ we have
$$
{\rm ev}_P :
\prod_{{\rm v} \in C_{0,{\rm int}}(\mathcal T)}  {\mathcal M}_{k_{\rm v}+1}(\beta({\rm v});P)
\to \prod_{{\rm e} \in C_{1,{\rm int}}(\mathcal T)} P^2.
$$
The fiber product (\ref{cornecomAinf122}) is by definition
\begin{equation}\label{fiberproducttreerev}
\aligned
& \left(\prod_{{\rm e} \in C_{1,{\rm int}}(\mathcal T)} P \times L\right)\\
&\,\, {}_{\prod_{{\rm e} \in C_{1,{\rm int}}(\mathcal T)} (P^2 \times L^2)}
\times_{({\rm ev}_P,{\rm ev})}
\left(\prod_{{\rm v} \in C_{0,{\rm int}}(\mathcal T)}  {\mathcal M}_{k_{\rm v}+1}(\beta({\rm v});P)\right).
\endaligned
\end{equation}
\begin{lem}
The fiber product (\ref{fiberproducttreerev}) is well-defined.
\end{lem}
\begin{proof}
We assumed that the map
$$
({\rm ev}_P,{\rm ev}_0) : \mathcal M_{k_{\rm v}+1}(\beta({\rm v});P) \to P \times L
$$
is weakly submersive stratumwisely. We also note that for each ${\rm e} \in C_{1,{\rm int}}(\mathcal T)$ there exists a unique vertex
${\rm v}$ such that ${\rm e}$ is its $0$-th vertex. These two facts imply the lemma immediately.
\end{proof}
\begin{defn}\label{defn2013}
A {\it $P$-parametrized $A_{\infty}$ correspondence}
\index{$A_{\infty}$ correspondence ! $P$-parametrized $A_{\infty}$ correspondence} over $L$
is a system of
(${\mathcal M}_{k+1}(\beta;P), {\rm ev}, \mu, E$)
satisfying Condition \ref{linAinfmaincondspara}.
\par
A {\it partial $P$-parametrized $A_{\infty}$ correspondence of energy cut level $E_0$
and minimal energy $e_0$ over $L$}
\index{$A_{\infty}$ correspondence ! partial $P$-parametrized $A_{\infty}$ correspondence}
\index{partial !  $P$-parametrized $A_{\infty}$ correspondence}
is defined in the same way as $P$-parametrized $A_{\infty}$ correspondence except the following:
\begin{enumerate}
\item
The moduli space of $P$-parametrized $A_{\infty}$ operations ${\mathcal M}_{k+1}(\beta;P)$
is defined only when $E(\beta) + ke_0 \le E_0$.
\item
The compatibility conditions that is Condition \ref{linAinfmaincondspara}
(IX)(X)(XI) are assumed only when  $E(\beta) + ke_0\le E_0$.
\item
We assume $\mathcal M_{k+1}(\beta,P) = \emptyset$ if $0 < E(\beta)  < e_0$.
\end{enumerate}
\end{defn}
\begin{lem}
Let ${\mathcal M}_{k+1}(\beta;P)$ be a $P$-parametrized $A_{\infty}$ correspondence
over $L$ and $\partial_iP$ a connected component of
the normalized boundary of $P$.
Then
$$
{\mathcal M}_{k+1}(\beta;\partial_i P)
=
\partial_i P \,{}_P \times_{{\rm ev}_P} {\mathcal M}_{k+1}(\beta;P)
$$
defines a $\partial_iP$-parametrized $A_{\infty}$ correspondence
over $L$.
The same holds for partial $P$-parametrized $A_{\infty}$ correspondence.
\end{lem}
The proof is obvious.
\begin{defn}\label{defn1913piso}
Suppose we are given two $A_{\infty}$ correspondences
over $L$ denoted by
${\mathcal M}^j_{k+1}(\beta)$ with $j=1,2$.
Then a {\it pseudo-isotopy}
\index{$A_{\infty}$ correspondence ! pseudo-isotopy between them}
between them is a
$P =[1,2]$ parametrized $A_{\infty}$ correspondence
${\mathcal M}_{k+1}(\beta;[1,2])$ such that
for $\{1\} \subset \partial [1,2]$
(resp. $\{2\} \subset \partial [1,2]$)
the  $A_{\infty}$ correspondence
${\mathcal M}_{k+1}(\beta;\{1\})$
(resp. ${\mathcal M}_{k+1}(\beta;\{2\})$)
is isomorphic to ${\mathcal M}^1_{k+1}(\beta)$
(resp ${\mathcal M}^2_{k+1}(\beta)$).
We define the notion of pseudo-isotopy of partial
$A_{\infty}$ correspondences in the same way.
\end{defn}

\begin{rem}
In this article we use the notion of pseudo-isotopy of
$A_{\infty}$ correspondences to
prove well-definedness of the filtered $A_{\infty}$
algebra induced by the $A_{\infty}$ correspondence.
(See Theorem \ref{theorem1934} (2).)
On the other hand,
we can define the notion of morphism of $A_{\infty}$ correspondences
and use it instead to prove well-definedness.
In other words, we are using the bifurcation method here
but not the cobordism method. (See \cite[Subsection 7.2.14]{fooo09} for these two methods.
In a slightly different formulation, a morphism of $A_{\infty}$ correspondences
is defined in\cite[Definition 8]{fooo010}.)
In \cite{fooo09} and \cite{fooo010}
we used the cobordism method.
In \cite{akahojoyce} and \cite{fooo091}
the bifurcation method was used. They will
give the same morphism
at the end of the day as explained in \cite[Remark 12.3]{fooo091}.
\end{rem}
\begin{defn}\label{inducAinfty}
An {\it inductive system of $A_{\infty}$ correspondences over $L$}
\index{$A_{\infty}$ correspondence ! inductive system of}
consists of the following objects.
\begin{enumerate}
\item
We are given a sequence $\{E^i\}_{i=1}^{\infty}$ of positive real numbers
such that $E^i < E^{i+1}$ and $\lim_{i\to \infty} E^i = \infty$.
We are also given $e_0 > 0$ independent of $i$.
\item
For each $i$, we are given a
partial
$A_{\infty}$ correspondence
${\mathcal M}^i_{k+1}(\beta)$ over $L$ of energy cut level $E^i$
and minimal energy $e_0$.
\item
For each $i$ we are given a pseudo-isotopy
${\mathcal M}_{k+1}(\beta;[i,i+1])$
between ${\mathcal M}^i_{k+1}(\beta)$ and ${\mathcal M}^{i+1}_{k+1}(\beta)$.
Here the energy cut level of ${\mathcal M}_{k+1}(\beta;[i,i+1])$
is $E^{i}$ and its minimal energy is $e_0$.
\item
We assume the following uniform Gromov compactness.
For each $E'>0$ the next set is of finite order.
\begin{equation}
\{\beta \in \frak G \mid \exists k \,\exists i\,\,
{\mathcal M}^i_{k+1}(\beta)
\ne \emptyset,
\,\, E(\beta) \le E'\}.
\end{equation}
The next set is also of finite order for each $E' >0$:
\begin{equation}
\{\beta \in \frak G \mid \exists k \,\exists i\,\,
{\mathcal M}^i_{k+1}(\beta;[i,i+1])
\ne \emptyset,
\,\, E(\beta) \le E'\}.
\end{equation}
\end{enumerate}
\end{defn}
\begin{rem}\label{rem1916}
In the situation of Definition \ref{inducAinfty},
suppose $E^{\prime i}$ is another sequence
so that $E^{\prime i} < E^{\prime i+1}$,
$\lim_{i\to \infty} E^{\prime i} = \infty$
and $E^{\prime i} < E^i$.
We forget the K-spaces
${\mathcal M}_{k+1}(\beta;[i,i+1])$ and
${\mathcal M}^i_{k+1}(\beta + ke_0)$ for
$E(\beta) > E^{\prime i}$.
Then we obtain another inductive system of $A_{\infty}$ correspondences over $L$.
\par
Let $e'_0 < e_0$. We replace $E^i$ by $E^{\prime i} = E^i e'_0/e_0$.
Then $ke'_0 + E \le E^{\prime i}$ implies $ke_0 + E \le E^i$.
Therefore any partial $A_{\infty}$ correspondence of energy
cut level $E^i$ and of minimal energy $e_0$ induces
one of energy cut level $E^{\prime i}$ and of minimal energy $e'_0$
by forgetting certain moduli spaces.
In this way when we compare two
inductive systems of $A_{\infty}$ correspondences over $L$
we may always assume that the numbers $E^i$, $e_0$ are common
without loss of generality.
We will assume it in the next definition.
\end{rem}
\begin{defn}\label{definition1918}
Suppose for $j=0,1$ we are given inductive systems
of $A_{\infty}$ correspondences over $L$, denoted by
${\mathcal M}^{ji}_{k+1}(\beta)$,
${\mathcal M}^j_{k+1}(\beta;[i,i+1])$.
(We take the same $E^i$ and $e_0$ for $j = 0,1$ as we explained in
Remark \ref{rem1916}.)
A {\it pseudo-isotopy between these two inductive systems}
\index{$A_{\infty}$ correspondence ! pseudo-isotopy between two inductive systems} consists
of the following objects.
\begin{enumerate}
\item
For each $i$ we are given a pseudo-isotopy
${\mathcal M}^i_{k+1}(\beta;[0,1])$ between
${\mathcal M}^{0i}_{k+1}(\beta)$ and
${\mathcal M}^{1i}_{k+1}(\beta)$.
The energy cut level and minimal energy of this pseudo-isotopy are
$E^i$ and $e_0$, respectively.
\item
For each $i$ we are given a $P = [0,1] \times [i,i+1]$ parametrized $A_{\infty}$ correspondence
${\mathcal M}_{k+1}(\beta;[0,1]\times [i,i+1])$ satisfying
the following properties.
\begin{enumerate}
\item
Its restriction to the boundary component
$\{j\} \times [i,i+1]$ is isomorphic to ${\mathcal M}^j_{k+1}(\beta;[i,i+1])$.
Here $j=0,1$.
\item
Its restriction to the boundary component
$[0,1]\times \{i\}$ is isomorphic to ${\mathcal M}^{1i}_{k+1}(\beta)$.
\item
Its restriction to the boundary component
$[0,1]\times \{i+1\}$ is isomorphic to ${\mathcal M}^{1i+1}_{k+1}(\beta)$.
\item
The isomorphisms in (a)(b)(c) are
consistent at
${\mathcal M}_{k+1}(\beta;\widehat S_2([0,1]\times [i,i+1]))$.
\item
The energy cut level of ${\mathcal M}_{k+1}(\beta;[0,1]\times [i,i+1])$
is $E^i$ and its minimal energy is $e_0$.
\item
The isomorphisms in Item (a)(b)(c)(d) satisfy appropriate corner
compatibility conditions.
\end{enumerate}
\item
We assume the following uniform Gromov compactness.
For each $E'>0$ the next set is of finite order.
\begin{equation}
\{\beta \in \frak G \mid \exists k ~\exists i\,\,
{\mathcal M}_{k+1}(\beta; [0,1] \times [i,i+1])
\ne \emptyset,
\,\, E(\beta) \le E'\}.
\end{equation}
\end{enumerate}
\end{defn}

\subsection{Filtered $A_{\infty}$ algebra and its pseudo-isotopy}
\label{subsec:19-3}

In this subsection we review certain algebraic material in
\cite{fooobook}, \cite{fooobook2}, \cite{fooo091}.

\begin{defn}\label{discmonoid}
We call a subset $G \subset \R_{\ge 0}\times 2\Z$ a
{\it discrete submonoid}\footnote{In the situation of oriented Lagrangian submanifolds, the Maslov index is even. Thus we assume
$G \subset \R_{\ge 0}\times 2\Z$, while we consider $G \subset \R_{\ge 0}\times \Z$
in Definition \ref{defn141212}.}
if the following holds.
We denote by $E : G  \to \R_{\ge 0}$ and $\mu : G  \to 2\Z$ the natural projections.
\begin{enumerate}
\item If $\beta_1,\beta_2 \in  G$, then $\beta_1 + \beta_2 \in G$.
$(0,0) \in G$.
\item The image $E(G) \subset \R_{\ge 0}$ is discrete.
\item For each $E_0 \in \R_{\ge 0}$ the inverse image
$G \cap E^{-1}([0,E_0])$ is a finite set.
\end{enumerate}
\end{defn}
Let $\Omega(L)$ be the de Rham complex of $L$.
We put
$\Omega(L)[1]^d = \Omega^{d+1}(L)$, where $\Omega^d(L)$
is the space of degree $d$ smooth forms.
We put
\begin{equation}
B_k(\Omega(L)[1])
=
 \underbrace{\Omega(L)[1] \otimes \dots \otimes \Omega(L)[1]}
_{\text{$k$ times}}.
\end{equation}
Let $G$ be a discrete submonoid as in Definition \ref{discmonoid}.
\begin{defn}\label{filAinfonOmega}{\rm (\cite[Definition 3.2.26, Definition 3.5.6, Remark 3.5.8]{fooobook})}
A {\it $G$-gapped filtered $A_{\infty}$ algebra structure}
\index{$A_{\infty}$ algebra ! $G$-gapped filtered $A_{\infty}$ algebra}
 on $\Omega(L)$ is a
sequence of multilinear maps
\begin{equation}\label{mkbeta}
\frak m_{k,\beta} : B_k(\Omega(L)[1]) \to \Omega(L)[1]
\end{equation}
for each $\beta\in G$ and $k\in \Z_{\ge 0}$ of degree $1-\mu(\beta)$
with the following properties:
\smallskip
\begin{enumerate}
\item
$\frak m_{k,\beta_0} = 0$ for $\beta_0 = (0,0)$, $k\ne 1,2$.
\item
$\frak m_{2,\beta_0}(h_1,h_2) = (-1)^* h_1 \wedge h_2$,
where $* = \deg h_1(\deg h_2+1)$.\footnote{See \cite[(3.5.9)]{fooobook} for the sign.}
\item
$\frak m_{1,\beta_0}(h) = (-1)^* dh$, where $* = n+1+\deg h$.
\footnote{The sign here is different from
the one given in \cite[(3.2.5)]{fooobook}.
The latter is a convention when we regard a DGA as an $A_{\infty}$
algebra with trivial higher multiplications. (Actually it does not matter if we change the sign of differential in DGA.) But the signs here and in \cite[Remark 3.5.8]{fooobook} provide the signs
in the filtered $A_{\infty}$ algebra which are compatible with the signs in its unfiltered DGA.}
\item
\begin{equation}\label{Ainfinityrelbeta}
\aligned
\sum_{k_1+k_2=k+1}&\sum_{\beta_1+\beta_2=\beta}\sum_{i=1}^{k-k_2+1} \\
&(-1)^*{\frak m}_{k_1,\beta_1}(h_1,\ldots,{\frak m}_{k_2,\beta_2}(h_i,\ldots,h_{i+k_2-1}),\ldots,h_{k}) = 0
\endaligned\end{equation}
holds for any $\beta \in G$ and $k$. The sign is given by
$
* = \deg' h_1 + \dots + \deg' h_{i-1}.
$
Here $\deg'$ is the shifted degree by $+1$.
\end{enumerate}
\end{defn}
\begin{defn}\label{defn:partial}
A {\it partial G-gapped filtered $A_{\infty}$ algebra structure
of energy cut level $E$ and minimal energy $e_0$}
\index{$A_{\infty}$ algebra ! partial G-gapped filtered $A_{\infty}$ algebra structure}
is a sequence of operators (\ref{mkbeta}) for $E(\beta) + ke_0\le E$
satisfying the same properties, except (\ref{Ainfinityrelbeta}) is
assumed only for $E(\beta) + ke_0\le E$.
We require $\frak m_{k,\beta} = 0$ if $0 < E(\beta) < e_0$.
\end{defn}
\begin{defn}
A multilinear map
$
F : B_k(\Omega(L)[1]) \to \Omega(L)[1]
$
is said to be {\it continuous in $C^{\infty}$ topology}
if the following holds.
Suppose $h_{i,a} \in \Omega(L)$ converges to $h_i \in \Omega(L)$
in  $C^{\infty}$ topology as $a \to \infty$, then
$F(h_{1,a},\dots,h_{k,a})$ converges to $F(h_{1},\dots,h_{k})$
in $C^{\infty}$ topology as $a \to \infty$.
\end{defn}
\begin{defn}
A (partial) $G$ gapped filtered
$A_{\infty}$ algebra structure on $\Omega(L)$ is said
to be {\it continuous}
\index{$A_{\infty}$ algebra ! continuous} if the operations
$\frak m_{k,\beta}$
is continuous in $C^{\infty}$ topology.
Hereafter we assume that all the
operations of (partial)
$A_{\infty}$ algebra structure on $\Omega(L)$
are continuous in $C^{\infty}$ topology.
\end{defn}
\begin{defn}\label{pisotopydef}(\cite[Definition 8.5]{fooo091})
For each $t\in [0,1]$,
$\beta\in G$ and $k\in \Z_{\ge 0}$,
let $\frak m_{k,\beta}^t$ be as in
Definition \ref{filAinfonOmega}
and
$\frak c_{k,\beta}^t$
a sequence of multilinear maps
\begin{equation}\label{ckbeta}
\frak c_{k,\beta}^t : B_k(\Omega(L)[1]) \to \Omega(L)[1]
\end{equation}
of degree $-\mu(\beta)$.
We say $(\{\frak m^t_{k,\beta}\},\{\frak c^t_{k,\beta}\})$
is a {\it pseudo-isotopy of $G$-gapped filtered $A_{\infty}$ algebra
structures}, or in short, {\it gapped pseudo-isotopy} on $\Omega(L)$ if the following holds:
\index{$A_{\infty}$ algebra ! pseudo-isotopy}
\index{pseudo-isotopy ! of $G$-gapped filtered $A_{\infty}$ algebras}
\index{pseudo-isotopy ! gapped}
\begin{enumerate}
\item  $\frak m^t_{k,\beta}$ and $\frak c^t_{k,\beta}$ are
continuous in $C^{\infty}$ topology.
The map sending $t$ to $\frak m^t_{k,\beta}$ or $\frak c^t_{k,\beta}$
is smooth. Here we use the operator topology with respect to the $C^{\infty}$
topology for $\frak m^t_{k,\beta}$ and $\frak c^t_{k,\beta}$ to define
their smoothness.
\item For each (but fixed) $t$, the set of operators $\{\frak m^t_{k,\beta}\}$ defines a $G$-gapped filtered $A_{\infty}$ algebra
structure on $\Omega(L)$.
\item  For each $h_i \in \Omega(L)[1]$ the following equality holds:
\begin{equation}\label{isotopymaineq}
\aligned
&\frac{d}{dt} \frak m_{k,\beta}^t(h_1,\ldots,h_k) \\
&+ \sum_{k_1+k_2=k+1}\sum_{\beta_1+\beta_2=\beta}\sum_{i=1}^{k-k_2+1}
(-1)^{*}\frak c^t_{k_1,\beta_1}(h_1,\ldots, \frak m_{k_2,\beta_2}^t(h_i,\ldots),\ldots,h_k) \\
&- \sum_{k_1+k_2=k+1}\sum_{\beta_1+\beta_2=\beta}\sum_{i=1}^{k-k_2+1}
\frak m^t_{k_1,\beta_1}(h_1,\ldots, \frak c_{k_2,\beta_2}^t(h_i,\ldots),\ldots,h_k)\\
&=0.
\endaligned
\end{equation}
Here $* = \deg' h_1 + \dots + \deg'h_{i-1}$.
\item
$\frak c_{k,\beta}^t = 0$ if $E(\beta) \le 0$.
\end{enumerate}
\par
Sometimes we say that $(\{\frak m^t_{k,\beta}\},\{\frak c^t_{k,\beta}\})$
is a pseudo-isotopy or {\it gapped pseudo-isotopy} between
$\{\frak m^0_{k,\beta}\}$ and $\{\frak m^1_{k,\beta}\}$.
\index{gapped ! pseudo-isotopy}
\end{defn}
\begin{defn}\label{pisotopydef2}
We define the notion of {\it pseudo-isotopy of partial
$G$-gapped filtered $A_{\infty}$ algebra
structures on $\Omega(L)$ of energy cut level $E$ and
minimal energy $e_0$}
\index{$A_{\infty}$ algebra ! pseudo-isotopy of partial
$G$-gapped filtered $A_{\infty}$ algebra}
in the same way,
except we replace (2) and (3) by the following (2)' and (3)'
and we further require (5) below.
\begin{enumerate}
\item[(2)']
We require that
the set of operators $\{\frak m^t_{k,\beta}\}$ defines a
partial $G$-gapped filtered $A_{\infty}$ algebra
structures on $\Omega(L)$ of energy cut level $E$ and of
minimal energy $e_0$.
\item[(3)']
We require (\ref{isotopymaineq}) only for $\beta, k$
with $E(\beta) + ke_0\le E$.
\item[(5)]
$\frak m^t_{k,\beta} = \frak c^t_{k,\beta} = 0$
if $0 < E(\beta) < e_0$.
\end{enumerate}
\end{defn}
We also use the notion of pesudo-isotopy of pseudo-isotopies.
It seems simpler to define a more general notion
of the $P$-parametrized family of $G$-gapped filtered $A_{\infty}$ algebra
structures on $\Omega(L)$ in general.
So we define this notion below.
\par
We consider the totality of smooth differential forms
on $P\times L$ which we write $\Omega(P\times L)$.
Here $\Omega(P\times L)[1]$ is its degree shift as before.
Let $t_1,\dots,t_d$ be local coordinates of $P$.
(In this article we only consider the case $P \subset \R^d$ for
$d= \dim P$. So we actually have canonical global coordinates.)
For $h \in \Omega(P\times L)$ we put
\begin{equation}
h = \sum_{I \subset \{1,\dots,d\}} dt_I \wedge h_I.
\end{equation}
Here $I=\{i_1,\dots,i_{\vert I\vert}\}$
($i_1 < \dots < i_{\vert I\vert}$),
$dt_I = dt_{i_1} \wedge \dots \wedge dt_{i_{\vert I\vert}}$
and $h_I$ does not contain $dt_i$.

\begin{defn}\label{defn1926}
A multilinear map
$
F : B_k(\Omega(P\times L)[1]) \to \Omega(P\times L)[1]
$
is said to be
{\it pointwise in $P$ direction}
\index{pointwise in $P$ direction}
if the following holds:
\par
For each $I \subset \{1,\dots,d\}$
and ${\bf t} \in P$ there exists a continuous map
$$
F^{\bf t}_{I;j_1,\dots,j_k} : B_k(\Omega(L)[1]) \to \Omega(L)[1]
$$
such that
\begin{equation}
\aligned
&F(dt_{j_1} \wedge h_1, \dots, dt_{j_k} \wedge h_k)\vert_{\{{\bf t}\} \times L} \\
&=
\sum_{I} dt_I \wedge dt_{j_1} \wedge \dots \wedge dt_{j_k}
\wedge F^{\bf t}_{I;j_1,\dots,j_k}(h^{\bf t}_1,\dots,h^{\bf t}_k),
\endaligned
\end{equation}
where $\vert_{\{{\bf t}\} \times L}$
means the restriction to $\{{\bf t}\} \times L$. Moreover $F^{\bf t}_{I;j_1,\dots,j_k}$
depends smoothly on ${\bf t}$ with respect to the operator topology.
Here $h^{\bf t}_i$ is the restriction of $h_i$ to $\{{\bf t}\} \times L$.
\end{defn}
\begin{rem}
This condition is equivalent to the following one.
\begin{enumerate}
\item[(*)]
For smooth differential forms $\sigma_i$ on $P$, we have
$$
F(\sigma_1 h_1,\dots, \sigma_k h_k)
=
\pm
\sigma_1 \wedge \dots \wedge \sigma_k \wedge
F(h_1,\dots, h_k).
$$
\end{enumerate}
\end{rem}
\begin{defn}\label{PparaAinfdef}
A {\it $P$-parametrized family of $G$-gapped filtered $A_{\infty}$ algebra
structures on $\Omega(L)$}
\index{$A_{\infty}$ algebra ! $P$-parametrized family of $G$-gapped filtered $A_{\infty}$ algebra}
is $\{\frak m^P_{k,\beta}\}$
satisfying the following properties:
\begin{enumerate}
\item
$
\frak m^P_{k,\beta} : B_k(\Omega(P\times L)[1]) \to \Omega(P\times L)[1]
$
is a multilinear map of degree $1$.
\item
$\frak m^P_{k,\beta}$ is
pointwise in $P$ direction if $\beta \ne \beta_0$.
\item
$\frak m^P_{k,\beta_0} = 0$ for  $k \ne 1,2$.
\item
$\frak m^P_{1,\beta_0}(h) = (-1)^*dh$. Here $d$ is the de Rham differential
and $* = n+1+\deg h$.
\item
$\frak m^P_{2,\beta_0}(h_1 \wedge h_2) = (-1)^* h_1 \wedge h_2$.
Here $\wedge$ is the wedge product and $* = \deg h_1(\deg h_2+1)$.
\item
$\frak m^P_{k,\beta}$ satisfies the following $A_{\infty}$ relation.
\begin{equation}\label{Ainfinityrelbeta2}
\aligned
\sum_{k_1+k_2=k+1}&\sum_{\beta_1+\beta_2=\beta}\sum_{i=1}^{k-k_2+1} \\
&(-1)^*{\frak m}^P_{k_1,\beta_1}(h_1,\ldots,{\frak m}^P_{k_2,\beta_2}(h_i,\ldots,h_{i+k_2-1}),\ldots,h_{k}) = 0,
\endaligned\end{equation}
where $* = \deg' h_1 + \ldots + \deg'h_{i-1}$.
\end{enumerate}
\end{defn}
\begin{defn}
A {\it partial $P$-parametrized family of $G$-gapped filtered $A_{\infty}$ algebra structures on $\Omega(L)$ of energy cut level $E$ and of minimal energy $e_0$}
\index{$A_{\infty}$ algebra ! partial $P$-parametrized family of $G$-gapped filtered $A_{\infty}$ algebras}
is $\{\frak m^P_{k,\beta}\}$
satisfying the same properties as above except the following points:
\begin{enumerate}
\item[(a)]
$\frak m^P_{k,\beta}$ is defined only for $\beta, k$ with $E(\beta) + ke_0\le E$.
\item[(b)]
We require the $A_{\infty}$ relation (\ref{Ainfinityrelbeta2})
only for $\beta, k$ with $E(\beta) + ke_0 \le E$.
\item[(c)]
$\frak m^P_{k,\beta} = 0$ if $0 < E(\beta) < e_0$.
\end{enumerate}
\end{defn}
\begin{lem}
The notion of pseudo-isotopy of $G$-gapped filtered $A_{\infty}$ algebra
structures on $\Omega(L)$ is the same as the notion of the
$P = [0,1]$ parametrized family of $G$-gapped filtered $A_{\infty}$ algebra
structures on $\Omega(L)$.
The same holds for the partial $G$-gapped filtered $A_{\infty}$ algebra
structure on $\Omega(L)$.
\end{lem}

\begin{proof}
Let $(\{\frak m^t_{k,\beta}\},\{\frak c^t_{k,\beta}\})$ be the objects as in Definition \ref{pisotopydef}.
We define $\frak m^P_{k,\beta}$ as follows.
It suffices to consider the case $\beta \ne \beta_0$.
\par
Suppose $h_i$ does not contain $dt$. Then we put
$$
\frak m^P_{k,\beta}(h_1,\dots,h_k)
=
\frak m^t_{k,\beta}(h_1,\dots,h_k)
+ dt \wedge \frak c^t_{k,\beta}(h_1,\dots,h_k).
$$
We also put
$$
\frak m^P_{k,\beta}(h_1,\dots,dt \wedge h_i, \dots, h_k)
=
(-1)^* dt \wedge \frak m^t_{k,\beta}(h_1,\dots,h_k)
$$
where $* = \deg' h_1 + \dots + \deg' h_{i-1}$.
If at least two of $\hat h_1, \dots, \hat h_k$ contain $dt$, then
$$
\frak m^P_{k,\beta}(\hat h_1,\dots\hat h_k)
=
0.
$$
It is straightforward to check that
(\ref{isotopymaineq}) is equivalent to
(\ref{Ainfinityrelbeta2}).
\end{proof}
\begin{lemdef}
Let $Q$ be a connected component of the normalized corner
$\widehat S_k(P)$. Then
a $P$-parametrized family of $G$-gapped filtered $A_{\infty}$ algebra
structures on $\Omega(L)$ induces a
$Q$-parametrized family of $G$-gapped filtered $A_{\infty}$ algebra
structures on $\Omega(L)$.
\end{lemdef}
\begin{proof}
This is a consequence of the following fact.
If $
F : B_k(\Omega(P\times L)[1]) \to \Omega(P\times L)[1]
$
is pointwise in $P$ direction,
it induces
$: B_k(\Omega(Q\times L)[1]) \to \Omega(Q\times L)[1]
$
which is pointwise in $Q$ direction.
This fact is a consequence of the definition of pointwise-ness.
\end{proof}

\begin{rem}
If two $G$-gapped filtered $A_{\infty}$ algebra structures on $\Omega(L)$
are pseudo-isotopic, then
the two filtered $A_{\infty}$ algebras induced from those two structures
are homotopy equivalent in the sense of \cite[Definition 4.2.42]{fooobook}.
This fact is proved in \cite[Theorem 8.2]{fooo091}.
\end{rem}

\subsection{Statement of the results}
\label{subsec:19-4}

In this subsection and hereafter we say a filtered $A_{\infty}$ structure,
pseudo-isotopy or $P$-parametrised $A_{\infty}$ structure is
{\it gapped} when it is $G$-gapped for some discrete submonoid
$G \subset \R_{\ge 0} \times 2\Z$ in Definition \ref{discmonoid}.

\begin{shitu}\label{situ1933}
Let $\mu$, $E$ be as in Situation \ref{situ191}.
Let $L$ be a compact oriented smooth manifold without boundary.
In addition to those, we consider one of the following situations.
\begin{enumerate}
\item
We are given
$\mathcal{AF} = \{({\mathcal M}_{k+1}(\beta),{\rm ev}) \mid \beta,k\}$,
which defines an $A_{\infty}$ correspondence over $L$.
Hereafter we will not include ${\rm ev}$ in the notation for simplicity.
\item
We are given $\mathcal{AF} = \{{\mathcal M}_{k+1}(\beta) \mid \beta,k\}$,
which defines a partial $A_{\infty}$ correspondence over $L$
of energy cut level $E_0$ and minimal energy $e_0$.
\item
We are given
$$
\mathcal{AF}^j = \{{\mathcal M}^j_{k+1}(\beta) \mid \beta,k\} ~(j=0,1),
\quad
\mathcal{AF}^{[0,1]} = \{{\mathcal M}_{k+1}(\beta;[0,1]) \mid \beta,k\},
$$
which are two $A_{\infty}$ correspondences over $L$
and a pseudo-isotopy among them, respectively.
\item
We are given
$$
\mathcal{AF}^j = \{{\mathcal M}^j_{k+1}(\beta) \mid \beta,k\} ~(j=0,1),
\quad
\mathcal{AF}^{[0,1]} = \{{\mathcal M}_{k+1}(\beta;[0,1]) \mid \beta,k\},
$$
which are two partial $A_{\infty}$ correspondences over $L$
and a pseudo-isotopy among them, respectively.
Their energy cut level are $E_0$ and minimal energy are $e_0$.
\item
We are given
$$
\mathcal{AF}^i = \{{\mathcal M}^i_{k+1}(\beta)\mid \beta,k\}, \quad
\mathcal{AF}^{[i,i+1]} = \{{\mathcal M}_{k+1}(\beta;[i,i+1])\mid \beta,k\}
$$
for $i=1,2,\dots$ such that
$$
\mathcal{IAF} =
\left(\{\mathcal{AF}^i\mid i=1,2,\dots\},\{\mathcal{AF}^{[i,i+1]}\mid i=1,2,\dots\}\right)
$$
consists of an inductive system of $A_{\infty}$ correspondences over $L$.
\item
For $j = 0,1$, we are given
$$
\mathcal{AF}^{ji} = \{{\mathcal M}^{ji}_{k+1}(\beta)\mid \beta,k\}, \quad
\mathcal{AF}^{j,[i,i+1]} = \{{\mathcal M}_{k+1}^j(\beta;[i,i+1])\mid \beta,k\}
$$
for $i=1,2,\dots$
such that
$$
\mathcal{IAF}^j = \left(\{\mathcal{AF}^{ji}\mid i=1,2,\dots\},\{\mathcal{AF}^{j,[i,i+1]}\mid i=1,2,\dots\}\right)
$$
consist of two inductive systems of $A_{\infty}$ correspondences over $L$.
Moreover we are given
$\{{\mathcal M}^i_{k+1}(\beta;[0,1])\mid \beta,k\}$ and
$\{{\mathcal M}_{k+1}(\beta;[0,1]\times [i,i+1])\mid \beta,k\}$
which consist of a pseudo-isotopy of the inductive systems $\mathcal{IAF}^0$,
$\mathcal{IAF}^1$ of $A_{\infty}$ correspondences.
\item
Let $\mathcal{AF}^j = \{{\mathcal M}^j_{k+1}(\beta) \mid \beta,k\}$ define $A_{\infty}$ correspondences
over $L$ for $j=0,1$.
Let $\mathcal{AF}^{[0,1],\ell} = \{{\mathcal M}_{k+1}(\beta;[0,1]) \mid \beta,k\}$, $\ell =a,b$ be
two pseudo-isotopies from $\mathcal{AF}^1$ to $\mathcal{AF}^2$.
Let
$$
\mathcal{AF}^{[0,1]\times [1,2]}
$$
be a $[0,1]\times [1,2]$-parametrized family of $A_{\infty}$ correspondences over $L$ with the
following properties.
\begin{enumerate}
\item
On $[0,1] \times \{1\}$ it is isomorphic to $\mathcal{AF}^{[0,1],a}$.
\item
On $[0,1] \times \{2\}$ it is isomorphic to $\mathcal{AF}^{[0,1],b}$.
\item Let $j=1$ or $j=2$ then on
$\{j\} \times [1,2]$ it is isomorphic to the direct product
$\mathcal{AF}^j \times [1,2]$.
\item
At the corner $\{0,1\} \times \{1,2\}$ the isomorphisms (a)(b)(c)(d) and the various isomorphisms included in the
definitions of $A_{\infty}$ correspondence and its pseudo-isotopies are compatible in the same sense
as Condition \ref{linAinfmaincondspara} (X)(XI).
\end{enumerate}
\item
For $j = 0,1$, we are given
$$
\mathcal{AF}^{ji} = \{{\mathcal M}^{ji}_{k+1}(\beta)\mid \beta,k\}, \quad
\mathcal{AF}^{j,[i,i+1]} = \{{\mathcal M}_{k+1}^j(\beta;[i,i+1])\mid \beta,k\}
$$
for $i=1,2,\dots$
such that
$$
\mathcal{IAF}^j = \left(\{\mathcal{AF}^{ji}\mid i=1,2,\dots\},\{\mathcal{AF}^{j,[i,i+1]}\mid i=1,2,\dots\}\right)
$$
consist of inductive systems of $A_{\infty}$ correspondences over $L$.
Moreover for $\ell = a,b$, we are given
$\{{\mathcal M}^{i,\ell}_{k+1}(\beta;[0,1])\mid \beta,k\}$ and
$\{{\mathcal M}^{\ell}_{k+1}(\beta;[0,1]\times [i,i+1])\mid \beta,k\}$
which consist of two pseudo-isotopies of the inductive systems $\mathcal{IAF}^0$,
$\mathcal{IAF}^1$ of $A_{\infty}$ correspondences.
\par
Furthermore we assume that we have a pseudo-isotopy of pseudo-isotopies among
these two pseudo-isotopies in the following sense:
For each $i$ we have a $[0,1] \times [i,i+1] \times [1,2]$-parametrized family of partial
$A_{\infty}$ correspondence
$$
\mathcal{AF}^{[0,1] \times [i,i+1] \times [1,2]}
$$
over
$L$ of energy cut level $E^i$ and minimal energy $e_0$ with the following
properties.
\begin{enumerate}
\item[(i)] On $[0,1] \times\{i\} \times [1,2]$, it satisfies the same condition
as (7) (a)-(d) up to energy level $E^i$.
\item[(ii)]  On $[0,1] \times [i,i+1] \times\{ c\}$ with $c = 1$ (resp. $c=2$)
it is isomorphic to
$\{{\mathcal M}^{a}_{k+1}(\beta;[0,1]\times [i,i+1])\mid \beta,k\}$
(resp. $\{{\mathcal M}^{b}_{k+1}(\beta;[0,1]\times [i,i+1])\mid \beta,k\}$.)
\item[(iii)]
On $\{j\} \times [i,i+1] \times [1,2]$ with  $j =0$ or $j=1$, it
is isomorphic to the direct product $\mathcal{AF}^{j,[i,i+1]} \times [1,2]$.
\item[(iv)] Various isomorphisms in (i)(ii)(iii) above and those appearing
in the definitions of $A_{\infty}$ correspondences or its pseudo-isotopies
are compatible at $S_{m}([0,1] \times [i,i+1] \times [1,2])$
in the same sense as Condition \ref{linAinfmaincondspara} (X)(XI).
\item[(v)]
A similar uniform Gromov compactness as in Definition \ref{definition1918} (3) is satisfied.
\end{enumerate}
\end{enumerate}
$\blacksquare$
\end{shitu}

\begin{thm}\label{theorem1934}
Suppose we are given $\mu$, $E$ as in Situation \ref{situ191}.
Let $L$ be a compact oriented smooth manifold without boundary.
\begin{enumerate}
\item
Suppose we are in Situation \ref{situ1933} (1).
We can associate a gapped filtered $A_{\infty}$ structure on $\Omega(L)$.
This filtered $A_{\infty}$ structure is independent of the
choices made for its construction up to pseudo-isotopy.
\item
Suppose we are in Situation \ref{situ1933} (2).
We can associate a gapped partial filtered $A_{\infty}$ structure on $\Omega(L)$
of energy cut level $E_0$ and  of minimal energy $e_0$.
This partial filtered $A_{\infty}$ structure is independent of the
choices made for its construction up to gapped pseudo-isotopy.
\item
Suppose we are in Situation \ref{situ1933} (3).
We can associate a gapped pseudo-isotopy of
filtered $A_{\infty}$ structures on $\Omega(L)$
among the two gapped filtered $A_{\infty}$ algebras
which associate by (1) to ${\mathcal M}^j_{k+1}(\beta)$
for $j=0,1$.
\par
In particular, it induces a gapped homotopy equivalence between
those two gapped filtered $A_{\infty}$ algebras.
\item
Suppose we are in Situation \ref{situ1933} (4).
We can associate a gapped partial pseudo-isotopy of
filtered $A_{\infty}$ structures on $\Omega(L)$
of energy cut level $E_0$ and of minimal energy $e_0$
among the two gapped partial filtered $A_{\infty}$ algebras
which we associate by (2) to ${\mathcal M}^j_{k+1}(\beta)$
for $j=0,1$.
\item
Suppose we are in Situation \ref{situ1933} (5).
We can associate a gapped
filtered $A_{\infty}$ structure on $\Omega(L)$.
This gapped filtered $A_{\infty}$ structure is independent of the
choices made for its construction up to gapped pseudo-isotopy.
\item
Suppose we are in Situation \ref{situ1933} (6).
We can associate a gapped pseudo-isotopy of
gapped filtered $A_{\infty}$ structures on $\Omega(L)$
among the two gapped filtered $A_{\infty}$ algebras
which we associate in (5) to ${\mathcal M}^{ji}_{k+1}(\beta)$
and ${\mathcal M}^j_{k+1}(\beta;[i,i+1])$
for $j=0,1$.
\par
In particular, it induces a gapped homotopy equivalence  between
those two gapped filtered $A_{\infty}$ algebras.
\item
Suppose we are in Situation \ref{situ1933} (7).
By (1) we obtain two gapped filtered $A_{\infty}$ algebras
$(\Omega(L),\frak m^j_{k})$ for $j=1,2$.
By (3) we obtain two gapped homotopy equivalences $\varphi_a$ and $\varphi_b$ from $(\Omega(L),\frak m^1_{k})$ to $(\Omega(L),\frak m^2_{k})$,
corresponding to $\ell =a$ and $\ell = b$, respectively.
\par
Then the claim of (7) is that $\varphi_a$ is gapped homotopic to $\varphi_b$.
\item
We obtain the same conclusion as (7) in the Situation \ref{situ1933} (8).
\end{enumerate}
\end{thm}

\section{Tree-like K-system: $A_{\infty}$ structure II: proof}
\label{sec:systemtree2}

In this section we prove Theorem \ref{theorem1934}.
The proof is mostly parallel to the proofs given
in Section \ref{sec:systemline3}, and
also similar to those given in \cite[Subection 7.2]{fooobook2}, \cite{fooo091}.

\subsection{Existence of CF-perturbations}
\label{subsec:20-1}

\begin{defn}\label{defassoSubmono}
Let $\mathcal M_{k+1}(\beta)$ be the moduli spaces of the $A_{\infty}$
operations of a (partial) $A_{\infty}$ correspondence $\mathcal{AC}$.
We use the notation of Condition \ref{linAinfmaincondspara}.
Consider the submonoid of $\R_{\ge 0} \times 2\Z$ generated by
the subset
$$
\{(E(\beta),\mu(\beta)) \mid \beta \in \frak G,
\mathcal M_{k+1}(\beta) \ne \emptyset\}
$$
of $\R_{\ge 0} \times 2\Z$.
This submonoid is discrete by Condition \ref{linAinfmaincondspara} (VIII).
We call it {\it the discrete submonoid associated to $\mathcal{AC}$},
and denote it by $G(\mathcal{AC})$.
\par
When we have a partial $P$-parametrized family $\mathcal{AC}_P$
of $A_{\infty}$ correspondences whose moduli spaces of
$P$-parametrized $A_{\infty}$ operations are
$\mathcal M_{k+1}(\beta;P)$,
we define the {\it discrete submonoid $G(\mathcal{AC}_P)$ associated to
$\mathcal{AC}_P$}
as the submonoid generated by the subset
$$
\{(E(\beta),\mu(\beta)) \mid \beta \in \frak G,
\mathcal M_{k+1}(\beta;P) \ne \emptyset\}.
$$
\end{defn}
\begin{defn}
Consider a discrete submonoid $G \subset \R_{\ge 0}\times 2\Z$ in the sense of
Definition \ref{discmonoid}.
\begin{enumerate}
\item
We put
$$
e_{\rm min}(G) = \inf \{ E(\beta) \mid \beta \in G, E(\beta)>0\},
$$
if the image of $E : G \to \R$ is not $0$. Otherwise we put
$e_{\rm min}(G) = 1$.
\item
For $E_0, e_0 > 0$ with $e_0 \le e_{\rm min}(G)$ we define
$$
\aligned
\mathcal{GK}(G;E_0,e_0)
= \{(\beta,k) \mid \beta \in G, k \in \Z_{\ge 0},
\,\,& E(\beta) = 0 \Rightarrow k > 0  \\
& E(\beta) + ke_0\le E_0\}.
\endaligned
$$
\end{enumerate}
\end{defn}
Note that if $e_0$ is a minimal energy of a $G$-gapped (partial) $A_{\infty}$
correspondence $\mathcal{AC}$ and $G \supseteq G(\mathcal{AC})$, then $e_0 \le e_{\rm min}(G)$.
In case $G \ne G(\mathcal{AC})$ we put
$\mathcal M_{k+1}(\beta) = \emptyset$
for $\beta \notin G(\mathcal{AC})$ as convention.
\begin{prop}\label{prop203333}
Let $\mathcal{AC}$ be a partial $A_{\infty}$ correspondence  of energy cut level
$E_0$ and minimal energy $e_0$, and $G$ a discrete submonoid containing $G(\mathcal{AC})$.
Suppose $e_0 \le e_{\rm min}(G)$ and $0 < \tau < \tau_0=1$.
We can find a system of
$\tau$-collared Kuranishi structures and
CF-perturbations,
$\{(\widehat{\mathcal U^+_{k+1}}(\beta),\widehat{\frak S}_{k+1}(\beta)) \mid
(\beta,k) \in \mathcal{GK}(G;E_0,e_0)\}$, with the following properties:
\begin{enumerate}
\item
$\widehat{\mathcal U^+_{k+1}}(\beta)$ is a $\tau$-collared Kuranishi structure of
$\mathcal M_{k+1}(\beta)^{\boxplus\tau_0}$ and is a thickening of
the Kuranishi structure
obtained from one in Condition \ref{linAinfmainconds} (III)
by trivialization of the corner. (Lemma-Definition \ref{lemdef1522}).
Evaluation maps are extended to $\tau$-collared strongly smooth maps on this
$\tau$-collared Kuranishi structure and the associated evaluation map ${\rm ev}_0$ is weakly submersive.
\item
$\widehat{\frak S}_{k+1}(\beta)$ is a
$\tau$-collared CF-perturbation of the $\tau$-collared
Kuranishi structure $\widehat{\mathcal U^+_{k+1}}(\beta)$.
It is transversal to $0$ and ${\rm ev}_0$ is
strongly submersive with respect to $\widehat{\frak S}_{k+1}(\beta)$.
\item
There exists an isomorphism of $\tau$-collared K-spaces
\footnote{See Remark \ref{rem:FiberProdOrd} for the sign and the order of the fiber products.}
\begin{equation}\label{formula166plus}
\aligned
\partial ({\mathcal M}_{k+1}(\beta)^{\boxplus\tau_0},\widehat{\mathcal U^+_{k+1}}(\beta))
&\cong
\coprod_{\beta_1,\beta_2,k_1,k_2,i}
(-1)^{\epsilon}({\mathcal M}_{k_1+1}(\beta_1)^{\boxplus\tau_0},
\widehat{\mathcal U^+_{k_1+1}}(\beta_1)) \\
&\qquad\qquad
 \,\, {}_{{\rm ev}_i}\times_{{\rm ev}_0}
({\mathcal M}_{k_2+1}(\beta_2)^{\boxplus\tau_0},
\widehat{\mathcal U^+_{k_2+1}}(\beta_2)),
\endaligned
\end{equation}
where
\begin{equation}\label{eq:boundarysign3}
\epsilon = (k_1 -1)(k_2 -1) + \dim L + k_1 +
(i-1)\Big(1+(\mu(\beta_2)+k_2) \dim L\Big).
\end{equation}
\item
The restriction of $\widehat{\frak S}_{k+1}(\beta)$
to the boundary is equivalent (see \cite[Definition 7.5]{part11} for the definition of equivalence of CF-perturbations) to the
fiber product of $\widehat{\frak S}_{k_1+1}(\beta_1)$
and $\widehat{\frak S}_{k_2+1}(\beta_2)$ under the
isomorphism (\ref{formula166plus}).
\item
On the normalized corner $\widehat{S}_m( {\mathcal M}_{k+1}(\beta)^{\boxplus\tau_0})$, we put the Kuranishi structure that is the restriction of
$\widehat{\mathcal U^+_{k+1}}(\beta)$.
Then  it is isomorphic to the
disjoint union of the fiber products
\begin{equation}\label{cornecomAinf122sec20}
\prod_{(\mathcal T,\beta(\cdot))}
({\mathcal M}_{k_{\rm v}+1}(\beta({\rm v}))^{\boxplus\tau_0},
\widehat{\mathcal U^+_{k_{\rm v}+1}}(\beta({\rm v}))).
\end{equation}
Here the union is taken over all
$(\mathcal T,\beta(\cdot))
\in \mathcal G(k+1,\beta)$
with $\#  C_{1,{\rm int}}(\mathcal T) = m$.
The fiber product (\ref{cornecomAinf122sec20}) is defined as in
Definition \ref{defn1933}.
This isomorphism is compatible with the evaluation maps.
It is also compatible with the embedding of the
Kuranishi structures.
 (Here we mean the embedding of the
Kuranishi structure obtained from one in  Condition \ref{linAinfmainconds} (III)
by trivialization of the corner to $\widehat{\mathcal U^+_{k+1}}(\beta)$.)
\item
The restriction of $\widehat{\frak S}_{k+1}(\beta)$
to $\widehat{S}_m( {\mathcal M}_{k+1}(\beta)^{\boxplus\tau_0})$
is equivalent to the fiber product of
$\widehat{\frak S}_{k_{\rm v}+1}(\beta({\rm v}))$
under the isomorphism (\ref{cornecomAinf122sec20}).
\item
(5) implies that the restriction of $\widehat{\mathcal U^+_{k+1}}(\beta)$ to
$\widehat{S}_{\ell}(\widehat{S}_m( {\mathcal M}_{k+1}(\beta)^{\boxplus\tau_0}))$
is a disjoint union of copies of fiber products
\begin{equation}\label{cornecomAinf12233sec203}
\prod_{(\mathcal T,\beta(\cdot))}
({\mathcal M}_{k_{\rm v}+1}(\beta({\rm v}))^{\boxplus\tau_0},
\widehat{\mathcal U^+_{k_{\rm v}+1}}(\beta({\rm v})))
\end{equation}
where the union is taken over all
$(\mathcal T,\beta(\cdot))
\in \mathcal G(k+1,\beta)$
with $\#  C_{1,{\rm int}}(\mathcal T)= m+\ell$.
The fiber product (\ref{cornecomAinf12233sec203}) is defined
as in Definition \ref{defn1933}.
\par
The covering maps
$\widehat{S}_{\ell}(\widehat{S}_m( {\mathcal M}_{k+1}(\beta)^{\boxplus\tau_0}))
\to \widehat{S}_{m+\ell}( {\mathcal M}_{k+1}(\beta)^{\boxplus\tau_0})$
are the underlying map of the covering maps of the
K-spaces (where the Kuranishi structure is
$\widehat{\mathcal U^+_{k+1}}(\beta)$ etc.)
that is the identity map
on each of the components (\ref{cornecomAinf12233sec203}).
\end{enumerate}
\end{prop}
To prove Proposition \ref{prop203333} we use the following:

\begin{lem}\label{lem2044}
If $(\mathcal T,\beta(\cdot))
\in \mathcal G(k+1,\beta)$, $(k,\beta) \in \mathcal{GK}(G;E_0,e_0)$
and ${\rm v} \in C_{0,{\rm int}}(\mathcal T)$,
then $(\beta({\rm v}),k_{\rm v}) \in \mathcal{GK}(G;E_1,e_0)$ with $E_1 < E_0$.
\end{lem}
\begin{proof}
We consider $\mathcal T \setminus \{{\rm v}\}$.
It has $k_{\rm v} + 1$ connected components $\mathcal T_i$,
$i=0,\dots,k$.
Let $\ell$ be the number of its connected components
that do not contain exterior vertices.
If $\mathcal T_i$ does not contain exterior vertices,
then by Definition \ref{defn192} (7) we have
$$
\sum_{{\rm v'} \in C_{0,{\rm int}}(\mathcal T) \cap \mathcal T_i}
E(\beta_{\rm v'}) > 0.
$$
Therefore
$\ell e_0 + E(\beta_{\rm v}) \le E_0$.
The lemma follows immediately from this fact.
\end{proof}
\begin{proof}[Proof of Proposition \ref{prop203333}.]
The proof is given by induction on $E_0$.
Lemma \ref{lem2044} implies that we already obtained
$\widehat{\mathcal U^+_{k_{\rm v}+1}}(\beta({\rm v}))$
and $\widehat{\frak S}_{k_{\rm v}+1}(\beta_{\rm v})$
while we construct
$\widehat{\mathcal U^+_{k+1}}(\beta)$
and $\widehat{\frak S}_{k+1}(\beta)$.
Here  $k_{\rm v}$ and $\beta({\rm v})$ are as in
(\ref{cornecomAinf12233sec203}).
Therefore using Propositions \ref{prop1562}, \ref{prop529rev},
we can prove Proposition \ref{prop203333}
in the same way as the proof of Proposition \ref{prop161}.
\end{proof}
Next we consider the case of $P$-parametrized family.

\begin{shitu}\label{shitu20555}
Suppose we have a partial $P$-parametrized family $\mathcal{AC}_P$
of $A_{\infty}$ correspondences whose moduli spaces of
$P$-parametrized $A_{\infty}$ operations are
$\mathcal M_{k+1}(\beta;P)$.
Let $E_0$ be its energy cut level and $e_0$ its minimal energy.
We assume that we are given the family of the pairs
$\widehat{\mathcal U^+_{k+1}}(\beta;\widehat S_{m}(P))$,
$\widehat{\frak S}_{k+1}(\beta;\widehat S_{m}(P))$
for $m\ge 1$ and $0 < \tau < \tau_0=1$ with the following properties:
\begin{enumerate}
\item
$\widehat{\mathcal U^+_{k+1}}(\beta;\widehat S_{m}(P))$
is a $\tau$-collared Kuranishi structure of $\mathcal M_{k+1}(\beta;\widehat S_{m}(P))^{\boxplus\tau_0}$
and is a thickening of the $\tau$-collared Kuranishi structure of
$\mathcal M_{k+1}(\beta;\widehat S_{m}(P))^{\boxplus\tau_0}$
that is a trivialization of the corner of the
Kuranishi structure
given by Condition \ref{linAinfmaincondspara} (III).
\item
$\widehat{\frak S}_{k+1}(\beta;\widehat S_{m}(P))$
is a CF-perturbation of
$\widehat{\mathcal U^+_{k+1}}(\beta;\widehat S_{m}(P))$.
\item
$\widehat{\frak S}_{k+1}(\beta;\widehat S_{m}(P))$
is transversal to $0$.
The evaluation map
$({\rm ev}_0,{\rm ev}_{\widehat S_{m}(P)})$ is strongly submersive
with respect to $\widehat{\frak S}_{k+1}(\beta;\widehat S_{m}(P))$.
Here ${\rm ev}_{\widehat S_{m}(P)}$ is the restriction of ${\rm ev}_{P}$ to
$\widehat S_{m}(P)$.
\item
By Condition \ref{linAinfmaincondspara} (XI),
$\widehat S_{\ell}(\mathcal M_{k+1}(\beta;\widehat S_{m}(P))^{\boxplus\tau_0})$ is
isomorphic to the
disjoint union of copies of
\begin{equation}\label{cornecomAinf122204}
\prod_{(\mathcal T,\beta(\cdot))} {\mathcal M}_{k_{\rm v}+1}(\beta({\rm v});
\widehat S_{\ell'+m}(P))^{\boxplus\tau_0}.
\end{equation}
This fiber product is defined in the same way as
(\ref{cornecomAinf122333}), (\ref{fiberproducttreerev}).
Here the union is taken over all
$(\mathcal T,\beta(\cdot))
\in \mathcal G(k+1,\beta)$ and $\ell' \in \Z_{\ge 0}$
with $\#  C_{1,{\rm int}}(\mathcal T)
+ \ell'= \ell$.
\begin{enumerate}
\item
The restriction of $\widehat{\mathcal U^+_{k+1}}(\beta;\widehat S_{m}(P))$
to $\widehat S_{\ell}(\mathcal M_{k+1}(\beta;\widehat S_{m}(P)))^{\boxplus\tau_0}$
is isomorphic to the fiber product of
$\widehat{\mathcal U^+_{\ell' + k_{\rm v}+1}}(\beta({\rm v});\widehat S_{\ell'+m}(P))$
on (\ref{cornecomAinf122204}).
\item
The isomorphism in (a) is compatible with the covering map  from
$\widehat S_{\ell}(\mathcal M_{k+1}(\beta;\widehat S_{m}(P)))^{\boxplus\tau_0}$
to $\mathcal M_{k+1}(\beta;\widehat S_{\ell+m}(P))^{\boxplus\tau_0}$.
Namely it is induced by the covering map of the Kuranishi structures.
\item
The restriction of the CF-perturbation
$\widehat{\frak S}_{k+1}(\beta;\widehat S_{m}(P))$
to
$$
\widehat S_{\ell}(\mathcal M_{k+1}(\beta;\widehat S_{m}(P)))^{\boxplus\tau_0}
$$
is equivalent to the fiber product of the CF-perturbations
$$
\widehat{\frak S}_{k_{\rm v}+1}(\beta({\rm v});\widehat S_{\ell'+m}(P))
$$
under the isomorphism in (a).
\end{enumerate}
\end{enumerate}
$\blacksquare$
\end{shitu}
\begin{prop}\label{prop203333rev}
In Situation \ref{shitu20555} let $G$ be a discrete monoid containing $G(\mathcal{AC}_P)$.
Then we can find a system of Kuranishi structures and
CF-perturbations
$$
\{(\widehat{\mathcal U^+_{k+1}}(\beta;P),\widehat{\frak S}_{k+1}(\beta;P)) \mid
(\beta,k) \in \mathcal{GK}(G;E_0,e_0)\}
$$
with the following properties:
\begin{enumerate}
\item
$\widehat{\mathcal U^+_{k+1}}(\beta;P)$ is a $\tau$-collared Kuranishi
structure of
$\mathcal M_{k+1}(\beta;P)^{\boxplus\tau_0}$ and is a thickening of
the Kuranishi structure
obtained from one in Condition \ref{linAinfmaincondspara} (III)
by trivialization of the corner. (Lemma-Definition \ref{lemdef1522}).
The evaluation maps are extended to strongly smooth maps on this
Kuranishi structure and $({\rm ev}_0,{\rm ev}_P)$ is stratumwise weakly submersive.
\item
$\widehat{\frak S}_{k+1}(\beta;P)$ is a CF-perturbation of the Kuranishi structure $\widehat{\mathcal U^+_{k+1}}(\beta;P)$.
It is transversal to $0$ and $({\rm ev}_0,{\rm ev}_P)$ is stratumwise
strongly submersive with respect to $\widehat{\frak S}_{k+1}(\beta;P)$.
\item
There exists an isomorphism of K-spaces\footnote{
As we note in the footnote at Proposition \ref{prop161} (2), we
simply write ${\rm ev}_i$ in place of ${\rm ev}^{\boxplus\tau_0}_i$. Similarly
we write ${\rm ev}_P$ in place of ${\rm ev}^{\boxplus\tau_0}_P$.}
\begin{equation}\label{formula166plusrev}
\aligned
&\partial ({\mathcal M}_{k+1}(\beta;P)^{\boxplus\tau_0},\widehat{\mathcal U^+_{k+1}}(\beta;P))
\\
&\cong
\coprod_{\beta_1,\beta_2,k_1,k_2,i}
(-1)^{\epsilon}({\mathcal M}_{k_1+1}(\beta_1;P)^{\boxplus\tau_0},
\widehat{\mathcal U^+_{k_1+1}}(\beta_1;P)) \\
&\qquad\qquad
 \,\, {}_{({\rm ev}_P,{\rm ev}_i)}\times_{({\rm ev}_P,{\rm ev}_0)}
({\mathcal M}_{k_2+1}(\beta_2;P)^{\boxplus\tau_0},
\widehat{\mathcal U^+_{k_2+1}}(\beta_2;P)) \\
&\quad\sqcup
({\mathcal M}_{k+1}(\beta;\partial P)^{\boxplus\tau_0},\widehat{\mathcal U^+_{k+1}}(\beta;\partial P)),
\endaligned
\end{equation}
where
$$
\epsilon = (k_1 -1)(k_2 -1) + \dim L + k_1 +
(i-1)\Big(1+(\mu(\beta_2)+k_2 +\dim P) \dim L\Big).
$$
\item
The restriction of $\widehat{\frak S}_{k+1}(\beta;P)$
to the boundary is equivalent to the
fiber product of $\widehat{\frak S}_{k_1+1}(\beta_1;P)$
and $\widehat{\frak S}_{k_2+1}(\beta_2;P)$ on the first
summand of the right hand side of the
isomorphism (\ref{formula166plusrev}).
It is equivalent to $\widehat{\frak S}_{k+1}(\beta;\widehat S_1(P))
= \widehat{\frak S}_{k+1}(\beta;\partial P)$
on the second
summand of the right hand side of (\ref{formula166plusrev}).
\item
On the normalized corner $\widehat{S}_m( {\mathcal M}_{k+1}(\beta;P)^{\boxplus\tau_0})$, we put the Kuranishi structure that is the restriction of
$\widehat{\mathcal U^+_{k+1}}(\beta;P)$.
Then it is isomorphic to the
disjoint union of the fiber products
\begin{equation}\label{cornecomAinf122sec20rev}
\prod_{(\mathcal T,\beta(\cdot))}
({\mathcal M}_{k_{\rm v}+1}(\beta({\rm v});\widehat{S}_{m'}(P))^{\boxplus\tau_0},
\widehat{\mathcal U^+_{k_{\rm v}+1}}(\beta({\rm v});\widehat{S}_{m'}(P))).
\end{equation}
This fiber product is defined in the same way as
(\ref{cornecomAinf122333}), (\ref{fiberproducttreerev}).
Here the union is taken over all
$(\mathcal T,\beta(\cdot))
\in \mathcal G(k+1,\beta)$
with $\#  C_{1,{\rm int}}(\mathcal T) + m' = m$.
This isomorphism is compatible with the evaluation maps.
It is also compatible with the embedding of the
Kuranishi structures. (Here we mean the embedding of the
Kuranishi structures obtained from one in  Condition \ref{linAinfmaincondspara} (III)
by trivialization of the corner to $\widehat{\mathcal U^+_{k+1}}(\beta;P)$.)
\item
The restriction of $\widehat{\frak S}_{k+1}(\beta;P)$
to $\widehat{S}_m( {\mathcal M}_{k+1}(\beta;P)^{\boxplus\tau_0})$
is equivalent to the fiber product of
$\widehat{\frak S}_{k_{\rm v}+1}(\beta({\rm v}),\widehat{S}_{m'}(P))$
under the isomorphism (\ref{cornecomAinf122sec20rev}).
\item
(5) and Situation \ref{shitu20555}
imply that the restriction of $\widehat{\mathcal U^+_{k+1}}(\beta;P)$ to
the iterated normalized corner $\widehat{S}_{\ell}(\widehat{S}_m( {\mathcal M}_{k+1}(\beta;P)^{\boxplus\tau_0}))$
is a disjoint union of copies of the fiber products
\begin{equation}\label{cornecomAinf12233sec203rev}
\prod_{(\mathcal T,\beta(\cdot))}
({\mathcal M}_{k_{\rm v}+1}(\beta({\rm v});\widehat S_{m'}(P))^{\boxplus\tau_0},
\widehat{\mathcal U^+_{k_{\rm v}+1}}(\beta({\rm v});\widehat S_{m'}(P)))
\end{equation}
where the union is taken over all
$(\mathcal T,\beta(\cdot))
\in \mathcal G(k+1,\beta)$
with $\#  C_{1,{\rm int}}(\mathcal T) + m'= m+\ell$.
The fiber product (\ref{cornecomAinf12233sec203rev}) is defined in the same way as
(\ref{cornecomAinf122333}), (\ref{fiberproducttreerev}).
\par
The covering map
$\widehat{S}_{\ell}(\widehat{S}_m( {\mathcal M}_{k+1}(\beta;P)^{\boxplus\tau_0}))
\to \widehat{S}_{m+\ell}( {\mathcal M}_{k+1}(\beta;P)^{\boxplus\tau_0})$
is the underlying map of the covering map of the
K-spaces (where the Kuranishi structure is
$\widehat{\mathcal U^+_{k+1}}(\beta;P)$ etc.)
that is the identity map
on each component of (\ref{cornecomAinf12233sec203rev}).
\item If we have a uniform family of $\widehat{\frak S}_{k+1}(\beta;\widehat S_{m}(P))$ as  in
Situation \ref{shitu20555}, then we obtain a uniform family of
$\widehat{\frak S}_{k+1}(\beta;P)$.
\end{enumerate}
\end{prop}
\begin{proof}
Using Lemma \ref{lem2044}, Propositions \ref{prop1562}, \ref{prop529rev},
we can prove Proposition \ref{prop203333rev}
in the same way as the proof of Proposition \ref{prop161}.
\end{proof}
\begin{rem}
In Section \ref{sec:systemline3} we use the ${\frak C}^{h}$-partial trivialization of the corner,
so the $P$-parametrized family after trivialization of the corner remains to be
$P$-parametrized.
Here we use the trivialization of the corner.
Therefore after the trivialization of the corner
we will get a $P^{\boxplus\tau_0}$-parametrized family.
\par
Later in the proof of Proposition \ref{prop2013}, we use the collared-ness of
the parametrized family in an algebraic model.
To obtain the collared family, we also need to trivialize the corner in $P$ direction.
\end{rem}
\subsection{Algebraic lemmas: promotion lemmas via pseudo-isotopy}
\label{subsec:20-2}
\begin{defn}
Let $E_0' < E_0$.
\begin{enumerate}
\item
Suppose $\{\frak m_{k,\beta}\}$ is a partial $G$-gapped filtered $A_{\infty}$ structure on $\Omega(L)$
of energy cut level $E_0$ and of minimal energy $e_0$.
We forget all the $\frak m_{k,\beta}$'s with
$E(\beta) > E'_0$ and obtain a
partial $G$-gapped filtered $A_{\infty}$ structure on $\Omega(L)$ of energy cut level $E'_0$.
We call it the partial filtered $A_{\infty}$ structure on $\Omega(L)$
obtained by the {\it energy cut at $E'_0$}.
\item
Suppose $\{\frak m_{k,\beta}\}$
is obtained from $\{\frak m'_{k,\beta}\}$  by the energy cut at $E'_0$
and $\{\frak m'_{k,\beta}\}$
is a partial $G$-gapped filtered $A_{\infty}$ structure on $\Omega(L)$ of energy
cut level $E_0$,
we call $\{\frak m'_{k,\beta}\}$ a {\it promotion of $\{\frak m_{k,\beta}\}$
to the energy
cut level $E_0$}.
\item
We define an energy cut or a promotion of
pseudo-isotopy or a $P$-parametrized family
of partial $A_{\infty}$ structures in the same way.
\end{enumerate}
\end{defn}
We will use the next proposition which says that we can extend the promotion
of partial $A_{\infty}$ structures using pseudo-isotopy.

\begin{prop}\label{prop208}
Fix a discrete submonoid $G$ and $e_0 \le e_{\rm min}(G)$.
Let $E_0 < E_1$. For each $j=0,1$
let $\{\frak m^j_{k,\beta}\}$ be a $G$-gapped partial filtered $A_{\infty}$
structure of energy cut level $E_j$ and minimal energy $e_0$
on $\Omega(L)$. Suppose that we are given
a $G$-gapped partial filtered $A_{\infty}$ pseudo-isotopy
$(\{\frak m^t_{k,\beta}\},\{\frak c^t_{k,\beta}\})$
of energy cut level $E_0$ and  minimal energy $e_0$,
from $\{\frak m^0_{k,\beta}\}$ to the
energy cut of $\{\frak m^1_{k,\beta}\}$ at $E_0$.
Then we can promote $\{\frak m^0_{k,\beta}\}$
to energy cut level $E_1$ and
$(\{\frak m^t_{k,\beta}\},\{\frak c^t_{k,\beta}\})$
to energy cut level $E_1$.
\end{prop}
\begin{proof}
The proof is the same as the proof of \cite[Theorem 8.1]{fooo091}.
(The only difference is
the following point: In \cite[Theorem 8.1]{fooo091}
partial structures are ones where we take only finitely many
$\beta$'s: Here we take finitely many $(\beta,k)$'s.)
We repeat the proof for completeness.
\par
We consider the set
$$
\frak E = \{ E(\beta) + ke_0 \mid (\beta,k) \in G \times \Z_{\ge 0} \}.
$$
This is a discrete set. So
by applying an induction we may and will assume that
\begin{equation}
\# (\frak E  \cap (E_0,E_1]) = 1.
\end{equation}
Let $(\beta,k) \in G \times \Z_{\ge 0}$ such that
$E(\beta) + ke_0 \in (E_0,E_1]$.
We will define $\frak m^t_{k,\beta}$ and $\frak c^t_{k,\beta}$
for each such $(\beta,k)$.
\par
We put $\frak c^t_{k,\beta} = 0$.
Then there exists a unique $\frak m^t_{k,\beta}$ such that it satisfies
(\ref{isotopymaineq}) and $\frak m^t_{k,\beta}=\frak m^1_{k,\beta}$ for $t=1$.
Note that (\ref{isotopymaineq}) can be regarded as an ordinary differential equation
for each fixed $(h_1,\dots,h_k)$. Therefore
$\frak m^t_{k,\beta}$ depends smoothly on $t$ and is local in the $[0,1]$-direction.
\par
Next we check the $A_{\infty}$ relation for each fixed $t$ in the case of
$(\beta,k)$. We calculate
\begin{equation}\label{ainfityinduction}
\aligned
&\frac{d}{dt}
\sum_{k_1+k_2=k+1}\sum_{\beta_1+\beta_2=\beta}\sum_{i=1}^{k-k_2+1}
\frak m^t_{k_1,\beta_1}(h_1,\ldots, \frak m_{k_2,\beta_2}^t(h_i,\ldots),\ldots,h_k) \\
=& \sum_{k_1+k_2=k+1}\sum_{\beta_1+\beta_2=\beta}\sum_{i=1}^{k-k_2+1}
\frac{d \frak m^t_{k_1,\beta_1}}{dt}(h_1,\ldots, \frak m_{k_2,\beta_2}^t(h_i,\ldots),\ldots,h_k) \\
&+ \sum_{k_1+k_2=k+1}\sum_{\beta_1+\beta_2=\beta}\sum_{i=1}^{k-k_2+1}
\frak m^t_{k_1,\beta_1}(h_1,\ldots, \frac{d \frak m^t_{k_2,\beta_2}}{dt}(h_i,\ldots),\ldots,h_k).
\endaligned
\end{equation}
Using (\ref{isotopymaineq}) and the $A_{\infty}$ relation,
(that is, the induction hypothesis), it is easy to see that
(\ref{ainfityinduction}) is zero.
We define $\frak m^0_{k,\beta}$ as the case $t=0$ of $\frak m^t_{k,\beta}$.
The proof of Proposition \ref{prop208} is complete.
\end{proof}
We next study the promotion of pseudo-isotopy
using pseudo-isotopy of pseudo-isotopies.
The proof is similar to that of \cite[Theorem 14.1]{fooo091}.
We repeat the detail of the proof for completeness.
\par
We first define the notion of a $P^{\boxplus\tau}$-parametrized (partial) $A_{\infty}$
structure to be collared.
Here $P^{\boxplus\tau}$ is the trivialization of the corner of our
manifold with corner $P$.
(A similar assumption appeared in \cite[Assumption 14.1]{fooo091}.)
\begin{defn}\label{Pparamet}
Let $\{ \frak m^{P^{\boxplus\tau}}_{k,\beta}\}$
be a $P^{\boxplus\tau}$-parametrized partial
$A_{\infty}$-structure of energy cut level $E_0$ and minimal energy $e_0$.
We say it is {\it $\tau$-collared}
if the following conditions (1)(2) are satisfied.
\par
Let ${\bf t} \in \overset{\circ\circ}S_k(P^{\boxplus\tau})$.
Its $\tau$-collared neighborhood is identified with
$V \times [-\tau,0)^k$.
Let $(t'_1,\dots,t'_m)$ be a coordinate of $V$ and let $(t''_1,\dots,t''_k)$
be the standard coordinate of $[-\tau,0)^k$.
A differential form on $P$ in a neighborhood is written as
$\sum f_{I'I''} dt^{\prime}_{I'} \wedge dt^{\prime\prime}_{I''}$
where $dt^{\prime}_{I'}$ are wedge products of $dt'_i$'s, and
$dt^{\prime\prime}_{I''}$ are wedges products of $dt''_i$'s.
\par
By definition, $\frak m^{P^{\boxplus\tau}}_{k,\beta}$ is written
on this neighborhood as the form
$$
\frak m^{P^{\boxplus\tau}}_{k,\beta}(h_1,\dots,h_k)
=
\sum_{I,I'}
dt^{\prime}_{I'} \wedge dt^{\prime\prime}_{I''}
\wedge \frak m^{t',t''}_{k,\beta;I',I''}(h_1,\dots,h_k).
$$
Now we require:
\begin{enumerate}
\item
$\frak m^{t',t''}_{k,\beta;I',I''}(h_1,\dots,h_k) = 0$
unless $I'' = \emptyset$.
\item
If $I''=\emptyset$,
$\frak m^{t',t''}_{k,\beta;I',\emptyset}(h_1,\dots,h_k)$
is independent of $t'' \in [-\tau,0)^k$.
\end{enumerate}
We say a {\it $P$-parametrized partial
$A_{\infty}$ structure is collared}
\index{$\tau$-collared ! $P$-parametrized partial
$A_{\infty}$ structure}
\index{$A_{\infty}$ algebra ! collared $P$-parametrized partial
$A_{\infty}$ structure}
if
there exist $\tau >0$ and $P'$ such that $P = P^{\prime \boxplus \tau}$
\end{defn}
\begin{exm}
The case when $P=[0,1]$
in Definition \ref{Pparamet} is nothing but the case of
pseudo-isotopy in Definition \ref{pisotopydef}.
In this case, $P^{\boxplus\tau} = [-\tau,1+\tau]$ and
the $\tau$-collared-ness property (1), (2)
in Definition \ref{Pparamet} implies the following property
of pseudo-isotopy
$(\{\frak m^t_{k,\beta}\},\{\frak c^t_{k,\beta}\})$,
respectively:
\begin{enumerate}
\item
$\frak c^t_{k,\beta} = 0$ for $t \in [-\tau,0] \cup [1,1+\tau]$.
\item
$\frac{d}{dt}\frak m^t_{k,\beta} = 0$ for $t \in [-\tau,0] \cup [1,1+\tau]$.
\end{enumerate}
\end{exm}

\begin{shitu}\label{situ2011}
Let $P$ be a manifold with corner and $E_1 > E_0 \ge 0$, $e_0 > 0$.
We assume that we are given the following objects.
\begin{enumerate}
\item
A $P\times [0,1]$-parametrized collared
partial
$A_{\infty}$ structure $\{ \frak m^{P\times [0,1]}_{k,\beta}\}$
of energy cut level $E_0$ and of minimal energy $e_0$ on $\Omega(L)$.
\item
A collared promotion of the restriction of $\{ \frak m^{P\times [0,1]}_{k,\beta}\}$
to $P \times \{1\}$ to energy cut level $E_1$.
\item
Let $\partial P = \coprod \partial_i P$ be the decomposition of the normal
boundary of $P$ into the connected components.
Then we also assume that a collared promotion of the restriction of
$\{ \frak m^{P\times [0,1]}_{k,\beta}\}$ to $\partial_i P \times [0,1]$
 to energy cut level $E_1$ is
given for each $i$.
\item
We assume that the restriction of the promotion in (2)
coincides with the promotion in (3) on
 $\partial_i P \times \{1\}$.
\item
Suppose that the images of $\partial_i P$ and $\partial_j P$ intersect
each other in
$P$ at the component $\partial_{ij} P$ of the codimension 2 corner of $P$.
(Note that the case $i=j$ is included.
In this case, $\partial_{ii} P$ is the `self intersection' of
$\partial_i P$.)
Then we assume that the promotions of the restrictions on
$\partial_i P \times [0,1]$ and on $\partial_j P \times [0,1]$
in (3) coincide with each other on $\partial_{ij} P \times [0,1]$.
\end{enumerate}
$\blacksquare$
\end{shitu}
\begin{rem}
In Situation \ref{situ2011}
we assumed the compatibility of the promotion only at the
codimension 2 corners. In this situation it automatically
implies that they coincide at higher codimensional corners.
This is because our assumptions are their exact coincidence and not
coincidence up to certain equivalence relation.
\end{rem}
\begin{prop}\label{prop2013}
In Situation \ref{situ2011} there exists a promotion of
$\{ \frak m^{P\times [0,1]}_{k,\beta}\}$ to energy cut level $E_1$
such that the promotion coincides with those given in
Situation \ref{situ2011} (2) (resp. (3))
on  $P \times \{1\}$ (resp. $\partial_i P \times [0,1]$).
\end{prop}
\begin{proof}
We first prove Proposition \ref{prop2013} for the case $P = [-\tau,1+\tau]$.
We regard $P \times [0,1]$ as $P \times [-\tau,1+\tau]
= ([0,1]^2)^{\boxplus\tau}$ and assume that our structures
are $\tau$-collared.
\par
We change the corner structure of $([0,1]^2)^{\boxplus\tau}$ so that
we smooth the corners at two points $(-\tau, 1+\tau)$, $(1+\tau,1+\tau)$ and
make two points $(-\tau, -\tau/4)$, $(1+\tau, -\tau/4)$
into new corners instead.
We leave two other corners
$(-\tau,-\tau), (1+\tau, -\tau)$ as corners.
We then get a new cornered 2 manifold $Q$ diffeomorphic to
$[-\tau,1+\tau]^2$.
We denote this diffeomorphism by
$F : [-\tau,1+\tau]^2 \to Q$.
The diffeomorphism
$F$ is different from the set theoretical identity map
${\rm id} : [-\tau,1+\tau]^2 \to Q$.
In fact, we can take $F$ satisfying the following properties. See Figure \ref{Figure21-1}.
\begin{enumerate}
\item $F(-\tau,-\tau/4) = (-\tau,1+\tau)$ and
$F(1+\tau,-\tau/4) = (1+\tau , 1+\tau)$.
\item $F$ is identity on the edge $[-\tau, 1+\tau] \times \{ -\tau \}$.
\end{enumerate}
\begin{figure}[h]
\centering
\includegraphics[scale=0.5]{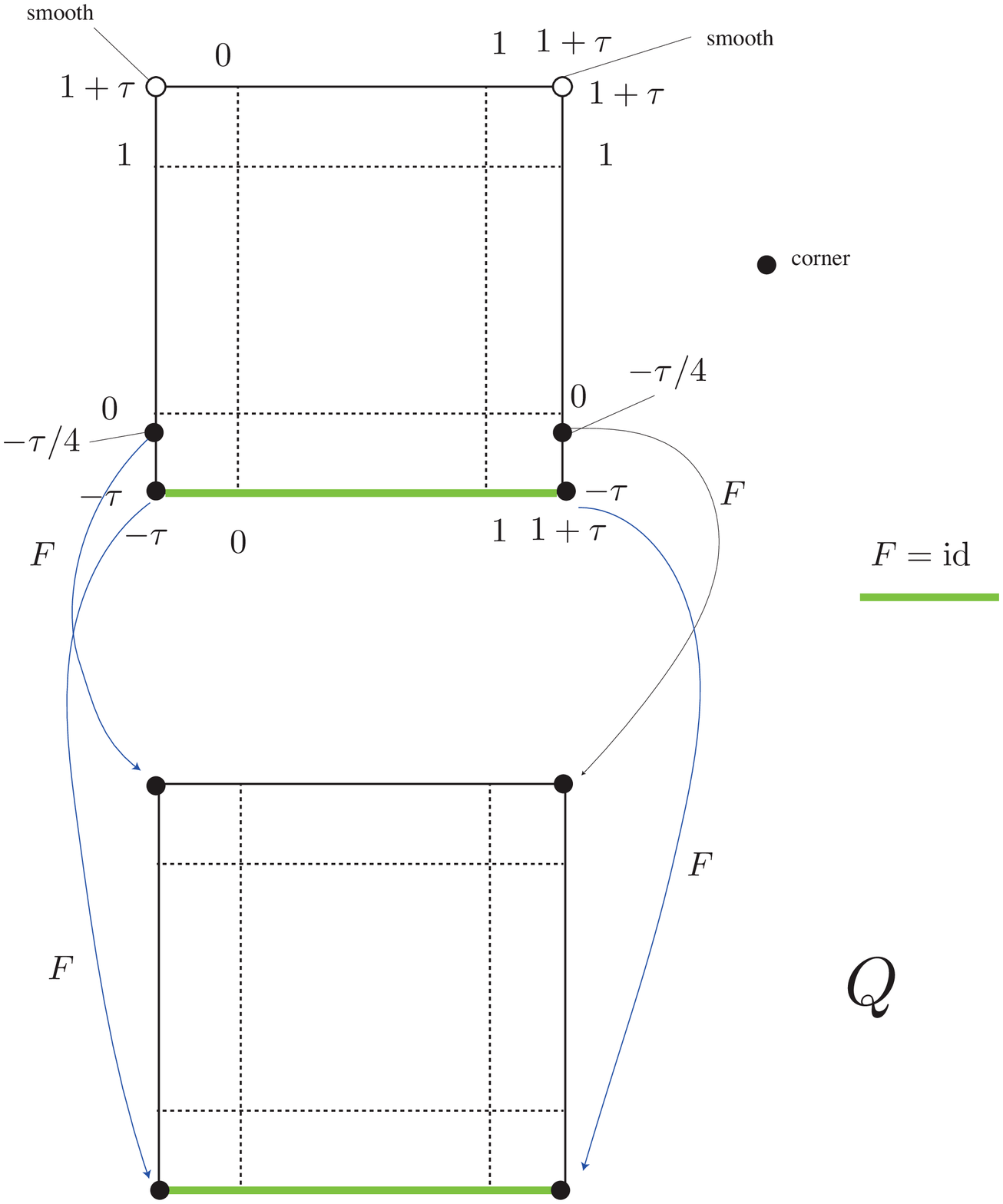}
\caption{$Q$ and $F$}
\label{Figure21-1}
\end{figure}
Using the $\tau$-collared-ness, our structures give a
$Q$-parametrized family of partial $A_{\infty}$ structures.
We regard it as a $[-\tau,1+\tau]^2$-parametrized family of partial
$A_{\infty}$ structures under the diffeomorphism $F$.
By assumption, its energy cut level is $E_0$ and
the energy cut level of its restriction to $[-\tau,1+\tau] \times \{1+\tau \}$ is $E_1$.
Therefore we can apply Proposition \ref{prop208}
to promote this $[-\tau,1+\tau]^2$-parametrized family to
energy cut level  $E_1$.
Using collared-ness again, we find that on $\partial [-\tau,1+\tau] \times [-\tau,1+\tau]
\subset  [-\tau,1+\tau]^2  \overset{F}\cong Q$ this promotion coincides with the
structure of energy cut level $E_1$ given at the beginning.
\par
Now we identify $Q \overset{\rm id}\cong [-\tau,1+\tau]^2$. At the place where
we smooth corners or make new corners, we can use the $\tau$-collared-ness to show that
the promotion coincides with the one originally given
at the beginning. Thus we obtain the required promotion.
Note the structure obtained is $\tau'$-collared for some $0< \tau' <\tau$ by construction.
\par
Thus we have proved Proposition \ref{prop2013} for the case $P = [0,1]$.
(We use only this case in this book.)
\par
The general case can be proved in a similar way.
Namely we smooth some of corners of $(P \times [0,1])^{\boxplus \tau}$ and make certain points
into new corners to obtain a cornered manifold $Q$
so that the following holds.
\begin{enumerate}
\item There exists a diffeomorphism
$F : (P \times [0,1])^{\boxplus \tau} \to Q$ which is
identity on $P^{\boxplus \tau} \times \{-\tau\}$.
\item
$\{ \frak m^{P\times [0,1]}_{k,\beta}\}$ induces a $Q$ parameter family
of partial $A_{\infty}$ structures.
\end{enumerate}
By the diffeomorphism in (2), the set
$P^{\boxplus \tau} \times\{1+\tau\} \subset (P\times [0,1])^{\boxplus \tau} = Q$,
where the last
equality is the set-theoretical one,
is mapped from a subset of
$(\partial P^{\boxplus \tau} \times [-\tau,1+\tau]) \cup
(P^{\boxplus \tau} \times \{1+\tau\})$.
Therefore the partial structure in (2) is one of energy cut level $E_1$ on
$P^{\boxplus \tau} \times\{1+\tau\}$.
On $(P \times [0,1])^{\boxplus \tau}$ it is of energy cut level $E_0$.
We apply Proposition \ref{prop208}  to promote it to the energy cut level $E_1$.
Using the diffeomorphism $F$, we regard it as a
$(P \times [0,1])^{\boxplus \tau}$-parametrized structure.
Using collared-ness, we find that it induces a
$(P \times[0,1])$-parametrized structure
via the identity map. (Note that identity map
$(P\times [0,1])^{\boxplus \tau} \cong Q$ is not a diffeomorphism.
However the structures are constant at the place where differentiability breaks down.)
\par
Thus we have obtained the required promotion.
\end{proof}

\subsection{Pointwise-ness of parametrized family of
smooth correspondences}
\label{subsec:20-3}

In this subsection we prove
Proposition \ref{prop2012} which reads that the operation
defined as a smooth correspondence associated to a $P$-parametrized family of $A_{\infty}$ correspondences is pointwise
in $P$ direction in the sense of Definition \ref{defn1926}.
To state the result in the way we can utilize in similar but different situations,
we slightly generalize Definition \ref{defn1926}.
We use the notation $t_J$ etc. of Definition \ref{defn1926} in the next definition.
\begin{defn}\label{defn1926rev}
Let $M_s$, $M_t$ be smooth manifolds (without boundary) and $P$
a smooth manifold with corner.
A linear map
$
F : \Omega(P\times M_s) \to \Omega(P\times M_t)
$
is said to be {\it pointwise in $P$ direction}
\index{pointwise in $P$ direction} if
the following holds:
\par
For each $I \subset \{1,\dots,d\}$
and ${\bf t} \in P$ there exists a linear and continuous map
$F^{\bf t}_{I;J} : \Omega(M_s) \to \Omega(M_t)$
such that
\begin{equation}\label{form208}
F(dt_{J} \wedge h)\vert_{\{{\bf t}\} \times M_t}
=
\sum_{I} dt_I \wedge dt_{J}
\wedge F^{\bf t}_{I;J}(h\vert_{{\{{\bf t}\}} \times M_t}).
\end{equation}
Moreover $F^{\bf t}_{I;J}$ depends smoothly on ${\bf t}$
(with respect to the operator topology) and is independent of $J$ up to sign.
\end{defn}
In case $M_s = L^k$ and $M_t = L$, Definition \ref{defn1926rev}
is nothing but Definition \ref{defn1926}.

\begin{shitu}\label{situ2010}
Let $(X,\widehat{\mathcal U})$ be a K-space, and let $M_s$, $M_t$
be smooth manifolds without boundary and $P$ a
smooth manifold with corners.
Let $f_s : (X,\widehat{\mathcal U}) \to M_s$, $f_t : (X,\widehat{\mathcal U}) \to M_t$
and $f_P : (X,\widehat{\mathcal U}) \to P$ be strongly smooth maps.
We assume that $(f_t,f_P) : (X,\widehat{\mathcal U}) \to P\times M$ is stratumwise
weakly submersive.
\par
Let $\widehat{\frak S}$ be a CF-perturbation of
$(X,\widehat{\mathcal U})$. We assume that $(f_t,f_P)$ is
stratumwise strongly submersive with respect to $\widehat{\frak S}$.
$\blacksquare$
\end{shitu}
\begin{defn}
We call $\frak X_P = ((X,\widehat{\mathcal U}),f_s,f_t,f_P)$ as in Situation \ref{situ2010}
a {\it $P$ parametrized family of smooth correspondences}.
\index{$P$-parametrized family of smooth correspondences}
\par
Let $\widehat{\frak S}$ be a CF-perturbation
such that $(f_t,f_P)$ is strongly submersive with respect to
$\widehat{\frak S}$.
Then for any $\epsilon >0$
we associate a linear map
$$
{\rm Corr}_{\frak X_P}(\cdot;\widehat{\frak S}^{\epsilon})
~:~ \Omega(P\times M_s) \longrightarrow \Omega(P\times M_t)
$$
by
\begin{equation}
{\rm Corr}_{\frak X_P}(h;\widehat{\frak S}^{\epsilon})
=
(f_P,f_t)!((f_P,f_s)^*h;\widehat{\frak S}^{\epsilon}).
\end{equation}
\end{defn}
Then we have
\begin{prop}\label{prop2012}
The map ${\rm Corr}_{\frak X_P}(\cdot;\widehat{\frak S}^{\epsilon})$ is pointwise in
$P$ direction.
\end{prop}
\begin{proof}
Let $h \in \Omega(M_s)$. We put
$$
(f_P,f_t)!((f_s)^*h;\widehat{\frak S}^{\epsilon})
= \sum_{I}dt_I
\wedge F_{I}(h).
$$
Let $F^{\bf t}_{I}(h)$ be the restriction of $F_{I}(h)$ to ${\{{\bf t}\}} \times M_t$.
Then it is easy to see that this $F^{\bf t}_{I}$ satisfies
(\ref{form208}) up to sign.
\end{proof}

\subsection{Proof of Theorem \ref{theorem1934}}
\label{subsec:20-4}

In this subsection, we complete the proof of Theorem \ref{theorem1934}.
\begin{proof}[Proof of Theorem \ref{theorem1934} (2)]
Suppose  $\mathcal{AF} = \{{\mathcal M}_{k+1}(\beta) \mid \beta,k\}$
defines a partial $A_{\infty}$ correspondence over $L$
of energy cut level $E_0$ and minimal energy $e_0$.
Let $G$ a discrete submonoid containing the discrete submonoid $G(\mathcal{AC})$
in Definition \ref{defassoSubmono}.
\begin{rem}
To prove Theorem \ref{theorem1934} (2)
itself, it suffices to take
$G=G(\mathcal{AC})$. However
we may also take $G$ which is strictly bigger than $G(\mathcal{AC})$.
We may replace $e_0$ by a smaller one.
\end{rem}
We apply Proposition \ref{prop203333} and
find a system of $\tau$-collared Kuranishi structures and CF-perturbations,
$\{(\widehat{\mathcal U^+_{k+1}}(\beta),\widehat{\frak S}_{k+1}(\beta)) \mid
(\beta,k) \in \mathcal{GK}(G;E_0,e_0)\}$.
We regard
$$
\left(\left(\mathcal M_{k+1}(\beta)^{\boxplus\tau_0},
\widehat{\mathcal U^+_{k+1}}(\beta)\right);
({\rm ev}_{1},\dots,{\rm ev}_{k}),{\rm ev}_0\right)
$$
as a smooth correspondence from
$L^k$ to $L$ and write it as
$\frak M_{k+1}(\beta)$.
We now define:
\begin{equation}\label{form2012}
\frak m^{\epsilon}_{k,\beta}(h_1,\dots,h_k)
:=
{\rm Corr}_{\frak M_{k+1}(\beta)}
\left(h_1\times\dots\times h_k;\widehat{\frak S}^{\epsilon}_{k+1}(\beta)
\right).
\end{equation}
Here and hereafter we denote
$$
h_1\times\dots\times h_k
:= \pi_1^*h_1 \wedge \dots \wedge \pi_k^*h_k
$$
where $\pi_i : L^k \to L$ is the $i$-th projection.
By Stokes' formula (\cite[Theorem 9.26]{part11})
we have
$$
\aligned
(d\circ \frak m^{\epsilon}_{k,\beta})(h_1,\dots,h_k)
&\pm
(\frak m^{\epsilon}_{k,\beta} \circ d)(h_1,\dots,h_k)  \\
&=
{\rm Corr}_{\partial\frak M_{k+1}(\beta)}
\left(h_1\times\dots\times h_k;\widehat{\frak S}^{\epsilon}_{k+1}(\beta)
\right).
\endaligned$$
We recall (\ref{formula166plus}), that is,
\begin{equation}\label{form2013}
\aligned
\partial ({\mathcal M}_{k+1}(\beta)^{\boxplus\tau_0},\widehat{\mathcal U^+_{k+1}}(\beta))
&\cong
\coprod_{\beta_1,\beta_2,k_1,k_2,i}
(-1)^{\epsilon}({\mathcal M}_{k_1+1}(\beta_1)^{\boxplus\tau_0},
\widehat{\mathcal U^+_{k_1+1}}(\beta_1)) \\
&\qquad\qquad
 \,\, {}_{{\rm ev}_i}\times_{{\rm ev}_0}
({\mathcal M}_{k_2+1}(\beta_2)^{\boxplus\tau_0},
\widehat{\mathcal U^+_{k_2+1}}(\beta_2)),
\endaligned
\end{equation}
where
\begin{equation}\label{eq:boundarysign33}
\epsilon = (k_1 -1)(k_2 -1) + \dim L + k_1 +
(i-1)\Big(1+(\mu(\beta_2)+k_2) \dim L\Big).
\end{equation}
We denote by $\frak M_{k_1,k_2,i}(\beta_1,\beta_2)$
the component corresponding to
$\beta_1,\beta_2,k_1,k_2,i$ in the right hand side
together with evaluation maps.
Note that the evaluation maps to the source of the left hand side
restrict to the evaluation maps to the source of either
the first or the second fiber product factor of the right hand side.
So our situation is (very slightly) different from one
of the composition formula \cite[Theorem 10.20]{part11}.
However we can apply \cite[Proposition 10.23]{part11}
instead by putting
$$
\aligned
&(X_1,\widehat{\mathcal U}_1,\widehat{\frak S}_1,\widehat
f_1)
= \left({\mathcal M}_{k_1+1}(\beta_1)^{\boxplus\tau_0},
\widehat{\mathcal U^+_{k_1+1}}(\beta_1),
\widehat{\frak S}^{\epsilon}_{k_1+1}(\beta_1),{\rm ev}_0
\right),
\\
&(X_2,\widehat{\mathcal U}_2,\widehat{\frak S}_2,\widehat
f_2)
= \left({\mathcal M}_{k_2+1}(\beta_2)^{\boxplus\tau_0},
\widehat{\mathcal U^+_{k_2+1}}(\beta_2),
\widehat{\frak S}^{\epsilon}_{k_2+1}(\beta_2),{\rm ev}_i
\right),
\\
&\widehat h_1 =
({\rm ev}_1,\dots,{\rm ev}_{k_1})^*(h_{i}\times\dots\times h_{i+k_1-1}),
\\
&\widehat h_2 =
({\rm ev}_1,\dots,{\rm ev}_{i-1},
{\rm ev}_{i+1}, \dots,{\rm ev}_{k_1},{\rm ev}_{0})^*
\\
&\qquad\qquad\qquad\qquad(h_1\times \dots
\times h_{i-1} \times
h_{i+k_1}\times \dots\times h_{k_2}\times h_0).
\endaligned
$$
Then \cite[(10.14)]{part11} and (\ref{formula166plus})=(\ref{form2013}) imply
\begin{equation}\label{2013form}
\aligned
&\int_L{\rm Corr}_{\partial\frak M_{k+1}(\beta)}
\left(h_1
\times\dots\times h_k;\widehat{\frak S}^{\epsilon}_{k+1}(\beta)
\right) \wedge h_0
\\
&=
\sum_{\beta_1,\beta_2,k_1,k_2,i}
{\rm Corr}_{\frak M_{k_2+1}(\beta_2)}
\left(\diamond;\widehat{\frak S}^{\epsilon}_{k_2+1}(\beta_2)
\right) \wedge h_0,
\endaligned
\end{equation}
where
$$
\aligned
\diamond
=
&h_1\times \dots
\times h_{i-1} \times \\
&\times
{\rm Corr}_{\frak M_{k_1+1}(\beta_1)}
\left(h_{i}\times\dots\times h_{i+k_1-1};
\widehat{\frak S}^{\epsilon}_{k_1+1}(\beta_1)
\right)
\times
h_{i+k_1}\times \dots\times h_{k_2}.
\endaligned
$$
(\ref{2013form}) implies that
$\{\frak m^{\epsilon}_{k,\beta}\}$ defines a $G$-gapped partial $A_{\infty}$-structure
of energy loss $E_0$ and minimal energy $e_0$.
\par
Thus we have constructed the required partial $A_{\infty}$-correspondence.
Its well-definedness up to
pseudo-isotopy will have been proved if Theorem \ref{theorem1934} (4) is proved.
\end{proof}
\begin{rem}
In the formulation of this article, we do not perturb $\frak m_{2,\beta_0}$.
(Recall $\beta_0 = 0$.)
Namely $\frak m_{2,\beta_0}$ coincides with the  wedge product up to sign and
$\frak m_{k,\beta_0} = 0$ for $k\ge 3$.
We take
$\widehat{\frak S}^{\epsilon}_{2+1}(\beta_0)$ as the trivial
perturbation.
Note that
${\mathcal M}_{2+1}(\beta_0) = L$ and
he evaluation map ${\rm ev}_0 : {\mathcal M}_{2+1}(\beta_0)
\to L$ is the identity map (that is a submersion).
So we do not need to perturb it.
We also take ${\mathcal M}_{k+1}(\beta_0) = L \times D^{k-2}$,
where we identify $D^{k-2}$ with the Stasheff cell.
Note that
${\rm ev}_0 : {\mathcal M}_{k+1}(\beta_0) \to L$
factors through the projection $L \times D^{k-2} \to L$
whose fiber is of positive dimension. Therefore this smooth correspondence
induces the zero map.
\par
We can proceed in a different way and
perturb ${\mathcal M}_{2+1}(\beta_0)$ so that  $\frak m_{2,\beta_0}$
has a smooth Schwartz kernel.
Then we necessarily include nonzero $\frak m_{k+1,\beta_0}$
for $k>2$.
We need to take such a choice of perturbation to generalize our story
to the case of bordered Riemann surfaces of higher genus
and/or those which have more than one boundary components,
because the corresponding moduli space of constant maps is not transversal.
\end{rem}

\begin{proof}[Proof of Theorem \ref{theorem1934} (4)]
We are given two
partial $A_{\infty}$ correspondences over $L$
$\mathcal{AF}^j = \{{\mathcal M}^j_{k+1}(\beta) \mid \beta,k\}$
$(j=0,1)$, and a pseudo-isotopy $\mathcal{AF}^{[0,1]} = \{{\mathcal M}_{k+1}(\beta;[0,1]) \mid \beta,k\}$ between them.
We assume that both of their energy cut levels are $E_0$ and minimal energies are $e_0$.
Let $G$ be a discrete submonoid containing both $G(\mathcal{AC}^j)$ for $j=0,1$.
We can make such a choice by Remark
\ref{rem1916}.
We also assume that $G$ contains $G(\mathcal{AF}^{[0,1]})$
and $e_0 \le e_{\rm min}(G)$.
\par
We assume that we have obtained partial filtered $A_{\infty}$ structures
$$
\{\frak m^{j,\epsilon_j}_{k,\beta} \mid (k,\beta) \in \mathcal{GK}(G;E_0,e_0)\}
$$
associated with the partial $A_{\infty}$ correspondences $\mathcal{AF}^j$
given in the proof of
Theorem \ref{theorem1934} (2) above.
It means that we have taken a system of
$$
\{(\widehat{\mathcal U^{j+}_{k+1}}(\beta),\widehat{\frak S}^j_{k+1}(\beta)) \mid
(\beta,k) \in \mathcal{GK}(G;E_0,e_0)\},
$$
where
$\widehat{\mathcal U^{j+}_{k+1}}(\beta)$ is a $\tau$-collared
Kuranishi structure on ${\mathcal M}^j_{k+1}(\beta)^{\boxplus\tau_0}$ and
$\widehat{\frak S}^j_{k+1}(\beta)$ is a
CF-perturbation of $\widehat{\mathcal U^{j+}_{k+1}}(\beta)$
such that they satisfy (\ref{form2013}).
Here recall that $\tau$ and $\tau_0$ satisfies the following inequality:
$$0 < \tau < \tau_0 =1.
$$
Now we apply Proposition \ref{prop203333rev} to
$P = [0,1]^{\boxplus(\tau_0 -\tau)}$
(then $P^{\boxplus\tau}=[0,1]^{\boxplus\tau_0}$)
to obtain objects
$$
\widehat{\mathcal U^{+}_{k+1}}(\beta;[0,1]), \quad
\widehat{\frak S}_{\rho,k+1}(\beta;[0,1])
$$
with $\rho \in (0,1]$ described below.
Firstly,
$\widehat{\mathcal U^{+}_{k+1}}(\beta;[0,1])$
is a $\tau$-collared Kuranishi structure
on ${\mathcal M}_{k+1}(\beta;[0,1])^{\boxplus\tau_0}$
with the following properties.
\begin{proper}
\begin{enumerate}
\item
Its restriction to
${\mathcal M}^j_{k+1}(\beta)^{\boxplus\tau_0} \subset
\partial ({\mathcal M}_{k+1}(\beta;[0,1])^{\boxplus\tau_0})$
is $\widehat{\mathcal U^{j+}_{k+1}}(\beta)$.
\item
Its restriction to
$
{\mathcal M}_{k_1+1}(\beta_1;[0,1])^{\boxplus\tau_0}
\, {}_{({\rm ev}_0,{\rm ev}_{[0,1]})}\times_{({\rm ev}_i,{\rm ev}_{[0,1]})}
{\mathcal M}_{k_2+1}(\beta_2;[0,1])^{\boxplus\tau_0}
$
is
$
\widehat{\mathcal U^{+}_{k_1+1}}(\beta_1;[0,1])
\, {}_{({\rm ev}_0,{\rm ev}_{[0,1]})}\times_{({\rm ev}_i,{\rm ev}_{[0,1]})}
\widehat{\mathcal U^{+}_{k_2+1}}(\beta_2;[0,1])
$.
\end{enumerate}
Note that $\widehat{\mathcal U^{j+}_{k+1}}(\beta)$ is a
$[0,1]^{\boxplus\tau_0} = [-\tau_0, 1+\tau_0]$-parametrized family.
\end{proper}
The $\tau$-collared
Kuranishi structure $\widehat{\mathcal U^{+}_{k+1}}(\beta;[0,1])$ also satisfies
the compatibility conditions at the corner. However, we do not
describe them here since they are special cases of the
statement of Proposition \ref{prop203333rev} and
we do not use them below directly.
\par
Secondly, $\widehat{\frak S}_{\rho,k+1}(\beta;[0,1])$
is a family of CF-perturbations
of $\widehat{\mathcal U^{+}_{k+1}}(\beta;[0,1])$
parametrized by $\rho \in (0,1]$
with the following properties.
\begin{proper}\label{prop2021}
\begin{enumerate}
\item
\begin{enumerate}
\item
Its restriction to
${\mathcal M}^j_{k+1}(\beta)^{\boxplus\tau_0} \subset
\partial {\mathcal M}_{k+1}(\beta;[0,1])^{\boxplus\tau_0}$ with $j=0$
is $\widehat{\frak S}^0_{k+1}(\beta)$.
\item
Its restriction to
${\mathcal M}^j_{k+1}(\beta)^{\boxplus\tau_0} \subset
\partial {\mathcal M}_{k+1}(\beta;[0,1])^{\boxplus\tau_0}$ with $j=1$
is $\epsilon \mapsto \widehat{\frak S}^{1 \rho\epsilon}_{k+1}(\beta)$.
\end{enumerate}
\item
Its restriction to
$
{\mathcal M}_{k_1+1}(\beta_1;[0,1])^{\boxplus\tau_0}
\, {}_{({\rm ev}_0,{\rm ev}_{[0,1]})}\times_{({\rm ev}_i,{\rm ev}_{[0,1]})}
{\mathcal M}_{k_2+1}(\beta_2;[0,1])^{\boxplus\tau_0}
$
is
$
{\widehat{\frak S}}_{\rho,k_1+1}(\beta_1;[0,1])
\, {}_{({\rm ev}_0,{\rm ev}_{[0,1]})}\times_{({\rm ev}_i,{\rm ev}_{[0,1]})}
{\widehat{\frak S}}_{\rho,k_2+1}(\beta_2;[0,1])
$.
\item
It is transversal to $0$.
The map $({\rm ev}_0,{\rm ev}_{[0,1]})$ is strongly
submersive with respect to $\widehat{\frak S}_{\rho,k+1}(\beta;[0,1])$.
\item
$\{\widehat{\frak S}_{\rho,k+1}(\beta;[0,1]) \mid \rho \in (0,1]\}$
is a uniform family.
\end{enumerate}
\end{proper}
We regard
$$
\left(\left(\mathcal M_{k+1}(\beta)^{\boxplus\tau_0},
\widehat{\mathcal U^+_{k+1}}(\beta;[0,1])\right);
\left(({\rm ev}_{1},{\rm ev}_{[0,1]}),\dots,
({\rm ev}_{k},{\rm ev}_{[0,1]})\right),
\left({\rm ev}_0,
{\rm ev}_{[0,1]}\right)\right)
$$
together with $\widehat{\frak S}_{\rho,k+1}(\beta;[0,1])$
as a smooth correspondence from
$(L \times [0,1]^{\boxplus\tau_0})^k$ to $L\times [0,1]^{\boxplus\tau_0}$ and write it as
$$
\frak M_{\rho,k+1}(\beta;[0,1]^{\boxplus\tau_0}).
$$
Now for differential forms $h_1,\dots,h_k$ on $L\times [0,1]^{\boxplus\tau_0}$ we put
$$
\frak m_{k,\beta}^{\rho,[0,1]^{\boxplus\tau_0}}(h_1,\dots,h_k)
=
{\rm Corr}_{\frak M_{\rho,k+1}(\beta;[0,1]^{\boxplus\tau_0})}
\left(h_1,\dots,h_k;{\widehat{\frak S}}^{\epsilon}_{\rho,k+1}(\beta;[0,1])\right).
$$
Using Property \ref{prop2021} (2)
and Stokes' formula (Theorem \ref{prop3210}) in the same way
as in the proof of Theorem \ref{theorem1934} (2),
we can prove
\begin{equation}\label{Ainfinityrelbeta2rev}
\aligned
&\sum_{k_1+k_2=k+1}\sum_{\beta_1+\beta_2=\beta}\sum_{i=1}^{k-k_2+1} \\
&(-1)^*{\frak m}^{\epsilon,\rho,[0,1]^{\boxplus\tau_0}}_{k_1,\beta_1}(h_1,\ldots,{\frak m}^{\epsilon,\rho,[0,1]^{\boxplus\tau_0}}_{k_2,\beta_2}(h_i,\ldots,h_{i+k_2-1}),\ldots,h_{k}) = 0,
\endaligned\end{equation}
where $* = \deg' h_1 + \ldots + \deg'h_{i-1}$.
\par
By Proposition \ref{prop2012},
$\frak m_{k,\beta}^{\epsilon,\rho,[0,1]^{\boxplus\tau_0}}$ is pointwise in
$[0,1]^{\boxplus\tau_0}$-direction.
Moreover,
Property \ref{prop2021} (1) (a) implies that
the restriction of
$\frak m_{k,\beta}^{\epsilon,\rho,[0,1]^{\boxplus\tau_0}}(h_1,\dots,h_k)$ to
$L \times \{-\tau\}$ is
$\frak m_{k,\beta}^{0,\epsilon}(h_1,\dots,h_k)$ and
Property \ref{prop2021} (1) (b) implies that
the restriction of
$\frak m_{k,\beta}^{\epsilon,\rho,[0,1]}(h_1,\dots,h_k)$ to
$L \times \{1+\tau\}$ is
$\frak m_{k,\beta}^{1,\epsilon\rho}(h_1,\dots,h_k)$.
Therefore $\frak m_{k,\beta}^{\epsilon,\rho,[0,1]^{\boxplus\tau_0}}$ gives a
partial $A_{\infty}$ pseudo-isotopy between
$\frak m_{k,\beta}^{0,\epsilon}$
and $\frak m_{k,\beta}^{1,\epsilon\rho}$.
\end{proof}
\begin{rem}
We introduced the parameter $\rho$ and
constructed a pseudo-isotopy between
$\frak m_{k,\beta}^{\epsilon,\rho,[0,1]^{\boxplus\tau_0}}$
and $\frak m_{k,\beta}^{1,\epsilon\rho}$
since we need it in the inductive construction
of a filtered $A_{\infty}$ structure associated to
an inductive system of $A_{\infty}$ correspondences.
See Remark \ref{rem171717}.
\end{rem}
\begin{proof}[Proof of Theorem \ref{theorem1934} (5)]
Recall from Definition \ref{inducAinfty} that we have a divergent sequence
$\{ E^i\}_i$ with
$$
0 < \cdots < E^i < E^{i+1} < \cdots \to +\infty
$$
of energy cut levels in the definition of the inductive system.
By (2) there exists a filtered $A_{\infty}$ structure of energy cut level
$E^i$ and minimal energy $e_0$ on $\Omega(L)$
induced by
$$
\mathcal{AF}^i = \{{\mathcal M}^i_{k+1}(\beta)\mid \beta,k\}.
$$
We denote it by $\{\frak m^i_k\}$.
By (4) there exists a partial pseudo-isotopy from  $\{\frak m^i_k\}$ to  $\{\frak m^{i+1}_k\}$.
Its energy cut level is $E^i$ and minimal energy is $e_0$.
It is induced by
$$
\mathcal{AF}^{[i,i+1]} = \{{\mathcal M}_{k+1}(\beta;[i,i+1])\mid \beta,k\}.
$$
We denote it by $\{\frak m^{[i,i+1]}_k\}$.
Then we can prove the following lemma by induction on $N$.
\begin{lem}\label{lem2024}
For each $n \le N$ we have the following.
\begin{enumerate}
\item
For $i \le n$,
there exists a promotion of $\{\frak m^i_k\}$ to the energy cut level $E^{n}$.
\item
For $i \le n-1$,
there exists a promotion of $\{\frak m^{[i,i+1]}_k\}$ to the energy cut level $E^n$.
It is a pseudo-isotopy between the promotions in (1).
\end{enumerate}
Moreover if $n' < n$ the structures in (1)(2) for $n'$ is the energy cut at $E^{n'}$ of
the structures in (1)(2) for $n$.
\end{lem}
\begin{proof}
This is immediate from Proposition \ref{prop208}.
\end{proof}
Now by mathematical induction we obtain the same conclusion as in Lemma \ref{lem2024} in the case $N=\infty$.
This implies Theorem \ref{theorem1934} (5).
Namely the filtered $A_{\infty}$ structure
associated to our inductive system of linear K-systems
is one $\{\frak m^{[0,1],1}_k \mid i=1,2,\dots\}$ obtained by promotion
to energy cut level $\infty$.
\end{proof}
\begin{proof}[Proof of Theorem \ref{theorem1934} (1)]
Using trivial pseudo-isotopy (the direct product), this is a special
case of Theorem \ref{theorem1934} (5).
\end{proof}
\begin{proof}[Proof of Theorem \ref{theorem1934} (6)]
Suppose we are in Situation \ref{situ1933} (6).
We apply Lemma \ref{lem2024} to each of the two inductive systems $\mathcal{IAF}^0$,
$\mathcal{IAF}^1$.
Namely we start with $\{\frak m^{ji}_k \mid i=1,2,\dots\}$ which are partial $A_{\infty}$
structures of energy cut level $E^i$
and minimal energy $e_0$ on $\Omega(L)$ (where $j=0,1$)
and with $\{\frak m^{j,[i,i+1]}_k \mid i=1,2,\dots\}$ which are partial pseudo-isotopies of energy cut level
$E^i$ and minimal energy $e_0$ among them.
Then by Lemma \ref{lem2024} and induction, we promote them to the energy cut level $\infty$.
\par
Next we consider $\{{\mathcal M}^i_{k+1}(\beta;[0,1])\mid \beta,k\}$ and
$\{{\mathcal M}_{k+1}(\beta;[0,1]\times [i,i+1])\mid \beta,k\}$
in Situation \ref{situ1933} (6).
\par
The former defines $\{\frak m^{[0,1],i}_k \mid i=1,2,\dots\}$ that is a pseudo-isotopy of  energy cut level $E^i$
and minimal energy $e_0$ from  $\{\frak m^{0i}_k \mid i=1,2,\dots\}$ to
$\{\frak m^{1i}_k \mid i=1,2,\dots\}$.
\par
The latter defines $\{\frak m^{[0,1],[i,i+1]}_k \mid i=1,2,\dots\}$
that is a pseudo-isotopy of pseudo isotopies.
In other words, it is a $([0,1] \times [i,i+1])$-parametrized family of partial
$A_{\infty}$ algebras of
energy cut level $E^i$ and minimal energy $e_0$ and
their restrictions to the normalized boundary are disjoint union of
$\{\frak m^{[0,1],i}_k \mid i=1,2,\dots\}$, $\{\frak m^{[0,1],i+1}_k \mid i=1,2,\dots\}$,
$\{\frak m^{0,[i,i+1]}_k \mid i=1,2,\dots\}$ and
$\{\frak m^{1,[i,i+1]}_k \mid i=1,2,\dots\}$.
\par
Now we apply Proposition \ref{prop208} and the same induction argument
as in the proof of Lemma \ref{lem2024}
to promote $\{\frak m^{[0,1],i}_k \mid i=1,2,\dots\}$ and $\{\frak m^{[0,1],[i,i+1]}_k \mid i=1,2,\dots\}$
to the energy cut level $\infty$.
Thus after promotion,
$\{\frak m^{[0,1],1}_k \mid i=1,2,\dots\}$ gives a
pseudo-isotopy from the promotion of $\{\frak m^{01}_k \mid i=1,2\dots\}$ to the promotion of $\{\frak m^{11}_k \mid i=1,2\dots\}$
of energy cut level $\infty$.
This is what we want to construct.
\par
We note that there exist $A_{\infty}$ homomorphisms
\begin{equation}\label{form2016}
\{\frak m^{[0,1],1}_k \mid i=1,2,\dots\} \longrightarrow
\{\frak m^{j1}_k \mid i=1,2\dots\}, \quad j=0,1,
\end{equation}
which are linear homotopy equivalences induced by the inclusion $L = L \times \{j\} \to L \times [0,1]$.
Therefore inverting one of them and using the Whitehead theorem for $A_{\infty}$ algebra
\cite[Theorem 4.2.45]{fooo08}, we obtain the
homotopy equivalence
$\{\frak m^{01}_k \mid i=1,2\dots\} \to \{\frak m^{11}_k \mid i=1,2\dots\}$.
\end{proof}
\begin{proof}[Proof of Theorem \ref{theorem1934} (3)]
This is a special case of Theorem \ref{theorem1934} (6).
\end{proof}
\begin{proof}[Proof of Theorem \ref{theorem1934} (8)]
For each $i$ we use $\mathcal{AF}^{[0,1] \times [i,i+1] \times [1,2]}$
to obtain a $[0,1] \times [i,i+1] \times [1,2]$ parametrized family of
partial $A_{\infty}$ structures of energy cut level $E^i$ and minimal energy $e_0$.
On $[0,1] \times [i,i+1] \times \{1\}$ and $[0,1] \times [i,i+1] \times \{2\}$
this family restricts to the family we obtain in the above proof of Theorem \ref{theorem1934} (6)
applied to $\{{\mathcal M}^{i,\ell}_{k+1}(\beta;[0,1])\mid \beta,k\}$,
$\{{\mathcal M}^{\ell}_{k+1}(\beta;[0,1]\times [i,i+1])\mid \beta,k\}$ with $\ell =a$ and $\ell = b$,
respectively.
Moreover it restricts
to the direct product on $\{j\} \times [i,i+1] \times [1,2]$ with  $j =0$ or $j=1$.
\par
Now applying Proposition \ref{prop208}, we use the same induction argument as in
the proof of
Lemma \ref{lem2024} to obtain at the part $i=1$ the
$[0,1]\times \{1\} \times [1,2]$ parametrized family of $A_{\infty}$ structures
of energy cut level $\infty$. We denote it by
$\{\frak m^{[0,1]\times \{1\} \times [1,2]}_k \mid i=1,2,\dots\} $.
We have a commutative diagram:
\begin{equation}\label{homohomohomo}
\begin{CD}
\{\frak m^{10}_k \}
@ <<< \{\frak m^{[0,1],1,b}_k\} @>>>  \{\frak m^{11}_k \}\\
@AAA   @AAA   @AAA \\
\{\frak m^{10}_k \}\times [1,2]
@<<<\{\frak m^{[0,1]\times \{1\} \times [1,2]}_k\} @>>>  \{\frak m^{11}_k \} \times [1,2]
\\
@VVV   @VVV   @VVV \\
\{\frak m^{10}_k \}
@<<<\{\frak m^{[0,1],1,a}_k\} @>>>  \{\frak m^{11}_k \}
\end{CD}
\end{equation}
The first and the third horizontal lines define homotopy equivalences we obtain
in Theorem \ref{theorem1934} (6) applied for $\ell =a$ and $\ell =b $ respectively.
(We invert the first arrow.)
\par
By the symbol $\{\frak m^{10}_k \}\times [1,2]$, we denote the direct product pseudo-isotopy,
that is, the isotopy $\{\frak m_k^t, \frak c_k^t\}$ such that $\frak m_k^t$ is
independent of $t$ and $\frak c_k^t$ are all zero.
Then by inverting one of the arrows in the first or third vertical lines we obtain identity maps.
Thus the commutativity of (\ref{homohomohomo}) implies that the two
homotopy equivalences obtained for  $\ell =a$ and $\ell =b $ are homotopic.
This is the conclusion of Theorem \ref{theorem1934} (8).
\end{proof}
\begin{proof}[Proof of Theorem \ref{theorem1934} (7)]
This is a special case of proof of Theorem \ref{theorem1934} (8).
\end{proof}
\begin{rem}
The proof we gave here using Diagram \eqref{homohomohomo} is similar to those
in \cite{akahojoyce}. It is based on a  way to define homotopy equivalence we took here, that is, to invert homotopy equivalence
(\ref{form2016}).
\par
There is an alternative way to define homotopy equivalence, given in \cite[Proof of Theorem 8.2]{fooo091},
where we take an appropriate integration and sum over trees to construct homotopy
equivalence from pseudo-isotopy directly.
This method has an advantage
that in case the pseudo-isotopy has cyclic symmetry the resulting homotopy equivalence is
also cyclic. (We do not know the version of \cite[Theorem 4.2.45]{fooo08} with cyclic symmetry.)
\par
We can use the same method to prove Theorem \ref{theorem1934} (7).
Namely we regard $\{\frak m^{[0,1]\times \{1\} \times [1,2]}_k\}$ as the
pseudo-isotopy from $\{\frak m^{10}_k \}\times [1,2]$ to
$\{\frak m^{11}_k \}\times [1,2]$ and apply the same formula
\cite[Definition 9.4]{fooo091}.
Using the fact that this is a pseudo-isotopy between direct product $A_{\infty}$ algebras,
we can easily check that the resulting filtered $A_{\infty}$ homomorphism becomes the
required homotopy.
\end{rem}
\begin{rem}
We stop here at the stage where we prove consistency up to homotopy of homotopies.
It is fairly obvious from the proof that we can prove as many higher consistency of the
homotopies as we want.
\end{rem}
\begin{rem}
We heard from some people that using iterated homotopies and homological algebra such as those we
developed in this section is cumbersome and had better be avoided.
Actually we do not think so at all.
\par
One of the origin of `homotopy everything structure' in algebraic topology is
to study `homotopy limit'.  So the language of $A_{\infty}$ structures which we are using here is
very much suitable for such a discussion.
Moreover, taking an inductive limit is necessary anyway to study, for example,
symplectic homology.
\par
Furthermore, the general strategy taken here does not use any of the special feature of the problem
and uses only the facts which are `intuitively obvious'. By this reason it should work in almost all the situations we meet and will
meet in the future. (Once we get used to it, applying this strategy becomes a routine.)
The general strategy is summarized as follows.
\par\medskip
\begin{enumerate}
\item
The problem is that it is hard and sometimes impossible to perturb infinitely many
moduli spaces simultaneously.
\item
Those moduli spaces are filtered by certain quantities, typically by energy.
\item
We fix certain `energy cut level' and perturb the (finitely many) moduli spaces only
up to that level.
\item
Then we obtain a partial structure, such as partial $A_{\infty}$ structure of energy cut level $E$.
\item
We may take $E$ as large as we want though $E$ should be finite.
\item
Let $E < E'$. we obtain the structures of energy cut level $E$ and of $E'$.
The later can be regarded as the structure of energy cut level $E$.
Those two structures are not the same but are `homotopy equivalent' in a sense of
`homotopy everything structure'.
\item
Then by a general method of homological algebra we can take homotopy limit to
obtain the desired structure.
\end{enumerate}
We note that a similar technique also appears in the renormalization theory.
Especially Costello \cite{Costello} uses a similar technique.
\end{rem}

\newpage
\part{Appendices}

\section{Orbifold and orbibundle by local coordinate}
\label{sec:ofd}

In this section we describe the story of orbifold as much as we
need in this article.
We restrict ourselves to effective orbifolds
and  regard only embeddings as morphisms.
The category $\mathcal{OB}_{\rm ef,em}$
where objects are effective
orbifolds and morphisms are embeddings among them
is naturally a $1$ category. Moreover
it has the following property.
We consider the forgetful map
$$
\frak{forget} : \mathcal{OB}_{\rm ef,em} \to \mathcal{TOP}
$$
where $\mathcal{TOP}$ is the category of topological spaces.
Then
$$
\frak{forget} : \mathcal{OB}_{\rm ef,em}(c,c')
\to \mathcal{TOP}(\frak{forget}(c),\frak{forget}(c'))
$$
is injective.
In other words, we can check the equality between morphisms
set-theoretically.
This is a nice property, which we use extensively in the main
body of this article.
If we go beyond this category, then we need to distinguish
carefully the two notions, two morphisms are equal,
two morphisms are isomorphic.
It will make the argument much more complicated.
\footnote{We need to use several maps between underlying topological
spaces of orbifolds, such as projection of bundles or covering maps.
In case we include those maps, we need to carefully examine whether
set-theoretical equality is enough to show their various properties
are preserved.}
\par
We emphasize that there is nothing new in this section.
The story of orbifold is  classical and is well-established.
It has been used in various branches of mathematics
since its invention by Satake \cite{satake} in more than 50 years ago.
Especially, if we restrict ourselves to
effective orbifolds, the story of orbifolds is nothing more than
a straightforward generalization of the theory of
smooth manifolds. The only important issue is the
observation that for effective orbifolds almost everything works
in the same way as manifolds.

\subsection{Orbifolds and embeddings between them}
\label{subsec;ofds}
\begin{defn}\label{2661}
Let $X$ be a paracompact Hausdorff space.
\begin{enumerate}
\item
An {\it orbifold chart of $X$ (as a topological space)}
\index{orbifold ! orbifold chart}
is a triple $(V,\Gamma,\phi)$ such that
$V$ is a manifold, $\Gamma$ is a finite group acting smoothly
and effectively on $V$ and $\phi : V \to X$
is a $\Gamma$ equivariant continuous map\footnote{The $\Gamma$ action
on $X$ is trivial.} which induces a
homeomorphism $\overline\phi : V/\Gamma \to X$ onto an open subset of $X$.
We assume that there exists $o \in V$ such that $\gamma o = o$
for all $\gamma \in \Gamma$.
We call $o$ the  {\it base point}.
\index{orbifold ! base point of an orbifold chart}
We say $(V,\Gamma,\phi)$ is an orbifold chart {\it at $x$} if
$x = \phi(o)$.
We call $\Gamma$ the {\it isotropy group},
\index{orbifold ! isotropy group of an orbifold chart}
$\phi$ the {\it local uniformization map}
\index{orbifold ! local uniformization map}  and
$\overline\phi$ the
{\it parametrization}
\index{orbifold ! parametrization}.
\item
Let $(V,\Gamma,\phi)$ be an orbifold chart and $p \in V$.
We put $\Gamma_p = \{ \gamma \in \Gamma \mid \gamma p = p\}$.
Let $V_p$ be a $\Gamma_p$ invariant open neighborhood of $p$
in $V$. We assume the map
$\overline\phi : V_p/\Gamma_p \to X$ is injective.
(In other words, we assume that
$\gamma V_{p} \cap V_p \ne \emptyset$ implies $\gamma \in \Gamma_p$.)
We call such a triple $(V_{p},\Gamma_p,\phi\vert_{V_p})$
a {\it subchart} \index{orbifold ! subchart} of $(V,\Gamma,\phi)$.
\item
Let $(V_i,\Gamma_i,\phi_i)$  $(i=1,2)$ be orbifold charts of $X$.
We say that they are {\it compatible}
\index{orbifold ! compatible (orbifold chart)} if the following holds
for each $p_1 \in V_1$ and $p_2 \in V_2$ with
$\phi_1(p_1) = \phi_2(p_2)$.
\begin{enumerate}
\item
There exists a group isomorphism $h : (\Gamma_1)_{p_1}
\to (\Gamma_2)_{p_2}$.
\item
There exists an $h$ equivariant diffeomorphism $\tilde\varphi : V_{1,p_1}
\to V_{2,p_2}$. Here $V_{i,p_i}$ is a $(\Gamma_i)_{p_i}$
equivariant subset of $V_i$ such that
 $(V_{i,p_i},(\Gamma_i)_{p_i},\phi\vert_{V_{i,p_i}})$
 is a subchart.
 \item
 $\phi_2 \circ \tilde\varphi = \phi_1$ on $V_{1,p_1}$.
\end{enumerate}
\item
A {\it representative of an orbifold structure}
\index{orbifold ! representative} on $X$
is a set of orbifold charts $\{(V_i,\Gamma_i,\phi_i) \mid i \in I\}$
such that each two of the charts are compatible in the sense of (3)
above and
$
\bigcup_{i\in I} \phi_i(V_i) = X,
$
is a locally finite open cover of $X$.
\end{enumerate}
\end{defn}
\begin{defn}\label{def262220}
Suppose that $X$, $Y$ are equipped with representatives
of orbifold structures $\{(V^X_i,\Gamma^X_i,\phi^X_i) \mid i \in I\}$
and $\{(V^Y_j,\Gamma^Y_j,\phi^Y_j) \mid j \in J\}$,
respectively.
A continuous map $f : X \to Y$ is said to be an
{\it embedding}
\index{orbifold ! embedding of orbifolds}\index{embedding ! of orbifolds}
if the following holds.
\begin{enumerate}
\item
$f$ is an embedding of  topological spaces.
\item
Let $p \in V^X_i$, $q \in V^Y_j$ with
$f(\phi_i(p)) = \phi_j(q)$. Then we have the following.
\begin{enumerate}
\item
There exists an isomorphism of groups
$h_{p;ji} : (\Gamma_i^X)_p \to (\Gamma_j^Y)_q$.
\item
There exist  $V^X_{i,p}$ and $V^Y_{j,q}$ such that
$(V^X_{i,p},(\Gamma^X_i)_p,\phi_i\vert_{V^X_{i,p}})$
is a subchart for $i=1,2$.
There exists an $h_{p;ji}$ equivariant embedding of manifolds
$\tilde f_{p;ji}: V^X_{i,p} \to V^Y_{j,q}$.
\item
The diagram below commutes.
\begin{equation}\label{diag2633}
\begin{CD}
V^X_{i,p} @ >{\tilde f_{p;ji}}>>
V^Y_{j,q}  \\
@ V{\phi_{i}}VV @ VV{\phi_{j}}V\\
X @ > {f} >> Y
\end{CD}
\end{equation}
\end{enumerate}
\end{enumerate}
Two orbifold embeddings are said to be {\it equal} if they coincide set-theoretically.
\end{defn}
\begin{rem}
Note that an embedding of effective orbifolds
is a continuous map $f : X \to Y$ of
underlying topological spaces, which has
the properties (2) above.
\par
When we study morphisms among ineffective
orbifolds or morphisms between effective orbifolds which is not
an embedding, such a morphism is a
continuous map $f : X \to Y$ of underlying topological spaces
{\it plus} certain additional data.
For example, if we consider an ineffective orbifold
that is a point with an action of a nontrivial finite group $\Gamma$,
then the morphism from this ineffective orbifold to itself
contains the datum of an automorphism of the group $\Gamma$.
(Two such morphisms $h_1,h_2$ are equivalent if there exists an
inner automorphism $h$ such that $h_1 = h \circ h_2$.)
\end{rem}

\begin{lem}\label{lem26444}
\begin{enumerate}
\item
The composition of embeddings is an embedding.
\item
The identity map is an embedding.
\item
If an embedding is a homeomorphism,
then its inverse is also an embedding.
\end{enumerate}
\end{lem}
The proof is easy and is left to the reader.

\begin{defn}\label{defn285555}
\begin{enumerate}
\item We call an embedding of orbifolds a {\it diffeomorphism}
\index{orbifold ! diffeomorphism} if it is a
homeomorphism in addition.
\item
We say that two representatives of orbifold structures on $X$ are
{\it equivalent} if the identity map regarded as a map between
$X$ equipped with those two representatives of orbifold structures
is a diffeomorphism. This is an equivalence relation by Lemma \ref{lem26444}.
\item
An equivalence class of the equivalence relation (2) is called
an {\it orbifold structure} on $X$.
\index{orbifold ! orbifold structure}
An {\it orbifold} is a pair of topological space and its
orbifold structure. \index{orbifold ! orbifold}
\item
The condition for a map $X \to Y$ to be an embedding
does not change if we replace representatives of orbifold structures
to equivalent ones. So we can define the notion of
an {\it embedding of orbifolds}.
\item
If $U$ is an open subset of an orbifold $X$, then
there exists a unique orbifold structure on $U$ such that
the inclusion $U \to X$ is an embedding.
We call $U$ with this orbifold structure an {\it open suborbifold}.
\index{orbifold ! open suborbifold}
\end{enumerate}
\end{defn}
\begin{defn}\label{defn26550}
\begin{enumerate}
\item
Let $X$ be an orbifold. An orbifold chart $(V,\Gamma,\phi)$
of underlying topological space $X$
in the sense of Definition \ref{2661} (1)
is called an {\it orbifold chart of an orbifold}  $X$ if
the map $\overline\phi : V/\Gamma \to X$ induced by $\phi$ is an
embedding of orbifolds.
\index{orbifold ! orbifold chart}
\item
Hereafter when $X$ is an orbifold,
an `orbifold chart' always means an orbifold chart of an orbifold in the sense of
(1).
\item
In case when an orbifold structure on $X$ is given,
a representative of its orbifold structure is
called an {\it orbifold atlas}.
\index{orbifold ! orbifold atlas}
\item
Two orbifold charts  $(V_i,\Gamma_i,\phi_i)$
are said to be {\it isomorphic}
if there exist a
group isomorphism $h : \Gamma_1 \to \Gamma_2$
and an $h$-equivariant diffeomorphism $\tilde\varphi : V_1 \to V_2$
such that $\phi_2 \circ \tilde\varphi = \phi_1$.
The pair $(h,\tilde\varphi)$ is called an
{\it isomorphism} or a {\it coordinate change}
\index{orbifold ! isomorphism of orbifold charts}
\index{orbifold ! coordinate change of orbifold charts}
between
two orbifold charts.
\end{enumerate}
\end{defn}
\begin{prop}\label{prop266}
In the situation of Definition \ref{defn26550} (4), suppose
$(h,\tilde\varphi)$ and $(h',\tilde\varphi')$ are both
isomorphisms between
two orbifold charts $(V_1,\Gamma_1,\phi_1)$ and $(V_2,\Gamma_2,\phi_2)$.
Then there exists $\mu \in \Gamma_2$
such that
\begin{equation}\label{262}
h'(\gamma) = \mu h(\gamma) \mu^{-1},
\qquad
\tilde\varphi'(x) = \mu \tilde\varphi(x).
\end{equation}
\par
Conversely, if $(h,\tilde\varphi)$ is an isomorphism between
orbifold charts, then $(h',\tilde\varphi')$ defined by (\ref{262})
is also an isomorphism between
orbifold charts.
In particular, any automorphism of an orbifold chart
$(h,\tilde\varphi)$ is given
by $h(\gamma) = \mu\gamma \mu^{-1}$, $\tilde\varphi(x) = \mu x$ for some element
$\mu \in \Gamma$.
\end{prop}
\begin{proof}
The proposition immediately follows from the next lemma.
\begin{lem}\label{lem21tenhatena}
Let $V_1$, $V_2$ be manifolds on which finite groups $\Gamma_1$,
$\Gamma_2$ act effectively and smoothly.
Assume that $V_1$ is connected.
Let $(h_i,\tilde\varphi_i)$ $(i=1,2)$ be pairs such that
$h_i : \Gamma_1 \to \Gamma_2$ are injective group homomorphisms
and $\tilde\varphi_i : V_1 \to V_2$ are $h_i$-equivariant embeddings of manifolds.
Moreover, we assume that the induced maps ${\varphi}_i :V_1/\Gamma_1
\to V_2/\Gamma_2$
are embeddings of orbifolds and
${\varphi}_1$ coincides with ${\varphi}_2$ set-theoretically.
Then there exists $\mu \in \Gamma_2$ such that
$$
\tilde\varphi_2(x) = \mu\tilde\varphi_1(x),
\qquad
h_2(\gamma) = \mu h_1(\gamma)\mu^{-1}.
$$
\end{lem}
\begin{proof}
For the sake of simplicity we prove only
the case when Condition \ref{convinv} below is satisfied.
Let $X$ be an orbifold.
For a point $x \in X$ we take its orbifold chart $(V_x,\Gamma_x,\psi_x)$.
We say $x \in {\rm Reg}(X)$ if
$\Gamma_x = \{1\}$,
and put ${\rm Sing}(X) = X \setminus {\rm Reg}(X)$.
\begin{conds}\label{convinv}
We assume that $\dim {\rm Sing}(X) \le \dim X -2$.
\end{conds}
This condition is satisfied if $X$ is oriented.  (In fact, Condition \ref{convinv} fails only
when there exists an element of $\Gamma_x$ (an isotropy group of some orbifold chart)
whose action is given by $(x_1,x_2,\dots,x_n) \mapsto (-x_1,x_2,\dots,x_n)$ for some coordinate $(x_1,\dots,x_n)$.
Therefore we can always assume Condition \ref{convinv}
in the study of Kuranishi structure, by adding a trivial factor which is acted by
the induced representation of $t \mapsto -t$ to both the obstruction bundle and to the Kuranishi neighborhood.)
\par
Let $x_0 \in V_1^0$.
By assumption there exists a unique $\mu \in \Gamma_2$
such that $\tilde\varphi_2(x_0) = \mu\tilde\varphi_1(x_0)$.
By Condition \ref{convinv} the subset $V^0_1$ is connected.
Therefore the above element $\mu$ is independent of $x_0 \in V_1^0$
by uniqueness.
Since $V_1^0$ is dense, we conclude $\tilde\varphi_2(x) = \mu\tilde\varphi_1(x)$ for any $x \in V_1$.
Now, for $\gamma \in \Gamma_1$, we calculate
$$
h_1(\gamma)\tilde\varphi_1(x_0)
=
\tilde\varphi_1(\gamma x_0)
=
\mu^{-1}\tilde\varphi_2( \gamma x_0)
=
\mu^{-1}h_2(\gamma)\tilde\varphi_2(x_0)
=
\mu^{-1}h_1(\gamma)\mu \tilde\varphi_1(x_0).
$$
Since the induced map is an embedding of orbifold,
it follows that the isotropy group of $\tilde\varphi_1(x_0)$
is trivial. Therefore
$h_1(\gamma) = \mu^{-1}h_2(\gamma)\mu$
as required.
\end{proof}
The proof of Proposition \ref{prop266} is complete.
\end{proof}
\begin{defn}\label{defn281010}
Let $X$ be an orbifold.
\begin{enumerate}
\item
A function $f : X \to \R$ is said to be a {\it smooth function} if
for any orbifold chart $(V,\Gamma,\phi)$ the composition
$f\circ \phi : V \to \R$ is smooth.
\index{orbifold ! smooth function on orbifold}
\item
A {\it differential form}
\index{orbifold ! differential form on orbifold}
\index{differential form ! on oribifold}
on an orbifold $X$ assigns a
$\Gamma$ invariant differential form $h_{\frak V}$ on $V$ to
each orbifold chart $\frak V = (V,\Gamma,\phi)$
such that the following holds.
\begin{enumerate}
\item
If $(V_1,\Gamma_1,\phi_1)$ is isomorphic to
$(V_2,\Gamma_2,\phi_2)$ and $(h,\tilde\varphi)$ is an isomorphism,
then $\tilde\varphi^*h_{\frak V_2} = h_{\frak V_1}$.
\item
If $\frak V_p = (V_p,\Gamma_p,\phi_p)$ is a subchart of $\frak V =(V,\Gamma,\phi)$,
then $h_{\frak V}\vert_{V_p} = h_{\frak V_p}$.
\end{enumerate}
\item
An $n$ dimensional orbifold $X$ is said to be {\it orientable} if there exists
a differential $n$-form $\omega$ such that $\omega_{\frak V}$
never vanishes.
\index{orbifold ! orientation}
\item
Let $\omega$ be an $n$-form as in (3) and
$\frak V = (V,\Gamma,\phi)$  an orbifold chart. Then
we give $V$ an orientation so that it is compatible with $\omega_{\frak V}$.
The $\Gamma$ action preserves the orientation.
We call such $(V,\Gamma,\phi)$ equipped with an orientation of $V$,
an {\it oriented orbifold chart}.
\index{orbifold ! oriented orbifold chart}
\item
Let $\bigcup_{i\in I}U_i = X$ be an open covering of an orbifold $X$.
A {\it smooth partition of unity subordinate to the covering} $\{U_i\}$
\index{orbifold ! partition of unity} is
a set of functions $\{\chi_i\mid i\in I\}$ such that:
\begin{enumerate}
\item
$\chi_i$ are smooth functions.
\item
The support of $\chi_i$ is contained in $U_i$.
\item
$\sum_{i\in I}\chi_i = 1$.
\end{enumerate}
\end{enumerate}
\end{defn}
\begin{lem}
For any locally finite open covering of an orbifold $X$ there exists a smooth
partition of unity subordinate thereto.
\end{lem}
We omit the proof, which is an obvious analogue of the standard proof
for the case of manifolds.

\begin{defn}\label{orbifolddefn}
An {\it orbifold with corner} is defined in the same way.
\index{orbifold ! orbifold with corner}
We require the following.
\begin{enumerate}
\item In Definition \ref{2661} (1) we assume that $V$ is a manifold with
corners.
\item
Let $S_k(V)$ be the set of points which lie on the codimension $k$ corner
and $\overset{\circ}S_k(V) = S_k(V) \setminus \bigcup_{k' > k}S_{k'}(V)$.
We require that $\Gamma$ action on each connected component of
$\overset{\circ}S_k(V)$ is effective.
(Compare \cite[Condition 4.14]{part11}.)
\item
For an embedding of orbifolds with corners
\index{orbifold ! embedding of orbifolds with corners}
\index{embedding ! of orbifolds with corners}
we require that the map
$\tilde f$ in
Definition \ref{def262220} (c) satisfies
$\tilde f(\overset{\circ}S_k(V_1)) \subset \overset{\circ}S_k(V_2)$.
\end{enumerate}
\end{defn}

\begin{lem}\label{lem26999}
Let $X_i$ $(i=1,2)$ be orbifolds and $\varphi_{21} : X_1 \to X_2$  an embedding.
Then we can find an orbifold atlas
$\{\frak V^i_{\frak r} =
\{ (V_{\frak r}^i, \Gamma_{\frak r}^i, \phi_{\frak r}^i )\} \mid \frak r \in \frak R_i\}$ with the following properties.
\begin{enumerate}
\item $\frak R_1 \subseteq \frak R_2$.
\item  $V^2_{\frak r} \cap \varphi_{21}(X_1) \ne \emptyset$
if and only if $\frak r \in \frak R_1$.
\item
If $\frak r \in \frak R_1$ then
$\varphi_{21}^{-1}(\phi_{\frak r}^2(V^2_{\frak r})) =
\phi_{\frak r}^1(V^1_{\frak r})$
and there exists
$(h_{\frak r,21},\tilde\varphi_{\frak r,21})$
such that:
\begin{enumerate}
\item
$h_{\frak r,21} : \Gamma^1_{\frak r} \to \Gamma^2_{\frak r}$
is a group isomorphism.
\item
$\tilde\varphi_{\frak r,21} : V^1_{\frak r} \to V^2_{\frak r}$
is an $h_{\frak r,21}$-equivariant embedding of smooth manifolds.
\item
The next diagram commutes.
\begin{equation}\label{diagin2611}
\begin{CD}
V^1_{\frak r} @ >{\tilde\varphi_{\frak r,21}}>>
V^2_{\frak r}  \\
@ V{\phi^{\frak r}_{1}}VV @ VV{\phi^{\frak r}_{2}}V\\
X_1 @ > {\varphi_{21}} >> X_2
\end{CD}
\end{equation}
\end{enumerate}
\item
In case $X_i$ has boundary or corners we may choose
our charts so that the following is satisfied.
\begin{enumerate}
\item
$V^i_{\frak r}$ is an open subset of $\overline V^i_{\frak r}
\times [0,1)^{d({\frak r})}$, where $d(\frak r)$ is independent of $i$
and $\overline V^i_{\frak r}$ is a manifold without boundary.
\item
There exists a point $o^i(\frak r)$ which is fixed by all
$\gamma \in \Gamma^i_{\frak r}$ such that
$[0,1)^{d({\frak r})}$ components of $o^i(\frak r)$ are all $0$.
\item If we write
$$
\varphi_{\frak r,21}(\overline y',(t'_1,\dots,t'_{d({\frak r})}))
=
(\overline y,(t_1,\dots,t_{d({\frak r})})),
$$
then $t_i = 0$ if and only if $t'_i = 0$.
\end{enumerate}
\end{enumerate}
We may take our atlas that are refinements of the given coverings of $X_1$ and $X_2$.
\end{lem}
\begin{proof}
For each $x \in X_1$ we can find  orbifold charts $\frak V^i_{x}$
for $i=1,2$, such that
$\varphi_{21}^{-1}(U^2_{x}) = U^1_{x}$,
$x \in U^1_{x}$ and that
 there exists a representative
$(h_{x,21},\tilde\varphi_{x,21})$
of embedding $U^1_{x} \to U^2_{x}$ that is a restriction
of $\varphi_{21}$.
In case $X_i$ has boundary or corners,
we choose them so that (4) is also satisfied.
\par
We cover $X_1$ by finitely many such $U_{x_j}^{1}$.
This is our choice of atlas
$\{\frak V^1_{\frak r} \mid \frak r \in \frak R_1\}$.
Then the associated $\{\frak V^2_{\frak r} \mid \frak r \in \frak R_1\}$ satisfies
(3)(4) and covers
$\varphi_{21}(X_1)$. We can extend it to
$\{\frak V^2_{\frak r} \mid \frak r \in \frak R_2\}$
so that (1)(2) are also satisfied.
\end{proof}
\begin{defn}
We call $(h_{\frak r,21},\tilde\varphi_{\frak r,21})$
a {\it local representative of embedding}  $\varphi_{\frak r,21}$ on the charts
$\frak V^1_{\frak r}$, $\frak V^2_{\frak r}$.
\index{orbifold ! local representative of embedding}
\index{embedding ! local representative of embedding of orbifolds}
\end{defn}
\begin{lem}
If $(h_{\frak r,21},\tilde\varphi_{\frak r,21})$, $(h'_{\frak r,21},\tilde\varphi'_{\frak r,21})$
are local representatives of an embedding of the same
charts
$\frak V^1_{\frak r}$, $\frak V^2_{\frak r}$,
then there exists $\mu \in \Gamma_2$
such that
$$
\tilde\varphi'_{\frak r,21}(x) = \mu\tilde\varphi_{\frak r,21}(x),
\qquad
h_{\frak r,21}^{\prime}(\gamma) = \mu h_{\frak r,21}(\gamma)\mu^{-1}.
$$
\end{lem}
This is a consequence of Lemma \ref{lem21tenhatena}.
\begin{lem}\label{smoothstruemb}
Let $X$ be a topological space, $Y$  an orbifold,
and $f : X \to Y$  an embedding of topological spaces.
Then the orbifold structure on $X$ by which $f$ becomes an
embedding of orbifolds is unique if there exists one.
\end{lem}
\begin{proof}
Let $X_1$, $X_2$ be orbifolds whose underlying topological spaces
are both $X$ and satisfy that $f_i : X_i \to Y$ are embeddings of orbifolds for
$i=1,2$. We will prove that the identity map
$\rm{id} : X_1 \to X_2$ is a diffeomorphism
of orbifods.
Since the condition for a homeomorphism
to be a diffeomorphism of orbifolds is a local condition,
it suffices to check it on a neighborhood of each point.
Let $p \in X$ and $q = f(p)$.
We take a representative $(h_i,\tilde\varphi_i)$ of the orbifold embeddings
$f_i : X_i \to Y$ using the orbifold charts
$\frak V^i_p = (V^i_p,\Gamma^i_p,\phi^i_p)$ of $X$ and
$\frak V_q = (V_q,\Gamma_q,\phi_q)$ of $Y$.
The maps $h_i : \Gamma^i_p \to \Gamma^i_q$ are group isomorphisms.
So we have a group isomorphism
$
h = h_2^{-1}\circ h_1 : \Gamma^1_p \to \Gamma^2_p.
$
Since $\tilde\varphi_1(V^1)/\Gamma_p = \tilde\varphi_2(V^2)/\Gamma_p$
set-theoretically, we have
$\tilde\varphi_1(V^1_p) = \tilde\varphi_2(V^2_p) \subset V_q$.
They are smooth submanifolds since $f_i$ are  embeddings of
orbifolds. Therefore
$
\varphi = \tilde\varphi_2^{-1}\circ \tilde\varphi_1
$
is defined in a neighborhood of the base point $o^i_p$
and is a diffeomorphism.
Then $(h,\tilde\varphi)$ is a local representative of ${\rm id}$.
\end{proof}

\subsection{Vector bundle on orbifold}
\label{subsec:vectorbundle}

\begin{defn}\label{defn2613}
Let $(X,\mathcal E,\pi )$ be a pair of orbifolds $X$ and $\mathcal E$ with
a continuous map $\pi : \mathcal E \to X$
between their underlying topological spaces.
Hereafter we write $(X,\mathcal E)$ in place of $(X,\mathcal E,\pi)$.
\begin{enumerate}
\item
An {\it orbifold chart} of $(X,\mathcal E)$ is
a quintuple $(V,E,\Gamma,\phi,\widehat\phi)$
with the following properties:
\index{orbifold ! orbifold chart of a vector bundle}
\begin{enumerate}
\item
$\frak V = (V,\Gamma,\phi)$ is an orbifold chart of the orbifold $X$.
\item $E$ is a finite dimensional vector space equipped with
a linear $\Gamma$
action.
\item
$(V \times E,\Gamma,\widehat\phi)$ is an orbifold chart of the
orbifold $\mathcal E$.
\item
The diagram below commutes set-theoretically.
\begin{equation}\label{diag26399}
\begin{CD}
V \times E @ >{\widehat\phi}>>
\mathcal E  \\
@ V{}VV @ VV{{\pi}}V\\
V @ > {\phi} >> X
\end{CD}
\end{equation}
Here the left vertical arrow is the projection to the
first factor.
\end{enumerate}
\item
In the situation of (1), let $p\in V$ and $( V_{p},\Gamma_p,\phi\vert_{V_p})$
be a subchart of $(V,\Gamma,\phi)$ in the sense of
Definition \ref{2661} (2).
Then  $(V_{p},E,\Gamma_p,\phi\vert_{V_p},\widehat\phi\vert_{V_p \times E})$
is an orbifold chart of $(X,\mathcal E)$.
We call it a {\it subchart} of $(V,E,\Gamma,\phi,\widehat\phi)$.
\index{orbifold ! subchart of an orbifold chart of vector bundle}
\item
Let $(V^i,E^i,\Gamma^i,\phi^i,\widehat{\phi^i})$  $(i=1,2)$ be orbifold charts of $(X,\mathcal E)$.
We say that they are {\it compatible} if the following holds
for each $p_1 \in V^1$ and $p_2 \in V^2$ with
$\phi^1(p_1) = \phi^2(p_2)$:
There exist open neighborhoods $V^i_{p_i}$ of $p_i \in V^i$ such that:
\par
\begin{enumerate}
\item
There exists an isomorphism
$(h,\tilde\varphi) : (V^1,\Gamma^1,\phi^1)\vert_{V^1_{p_1}} \to (V^2,\Gamma^2,\phi^2)\vert_{V^2_{p_2}}$
between orbifold charts of $X$, which are subcharts.
\item
There exists an isomorphism
$(h,\tilde{\hat{\varphi}}) : (V^1\times E^1,\Gamma^1,\phi^1)\vert_{V^1_{p_1}\times E^1} \to
(V^2\times E^2,\Gamma^2,\phi^2)\vert_{V^2_{p_2} \times E^2}$
between orbifold charts of $\mathcal E$, which are subcharts.
\item
For each $y \in V^1_{p_1}$ the map $ E^1 \to E^2$ given by
$\xi \to \pi_{E^2}\tilde{\hat{\varphi}}(y,\xi)$ is a linear isomorphism.
Here $\pi_{E^2} : V^2 \times E^2 \to E^2$ is the projection.
\end{enumerate}
\item
A {\it representative of a vector bundle structure} on $(X,\mathcal E)$
is a set of orbifold charts $\{(V_i,E_i,\Gamma_i,\phi_i,\widehat\phi_i) \mid i \in I\}$
such that any two of the charts are compatible in the sense of (3)
above and
$$
\bigcup_{i\in I} \phi_i(V_i) = X,
\quad
\bigcup_{i\in I} \widehat\phi_i(V_i \times E_i) = \mathcal E,
$$
are locally finite open covers.
\index{orbifold ! representative of a vector bundle structure on orbifold}
\end{enumerate}
\end{defn}
\begin{defn}\label{def26222}
Suppose $(X^*,\mathcal E^*)$ $(* = a,b)$ have representatives
of vector bundle structures $\{(V^*_i,E^*_i,\Gamma^*_i,\phi^*_i,\widehat\phi^*_i) \mid i \in I^*\}$,
respectively.
A pair of orbifold embeddings $(f,\widehat f)$,
$f : X^a \to X^b$, $\widehat f : \mathcal E^a \to \mathcal E^b$ is said to be an
{\it embedding of vector bundles} if the following holds.
\index{orbifold ! embedding of vector bundles on orbifolds}
\index{embedding ! of vector bundles on orbifolds}
\begin{enumerate}
\item
Let
$p \in V^a_i$, $q \in V^b_j$ with
$f(\phi^a_i(p)) = \phi^b_j(q)$.
Then there exist open subcharts
$(V^a_{i,p}\times E^a_{i,p},\Gamma^a_{i,p},\widehat\phi^a_{i,p})$
and
$(V^b_{j,q}\times E^b_{j,q},\Gamma^b_{j,q}\widehat\phi^b_{j,q})$
and a local representative
$(h_{p;i,j},f_{p;i,j},\widehat f_{p;i,j})$ of the embeddings $f$ and $\widehat f$
such that
for each $y \in V^a_i$ the map
$\xi \mapsto \pi_{E^b}(\widehat f_{p;i,j}(y,\xi))$,
$E^a_{i,p} \to E^b_{j,q}$ is a linear embedding. Here $\pi_{E^b} : V^b \times E^b \to E^b$ is the projection.
\item
The diagram below commutes set-theoretically.
\begin{equation}\label{diag2633}
\begin{CD}
\mathcal E^a @ >{\widehat f}>>
\mathcal E^b  \\
@ V{\pi_{E^a}}VV @ VV{\pi_{E^b}}V\\
X^a @ > {f} >> X^b
\end{CD}
\end{equation}
\end{enumerate}
Two orbifold embeddings of vector bundles are said to be {\it equal}
if they coincide set-theoretically
as pairs of maps.
\end{defn}
\begin{lem}\label{lem26444AA1}
\begin{enumerate}
\item
A composition of embeddings of vector bundles is an embedding.
\item
The pair of identity maps $({\rm id}, \widehat{\rm id})$
is an embedding.
\item
If an embedding of vector bundles is a pair of homeomorphisms,
then the pair of their inverses is also an embedding.
\end{enumerate}
\end{lem}
The proof is easy and is omitted.
\begin{defn}\label{defn2820}
Let
$(X,\mathcal E)$ be as in Definition \ref{defn2613}.
\begin{enumerate}
\item
An embedding of vector bundles is said to be an {\it isomorphism}
\index{orbifold ! isomorphism of vector bundles on orbifolds}
if it is a pair of diffeomorphisms of orbifolds.
\item
We say that two representatives of a vector bundle structure on $(X,\mathcal E)$
are
{\it equivalent} if the pair of identity maps regarded as a self-map of
vector bundle
$(X,\mathcal E)$ equipped with those two representatives of vector bundle
structure
is an isomorphism.
This is an equivalence relation by Lemma \ref{lem26444AA1}.
\item
An equivalence class of the equivalence relation (2) is called
a {\it vector bundle structure} on $(X,\mathcal E)$.
\index{orbifold ! vector bundle structure}
\item
A pair $(X,\mathcal E)$ together with its vector bundle
structure is called a {\it vector bundle} on $X$.
\index{orbifold ! vector bundle on orbifold}
\index{orbifold ! total space of vector bundle}
\index{orbifold ! base space of vector bundle}
\index{orbifold ! projection of vector bundle}
We call $\mathcal E$ the {\it total space}, $X$ the
{\it base space}, and $\pi : \mathcal E \to X$ the {\it projection}.
\item
The condition for the pair $(f,\widehat f) : (X^a,\mathcal E^a) \to (X^b,\mathcal E^b)$
to be an embedding
depends only on the equivalence class of vector bundle structures
independent of its representatives.
This enable us to define the notion of
an {\it embedding of vector bundles}.
%
\item
We say $(f,\widehat f)$ is an embedding {\it over the orbifold embedding $f$.}
\end{enumerate}
\end{defn}
\begin{rem}
\begin{enumerate}
\item
We may use the terminology `orbibundle' in place of vector bundle.
We use this terminology in case we emphasize that it is different
from the vector bundle over the underlying topological space.
\index{orbifold ! orbibundle}
\item
We sometimes simply say $\mathcal E$ is a vector bundle on an
orbifold $X$.
\end{enumerate}
\end{rem}
\begin{defn}\label{defn2655}
\begin{enumerate}
\item
Let
$(X,\mathcal E)$ be a vector bundle.
We call an orbifold chart $(V,E,\Gamma,\phi,\widehat\phi)$
in the sense of Definition \ref{defn2613} (1)
of underlying pair of topological spaces $(X,\mathcal E)$
an {\it orbifold chart of a vector bundle} $(X,\mathcal E)$ if
the pair of maps $(\overline\phi,\overline{\widehat\phi}) : (V/\Gamma,(V\times E)/\Gamma) \to (X,\mathcal E)$ induced from $(\phi,\widehat\phi)$ is an
embedding of vector bundles.
\item
If $(V,E,\Gamma,\phi,\widehat\phi)$
is an orbifold chart of a vector bundle,
we call a pair $(E,\widehat\phi)$
a {\it trivialization} \index{orbifold ! trivialization of vector bundle} of our vector bundle on $V/\Gamma$.
\item
Hereafter when $(X,\mathcal E)$ is a vector bundle,
its `orbifold chart' always means an orbifold chart of a vector bundles in the
sense of (1).
\item
In case when a vector bundle structure on $(X,\mathcal E)$ is given,
a representative of this vector bundle structure is
called an {\it orbifold atlas} of $(X,\mathcal E)$.
\index{orbifold ! orbifold atlas}
\item
Two orbifold charts $(V_i,E_i,\Gamma_i,\phi_i,\widehat \phi_i)$ of
a vector bundle
are said to be {\it isomorphic} if there exist an
isomorphism $(h,\tilde\varphi)$ of orbifold charts
$(V_1,\Gamma_1,\phi_1) \to (V_2,\Gamma_2,\phi_2)$
and an isomorphism
$(h,\tilde{\hat{\varphi}})$ of orbifold charts
$(V_1\times E_1,\Gamma_1,\widehat\phi_1) \to (V_2\times E_2,\Gamma_2,\widehat\phi_2)$
such that they induce an embedding of vector bundles
$(\varphi,\hat\varphi) : (V_1/\Gamma_1,(V_1\times E_1)/\Gamma_1) \to (V_2/\Gamma_2,(V_2\times E_2)/\Gamma_2)$.
The triple
$(h,\tilde\varphi,\tilde{\hat{\varphi}})$ is called an {\it isomorphism}
or a {\it coordinate change} between
orbifold charts of the vector bundle.
\index{orbifold ! isomorphism of orbifold charts of a  vector bundle}
\index{orbifold ! coordinate change of orbifold charts of a  vector bundle}
\end{enumerate}
\end{defn}
\begin{lem}\label{lem2619}
Let $(X^b,\mathcal E^b)$ be a vector bundle over an orbifold $X^b$
and $f : X^a \to X^b$  an embedding of orbifolds.
Let $\mathcal E^a = X^a \times_{X^b} \mathcal E^b$ be the fiber product
in the category of topological space.
By definition of the fiber product, we have maps
$\pi : \mathcal E^a \to X^a$ and $\widehat f : \mathcal E^a \to \mathcal E^b$.
Then the exists a unique structure of vector bundle on $(X^a,\mathcal E^a)$
such that the projection is given the above map $\pi$  and
$(f,\widehat f)$ is an embedding of vector bundles.
\end{lem}
\begin{proof}
Let
$\{\frak V^*_{\frak r} \mid \frak r \in \frak R_*\}$, $*=a,b$
be orbifold atlases where
$\frak V^*_{\frak r} = (V^*_{\frak r},\Gamma^*_{\frak r},\phi^*_{\frak r})$.
Let
$(V^b_{\frak r},E^b_{\frak r},\Gamma^b_{\frak r},\phi^b_{\frak r},\widehat\phi^b_{\frak r})$
be orbifold atlas of the vector bundle $(X^b,\mathcal E^b)$.
Let
$(h_{\frak r,ba},\tilde\varphi_{\frak r,ba})$ be
a local representative of the embedding $f$ on the charts
$\frak V^a_{\frak r}$, $\frak V^b_{\frak r}$.
We put $E^a_{\frak r} = E^b_{\frak r}$, on which $\Gamma^a_{\frak r}$ acts
by the isomorphism $h_{\frak r,ba}$.
By definition of fiber product, there exists uniquely a map
$\widehat\phi^a_{\frak r} : V_{\frak r}^b \times E_{\frak r}^b \to \mathcal E^a$
such that the next diagram commutes.
\begin{equation}\label{diag2619diag}
\begin{CD}
V_{\frak r}^a @ <{\pi}<<
V_{\frak r}^a \times E_{\frak r}^b   @>{\tilde\varphi_{\frak r,ba}\times id}>> V_{\frak r}^b \times E_{\frak r}^b\\
@ V{\phi^a_{\frak r}}VV @ VV{\widehat\phi^a_{\frak r}}V
@VV{\widehat\phi^b_{\frak r}}V\\
X^a @ < {\pi} <<\mathcal E^a @>{\hat f}>>\mathcal E^b
\end{CD}
\end{equation}
In fact,
$$
f\circ \phi^a_{\frak r}\circ \pi
= \phi^b_{\frak r}   \circ \varphi_{\frak r,ba} \circ \pi
= \phi^b_{\frak r} \circ \pi  \circ (\tilde\varphi_{\frak r,ba}\times id)
= \pi \circ \widehat\phi^b_{\frak r} \circ (\tilde\varphi_{\frak r,ba}\times id).
$$
Thus $\{(V^a_{\frak r},E^a_{\frak r},\Gamma^a_{\frak r},\phi^a_{\frak r},\widehat\phi^a_{\frak r}) \mid \frak r \in \frak R\}$
is an atlas of the vector bundle $(X^a,\mathcal E^a)$.
\end{proof}
\begin{defn}
We call the vector bundle in Lemma \ref{lem2619} the {\it pull-back}
\index{pull-back ! of vector bundle on orbifold}
and write $f^*(X^b,\mathcal E^b)$.
(Sometimes we write $f^*\mathcal E^b$ by an abuse of notation.)
\par
In case $X^a$ is an open subset of $X^b$ equipped with
open substructure we call {\it restriction}
\index{orbifold ! restriction of vector bundle}
in place of pull-back of $\mathcal E^b$
and write $\mathcal E^b\vert_{X^a}$ in place of $f^*\mathcal E^b$.
\end{defn}
\begin{lem}\label{lem2622}
In the situation of Lemma \ref{lem26999}
suppose in addition that $\mathcal E^i$ is a vector bundle over $X^i$
and $\widehat\varphi_{21} : \mathcal E^1 \to \mathcal E^2$
is an embedding of vector bundles over $\varphi_{21}$.
Then in addition to the conclusion of Lemma \ref{lem26999},
there exists $\tilde{\hat{\varphi}}_{\frak r;21} :
V^{1}_{\frak r} \times E_{\frak r}^1
\to V^{2}_{\frak r} \times E_{\frak r}^2$
that is an $h_{\frak r;21}$ equivariant embedding of manifolds
with the following properties:
\begin{enumerate}
\item
The next diagram commutes.
\begin{equation}\label{diagin26777}
\begin{CD}
V^1_{\frak r}\times E^1_{\frak r} @ >{\tilde{\hat{\varphi}}_{\frak r,21}}>>
V^2_{\frak r} \times E^2_{\frak r} \\
@ V{\widehat\phi_{\frak r}^{1}}VV @ VV{\widehat\phi_{\frak r}^{2}}V\\
\mathcal E^1 @ > {\widehat\varphi_{21}} >> \mathcal E^2
\end{CD}
\end{equation}
\item
For each $y \in V^1_{\frak r}$ the map
$\xi \mapsto \pi_2(\tilde{\hat{\varphi}}_{\frak r,21}(y,\xi))$
$: E^1_{\frak r} \to E^2_{\frak r}$ is a linear embedding.
\end{enumerate}
\end{lem}
The proof is similar to the proof of Lemma \ref{lem26999}
and is omitted.
\begin{defn}\label{def28262826}
We call $(h_{\frak r,21},\tilde\varphi_{\frak r,21},\tilde{\hat{\varphi}}_{\frak r,21})$
a {\it local representative of embedding}  $(\varphi_{21},\widehat\varphi_{21})$
on the charts
$(V^1\times E^1,\Gamma^1,\widehat\phi^1)$, $ (V^2\times E^2,\Gamma^2,\widehat\phi^2)$.
\index{orbifold ! local representative of embedding of vector bundle on orbifold}
\index{embedding ! local representative of embedding of vector bundle on orbifold}
\end{defn}
\begin{lem}\label{lem2715}
If $(h_{\frak r,21},\tilde\varphi_{\frak r,21},\tilde{\hat{\varphi}}_{\frak r,21})$,
$(h'_{\frak r,21},\tilde\varphi'_{\frak r,21},\tilde{\hat{\varphi}}'_{\frak r,21})$
are local representatives of an embedding of vector bundles of the same
charts
$(V^1\times E^1,\Gamma^1,\widehat\phi^1)$, $ (V^2\times E^2,\Gamma^2,\widehat\phi^2)$,
then there exists $\mu \in \Gamma^2$
such that
$$
\tilde\varphi_{\frak r,21}'(x) = \mu\tilde\varphi_{\frak r,21}(x),
\quad
\tilde{\hat{\varphi}}_{\frak r,21}'(x,\xi) = \mu\tilde{\hat{\varphi}}_{\frak r,21}(x,\xi)
\quad
h_{\frak r,21}'(\gamma) = \mu h_{\frak r,21}(\gamma)\mu^{-1}.
$$
\end{lem}
\begin{proof}
This is a consequence of Lemma \ref{lem21tenhatena}.
\end{proof}
\begin{rem}
In \cite[Situation 6.3]{part11} we introduced the
notation
$(h_{\frak r,21},\tilde\varphi_{\frak r,21},\breve{{\varphi}}_{\frak r,21})$
where $\breve{{\varphi}}_{\frak r,21}$ is related to
$\tilde{\hat{\varphi}}_{\frak r,21}$ by the formula
$$
\tilde{\hat{\varphi}}_{\frak r,21}(y,\xi)
=(\tilde{{\varphi}}_{\frak r,21}(y),\breve{{\varphi}}_{\frak r,21}(y,\xi)).
$$
\end{rem}
We use the pull-back of vector bundles in a different situation.
Let $\mathcal E^i$, $i=1,2$, be vector bundles over an orbifold $X$.
We take the Whitney sum  bundle $\mathcal E^1 \oplus \mathcal E^2$ and
denote by $\vert\mathcal E^1 \oplus \mathcal E^2\vert$ its total
space. There exists a projection
\begin{equation}\label{projfromWhe}
\vert\mathcal E^1 \oplus \mathcal E^2\vert
\to \vert\mathcal E^2\vert.
\end{equation}

\begin{defnlem}\label{pullbackbyproj}
The total space $\vert\mathcal E^1 \oplus \mathcal E^2\vert$
has a structure of vector bundle over  $\vert\mathcal E^2\vert$
such that (\ref{projfromWhe}) is the projection.
We write it as
$\pi_{\mathcal E^2}^*\mathcal E^1$
and call the {\it pull-back} of $\mathcal E^1$ by the projection
$\pi_{\mathcal E^2} : \vert\mathcal E^2\vert \to X$.
\index{pull-back ! of vector bundle on orbifold}
\par
When $U$ is an open subset of $\vert\mathcal E^2\vert$ and
$\pi : U \to X$ is the restriction of $\pi_{\mathcal E^2}$ to $U$,
the pull-back $\pi^*_{\mathcal E^2}\mathcal E^1$ is by definition
the restriction of $\pi_{\mathcal E^2}^*\mathcal E^1$ to $U$.
\end{defnlem}
The proof is immediate from definition.
\begin{rem}
We note that the total space $\vert\mathcal E^1 \oplus \mathcal E^2\vert$
is {\it not} a fiber product
$\vert\mathcal E^1 \vert \times_X \vert\mathcal E^2\vert$.
In fact, if $X$ is a point and $\mathcal E^1
= \mathcal E^2= \R^n/\Gamma$ with linear $\Gamma$ action,
then the fiber of $\vert\pi_{\mathcal E^2}^*\mathcal E^1 \vert \to
\vert\mathcal E^2\vert \to X$ at $[0]$ is
$(E^1 \times E^2)/\Gamma$.
The fiber of the map
$\vert\mathcal E^1 \vert \times_X \vert\mathcal E^2\vert \to X$
at $[0]$ is $(E^1/\Gamma) \times (E^2/\Gamma)$.
\end{rem}
\begin{defn}
Let $(X,\mathcal E)$ be a vector bundle. A {\it section}
of $(X,\mathcal E)$ is an embedding of orbifolds $s : X \to \mathcal E$
such that the composition of $s$ and the projection
is the identity map set-theoretically.
\index{orbifold ! section of vector bundle}
\end{defn}

\begin{lem}\label{lem2626}
Let $\{(V_{\frak r},E_{\frak r},\Gamma_{\frak r},
\psi_{\frak r},\widehat\psi_{\frak r})\mid \frak r \in \frak R\}$
be an atlas of $(X,\mathcal E)$.
Then a section of  $(X,\mathcal E)$
corresponds one to one to the following object.
\begin{enumerate}
\item
For each $\frak r$ we have a $\Gamma_r$ equivariant
map $s_{\frak r} : V_{\frak r} \to E_{\frak r}$,
which is compatible in the sense of (2) below.
\item
Suppose $\phi_{\frak r_1}(x_1) = \phi_{\frak r_2}(x_2)$.
Then the definition of orbifold atlas implies that
there exist subcharts
$(V_{\frak r_i,x_i},E_{\frak r_i,x_i},\Gamma_{\frak r_i,x_i},
\phi_{\frak r_i,x_i},\widehat\phi_{\frak r_i})$
of the orbifold charts
$(V_{\frak r_i},E_{\frak r_i},\Gamma_{\frak r_i},
\phi_{\frak r_i},\widehat\phi_{\frak r_i})$
at $x_i \in V_{\frak r_i}$ for $i = 1,2$ and
an isomorphism of charts
$$
\aligned
(h^{\frak r,p}_{12},\tilde\varphi^{\frak r,p}_{12},
\tilde{\hat{\varphi}}^{\frak r,p}_{12}) ~:~
&(V_{\frak r_2,x_2},E_{\frak r_2},\Gamma_{\frak r_2,x_2},
\phi_{\frak r_2,x_2},\widehat\phi_{\frak r_2})
\\
&\to (V_{\frak r_1,x_1},E_{\frak r_1,x_1},\Gamma_{\frak r_1,x_1},
\phi_{\frak r_1,x_1},\widehat\phi_{\frak r_1}).
\endaligned
$$
\par
Now we require the following equality:
\begin{equation}\label{sectioncompati}
\tilde{\hat{\varphi}}^{\frak r,p}_{12}(s_{\frak r_1}(y,\xi)) =
s_{\frak r_2}(\tilde\varphi^{\frak r,p}_{12}(y),\xi).
\end{equation}
\end{enumerate}
\end{lem}
\begin{proof}
The proof is mostly the same as the corresponding
standard result
for the case of vector bundle on a manifold or on a topological space.
Let $s : X \to \mathcal E$ be a section, which is an orbifold embedding.
Let $p \in \phi_{\frak r}(V_{\frak r})$.
Then there exist a subchart $(V_{\frak r,p},\Gamma_{\frak r,p},\phi_{\frak r,p})$
of $\frak V_{\frak r}$ and
a subchart $(\widehat V_{\frak r,\tilde p},
\Gamma_{\frak r,\tilde p},\phi_{\frak r,\tilde p})$
of $(V_{\frak r}\times E_{\frak r},\Gamma_{\frak r},\phi_{\frak r,\tilde p})$
such that a representative  $(h',\tilde\varphi')$
of $s$ exists on the subcharts.
Since $\pi \circ s = $identity set-theoretically,
it follows that $\pi_1(\tilde\varphi(y)) \equiv y \mod \Gamma_p$
for any $y \in V_{\frak r,p}$.
We take $y$ such that $\Gamma_{y} = \{1\}$.
Then, there exists a {\it unique} $\mu \in \Gamma_p$
such that $\pi_1(\tilde\varphi'(y)) \equiv \mu y$.
By continuity this $\mu$ is independent of $y$.
(We use Condition \ref{convinv} here.)
\par
We replace $\tilde p$ by $\mu^{-1}\tilde p$
and $(\widehat V_{\frak r,\tilde p},\Gamma_{\frak r,\tilde p},\phi_{\frak r,\tilde p})$
by $(\mu^{-1}\widehat V_{\frak r,\tilde p},\mu^{-1}\Gamma_{\tau,\tilde p} \mu,\phi_{\frak r,\tilde p}\circ \mu)$
and $(h',\tilde\varphi')$ by $(h'\circ {\rm conj}_{\mu},\tilde\varphi' \circ \mu^{-1})$.
(Here ${\rm conj}_{\mu}(\gamma) = \mu \gamma \mu^{-1}$.)
Therefore we may assume $\pi_1(\tilde\varphi'(y)) = y$.
Note that $\tilde\varphi'$ is $h'$-equivariant and $\pi_1$ is ${\rm id}$-equivariant.
Here ${\rm id}$ is the identity map $\Gamma_{\frak r,y} \to \Gamma_{\frak r,y}$.
Therefore the identity map $V_{\frak r,p} \to V_{\frak r,p}$ is
$h'$ equivariant.
Hence $h' = {\rm id}$.
\par
In sum, we have the following.
(We put $s_{\frak r,p} = \tilde\varphi'$.)
For a sufficiently small $\frak V_{\frak r,p}$
there exists uniquely a map $
s_{\frak r,p} : V_{\frak r,p} \to V_{\frak r,p} \times E_{\frak r}$
such that
\begin{enumerate}
\item[(a)]
$
\pi_1(s_{\frak r,p}(x)) =x
$
\item[(b)]
$s_{\frak r,p}$ is equivariant with respect to the embedding
$\Gamma_{\frak r,p} \to \Gamma_{\frak r}$.
(Recall $\Gamma_{\frak r,p} =
\{\gamma \in \Gamma_{\frak r} \mid \gamma p = p\}$.)
\item[(c)]
$({\rm id},s_{\frak r,p})$ is a local representative of $s$.
\end{enumerate}
We can use uniqueness of such $s_{\frak r,p}$ to glue them
to obtain a map $V_{\frak r} \to V_{\frak r} \times E_{\frak r}$.
By (a)  this map is of the form
$x \mapsto (x,\frak s_{\frak r}(x))$
for some map $\frak s_{\frak r} : V_{\frak r} \to E_{\frak r}$.
This is the map $\frak s_{\frak r}$ required in (1).
Since $x \mapsto \gamma^{-1}\frak s_{\frak r}(\gamma x)$
also has the same property, the uniqueness implies that
$\frak s_{\frak r}$ is $\Gamma_{\frak r}$
equivariant.
(\ref{sectioncompati}) is also a consequence of
the uniqueness.
\par
Thus we find a map from the set of sections to the
set of $(s_{\frak r})_{\frak r \in \frak R}$ satisfying (1)(2).
The construction of the converse map is obvious.
\end{proof}
The next lemma is proved during the proof of Lemma \ref{lem2626}.

\begin{lem}\label{lem2627}
Let $(V_{\frak r},E_{\frak r},\Gamma_{\frak r},\phi_{\frak r},\widehat\phi_{\frak r})$ be an orbifold chart of
$(X,\mathcal E)$ and $s$ a section of $(X,\mathcal E)$.
Then there exists uniquely a $\Gamma$ equivariant map
$s_{\frak r} : V_{\frak r} \to E_{\frak r}$ such that the following diagram commutes.
\begin{equation}\label{diagin26277XXrev3}
\begin{CD}
V_{\frak r}\times E_{\frak r} @ >{\widehat\phi_{\frak r}}>>
\mathcal E_{\frak r} \\
@ A{{\rm id} \times s_{\frak r}}AA @ AA{s}A\\
V_{\frak r} @ > {\phi_{\frak r}} >> X
\end{CD}
\end{equation}
\end{lem}
\begin{defn}\label{defnlocex}
We call the system of maps $s_{\frak r}$ the {\it local expression} of $s$ in the orbifold
chart $(V_{\frak r},E_{\frak r},\Gamma_{\frak r},\phi_{\frak r},\widehat\phi_{\frak r})$.
\index{orbifold ! local expression of a section on orbifold chart}
\end{defn}
Next we review the proofs of a few well-known facts on
pull-back bundle etc.. Those proofs are straightforward generalization of the
corresponding proofs for the pull-back in the manifold theory. We include them
only for completeness' sake.
\begin{prop}\label{homotopicpulback}
Let $\mathcal E$ be a vector bundle on $X \times [0,1]$, where
$X$ is an orbifold. We identify $X \times \{0\}$, $X \times \{1\}$
with $X$ in an obvious way. Then
there exists an isomorphism of vector bundles
$$
I : \mathcal E\vert_{X \times \{0\}} \cong \mathcal E\vert_{X \times \{1\}}.
$$
Suppose in addition that
we are given a compact set $K \subset X$,
its neighborhood $V$ and an isomorphism
$$
I_0 : \mathcal E\vert_{V \times [0,1]} \cong \mathcal E\vert_{V \times \{0\}} \times [0,1].
$$
Then we may choose $I$ so that it coincides with the isomorphism induced by
$I_0$ on $K$.
If $K$ is a submanifold, we may take $K=V$.
\end{prop}
For the proof of the proposition we use the notion of connection on vector
bundle on orbifolds, which we now recall here.
Note that a vector field on an orbifold is a section of the tangent bundle.
\begin{defn}
A {\it connection}\index{connection} on a vector bundle $(X,\mathcal E)$ is
an $\R$
linear map
$$
\nabla : C^{\infty}(TX) \otimes_{\R} C^{\infty}(\mathcal E)
\to C^{\infty}(\mathcal E)
$$
such that $\nabla_X(V) = \nabla(X,V)$ satisfies
$$
\nabla_{fX}(V) = f\nabla_X(V), \qquad
\nabla_{X}(fV) = f\nabla_X(V) + X(f)V.
$$
\end{defn}
Here $C^{\infty}(\mathcal E)$ is the vector space consisting of
all smooth sections of $\mathcal E$.
\par
For any connection $\nabla$ and piecewise smooth map $
\ell : [a,b] \to X$ we obtain parallel transport
$$
{\rm Pal}^{\nabla} : \mathcal E_{\ell(a)}  \to \mathcal E_{\ell(b)}
$$
in the same way as in the case of manifolds.
\begin{rem}
Here $\mathcal E_{\ell(a)}$ is the fiber of $\mathcal E$ at $\ell(a) \in  X$
and is defined as follows.
We take a chart $(V_{\frak r},E_{\frak r},\Gamma_{\frak r},\psi_{\frak r},\widehat\psi_{\frak r})$ of $(\mathcal E,X)$
at $\ell(a)$. Then $\mathcal E_{\ell(a)} = E_{\frak r}$.
If $(V_{\frak r'},E_{\frak r'},\Gamma_{\frak r'},\psi_{\frak r'},\widehat\psi_{\frak r'})$ is another chart,
we can identify $E_{\frak r}$ with $E_{\frak r'}$ by
$\xi \mapsto \breve{\varphi}_{\frak r'\frak r}(\xi,y)$ where $\psi_{\frak r}(y) = \ell(a)$ and
$ \breve{\varphi}_{\frak r'\frak r} : V_{\frak r} \times E_{\frak r} \to E_{\frak r'}$ is a part of the
coordinate change.
(\cite[Situaion 6.3]{part11}.)
Note that the identification $\xi \mapsto \breve{\varphi}_{\frak r'\frak r}(\xi,y)$ is well-defined up to the $\Gamma_{\ell(a)}
= \{ \gamma \in \Gamma_{\frak r} \mid \gamma(y) = y\}$ action.
Thus
the parallel transport
$
{\rm Pal}^{\nabla} : \mathcal E_{\ell(a)}  \to \mathcal E_{\ell(b)}
$ is well-defined up to the $\Gamma_{\ell(a)} \times \Gamma_{\ell(b)}$ action.
\end{rem}
\begin{lem}\label{lem2328}
Any vector bundle $\mathcal E$ over orbifold $X$ has a connection.
Moreover if a connection is given for $\mathcal E\vert_{V}$
where $V$ is an open neighborhood of a compact subset $K$ of $X$,
then we can extend it without changing it on $K$.
If $K$ is a submanifold, we may take $K=V$.
\end{lem}
The proof is an obvious analogue of the proof of the existence of
connection on a vector bundle over a manifold, which uses a partition of unity.
\par
We are now ready to give the proof of Proposition \ref{homotopicpulback}.
\begin{proof}[Proof of Proposition \ref{homotopicpulback}]
We take a connection on $\mathcal E\vert_V$.
We then take direct product connection on
$\mathcal E\vert_{V \times \{0\}} \times [0,1]$,
and use $I_0$ to obtain a connection on
$ \mathcal E\vert_{V \times [0,1]}$.
We extend it to a connection on $\mathcal E$
without changing it on $K \times [0,1]$.
For each fixed $x \in X$, we can use parallel transportation
along the path $t \mapsto (x,t)$ to get an isomorphism
$\mathcal E_{(x,0)} \cong \mathcal E_{(x,1)}$.
We have thus obtained a set theoretical map
$$
\vert\mathcal E\vert_{X \times \{0\}}\vert \cong \vert\mathcal E\vert_{X \times \{1\}}\vert.
$$
It is easy to see that it induces an isomorphism of vector bundles.
Using the fact that our connection is direct product on $K \times [0,1]$,
we can check the second half of the statement.
\end{proof}
\begin{defn}
We say two embeddings of orbifold $f_i : X \to Y$
($i=1,2$) are {\it isotopic} to
\index{isotopic} each other
if there exists an embedding of orbifolds
$H : X\times [0,1] \to Y \times [0,1]$ such that
the second factor of $H(x,t)$ is $t$ and that
$$
H(x,0) = (f_1(x),0) \qquad
H(x,1) = (f_2(x),1).
$$
Suppose $V \subset X$ and $f_1 = f_2$ on a neighborhood $V$ of $K$. We say
$f_1$ is {\it isotopic to $f_2$ relative to $K$}
\index{isotopic relative to a compact subset} if we may take $H$ such that
\begin{equation}\label{homoisidentity}
H(x,t) = (f_1(x),t) = (f_2(x),t)
\end{equation}
for $x$ in a neighborhood of $K$.
In case $K$ is a submanifold, we may take $K=V$
and then (\ref{homoisidentity}) holds for $x \in K$.
\end{defn}
\begin{cor}\label{cor2939}
Let $f_i : X \to Y$  be two isotopic embeddings and
$\mathcal E$ a vector bundle on $Y$. Then
the pull-back bundle $f_1^*\mathcal E$ is isomorphic
to $f_2^*\mathcal E$.
If $f_1 = f_2$ on a neighborhood of $K \subset X$ and
$f_1$ is isotopic to $f_2$ relative to $K$,
then we may choose the isomorphism $f_1^*\mathcal E \cong f_2^*\mathcal E$
so that its restriction to $K$ is the identity map.
\end{cor}
\begin{proof}
This is immediate from
Proposition \ref{homotopicpulback} and the definition.
\end{proof}
We next recall \cite[Definition 12.23]{part11} which we re-state here.
\begin{defn}\label{lem1230002}
Let $f : X \to Y$ be an embedding of orbifolds and $K\subset X$  a
compact subset and $U$ an open neighborhood of $K$ in
$Y$.\footnote{Here we regard the image $f(K)$ as a subset of $Y$ via the embedding $f$, and
write $f(K)$ as $K$ to simply the notation.}
We say that
a continuous map $\pi : U \to X$ is
{\it diffeomorphic to the projection of the normal bundle}
\index{diffeomorphic to the projection of the normal bundle} if the following holds.
\par
Let ${\rm pr} : N_XY \to X$ be the normal bundle. Then there exists
a neighborhood $U'$  of $K$ in $ N_XY$, (Note $K\subset X \subset N_XY$.) and a diffeomorphism
of orbifolds $h : U' \to U$ such that $\pi\circ h = {\rm pr}$.
We also require that $h(x) = x$ for $x$ in a neighborhood of $K$ in  $X$.
\end{defn}
\begin{defnlem}\label{defpullbackbundlenbd}
Let $\pi : U \to X$ be diffeomorphic to the
projection of the normal bundle as in Definition \ref{lem1230002}
and $\mathcal E$ a vector bundle on $X$. We define
$\pi^*\mathcal E$, the pull-back bundle as follows.
\par
Let $h$, $U'$ be as in Definition \ref{lem1230002}.
We defined a pull-back bundle ${\rm pr}^*\mathcal E$ on $N_XY$
in Definition \ref{pullbackbyproj}.
We put
$$
\pi^*\mathcal E = (h^{-1})^*{\rm pr}^*\mathcal E\vert_{U'}.
$$
This is independent of the choice of $(U',h)$ in the following sense.
Let $U'_i$, $h_i$ ($i=1,2$) be two choices.
Then we can shrink $U$ and $U'_i$ so that for each $i$ the restriction of $h_i$
becomes an isomorphism between them.
Then
\begin{equation}\label{294949}
(h_1^{-1})^*{\rm pr}^*\mathcal E\vert_{U'_1}
\cong
(h_2^{-1})^*{\rm pr}^*\mathcal E\vert_{U'_2}.
\end{equation}
Moreover the isomorphism (\ref{294949})
can be taken so that the following holds in addition.
We regard $K \subset U$.
Then by definition it is easy to see that
the restriction of both sides of (\ref{294949})
is canonically identified with
the restriction of $\mathcal E$ to $K \subset X$.
The isomorphism  (\ref{294949}) becomes the identity map
on $K$ under this isomorphism.
\end{defnlem}
\begin{proof}
We can replace $U$ by a smaller open neighborhood so that
$h_1^{-1} : U \to N_XY$ is isotopic to $h_2^{-1} : U \to N_XY$.
(See the proof of Proposition \ref{prop2942} below.)
Then (\ref{294949}) follows from Corollary \ref{cor2939}.
The second half of the claim also follows from the second half of
Corollary \ref{cor2939}.
\end{proof}
The pull-back bundle is independent of the projection $\pi$ but depends only
on the neighborhood $U$ in the situation of Definition \ref{lem1230002}.
In fact, we have
\begin{prop}\label{prop2942}
Let $\pi_i : U \to X$ be as in Definition \ref{lem1230002} for $i=1,2$.
Then there exist a neighborhood $U_0$ of $X$ in $Y$ and a map
$f : U_0 \to U$ such that
\begin{enumerate}
\item
$
\pi_2 \circ f = \pi_1
$
\item
$f : U_0 \to U$ is isotopic to the inclusion map $U_0 \hookrightarrow U$
relative to $X$.
\end{enumerate}
\end{prop}
\begin{proof}
Let $h_i : U'_i \to U_i$ be as in  Definition \ref{lem1230002}.
We put $f = h_2 \circ h^{-1}_1$ which is defined for sufficiently
small  $U_0$. If suffices to show that $f$ is isotopic to
the inclusion map.
We first prove the proposition in the case when the following additional assumption
is satisfied. (We will  remove this assumption later.)
\begin{assump}\label{assym2929}
For any $x \in K \subset N_XY$ the first derivative at $x$, $D_xf : T_x(N_XY) \to T_x(N_XY)$ is the
identity map.
\end{assump}
\par
For the case of manifolds, this assumption enables us to prove
Proposition \ref{prop2942} by observing that $f$ is $C^1$-close to the inclusion map.
Then, for example, using minimal geodesic, we can show that $f$ is isotopic to the inclusion map.
\par
For the case of orbifolds, we need to work out this last step a bit more
carefully since the number
$$
\inf \{ r \mid \text{if $d(x,y) < r$, the minimal geodesic joining $x$ and $y$ is unique}\}
$$
can be $0$ in general unlike the case of manifolds.
\par
To clarify this point,
we need to prepare certain lemmas whose statements require some digression.
We can define the notion of Riemannian metric of orbifold $X$ in a straightforward way
similarly as in the manifold case.
For $p\in X$ we have a geodesic coordinate
$(TB_p(c_p),\Gamma_p,\psi_p)$ where
$$
TB_p(c_p) = \{ \xi \in T_pX \mid \Vert \xi\Vert < c_p\}
$$
and the group $\Gamma_p$ is the isotropy group of the orbifold chart of $X$ at $p$.
The uniformization map $\psi : TB_p(c_p) \to X$ is defined by using
minimal geodesic in the same way as the construction of the exponential map
in the standard Riemannian geometry.
We note that this map is well defined up to the action of $\Gamma_p$.
We need to take the number $c_p$ small
so that the exponential map $\psi$ induces a homeomorphism
$TB_p(c_p)/\Gamma_p \to X$. We may not be able to choose $c_p$ uniformly away from $0$
even on the compact subset of the given orbifold in general.
(This is because $d(p,q) < c_p$ must imply
$\#\Gamma_q \le \# \Gamma_p$, when $c_p$ is sufficiently small.)
However we can prove the following.
Let $X$ be an orbifold and $Z$ a compact set. Suppose
$B_{c_0}(Z) = \{x \mid d(x,Z) \le c_0\}$ is complete with respect to the metric induced by
the Riemannian metric.
\begin{lem}
Let $Z \subset X$ be a compact subset.
Then there exists a finite set $\{p_i \mid j \in J\} \subset Z$ and $0 < c_j < c_0$
such that
\begin{enumerate}
\item The geodesic coordinate $(TB_{p_j}(c_j),\Gamma_{p_j},\psi_{p_j})$
exists.
\item
$$
Z \subset \bigcup_{j} \psi_{p_j}(TB_{p_j}(c_j/2)).
$$
\end{enumerate}
\end{lem}
The proof is immediate from the compactness of $Z$.
We call such $\{(TB_{p_j}(c_j),\Gamma_{p_j},\psi_{p_j}) \mid j\}$
a {\it geodesic coordinate system}\index{geodesic coordinate system} of $(X,Z)$.
We put $P = \{p_j \mid j =1,\dots, J\}$.
\begin{defn}\label{lem2945}
We fix a geodesic coordinate system of $(X,Z)$.
Let $Z_0 \subset Z$ be a compact subset containing $P$ and
$F : U \to X$ be an
embedding of orbifolds where $U \supset Z$ is an open neighborhood of $Z$.
We say $F$ is {\it $C^1$ $\epsilon$-close to the identity}
\index{$C^1$ $\epsilon$-close to the identity} on $Z_0$ if the following holds.
\begin{enumerate}
\item
$F(B_{p_j}(c_j/2)) \subset B_{p_j}(c_j)$.
\item
There exists
$\tilde F_{j} : B_{p_j}(c_j/2) \to B_{p_j}(c_j)$
such that:
\begin{enumerate}
\item
$
\psi_{p_j} \circ  \tilde F_{j} = F \circ \psi_{p_j}$.
\item
$d(x,\tilde F_{j}(x)) < \epsilon$ for $x \in TB_{p_j}(c_j/2) \cap \psi_{p_j}^{-1}(Z_0)$.
\item
$d(D_x\tilde F_{j},id) < \epsilon$ for $x \in TB_{p_j}(c_j/2) \cap \psi_{p_j}^{-1}(Z_0)$.
\end{enumerate}
\smallskip
Here $d$ in Item (b) is the standard metric on Euclidean space $T_{p_j}X$
(together with the metric induced by our Riemannian metric),
$d$ in Item (c) is a distance in the space of $n\times n$ matrices.
(Here $n = \dim X$. We use our Riemannian metric to define a metric
on this space of matrices, which is a vector space of dimension $n^2$ with metric.)
\end{enumerate}
\end{defn}
\begin{lem}\label{lem2946}
For each $Z$ and a geodesic coordinate system of $(X,Z)$,
there exists $\epsilon >0$ such that the following
holds for any given $Z_0 \subset Z$ and $F : U \to X$:
\par
If $F$ is $C^1$ $\epsilon$-close to the identity on $Z_0$,
then $F$ is isotopic to the identity on $Z_0$.
Moreover for any $\delta >0$ there exists $\epsilon(\delta) >0$
such that if $F$ is $C^1$ $\epsilon(\delta)$-close to the identity on $Z_0$,
then the isotopy from $F$ to the identity map is taken
to be  $C^1$ $\delta$-close to the identity on $Z_0$.
\end{lem}
\begin{proof}
We first observe that if $\epsilon >0$ is sufficiently small, then
for each $j$ the map
$\tilde F_j$ satisfying Definition \ref{lem2945} (2) (a),(b) and (c) is unique.
Indeed, it easily follows that
such $\tilde F_j$ is unique up to the action of
$\Gamma_{p_j}$.
Since the $\Gamma_{p_j}$ action is effective,
$\Gamma_{p_j}$ is a finite group and $p_j \in Z_0$,
we find that at most one such $\tilde F_j$ can satisfy (c).
(Note that
the map $\Gamma_{p_j} \to O(n)$ taking the linear part at $p_j$ of the action
is injective since the $\Gamma_{p_j}$ action has a fixed point and is effective.)
\par
Next for each $t \in [0,1]$, we define a map
$$
\tilde F_{t,j} : V_j \to TB_{p_j}(c_j/2)
$$
as follows. Here $V_j$ is a sufficiently small neighborhood of
$TB_{p_j}(c_j/2) \cap \psi_j^{-1}(Z_0)$.
We take a Riemannian metric on $V_j$ that is the pull-back
of the given Riemannian metric on $X$ by the map $\psi_{p_j}$.
By choosing $\epsilon$ sufficiently small and using (b),
we find
a unique minimal geodesic $\ell_{x,j} : [0,1] \to TB_{p_j}(c_j/2+2\epsilon)$ joining
$x$ to $\tilde F_{j}(x)$. We put
$$
\tilde F_{t,j}(x) = \ell_{x,j}(t).
$$
In the same way as in the proof of the uniqueness of $\tilde F_j$ we can show that there exists
a map $F_t$ such that
$
\psi_{p_j} \circ  \tilde F_{t,j} = F_t \circ \psi_{p_j}$.
Using Definition \ref{lem2945} (2)(b) and (c), we can show that $F_t$ is
$C^1$ $\epsilon$-close to the identity map.
By choosing $\epsilon >0$ sufficiently small,
we derive that $F_t$ is a diffeomorphism to
its image.
Thus
$F_t$ is the required isotopy from $F$ to the identity map.
\end{proof}
Now we use Lemma \ref{lem2946} to complete the proof of Proposition \ref{prop2942}
under Assumption \ref{assym2929}.
Recall from the beginning of the proof of Proposition \ref{prop2942} that
we put
$f = h_2 \circ h^{-1}_1$ which we want to prove is isotopic to
the identity map in a neighborhood of $K \subset X \subset Y$.
\par
We take a finite set of points $p_j \in K$ so that
$
K \subset \bigcup_{j} \psi_{p_j}(TB_{p_j}(c_j/2))
$.
Then we take a compact neighborhood $Z \supset K$ such that
$Z \subset \bigcup_{j} \psi_{p_j}(TB_{p_j}(c_j/2))$.
We apply Lemma \ref{lem2946} to obtain $\epsilon$.
Note that $f$ is the identity map on $K$ and its first derivative
at any point in $K$ is also the identity by Assumption \ref{assym2929}.
Therefore we can find a sufficiently small compact neighborhood $Z_0$ of $K$
so that $f$ is $C^1$ $\epsilon$-close to the identity on $Z_0$.
Thus Lemma \ref{lem2946} implies
that $f$ is isotopic to the identity map.
The proof of Proposition \ref{prop2942} under the additional Assumption \ref{assym2929}
is now complete.
\par
To remove Assumption \ref{assym2929},
we use the following lemma.
\begin{lem}\label{lem292333}
Let $U$ be an open neighborhood of $K$ in $N_XY$
and $F : U \to N_XY$ be an open embedding of orbifolds.
Assume $F ={\rm id}$ on a neighborhood of $K$ in $X$
and $D_xF(V) \equiv V \mod T_xX$ at any $x$ contained in a neighborhood of $K$.
\par
Then there exists a smaller neighborhood $U'$ of $K$ such that
the restriction of $F$ to $U'$ is isotopic to the embedding satisfying Assumption \ref{assym2929}.
\end{lem}
\begin{proof}\label{lem2947}
We take the first derivative of $F$
$$
D_x F : T_x N_XY \to N_XY
$$
at $x$ contained in a neighborhood of $K$ in $X$.
Note that $T_xX \subset N_XY$ is preserved under this map
and we have a decomposition $T_x N_XY = T_xX \oplus (N_XY)_x$.
Therefore there exists a linear bundle map
$$
H : N_XY \to TX
$$
on a neighborhood of $K$ such that
$$
D_xF (\xi,\eta) = (\xi + H_x(\eta), \eta),
$$
where $\xi \in  T_xX$ and $\eta \in (N_XY)_x$.
Now we define $G_t : U' \to N_XY$ as follows.
(Here $U'$ is a small neighborhood of $K$ in $N_XY$.)
Let $(x,\eta) \in U'$ where $x$ is in a neighborhood of $K$ in $X$ and $\eta \in (N_XY)_x$.
We take a geodesic $\ell : [0,1] \to X$ of constant speed with $\ell(0) = x$ and
$D\ell/dt(0) = H_x(\eta)$.
Let $\ell_{\le t}$ be its restriction to $[0,t]$.
Then $G_t(x,\eta) = (\ell(t),{\rm Pal}_{\ell_{\le t}}(\eta))$,
where ${\rm Pal}_{\ell_{\le t}}(\eta) \in (N_XY)_{\ell(t)}$ is the
parallel transport along $\ell_{\le t}$.
By construction the first derivative of $G_t$ at a point in $K$ is given by
$$(\xi,\eta)
\mapsto
(\xi + tH_x(\eta), \eta),
$$
which is invertible.
Therefore,
if $V'$ is a sufficiently small neighborhood of $K$, then the restriction of $G_t$ is an embedding $V' \to N_XY$.
Note that $F\circ G_1^{-1}$ satisfies Assumption \ref{assym2929}.
Thus the proof of Lemma \ref{lem292333} is complete.
\end{proof}
Using Lemma \ref{lem292333}, we can reduce the general case of
Proposition \ref{prop2942} to the case when Assumption \ref{assym2929}
is satisfied.
The proof of Proposition \ref{prop2942} is complete.
\end{proof}

We used the next result in \cite[Subsection 13.2]{part11}.
\begin{prop}\label{prop2949}
Let $f : X \to Y$ be an embedding of orbifolds and $K_i$  compact subsets of $X$
for $i=1,2$ such that $K_1 \subset {\rm Int}K_2$.
\footnote{Hereafter we simply write $f(X), f(K_i)$ as $X, K_i$ and regard them as subsets of $Y$
via the embedding $f$, when no confusion can occur.}
Suppose $U_1$ is an open neighborhood $K_1$ in $Y$ and $\pi_1 : U_1 \to X$
is diffeomorphic to the projection of the normal bundle
in the sense of Definition \ref{lem1230002}
\par
Then there exists an open neighborhood $U_2$ of $K_2$ in $Y$ and $\pi_2 : U_2 \to X$
such that it is diffeomorphic to the projection of the normal bundle
and $\pi_1 = \pi_2$ on an open neighborhood of $K_1$.
\end{prop}
\begin{proof}
For the case of manifolds, the corresponding result is standard
and can be proved by applying the isotopy extension lemma.
We can substitute the isotopy extension lemma by
Lemma \ref{lem2946} in the same way to handle the orbifold case.
For completeness' sake we give detail of the proof below.
\par
We first apply \cite[Lemma 6.5]{fooooverZ} to obtain an open
neighborhood $U'_2$ of $K_2$ in $Y$ and a submersion $\pi'_2 : U'_2 \to X$ such that
$(U'_2,\pi'_2)$ is diffeomorphic to the normal bundle $N_XY$.
We may adjust the map $\pi'_2$ to $\pi_2$
so that $\pi_1 = \pi_2$ on an open neighborhood of $K_1$ whose detail is in order.
\par
Let $W_1^{(i)}$ be a neighborhood of $K_1$ in $X$ such that
$$
\overline{W_1^{(1)}} \subset W_1^{(2)} \subset \overline{W_1^{(2)}}
\subset U_1 \cap X.
$$
Let $\Omega$ be an open subset of $U_1$ with
$$
\overline{W_1^{(2)}} \subset \Omega
\subset \overline\Omega \subset U_1.
$$
Later on we will choose $\Omega$ so that it is contained in a sufficiently small
neighborhood of $X$ in $Y$.
We put
$$
V_1^{(i)} = \pi_1^{-1}(W_1^{(i)}) \cap  \Omega.
$$
Let $\tilde\chi : {\rm Int}(K_2) \cup \Omega \to [0,1]$ be a smooth function such that
$$
\tilde\chi =
\begin{cases}
1  &\text{on $W_1^{(1)}$} \\
0  &\text{on the complement of $W_1^{(2)}$},
\end{cases}
$$
and put $\chi = \tilde\chi \circ \pi'_2$.
We may choose $\Omega$ so small that the following holds.
$$
\chi =
\begin{cases}
1  &\text{on $V_1^{(1)}$} \\
0  &\text{on the complement of $V_1^{(2)}$}.
\end{cases}
$$
\par
Let $Z =  \overline{W_1^{(2)}} \setminus  {W_1^{(1)}}$.
We take a neighborhood $U'$ of $Z$ and restrict
$\pi_1$ and $\pi'_2$ there.
Then we can apply  Proposition \ref{prop2942}
to prove that there exists an isotopy $F_t : U' \to X$ such that
$\pi'_2 \circ F_1 = \pi_1$ and $F_0$ is the identity map.
Now we put
$$
\pi_2(x)
=
\begin{cases}
(\pi_2 \circ F_{\chi(x)})(x) &\text{on $U'$} \\
\pi_1 &\text{ on $V_1^{(1)}$} \\
\pi'_2(x) &\text{elsewhere on $U'_2$}.
\end{cases}
$$
It is easy to see that they are glued to define a map.
To complete the proof it suffices to show that
$x \mapsto F_{\chi(x)}(x)$ is an embedding $: U' \to X$.
We will prove this below.
\par
We first consider the case that  Assumption \ref{assym2929} is satisfied for the map $f : U' \to X$ with
$\pi'_2 \circ f = \pi_1$.
In this case we may choose the isotopy $F_t$ arbitrarily close to the identity
map in $C^1$ sense by taking $\Omega$ sufficiently small.
(This is the consequence of the second half of Lemma \ref{lem2946}.)
Therefore the first derivative of $x \mapsto F_{\chi(x)}(x)$ is close to
the identity. It follows that this map is an embedding.
\par
We finally show that we can choose $(U'_2,\pi'_2)$ so that
Assumption \ref{assym2929} is satisfied for the map $f : U' \to X$ with
$\pi'_2 \circ f = \pi_1$.
We consider the fiber $\pi_1^{-1}(x)$ of $\pi_1$.
We may choose a Riemannian metric of $X$ in a neighborhood of $K_1$
so that this fiber is perpendicular to $X$ at any $x$ in a neighborhood of $K_1$.
We now extend this Riemannian metric to the whole $X$.
We use this Riemannian metric and the associated exponential map in the
normal direction to identify a neighborhood of $K_2$ with its
normal bundle and obtain $U'_2$ and $\pi'_2$.  Then Assumption \ref{assym2929} is satisfied.
The proof of Proposition \ref{prop2949} is now complete.
\end{proof}
\section{Covering space of effective orbifold and K-space}
\label{sec:cover}

\subsection{Covering space of orbifold}
\label{subsec:cover}

We first define the notion of a covering space of an orbifold.
Let $U_1, U_2$ be orbifolds and let $\pi : U_1 \to U_2$ be a
continuous map between their underlying topological spaces.

\begin{defn}
For $i=1,2$ let $x_i \in U_i$ with $\pi(x_1) = x_2$ and
$\frak V_i = (V_i,\Gamma_i,\phi_i)$ be orbifold charts of $U_i$ at $x_i$.
We say that
$(\frak V_1,\frak V_2)$ is a {\it covering chart}
\index{orbifold ! covering chart}
if the following holds:
\begin{enumerate}
\item
There exists an injective group homomorphism
$h_{21} : \Gamma_1 \to \Gamma_2$.
\item
There exists an $h_{21}$-equivariant diffeomorphism $\varphi_{21}
: V_1 \to V_2$.
\item
$\phi_2 \circ \varphi_{21} = \phi_1$.
\end{enumerate}
The index $[\Gamma_2:h_{21}(\Gamma_1)]$
is called the {\it covering index}
of the covering chart $(\frak V_1,\frak V_2)$.
\end{defn}

\begin{defn}\label{coveringdef}
The map $\pi : U_1 \to U_2$ is called a
{\it covering map}
\index{orbifold ! covering map}
if the following holds at each $x \in U_2$.
\begin{enumerate}
\item
The set $\pi^{-1}(x)$
is a finite set, which we write
$
\{\tilde x_1,\dots,\tilde x_{m_x}\}.
$
\item
There exist an orbifold chart $\frak V_x$ of $U_2$ at $x$
and orbifold charts $\frak V_{\tilde x_j}$ of $U_1$ at $\tilde x_j$
respectively for $j=1,\dots,m_x$ such that
$(\frak V_{\tilde x_j}, \frak V_x)$ is a covering chart.
(Here $\frak V_x$ is independent of $j$.)
We write its covering index $n_j(x)$.
\item
$\sum_{j=1}^{m_x} n_j(x)$ is independent of $x$.
\end{enumerate}
We call $\sum_{j=1}^{m_x} n_j(x)$
the {\it covering index}\index{orbifold ! covering index}  of $\pi$.
\end{defn}
\begin{rem}
\begin{enumerate}
\item
We only define a finite covering here since we do not use an infinite covering
in the present article.
\item
If Definition \ref{coveringdef} (1)-(2) is satisfied
and $U_2$ is connected, then
(3) is equivalent to the following condition:
\par\medskip
\begin{enumerate}
\item[(3)']
$
\pi^{-1}(U_x) = \bigcup_{j=1}^{m_x} U_{\tilde x_j}
$
and the right hand side is a disjoint union.
\end{enumerate}
\par\medskip\noindent
In fact, (3)' implies that $\sum_{j=1}^{m_x} n_j(x)$
is a locally constant function.
We can replace (3) by (3)' and define an infinite covering in the same way.
\item
The composition of covering maps is a covering map.
\end{enumerate}
\end{rem}
\begin{lem}\label{lem274}
Let $\varphi_{21} : U_1 \to U_2$ be an embedding of orbifolds and
$\pi_2 : \widetilde U_2 \to U_2$ a covering map of orbifolds.
We consider the fiber product
$U_1 \times_{U_2} \widetilde U_2$ in the category of topological spaces.
It comes with continuous maps
$\pi_1 : U_1 \times_{U_2} \widetilde U_2 \to U_1$
and
$\widetilde\varphi_{21} : U_1 \times_{U_2} \widetilde U_2 \to \widetilde U_2$.
Then $U_1 \times_{U_2} \widetilde U_2$ has a structure of orbifolds such that:
\begin{enumerate}
\item
$\pi_1$ is a covering map.
\item
$\widetilde\varphi_{21}$ is an embedding of orbifolds.
\end{enumerate}
The conditions (1) (2) uniquely determine the orbifold structure on
$U_1 \times_{U_2} \widetilde U_2$.
\end{lem}
\begin{proof}
We take the atlas $\{\frak V^i_{\frak r} \mid \frak r \in \frak R_i\}$
as in Lemma \ref{lem26999}.
We may choose it sufficiently fine so that
for each $\frak r \in \frak R_2$ there exist
$\frak V^{2,\frak r}_{j}$, $j=1,\dots, m_{\frak r}$
such that
$(\frak V^{2,\frak r}_{j}, \frak V^2_{\frak r})$
is a covering atlas and satisfies
$
\pi^{-1}(U^2_{\frak r}) = \bigcup_{j=1}^{m_{\frak r}} U^2_{\frak r,j}.
$
For each $j=1,\dots,m_{\frak r}$, the map $\pi$ determines a finite index subgroup
$\Gamma^2_{j,\frak r}$ of
$\Gamma^2_{\frak r}$ for $\frak r \in \frak R_2$.
Note that in case $\frak r \in \frak R_1$  the group
$\Gamma^2_{\frak r}$ is isomorphic to $\Gamma^1_{\frak r}$.
Then a finite subgroup $\Gamma^1_{j,\frak r}$ of $\Gamma^1_{\frak r}$
determines $\Gamma^2_{j,\frak r}
\subset \Gamma^2_{\frak r}$.
Therefore
the collection
$(V^1_{\frak r},\Gamma^1_{\frak r,j},\phi^1_{\frak r,j})$
determines an orbifold chart. Here the map
$\phi^1_{\frak r,j}$ is defined by
$$
\phi^1_{\frak r,j}(y) = (\phi^1_{\frak r}(y),\phi^2_{\frak r,j}(
\varphi^{\frak r}_{21}(y))),
$$
where $\varphi^{\frak r}_{21} : V^{\frak r}_1 \to V^{\frak r}_2$
is determined by the orbifold embedding $\varphi_{21}$
and $\phi^2_{\frak r,j} : V^{\frak r}_2 \to \widetilde U_2$
is a part of the orbifold chart $\frak V^{2,\frak r}_{j}$.
It is easy to check that the charts
$(V^1_{\frak r},\Gamma^1_{\frak r,j},\phi^1_{\frak r,j})$ for various
$\frak r$ and $j$ determine an orbifold structure on the fiber product.
We can easily check (1), (2). The uniqueness is also easy to check.
\end{proof}
\begin{lemdef}\label{lemdef2425}
We can pull back a vector bundle and its section by a covering map of an orbifold.
\end{lemdef}
\begin{proof}
Let $\widetilde U \to U$ be a covering and $\mathcal E \to U$ a vector bundle.
We take a local coordinate $(V\times E,\Gamma,\widehat\phi)$ of $\mathcal E \to U$
at $p \in U$ where
$(V,\Gamma,\phi)$ is the corresponding chart of $U$.
We may shrink $V$ and may assume that there is a covering
chart $(V,\Gamma_i,\phi_i)$ for each $p_i \in \tilde U$ with $\pi(p_i) = p$.
The two maps $\phi_i \circ {\rm pr}_V : V \times E \to \tilde U$
and $\widehat\phi : V \times E \to \vert\mathcal E\vert$
give rise to a map
$$
\widehat{\phi}_i : V \times E \to \widetilde U \times_{U} \vert\mathcal E\vert.
$$
It is easy to see that $(V\times E,\Gamma_i,\widehat\phi_i)$ gives
a structure of vector bundle on $\widetilde U \times_{U} \vert\mathcal E\vert
\to \widetilde U$.
\end{proof}
\begin{rem}
Our covering map in the sense of Definition \ref{coveringdef} defines a {\it good map} in the sense of \cite{ofdruan}. Then the above lemma is a special case of
\cite[Theorem 2.43]{ofdruan}.
\end{rem}
\begin{lem}\label{smoothstrucov}
Let $X$ be a topological space, $Y$ an orbifold,
and let $f : X \to Y$ be a continuous map.
Then the orbifold structure on $X$, under which $f$ becomes a
covering map of orbifolds, is unique if it exists.
\end{lem}
\begin{proof}
Suppose we have two such structures. It suffices to show that
the identity map is a diffeomorphism of orbifold.
Let $p \in Y$ and $(V,\Gamma,\phi)$ be an orbifold chart of $X$ at $p$.
We may shrink $V$ if necessary and assume that there are
covering charts $(V_i,\Gamma_i,\phi_i)$ ($i=1,2$)
for two orbifold structures on $X$ such that $f : X \to Y$ is a covering map.
We put
$$
S(V) = \{x \in V \mid \exists \gamma \in \Gamma \setminus \{1\},\,\,
\gamma x = x\}, \quad U^{\rm reg} = (V \setminus S(V))/\Gamma.
$$
Note that $f^{-1}(U^{\rm reg}) \to U^{\rm reg}$ is a covering space.
An element of $\gamma \in \Gamma$
is contained in the image of $\Gamma_i$
if and only if the corresponding loop
in $U^{\rm reg}$ lifts to a loop in $f^{-1}(U^{\rm reg})$.
Therefore $\Gamma_1 = \Gamma_2$ as subgroups of $\Gamma$.
It is now easy to see that the identity map is a
diffeomorphism at $p_i$.
\end{proof}

\subsection{Covering space of K-space}
\label{subsec:coverkura}
Using Lemma \ref{lem274} and Lemma-Definition \ref{lemdef2425},
it is fairly straightforward to define the notion of a covering space of a K-space
and show that various objects are pulled back by a covering map.
We will describe them below for completeness.
\par
\begin{shitu}\label{shitu23.8}
Let $\mathcal U = (U,\mathcal E,s,\psi)$
be a Kuranishi chart of $X$.
Suppose we are given a topological space $\widetilde X$ and a continuous map
$\pi : \widetilde X \to X$.
$\blacksquare$
\end{shitu}
\begin{defn}
We call $\widetilde{\mathcal U} = (\widetilde U,\widetilde{\mathcal E},\widetilde s,\widetilde \psi)$
a {\it covering chart} if
$\widetilde U$ is a covering space of the orbifold $U$,
$\widetilde{\mathcal E}$ is the pull-back of $\mathcal E$ by $\widetilde U \to U$,
$\widetilde s$ is obtained from $s$ by the pull-back
and $\widetilde \psi$ is a homeomorphism from $\widetilde s^{-1}(0)$ to  $\widetilde X$.
We require
$
\pi \circ \widetilde \psi = \psi.
$
\end{defn}

\begin{lemdef}\label{lemdef2310}
Suppose we are in Situation \ref{shitu23.8}.
Let $\mathcal U' = (U',\mathcal E',s',\psi')$
be another Kuranishi chart and
$\Phi : \mathcal U' \to \mathcal U$ an embedding of Kuranishi chart.
Let $\widetilde{\mathcal U}$ be a covering chart of $\mathcal U$.
Then we can define a covering chart $\widetilde{\mathcal U'}
= (\widetilde U',\widetilde{\mathcal E'},\widetilde s',\widetilde \psi')$
of $\mathcal U'$ and an embedding of
Kuranishi chart $\widetilde{\Phi} : \widetilde{\mathcal U'}
\to \widetilde{\mathcal U}$ such that
\begin{equation}\label{formula2344}
\widetilde U' = U' {}_{\Phi}\!\times_U \tilde U
\end{equation}
with the orbifold structure given by Lemma \ref{lem274}.
We call $\widetilde{\mathcal U'}$ the {\rm pull-back covering chart}
of $\widetilde{\mathcal U}$ by $\Phi$.\index{pull-back ! covering chart}
In particular, if $U_0$ is an open subset of $U$, then
we can define a restriction $\widetilde{\mathcal U'}\vert_{U_0}$
\index{restriction of a covering} of
$\widetilde{\mathcal U'}$.
\end{lemdef}
\begin{proof}
We define $\widetilde U'$ by (\ref{formula2344}).
Then $\widetilde U' \to U'$ is a covering space.
Therefore by Lemma \ref{lemdef2425} we can pull back $E$ to $\widetilde U'$.
We define $E'$ by the pull-back. We can then define $\widetilde s'$, $\widetilde \psi'$
from definition. The existence of $\widetilde{\Phi}$ is immediate from
construction.
\end{proof}

\begin{defn}\label{def2788}
Let $\widehat{\mathcal U}$ be a Kuranishi structure on $X$,
$\pi : \widetilde X \to X$ a finite to one continuous map. Let
$\widehat{\widetilde{\mathcal U}}$ be a Kuranishi structure on $\widetilde X$.
We say that $(\widetilde X,\widehat{\widetilde{\mathcal U}})$ is a {\it covering space}
of $(X,\widehat{\mathcal U})$ if the following holds.
\begin{enumerate}
\item
For each $p \in X$ we are given a covering chart $\widetilde{\mathcal U_p}$
of the Kuranishi chart $\mathcal U_p$, which is a part of the data
consisting $\widehat{\mathcal U}$.
\item
If $q \in \psi_p(s_p^{-1}(0))$, then the restriction
$\widetilde{\mathcal U_q}\vert_{U_{pq}}$ of
the covering chart $\widetilde{\mathcal U_q}$
given in Item (1)
is an open subchart of the pull-back covering chart
of $\widetilde{\mathcal U_p}$ by the coordinate change $\Phi_{pq}$.
\item
Let $p_j \in \widetilde X$ and $p \in X$ with
$\pi(p_j) = p$ and $\mathcal U_{p_j}
=  (U_{p_j},\mathcal E_{p_j},s_{p_j},\psi_{p_j})$
a Kuranishi chart of $\tilde X$ at $p_j$ which is a
part of the data of $(\widetilde X,\widehat{\widetilde{\mathcal U}})$.
Then there exists an open embedding $\Phi_{p_j}$
of the Kuranishi chart $\mathcal U_{p_j}$
to the covering chart
$\widetilde{\mathcal U}_{p}$ given in Item (1).
\item
For any $q \in \psi_p(s_p^{-1}(0))$, $q_j \in \psi_{p_j}(s_{p_j}^{-1}(0))$,
the following diagram commutes.
\begin{equation}\label{diag27XX}
\begin{CD}
\mathcal U_{q_j}\vert_{U_{p_jq_j}} @ > {\Phi_{p_j q_j}} >>
{\mathcal U}_{p_j}  \\
@ V{\Phi_{q_j}}VV @ VV{\Phi_{p_j}}V\\
\widetilde{\mathcal U}_{q}\vert_{U_{pq}} @ > {\widetilde\Phi_{pq}} >>\widetilde{\mathcal U}_{p}
\end{CD}
\end{equation}
Here $\Phi_{p_j q_j}$ is the coordinate change
of the Kuranishi structure $(\widetilde X,\widehat{\widetilde{\mathcal U}})$
and $\widetilde\Phi_{pq}$ is
obtained by Lemma-Definition \ref{lemdef2310}.
\item
Let $n_{p,j}$ be the covering index of the covering
$\widetilde U_{p,j} \to U_{p}$. Then the number
\begin{equation}\label{coveringindex}
\sum_{p_j \in \widetilde X: \pi(p_j) = p} n_{p,j}
\end{equation}
is independent of $p$.
We call it   the {\it covering index} of
the covering
$(\widetilde X,\widehat{\widetilde{\mathcal U}})$ of $(X,\widehat{\mathcal U})$.
\end{enumerate}
\end{defn}
\begin{rem}
\begin{enumerate}
\item
If $X$ is connected, Condition (5) follows from
Conditions (1) (2).
In fact, Conditions (1) (2) imply that the covering index of $\tilde U_p \to U_p$
is locally constant. However it does not seem to be a good idea
to assume connectivity of $X$ in our situation,
since the topology of $X$ can be pathological.
\item
The commutativity of Diagram \ref{diag27XX}
means the set theoretical equalities
of the underlying topological spaces of orbifolds and of
the total spaces of obstruction bundles.
We can safely do so since all
the maps involved are embeddings.
\item
We can define the notion of a covering space of a space equipped with a good coordinate
system and can prove that
for a given covering space of K-space
we can construct a covering space equipped with a good coordinate system
in the same way.
\end{enumerate}
\end{rem}

\begin{lem}
Let $(\widetilde X,\widehat{\widetilde{\mathcal U}})$ be a covering space
of $(X,\widehat{\mathcal U})$ in the sense of Definition \ref{def2788}.
\begin{enumerate}
\item
If $\widehat{\frak S}$ is a CF-perturbation of
$\widehat{\mathcal U}$, it induces a CF-perturbation
$\widehat{\widetilde{\frak S}}$ of
$\widehat{\widetilde{\mathcal U}}$.
\item
A strongly continuous (strongly smooth) map $\widehat f$
from $(X,\widehat{\mathcal U})$ can be pulled back to a strongly continuous (strongly smooth) map
$\widehat{\widetilde f}$ from
$(\widetilde X,\widehat{\widetilde{\mathcal U}})$.
\item
A differential form $\widehat h$ on $(X,\widehat{\mathcal U})$ can be pulled back to
$\widehat{\widetilde h}$ on $(\widetilde X,\widehat{\widetilde{\mathcal U}})$.
\item
A strongly smooth map $\widehat f$ is strongly submersive with respect to
$\widehat{\frak S}$ if and only if $\widehat{\widetilde f}$ in (2)
is strongly submersive with respect to
$\widehat{\widetilde{\frak S}}$.
\item
The statements such that
`CF-perturbations' in (2) (4) are replaced by
`multivalued perturbations' also hold.
\item
In the situation of (4) we have
$$
n\widehat f! (\widehat h;\widehat{\frak S})
=
\widehat{\widetilde f}(\widehat{\widetilde h};\widehat{\widetilde{\frak S}}).
$$
Here $n$ is the covering index.
\end{enumerate}
\end{lem}
The proof is straightforward and so omitted.
\subsection{Covering spaces associated to the corner structure
stratification}
\label{subsec:coverconer}
One of the main reasons we introduced the notion of covering space
of a K-space is that we use the next lemma to
clarify the discussion of normalized boundary.
\begin{prop}\label{lem27111}
If $(X,\widehat{\mathcal U})$ is a K-space with corners,
then $\overset{\circ} S_{k-1}(\partial(X,\widehat{\mathcal U}))$
is a covering space of  $\overset{\circ} S_{k}(X,\widehat{\mathcal U})$
with covering index $k$.
\end{prop}
\begin{proof}
We first prove the proposition for the case of orbifolds.
\begin{lem}\label{lem2712}
Let $U$ be an orbifold with corner.
The map $\overset{\circ}S_{k-1}(\partial U) \to \overset{\circ}S_k U$
which is the restriction of the map $\pi$ in \cite[Lemma 8.8 (2)]{part11}
is a $k$-fold covering of orbifolds.
If $\mathcal E$ is a vector bundle on $U$, then the
bundle induced on $\overset{\circ}S_{k-1}(\partial U)$
is canonically isomorphic to the pull-back of
the restriction of $\mathcal E$ to $\overset{\circ}S_k U$.
\par
\end{lem}
\begin{proof}
Let $x \in \overset{\circ}S_k U$ and $\frak V_x =
(V_x,\Gamma_x,\phi_x)$ be an orbifold chart of $U$ at $x$.
We may assume $V_x \subset \overline{V_x} \times [0,1)^{k}$ and
$o_x = (\overline o_x,(0,\dots,0))$, where
$\overline{V_x}$ is a manifold without boundary.
There exists a group homomorphism
$\sigma : \Gamma_x \to {\rm Perm}(k)$
such that if
$\gamma(\overline y,(t_1,\dots,t_k))
= (\overline y',(t'_1,\dots,t'_k))$ then
$t'_k = 0$ if and only if $t_{\sigma(\gamma)(k)}=0$.
\par
We consider $I \subset \{1,\dots,k\}$ a complete set of representatives
of $\{1,\dots,k\}/\Gamma_x$.
For each $i \in I$, we put
$$
\aligned
\Gamma_{x,i} &= \{\gamma \in \Gamma_x \mid \sigma(\gamma) i = i\},
\\
\partial_iV_x
&=
\{(\overline y,(t_1,\dots,t_k)) \in V_x \mid t_i = 0\}.
\endaligned$$
The given $\Gamma_x$ action on $V_x$
induces a $\Gamma_{x,i}$ action on $\partial_iV_x$.
By definition of normalized boundary, there exists a map
$\phi_{x,i} : \partial_iV_x \to \partial U$
such that $(\partial_iV_x,\Gamma_{x,i},\phi_{x,i})$
is an orbifold chart of $\partial U$ at $\tilde x_i$.
Here $\{\tilde x_i \mid i \in I\} = \pi^{-1}(x) \subset \overset{\circ}S_{k-1}(\partial U)$.
An orbifold chart of $\overset{\circ}S_{k-1}(\partial U)$
at $\tilde x_i$ is
$$
(\overline{V_x},\Gamma_{x,i},\phi_{x,i}\vert_{\overline{V_x}}),
$$
where $\overline{V_x}$ is identified with the subset
$\overline{V_x} \times \{0\}$ of $\partial_iV_x$.
\par
On the other hand,
an orbifold chart of $\overset{\circ}S_k U$ at $x$  is
$(\overline{V_x},\Gamma_{x},\phi_{x}\vert_{\overline{V_x}})$.
Since $\#(\Gamma_{x}(i)) = \# (\Gamma_x/\Gamma_{x,i})$,
we have
$$
\sum_{i \in I}\# (\Gamma_x/\Gamma_{x,i}) = k.
$$
We have thus proved that $\pi : \overset{\circ}S_{k-1}(\partial U) \to
\overset{\circ}S_k (U)$ is a $k$-fold covering of orbifolds.
We can prove the second half of the lemma using the above
description of the orbifold charts.
\end{proof}
Now we consider the case of Kuranishi structure.
Let $\mathcal U_p = (U_p,\mathcal E_p,\psi_p,s_p)$ be a Kuranishi chart of
$\widehat{\mathcal U}$ at $p \in \overset{\circ} S_{k}(X,\widehat{\mathcal U})$.
We put
$$
\overset{\circ} S_{k}(\mathcal U_p) :=
(\overset{\circ} S_{k}(U_p),\mathcal E_p\vert_{\overset{\circ} S_{k}},\psi_p\vert_{{\overset{\circ} S_{k}}},s_p\vert_{{\overset{\circ} S_{k}}}),
$$
which is a Kuranishi chart at $p$ of the K-space
$\overset{\circ} S_{k}(X,\widehat{\mathcal U})$.
Let $\pi$ be the map from the underlying topological space
of $\overset{\circ} S_{k-1}(\partial(X,\widehat{\mathcal U}))$
to that of $\overset{\circ} S_{k}(X,\widehat{\mathcal U})$.
We use the notation in the proof of Lemma \ref{lem2712}
by putting $x = o_p \in U_p$, $U = U_p$ and $\mathcal E = \mathcal E_p$.
Then $\pi^{-1}(p)$ consists of $\# I$ points which
we write $\widetilde p_{i}$, $i\in I$.
Then $(\overline{V_x},\Gamma_{x,i},\psi_{x,i}\vert_{\overline{V_x}})$
is an orbifold chart of $\overset{\circ} S_{k-1}(\partial(X,\widehat{\mathcal U}))$
at $\widetilde p_i$.
We restrict $\overset{\circ} S_{k}(\mathcal U_p)$
to $\overline{V_x}/\Gamma_x \subset U_p$ to get
a Kuranishi chart of
$\overset{\circ} S_{k-1}(\partial(X,\widehat{\mathcal U}))$.
\par
The second half of Lemma \ref{lem2712}
implies that the pull-back of the restriction of the obstruction bundle $\mathcal E_p$
to $\overline{V_x}/\Gamma_x \subset U_p$
defines a vector bundle on $\overline{V_x}/\Gamma_{x,i}$
which is isomorphic to the obstruction bundle
of the Kuranishi structure $\overset{\circ} S_{k-1}(\partial(X,\widehat{\mathcal U}))$.
It is easy to see from construction
that Kuranishi maps are preserved by this isomorphism.
\par
We have thus constructed open substructures of
$\overset{\circ} S_{k-1}(\partial(X,\widehat{\mathcal U}))$
and
$\overset{\circ} S_{k}(X,\widehat{\mathcal U})$
so that the Kuranishi chart of the former is a $k$ fold covering
of the Kuranishi chart of the latter.
It is easy to see that this isomorphism is compatible with
the coordinate change. Hence
the proof of Proposition \ref{lem27111} is complete.
\end{proof}
We next generalize Proposition  \ref{lem27111}
to the corners of arbitrary codimension.
\begin{prop}\label{prop2813}
Let $(X,\widehat{\mathcal U})$ be an $n$-dimensional K-space.
Then for each $k$ there exists an
$(n-k)$-dimensional K-space $\widehat S_k(X,\widehat{\mathcal U})$ with
corners and maps
$\pi_k :  \widehat S_k(X,\widehat{\mathcal U})  \to
S_k(X,\widehat{\mathcal U})$,
$\pi_{\ell,k} : \widehat S_{\ell}(\widehat S_k(X,\widehat{\mathcal U}))
\to \widehat S_{k+\ell}(X,\widehat{\mathcal U})$
with the following properties:
\begin{enumerate}
\item
$\pi_k$ is a continuous map between underlying topological spaces.
\item
$\widehat S_1(X,\widehat{\mathcal U})$ is the normalized boundary
$\partial(X,\widehat{\mathcal U})$.
\item The interior of $\widehat S_k(X,\widehat{\mathcal U})$
is isomorphic to $\overset{\circ}S_k(X,\widehat{\mathcal U})$.
The underlying homeomorphism of this isomorphism is the
restriction of $\pi_k$.
\item
$\pi_{\ell,k}$ is an $(\ell+k)!/\ell!k!$ fold covering map of K-spaces.
\item
The following objects on $(X,\widehat{\mathcal U})$ induce the corresponding ones on
$\widehat S_k(X,\widehat{\mathcal U})$.
The induced objects are compatible with the covering maps $\pi_{\ell,k}$.
\begin{enumerate}
\item
CF-perturbation.
\item
Multivalued perturbation.
\item
Differential form.
\item Strongly continuous map. Strongly smooth map.
\item Covering map.
\end{enumerate}
\item
The following diagram commutes.
\begin{equation}\label{diagin26277XX2}
\begin{CD}
\widehat S_{k_1}(\widehat S_{k_2}(\widehat S_{k_3}(X,\widehat{\mathcal U}))) @ >{\pi_{k_1,k_2}}>>
\widehat S_{k_1+k_2}(\widehat S_{k_3}(X,\widehat{\mathcal U})) \\
@ V{\widehat S_{k_1}(\pi_{k_2,k_3})}VV @ VV{\pi_{k_1+k_2,k_3}}V\\
\widehat S_{k_1}(\widehat S_{k_2+k_3}(X,\widehat{\mathcal U})) @ > {\pi_{k_1,k_2+k_3}} >> \widehat S_{k_1+k_2+k_3}(X,\widehat{\mathcal U})
\end{CD}
\end{equation}
Here $\widehat S_{k_1}(\pi_{k_2,k_3})$ is the covering map induced from $\pi_{k_2,k_3}$.
\item
Let $f_i : (X_i,\widehat{\mathcal U_i}) \to M$ be a strongly smooth map
and assume that $f_1$ is transversal to $f_2$. Then
$$
\widehat S_{k}\left((X_1,\widehat{\mathcal U_1}) \times_M (X_2,\widehat{\mathcal U_2})\right)
\cong
\coprod_{k_1+k_2=k}
\widehat S_{k_1}(X_1,\widehat{\mathcal U_1}) \times_M \widehat S_{k_2}(X_2,\widehat{\mathcal U_2}).
$$
Here the right hand side is the disjoint union.
\item
(1)-(6) also hold when we replace `Kuranishi structure' by `good coordinate system'.
\item
Various kinds of embeddings of Kuranishi structures and/or good coordinate systems induce
the corresponding ones of
$\widehat S_{k}(X,\widehat{\mathcal U})$.
\end{enumerate}
\end{prop}
We note that
Proposition \ref{prop2813} (7) implies
\begin{equation}\label{form3077}
\widehat S_{\ell}(X_1 \times_{M_1} \dots \times_{M_{n-1}} X_n)
=
\coprod_{\ell_1+\dots+\ell_n = \ell}
(\widehat S_{\ell_1}(X_1) \times_{M_1} \dots \times_{M_{n-1}} \widehat S_{\ell_n}(X_n)).
\end{equation}
\begin{defn}\label{norcor}
We call $\widehat S_{k}(X,\widehat{\mathcal U})$ the
{\it normalized (codimension $k$) corner} of $(X,\widehat{\mathcal U})$.
\index{normalized codimension $k$ corner}
\index{corner ! normalized codimension $k$ corner}
\end{defn}
\begin{proof}[Proof of Proposition \ref{prop2813}]
Let $M$ be a manifold with corners. We first define $\widehat{S}_k(M)$.
Let $x \in \overset{\circ}{S}_m(M)$.
We take a chart $\frak V_x = (V_x ,\psi_x)$ so that
$V_x \cong \overline V_x \times [0,1)^m$.
Let $A \subset \{1,\dots,m\}$ with $\# A = k$.
A pair $(x,A)$ becomes an element of $\widehat{S}_k(M)$.
\par
We next define a topology on $\widehat{S}_k(M)$.
Let $y = \psi(\tilde y)$ with $\tilde y \in V_x$.
We write $\tilde y = (\tilde y_0,(t_1,\dots,t_m))$.
If $t_i = 0$ for all $i\in A$, we consider elements
$y_A = (y,A) \in \widehat{S}_k(M)$ as follows.
Suppose $B = \{ i \mid t_i = 0\} \supset A$.
Let $W$ be a neighborhood of $(t_i)_{i\notin B}$ in $(0,1)^{\{1,\dots,m\} \setminus B}$.
Then $\overline V \times W \times [0,1)^B$
together with the restriction of $\psi_x$ is a chart of $y$.
Thus we have $(y,A) \in \widehat S_{k}(U)$.
We say $(y^a,A)$ above converges to $(x,A)$ if $\tilde y^a$ converges to $o_x$,
where $o_x$ is the point such that $\psi_x(o_x) = x$.
\par
It is easy to see that $\widehat{S}_k(M)$ with this topology becomes
a manifold with corners.
This construction is canonical so that it induces one of orbifolds and
of Kuranishi structure.
(The proof of this part is entirely similar to the case of normalized boundary
and so is omitted.)
\par
We next construct the covering map $\pi_{\ell,k}$.
We consider the case of manifolds.
Let $x \in \overset{\circ}{S}_m(M)$ and let $\frak V_x$, $A$ be as above.
For simplicity of notation, we put $A=\{1,\dots,k\}$.
Suppose $(x,A) \in S_{\ell}(\widehat S_k(M))$.
It implies $m \ge k+\ell$.
\par
By definition the neighborhood of $(x,A)$ in $\widehat S_k(M)$
is $(y,A)$ where
$y \in \overline V \times \{(0,\dots,0)\} \times [0,1)^{m-k}$.
Therefore
a point $\tilde x$ in $\widehat S_{\ell}(\widehat S_k(M))$ such that
$\pi_{\ell}(\tilde x) = (x,A)$
corresponds one to one to the set $A^+ \supset A$ with $\#A^+ = \ell+k$.
We put $B = A^+ \setminus A$.
We thus may regard $(x,A,B) \in \widehat S_{\ell}(\widehat S_k(M))$.
(Here $\# B = \ell$.)
Now we define the map $\pi_{\ell,k} : \widehat S_{\ell}(\widehat S_k(M)) \to \widehat S_{\ell+k}(M)$ by
$$
\pi_{\ell,k}(x,A,B) = (x,A\cup B).
$$
Given $(x,C) \in \widehat S_{\ell+k}(M)$, the element in the fiber of $\pi_{\ell,k}$ corresponds one to one
to the partition of $C$ into $A\cup B$ where $\#A = k$ and $\#B = \ell$.
We can use this fact to show that $\pi_{\ell,k}$ is a covering map of covering index $(k+\ell)!/k!\ell!$.
\par
We have thus constructed the covering map $\pi_{\ell,k}$ in the case of manifolds.
To prove the case of
orbifolds and of K-spaces, it suffices to observe that this construction
is canonical and so is compatible with various kinds of coordinate changes.
\par
Once the K-space $\widehat{S}_k(X,\widehat{\mathcal U})$ and the covering map
$\pi_{\ell,k}$ are defined as above, it is very easy to check the properties (1)-(9).
\end{proof}
\begin{rem}
As well as other parts of this article,
Proposition \ref{prop2813} is {\it not} new.
Especially we would like to mention that mostly the same construction
appeared in D. Joyce's paper \cite{joyce3}.
(The article \cite{joyce3} discusses the case of manifolds.
However it is straightforward to generalize the story to the
case of K-spaces.)
In \cite{joyce3} the notion of the boundary $\partial X$ is defined,
which is the same as our definition of normalized boundary.
Then the action of ${\rm Perm}(k)$ on
$\underbrace{\partial \cdots \partial}_{\text{$k$ times}}X$ is introduced.
The quotient space $\underbrace{\partial \cdots
\partial}_{\text{$k$ times}}X/{\rm Perm}(k)$
(which is denoted by $C_kX$ in \cite{joyce3}),
coincides with our $\widehat S_k(X)$.
\end{rem}

\subsection{Finite group action on K-space}
\label{subsec:finitegroupact}
\begin{defn}
Let $X$ be an orbifold and $G$ a finite group.
A $G$ action on $X$ as a topological space
is said to be a {\it $G$ action on the orbifold} $X$ if
the homeomorphism $X \to X$ induced by each element of $G$ is a
diffeomorphism of orbifold.
\par
Two actions are said to be {\it the same} if they are the same as maps $G \times X \to X$,
set-theoretically.
\end{defn}
\begin{lem}\label{2918lem}
Let $X$ be an orbifold on which a finite group $G$ acts (as orbifold).
We assume that the action is effective on each connected component. Then
there exists a unique orbifold structure on $X/G$ such that
$X$ is a covering space of $X/G$ and the natural
map $X \to X/G$ is a covering map.
\end{lem}
\begin{proof}
Let $x \in X$ and $G_x = \{g \in G \mid gx = x\}$.
We take an orbifold chart $\frak V_x = (V_x,\Gamma_x,\psi_x)$
such that $U_x$ is $G_x$ invariant (by using a $G$ invariant Riemannian
metric, for example).
Using the effectivity of $G$ action on each connected component,
we can easily show that the $G_x$ action on $U_x$ is effective.
For each $g \in G_x$ we obtain a map $\varphi_g : V_x \to V_x$ and
a group homomorphism $h_g : \Gamma_x \to \Gamma_x$.
Since $\varphi_{g_1}\varphi_{g_2}$ induces the same continuous map as
$\varphi_{g_1 g_2}$ between the underlying topological spaces, there exists a unique
element $\gamma_{g_1g_2g_3} \in \Gamma_x$
such that
$$
\varphi_{g_1}\varphi_{g_2} = \gamma_{g_1 g_2}\varphi_{g_1 g_2}.
$$
Moreover we have
$$
h_{g_1}h_{g_2} = {\rm conj}_{\gamma_{g_1 g_2}}h_{g_1 g_2}.
$$
Note that $\varphi_g$ is $h_g$ equivariant.
Then we define a group structure on the direct product set
$\Gamma_x \times G_x$ by
\begin{equation}
(\gamma_1,g_1)\circ (\gamma_2,g_2)
=
(\gamma_1 h_{g_1}(\gamma_2) \gamma_{g_1,g_2}, g_1g_2).
\end{equation}
We define $\cdot : (\Gamma_x \times G_x) \times V_x \to V_x$
by
$$
(\gamma,g)\cdot x = \gamma(\varphi_{g}(x)).
$$
Then we observe
$$
\aligned
&(\gamma_1,g_1)\cdot((\gamma_2,g_2))\cdot x)
=
(\gamma_1,g_1)\cdot \gamma_2(\varphi_{g_2}(x))
=
\gamma_1(\varphi_{g_1}(\gamma_2(\varphi_{g_2}(x)))) \\
&= \gamma_1 h_{g_1}(\gamma_2) \varphi_{g_1}(\varphi_{g_2}(x))
=
\gamma_1 h_{g_1}(\gamma_2) \gamma_{g_1g_2}(\varphi_{g_1g_2}(x))
=
((\gamma_1,g_1)\circ (\gamma_2,g_2))\cdot x.
\endaligned
$$
It follows from effectivity that $\circ$ defines a group structure.
We denote this group by $\Gamma_x \tilde\times G_x$.
\par
Now we define $\overline\phi_x : V_x \to X/G$ by the composition of $\phi_x :
V_x \to X$ and the projection $X \to X/G$.
Then
it is easy to see that $(V_x,\Gamma_x \tilde\times G_x,\overline \phi_x)$
defines an orbifold structure on $X/G$.
The rest of the proof is obvious.
\end{proof}
Now we define the definition of the action of finite group on K-space.
\begin{defn}
Let $(X,\widehat{\mathcal U})$ be a K-space.
An {\it automorphism}
\index{K-space ! automorphism}
$\Phi$ consists of a pair $(\vert\Phi\vert,\{\Phi_p\})$
of a homeomorphism $\vert\Phi\vert : X \to X$ and
an assignment $X \ni p \mapsto \Phi_p =(\varphi_p,\widehat\varphi_p)$
with the following properties:
\begin{enumerate}
\item $\varphi_p : U_p \to U_{\vert\Phi\vert(p)}$ is a diffeomorphism of orbifolds.
\item $\widehat\varphi_p : E_p \to E_{\vert\Phi\vert(p)}$ is an
isomorphism of vector bundles over $\varphi_p$.
\item
$\widehat\varphi_p \circ s_p = s_{\vert\Phi\vert(p)} \circ \varphi_p$
holds on $U_p$.
\item
$\vert\Phi\vert \circ \varphi_p
= \varphi_{\vert\Phi\vert(p)} \circ \varphi_p$ holds
on $s_p^{-1}(0)$.
\item
Let $q \in \psi_(s_p^{-1}(0))$. Suppose $\Phi_{pq} =
(U_{pq},\varphi_{pq},\widehat\varphi_{pq})$
and
$$
\Phi_{\vert\Phi\vert(p)\vert\Phi\vert(q)} =
(U_{\vert\Phi\vert(p)\vert\Phi\vert(q)},\varphi_{\vert\Phi\vert(p)\vert\Phi\vert(q)},\widehat\varphi_{\vert\Phi\vert(p)\vert\Phi\vert(q)})
$$
are the coordinate
changes.
Then we have the following:
\begin{enumerate}
\item
$\varphi_q(U_{pq}) = U_{\vert\Phi\vert(p)\vert\Phi\vert(q)}$.
\item
$\varphi_{\vert\Phi\vert(p)\vert\Phi\vert(q)} \circ \varphi_q
= \varphi_p \circ \varphi_{pq}$.
\item
$\widehat\varphi_{\vert\Phi\vert(p)\vert\Phi\vert(q)} \circ \widehat\varphi_q
= \widehat\varphi_p \circ \widehat\varphi_{pq}$.
\end{enumerate}
\end{enumerate}
We say $(\vert\Phi\vert,\{\Phi_p\})$ is {\it the same} as
$(\vert\Phi'\vert,\{\Phi'_p\})$  if
$\vert\Phi\vert = \vert\Phi'\vert$ and
$\Phi_p = \Phi'_p$ for all $p$.
\end{defn}
\begin{rem}
\begin{enumerate}
\item
The equality $\Phi_p = \Phi'_p$ has an obvious meaning. Namely
we defined the notion of two diffeomorphisms and bundle isomorphisms
of orbifolds to be the same. (That is, they are the same set-theoretically.)
\item
It happens that two automorphisms $\Phi$ and $\Phi'$ with the
same underlying homeomorphisms $\vert\Phi\vert$ and $\vert\Phi'\vert$
could be different.
\end{enumerate}
\end{rem}
\begin{defnlem}
\begin{enumerate}
\item
We can compose two automorphisms of K-spaces.
The composition is again an automorphism.
\item
The set of automorphisms of a given K-space $(X,\widehat{\mathcal U})$
is a group whose
product is the composition of automorphisms.
We denote this group by ${\rm Aut}(X,\widehat{\mathcal U})$.
\item An action of a finite group $G$ on
$(X,\widehat{\mathcal U})$ is, by definition, a
group homomorphism $G \to {\rm Aut}(X,\widehat{\mathcal U})$.
\item
A $G$ action on $(X,\widehat{\mathcal U})$ induces a $G$ action
on the underlying topological space $X$.
\end{enumerate}
\end{defnlem}
\begin{defn}
An action of a finite group $G$ on a K-space $(X,\widehat{\mathcal U})$
is said to be {\it effective}
\index{K-space ! effective action on K-space}
if the following holds for each $p \in X$.
\par
We put $G_p = \{g \in G \mid gp = p\}$.
Let $U_p$ be the Kuranishi neighborhood of $p$.
By definition  $G_p$ acts on
$U_p$. We require that this action is effective on each connected component of
$U_p$.
\end{defn}
\begin{lem}
Suppose a finite group $G$ acts effectively on a K-space
$(X,\widehat{\mathcal U})$.
Then  there exists a unique Kuranishi structure on $X/G$
such that the projection $X \to X/G$ is an underlying map of the covering map and
each ${\mathcal U}_p$ gives a covering chart of this covering.
\end{lem}
\begin{proof}
This follows from Lemma \ref{2918lem}.
\end{proof}

\section{Admissible Kuranishi structure}\label{sec:admKura}

In this section we introduce the notion of admissible
Kuranishi structure. For this purpose we introduce the notion of
admissible orbifold, admissible vector bundle, and various admissible objects associated to them, like
admissible section, admissible Riemannian metric, etc., and
provide their fundamental properties.
Roughly speaking, `admissibility' in this section
is some condition that objects in question obey
certain exponential decay estimates at asymptotic ends.
Here we regard boundary or corner points as the end points, so
the `exponential decay estimates at asymptotic ends' means
the exponential decay estimates in the direction normal to boundary or corners.
\subsection{Admissible orbifold}
\label{subsec:admord}
Firstly, we discuss admissibility for the case of manifold with corner before
going to the case of orbifold with corner, because the key idea can be seen
even for the case of manifold.
\begin{shitu}\label{Situationcorner}
Let $V \subset \overline{V} \times [0,1)^k$ be an open subset where $\overline V$ is a
manifold without corner.
We denote ${\bf t} = (t_1,\dots,t_k) \in [0,1)^k$.
$\blacksquare$
\end{shitu}

\begin{conven}
\begin{enumerate}
\item We put
\begin{equation}\label{Tandt}
T_i = e^{1/t_i}, \qquad (\text{i.e.}\,\, t_i = \frac{1}{\log T_i}).
\end{equation}
We consider the corner structure stratification of $V$.
Then each connected component of open stratum $\overset{\circ}S_{\ell}V$
has coordinates that are union of the coordinate of $\overline V$ and
$k-\ell$ of $T_i$'s. Here $T_i \in [0,\infty)$.
We consider the $C^m$ norm of a function $f : V  \to \R$
stratumwisely (i.e., the norm of the differential in the stratum direction) using the above coordinate $T_i$.
(Namely we use $T_i$ and not $t_i$ to define the $C^m$ norm.)
\item
For a function $f : V  \to \R$, we denote by $\vert f\vert_{C^m}$
the pointwise $C^m$ norm in the above sense.
Thus $\vert f\vert_{C^m}$ can be regarded as a non negative real valued
function on $V$.
To define it we use a certain Riemannian metric on
$\overline V$.
(We use the standard metric for $T_i \in [0,\infty)$.)
Since we only consider its value on a compact subset, the difference of
the metric
affects only by a bounded ratio. So we do not need to care about the difference of the metric.
\end{enumerate}
\end{conven}
\begin{defn}\label{defn31333}
\begin{enumerate}
\item
We say a function $f : V \to \R$ is {\it admissible}\index{admissible ! function}
if for each compact subset $K$ and $m>0$, there exist
$\sigma(m,K)>0$ and $C(m,K)>0$ such that the following holds for each $i$.
\begin{equation}
\left\vert
\frac{\partial f}{\partial T_{i}}
\right\vert_{C^m} \le C(m,K) e^{-\sigma(m,K) T_i}.
\end{equation}
\item
We say a function $f : V \to \R$ is {\it exponentially small near the boundary}
if for each each compact subset $K$ and $m>0$, there exist
$\sigma(m,K)>0$ and $C(m,K)>0$ such that the following holds for each $i$.
\begin{equation}
\left\vert
f
\right\vert_{C^m} \le C(m,K) e^{-\sigma(m,K) T_i}.
\end{equation}
\end{enumerate}
%
\end{defn}
\begin{exm}\label{exam294}
If $k=1$, an admissible function $f(\overline y,t)$
is written as the form
$$
f(\overline y,t) = f_0(\overline y)
+ f_1(\overline y,t)
$$
such that
$f_1(\overline y,1/\log T)$ decays in an exponential order
in $T$. (Here $t = 1/\log T$.)
\end{exm}
In a way similar to the above example, we can prove the following.
\begin{lem}\label{canonicalformforadmi}
On a subset $K \times [0,c)^k$, any admissible function $f$ is written
in the following form:
\begin{equation}\label{decomp}
f(\overline y,(t_1,\dots,t_k))
=
\sum_{I \subseteq \{1,\dots,k\}}
f_I(\overline y,t_{I}).
\end{equation}
Here $t_{I} = (t_i)_{i\in I}$ and
$f_I$ is a function on $K \times [0,c)^I$
which is exponentially small near the boundary.
\end{lem}
\begin{proof}
Firstly we observe that the set of functions of the form (\ref{decomp})
forms an $\R$ vector space.
We also note that
if $\overline V \times [0,c)^k \to \overline V \times [0,c)^I$
is a projection then the pull-back of admissible
functions are admissible.
The same holds for
`exponentially small near the boundary'
and `of the form (\ref{decomp})'.
\begin{sublem}
If an admissible function is zero on the boundary, then
it is exponentially small near the boundary.
\end{sublem}
\begin{proof}
We have
$$
f(\overline y;(t_1,\dots,t_k))
=
-\int_{T_i = e^{1/t_i}}^{\infty}
\frac{\partial f}{\partial T_i}dT_i.
$$
This implies the sublemma.
\end{proof}
We are going to construct $f_{I}$ in \eqref{decomp} by an upward induction on $\#I$
such that
$$
f(\overline y,(t_1,\dots,t_k))
-
\sum_{I \subseteq \{1,\dots,k\}, \# I\le m}
f_I(\overline y,t_{I})
= 0
$$
on $S_{k-m}(K \times [0,c)^k)$.
\par
For $m=0$ we put $f_{\emptyset}(\overline y)
= f(\overline y;(0,\dots,0))$.
Suppose we have $f_{I'}$ for $I'$ with $\#I' < m$.
Let $I \subset \{1,\dots,k\}$ with $\# I = m$.
By induction hypothesis,
we may assume that $f$ is zero on $S_{k-m+1}(K \times [0,c)^k)$.
We embed $K \times [0,c)^I$ into
$K \times [0,c)^k$ by putting $t_i = 0$ for $i \notin I$.
Restricting $f$ to its image we obtain an
admissible function on $K \times [0,c)^I$.
Since we assumed $f$ is zero on $S_{k-m+1}(K \times [0,c)^k)$,
then $f =0$ on $\partial(K \times [0,c)^I)$.
We define $f_I = f\vert_{\partial(K \times [0,c)^I)}$.
Then $f_I$ is exponentially small near the boundary
and therefore its pull-back to $K \times [0,c)^k$
is of the form (\ref{decomp}).
By taking
$$
S_{k-m}(K \times [0,c)^k) =
\bigcup_{I \subseteq \{1,\dots,k\}, \# I= m}K \times [0,c)^I
$$
into account, we can see that
$$
f - \sum_{I \subseteq \{1,\dots,k\}, \# I= m}f_I
$$
is zero on $S_{k-m}(K \times [0,c)^k)$.
The proof of the lemma is complete by induction.
\end{proof}
\begin{defn}\label{defn297}
\begin{enumerate}
\item
Let $V_i \subset \overline{V}_i \times [0,1)^k$ $(i=1,2)$ be open subsets
as in Situation \ref{Situationcorner}
and let $\varphi_{21} : V_1 \to V_2$ be an embedding of manifolds.
We say that $\varphi_{21}$ is an
{\it admissible embedding}\index{admissible ! embedding} if there exists
a permutation $\sigma : \{1,\dots,k\}
\to  \{1,\dots,k\}$ such that
$
\varphi(\overline y,(t_1,\dots,t_k))
=
(\overline y',(t'_1,\dots,t'_k))
$
has the following properties.
\begin{enumerate}
\item
The coordinates of $\overline y'$ are admissible
in the sense of Definition \ref{defn31333}.
\item
If we put $T'_i = e^{1/t'_{i}}$,
then for each $i$,
$
T'_{\sigma(i)} - T_i
$
is admissible in the sense of Definition \ref{defn31333}
near the boundary.
\end{enumerate}
\item
An admissible embedding is said to be an {\it admissible diffeomorphism}
\index{admissible ! diffeomorphism}
if it is also a homeomorphism.
\item
An action of a finite group $\Gamma$ on $V \subset \overline{V} \times [0,1)^k$
is said to be an {\it admissible action}\index{admissible ! action}
if each element of $\Gamma$ induces an admissible diffeomorphism.
\item
An orbifold chart in the sense of Definition \ref{2661} (1) is said to be an
{\it admissible chart}\index{admissible ! chart}
if the $\Gamma$ action is admissible.
\end{enumerate}
\end{defn}
\begin{rem}
\begin{enumerate}
\item
In the geometric context of pseudo-holomorphic curves
we took $T$ to be the `length' of the neck region
in \cite{fooobook2} and etc..
For this choice, Definition \ref{defn297} (1)(b) is satisfied.
See Subsection \ref{subset:admphc}.
\item
The choice of the coordinate $t = 1/T$ used in \cite{fooobook2} is different from
that of (\ref{Tandt}). See also Subsection \ref{subset:admphc} for this point.
\end{enumerate}
\end{rem}
\begin{lem}\label{newlem319}
Let $\varphi_{21} : V_1 \to V_2$ be an admissible embedding
as in Definition \ref{defn297}.
\begin{enumerate}
\item
An admissible embedding $\varphi_{21}$ induces a smooth embedding
$$
\varphi_{21} : \overline{V}_1 \times [0,1)^k
\to
\overline{V}_2 \times [0,1)^k.
$$
\item
Denote by $\varphi_{21}^{j}$
the $j$-th component of the
$[0,1)^k$ factor of $\varphi_{21}$.
If we put
\begin{equation}\label{eq:S/t}
S_i = \frac{1}{t_i},
\end{equation}
then for each compact set $K \subset V$ and $m\ge 0$ there exist
$C(m,K)>0$ and $\sigma(m,K) > 0$ such that
\begin{equation}\label{eq:CCexpS}
\left\Vert
\varphi_{21}^{\sigma(i)} - t_i
\right\Vert_{C^{m,S}_K}
\le
C(m,K) e^{-\sigma(m,K) S_i}.
\end{equation}
Here $\sigma(i) \in \{1,\dots ,k\}$ and we identify $V \cong \overline V \times (1,\infty]^k$
using $S_i$ as the coordinates of the second factor and
$C^{m,S}_K$ stands for the $C^m$ norm on $K$ with respect to the $S_i$ coordinates.
\end{enumerate}
\end{lem}
\begin{proof}
If $T'_{\sigma(i)} = T_i$
for all $i$ in addition,
then it is easy to see that an admissible embedding induces a smooth embedding.
So it suffices to consider the case when
$\overline y' = \overline y$ and $\sigma (i) = i$.
We also note that (1) follows from (2).
Then by Definition \ref{defn297} (1)(b) and \eqref{Tandt} we have
$$
t'_i = (\log (e^{\frac{1}{t_i}}+f_i(x,(t_1,\dots,t_k))))^{-1}
$$
for some admissible function $f_i$. The right hand side is equal to
$$
\frac{t_i}{1 + t_i\log (1 + e^{-S_i}f_i(x,(t_1,\dots,t_k)))}.
$$
Statement (2) easily follows from
this formula.\footnote{Recall $T_i = e^{S_i}$. Therefore a function which decays in exponential
order in $T_i$ coordinates also decays in exponential order in $S_i$ coordinates.}
\end{proof}
\begin{rem}\label{rem3110}
The admissible function $f_i$ above may not be zero at $t_i=0$.
However, $\log (1 + e^{-S_i}f_i(x,(t_1,\dots,t_k)))$ goes to
$0$ in exponential order as $S_i \to \infty$.
Therefore we can smoothly
extend $\varphi_{21}$ to a {\it collared} neighborhood.
Here is the key point that
coordinate changes of Kuranishi structure
can be extended smoothly to a collared neighborhood,
once we establish the exponential decay estimate
\eqref{eq:CCexpS} of the {\it coordinate changes} with respect to
the {\it $S_i$ coordinates}.
See also Remark \ref{rem:3140}.
On the other hand, we recall from Definition \ref{defn31333} that
admissible {\it functions} are required
that their derivatives\footnote{We assume $m>0$ in Definition \ref{defn31333},
while we assume $m\ge 0$ in \eqref{eq:CCexpS}.} with respect to
the {\it $T_i$ coordinates}
satisfy the
exponential decay estimate.
\end{rem}
Next, we go to the case of orbifold with corner.
\begin{defn}\label{def:3112}
\begin{enumerate}
\item
Two admissible charts are said to be {\it compatible as admissible charts}
if the diffeomorphism $\tilde\varphi$ in Definition \ref{2661} (3)
is admissible.
\item
A representative of an orbifold structure (with boundary or corner)
is said to be a {\it representative of an admissible orbifold structure}
if
\begin{enumerate}
\item
Each member is an admissible chart.
\item
Two of them are compatible as admissible charts.
\end{enumerate}
\item
In Definition \ref{def262220}, suppose $\{(V^X_i,\Gamma^X_i,\phi^X_i) \mid i \in I\}$
and $\{(V^Y_j,\Gamma^Y_j,\phi^Y_j) \mid j \in J\}$
are representatives of admissible orbifold structures.
Then the embedding $f$ in Definition \ref{def262220} (2) is said to be an {\it admissible embedding}\index{admissible ! embedding}
if $\tilde{f}_{p;ji}$ in Definition \ref{def262220} (2) is an admissible embedding
in the sense of Definition \ref{defn297} (1).
\end{enumerate}
\end{defn}
\begin{lem}\label{lem26444rev}
\begin{enumerate}
\item
Composition of admissible embeddings is an admissible embedding.
\item
The identity map is an admissible embedding.
\item
If an admissible embedding is a homeomorphism,
the inverse is also an admissible embedding.
\end{enumerate}
\end{lem}
The proof is obvious.
\begin{defn}\label{def:3113}
\begin{enumerate}
\item We say an admissible embedding of orbifolds is an
{\it admissible diffeomorphism}\index{admissible ! diffeomorphism}
if it is a
homeomorphism in addition.
\item
We say that two representatives of admissible orbifold structures on $X$ are
{\it equivalent} if the identity map regarded as a map between
$X$ equipped with those two representatives of admissible orbifold structures
is an admissible diffeomorphism.
This is an equivalence relation by Lemma \ref{lem26444rev}.
\item
An equivalence class of the equivalence relation (2) is called
an {\it admissible orbifold structure}
\index{admissible ! orbifold structure}
\index{orbifold ! admissible orbifold structure} on $X$.
\item
An orbifold $X$ with an admissible orbifold structure is called
an {\it admissible orbifold}.
\index{admissible ! orbifold}
\index{orbifold ! admissible}
\item
The condition for a map $X \to Y$ to be an admissible embedding
does not change if we replace representatives of admissible orbifold structures
to equivalent ones. So we can define the notion of
an {\it admissible embedding of admissible orbifolds}.
\item
If $U$ is an open subset of an admissible orbifold $X$, then
there exists a unique admissible orbifold structure on $U$ such that
the inclusion $U \to X$ is an admissible embedding.
We call $U$ with this admissible orbifold structure an {\it open admissible suborbifold}.
\end{enumerate}
\end{defn}
\begin{defn}\label{defn265rev5}
\begin{enumerate}
\item
Let $X$ be an admissible orbifold. An admissible orbifold chart $(V,\Gamma,\phi)$
of underlying topological space $X$
is an {\it admissible orbifold chart of orbifold}
\index{orbifold ! admissible orbifold chart}
\index{admissible ! admissible orbifold chart}
$X$ if
the map $V/\Gamma \to X$ induced by $\phi$ is an
admissible embedding of orbifolds.
\item
Hereafter when $X$ is an admissible orbifold,
an admissible orbifold chart always means an admissible orbifold chart of orbifold in the sense of (1).
\item
In case when an admissible orbifold structure on $X$ is given,
a representative of its admissible orbifold structure is
called an {\it admissible orbifold atlas}.
\item
Two admissible orbifold charts  $(V_i,\Gamma_i,\phi_i)$
are said to be {\it isomorphic} if there exists a
group isomorphism $h : \Gamma_1 \to \Gamma_2$
and an $h$ equivariant admissible diffeomorphism $\varphi : V_1 \to V_2$
such that $\phi_2 \circ \varphi = \phi_1$.
The pair $(h,\varphi)$ is called an {\it admissible isomorphism} or
{\it admissible coordinate change} between two admissible orbifold charts.
\index{admissible ! admissible coordinate change of admissible orbifold charts}
\index{admissible ! admissible isomorphism of admissible orbifold charts}
\end{enumerate}
\end{defn}
The proofs of the following lemmas are obvious from definition.
\begin{lem}\label{lem2810}
Suppose $V, V_1, V_2$ are as in Situation \ref{Situationcorner}.
\begin{enumerate}
\item A restriction of an admissible function on $V$ to its open subset is also admissible.
\item
Let $f : V_2 \to \R$ be an admissible function
and $\varphi_{21} : V_{21} \to V_1$ an admissible embedding,
then the composition $f \circ \varphi_{21} : V_1 \to \R$ is
admissible.
\end{enumerate}
\end{lem}
\begin{lem}
A subchart (in the sense of  Definition \ref{2661} (2)) of an admissible chart is
also admissible.
\end{lem}
\begin{defn}\label{def2810}
Let $X$ be an admissible orbifold.
A function $f : X \to \R$ is said to be
{\it admissible}\index{admissible ! function}
if for all the orbifold charts $(V,\Gamma,\phi)$ of $X$ the composition
$f\circ \phi : V \to \R$ is admissible in the sense of
Definition \ref{defn31333}.
\end{defn}
\begin{lem}
Let $X$ be an admissible orbifold and $f : X \to \R$ a function.
\begin{enumerate}
\item
If there exists a representative
$\{ (V_i, \Gamma_i, \phi_i) ~\vert~ i \in I \}$ of
the orbifold structure on $X$ such that
$f\circ \phi_i : V_i \to \R$ is admissible
for any $i$, then $f$ is admissible.
\item The composition of an admissible embedding and an admissible function
(resp. a function of exponential decay)
is again admissible (resp. a function of exponential decay).
\end{enumerate}
\end{lem}
This is a consequence of Lemma \ref{lem2810}.
\subsection{Admissible tensor}
\label{subsec:admtens}
Next we discuss tensor calculus etc.
on an admissible orbifold.
\begin{conven}\label{convNormalcoord}
Let $X$ be an admissible orbifold and let $V$ be as in Situation
\ref{Situationcorner}.
The coordinates of $V$ consist of those in $\overline V$ direction and
in ${\bf t}$ direction.
We call the coordinates
in ${\bf t}$ direction the
{\it normal coordinates}\index{normal coordinates}, and
the coordinates in $\overline V$ direction
the {\it horizontal coordinates}.
\end{conven}
\begin{defnlem}\label{def:admtensor}
\begin{enumerate}
\item
Let $\mathcal T$ be a tensor on an admissible orbifold and
$\mathcal T^{i_1,\dots,i_m}_{j_1,\dots,j_{m'}}$ be its expression
by admissible local coordinates.
We call it  {\it admissible}\index{admissible ! tensor}
if all $\mathcal T^{i_1,\dots,i_m}_{j_1,\dots,j_{m'}}$ are
admissible functions.
This notion is invariant under the
admissible coordinate change defined in Definition \ref{defn265rev5} (4).
\item
A tensor $\mathcal T$ is called {\it strongly admissible}\index{admissible ! strongly admissible tensor}
if the following holds.
\begin{enumerate}
\item
We rewrite $\mathcal T^{i_1,\dots,i_m}_{j_1,\dots,j_{m'}}$ using $T_i = e^{1/t_i}$
in place of $t_i$ as coordinates to obtain
$\widehat{\mathcal T}^{i_1,\dots,i_m}_{j_1,\dots,j_{m'}}$.
Then if $i_1,\dots,i_m$ and $j_1,\dots,j_m$ contain
the normal coordinates, $\widehat{\mathcal T}^{i_1,\dots,i_m}_{j_1,\dots,j_{m'}}$ is exponentially
small near the boundary.\footnote{In particular it is zero if the corresponding $t_i$ coordinate
is zero.}
\item
$\mathcal T$ is admissible.
\end{enumerate}
This notion is also invariant under the
admissible coordinate change.
\item
We define an {\it admissibile differential form}
as a special case of an admissible tensor.
\item
Various algebraic operations of tensors and differential forms
preserve the admissibility.
\item
The exterior derivative of an admissible differential form
is admissible.
Moreover the exterior derivative of an admissible function
is strongly admissible.
\item
The Lie derivative by an admissible vector field preserves admissibility.
The bracket of admissible vector fields
(resp. strongly admissible vector fieles) is
admissible (resp. strongly admissible).
\item
An {\it admissible Riemannian metric}
\index{admissible ! Riemannian metric} of an admissible
orbifold is a stratumwise Riemannian metric with the following properties.
(Here we consider the corner structure stratification.)
\par
Let $\overline V \times [0,1)^k$ be an
admissible chart. We write the $[0,1)^k$ coordinates as $t_i$
and put $T_i = e^{1/t_i}$.
Then there exists a symmetric 2-tensor $g'$
on $\overline V \times [0,1)^k$
with the following properties:
\begin{enumerate}
\item
$g'$ is strongly admissible as a 2-tensor.
\item
For $A \subseteq \{1,\dots,k\}$ we put
$$
V^A = \{(x,(t_1,\dots,t_k)) \mid x \in \overline V,
\,\, \text{$t_i = 0$ if $i \in A$}, \text{$t_i \ne 0$ if $i \notin A$} \},
$$
which is an open subset of a stratum.
We require that on $V^A$ the stratumwise Riemannian metric is given by
\begin{equation}\label{adRiemann}
g' + \sum_{i \notin A} (dT_i)^{2\otimes}.
\end{equation}
\end{enumerate}
In particular, this condition implies that the Riemannian
metric is stratumwise complete.
Note that we use the $T_i$ coordinates and not the $t_i$ coordinates in \eqref{adRiemann}.
We also note that $g'$ $C^{\infty}$-converges
to the stratumwise metric on $\overline{V}$ as $T_i \to \infty$.
\item
Let $\bigcup_{i\in I}U_i = X$ be an open covering of an orbifold $X$.
A smooth partition of unity  $\{\chi_i\mid i\in I\}$ subordinate to this covering is
called an {\it admissible partition of unity} if
each $\chi_i$ is an admissible function.
\end{enumerate}
\end{defnlem}
The proof is straightforward by definition, so is omitted.
\begin{lem} Let $V$be as in Situation \ref{Situationcorner}.
\begin{enumerate}
\item
Pull-back of an admissible differential form on $V$ by an admissible
embedding is admissible.
\item
Pull-back of a differential form on a smooth manifold $M$ by an admissible
map from $V$ to $M$ is admissible.
\item
Any locally finite covering of an admissible orbifold admits an
admissible partition of unity.
\item
An admissible Riemannian metric exists.
\end{enumerate}
\end{lem}
\begin{proof}
We can prove the lemma by modifying the standard proof of the corresponding results
in manifold theory in an obvious way, so omit it.
\end{proof}
\par
When we apply certain operations on tensors, admissibility
is mostly preserved under the operations.
In certain case if
we take derivative on the normal direction,
an admissible object changes to a strongly admissible
one.
Since it is easy to see when it happens, we do not
provide a detailed account thereon here.
We will state that kinds of facts in case we need it.

\subsection{Admissible vector bundle}
\label{subsec:admvecbund}
We next describe the admissible version of Definition \ref{defn2613}.
\begin{defn}\label{defn2613rev}
Let $(X,\mathcal E,\pi)$ be a pair of admissible orbifolds $X$ and $\mathcal E$
with a continuous map $\pi : \mathcal E \to X$
between their underlying topological spaces.
Hereafter we write $(X,\mathcal E)$ in place of $(X,\mathcal E,\pi)$.
\begin{enumerate}
\item
An {\it admissible orbifold chart} of $(X,\mathcal E)$ is
a quintuple $(V,E,\Gamma,\phi,\widehat\phi)$
with the following properties:
\begin{enumerate}
\item
$\frak V = (V,\Gamma,\phi)$ is an admissible orbifold chart of $X$.
\item $E$ is a finite dimensional vector space equipped with a linear $\Gamma$
action.
\item
$(V \times E,\Gamma,\widehat\phi)$ is an admissible orbifold chart of $\mathcal E$.
\item
The diagram below commutes set-theoretically.
\begin{equation}\label{diag26399}
\begin{CD}
V \times E @ >{\widehat\phi}>>
\mathcal E  \\
@ V{}VV @ VV{\pi}V\\
V @ > {\phi} >> X
\end{CD}
\end{equation}
Here the left vertical arrow is the projection to the
first factor.
\end{enumerate}
\item
In the situation of (1), let $p\in V$ and $( V_{p},\Gamma_p,\phi\vert_{V_p})$
be a subchart of $(V,\Gamma,\phi)$ in the sense of
Definition \ref{2661} (2).
Then  $(V_{p},E,\Gamma_p,\phi\vert_{V_p},\widehat\phi\vert_{V_p \times E})$
is an admissible orbifold chart of $(X,\mathcal E)$.
We call it a {\it subchart} of $(V,E,\Gamma,\phi,\widehat\phi)$.
\item
Let $(V_i,E_i,\Gamma_i,\phi_i,\widehat\phi_i)$  $(i=1,2)$ be
admissible orbifold charts of $(X,\mathcal E)$.
We say that they are {\it compatible as admissible charts} if the following holds
for each $p_1 \in V_1$ and $p_2 \in V_2$ with
$\phi_1(p_1) = \phi_2(p_2)$.
\begin{enumerate}
\item
There exists an isomorphism
$(h,\varphi) : (V_1,\Gamma_1,\phi_1) \to (V_2,\Gamma_2,\phi_2)$
of admissible orbifold charts of $X$.
\item
There exists an isomorphism
$(h,\widehat\varphi) : (V_1\times E_1,\Gamma_1,\widehat\phi_1) \to
(V_2\times E_2,\Gamma_2,\widehat\phi_2)$
of admissible orbifold charts of $\mathcal E$.
\item
For each $y \in V_1$ the map $ E_1 \to E_2$ given by
$\xi \to \pi_{E_2}\widehat{\varphi}(y,\xi)$ is a linear isomorphism.
Here $\pi_{E_2} : V_2 \times E_2 \to E_2$ is the projection.
\item
Each component of the map $V_1\times E_1 \to E_2$ that is a composition of
$\widehat\varphi$ and the projection $V_2 \times E_2
\to E_2$ is an admissible function.
\end{enumerate}
\item
A {\it representative of an admissible vector bundle structure} on $(X,\mathcal E)$
is a set of admissible orbifold charts $\{(V_i,E_i,\Gamma_i,\phi_i,\widehat\phi_i) \mid i \in I\}$
such that any two of the charts are compatible in the sense of (3)
above and
$$
\bigcup_{i\in I} \phi_i(V_i) = X,
\quad
\bigcup_{i\in I} \widehat\phi_i(V_i \times E_i) = \mathcal E,
$$
are locally finite open coverings.
\end{enumerate}
\end{defn}
\begin{defn}\label{def26222rev}
Let $(X^*,\mathcal E^*)$ $(* = a,b)$ have representatives
of vector bundle structures $\{(V^*_i,E^*_i,\Gamma^*_i,\phi^*_i,\widehat\phi^*_i) \mid i \in I^*\}$,
respectively.
A pair of admissible orbifold embeddings $(f,\widehat f)$,
$f : X^a \to X^b$, $\widehat f : \mathcal E^a \to \mathcal E^b$ is said to be an
{\it admissible embedding of admissible vector bundles} if the following holds.
\begin{enumerate}
\item
Let
$p \in V^a_i$, $q \in V^b_j$ with
$f(\phi^a_i(p)) = \phi^b_j(q)$.
Then there exist admissible open subcharts
$(V^a_{i,p}\times E^a_{i,p},\Gamma^a_{i,p},\widehat\phi^a_{i,p})$
and
$(V^b_{j,q}\times E^b_{j,q},\Gamma^b_{j,q}\widehat\phi^b_{j,q})$
and a local representative
$(h_{p;ji},f_{p;ji},\widehat f_{p;ji})$ of the embeddings $f$ and $\widehat f$
such that the following holds.
For each $y \in V^a_i$ the map
$\xi \mapsto \pi_2(\widehat f_{p;ji}(y,\xi))$,
$E^a_1 \to E^b_2$ is a linear embedding.
\item
Each component of the map $\pi_2\circ \widehat f_{p;ji}:
V_1^a \times E_1^a  \to E_2^b$ is admissible.
\item
The diagram below commutes set-theoretically.
\begin{equation}\label{diag2633rev}
\begin{CD}
\mathcal E^a @ >{\widehat f}>>
\mathcal E^b  \\
@ V{\pi}VV @ VV{\pi}V\\
X^a @ > {f} >> X^b
\end{CD}
\end{equation}
\end{enumerate}
Two orbifold embeddings are said to be {\it equal} if they coincide set-theoretically
as pairs of maps.
\end{defn}
\begin{lem}\label{lem26444AA}
\begin{enumerate}
\item
The composition of admissible embeddings of vector bundles is an
admissible embedding.
\item
The pair of identity maps $({\rm id}, \widehat{\rm id})$ is an admissible embedding.
\item
If an admissible embedding of vector bundles is a pair of homeomorphisms,
then the pair of their inverses is also an admissible embedding.
\end{enumerate}
\end{lem}
The proof is obvious.
\begin{defn}\label{def:3125}
Let
$(X,\mathcal E)$ be as in Definition \ref{defn2613}.
\begin{enumerate}
\item
An admissible embedding of vector bundles is said to be an {\it isomorphism}
\index{admissible ! isomorphism of admissible vector bundles}
if it is a pair of admissible diffeomorphisms of admissible orbifolds.
\item
We say that two representatives of an admissible vector bundle structure on $(X,\mathcal E)$ are
{\it equivalent} if the pair of identity maps regarded as a map between
$(X,\mathcal E)$ equipped with those two representatives of vector bundle structures
is an admissible embedding. This is an equivalence relation by Lemma \ref{lem26444AA}.
\item
An equivalence class of the equivalence relation (1) is called
an {\it admissible vector bundle structure} on $(X,\mathcal E)$.
\item
A pair $(X,\mathcal E)$ together with its admissible vector bundle
structure is called an
{\it admissible vector bundle}\index{admissible ! vector bundle}
on $X$.
We call $\mathcal E$ the {\it total space}, $X$ the
{\it base space}, and $\pi : \mathcal E \to X$ the {\it projection}.
\item
The condition for $(f,\widehat f) : (X^a,\mathcal E^a) \to (X^b,\mathcal E^b)$
to be an admissible embedding
does not change if we replace representatives of admissible vector bundle
structures
to equivalent ones. So we can define the notion of
an {\it admissible embedding of vector bundles}.
\index{admissible ! embedding of vector bundles}
\item
We say $(f,\widehat f)$ is an admissible embedding {\it over the admissible orbifold embedding $f$.}
\end{enumerate}
\end{defn}
\begin{defn}\label{defnAdm2655}
\begin{enumerate}
\item
Let $(X,\mathcal E)$ be an admissible vector bundle.
We call an admissible orbifold chart $(V,E,\Gamma,\phi,\widehat\phi)$
in the sense of Definition \ref{defn2613rev} (1)
of underlying pair of topological spaces $(X,\mathcal E)$
an {\it admissible orbifold chart of an admissible vector bundle} $(X,\mathcal E)$
\index{admissible ! admissible orbifold chart of admissible vector bundle}  if
the pair of maps $(\overline\phi,\overline{\widehat\phi}) : (V/\Gamma,(V\times E)/\Gamma) \to (X,\mathcal E)$ induced from $(\phi,\widehat\phi)$ is an
admissible embedding of admissible vector bundles.
\item
If $(V,E,\Gamma,\phi,\widehat\phi)$
is an admissible orbifold chart of an admissible vector bundle,
we call a pair $(E,\widehat\phi)$
a {\it trivialization}\index{orbifold ! trivialization of vector bundle}
\index{admissible ! trivialization of admissible vector bundle}
of our admissible vector bundle on $V/\Gamma$.
\item
Hereafter when $(X,\mathcal E)$ is an admissible vector bundle,
its `admissible orbifold chart' always means an admissible orbifold chart
of an admissible vector bundle in the sense of (1).
\item
In case when an admissible vector bundle structure on $(X,\mathcal E)$ is given,
a representative of this admissible vector bundle structure is
called an {\it admissible orbifold atlas} of $(X,\mathcal E)$.
\index{orbifold ! admissible orbifold atlas of admissible vector bundle}
\index{admissible ! admissible orbifold atlas of admissible vector bundle}
\item
Two admissible orbifold charts $(V_i,E_i,\Gamma_i,\phi_i,\widehat \phi_i)$
of an admissible vector bundle
are said to be {\it isomorphic} if there exist an
isomorphism $(h,\tilde\varphi)$ of admissible orbifold charts
$(V_1,\Gamma_1,\phi_1) \to (V_2,\Gamma_2,\phi_2)$
and an admissible isomorphism
$(h,\tilde{\hat{\varphi}})$ of admissible orbifold charts
$(V_1\times E_1,\Gamma_1,\widehat\phi_1) \to (V_2\times E_2,\Gamma_2,\widehat\phi_2)$
such that they induce an admissible embedding of admissible vector bundles
$(\varphi,\hat\varphi) : (V_1/\Gamma_1,(V_1\times E_1)/\Gamma_1) \to (V_2/\Gamma_2,(V_2\times E_2)/\Gamma_2)$.
The triple
$(h,\tilde\varphi,\tilde{\hat{\varphi}})$ is called an {\it admissible isomorphism} or
{\it admissible coordinate change} between
admissible orbifold charts of the admissible vector bundle.
\index{admissible ! admissible coordinate change of admissible vector bundle}
\index{admissible ! admissible isomorphism of admissible orbifold charts}
\end{enumerate}
\end{defn}
Once we have established these basic notions related to admissible vector bundles
as above, the next lemma obviously holds.
\begin{lem}
The tangent bundle of an admissible orbifold has a canonical structure
of admissible vector bundle.
\par
Taking Whitney sum, tensor product, dual, exterior power, quotient of
admissible vector bundles preserve admissibility.
\end{lem}
\begin{lemdef}
If $(X^b,\mathcal E^b)$ is an admissible vector bundle and
$f : X^a \to X^b$ is an admissible embedding,
then the pull-back $f^*\mathcal E^b$ has a unique structure
of admissible vector bundle such that the embedding
of vector bundles $(f,\widehat f) : (X^a,f^*\mathcal E^b)
\to (X^b,\mathcal E^b)$ becomes an admissible embedding.
\par
We call $f^*\mathcal E^b$ equipped
with this admissible vector bundle structure the
{\rm pull-back in the sense of admissible vector bundles}.
\index{admissible ! pull-back of admissible vector bundle}
\index{pull-back ! of admissible vector bundle}
\end{lemdef}
The proof is straightforward, so is omitted.
The following lemmas are also straightforward consequences from definitions.
\begin{lem}
A covering space $\tilde X$ of an admissible orbifold $X$ has a
canonical structure of an admissible orbifold such that
admissible functions of $X$ are pulled back to admissible
functions.
\end{lem}
\begin{lem}
The normalized boundary or corner of an admissible orbifold
is admissible.
The covering maps in Lemma \ref{lem2712} or
Proposition \ref{prop2813} are admissible.
\end{lem}
\begin{defn}\label{def:3128}
An {\it admissible
section}\index{admissible ! section}
of an admissible vector bundle $(X,\mathcal E)$ is an admissible embedding of orbifolds $s : X \to \mathcal E$
such that the composition of $s$ and the projection
is the identity map set-theoretically.
\end{defn}
The next lemma obviously follows from definition.
\begin{lem}
An admissible tensor of an admissible orbifold
is regarded as an admissible section of an appropriate tensor product
bundle of the tangent bundle or its dual.
\end{lem}
\begin{defnlem}\label{defnlem:admconne}
A connection $\nabla$ on an admissible vector bundle $(X,\mathcal E)$ is called
{\it strongly admissible}\index{admissible ! strongly admissible connection}
if for an admissible orbifold chart $(V,E,\Gamma,\phi,\widehat\phi)$
of $(X,\mathcal E)$
with $V\subset \overline{V}\times [0,1)^k$ as in Situation \ref{Situationcorner},
it is locally expressed by a 1-form
$\sum_{\alpha=1}^{\overline{v}} A^a_{b,\alpha}dx_{\alpha} +
\sum_{i=1}^{k} A^a_{b,\overline{v}+i}dT_{i}$
satisfying the following:
 (Here $x_1, \dots, x_{\overline{v}}$ are the
local coordinates of $\overline{V}$ and $A^a_{b, \ast}$ is the $(a,b)$-component of a matrix $(A^a_{b,\ast})$ defined by an endmorphism of $E$, and
$T_i=e^{1/t_i}$ $(i=1,\dots,k)$.)
\begin{enumerate}
\item the function $A^a_{b,\alpha}$ is admissible, and
\item the function $A^a_{b,\overline{v}+i}$ is exponentially small near the boundary.
\end{enumerate}
This notion is independent of
the choices of the admissible orbifold charts of
$(X,\mathcal E)$ up to admissible coordinate changes in the sense of Definition \ref{defnAdm2655} (5).
(Note that the exterior derivative of an admissible function is strongly admissible
by Definition-Lemma \ref{def:admtensor} (5).)
\end{defnlem}
The next lemma is used in the proof of Lemma \ref{admissiblenormalbundle}.
\begin{lem}\label{lem:admiLC}
The Levi-Civita connection of an admissible Riemannian metric is strongly admissible.
Moreover, the Christoffel symbol $\Gamma^{k}_{ij}$ of the Levi-Civita connection
enjoys the following property:
If
$i,j,k$ contain the horizontal coordinates (see
Convention \ref{convNormalcoord}), $\Gamma^{k}_{ij}$ is admissible,
and if $i,j,k$ contain the normal coordinates, $\Gamma^{k}_{ij}$ is exponentially small
near the boundary
(under the coordinates $T_{\ast}=e^{1/t_{\ast}}$).
\end{lem}
\par
The next is an analog of Proposition \ref{homotopicpulback}
in admissible vector bundles.
\begin{prop}\label{homotopicpulbackAdm}
Let $\mathcal E$ be an admissible vector bundle on $X\times [0,1]$, where
$X$ is an admissible orbifold.
We identify $X \times \{0\}$, $X \times \{1\}$
with $X$ in an obvious way. Then
there exists an isomorphism of admissible vector bundles
$$
I : \mathcal E\vert_{X \times \{0\}} \cong \mathcal E\vert_{X \times \{1\}}.
$$
Suppose in addition that
we are given a compact set $K \subset X$,
its neighborhood $V$ and an isomorphism
$$
I_0 : \mathcal E\vert_{V \times [0,1]} \cong \mathcal E\vert_{V \times \{0\}} \times [0,1].
$$
Then we may choose $I$ so that it coincides with the isomorphism induced by
$I_0$ on $K$.
If $K$ is a submanifold, we may take $K=V$.
\end{prop}
\begin{proof}
Using an admissible connection defined in Definition-Lemma \ref{defnlem:admconne},
we can define a parallel transform for admissible vector bundles.
Then the proof goes in a way similar to that of Proposition \ref{homotopicpulback}.
See also the proof of Lemma \ref{admissibletrivializationofbundle} below.
\end{proof}
In this way, we can translate various stories for manifolds to the admissible
realm.
Especially we can define the notion of admissibility of Kuranishi structure.
In fact, the whole story works just by adding the word admissible
to various constructions.
\begin{defn}\label{def:adKura}
A Kuranishi structure
$$
\widehat{\mathcal U} =
\left(\left\{ \mathcal U_p = (U_p, \mathcal E_p, \psi_p, s_p) \right\},
\left\{\Phi_{pq} = (U_{pq}, \varphi_{pq}, \widehat{\varphi}_{pq}) \right\}
\right)
$$
is {\it admissible}\index{Kuranishi structure ! admissible}
\index{admissible ! Kuranishi structure}
if
\begin{enumerate}
\item[$\bullet$] $U_p$ is an admissible orbifold in the sense of
Definition \ref{def:3113}.
\item[$\bullet$] $\mathcal E_p$ is an admissible vector bundle in the sense of
Definition \ref{def:3125}.
\item[$\bullet$] $s_p$ is an admissible section of $\mathcal E_p$ in the sense of
Definition \ref{def:3128}.
\item[$\bullet$] $U_{pq}$ is an open admissible suborbifold of $U_q$ in the sense of Definition \ref{def:3113}.
\item[$\bullet$] $(\varphi_{pq}, \widehat{\varphi}_{pq})$ is an admissible
embedding of admissible vector bundles in the sense of
Definition \ref{def:3125}.
\end{enumerate}
\end{defn}

\subsection{Admissibility of bundle extension data}
\label{subset:adbedata}
In \cite[Subsections 12.2 and 13.3]{part11} we used
a pair of the projection of the normal bundle and
an isomorphism of bundles, which we called a bundle extension data.
In this subsection we introduce the corresponding notion in the
admissible category and prove its existence.
\begin{defn}
In \cite[Situation 12.21]{part11}, we assume
that Kuranishi charts $\mathcal U_i = (U_i, \mathcal E_i, \psi_i, s_i)$ ($i=1,2$) are admissible and the embedding
$\Phi_{21} : \mathcal U_1 \to \mathcal U_2$ is admissible.
Then we say the bundle extension data
$(\pi_{12},\tilde\varphi_{21},\Omega_{12})$ associated with $\Phi_{21}$
defined as in Definition \ref{coorchangebed}
is an {\it admissible bundle extension data}
\index{bundle extension data ! admissible bundle extension data}
\index{admissible ! bundle extension data}
if both of $\pi_{12}:\Omega_{12} \to U_1$ and
$\tilde\varphi_{21}: \pi_{12}^{\ast}\mathcal E_1 \to \mathcal E_2$ are admissible.
\end{defn}
\begin{prop}\label{admbdlextdat}
An admissible bundle extension data exists.
\end{prop}
\begin{proof}
We first construct a tubular neighborhood in the admissible
category.
\par
Let $f : X \to Y$ be an admissible embedding of admissible orbifolds.
Its normal bundle $N_XY$ is defined and is an admissible vector bundle
over $X$. We take an admissible Riemannian metric on $Y$.
We will use it to define an exponential map
$$
{\rm Exp} : BN_XY \to Y
$$
that is a map from the neighborhood of the zero section
of the total space of the normal bundle $N_XY$ to $Y$.

\begin{lem}\label{admissiblenormalbundle}
The exponential map ${\rm Exp}$ is an admissible diffeomorphism
between admissible orbifolds.
\end{lem}
\begin{proof}
Let $Y$ be an admissible suborbifold.
Using local coordinates we can represent them as follows.
Let $p \in X \subset Y$ and $\frak V_{X,p} = (V_{X,p}, \Gamma_{X,p}, \phi_{X,p})$
(resp. $\frak V_{Y,p}= (V_{Y,p}, \Gamma_{Y,p}, \phi_{Y,p})$)
be an orbifold chart
of $X$ (resp. $Y$) at $p$.
We may take $V_{X,p} = \overline V_{X,p} \times [0,1)^k$
and $V_{Y,p} = \overline V_{Y,p} \times [0,1)^k$.
Here  $\overline V_{X,p}$, $\overline V_{Y,p}$ are open subsets of Euclidean space.
\par
We change the coordinates from $(t_1,\dots,t_k) \in [0,1)^k$
to $(T_1,\dots,T_k) \in (1,\infty]$ by
$$
T_i = e^{1/t_i}.
$$
Then the embedding $X \to Y$ in the coordinates is
given by
$$
\varphi :  \overline V_{X,p} \times (0,\infty]^k
\to  \overline V_{Y,p} \times (0,\infty]^k
$$
such that
$$
\varphi(\overline y,(T_1,\dots,T_k))
=
\left(\varphi_0(\overline y,(T_1,\dots,T_k)),~
(T_i+\varphi_i(\overline y,(T_1,\dots,T_k))_{i=1}^k)\right),
$$
where $\varphi_0$, $\varphi_i$ are admissible.
Therefore the normal bundle has an orthonormal frame of the form
\begin{equation}\label{normalframeeq}
e^a = (e^a_0,(e^a_1,\dots,e^a_k)),
\end{equation}
where $e^a_0$ is admissible and $e^a_i$ are exponentially small at the boundary.
(We also note that we use an admissible Riemannian metric
of the form (\ref{adRiemann}).)
\par
Now let $\ell(s)$ be a geodesic of $Y$ such that $\ell(0) \in\varphi(X)$
and
$$
\overset{\cdot}\ell(0)
= \sum_{a=1}^{\dim X - \dim Y} c_ae^a.
$$
Let us write down the equation that $\ell$ is a geodesic by local coordinates.
We use Greek letters for the indices of the horizontal coordinates
and $h,i,j$ for the indices of the normal coordinates in the sense of Convention \ref{convNormalcoord}.
(Recall that we use $T_i = e^{1/t_i}$ for the normal coordinates.)
Then the equation of the geodesic is:
\begin{equation}\label{geodesiceq}
\aligned
&\frac{d^2\ell_{\alpha}}{ds^2}
+ \sum_{i,j}\Gamma^{\alpha}_{ij}\frac{d\ell_{i}}{ds}\frac{d\ell_{j}}{ds}
+ \sum_{i,\beta}\Gamma^{\alpha}_{i\beta}\frac{d\ell_{i}}{ds}\frac{d\ell_{\beta}}{ds}
+ \sum_{\beta,\gamma}\Gamma^{\alpha}_{\beta\gamma}\frac{d\ell_{\beta}}{ds}\frac{d\ell_{\gamma}}{ds} = 0,
\\
&\frac{d^2\ell_{h}}{ds^2}
+ \sum_{i,j}\Gamma^{h}_{ij}\frac{d\ell_{i}}{ds}\frac{d\ell_{j}}{ds}
+ \sum_{i,\beta}\Gamma^{h}_{i\beta}\frac{d\ell_{i}}{ds}\frac{d\ell_{\beta}}{ds}
+ \sum_{\beta,\gamma}\Gamma^{h}_{\beta\gamma}\frac{d\ell_{\beta}}{ds}\frac{d\ell_{\gamma}}{ds} = 0.
\endaligned
\end{equation}
Since $\overset{\cdot}\ell(0) \perp T\varphi(X)$, the functions
$
\frac{d\ell_{\alpha}}{ds}(0)
$
are admissible and
$
\frac{d\ell_{h}}{ds}(0)
$
are exponentially small.
We also note that the Christoffel symbols
$\Gamma^{*}_{**}$ are admissible
and $\Gamma^{*}_{i*}$, $\Gamma^{h}_{**}$ are exponentially small.
(See Lemma \ref{lem:admiLC}.)
Therefore the equation (\ref{geodesiceq})
implies that
the map which associates
$\ell(s)$ to the initial value $(\ell(0),(c_1,\dots,c_a))$
is admissible. The proof of Lemma \ref{admissiblenormalbundle} is complete.
\end{proof}
\begin{rem}
The above proof shows that the geodesic $s \mapsto \ell(s)$ such that
$\ell(0) \in\varphi(X)$
and $\overset{\cdot}\ell(0) \perp T\varphi(X)$, $\Vert
\overset{\cdot}\ell(0)\Vert =1$
is defined on $s \in [-S,S]$, where $S$ is independent of the
initial condition
$\overset{\cdot}\ell(0)$ and $\ell(0)$.
\end{rem}

In the situation of Lemma \ref{admissiblenormalbundle}
let $\mathcal E$ be an admissible vector bundle on $X$.
Its pull-back $f^*\mathcal E$ is an admissible vector
bundle on $Y$.
We identity a neighborhood of $f(X)$ in $Y$ with
an open subset $BN_XY$of $N_XY$ that is a neighborhood of the origin.
The pull-back $\pi^*f^*\mathcal E$ is an
admissible vector bundle.
We define a map
$$
\tilde f : \pi^*f^*\mathcal E \to \mathcal E
$$ as follows.
We describe the construction of $\tilde f$ for the case when $X,Y$ are manifolds.
For the case of orbifolds we construct the map locally in a way equivariant under the action
of the isotropy group.
Then we can globalize the construction using
the invariance under the coordinate change.
\par
Let $p \in BN_XY$. The fiber $\pi^*f^*\mathcal E_p$
is identified with $\mathcal E_{f(p)}$.
We assume that our neighborhood $BN_XY$ of $X$ is sufficiently
small. Then there exists a unique geodesic joining
$p=f(q)$ and $q$.
Taking an admissible connection on $\mathcal E$,
we use the parallel transport along
this geodesic to send elements of $\mathcal E_{f(p)}$
to $\mathcal E_{p}$.
Thus we have obtained the map
$
\tilde f : \pi^*f^*\mathcal E \to \mathcal E
$.
\begin{lem}\label{admissibletrivializationofbundle}
The map $\tilde f$ is an admissible isomorphism of admissible vector bundles.
\end{lem}
\begin{proof}
Let $\ell(s)$ be the geodesic joining $p$ to $q$.
Note $\ell(0) \in\varphi(X)$
and $\overset{\cdot}\ell(0) \perp T\varphi(X)$.
Let $\xi(s) \in \mathcal E_{\ell(s)}$ is a parallel section
along $\ell$. We take an admissible orbifold chart
of $(X,\mathcal E)$ and write  $\xi_a$ be the expression of $\xi$.
We use the letter $a,b$ for the indices of $\xi$.
The equation that $\xi(s)$ is parallel is written as
$$
\frac{d\xi_a}{ds}
+
\sum_{b,i} A^a_{b,i} \frac{d\ell_i}{ds}\xi_b
+
\sum_{b,\alpha} A^a_{b,\alpha} \frac{d\ell_\alpha}{ds}\xi_b
=0.
$$
Here $A^*_{**}$ is an admissible connection on our admissible vector bundle.
Note that $A^a_{b,i}$ and $\frac{d\ell_i}{ds}$ are exponentially small
near the boundary and
$ A^a_{b,\alpha}$, $\frac{d\ell_\alpha}{ds}$ are admissible.
Therefore the map which associates $\xi_a(s)$ to $\xi_a(0)$
is admissible.
\end{proof}
We put $\pi = {\rm pr} \circ {\rm Exp}^{-1} : Y \to X$, where
$ {\rm Exp}$ is as in Lemma \ref{admissiblenormalbundle} and
$ {\rm pr} : BN_XY \to X$ is the projection of the normal bundle.
The pair $(\pi,\tilde f)$ where $\tilde f$ is as in Lemma
\ref{admissibletrivializationofbundle}
is an admissible bundle extension data.
\end{proof}
\begin{prop}\label{prop12333rev}
In the situation of \cite[Proposition 13.9]{part11}
we assume that $(\widetriangle{\mathcal U},\mathcal K)$
is admissible.
Then there
exists a system of bundle extension data of $(\widetriangle{\mathcal U},\mathcal K)$
consisting of admissible bundle extension data.
\end{prop}
\begin{proof}
The proof is the same as that of \cite[Proposition 13.9]{part11},
using Proposition \ref{admbdlextdat}.
In the course of the proof we need to extend the given
admissible extension data without changing it at the compact set
where the bundle extension data is already defined.
(See the last part of the proof of \cite[Lemma 13.10]{part11}.)
We can do the procedure in the admissible category using an admissible
partition of unity.
\end{proof}

\subsection{Admissibility of the
moduli spaces of pseudo-holomorphic curves}
\label{subset:admphc}

In \cite[Theorem 6.4, Proposition 8.27, Proposition 8.32]{foooexp} (see also
\cite[Proposition 16.11]{foootech}), we proved the exponential
decay of the coordinate change. It implies
\begin{equation}\label{formula308}
\left\Vert \frac{\partial}{\partial T_{\rm e}'}(T_{\rm e}' - T_{\rm e}) \right\Vert
\le C e^{-\delta T_{\rm e}'}.
\end{equation}
Here $T_{\rm e}$ is the coordinate corresponding to a singular point
which is
resolved and $T'_{\rm e}$ is the  coordinate corresponding to the same
singular point after coordinate change.
\par
In this subsection, we axiomatize the properties of the coordinate change
proved in \cite{foootech} \cite{foooexp}
under an abstract setting described in Situation \ref{situ3136}.
Then we can directly see that the results proved in \cite{foootech} \cite{foooexp} actually
imply that the Kuranishi structure of the moduli space of pseudo-holomorphic curves
is admissible in the sense of Definition \ref{def:adKura}.
This is tautological,
but such an axiomatization will be useful when we prove admissibility of
coordinate change of a Kuranishi structure constructed in other situation.
\begin{shitu}\label{situ3136}
We consider open subsets $V_1 \subset \overline V_1 \times (0,\infty]^k$ and
$V_2 \subset \overline V_2 \times (0,\infty]^k$,
where $\overline V_i$ are open subsets of $\R^{n_i}$.
Let  $\varphi : V_1 \to V_2$ be an embedding of topological
spaces such that the following holds:
We consider the stratification of $V_i$
such that $S_mV_i$ consists of the points where at least $m$
of the coordinates of the $(0,\infty]^k$ factor are $\infty$.
We put $\overset{\circ}S_mV_i = S_mV_i \setminus S_{m+1}V_i$.
We assume the following:
\begin{enumerate}
\item $\varphi(p) \in S_mV_2$ if and only if $p \in S_mV_1$.
\item  The restriction of $\varphi$ to $\overset{\circ}S_mV_1$
is a smooth embedding $\overset{\circ}S_mV_1 \to
\overset{\circ}S_mV_2$.
\item
We write
$$
\varphi(x;T_1,\dots,T_k)
= (\overline{\varphi}(x;T_1,\dots,T_k);
T'_1(x;T_1,\dots,T_k),\dots,T'_k(x;T_1,\dots,T_k)).
$$
Then the following holds:
\begin{enumerate}
\item
$$
\left\Vert \frac{\partial \overline{\varphi}}{\partial T_i}\right\Vert_{C^k}
\le C_k e^{-c_k T_i}.
$$
Here $C_k, c_k$ are positive numbers depending on $k$ and $\varphi$,
and $\Vert \cdot \Vert_{C^k}$ is the $C^k$ norm with respect to all
of $x, T_i$.
\item
\begin{equation}\label{expTTTT}
\left\Vert \frac{\partial T'_j}{\partial T_i}- \delta_{ij}\right\Vert_{C^k}
\le C_k e^{-c_k T_i}.
\end{equation}
Here $C_k,c_k,\Vert \cdot\Vert_{C^k}$ are the same as above and
$\delta_{ij} = 0$ if $i\ne j$ and $\delta_{ii} = 1$.
\end{enumerate}
We assume the above inequality (a)(b) stratumwise.
Namely in case when certain coordinate $T_i$ is $\infty$, we
require the inequality only for $T_j$ derivatives which are not $\infty$.
$\blacksquare$
\end{enumerate}
\end{shitu}

\begin{lem}\label{3137}
In Situation \ref{situ3136}
the coordinate change $\varphi$ is admissible in the sense of
Definition \ref{defn297}.
\end{lem}
\begin{proof}
This is immediate from the definition.
\end{proof}
We can use this lemma to obtain an admissible coordinate
in the geometric situation appearing in the moduli space of
pseudo-holomorphic curves.
\begin{rem}\label{rem:3140}
As we explained in \cite[Remark A1.63]{fooobook2},
the coordinate appearing in algebraic geometry
is $e^{-cT_{\rm e}}$ which decays faster than $1/T_{\rm e}$.
On the other hand,
$1/T_{\rm e}$ is the coordinate used in \cite{fooo09} \cite{foootech} etc.
Here the coordinate change is smooth with respect to this coordinate $1/T_{\rm e}$.
In this article we take an even more slower coordinate $1/\log T_{\rm e}$
than $1/T_{\rm e}$
(see (\ref{Tandt}))
so that the coordinate change is admissible.
\par
When we use the coordinate $s_{\rm e} = 1/T_{\rm e}$
in place of $1/\log T_{\rm e}$,
Lemma \ref{newlem319} (1) still holds while Lemma \ref{newlem319}
(2) does not.
In fact,
$
f(x,T) = 1
$
is an admissible coordinate and so
\begin{equation}\label{312312}
T' = T+1
\end{equation}
is an admissible coordinate change.
Then putting $s' = 1/T'$ and $s = 1/T$, we have
$
s' = \frac{1}{1+1/s} = \frac{s}{s+1}.
$
Hence
$
\frac{d^2s'}{ds^2}(0)  = -2 \ne 0.
$
As we mentioned in Remark \ref{rem3110},
Lemma \ref{newlem319} (2) is necessary to extend the coordinate change
to the collared neighborhood.
Indeed, in Section \ref{sec:triboundary} of this article
we extended the coordinate change
from $\overline V \times [0,1)$ to
$\overline V \times (-1,1)$ by taking the $(-1,1)$
component to be the identity map on $(-1,0)$.
Lemma \ref{newlem319}
(2) implies that this extended coordinate change is smooth
at $t = 1/\log T = 0$.
\par
On the other hand,
in \cite{fooobook2} we used the fact that the coordinate change
is smooth in the coordinate $s = 1/T$.
However, we note that we did not need Lemma \ref{newlem319} (2)
in \cite{fooobook2} because we did not use
the {\it collared} Kuranishi neighborhood
and did not need to extend the coordinate change.
\par
The fact that the extended coordinate change is smooth
will be used when we make Kuranishi structures
on the moduli spaces of pseudo-holomorphic curves
compatible with the
{\it forgetful map of marked points}\index{forgetful map of marked points},
by adopting the framework of {\it collared} Kuranishi structure developed in
Section \ref{sec:triboundary}.
We will discuss this point in detail in the forthcoming paper \cite{foootech3}.
\end{rem}
\begin{rem}\label{rem:3141}
In the geometric situation, the coordinate change where
$T' - T$ is positive at $T=\infty$ actually occurs.
So the coordinate change of the form (\ref{312312})
should be considered.
In fact, the parameter $T$ corresponds to the `length of the neck' region.
Namely if the neck of the source curve is $[0,1] \times [-5T,5T]$,
the corresponding element in the (thickened) moduli space has the coordinate
$T$. However the value of $T$ depends on the choice of the coordinate at the
nodes. Actually when nodes have coordinates $z$ and $w$,
we identify $zw = -r$ and
$$
r = e^{10\pi T}.
$$
See \cite[Section 8]{foooexp} right above Figure 12.
If we take a different choice of $z$, say $z' = e z$, then
$r$ becomes $r' =  er$. So $T' = T+1$.
\end{rem}
\begin{rem}\label{rem:geoKuraadmi}
(1) In \cite[Section 8]{foooexp} it is also proved that the
coordinate change of the obstruction bundle is admissible
and the Kuranishi map is also admissible.
(The former is a consequence of \cite[Proposition 8.19]{foooexp}
and the latter is proved in the course of the proof of
\cite[Proposition 8.31]{foooexp}.)
Therefore it turns out to be that the Kuranishi structure we constructed
on the moduli space of pseudo-holomorphic curves
is admissible in the sense of Definition \ref{def:adKura}.
\par
(2) The argument of \cite[Section 8]{foooexp} quoted above
is the cases of the Kuranishi chart of a stable map
$((\Sigma,\vec z),u)$ when the marked source curve
$(\Sigma,\vec z)$ is stable. There are cases when
the pair of a marked source curve
$(\Sigma,\vec z)$ and a map $u : \Sigma \to X$ is stable
but $(\Sigma,\vec z)$ is not stable.
We can also prove admissibility of the Kuranishi chart
in such cases. See \cite[Part 4 especially Section 21]{foootech}.
\end{rem}
\begin{rem}
There are different kind of boundaries or corners appearing in
applications. For example, to prove independence of
the filtered $A_{\infty}$ structure associated to a Lagrangian
submanifold under the change of compatible almost complex
structures chosen in the course of the construction,
we take a one parameter family of almost complex structures
$\{ J_{s}\}$ joining two almost complex structures $J_0$ and $J_1$
which we choose for the construction.
Then we consider the union of moduli spaces of pseudo-holomorphic discs
bounding our Lagrangian submanifold with respect to the almost complex structures
$J_{s}$ for $s \in [0,1]$.
In this case the part $s=0,1$ becomes a boundary.
To prove that our Kuranishi structure is admissible at this boundary,
we choose our family
$\{J_s \}$ so that $J_s = J_0$ (resp. $J_s = J_1$) for $s \in [0,\epsilon]$
(resp. for $s \in [1-\epsilon,1]$).
Then the boundary corresponding to $s=0,1$ has a canonical collar.
Therefore admissibility is obviously satisfied for this collar.
\par
We can study the case when we consider a homotopy of Hamiltonians
or homotopy of homotopies of almost complex structures
(or Hamiltonians) in the same way.
The study of such boundaries is much easier than that of the
boundary corresponding to the boundary node.
\end{rem}

\section{Stratified submersion to a manifold with corners.}
\label{subsec:ssmaptocorners}
In \cite[Definition 3.39]{part11} we defined the notions of strongly smooth map and
weekly submersive map from a K-space to a manifold {\it without boundary or corners}.
In this section we give the corresponding definitions
for the case when the target manifold
$P$ has boundary or corner.
In Sections \ref{sec:systemline1} and \ref{sec:systemline3},
we used them to define and study homotopy and/or higher
homotopy of morphisms of linear K-systems.
In Sections \ref{sec:systemtree1} and \ref{sec:systemtree2},
we also used them to define and study pseudo-isotopy of
filtered $A_{\infty}$ structure associated to a Lagrangian
submanifold.
\par
Let $P$ be a manifold with corners (cornered manifold).
For $p \in P$ we take a coordinate $(V_p,\phi_p)$
such that
$
V_p = \overline V_p \times [0,1)^k
$
where $\overline V_p$ is an open set of $\R^{\dim P - k}$
and $\phi_p : V_p \to P$ is a parametrization.
Here $p \in \overset{\circ}S_k(P)$ and
$p = \phi(y_p,(0,\dots,0))$.
In this section we take the coordinate of $P$ in this form.
\par
Let $X$ be an orbifold with corners (cornered orbifold),
and $f : X \to P$ a continuous map.

\begin{defn}\label{def301}
Under the situation above,
we say $f : X \to P$ is a
{\it corner stratified smooth map}
\index{corner ! corner stratified smooth map}
if for each $q \in X$ and $p = f(q)$ we can
choose coordinates $\frak V_q =(V_q, \Gamma_q, \phi_q)$ and
$\frak V_p=(V_p, \phi_p)$ respectively
with the following properties.
(Note that since $P$ is a smooth manifold, $\Gamma_p =\{ {\rm id} \}$.)
\begin{enumerate}
\item
$V_p = \overline V_p \times [0,1)^{k}$ is as above and
$V_q = \overline V_q \times [0,1)^{\ell+k}$.
Here $p \in \overset{\circ}S_k(P)$,
$q \in \overset{\circ}S_{k+\ell}(X)$ and
$q = \phi(y_q,(0,\dots,0))$.
\item
There exists a map $f_q : V_q \to V_p$
of the form
$$
f_q(y;(s_1,\dots,s_{\ell},t_1,\dots,t_k))
=
(\overline f_q(y;(s_1,\dots,s_{\ell},t_1,\dots,t_k)),(t_1,\dots,t_k))
$$
such that $\overline f_q : V_q \to \overline V_p$ is admissible.
\item
The following diagram commutes.
\begin{equation}\label{diag3011}
\begin{CD}
V_q @ >{\phi_q}>>
X  \\
@ V{f_p}VV @ VV{f}V\\
V_p @ > {\phi_p} >> P
\end{CD}
\end{equation}
\end{enumerate}
See Figure \ref{Figure25-1}.
In the case $X, P$ are admissible, we require
$\frak V_p$, $\frak V_q$ are admissible charts.
\end{defn}
\begin{figure}[h]
\centering
\includegraphics[scale=0.3]{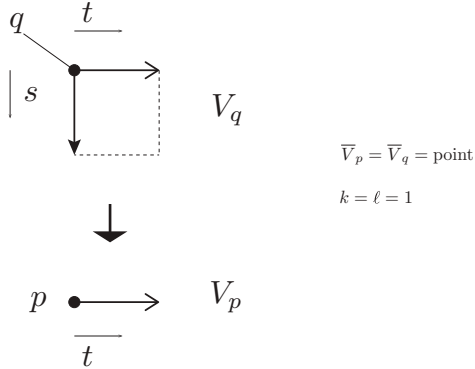}
\caption{Figure of Definition \ref{def301}}
\label{Figure25-1}
\end{figure}
\begin{rem}
Throughout this section, we can work either in the category
of smooth manifold (or orbifold) with corners,
or in the admissible category.
We do not mention about admissibility from now on in this
section.
\end{rem}
\begin{defn}
In the situation of Definition \ref{def301},
we say $f: X \to P$ is a
{\it corner stratified submersion}
\index{corner ! corner stratified submersion (from orbifold)}
\index{submersion ! corner stratified submersion (from orbifold)}
if $\overline f_q : V_q \to \overline V_p$ is a submersion for any $q \in X$.
\end{defn}

\begin{lem}\label{lem32333}
Let $X_1, X_2$ be cornered orbifolds and let
$P_1,P_2,P$ be cornered manifolds, and
$R$ a smooth manifold without boundary.
\begin{enumerate}
\item
Let $f_i : X_i \to P \times R$
be smooth maps.
Suppose $f_1$ is a corner stratified submersion.
Then the fiber product $X_1 \times_{P\times R} X_2$ carries a
structure of cornered orbifold.
In addition, if $\pi_P \circ f_2 :  X_2 \to P$ is
also a corner stratified submersion, then the map
$X_1 \times_P X_2 \to P$ induced from $f_1$ and $f_2$
in an obvious way is a corner stratified submersion.
\item
Let $f_i : X_i \to P_i \times R$
be smooth maps.
Suppose $f_1$
is a corner stratified submersion.
Then the fiber product $X_1 \times_{R} X_2$ carries a
structure of cornered orbifold.
In addition, if $\pi_{P_2} \circ f_2 :  X_2 \to P_2$
is also a corner stratified submersion, then
the map
$X_1 \times_R X_2 \to P_1 \times P_2$ induced from $f_1$ and $f_2$
in an obvious way is a corner stratified submersion.
\end{enumerate}
\end{lem}
The proof is obvious.

\begin{lem}\label{lem304}
Let $f : X \to P$ be a corner stratified smooth map from
a cornered orbifold to a cornered manifold.
Let $\widehat S_k(P)$ be the normalized corner of $P$
and $\pi : \widehat S_k(P) \to S_k(P) \subset P$ the projection.
\begin{enumerate}
\item
The fiber product
$
\widehat S_k(P) \times_P X
$
as topological space carries a structure of cornered orbifold.
The projection
$\widehat S_k(P) \times_P X \to \widehat S_k(P)$
is a corner stratified smooth map.
\item
The projection $\widehat S_k(P) \times_P X \to \widehat S_k(P)$
is a corner stratified submersion if
$f : X \to P$
is a corner stratified submersion.
\item
The map
$
\widehat S_{\ell}(\widehat S_k(P)) \times_P X
\to \widehat S_{\ell+k}(P) \times_P X
$
is a $(k+\ell)!/k!\ell!$ fold covering map.
\end{enumerate}
\end{lem}
The proof is again obvious.
\par
Now it is straightforward to generalize the story for an orbifold $X$ to the case when $X$
is a K-space.
\begin{defn}\label{defn3055}
Let $(X,\widehat{\mathcal U})$ be a K-space and $P$  a
manifold with corner.
\begin{enumerate}
\item A strongly continuous map $\widehat f : (X,\widehat{\mathcal U})
\to P$ is said to be a
{\it corner stratified smooth map}
\index{corner ! corner stratified smooth map (for K-space)} if
$f_p : U_p \to P$ is a corner stratified smooth map for any $p \in X$.
\item
A corner stratified smooth map $\widehat f : (X,\widehat{\mathcal U})
\to P$ is said to be a {\it corner stratified weak submersion}
\index{submersion ! corner stratified weak submersion}
\index{corner ! corner stratified weak submersion}
if $f_p : U_p \to P$ is a corner stratified submersion for any $p \in X$.
\item
Let $\widehat{\frak S}$ be a CF-perturbation of $X$.
We say that a corner stratified smooth map $\widehat f : (X,\widehat{\mathcal U}) \to P$ is
a {\it corner stratified strong submersion with respect to}
\index{strongly submersive ! corner stratified strong submersion}
\index{submersion ! corner stratified strong submersion}
\index{corner ! corner stratified strong submersion}
$\widehat{\frak S}$
if the following holds.
Let $p \in X$ and $(U_p,E_p,\psi_p, s_p)$ be a Kuranishi chart at $p \in X$.
Let $(V_{\frak r},\Gamma_{\frak r},\phi_{\frak r})$
be an orbifold chart of $U_p$ at some point and
$(W_{\frak r},\omega_{\frak r},\frak s_{\frak r}^{\epsilon})$
be a representative of $\widehat{\frak S}$ in this orbifold chart.
Then
$$
f \circ \psi_p \circ \phi_{\frak r} \circ {\rm pr} :
(\frak s_{\frak r}^{\epsilon})^{-1}(0) \to P
$$
is a corner stratified submersion.
Here  ${\rm pr} : V_{\frak r} \times W_{\frak r} \to V_{\frak r}$
is the projection.
\end{enumerate}
We can define a corner stratified smooth map and a corner stratified weak submersion
from a space equipped with a good coordinate system in the same way.
\end{defn}
\begin{lemdef}
Let $P$ be a cornered manifold and $R$ a manifold with boundary.
Let $\widehat f : (X,\widehat{\mathcal U}) \to P \times R$
be a corner stratified strong submersion with respect to $\widehat{\frak S}$.
Then for any differential form $h$ on $(X,\widehat{\mathcal U})$ and
for each sufficiently small $\epsilon >0$,
we can define the push out
$$
\widehat{f}!(h;\widehat{\frak S}^{\epsilon})
$$
which is a smooth differential form on $P \times R$,
in the same way as in \cite[Theorem 9.14]{part11}
using a good coordinate system compatible with the given
Kuranishi structure.
It is independent of the choice of the compatible good coordinate system if
$\epsilon >0$ is sufficiently small.
\end{lemdef}
The proof is the same as
that of \cite[Theorem 9.14]{part11} and so is omitted.
\begin{lem}
For each $i=1,2$ let $(X_i,\widehat{\mathcal U_i})$ be a K-space and
$\widehat f_i : (X_i,\widehat{\mathcal U_i})\to P \times R$ a corner stratified smooth map,
where $P$ is a cornered manifold and $R$ is a manifold
without boundary.
We assume that $\widehat f_1$ is a corner stratified weak submersion
and $\pi\circ \widehat f_2 : (X_2,\widehat{\mathcal U_2})\to P$
is a corner stratified weak submersion.
\begin{enumerate}
\item
The fiber product $X_1 \times_{P\times R} X_2$ carries a Kuranishi
structure and the map
$X_1 \times_{P\times R} X_2 \to P$ induced from $\widehat f_1$ and $\widehat f_2$
in an obvious way is a corner stratified weak submersion.
\item
Let $\widehat{\frak S_i}$ be a CF-perturbation
of $(X_i,\widehat{\mathcal U_i})$.
We assume that $\widehat f_1$ is a corner stratified strong submersion
with respect to $\widehat{\frak S_1}$
and $\pi_P\circ \widehat f_2 : (X_2,\widehat{\mathcal U_2})\to P$
is a corner stratified strong submersion
with respect to $\widehat{\frak S_2}$.
Then we can define the fiber product
$\widehat{\frak S_1} \times_{P \times R} \widehat{\frak S_2}$
of CF-perturbations.
The map
$X_1 \times_{P\times R} X_2
\to P$ is a corner stratified strong submersion
with respect to $\widehat{\frak S_1} \times_{P \times R} \widehat{\frak S_2}$.
\end{enumerate}
\end{lem}
\begin{proof}
This follows from Lemma \ref{lem32333} (1).
\end{proof}
\par
There is a slightly different situation we take the fiber product as follows:
\begin{lem}\label{lem32777}
For each $i=1,2$ let $(X_i,\widehat{\mathcal U_i})$ be a K-space and
$\widehat f_i : (X_i,\widehat{\mathcal U_i})\to R \times P_i$
a corner stratified smooth map, where $P_i$ is a cornered manifold and $R$ is a manifold
without boundary.
We assume that $\widehat f_1$ is a corner stratified weak submersion
and $\pi_{P_2} \circ \widehat f_2$ is a corner stratified weak submersion.
\begin{enumerate}
\item
The fiber product $X_1 \times_{R} X_2$ carries a Kuranishi
structure and the map
$X_1 \times_{R} X_2 \to P_1 \times P_2$ induced from $\widehat f_1$ and $\widehat f_2$
in an obvious way is a corner stratified weak submersion.
\item
Let $\widehat{\frak S_i}$ be a CF-perturbation
of $(X_i,\widehat{\mathcal U_i})$.
We assume that $\widehat f_i$ is a corner stratified strong submersion
with respect to $\widehat{\frak S_i}$.
Then we can define the fiber product
$\widehat{\frak S_1} \times_{R} \widehat{\frak S_2}$ of
CF-perturbations. The map
$X_1 \times_{R} X_2
\to P_1 \times P_2$ is a corner stratified strong submersion
with respect to $\widehat{\frak S_1} \times_{R} \widehat{\frak S_2}$.
\end{enumerate}
\end{lem}
\begin{proof}
This follows from Lemma \ref{lem32333} (2).
\end{proof}
\par
Next we discuss Stokes' formula and the composition formula under the correspondence.
\begin{defn}\label{defn3288}
Let $\widehat f : (X,\widehat{\mathcal U})\to P$
be a corner stratified weak submersion.
We divide the boundary of $X$ into two components:
\begin{equation}\label{3222foirm}
\partial_{\frak C^v} (X,\widehat{\mathcal U})
= \widehat f^{-1}(\partial P)
\end{equation}
and
\begin{equation}\label{3223foirm}
\partial_{\frak C^h} (X,\widehat{\mathcal U})
= \partial(X,\widehat{\mathcal U}) \setminus
\partial_{\frak C^v} (X,\widehat{\mathcal U}).
\end{equation}
We call (\ref{3222foirm}) the {\it vertical boundary}
\index{boundary ! vertical boundary}
and (\ref{3223foirm}) the {\it horizontal boundary}.
\index{boundary ! horizontal boundary}
See Figure \ref{Figure25-2}.
They are induced by the decomposition of the boundary
satisfying (\ref{eq16888tugi}) in Situation \ref{decomporbbdrkura}.
\end{defn}
\begin{figure}[h]
\centering
\includegraphics[scale=0.3]{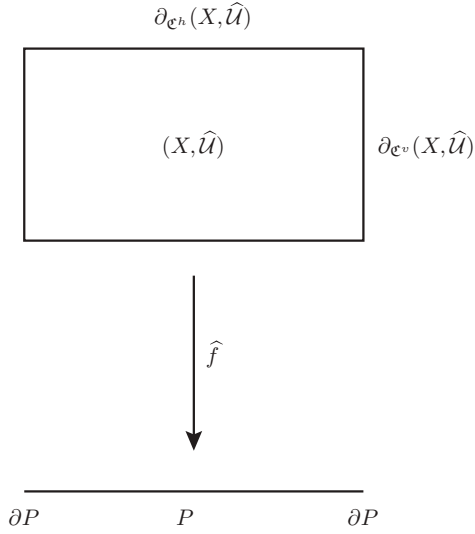}
\caption{vertical/horizontal boundary}
\label{Figure25-2}
\end{figure}
\begin{lem}\label{lem2511}
In the situation of Definition \ref{defn3288}, the restriction of $\widehat f$
to the horizontal boundary induces a corner stratified weak submersion
$\widehat f\vert_{\partial_{\frak C^h} (X,\widehat{\mathcal U})}
: \partial_{\frak C^h} (X,\widehat{\mathcal U}) \to P$.
If $\widehat f$ is a corner stratified strong submersion with
respect to a CF-perturbation $\widehat{\frak S}$, then
the restriction $\widehat f\vert_{\partial_{\frak C^h} (X,\widehat{\mathcal U})}$
is a corner stratified strong submersion with respect to
$\widehat{\frak S}\vert_{\partial_{\frak C^h} (X,\widehat{\mathcal U})}$.
\end{lem}
The proof is obvious.
\begin{thm}\label{prop3210}
In the  situation of Lemma \ref{lem2511},
let $h$ be a differential form on $(X,\widehat{\mathcal U})$.
Then for each sufficiently small $\epsilon >0$ we have
\begin{equation}\label{form342}
d \widehat f!(h;\widehat{\frak S}^{\epsilon})
=
\widehat f!(dh;\widehat{\frak S}^{\epsilon})
+
\widehat f!(h\vert_{\partial_{\frak C^h} (X,\widehat{\mathcal U})};
\widehat{\frak S}^{\epsilon}\vert_{\partial_{\frak C^h} (X,\widehat{\mathcal U})}).
\end{equation}
\end{thm}
\begin{proof}
Since both hand sides are smooth forms, it suffices to prove the formula
pointwise on ${\rm Int} \,P$. Let $p \in {\rm Int} \,P$
and take a compact set $K \subset {\rm Int} \,P$ containing an open neighborhood of $p$.
Using a partition of unity, we may assume without loss of generality
that $h$ is supported in $\widehat f^{-1}(K)$.
Then we can apply \cite[Theorem 9.26]{part11} to
$X \setminus \partial_{\frak C^v} (X,\widehat{\mathcal U})$
to prove the equality (\ref{form342}) at $p$.
\end{proof}
\begin{defn}
Let $P$ be a manifold with corner and let $M_s, M_t$ be manifolds
without boundary.
\begin{enumerate}
\item
A {\it $P$-parametrized
smooth correspondence}\index{smooth correspondence ! $P$-parametrized}
is a system
$$\frak X = \left((X,\widehat{\mathcal U}), \widehat f_s, \widehat f_t, \pi_{s,P}, \pi_{t,P}
\right)
$$
where $(X,\widehat{\mathcal U})$ is a K-space
and $\widehat f_s : (X,\widehat{\mathcal U})
\to P \times M_s $ and $\widehat f_t : (X,\widehat{\mathcal U})
\to P \times M_t$ are strongly smooth maps.
We assume that $\widehat f_t$ is a corner stratified weak submersion
and satisfies
\begin{equation}\label{form32555}
\pi_{s,P} \circ \widehat f_s = \pi_{t,P} \circ \widehat f_t.
\end{equation}
Here $\pi_{s,P}$ and $\pi_{t,P}$ are the projections to the $P$ factor.
\item
A {\it $P$-parametrized perturbed smooth correspondence}\index{smooth correspondence ! $P$-parametrized perturbed} is
$(\frak X,\widehat{\frak S})$ where  $\frak X = ((X,\widehat{\mathcal U}),
\widehat f_s,\widehat f_t)$ is a $P$-parametrized
smooth correspondence and $\widehat{\frak S}^{\epsilon}$ is a
CF-perturbation
of $(X,\widehat{\mathcal U})$ such that
$\widehat f_t$ is a corner stratified strong submersion with respect to
$\widehat{\frak S}$.
\item
Let  $(\frak X,\widehat{\frak S})$
($\frak X = ((X,\widehat{\mathcal U}),
\widehat f_s,\widehat f_t, \pi_{s,P}, \pi_{t,P})$) be a
$P$-parametrized perturbed smooth correspondence.
The restrictions of $\widehat f_s, \widehat f_t$, $\widehat{\frak S}$ to
the horizontal boundary
$\partial_{\frak C^h}(X,\widehat{\mathcal U})$
define a $P$-parametrized perturbed smooth correspondence.
We call it the {\it boundary} of $(\frak X,\widehat{\frak S})$
and denote it by $\partial(\frak X,\widehat{\frak S})$.
\end{enumerate}
\end{defn}
\begin{defn}\label{defn3213}
A $P$ parametrized perturbed smooth correspondence
$(\frak X,\widehat{\frak S})$ from $M_s$ to $M_t$ induces a map
$$
{\rm Corr}_{(\frak X,\widehat{\frak S}^{\epsilon})} : \Omega^k(P\times M_s) \to
\Omega^{k+\ell}(P\times M_t)
$$
by
\begin{equation}
{\rm Corr}_{(\frak X,\widehat{\frak S}^{\epsilon})}(h)
= \widehat f_t!(\widehat f_s^*h;\widehat{\frak S}^{\epsilon})
\end{equation}
for each sufficiently small $\epsilon >0$.
\end{defn}
\begin{lem}
Suppose we are in the situation of Definition \ref{defn3213}.
If $\rho$ is a differential form on $P$,
then for each sufficiently small $\epsilon >0$ we have
$$
{\rm Corr}_{(\frak X,\widehat{\frak S}^{\epsilon})}(\rho \wedge h)
= \pm \rho \wedge  {\rm Corr}_{(\frak X,\widehat{\frak S}^{\epsilon})}(h).
$$
\end{lem}
\begin{proof}
This is a consequence of (\ref{form32555}).
\end{proof}
Theorem \ref{prop3210} immediately implies the following.
\begin{prop}
For each sufficiently small $\epsilon >0$ we have
$$
d \circ {\rm Corr}_{(\frak X,\widehat{\frak S}^{\epsilon})} = {\rm Corr}_{(\frak X,\widehat{\frak S}^{\epsilon})} \circ d
=
{\rm Corr}_{\partial (\frak X,\widehat{\frak S}^{\epsilon})}.
$$
\end{prop}
Next we discuss the composition formula.
\begin{defnlem}\label{defnlem3216}
Let $M_1,M_2,M_3$ be smooth manifolds and $P, P_1, P_2$
manifolds with corner.
\begin{enumerate}
\item
Let $\frak X_{i i+1} = (X_{i i+1},\widehat{\mathcal U}_{i i+1},
\widehat f_{ii+1;i},\widehat f_{ii+1;i+1})$ be
$P$-parametrized smooth correspondences from $M_i$ to $M_{i+1}$
for $i=1,2$.
\begin{enumerate}
\item
The {\it composition} $\frak X_{13} = \frak X_{23} \circ \frak X_{12}$
is
$$
\left((X_{12},\widehat{\mathcal U}_{12}) \times_{P\times M_2}
(X_{23},\widehat{\mathcal U}_{23}),\widehat f_{12;1}\circ \pi,\widehat f_{12,2}\circ \pi
\right),
$$
which is a $P$-parametrized smooth correspondence from $M_1$ to $M_{3}$.
\item
In addition, if $(\frak X_{i i+1},\widehat{\frak S}_{i i+1})$
is a $P$-parametrized perturbed smooth correspondence,
then together with
$\widehat{\frak S}_{23} =
\widehat{\frak S}_{12} \times_{P\times M_2} \widehat{\frak S}_{23}$,
the composition
$\frak X_{13} = \frak X_{23} \circ \frak X_{12}$
defines a $P$-parametrized perturbed smooth correspondence
from $M_1$ to $M_3$.
We say $(\frak X_{13},\widehat{\frak S}_{13})$ is the {\it composition} of
$(\frak X_{12},\widehat{\frak S}_{12})$ and
$(\frak X_{23},\widehat{\frak S}_{23})$.
\end{enumerate}
\item
Let $\Xi_{i i+1} = (X_{i i+1},\widehat{\mathcal U}_{i i+1})$ be
$P_i$-parametrized smooth correspondences from $M_i$ to $M_{i+1}$
for $i=1,2$.
\begin{enumerate}
\item
The {\it composition}  $\frak X_{13} = \frak X_{23} \circ \frak X_{12}$
is defined by $(X_{12},\widehat{\mathcal U}_{12}) \times_{M_2}
(X_{23},\widehat{\mathcal U}_{23})$
which is a $(P_1 \times P_2)$-parametrized smooth correspondence from $M_1$ to $M_{3}$.
\item
In addition, if $(\frak X_{i i+1},\widehat{\frak S}_{i i+1})$
is a $P_i$-parametrized perturbed smooth correspondence, then
together with
$\widehat{\frak S}_{23} =
\widehat{\frak S}_{12} \times_{M_2} \widehat{\frak S}_{23}$,
the composition
$\frak X_{13}$ defines a $(P_1 \times P_2)$-parametrized perturbed smooth correspondence
from $M_1$ to $M_3$.
We say $(\Xi_{13},\widehat{\frak S}_{13})$ is the {\it composition} of
$(\frak X_{12},\widehat{\frak S}_{12})$ and
$(\frak X_{23},\widehat{\frak S}_{23})$.
\end{enumerate}
\end{enumerate}
\end{defnlem}
\begin{proof}
This is a consequence of Lemma \ref{lem32777}.
\end{proof}
\begin{prop}\label{prop328}
In the situation of Definition-Lemma \ref{defnlem3216} (1)
we have
\begin{equation}
{\rm Corr}_{(\frak X_{23},\widehat{\frak S}_{23}^{\epsilon})}
\circ {\rm Corr}_{(\frak X_{12},\widehat{\frak S}_{12}^{\epsilon})}
=
{\rm Corr}_{(\frak X_{13},\widehat{\frak S}_{13}^{\epsilon})}
\end{equation}
for each sufficiently small $\epsilon >0$.
\end{prop}
We will discuss the situation of Definition-Lemma \ref{defnlem3216} (2)
later.
\begin{proof}
Let us consider the following situation.
\begin{shitu}\label{situ3210}
For $i=1,2$, let $(X_i,\widehat{\mathcal U_i})$ be K-spaces,
$P,P_i$ manifolds with corner, and $R$ a manifold
without boundary.
Let $\widehat{\frak S_i}$ be CF-perturbations
of $(X_i,\widehat{\mathcal U_i})$.
\begin{enumerate}
\item
$\widehat f_i : (X_i,\widehat{\mathcal U_i}) \to P \times R$
are corner stratified strongly smooth maps for $i=1,2$.
We assume that $\widehat f_1 : (X_1,\widehat{\mathcal U_1}) \to P \times R$
and $\pi_P \circ \widehat f_2 :
(X_2,\widehat{\mathcal U_2}) \to P$ are corner stratified
weakly submersive and corner stratified strongly
submersive with respect to $\widehat{\frak S_1}$, $\widehat{\frak S_2}$, respectively.
\item
$\widehat f_i : (X_i,\widehat{\mathcal U_i}) \to P_i \times R$
are corner stratified strongly smooth maps.
We assume $\widehat f_1$ and $\pi_{P_2} \circ \widehat f_2$
are corner stratified weakly
submersive and corner stratified strongly
submersive
with respect to $\widehat{\frak S_1}$, $\widehat{\frak S_2}$, respectively.
$\blacksquare$
\end{enumerate}
\end{shitu}
\begin{lem}
In Situation \ref{situ3210} (1) we consider differential
forms $h_i$ on $(X_i,\widehat{\mathcal U_i})$.
They define a differential form $h_1 \wedge h_2$ on the fiber product
$(X_1,\widehat{\mathcal U_1}) \times_{P\times R}
(X_2,\widehat{\mathcal U_2})$.
Then for each sufficiently small $\epsilon >0$ we have
\begin{equation}\label{form3255}
\aligned
&\int_{((X_1,\widehat{\mathcal U_1}) \times_{R\times P}
(X_2,\widehat{\mathcal U_2}),(\widehat{\frak S_1}
\times_{R\times P}\widehat{\frak S_2})^{\epsilon})}
h_1 \wedge h_2 \\
&=
\int_{((X_2,\widehat{\mathcal U_2}),\widehat{\frak S_2})}
\widehat f_2^* \widehat f_1!(h_1;(\widehat{\frak S_1})^{\epsilon})
\wedge h_2.
\endaligned
\end{equation}
\end{lem}
\begin{proof}
We can use Sublemma \ref{lem321113} below in place of \cite[Lemma 10.27]{part11}.
Then the proof is the same as that of \cite[Proposition 10.23]{part11}.
\end{proof}
\begin{sublem}\label{lem321113}
For $i=1,2$ let $N_i$, $P$ be smooth manifolds with corner, and
$f_i : N_i \to P \times M$ smooth maps,
and $h_i$ smooth differential forms on $N_i$ of compact support.
Suppose that $f_1$ is a corner stratified submersion.
Then we have
\begin{equation}\label{form821}
\int_{N_1\,{}_{f_1}\times_{f_2} N_2} h_1 \wedge h_2
= \pm
\int_{N_2} f_2^* (f_1!(h_1)) \wedge h_2.
\end{equation}
\end{sublem}
\begin{proof}
Using a partition of unity, it suffices to prove \eqref{form821}
when $P =\overline P \times [0,1)^{b}$,
$N_1 = P\times M \times \R^{a_1} \times [0,1)^{a_2}$ and
$f_i : N_i \to P \times M$ is the obvious projection.
We can prove this by Fubini's theorem in the same way as in
\cite[Lemma 10.27]{part11}
\end{proof}
\begin{lem}
In Situation \ref{situ3210} (2) we consider differential
forms $h_i$ on $(X_i,\widehat{\mathcal U_i})$ ($i=1,2$) and
the wedge product $h_1 \wedge h_2$ on the fiber product
$(X_1,\widehat{\mathcal U_1}) \times_{R}
(X_2,\widehat{\mathcal U_2})$.
Then for each sufficiently small $\epsilon >0$ we have
\begin{equation}\label{form3255}
\aligned
&\int_{((X_1,\widehat{\mathcal U_1}) \times_{R\times P}
(X_2,\widehat{\mathcal U_2}),(\widehat{\frak S_1}
\times_{R\times P}\widehat{\frak S_2})^{\epsilon})}
h_1 \wedge h_2 \\
&=
\int_{((X_2,\widehat{\mathcal U_2}),\widehat{\frak S_2})}
(\pi_{R} \circ \widehat f_2)^*(\pi_{R} \circ \widehat f_1)!(h_1;(\widehat{\frak S_1})^{\epsilon})
\wedge h_2.
\endaligned
\end{equation}
\end{lem}
\begin{proof}
This is a consequence of
\cite[Proposition 10.23]{part11}.
\end{proof}
The proof of Proposition \ref{prop328} is now complete.
\end{proof}
To discuss the situation of Situation \ref{situ3210} (2)
we slightly generalize the notion of correspondence.
\begin{defn}
Let
$(\frak X,\widehat{\frak S}) = ((X,\widehat{\mathcal U},
\widehat f_{s},\widehat f_{t}, \pi_{s,P_1}, \pi_{t,P_1}),\widehat{\frak S})$ be a $P_1$-parametrized
perturbed smooth correspondence from $M_s$ to $M_t$.
Let $P_2$ be a manifold with corner.
We regard
$$
\left(P_2 \times (X,\widehat{\mathcal U}),
{\rm id}_{P_2}\times \widehat f_{s},{\rm id}_{P_2}\times \widehat f_{t},
\pi_{P_2}\circ ({\rm id}_{P_2} \times\pi_{s,P_1}),
\pi_{P_2}\circ ({\rm id}_{P_2} \times\pi_{t,P_1})\right)
$$
as a $(P_2\times P_1)$-parametrized smooth correspondence from $M_s$ to
$M_t$.
Here $\pi_{P_2} : P_2 \times P_1 \to P_2$ is the projection.
Then this defines a map which we denote by
\begin{equation}\label{form32111}
{\rm Corr}_{(\frak X,\widehat{\frak S}^{
\epsilon}),P_2} ~:~
\Omega^k(P_2 \times P_1 \times M_s)
\to \Omega^{k+\ell}(P_2 \times P_1 \times M_t)
\end{equation}
for each sufficiently small $\epsilon >0$.
Here $\ell =\dim M_t - \dim (X,\widehat{\mathcal U})$.
Similarly we define
\begin{equation}\label{form321111}
{\rm Corr}_{(\frak X,\widehat{\frak S}^{
\epsilon}),P_1} ~:~
\Omega^k(P_1 \times P_2 \times M_s)
\to \Omega^{k+\ell}(P_1 \times P_2 \times M_t).
\end{equation}
We note that when we define these maps, we do not use the orientations on $P_1, P_2$.
So the order of factors in the direct product does not cause the sign problem.
Thus we may write as
$$
{\rm Corr}_{(\frak X,\widehat{\frak S}^{
\epsilon}),P_2} ~:~
\Omega^k(P_1 \times P_2 \times M_s)
\to \Omega^{k+\ell}(P_1 \times P_2 \times M_t).
$$
\end{defn}
\begin{prop}\label{prop3223}
In the situation of Definition-Lemma \ref{defnlem3216} (2) we have
\begin{equation}
{\rm Corr}_{(\frak X_{23},\widehat{\frak S}^{
\epsilon}_{23}),P_1} \circ
{\rm Corr}_{(\frak X_{12},\widehat{\frak S}^{
\epsilon}_{12}),P_2}
=
{\rm Corr}_{(\frak X_{13},\widehat{\frak S}^{
\epsilon}_{13})}
\end{equation}
for each sufficiently small $\epsilon >0$.
\end{prop}
\begin{proof}
\begin{lem}\label{lem3224}
Suppose we are in Situation \ref{situ3210} (2).
Let $h_i$ be differential forms on $(X_i,\widehat{\mathcal U}_i)$ and
$\rho_i$ differential forms on $P_i$ for $i=1,2$. Then we obtain
a differential form
$$
h_1 \wedge (\pi_{P_1} \circ f_1)^*\rho_1
\wedge
h_2  \wedge (\pi_{P_2} \circ f_2)^*\rho_2
$$
on  $(X_1,\widehat{\mathcal U}_1) \times_R (X_2,\widehat{\mathcal U}_2)$.
Moreover we have the following equality:
\begin{equation}\label{form3255revrev}
\aligned
&\int_{(X_1,\widehat{\mathcal U}_1) \times_R (X_2,\widehat{\mathcal U}_2)}
h_1 \wedge (\pi_{P_1} \circ f_1)^*\rho_1 \wedge h_2  \wedge
(\pi_{P_2} \circ f_2)^*\rho_2 \\
&=
\int_{((X_2,\widehat{\mathcal U_2}),\widehat{\frak S_2})}
(\pi_R \circ f_2)^*(\pi_R  \circ f_1)!(h_1\wedge (\pi_{P_1} \circ f_1)^*\rho_1;(\widehat{\frak S_1})^{\epsilon})
\\
&\qquad\qquad\qquad\qquad\qquad\qquad\qquad\qquad\qquad
\wedge h_2 \wedge (\pi_{P_2} \circ f_2)^*\rho_2.
\endaligned
\end{equation}
\end{lem}
\begin{proof}
Applying \cite[Proposition 10.23]{part11} to
$h_1 \wedge (\pi_{P_1} \circ f_1)^*\rho_1$ and
$ h_2  \wedge
(\pi_{P_2} \circ f_2)^*\rho_2$, we obtain Lemma \ref{lem3224}.
\end{proof}
Proposition \ref{prop3223} is a consequence of Lemma \ref{lem3224}.
\end{proof}
Next we rewrite Lemma \ref{lem304} to the case when $X$ is a K-space.
\begin{lem}
Let $\widehat f : (X,\widehat{\mathcal U}) \to P$ be a corner stratified smooth map from
a K-space to a cornered manifold.
Let $\widehat S_k(P)$ be the normalized corner of $P$
and $\pi : \widehat S_k(P) \to S_k(P) \subset P$ the projection.
\begin{enumerate}
\item
The fiber product
$
\widehat S_k(P) \times_P X
$
as topological space carries a Kuranishi structure.
The projection
$\widehat S_k(P) \times_P X  \to \widehat S_k(P)$
is a corner stratified smooth map.
\item
The projection $\widehat S_k(P) \times_P X  \to \widehat S_k(P)$
is a corner stratified submersion if
$\widehat f : X \to P$
is a corner stratified submersion.
\item
The map
$
\widehat S_{\ell}(\widehat S_k(P)) \times_P X
\to \widehat S_{\ell+k}(P) \times_P X
$
is induced by a $(k+\ell)!/k!\ell!$ fold covering map of
K-spaces.
\end{enumerate}
\end{lem}
The proof is again obvious.
\par
Next we mention the relation to (partial) trivialization of the corners.
(See Sections \ref{sec:triboundary} and \ref{section:compomorphis}.)
The proof of the next lemma is straightforward so omitted.
\begin{lem}
Suppose $\widehat f : (X,\widehat {\mathcal U}) \to P$ is a corner stratified submersion from a
K-space to an (admissible) manifold with corner.
Then
$\widehat f$ induces a map
$$
\widehat f^{\boxplus\tau_0} : (X,\widehat {\mathcal U})^{\boxplus\tau_0} \to P^{\boxplus\tau_0}.
$$
Let $\frak C^{\rm vert}$ be a component of the corner of
$(X,\widehat {\mathcal U})$. Then we obtain a map
$$
\widehat f^{\frak C^{\rm vert}\boxplus\tau_0} :
 (X,\widehat {\mathcal U})^{\frak C^{\rm vert}\boxplus\tau_0} \to P.
$$
In both cases, if $\widehat f$ is corner stratified weakly submersive (resp.
corner stratified strongly submersive with respect to $\widehat{\frak S}$),
then $\widehat f^{\boxplus\tau_0}$ and $\widehat f^{\frak C^{\rm vert}\boxplus\tau_0}$
are corner stratified weakly submersive
(resp. corner stratified strongly submersive with respect to $\widehat{\frak S}$).
\end{lem}
Most of the stories of Kuranishi structure, CF-perturbation, push out etc. can be generalized to the case when the target space
has a corner, in a straightforward way. We will describe them when we need them.
\begin{rem}
In Section \ref{sec:systemline3} etc.
we are using corner stratified submersions (Definition \ref{defn3055})
to define and study homotopy and higher homotopy
of morphisms of linear K-systems.
On the other hand, we like to mention that there is another way to define homotopy and/or higher
homotopy of morphisms
etc. without using corner stratified submersions to manifolds
with corner.
Indeed, while we were writing
\cite{fooo06} we sometimes took this way.
For example, in \cite[Subsection 19.2]{fooo06}
we take a small number $\epsilon >0$ and
consider $P = (-\epsilon,1+\epsilon)$ instead of $P =[0,1]$.
When we construct $\mathcal N(\alpha_1,\alpha_2;P)$ in \cite{fooo06},
we consider the $P$-parametrized version of the moduli space so that
it is constant on $ (-\epsilon,0)$ and $(1,1+\epsilon)$.
This method also works rigorously.
However choosing $P =[0,1]$
and using corner stratified submersions seem more natural.
\end{rem}

\section{Local system and
smooth correspondence in de Rham theory with twisted coefficients}
\label{sec:lcalsystemtwist}
Let $\mathcal L$ be a local system, i.e., a flat vector bundle, on a manifold $M$.
We denote by $(\Omega^{\bullet}(M;{\mathcal L}), d=d_{\mathcal L})$ the de Rham complex with coefficients in $\mathcal L$.
We recall some basic operations on the de Rham complex with  twisted coefficients.
\par
\noindent
{\bf 1.} (Pull-back.)
Let $f:N \to M$ be a smooth map.  Clearly, the pull-back $f^*{\mathcal L}$ is a flat vector bundle and
we have a cochain homomorphisms
$$f^* : \Omega^{\bullet}(M;{\mathcal L}) \to \Omega^{\bullet}(N;f^*{\mathcal L}).$$
As in the usual de Rham theory, we have $d \circ f^*=f^* \circ d$.

\noindent
{\bf 2.} (Wedge product.)
Let ${\mathcal L}_1, {\mathcal L}_2$ be flat vector bundles.
Then ${\mathcal L}_1 \otimes {\mathcal L}_2$
is a flat vector bundle and we have the product
$$
\wedge : \Omega^{\bullet}(M;{\mathcal L}_1) \otimes
\Omega^{\bullet}(M;{\mathcal L}_2) \to
\Omega^{\bullet}(M;{\mathcal L}_1 \otimes {\mathcal L}_2).
$$
The wedge product and the differential enjoy the Leibniz' rule.

\noindent
{\bf 3.} (Integration along fibers.)
Let $\pi:N \to M$ be a proper submersion and let $O_M$ (resp. $O_N$) be the flat real line bundle associated with the orientation $O(1)$-bundle
of $M$ (resp. $N$).
We denote by $O_{\pi}=O_N \otimes \pi^*O_M$ be the relative orientation bundle of the submersion $\pi$.
Then we have the integration along fibers
$$\pi !:\Omega^{\bullet}(N;O_{\pi}) \to \Omega^{\bullet}(M).$$
For a flat vector bundle $\mathcal L$ on $M$, we have the integration along fibers with twisted coefficients
$$\pi !: \Omega^{\bullet}(N; \pi^*{\mathcal L}) \to \Omega^{\bullet}(M;{\mathcal L}).$$
Suppose that the boundary $\partial N$ of $N$ is not empty and the restriction $\pi\vert_{\partial N}: \partial N \to M$ is also submersion, then we have
$$d \circ \pi ! = \pi ! \circ d + (\pi\vert_{\partial N}) !.$$

Now consider the following situation.  Let $f_s:X \to M_s$ is a smooth map and $f_t:X \to M_t$ is a proper submersion.
\begin{equation}
\xymatrix{
&&  X \ar[ld]\ar[rd] \\
& M_s && M_t}
\nonumber
\end{equation}

Let ${\mathcal L}_s$ (resp. ${\mathcal L}_t$) be a flat vector bundle such that
$(f_s)^* {\mathcal L}_s \cong (f_t)^* {\mathcal L}_t \otimes O_{f_s}$.
By composing the pull-back operation by $f_s$ and the integration along fibers of $f_t$, we obtain the correspondence
$$f_t! \circ f_s^*:\Omega^{\bullet}(M_s;{\mathcal L}_s) \to \Omega^{\bullet}(M_t;{\mathcal L}_t).$$
Taking these arguments into account, we can obtain the following.
Now we consider \cite[Situation 7.1]{part11}.
Let $\mathfrak X=  ((X;\widehat{\mathcal U});\widehat f_s,\widehat f_t)$ be a smooth correspondence from $M_s$ to $M_t$.
Namely, $(X, \widehat{\mathcal U})$ is a K-space, $f_s:(X, \widehat{\mathcal U}) \to M_s$ is a weakly smooth map and
$f_t:(X,\widehat{\mathcal U}) \to M_t$ is a weak submersion.
\begin{equation}
\xymatrix{
&&  (X,\widehat{\mathcal U}) \ar[ld]\ar[rd] \\
& M_s && M_t}
\nonumber
\end{equation}
\begin{thm}\label{thm:261}
Let  $\widehat{\mathfrak S}= \{\widehat{\mathfrak S}^{\epsilon}\}$ be a CF-perturbation with respect to which $f_t$ is a strong submersion.
Let ${\mathcal L}_s$ (resp. ${\mathcal L}_t$) be a flat vector bundle on $M_s$ (resp. $M_t$) such that
\begin{equation}\label{eq:261}
(f_s)^*{\mathcal L}_s \cong (f_t)^*{\mathcal L}_t \otimes O_{f_t}.
\end{equation}
Here $O_{f_t}$ is the flat real line bundle associated with the relative orientation $O(1)$-bundle, i.e.,
\begin{equation}\label{eq:262}
O_{f_t} = (f_t)^*O_{M_t}^* \otimes O_{X}.
\end{equation}
Then
for $\widehat{\mathfrak X}=(\mathfrak X, \widehat{\mathfrak S})$ we have the map
$$
{\rm Corr}^{\epsilon}_{\widehat{\mathfrak X}} :
\Omega^{\bullet}(M_s; {\mathcal L}_s) \to
\Omega^{\bullet}(M_t;{\mathcal L}_t).
$$
\end{thm}
The following properties are fundamental.
\begin{thm}\label{thm:262}
If the restriction of $f_t$ to $(\partial X, \partial \widehat U)$ is strongly submersive with respect to $\widehat{\mathfrak S}\vert_{\partial X}$,
we have
$$d ~{\rm Corr}^{\epsilon}_{\widehat{\mathfrak X}} (h)
= {\rm Corr}^{\epsilon}_{\widehat{\mathfrak X}} (dh) + {\rm Corr}^{\epsilon}_{\partial\widehat{\mathfrak X}} (h)$$
for $h \in \Omega^{\bullet}(M_s;{\mathcal L}_s)$.
\end{thm}

In addition to \cite[Situation 10.15]{part11},
let ${\mathcal L}_i$ be a flat vector bundle on $M_i$, $i=1,2,3$ such that
$$
f_{1, 21}^* {\mathcal L}_1 \cong f_{2,21}^* {\mathcal L}_2 \otimes O_{f_{2,21}}, \quad  f_{2,32}^*{\mathcal L}_2 \cong f_{3,32}^* {\mathcal L}_3 \otimes O_{f_{3,32}}.
$$
Note that $O_{f_{3,31}} \cong g_{32,31}^*O_{f_{3,32}} \otimes g_{21,31}^*O_{f_{2,21}}$.
Here $g_{32,321}:\mathfrak X_{31} \to \mathfrak X_{32}$ and $g_{21,31}:\mathfrak X_{31} \to \mathfrak X_{21}$ are natural projections from the fiber product.
Hence we have
$$
f_{1,31}^* {\mathcal L}_1 \cong f_{3,31}^* {\mathcal L}_3 \otimes O_{f_{3,31}}.
$$
Then we have the following composition formula.
\begin{thm}\label{thm:263}
$$
{\rm Corr}^{\epsilon}_{\widehat{{\mathfrak X}_{32}} \circ \widehat{{\mathfrak X}_{21}}}
= {\rm Corr}^{\epsilon}_{\widehat{{\mathfrak X}_{32}}} \circ {\rm Corr}^{\epsilon}_{\widehat{{\mathfrak X}_{21}}}.
$$
\end{thm}

\newpage

\bibliographystyle{amsalpha}

\include{index}
\printindex
\end{document}